А.М. Гальмак

# n-АРНЫЕ ГРУППЫ

А.М. Гальмак

# n-АРНЫЕ ГРУППЫ

## ЧАСТЬ I

Приведены основные понятия и результаты из общей теории n-арных групп об n-арных подгруппах, порождающих множествах, смежных классах, различных обобщениях нормальных и сопряженных подгрупп в группе, а также о произведениях n-арных групп. Много внимания уделено изучению связи между n-арными подгруппами n-арной группы и подгруппами группы, к которой она приводима согласно теоремам Поста и Глускина-Хоссу. Рассмотрены вопросы, связанные со строением абелевых, полуабелевых, циклических и полуциклических n-арных групп.

Библиогр.: 47 назв.

# ВВЕДЕНИЕ

Начало развитию теории n-арных групп положила опубликованная в 1928 году в журнале "Mathematische Zeitschrift" статья В. Дёрнте "Untersuchungen über einen verallgemeinerten Gruppenbegriff" [1], в которой впервые было введено понятие n-группы, называемой также n-арной или полиадической группой. Уже из названия статьи видно, что истоки теории n-арных групп лежат в теории групп. Непосредственное отношение к возникновению новой теории имела также Эмми Нётер, по инициативе которой Дёрнте и занялся реализацией лежащей почти на поверхности идеи о замене в определении группы ассоциативной и однозначно обратимой слева и справа бинарной операции на ассоциативную и однозначно обратимую на каждом месте n-арную операцию. До Дёрнте такие тернарные, т. е. 3-арные операции, удовлетворяющие некоторым дополнительным условиям, изучал Х. Прюфер, опубликовавший в 1924 году в том же "Mathematische Zeitschrift" статью [2], в которой применял введенные им тернарные операции для исследования бесконечных абелевых групп. Впоследствии алгебры с такими операциями стали называть грудами Прюфера. Дёрнте установил, что груды Прюфера являются частным случаем n-арных групп, а именно, полуабелевыми тернарными группами, все элементы которых являются идемпотентами.

Первым, кто обратил серьезное внимание на статью В. Дёрнте, был, по-видимому, Э. Пост, сумевший разглядеть в небольшой статье зачатки многообещающей теории с широкими возможностями и блестящими перспективами. В 1940 году Э. Пост опубликовал в "Trans. Amer. Math. Soc."

объемную статью "Polyadic groups" [3], которая по важности полученных результатов и предложенных идей является одним из краеугольных камней теории n-арных групп и во многом предопределила тематику современных исследований по n-арным группам. Авторитет Поста имел решающее значение для привлечения свежих сил к изучению n-арных групп. Число изучающих n-арные группы стало постепенно возрастать, хотя тематика исследований расширялась незначительно, группируясь в основном вокруг аксиоматики n-арных групп и приводимости n-арных групп к группам. После Поста наибольший вклад в теорию n-арных групп внес С.А. Русаков, многочисленные результаты которого по n-арным группам, посвященные в основном силовскому строению n-арных групп и приложениям n-арных групп, систематизированы в двух его монографиях [4, 5]. Информация по n-арным группам имеется в книгах [6 – 9], а также в обзорах [10, 11].

К настоящему времени теория n-арных групп, несмотря на свой довольно почтенный возраст, остается для широкой математической общественности малоизвестной областью современной алгебры, значительно уступающей в своем развитии теории групп. Одной из основных причин сложившегося положения является широко распространенное заблуждение об отсутствии принципиальных различий между теорией групп и теорией n-арных групп при n ≥ 3. На самом деле это не так. В теории n-арных групп наряду со свойствами, общими для групп и n-арных групп, систематически изучаются и свойства n-арных групп, отсутствующие у групп. Изучение таких специфических свойств, среди которых встречаются и довольно экзотические, является одной из главных задач теории n-арных групп, причем не менее важной, чем получение n-арных аналогов известных групповых результатов.

Еще одной причиной замедленного развития теории n-арных групп является, на наш взгляд, очевидный дефицит учебной и монографической литературы по n-арным группам, без ликвидации которого невозможен сколько-нибудь значительный прогресс в изучении n-арных групп. В этой связи, кроме упомянутой выше статьи Поста [3], давно ставшей библиографической редкостью, и книг С.А. Русакова [4, 5], посвященных в основном его собственным результатам, можно указать еще книгу автора [12].

Предлагаемая книга пополнит небольшой список книг, посвященных n-арным группам. Она предназначена для первоначального знакомства с теорией n-арных групп и поэтому будет полезна в первую очередь тем, кто, несмотря на все предубеждения, осмелился ступить на давно открытый, но до сих пор малоизученный материк n-арных групп, остающийся, по сути дела, "terra incognita" на карте современной алгебры. Для чтения книги желательно знакомство с основами теории групп.

# Г Л А В А  1

# ТЕОРЕМЫ ПОСТА И ГЛУСКИНА-ХОССУ

Теоремы Поста и Глускина-Хоссу, играющие фундаментальную роль в теории полиадических групп, относятся к одному из важнейших её направлений, в рамках которого полиадические группы одной арности изучаются с помощью полиадических групп другой, в частности, меньшей арности. Значение теорем Поста и Глускина-Хоссу для теории полиадических групп заключается прежде всего в том, что они дают хороший инструмент для их изучения, позволяя сводить его к исследованию групп. По существу, теоремы Поста и Глускина-Хоссу позволяют вначале, погрузившись в теорию групп, воспользоваться её глубоко разработанным аппаратом, а затем подняться на поверхность теории полиадических групп уже с новыми для неё результатами.

## §1.1. КЛАССИЧЕСКИЕ ОПРЕДЕЛЕНИЯ n-АРНОЙ ГРУППЫ. ПРИМЕРЫ

Как уже отмечалось во введении, идея изучения ассоциативных и однозначно обратимых на каждом месте n-арных операций восходит к Э. Нётер. Реализуя эту идею, В. Дёрнте в 1928 году впервые ввёл [1] понятие n-группы, называемой также n-арной или полиадической группой. Статья В. Дёрнте привлекла внимание Э. Поста, опубликовавшего в 1940 году фундаментальную работу [3], в которой он привёл два новых определения n-арной группы, которые вместе с определением В. Дёрнте стали классическими. Эти определения являются обобщениями определения группы, как полугруппы, в которой разрешимы уравнения

$$xa = b, \qquad ay = b.$$

**1.1.1. Определение** [1, Дёрнте]. Универсальная алгебра $< A, [\ ] >$ с одной n-арной ($n \geq 2$) операцией $[\ ] : A^n \to A$ называется *n-арной группой*, если выполняются следующие условия:

1) n-арная операция $[\ ]$ на множестве A *ассоциативна*, т. е.
$$[[a_1...a_n]a_{n+1}...a_{2n-1}] = [a_1...a_i[a_{i+1}...a_{i+n}]a_{i+n+1}...a_{2n-1}]$$

для всех $i = 1, 2, ..., n - 1$ и всех $a_1, a_2, ..., a_{2n-1} \in A$;

2) каждое из уравнений

$$[a_1...a_{i-1} x_i a_{i+1}...a_n] = b, \; i = 1, 2, ..., n$$

однозначно разрешимо в A относительно $x_i$ для всех $a_1, ..., a_{i-1}, a_{i+1}, ..., a_n, b \in A$.

Полагая в определении 1.1.1 n = 2, получаем определение бинарной группы.

Если алгебра $< A, [\ ] >$ удовлетворяет условию 1) определения 1.1.1, то она называется *n-арной полугруппой*. Алгебра $< A, [\ ] >$, удовлетворяющая условию 2) того же определения, называется *n-арной квазигруппой*.

Пост заметил, что требование однозначной разрешимости уравнений в определении Дёрнте можно ослабить, потребовав только их разрешимость, а число уравнений уменьшить с n до двух, а при $n \geq 3$ даже до одного.

**1.1.2. Определение** [3, Пост]. n-Арная полугруппа $< A, [\ ] >$ называется n-арной группой, если в A разрешимы уравнения
$$[xa_2...a_n] = b, \quad [a_1...a_{n-1} y] = b$$

для всех $a_1, ..., a_n, b \in A$.

**1.1.3. Определение** [3, Пост]. n-Арная полугруппа $<A, [\;]>$ называется n-арной группой ($n \geq 3$), если в A разрешимо уравнение

$$[a_1...a_{i-1}\, x a_{i+1}...a_n] = b$$

для всех $a_1, ..., a_{i-1}, a_{i+1}, ..., a_n, b \in A$ и некоторого $i \in \{2, ..., n-1\}$.

В дальнейшем для сокращения записей будем использовать общепринятые в теории n-арных групп обозначения:

$$a_m^k = \begin{cases} a_m a_{m+1}...a_k, & m \leq k, \\ \varnothing, & m > k; \end{cases} \qquad \overset{k}{a} = \begin{cases} \underbrace{a...a}_{k}, & k > 0, \\ \varnothing, & k = 0. \end{cases}$$

Для всякого $m = k(n-1) + 1$, где $k \geq 1$, положим

$$[\,a_1^m\,] = [\,a_1^{k(n-1)+1}\,] = [[...[[\,a_1^n\,]\,a_{n+1}^{2n-1}\,]...]\,a_{(k-1)(n-1)+2}^{k(n-1)+1}\,].$$

Имеет место

**1.1.4. Теорема** [4, с. 9]. Пусть $<A, [\;]>$ – m-арная полугруппа, $m = k(n-1) + 1$, $r = t(n-1) + 1$, $1 \leq t \leq k$. Тогда

$$[\,a_1^m\,] = [\,a_1^j\,[\,a_{j+1}^{j+r}\,]\,a_{j+r+1}^m\,]$$

для всех $a_1, ..., a_m \in A$, где $j = 0, 1, ..., m - r$.

**1.1.5. Следствие.** Пусть $<A, [\;]>$ – универсальная алгебра, удовлетворяющая определению 1.1.1. Тогда для любого натурального k и любых $a_1, ..., a_{i-1}, a_{i+1}, ..., a_{k(n-1)}, a \in A$ уравнение

$$[\,a_1^{i-1}\, x_i\, a_{i+1}^{k(n-1)+1}\,] = a, \quad i \in \{1, 2, ..., k(n-1) + 1\}$$

разрешимо в A и его решение единственно.

**1.1.6. Пример.** Определим на группе A n-арную операцию

$$[a_1 a_2 ... a_n] = a_1 a_2 ... a_n a,$$

где а – элемент из центра Z(A) группы A. Так как

$$[[a_1 \ldots a_n]a_{n+1} \ldots a_{2n-1}] = (a_1 \ldots a_n a)a_{n+1} \ldots a_{2n-1}a =$$

$$= a_1 \ldots a_i(a_{i+1} \ldots a_{i+n}a)a_{i+n+1} \ldots a_{2n-1}a =$$

$$=[a_1 \ldots a_i[a_{i+1} \ldots a_{i+n}]a_{i+n+1} \ldots a_{2n-1}]$$

для всех $i = 1, 2, \ldots, n$ и всех $a_1, a_2, \ldots, a_{2n-1} \in A$, то $< A, [\ ] >$ – n-арная полугруппа.

Разрешимость в A уравнений

$$[x\, a_2 \ldots a_n] = b, \quad [a_1 \ldots a_{n-1}y] = b.$$

вытекает из разрешимости в A уравнений

$$xa_2 \ldots a_n a = b, \quad a_1 \ldots a_{n-1}ya = b.$$

Следовательно, $< A, [\ ] >$ – n-арная группа.

**1.1.7. Пример.** Положив в примере 1.1.6 $a = 1$, получим n-арную группу $< A, [\ ] >$ с n-арной операцией

$$[a_1 a_2 \ldots a_n] = a_1 a_2 \ldots a_n,$$

которая называется *производной* n-арной группой от группы A.

В предыдущих примерах мы строили n-арную групповую операцию при помощи групповой операции на множестве, которое относительно групповой операции было группой. Следующий пример показывает, что n-арную групповую операцию можно построить при помощи групповой операции на множестве, которое относительно групповой операции не является группой.

**1.1.8. Пример.** Определим на множестве $T_n$ всех нечетных подстановок степени n тернарную операцию $[\alpha\beta\gamma] = \alpha*\beta*\gamma$, где $*$ – умножение подстановок. Так как произведение трех нечетных подстановок является нечетной подстановкой, то множество $T_n$ замкнуто относительно тернарной операции $[\ ]$. Ассоциативность тернарной операции $[\ ]$ следует из ассоциативности бинарной операции $*$ в $S_n$. Ясно, что в $T_n$ однозначно разрешимы уравнения

$$[x\alpha_2\alpha_3] = \alpha, \quad [\alpha_1 y\alpha_3] = \alpha, \quad [\alpha_1\alpha_2 z] = \alpha.$$

Следовательно, $< T_n, [\ ] >$ – тернарная группа.

Пример 1.1.8 обобщается следующим предложением.

**1.1.9. Предложение.** Пусть B – подмножество группы A, удовлетворяющее следующим условиям:
  1) если $b_1, b_2, \ldots, b_n \in B$, то $b_1 b_2 \ldots b_n \in B$;
  2) если $b \in B$, то $b^{-1} \in B$.
Тогда $<B, [\ ]>$ – n-арная группа с n-арной операцией
$$[b_1 b_2 \ldots b_n] = b_1 b_2 \ldots b_n.$$

*Доказательство.* Если $b_1, b_2, \ldots, b_n \in B$, то
$$[b_1 b_2 \ldots b_n] = b_1 b_2 \ldots b_n \in B.$$

Ассоциативность n-арной операции $[\ ]$ является следствием ассоциативности операции в группе A.

Для произвольных $b_1, \ldots, b_{n-1}, b \in B$ положим
$$c = b\, b_{n-1}^{-1} \ldots b_1^{-1}, \qquad d = b_{n-1}^{-1} \ldots b_1^{-1}\, b.$$

Согласно 2), $b_1^{-1}, \ldots, b_{n-1}^{-1} \in B$, а согласно 1), $c, d \in B$. Так как
$$[cb_1 \ldots b_{n-1}] = cb_1 \ldots b_{n-1} = b\, b_{n-1}^{-1} \ldots b_1^{-1} b_1 \ldots b_{n-1} = b,$$
$$[b_1 \ldots b_{n-1} d] = b_1 \ldots b_{n-1} d = b_1 \ldots b_{n-1}\, b_{n-1}^{-1} \ldots b_1^{-1} b = b,$$

то в B разрешимы уравнения
$$[xb_1 \ldots b_{n-1}] = b, \qquad [b_1 \ldots b_{n-1} y] = b.$$

Следовательно, $<B, [\ ]>$ – n-арная группа. ∎

**1.1.10. Следствие.** Если B – подмножество группы A, удовлетворяющее условиям 1) и 2) предложения 1.1.9, и n – четное, то n-арная группа $<B, [\ ]>$ – производная от группы.

*Доказательство.* Если $n = 2k$, то согласно 2) предложения 1.1.9, $b^{-1} \in B$ для любого $b \in B$. Тогда согласно 1) того же предложения,
$$\underbrace{bb^{-1} \ldots bb^{-1}}_{k} = 1 \in B.$$

Кроме того, для любых $b_1, b_2 \in B$ имеем

$$b_1 b_2 = b_1 b_2 \underbrace{1\ldots 1}_{n-2} \in B.$$

Таким образом, B замкнуто относительно бинарной операции, содержит единицу и все свои обратные, то есть является группой. ∎

Так как $b = b^{-1}$ для всякой инволюции b группы A, то справедливо

**1.1.11. Следствие.** Если B – множество инволюций группы A, удовлетворяющее условию 1) предложения 1.1.9, то $< B, [\ ] >$ – n-арная группа с n-арной операцией из того же предложения.

**1.1.12. Пример.** Пусть b – элемент группы A, удовлетворяющий условию $b^{n-1} = 1$, [ ] – n-арная операция, производная от операции в группе. Так как $[\underbrace{b\ldots b}_{n}] = b^{n-1}b = b$, то $< \{b\}, [\ ] >$ – n-арная группа.

**1.1.13. Пример.** Если b – инволюция группы A, то есть $b^2 = 1$, то $< \{b\}, [\ ] >$ – тернарная группа с тернарной операцией, производной от операции в группе.

Если $< \{b\}, [\ ] >$ – n-арная группа из примера 1.1.12 и элемент b не является инволюцией, то $b^{-1} \notin \{b\}$, и поэтому обратное утверждение к предложению 1.1.9 в общем случае неверно.

Большое число примеров тернарных групп доставляют группы движений.

**1.1.14. Пример.** Любое движение плоскости является либо параллельным переносом $T_{\bar{a}}$ на некоторый вектор $\bar{a}$, либо поворотом $R_O^{\alpha}$ вокруг некоторой точки O на угол α, либо скользящим отражением $S_l^a = T_{\bar{a}} S_l = S_l T_{\bar{a}}$ относительно некоторого вектора $\bar{a}$ и некоторой прямой $l$, где $S_l$ – отражение относительно прямой $l$. Ясно, что $S_l^a = S_l$ при $\bar{a} = 0$. Параллельные переносы и повороты называют движениями первого рода, а скользящие отражения – движениями второго рода. Известно также, что: произведение двух движений пер-вого рода является движением первого рода; произведение движений первого и второго рода – движением второго рода; произведение двух движений второго рода – движением первого рода. Обозначим через $E_2(2)$ – множество всех движений второго рода плоскости и определим на $E_2(2)$ тернарную операцию $[s_1 s_2 s_3] = s_1 s_2 s_3$.

Используя приведенные свойства произведений движений и рассуждая так же, как в примере 1.1.8, можно показать, что $< E_2(2), [\ ] > -$ тернарная группа.

Аналогично устанавливается, что множество $E_2(3)$ всех движений второго рода пространства также является тернарной группой.

**1.1.15. Пример.** Всякий поворот прямой в некоторой плоскости на угол 180° вокруг любой точки этой прямой является инволюцией в группе всех самосовмещений прямой в выбранной плоскости. Кроме того, произведение трех таких поворотов снова является поворотом на 180°. Поэтому, согласно следствию 1.1.11, множество всех поворотов прямой на 180° в фиксированной плоскости является тернарной группой с тернарной операцией, производной от операции в группе.

**1.1.16. Пример.** Известно, что для любых трех прямых a, b, c, лежащих в одной плоскости и проходящих через точку O, существует прямая d, лежащая в той же плоскости и проходящая через точку O, и такая, что $S_a S_b S_c = S_d$. А так как, кроме того, всякое отражение вида $S_a$ является инволюцией в группе всех движений плоскости, то, согласно следствию 1.1.11, множество всех отражений вида $S_a$ относительно прямых, лежащих в выбранной плоскости и проходящих через общую точку, является тернарной группой с тернарной операцией, производной от операции в группе.

В следующем примере все рассматриваемые прямые также лежат в одной плоскости.

**1.1.17. Пример.** Известно, что если прямые a, b, c перпендикулярны прямой $l$, то существует прямая d, перпендикулярная $l$, и такая, что $S_a S_b S_c = S_d$. Таким образом, снова, согласно следствию 1.1.11, заключаем, что множество всех отражений вида $S_a$ относительно прямых, перпендикулярных одной прямой, является тернарной группой с тернарной операцией, производной от операции в группе.

**1.1.18. Пример.** Пусть $1 + (n-1)Z$ – класс вычетов по модулю $n-1$, где $n \geq 3$. Так как для любых

$$a_1 = 1+(n-1)z_1,\ a_2 = 1+(n-1)z_2, \ldots, a_n = 1+(n-1)z_n \in 1+(n-1)Z$$

верно

$$a_1+a_2+...+a_n = 1+(n-1)z_1+1+(n-1)z_2+...+1+(n-1)z_n =$$
$$= 1+n-1+(n-1)z_1+(n-1)z_2+...+(n-1)z_n =$$
$$= 1+(n-1)(1+z_1+z_2+...+z_n) \in 1+(n-1)Z,$$

то класс вычетов $1+(n-1)Z$ замкнут относительно n-арной операции

$$[a_1a_2...a_n] = a_1+a_2+...+a_n.$$

Ассоциативность n-арной операции [ ] является следствием ассоциативности операции в группе Z. Легко убедиться, что

$$x = y = 1+(n-1)(z_n-z_1-...-z_{n-1}-1) \in 1+(n-1)Z$$

является решением уравнений

$$[xa_1...a_{n-1}] = a_n, \quad [a_1...a_{n-1}y] = a_n.$$

Таким образом, $< 1+(n-1)Z, [\ ] > $ – n-арная группа.

n-Арная группа $< A, [\ ] >$ называется *конечной*, если множество A конечно. В этом случае число элементов $|A|$ множества A называется *порядком* n-арной группы $< A, [\ ] >$. Если A – бесконечное множество, то говорят, что она имеет *бесконечный порядок*.

n-Арные группы из примеров 1.1.8 и 1.1.12 являются конечными, причем их порядки равны соответственно n!/2 и 1. Примеры 1.1.14 – 1.1.18 доставляют примеры бесконечных n-арных групп.

## §1.2. АНАЛОГИ ЕДИНИЦЫ И ОБРАТНОГО ЭЛЕМЕНТА

Следующие три определения обобщают на n-арный случай определение единицы группы A, как элемента $e \in A$ такого, что $ea = ae = a$ для любого $a \in A$.

**1.2.1. Определение.** Элемент $e \in A$ n-арной группы $< A, [\ ] >$ называется *единицей* этой n-арной группы, если

$$[\underbrace{e...e}_{i-1}a\underbrace{e...e}_{n-i}] = a$$

для любого a ∈ A и любого i = 1, 2, ..., n.

**1.2.2. Определение.** Элемент $\varepsilon \in A$ n-арной группы < A, [ ] > называется её *идемпотентом* если

$$[\underbrace{\varepsilon...\varepsilon}_{n-1}a] = [a\underbrace{\varepsilon...\varepsilon}_{n-1}] = a$$

для любого a ∈ A.

Ясно, что единица n-арной группы является и ее идемпотентом.

**1.2.3. Определение.** Последовательность $e_1...e_{k(n-1)}$ ($k \geq 1$) элементов n-арной группы < A, [ ] > называется *нейтральной*, если

$$[e_1...e_{k(n-1)}a] = [ae_1...e_{k(n-1)}] = a$$

для любого a ∈ A.

Ясно, что если $\varepsilon$ – идемпотент, в частности, единица n-арной группы, то последовательности

$$\underbrace{\varepsilon...\varepsilon}_{n-1}, \underbrace{\varepsilon...\varepsilon}_{2(n-1)}, ..., \underbrace{\varepsilon...\varepsilon}_{k(n-1)}, ...$$

являются нейтральными.

**1.2.4. Пример.** Пусть < A, [ ] > – n-арная группа, производная от бинарной группы A (пример 1.1.7). Если e – единица группы, то для любого a ∈ A и любого i = 1, 2, ..., n верно

$$[\underbrace{e...e}_{i-1}a\underbrace{e...e}_{n-i}] = e^{i-1}a\, e^{n-i} = a,$$

то есть e – единица n-арной группы < A, [ ] >.

Имеет место и обратное утверждение [1].

**1.2.5. Предложение** [1]. n-Арная группа, обладающая единицей, является производной от группы.

В n-арной группе при n > 2, в отличие от групп, может быть несколько единиц. Более того существуют n-арные группы, в которых все элементы являются единицами.

Следующий пример нам понадобится для того, чтобы показать, что существуют n-арные группы (n > 2) любого конечного порядка, в которых вообще нет единиц.

**1.2.6. Пример.** Пусть $D_n$ – диэдральная группа, т. е. полная группа преобразований симметрии правильного n-угольника. Поворот с n-угольника в его плоскости на угол $2\pi/n$ вокруг центра n-угольника порождает циклическую подгруппу

$$C_n = <c> = \{e, c, c^2, \ldots, c^{n-1}\}$$

поворотов. Диэдральная группа содержит еще n-отражений. Если b – отражение, то
$$B_n = \{b, bc, \ldots, bc^{n-1}\} = \{b_1, b_2, \ldots, b_n\}$$
есть множество всех отражений.

Определим на $B_n$ тернарную операцию $[\varphi\psi\theta] = \varphi\psi\theta$. Так как произведение двух отражений является поворотом, то $\varphi\psi$ – поворот. А так как произведение поворота на отражение является отражением, то $[\varphi\psi\theta] = \varphi\psi\theta$ – отражение. Следовательно, множество $B_n$ замкнуто относительно тернарной операции [ ].

Ассоциативность тернарной операции [ ] вытекает из ассоциативности бинарной операции в диэдральной группе.

Рассмотрим в $B_n$ уравнение $[x\psi\theta] = \tau$, которое равносильно уравнению $x\psi\theta = \tau$. Последнее уравнение имеет в $D_n$ решение $x = \varphi$. Если $\varphi$ – поворот, то $\varphi\psi\theta = \tau$ – поворот, что противоречит выбору $\tau \in B_n$. Аналогично доказывается разрешимость в $B_n$ уравнений

$$[\varphi y\theta] = \tau, \qquad [\varphi\psi z] = \tau.$$

Мы показали, что $<B_n, [\ ]>$ – тернарная группа.

**1.2.7. Предложение.** В тернарной группе $<B_n, [\ ]>$ все элементы являются идемпотентами, среди которых нет единицы.

*Доказательство.* Так как произведение любого отражения на себя является тождественным преобразованием, то

$$[\varphi\varphi\psi] = \varphi\varphi\psi = \psi = \psi\varphi\varphi = [\psi\varphi\varphi]$$

для любых $\varphi, \psi \in B_n$, то есть все элементы в $< B_n, [\ ] >$ являются идемпотентами.

Для того, чтобы установить, что в $< B_n, [\ ] >$ нет единицы, покажем, что для любого $\varphi \in B_n$ $(n > 2)$ существует $\psi \in B_n$ такой, что $[\varphi\psi\varphi] \neq \psi$. Если $\varphi = b$, то, положив $\psi = bc$ и используя равенство $bc^i = c^{n-i}b$, получим

$$[\varphi\psi\varphi] = \varphi\psi\varphi = bbcb = cb = bc^{n-1} \neq bc = \psi.$$

Если же $\varphi = bc^i$ $(i = 1, 2, \ldots, n-1)$, то, положив $\psi = bc^{i+1}$, получим

$$[\varphi\psi\varphi] = \varphi\psi\varphi = bc^i bc^{i+1} bc^i = bc^i c^{n-i-1} bbc^i =$$
$$= bc^{n-1}c^i = bc^{i-1} \neq bc^{i+1} = \psi.$$

Следовательно, в $< B_n, [\ ] >$ нет единицы. ∎

Существуют n-арные группы, в которых нет не только единиц, но и идемпотентов.

**1.2.8. Пример.** Пусть $R = \{1, a, a^2, a^3, b, ba, ba^2, ba^3\}$ – группа кватернионов. На множестве R определим 5-арную операцию [ ] через бинарную операцию группы R следующим образом:

$$[x_1 x_2 \ldots x_5] = x_1 x_2 \ldots x_5 a^2.$$

Так как $Z(R) = \{1, a^2\}$ – центр группы R, то $< R, [\ ] >$ – 5-арная группа (см. пример 1.1.6).

Используя выполняющиеся в группе кватернионов тождества

$$a^4 = 1, \ a^2 = b^2, \ ab = ba^3,$$

можно показать, что 5-арная группа $< R, [\ ] >$ не обладает идемпотентами. Отметим, что 5-арная группа $< R, [\ ] >$ была построена С.А. Русаковым [4]. Он первым установил существование гамильтоновых n-арных групп без единицы, где $n = 4k + 1$.

Пример 1.2.8 можно обобщить.

**1.2.9. Пример.** Пусть A – группа экспоненты $n - 1$, $a \in Z(A)$, $a \neq 1$. В примере 1.1.6 установлено, что $< A, [\ ] >$ – n-арная группа с n-арной операцией

$$[a_1a_2...a_n] = a_1a_2...a_na.$$

Если x – произвольный элемент из A, то

$$[\underbrace{x...x}_{n}] = \underbrace{x...x}_{n}\,a = x^{n-1}xa = xa \neq x,$$

так как $a \neq 1$. Таким образом,

$$[\underbrace{x...x}_{n}] \neq x$$

для любого $x \in A$ и, следовательно, $< A, [\ ] >$ не содержит идемпотентов. Отметим, что в предыдущем примере $\operatorname{Exp} R = 4 = 5 - 1$.

**1.2.10. Предложение.** Элемент e n-арной группы $< A, [\ ] >$ является единицей, если

$$[a\underbrace{e...e}_{n-1}] = [ea\underbrace{e...e}_{n-2}] = [\underbrace{e...e}_{n-1}a] = a$$

для любого $a \in A$.

***Доказательство.*** Применив последовательно $n - 2$ раза равенства

$$[ea\underbrace{e...e}_{n-2}] = a, \quad [a\underbrace{e...e}_{n-1}] = a,$$

получим

$$[a\underbrace{e...e}_{n-1}] = [[ea\underbrace{e...e}_{n-2}]\underbrace{e...e}_{n-1}] =$$

$$= [e\,[a\underbrace{e...e}_{n-1}]\underbrace{e...e}_{n-2}] = [ea\underbrace{e...e}_{n-2}] =...= [\underbrace{e...e}_{n-2}ae],$$

то есть

$$[\underbrace{e...e}_{i-1}\,a\,\underbrace{e...e}_{n-i}] = a$$

для любого $a \in A$ и любого $i = 1, 2, …, n$. ∎

Доказанное предложение позволяет дать еще одно определение единицы n-арной группы, эквивалентное определению 1.2.1.

**1.2.11. Определение.** Элемент e n-арной группы $< A, [\ ] >$ называется ее *единицей*, если

$$[a\underbrace{e...e}_{n-1}] = [ea\underbrace{e...e}_{n-2}] = [\underbrace{e...e}_{n-1}a] = a$$

для любого $a \in A$.

Справедливо следующее

**1.2.12. Предложение.** Если $< A, [\ ] >$ – n-арная группа, $e_1, ..., e_{k(n-1)} \in A$, $k \geq 1$, то следующие утверждения эквивалентны:

1) последовательность $e_1...e_{k(n-1)}$ – нейтральная;
2) существует элемент $a \in A$ такой, что

$$[e_1...e_{k(n-1)}a] = a;$$

3) существует элемент $a \in A$ такой, что

$$[a\, e_1...e_{k(n-1)}] = a.$$

Предложение 1.2.12 позволяет дать еще одно определение идемпотента, эквивалентное определению 1.2.2.

**1.2.13. Определение.** Элемент $\varepsilon$ n-арной группы $< A, [\ ] >$ называется *идемпотентом*, если $[\underbrace{\varepsilon...\varepsilon}_{n}] = \varepsilon$.

В n-арной группе всегда разрешимо уравнение

$$[x\, e_2...e_{k(n-1)}a] = a.$$

Поэтому с учетом утверждения 2) предложения 1.2.12 справедливо

**1.2.14. Предложение.** В любой n-арной группе существуют нейтральные последовательности.

Существование в n-арной группе нейтральных последовательностей является также и следствием разрешимости в n-арной группе уравнения

$$[a\, e_1...e_{k(n-1)-1}y] = a$$

и утверждения 3) предложения 1.2.12.

Нейтральные последовательности n-арной группы определяются неоднозначно.

Иногда, для сокращения записей, последовательности элементов будем обозначать малыми греческими буквами: $a_1...a_i = \alpha$. При этом число $l(\alpha) = i$ будем называть длиной последовательности $\alpha$. Для пустой последовательности $\varnothing$ считают $l(\varnothing) = 0$.

**1.2.15. Предложение.** Если $\alpha\beta$ – нейтральная последовательность n-арной группы $< A, [\,] >$, то $\beta\alpha$ – также нейтральная последовательность $< A, [\,] >$.

**1.2.16. Предложение.** Если $\alpha$ и $\beta$ – нейтральные последовательности n-арной группы, то $\alpha\beta$ – также нейтральная последовательность.

Следующее определение обобщает на n-арный случай понятие обратного элемента группы.

**1.2.17. Определение.** Последовательность $\beta$ элементов n-арной группы $< A, [\,] >$ называется *обратной* к последовательности $\alpha$ элементов из A, если последовательности $\alpha\beta$ и $\beta\alpha$ являются нейтральными.

Ясно, что если $\beta$ – обратная к $\alpha$, то $\alpha$ – обратная к $\beta$.

**1.2.18. Предложение.** Для любой последовательности $\alpha$ элементов n-арной группы $< A, [\,] >$ существует обратная последовательность $\beta$.

Отметим, что обратная последовательность, длина которой больше единицы, определяется неоднозначно.

Следствием предложения 1.2.15 является

**1.2.19. Предложение.** Если $\alpha$ и $\beta$ – последовательности элементов n-арной группы, то следующие утверждения равносильны:

1) $\beta$ – обратная к $\alpha$;
2) $\alpha\beta$ – нейтральная;

3) $\beta\alpha$ – нейтральная.

**1.2.20. Предложение.** Пусть $\alpha_1, \ldots, \alpha_r$ – последовательности, составленные из элементов n-арной группы $<A, [\ ]>$, и пусть $\beta_1, \ldots, \beta_r$ – последовательности, обратные соответственно данным. Тогда $\beta_r \ldots \beta_1$ – обратная последовательность для последовательности $\alpha_1 \ldots \alpha_r$.

**1.2.21. Определение.** Элемент b n-арной группы $<A, [\ ]>$ называется *косым* элементом для элемента $a \in A$, если

$$[\underbrace{a \ldots a}_{i-1} b \underbrace{a \ldots a}_{n-i}] = a$$

для любого $i = 1, 2, \ldots, n$.

Если b – косой элемент для a, то употребляют обозначение $b = \bar{a}$. Таким образом, по определению имеем

$$[\underbrace{a \ldots a}_{i-1} \bar{a} \underbrace{a \ldots a}_{n-i}] = a.$$

Из определения 1.2.21 вытекает, что для всякого элемента a n-арной группы, косой элемент $\bar{a}$ определяется однозначно.

Следствием определения 1.2.21 и предложения 1.2.12 является следующее

**1.2.22. Предложение.** Для любого элемента a n-арной группы $<A, [\ ]>$ и любого $i = 1, 2, \ldots, n-1$, последовательность

$$\underbrace{a \ldots a}_{i-1} \bar{a} \underbrace{a \ldots a}_{n-i-1}$$

является нейтральной.

**1.2.23. Предложение.** Решение уравнения

$$[\underbrace{a \ldots a}_{i-1} x \underbrace{a \ldots a}_{n-i}] = a$$

для фиксированного i = 1, 2, …, n является косым элементом для a.

**Доказательство.** Пусть x = $\bar{a}$ – решение данного уравнения, то есть

$$[\underbrace{a\ldots a}_{i-1}\,\bar{a}\,\underbrace{a\ldots a}_{n-i}] = a.$$

Обе последовательности

$$\underbrace{a\ldots a}_{i-1},\quad \underbrace{a\ldots a}_{n-i}$$

не могут одновременно быть пустыми. Поэтому пусть для определенности $\underbrace{a\ldots a}_{n-i} \neq \varnothing$. Тогда, согласно предложению 1.2.12, последовательность

$$\underbrace{a\ldots a}_{i-1}\,\bar{a}\,\underbrace{a\ldots a}_{n-i-1}$$

является нейтральной, а по предложению 1.2.15 нейтральными будут и последовательности

$$\underbrace{a\ldots a}_{i-2}\,\bar{a}\,\underbrace{a\ldots a}_{n-i},\ \ldots,\ \bar{a}\,\underbrace{a\ldots a}_{n-2},$$

$$\underbrace{a\ldots a}_{n-2}\,\bar{a},\ \ldots,\ \underbrace{a\ldots a}_{i}\,\bar{a}\,\underbrace{a\ldots a}_{n-i-2},$$

то есть для любого i = 1, 2, …, n, последовательность

$$\underbrace{a\ldots a}_{i-1}\,\bar{a}\,\underbrace{a\ldots a}_{n-i-1}$$

является нейтральной. Следовательно,

$$[\underbrace{a\ldots a}_{i-1}\,\bar{a}\,\underbrace{a\ldots a}_{n-i}] = a$$

для любого i = 1, 2, …, n. ∎

**1.2.24. Пример.** Укажем косые элементы для каждого элемента 5-арной группы из примера 1.2.8:

$$\bar{1} = a^2,\ \bar{a} = a^3,\ \overline{a^2} = 1,\ \overline{a^3} = a,$$

$$\overline{b} = ba^2,\ \overline{ba} = ba^3,\ \overline{b\,a^2} = b,\ \overline{b\,a^3} = ba.$$

Следствием определений является

**1.2.25. Предложение.** Всякий идемпотент n-арной группы совпадает со своим косым.

**1.2.26. Предложение.** Если $< A, [\ ] > -$ n-арная группа $(n \geq 3)$, то

$$\overline{[a_1 \ldots a_n]} = [\underbrace{\overline{a}_n \overset{n-3}{a_n} \ldots \overline{a}_1 \overset{n-3}{a_1} \ldots \overline{a}_n \overset{n-3}{a_n} \ldots \overline{a}_1 \overset{n-3}{a_1}}_{n-2}]$$

для любых $a_1, \ldots, a_n \in A$. В частности, $\overline{[abc]} = [\,\overline{c}\,\overline{b}\,\overline{a}\,]$.

***Доказательство.*** Так как

$$n(n{-}2)(n{-}2) = n(n{-}2)(n{-}1{-}1) = n(n{-}2)(n{-}1) - n(n{-}2) =$$
$$= (n^2{-}2n)(n{-}1) - n(n{-}1{-}1) = (n^2{-}2n)(n{-}1) - n(n{-}1) + n =$$
$$= (n^2{-}3n)(n{-}1) + n - 1 + 1 = (n^2{-}3n{-}1)(n{-}1) + 1,$$

то правая часть равенства из условия леммы имеет смысл.

Кроме того,

$$[\underbrace{[a_1^n]\ldots[a_1^n]}_{n-2}\,[\underbrace{\overline{a}_n \overset{n-3}{a_n} \ldots \overline{a}_1 \overset{n-3}{a_1} \ldots \overline{a}_n \overset{n-3}{a_n} \ldots \overline{a}_1 \overset{n-3}{a_1}}_{n-2}]\,[a_1^n]] =$$

$$= [[\underbrace{[a_1^n]\ldots[a_1^n]}_{n-3}]\,a_1\ldots \underbrace{a_n\,\overline{a}_n\,\overset{n-3}{a_n}}_{\text{нейтр.}}\ldots \overline{a}_1\,\overset{n-3}{a_1}$$

$$\underbrace{\phantom{a_n\,\overline{a}_n\,a_n \ldots \overline{a}_1\,a_1}}_{\text{нейтр.}}$$

$$\underbrace{\overline{a}_n \overset{n-3}{a_n} \ldots \overline{a}_1 \overset{n-3}{a_1} \ldots \overline{a}_n \overset{n-3}{a_n} \ldots \overline{a}_1 \overset{n-3}{a_1}}_{n-3}\,[a_1^n]] =$$

$$= [[\underbrace{[a_1^n]\ldots[a_1^n]}_{n-3}\overline{a}_n \overset{n-3}{a_n}\ldots \overline{a}_1 \overset{n-3}{a_1}\ldots \overline{a}_n \overset{n-3}{a_n}\ldots \overline{a}_1 \overset{n-3}{a_1} [a_1^n]]] = \ldots = [a_1^n],$$

т. е.

$$[\underbrace{[a_1^n]\ldots[a_1^n]}_{n-2}\,[\underbrace{\overline{a}_n \overset{n-3}{a_n}\ldots \overline{a}_1 \overset{n-3}{a_1}\ldots \overline{a}_n \overset{n-3}{a_n}\ldots \overline{a}_1 \overset{n-3}{a_1}}_{n-2}]\,[a^n{}_1]] = [a_1^n].$$

Применяя теперь предложение 1.2.23, заключаем, что правая часть равенства из условия леммы действительно является косым элементом для элемента $[a_1\ldots a_n]$. ∎

**1.2.27. Предложение.** Если $<A, [\ ]>$ – n-арная группа, то для любого $a \in A$ верно

$$\overline{\overline{a}} = [\underbrace{a\ldots a}_{(n-3)(n-1)+1}].$$

*Доказательство.* Так как

$$[\underbrace{\overline{a}\ldots\overline{a}}_{n-1}[\underbrace{a\ldots a}_{(n-3)(n-1)+1}]] = [\underbrace{\overline{a}\ldots\overline{a}}_{n-2}[\overline{a}\underbrace{a\ldots a}_{n-1}]\underbrace{a\ldots a}_{(n-4)(n-1)+1}] =$$

$$= [\underbrace{\overline{a}\ldots\overline{a}}_{n-2}aa\underbrace{a\ldots a}_{(n-4)(n-1)+1}] = [\underbrace{\overline{a}\ldots\overline{a}}_{n-2}\underbrace{a\ldots a}_{(n-4)(n-1)}aa] = \ldots =$$

$$= [\overline{a}\,\overline{a}\underbrace{a\ldots a}_{n-4}aa] = \overline{a},$$

т. е.

$$[\underbrace{\overline{a}\ldots\overline{a}}_{n-1}[\underbrace{a\ldots a}_{(n-3)(n-1)+1}]] = \overline{a},$$

то по предложению 1.2.23 элемент

$$[\underbrace{a\ldots a}_{(n-3)(n-1)+1}]$$

является косым для $\overline{a}$. ∎

**1.2.28. Следствие.** Если $<A, [\ ]>$ – тернарная группа, то $\bar{\bar{a}} = a$.

## §1.3. ЭКВИВАЛЕНТНЫЕ ПОСЛЕДОВАТЕЛЬНОСТИ

Если $<A, \circ>$ – группа, то всякой последовательности $a_1 a_2 ... a_i$ ($i \geq 2$) элементов этой группы можно поставить в соответствие элемент $a_1 \circ a_2 \circ ... a_i$ этой же группы. В n-арной группе $<A, [\ ]>$ при $n > 2$ последовательности $a_1 a_2 ... a_i$ ($i \geq 2$) ее элементов подобным образом с помощью n-арной операции $[\ ]$ можно поставить в соответствие элемент из A только в том случае, если $i = k(n-1) + 1$, где $k \geq 1$:

$$a_1 a_2 ... a_{k(n-1)+1} \to [a_1 a_2 ... a_{k(n-1)+1}].$$

Если же $i \neq k(n-1) + 1$, то такое соответствие отсутствует. В связи с этим в теории n-арных групп вводится понятие эквивалентности последовательностей элементов n-арной группы, которое является отношением эквивалентности в обычном смысле, а при $n = 2$ совпадает с отношением равенства элементов группы.

**1.3.1. Определение.** Последовательности $\alpha$ и $\beta$ элементов n-арной группы $<A, [\ ]>$ называются *эквивалентными* в ней, если существуют последовательности $\gamma$ и $\delta$ элементов этой же n-арной группы такие, что

$$[\gamma \alpha \delta] = [\gamma \beta \delta]. \qquad (*)$$

Если последовательности $\alpha$ и $\beta$ элементов n-арной группы эквивалентны в ней, то для сокращения записей будем употреблять обозначение $\alpha \theta \beta$, указывая в необходимых случаях саму n-арную группу.

**1.3.2. Предложение.** Если $\alpha \theta \beta$, то $l(\alpha) \equiv l(\beta) \pmod{n-1}$.

*Доказательство.* Так как $\alpha \theta \beta$, то верно $(*)$, откуда

$$l(\gamma)+l(\alpha)+l(\delta) \equiv 1 \ (\mathrm{mod}\ n-1),$$
$$l(\gamma)+l(\beta)+l(\delta) \equiv 1 \ (\mathrm{mod}\ n-1).$$

Из этих двух сравнений получаем

$$l(\alpha)-l(\beta) \equiv 0 \ (\mathrm{mod}\ n-1),$$
$$l(\alpha) \equiv l(\beta) \ (\mathrm{mod}\ n-1). \blacksquare$$

Понятие эквивалентности последовательностей в n-арной группе было введено Э. Постом. Им же доказана следующая

**1.3.3. Теорема** [3]. Если $\alpha\theta\beta$ в n-арной группе $< A, [\ ] >$, то

$$[\nu\alpha\mu] = [\nu\beta\mu]$$

для любых последовательностей $\nu$ и $\mu$ элементов из A и таких, что

$$l(\nu) + l(\alpha) + l(\mu) \equiv 1 \ (\mathrm{mod}\ n-1).$$

*Доказательство.* Так как $\alpha\theta\beta$, то для них верно (∗). Пусть $\gamma$ и $\delta$ – непустые последовательности. По предложению 1.2.18 существуют последовательности $\gamma'$ и $\delta'$ – обратные соответственно к $\gamma$ и $\delta$, что влечет нейтральность последовательностей $\gamma'\gamma$ и $\delta\delta'$. Так как

$$l(\gamma) + l(\alpha) + l(\delta) \equiv 1 \ (\mathrm{mod}\ n-1),$$
$$l(\nu) + l(\alpha) + l(\mu) \equiv 1 \ (\mathrm{mod}\ n-1),$$

то

$$l(\nu) - l(\gamma) + l(\mu) - l(\delta) \equiv 0 \ (\mathrm{mod}\ n-1).$$

Кроме того,

$$l(\gamma') + l(\gamma) \equiv 0 \ (\mathrm{mod}\ n-1),$$
$$l(\delta) + l(\delta') \equiv 0 \ (\mathrm{mod}\ n-1).$$

Складывая почленно три последние сравнения, получим

$$l(\nu) + l(\gamma') + l(\delta') + l(\mu) \equiv 0 \ (\text{mod } n - 1).$$

Поэтому, приписывая к обеим частям равенства (∗) соответствующие последовательности, получим

$$[\nu\gamma'[\gamma\alpha\delta]\delta'\mu] = [\nu\gamma'[\gamma\beta\delta]\delta'\mu],$$

откуда, используя нейтральность последовательностей, получаем

$$[\nu\gamma'\gamma\alpha\delta\delta'\mu] = [\nu\gamma'\gamma\beta\delta\delta'\mu],$$

$$[\nu\alpha\mu] = [\nu\beta\mu].$$

Если последовательность $\gamma$ – пустая, то (∗) принимает вид $[\alpha\delta] = [\beta\delta]$, откуда

$$[\nu[\alpha\delta]\delta'\mu] = [\nu[\beta\delta]\delta'\mu],$$

$$[\nu\alpha\delta\delta'\mu] = [\nu\beta\delta\delta'\mu],$$

$$[\nu\alpha\mu] = [\nu\beta\mu].$$

Если же $\delta$ – пустая последовательность, то (∗) принимает вид $[\gamma\alpha] = [\gamma\beta]$, откуда

$$[\nu\gamma'[\gamma\alpha]\mu] = [\nu\gamma'[\gamma\beta]\mu],$$

$$[\nu\gamma'\gamma\alpha\mu] = [\nu\gamma'\gamma\beta\mu],$$

$$[\nu\alpha\mu] = [\nu\beta\mu].$$

Если $\gamma$ и $\delta$ – пустые последовательности, то (∗) принимает вид $[\alpha] = [\beta]$, откуда

$$[\nu[\alpha]\mu] = [\nu[\beta]\mu],$$

$$[\nu\alpha\mu] = [\nu\beta\mu]. \qquad \blacksquare$$

**1.3.4. Следствие.** Отношение $\theta$ является эквивалентностью на множестве всех последовательностей элементов n-арной группы.

*Доказательство.* Рефлексивность и симметричность являются простыми следствиями определения.

Если теперь $\alpha\theta\beta$, то существуют последовательности $\gamma$ и $\delta$ такие, что
$$[\gamma\alpha\delta] = [\gamma\beta\delta], \qquad (1)$$
откуда
$$l(\gamma) + l(\beta) + \alpha(\delta) \equiv 1 \pmod{n-1}.$$

Если же теперь $\beta\theta\tau$, то по теореме 1.3.3,
$$[\gamma\beta\delta] = [\gamma\tau\delta]. \qquad (2)$$

Из (1) и (2) получаем $[\gamma\alpha\delta] = [\gamma\tau\delta]$, т. е. $\alpha\theta\tau$ и значит отношение $\theta$ – транзитивно. ∎

Класс эквивалентности, содержащий последовательность $\alpha$, обозначим через $\theta(\alpha)$.

**1.3.5. Предложение.** Пусть $< A, [\ ] >$ – n-арная группа, $a, b \in A$. Тогда и только тогда $\theta(a\alpha) = \theta(b\alpha)$, когда $a = b$.

*Доказательство.* Достаточность очевидна.

*Необходимость.* Если $\theta(a\alpha) = \theta(b\alpha)$, то $a\alpha\theta b\alpha$, т. е. существуют последовательности $\gamma$ и $\delta$ такие, что
$$[\gamma a\alpha\delta] = [\gamma b\alpha\delta].$$

Из однозначной разрешимости в n-арной группе соответствующих уравнений вытекает $a = b$. ∎

Обозначим, как обычно, через $F_A$ свободную полугруппу над алфавитом $A$, т. е. множество всех последовательностей, составленных из элементов множества $A$ с бинарной операцией «склеивания» последовательностей.

**1.3.6. Определение.** Для всякой n-арной группы $< A, [\ ] >$ положим
$$\mathscr{A} = F_A / \theta.$$

Для всякого $i = 1, ..., n-1$ определим также множество

$$A^{(i)} = \{\theta(\alpha) \mid \theta(\alpha) \in \mathscr{A},\ l(\alpha) = i\},$$

в частности,

$$A' = \{\theta(a) \mid a \in A\}, \quad A'' = \{\theta(ab) \mid a, b \in A\}.$$

**1.3.7. Предложение.** Справедливы следующие утверждения:

1) $A^{(i)} = \{\theta(\alpha a) \mid a \in A\} = \{\theta(a\alpha) \mid a \in A\}$,

где $\alpha$ – фиксированная последовательность длины $i - 1$;

2) $|A^{(i)}| = |A|$ для любого $i = 1, ..., n - 1$;

3) $A^{(i)} \cap A^{(j)} = \varnothing$, где $i, j \in \{1, ..., n - 1\}, i \neq j$;

4) $\mathscr{A} = \bigcup_{i=1}^{n-1} A^{(i)}$;

5) $\theta$ – конгруэнция на $F_A$.

*Доказательство.* 1) Включение

$$\{\theta(\alpha a) \mid a \in A\} \subseteq A^{(i)}$$

очевидно.

Пусть теперь $\theta(a_1...a_i)$ – произвольный элемент из $A^{(i)}$. Зафиксируем элементы $a_{i+1}, ..., a_n \in A$ и рассмотрим в $A$ уравнение

$$[a_1...a_i a_{i+1}...a_n] = [\alpha x a_{i+1}...a_n],$$

которое имеет решение $x = a$, т. е.

$$[a_1...a_i a_{i+1}...a_n] = [\alpha a a_{i+1}...a_n].$$

Из последнего равенства вытекает эквивалентность последовательностей $a_1...a_i$ и $\alpha a$. Следовательно, $\theta(a_1...a_i) = \theta(\alpha a)$, откуда с учетом произвольного выбора $\theta(a_1...a_i)$, получаем

$$A^{(i)} \subseteq \{\theta(\alpha a) \mid a \in A\}.$$

Этим доказано равенство

$$A^{(i)} = \{\theta(\alpha a) \mid a \in A\}.$$

Равенство
$$A^{(i)} = \{\theta(a\alpha) \mid a \in A\}$$

доказывается аналогично.

2) Отображение $\varphi_i : A \to A^{(i)}$ по правилу

$$\varphi_i : a \mapsto \theta(\alpha a), \; i = 1, ..., n-1$$

является биекцией.

3) Очевидно.

4) Пусть $\theta(\alpha)$ – произвольный элемент из $\mathscr{A}$. Если $1 \leq l(\alpha) = j \leq n-1$, то

$$\theta(\alpha) \in A^{(j)} \subseteq \bigcup_{i=1}^{n-1} A^{(i)}.$$

Если же $l(\alpha) > n-1$, то $l(\alpha) = k(n-1)+j$, где $k \geq 1$, $1 \leq j \leq n-1$. Пусть для определенности $\alpha = a_1...a_{k(n-1)+j}$ и зафиксируем элементы $b_1, ..., b_{n-j} \in A$. Тогда

$$[a_1...a_{k(n-1)+j}b_1...b_{n-j}] = [a_1...a_{j-1}[a_j...a_{k(n-1)+j}]b_1...b_{n-j}],$$

что означает эквивалентность последовательностей

$$\alpha = a_1...a_{k(n-1)+j}, \quad a_1...a_{j-1}[a_j...a_{k(n-1)+j}].$$

Следовательно,

$$\theta(\alpha) = \theta(a_1...a_{j-1}[a_j...a_{k(n-1)+j}]) \in A^{(j)} \subseteq \bigcup_{i=1}^{n-1} A^{(i)}.$$

Так как класс $\theta(\alpha)$ выбран из $\mathscr{A}$ произвольно, то доказано включение

$$\mathscr{A} \subseteq \bigcup_{i=1}^{n-1} A^{(i)}.$$

Включение

$$\bigcup_{i=1}^{n-1} A^{(i)} \subseteq \mathscr{A}$$

очевидно. Из последних двух включений получаем равенство

$$\mathscr{A} = \bigcup_{i=1}^{n-1} A^{(i)}.$$

5) Если $\alpha\theta\beta$, $\alpha'\theta\beta'$, а $\gamma$ и $\delta$ такие, что

$$l(\gamma) + l(\alpha\alpha') + l(\delta) \equiv 1 \pmod{n-1},$$

то, применяя дважды теорему 1.3.3, получим

$$[\gamma\alpha\alpha'\delta] = [\gamma\beta\alpha'\delta] = [\gamma\beta\beta'\delta],$$

т. е.

$$[\gamma\alpha\alpha'\delta] = [\gamma\beta\beta'\delta],$$

откуда, согласно определению, $\alpha\alpha'\theta\beta\beta'$. ∎

Из 5) предыдущего предложения следует, что $\mathscr{A}$ – полугруппа, операцию которой будем обозначать через $*$. Во многих случаях для сокращения записей будем писать $\theta(\alpha)\theta(\beta)$ вместо $\theta(\alpha)*\theta(\beta)$.

**1.3.8. Предложение.** В $A^{(n-1)}$ существует класс, содержащий все нейтральные последовательности n-арной группы, и любая последовательность из которого является нейтральной.

*Доказательство.* Рассмотрим класс $\theta(e_1...e_{n-1}) \subseteq A^{(n-1)}$, где $e_1...e_{n-1}$ – нейтральная последовательность длины $n-1$ n-арной группы $<A, [\ ]>$. Если $\beta$ – другая нейтральная последовательность, то

$$[e_1...e_{n-1}a] = a = [\beta a],$$

т. е. последовательности $e_1...e_{n-1}$ и $\beta$ эквивалентны, и поэтому $\beta \in \theta(e_1...e_{n-1})$.

Если теперь $\gamma \in \theta(e_1...e_{n-1})$, то последовательности $\gamma$ и $e_1...e_{n-1}$ – эквивалентны. По теореме 1.3.3

$$[\gamma a] = [e_1...e_{n-1}a]$$

для любого $a \in A$. А так как $[e_1...e_{n-1}a] = a$, то $[\gamma a] = a$, что означает нейтральность последовательности $\gamma$. ∎

**1.3.9. Предложение.** Если β – обратная последовательность для последовательности α, то класс θ(β) содержит все обратные последовательности для последовательности α, и любая последовательность из θ(β) является обратной к α.

*Доказательство.* Последовательность αβ – нейтральная. Если β′ – еще одна обратная к α, то αβ′ – также нейтральная. Поэтому
$$[\alpha\beta a] = a, \ [\alpha\beta' a] = a$$
для некоторого a ∈ A, откуда
$$[\alpha\beta a] = [\alpha\beta' a],$$
что означает эквивалентность последовательностей β и β′. Следовательно, β′ ∈ θ(β).

Если теперь γ ∈ θ(β), то последовательности γ и β эквивалентны. По теореме 1.3.3 для любого a ∈ A верно
$$[\alpha\gamma a] = [\alpha\beta a].$$
А так как αβ – нейтральная последовательность, то [αβa] = a, откуда [αγa] = a, что означает нейтральность последовательности αγ. Следовательно, γ – обратная последовательность для α. ∎

## §1.4. ТЕОРЕМА ПОСТА О СМЕЖНЫХ КЛАССАХ

**1.4.1. Лемма.** Если ε – нейтральная последовательность из $F_A$, α – произвольная последовательность из $F_A$, то εαθαθαε.

*Доказательство.* Зафиксируем последовательность β ∈ $F_A$ такую, что
$$l(\alpha) + l(\beta) \equiv 1 \ (\text{mod } n-1).$$
Тогда, учитывая, что $l(\varepsilon) \equiv 0 \ (\text{mod } n-1)$, будем иметь

$$l(\varepsilon) + l(\alpha) + l(\beta) \equiv 1 \pmod{n-1}.$$

Таким образом, к последовательностям $\varepsilon\alpha\beta$ и $\alpha\beta$ применима n-арная операция [ ]. А так как $\varepsilon$ – нейтральная последовательность, то

$$[\varepsilon\alpha\beta] = [\varepsilon[\alpha\beta]] = [\alpha\beta],$$

откуда $\varepsilon\alpha\theta\alpha$. Аналогично доказывается $\alpha\theta\alpha\varepsilon$. ∎

**1.4.2. Теорема Поста о смежных классах [3].** Для всякой n-арной группы $<A, [\ ]>$ справедливы следующие утверждения:

1) $<\mathcal{A}, *>$ – группа;

2) $\theta([a_1...a_n]) = \theta(a_1)...\theta(a_n)$;

3) $\mathcal{A} = <A'>$;

4) $\theta(\beta)A^{(n-1)} = A^{(n-1)}\theta(\beta) = A^{(i)}$ для любого $\theta(\beta) \in A^{(i)}$, где $i \in \{1, ..., n-1\}$;

5) $A^{(n-1)}$ – инвариантная подгруппа группы $\mathcal{A}$;

6) $\mathcal{A}/A^{(n-1)} = \{A', A'', ..., A^{(n-1)}\} = <\{A'\}>$ – циклическая группа порядка $n-1$.

*Доказательство.* 1) Согласно 5) предложения 1.3.7, $<\mathcal{A}, *>$ – полугруппа. Обозначим через E – множество всех нейтральных последовательностей n-арной группы $<A, [\ ]>$, которое по предложению 1.3.8 совпадает с классом $\theta(e_1...e_{n-1}) \subseteq A^{(n-1)}$, где $e_1...e_{n-1}$ – нейтральная последовательность. Если $\theta(\alpha)$ – произвольный элемент из $\mathcal{A}$, то, применяя лемму 1.4.1, получим

$$E\theta(\alpha) = \theta(e_1...e_{n-1})\theta(\alpha) = \theta(e_1...e_{n-1}\alpha) = \theta(\alpha) =$$
$$= \theta(\alpha e_1...e_{n-1}) = \theta(\alpha)\theta(e_1...e_{n-1}) = \theta(\alpha)E,$$

т. е.

$$E\theta(\alpha) = \theta(\alpha) = \theta(\alpha)E.$$

Этим показано, что $\mathcal{A}$ – полугруппа с единицей $E = \theta(e_1...e_{n-1})$.

Пусть снова $\theta(\alpha)$ – произвольный элемент из $\mathcal{A}$ и $\beta$ – обратная последовательность для последовательности $\alpha$. Так как

$$\theta(\alpha)\theta(\beta) = \theta(\alpha\beta) = E = \theta(\beta\alpha) = \theta(\beta)\theta(\alpha),$$

то $\theta(\beta)$ – обратный элемент для $\theta(\alpha)$. Мы показали, что $\mathcal{A}$ – группа.

2) Элемент $[a_1...a_n]$, как последовательность длины 1 и последовательность $a_1...a_n$ эквивалентны. Поэтому

$$\theta([a_1...a_n]) = \theta(a_1...a_n),$$

а так как

$$\theta(a_1...a_n) = \theta(a_1)...\theta(a_n),$$

то

$$\theta([a_1...a_n]) = \theta(a_1)...\theta(a_n).$$

3) Пусть $\theta(\alpha)$ – произвольный элемент из $\mathcal{A}$. Согласно 4) предложения 1.3.7, можно считать $\alpha = a_1...a_i$, где $1 \le i \le n-1$. Тогда

$$\theta(\alpha) = \theta(a_1...a_i) = \theta(a_1)...\theta(a_i),$$

где $\theta(a_1), ..., \theta(a_i) \in A'$.

4) Можно считать $\beta = b_1...b_i$. Применяя 1) предложения 1.3.7, где $\alpha = a_1...a_{n-2}$, получим

$$\theta(\beta)A^{(n-1)} = \{\theta(\beta)\theta(a_1...a_{n-2}a) \mid a \in A\} =$$
$$= \{\theta(b_1...b_i a_1...a_{n-2}a) \mid a \in A\} =$$
$$= \{\theta(b_1...b_{i-1}[b_i a_1...a_{n-2}a]) \mid a \in A\} \subseteq A^{(i)},$$

т. е.

$$\theta(\beta)A^{(n-1)} \subseteq A^{(i)}.$$

Если теперь $\theta(\beta) = \theta(\delta b)$ – произвольный класс из $A^{(i)}$, то для любого другого класса $\theta(\alpha a) \in A^{(i)}$ существуют $a_1(a), ..., a_{n-1}(a) \in A$ такие, что

$$\alpha a \theta \delta[b a_1(a)...a_{n-1}(a)]\theta \delta b a_1(a)...a_{n-1}(a).$$

Тогда

$$A^{(i)} = \{\theta(\alpha a) \mid a \in A\} = \{\theta(\delta[ba_1(a)...a_{n-1}(a)]) \mid a \in A\} =$$
$$= \{\theta(\delta b a_1(a)...a_{n-1}(a)) \mid a \in A\} =$$
$$= \{\theta(\beta)\theta(a_1(a)...a_{n-1}(a)) \mid a \in A\} \subseteq \theta(\beta)A^{(n-1)},$$

т. е.

$$A^{(i)} \subseteq \theta(\beta)A^{(n-1)}.$$

Из обоих доказанных включений, получаем

$$\theta(\beta)A^{(n-1)} = A^{(i)}.$$

Второе равенство доказывается аналогично.

5) Так как для любых $\theta(\alpha a)$ и $\theta(\alpha b)$ из $A^{(n-1)}$ имеем

$$\theta(\alpha a)\theta(\alpha b) = \theta(\alpha a\alpha b) = \theta(\alpha[a\alpha b]) \in A^{(n-1)},$$

то групповая операция замкнута на $A^{(n-1)}$. Согласно предложению 1.3.8, единица $E = \theta(e_1...e_{n-1})$ принадлежит $A^{(n-1)}$. Обратный элемент для $\theta(\alpha a) \in A^{(n-1)}$ также принадлежит $A^{(n-1)}$. Таким образом, $A^{(n-1)}$ – подгруппа группы $\mathscr{A}$.

Инвариантность $A^{(n-1)}$ в $\mathscr{A}$ вытекает из предыдущего пункта и 4) предложения 1.3.7.

6) Из 4) с учетом предложения 1.3.7, получаем

$$\mathscr{A}/A^{(n-1)} = \{A', A'', ..., A^{(n-1)}\}.$$

Кроме того,

$$\underbrace{A'...A'}_{i} = \{\theta(a_1) \mid a_1 \in A\}...\{\theta(a_i) \mid a_i \in A\} =$$
$$= \{\theta(a_1)...\theta(a_i) \mid a_1, ..., a_i \in A\} =$$
$$= \{\theta(a_1...a_i) \mid a_1, ..., a_i \in A\} = A^{(i)},$$

т. е.

$$\underbrace{A'...A'}_{i} = A^{(i)}.$$

Следовательно, $\mathscr{A}/A^{(n-1)}$ – циклическая группа порядка $n-1$, порожденная элементом $A'$. ∎

В дальнейшем для сокращения записей подгруппу $A^{(n-1)}$ будем обозначать распространенным в литературе по n-арным группам символом $A_o$, т. е. $A_o = A^{(n-1)}$.

Изоморфизм n-арных групп определяется как изоморфизм универсальных алгебр.

**1.4.3. Следствие.** Всякая n-арная группа изоморфно вкладывается в n-арную группу, производную от группы.

*Доказательство.* Для всякой n-арной группы $< A, [\;] >$ определим на группе $< \mathcal{A}, * >$ n-арную операцию $\lfloor\;\rfloor$ по правилу

$$\lfloor \theta(\alpha_1)\theta(\alpha_2)...\theta(\alpha_n) \rfloor = \theta(\alpha_1)*\theta(\alpha_2)*...\theta(\alpha_n) =$$
$$= \theta(\alpha_1)\theta(\alpha_2)...\theta(\alpha_n).$$

Тогда $< \mathcal{A}, \lfloor\;\rfloor >$ – n-арная группа, производная от группы $< \mathcal{A}, * >$.

Так как $< \mathcal{A}, \lfloor\;\rfloor >$ – n-арная группа, то для любых $\theta(a)$, $\theta(a_1), ..., \theta(a_{n-1}) \in A' \subseteq \mathcal{A}$ в $\mathcal{A}$ разрешимы уравнения

$$\lfloor x\theta(a_1)...\theta(a_{n-1}) \rfloor = \theta(a),$$
$$\lfloor \theta(a_1)...\theta(a_{n-1})y \rfloor = \theta(a),$$

решения которых, очевидно, принадлежат $A'$. Поэтому $< A', \lfloor\;\rfloor >$ – n-арная подгруппа n-арной группы $< \mathcal{A}, \lfloor\;\rfloor >$.

Отображение $\varphi : a \mapsto \theta(a)$ устанавливает изоморфизм n-арной группы $< A, [\;] >$ на n-арную группу $< A', \lfloor\;\rfloor >$. ∎

В силу изоморфизма из следствия 1.4.3, n-арные группы $< A, [\;] >$ и $< A', \lfloor\;\rfloor >$ можно отождествлять, что позволяет сформулировать

**1.4.4. Следствие.** Всякая n-арная группа является смежным классом группы по нормальной подгруппе.

Из 2) предложения 1.3.7 и 6) теоремы 1.4.2. вытекает

**1.4.5. Следствие.** Если $<A, [\,]>$ – конечная n-арная группа, то $|\mathcal{A}| = |A|(n-1)$.

Если в n-арной группе $<A, [\,]>$ зафиксировать элементы $b_1, ..., b_{n-2} \in A$, то, согласно предложению 1.3.7, группу $\mathcal{A}$ можно представить в виде

$$\mathcal{A} = \bigcup_{j=1}^{n-1} \{\theta(ab_1...b_{j-1}) \mid a \in A\}.$$

Поэтому, если $<A, [\,]>$ – конечная n-арная группа, где $A = \{a_1,...,a_m\}$, то группу $\mathcal{A}$ можно наглядно представить в виде таблицы на рис.1.

| $\theta(a_1)$ | $\theta(a_1b_1)$ | $\theta(a_1b_1b_2)$ | ... | $\theta(a_1b_1...b_{j-1})$ | ... | $\theta(a_1b_1...b_{n-2})$ |
|---|---|---|---|---|---|---|
| $\theta(a_2)$ | $\theta(a_2b_1)$ | $\theta(a_2b_1b_2)$ | ... | $\theta(a_2b_1...b_{j-1})$ | ... | $\theta(a_2b_1...b_{n-2})$ |
| ... | ... | ... | | ... | | ... |
| $\theta(a_i)$ | $\theta(a_ib_1)$ | $\theta(a_ib_1b_2)$ | ... | $\theta(a_ib_1...b_{j-1})$ | ... | $\theta(a_2b_1...b_{n-2})$ |
| ... | ... | ... | | ... | | ... |
| $\theta(a_m)$ | $\theta(a_mb_1)$ | $\theta(a_mb_1b_2)$ | ... | $\theta(a_mb_1...b_{j-1})$ | ... | $\theta(a_mb_1...b_{n-2})$ |
| $A'$ | $A''$ | $A'''$ | | $A^{(j)}$ | | $A_o = A^{(n-1)}$ |

Рис. 1.

Группа $\mathcal{A}$ совпадает с множеством всех $m(n-1)$ клеток таблицы, причем все клетки j-го столбца образуют множество $A^{(j)}$. Первый столбец $A'$ порождает группу $\mathcal{A}$ и является n-арной группой, изоморфной n-арной группе $<A, [\,]>$. Элементы последнего столбца образуют нормальную подгруппу $A_o$ группы $\mathcal{A}$, факторгруппа $\mathcal{A}/A_o$ – это все столбцы $A'$, $A''$, ..., $A^{(n-1)} = A_o$, причем первый столбец является образующим циклической группы $\mathcal{A}/A_o$.

Если в n-арной группе $<A, [\,]>$ зафиксировать элемент b, то группу $\mathcal{A}$ можно представить в виде

$$\mathcal{A} = \bigcup_{j=1}^{n-1} \{\theta(a\underbrace{b...b}_{j-1}) \mid a \in A\} = \bigcup_{j=1}^{n-1} \{\theta(\underbrace{b...b}_{j-1}a) \mid a \in A\}.$$

Поэтому, если $< A, [\ ] >$ – конечная n-арная группа, где $A = \{a_1, ..., a_m\}$, то группу $\mathcal{A}$ можно наглядно представить в виде таблицы на рис.2.

| $\theta(a_1)$ | $\theta(a_1b)$ | $\theta(a_1bb)$ | ... | $\theta(a_1\underbrace{b...b}_{j-1})$ | ... | $\theta(a_1\underbrace{b...b}_{n-2})$ |
|---|---|---|---|---|---|---|
| $\theta(a_2)$ | $\theta(a_2b)$ | $\theta(a_2bb)$ | ... | $\theta(a_2\underbrace{b...b}_{j-1})$ | ... | $\theta(a_2\underbrace{b...b}_{n-2})$ |
| ... | ... | ... |   | ... |   | ... |
| $\theta(a_i)$ | $\theta(a_ib)$ | $\theta(a_ibb)$ | ... | $\theta(a_i\underbrace{b...b}_{j-1})$ | ... | $\theta(a_i\underbrace{b...b}_{n-2})$ |
| ... | ... | ... |   | ... |   | ... |
| $\theta(a_m)$ | $\theta(a_mb)$ | $\theta(a_mbb)$ | ... | $\theta(a_m\underbrace{b...b}_{j-1})$ | ... | $\theta(a_m\underbrace{b...b}_{n-2})$ |
| $A'$ | $A''$ | $A'''$ |   | $A^{(j)}$ |   | $A_o = A^{(n-1)}$ |

Рис. 2.

Согласно 4) теоремы 1.4.2,

$$A^{(j)} = \theta^j(b)A_o = A_o\theta^j(b),\ j = 1, ..., n - 1.$$

Учитывая также 3) и 4) предложения 1.3.7, получаем

**1.4.6. Предложение.**

$$\mathcal{A} = A_o + \theta(b)A_o + \theta^2(b)A_o + ... + \theta^{n-2}(b)A_o =$$
$$= A_o + A_o\theta(b) + A_o\theta^2(b) + ... + A_o\theta^{n-2}(b).$$

Приведем определение обертывающей группы, впервые введенное Постом.

**1.4.7. Определение** [3]**.** Группа $< \widetilde{A}, \bullet >$ называется *обертывающей* для n-арной группы $< A, [\ ] >$, если:
1) группа $< \widetilde{A}, \bullet >$ порождается множеством A;
2) $[x_1x_2...x_n] = x_1 \bullet x_2 \bullet ... \bullet x_n$ для любых $x_1, x_2, ..., x_n \in A$.

По аналогии с определением множеств $A^{(i)}$, определим для любого $i = 1, ..., n - 1$ множества

$$\widetilde{A}^{(i)} = \{a_1 \bullet ... \bullet a_i \mid a_1, ..., a_i \in A\}.$$

В частности, $\widetilde{A}' = A$, $\widetilde{A}^{(n-1)} = \widetilde{A}_o$.

Следующее предложение является аналогом предложения 1.3.7.

**1.4.8. Предложение.** Если $< \widetilde{A}, \bullet > $ – обертывающая группа n-арной группы $< A, [\ ] >$, то справедливы следующие утверждения:

1) $\widetilde{A}^{(i)} = \{a_1 \bullet ... \bullet a_{i-1} \bullet a \mid a \in A\} = \{a \bullet a_1 \bullet ... \bullet a_{i-1} \mid a \in A\}$,

где $a_1, ..., a_{i-1}$ – фиксированные элементы из $A$ ($i = 1, ..., n - 1$), причем

$$a_1 \bullet ... \bullet a_{i-1} \bullet a = a_1 \bullet ... \bullet a_{i-1} \bullet b$$

тогда и только тогда, когда $a = b$;

2) $|\widetilde{A}^{(i)}| = |A|$;

3) для любых $i, j \in \{1, ..., n - 1\}$ множества $\widetilde{A}^{(i)}$ и $\widetilde{A}^{(j)}$ либо не пересекаются, либо совпадают;

4) $\widetilde{A} = \bigcup_{i=1}^{n-1} \widetilde{A}^{(i)}$.

*Доказательство.* 1) Включение

$$\{a_1 \bullet ... \bullet a_{i-1} \bullet a \mid a \in A\} \subseteq \widetilde{A}^{(i)}$$

очевидно.

Пусть теперь $b_1 \bullet ... \bullet b_i$ – произвольный элемент из $\widetilde{A}^{(i)}$. Зафиксируем элементы $b_{i+1}, ..., b_n \in A$ и рассмотрим в $A$ уравнение

$$[b_1...b_i b_{i+1}...b_n] = [a_1...a_{i-1} x b_{i+1}...b_n],$$

которое имеет решение $x = a \in A$, т. е.

$$[b_1...b_i b_{i+1}...b_n] = [a_1...a_{i-1} a b_{i+1}...b_n].$$

Учитывая условие 2) определения 1.4.7, получаем

$$b_1 \bullet ... \bullet b_i \bullet b_{i+1} \bullet ... \bullet b_n = a_1 \bullet ... \bullet a_{i-1} \bullet a \bullet b_{i+1} \bullet ... \bullet b_n,$$

откуда

$$b_1 \bullet ... \bullet b_i = a_1 \bullet ... \bullet a_{i-1} \bullet a.$$

Следовательно,

$$\widetilde{A}^{(i)} \subseteq \{ a_1 \bullet ... \bullet a_{i-1} \bullet a \mid a \in A\},$$

и равенство

$$\widetilde{A}^{(i)} = \{ a_1 \bullet ... \bullet a_{i-1} \bullet a \mid a \in A\}$$

доказано. Второе равенство доказывается аналогично.

Если

$$a_1 \bullet ... \bullet a_{i-1} \bullet a = a_1 \bullet ... \bullet a_{i-1} \bullet b,$$

то $a = b$ в силу того, что $< \widetilde{A}, \bullet >$ – группа.

2) Отображение $f_i : A \to \widetilde{A}^{(i)}$ по правилу

$$f_i : a \mapsto a_1 \bullet ... \bullet a_{i-1} \bullet a, \; i = 1, ..., n - 1$$

является биекцией.

3) Предположим, что существуют

$$a_1 \bullet ... \bullet a_i \in \widetilde{A}^{(i)}, \; b_1 \bullet ... \bullet b_j \in \widetilde{A}^{(j)}$$

такие, что

$$a_1 \bullet ... \bullet a_i = b_1 \bullet ... \bullet b_j$$

Пусть для определенности $i < j$. Тогда, зафиксировав элементы $b_1, ..., b_{i-1} \in A$, получим согласно 1)

$$a_1 \bullet ... \bullet a_i = a \bullet b_1 \bullet ... \bullet b_{i-1}, \; a \in A.$$

Зафиксировав элементы $b_1, ..., b_{j-2}$, и снова применяя 1), получим

$$b_1 \bullet ... \bullet b_j = a \bullet b_1 \bullet ... \bullet b_{i-1} \bullet b_i \bullet ... \bullet b_{j-2} \bullet c, \; c \in A.$$

Таким образом,

$$a \bullet b_1 \bullet ... \bullet b_{i-1} = a \bullet b_1 \bullet ... \bullet b_{i-1} \bullet b_i \bullet ... \bullet b_{j-2} \bullet c,$$

т. е. $b_i \bullet ... \bullet b_{j-2} \bullet c = e$ – единица группы $< \widetilde{A}, \bullet >$.

Если теперь $c_1 \bullet ... \bullet c_i$ – произвольный элемент из $\widetilde{A}^{(i)}$, то

$$c_1 \bullet ... \bullet c_i = c_1 \bullet ... \bullet c_i \bullet e = c_1 \bullet ... \bullet c_i \bullet b_i \bullet ... \bullet b_{j-2} \bullet c \in \widetilde{A}^{(j)},$$

т. е.
$$\widetilde{A}^{(i)} \subseteq \widetilde{A}^{(j)}.$$

Произвольный элемент $u \in \widetilde{A}^{(j)}$ можно согласно 1) представить в виде
$$u = a' \bullet b_1 \bullet \ldots \bullet b_{i-1} \bullet b_i \bullet \ldots \bullet b_{i-2} \bullet c, \ a' \in A,$$
откуда
$$u = a' \bullet b_1 \bullet \ldots \bullet b_{i-1} \bullet e = a' \bullet b_1 \bullet \ldots \bullet b_{i-1} \in \widetilde{A}^{(i)},$$
т. е.
$$\widetilde{A}^{(j)} \subseteq \widetilde{A}^{(i)}.$$

Мы доказали равенство
$$\widetilde{A}^{(j)} = \widetilde{A}^{(i)}.$$

4) Включение
$$\bigcup_{i=1}^{n-1} \widetilde{A}^{(i)} \subseteq \widetilde{A}$$
очевидно.

Так как $\widetilde{A} = <A>$, то произвольный элемент $v \in \widetilde{A}$ имеет вид
$$v = a_1 \bullet \ldots \bullet a_t, \ t \in Z, \ a_1,\ldots,a_t \in A.$$

Если $t \leq n - 1$, то
$$v = a_1 \bullet \ldots \bullet a_t \in \bigcup_{i=1}^{n-1} \widetilde{A}^{(i)}.$$

Если же $t \geq n$, то $t = k(n - 1) + i$, где $k \geq 1$, $1 \leq i \leq n - 1$. Поэтому
$$v = a_1 \bullet \ldots \bullet a_t = a_1 \bullet \ldots \bullet a_{i-1} \bullet (a_i \bullet \ldots \bullet a_{k(n-1)+i}) =$$
$$= a_1 \bullet \ldots \bullet a_{i-1} \bullet [a_i \ldots a_{k(n-1)+i}] \in \widetilde{A}^{(i)},$$
откуда, в силу произвольного выбора $v \in \widetilde{A}$, получаем
$$\widetilde{A} \subseteq \bigcup_{i=1}^{n-1} \widetilde{A}^{(i)}$$

Таким образом, $\widetilde{A} = \bigcup_{i=1}^{n-1} \widetilde{A}^{(i)}.$ ∎

**1.4.9. Теорема.** Если $<\widetilde{A},\bullet>$ – обертывающая группа n-арной группы $<A,[\ ]>$, то справедливы следующие утверждения:

1) существует гомоморфизм $\psi$ группы $<\mathcal{A},*>$ на группу $<\widetilde{A},\bullet>$;

2) сужение гомоморфизма $\psi$ на $A_o$ является изоморфизмом групп $<A_o,*>$ и $<\widetilde{A}_o,\bullet>$, причем $\widetilde{A}_o = \psi(A_o)$ – инвариантная подгруппа группы $<\widetilde{A},\bullet>$;

3) существует гомоморфизм $\widetilde{\psi}$ факторгруппы $\mathcal{A}/A_o$ на факторгруппу $\widetilde{A}/\widetilde{A}_o$;

4) $\widetilde{A}/\widetilde{A}_o = <\{A\}>$ – циклическая группа с образующим элементом $A$, имеющая порядок, делящий $n-1$.

***Доказательство.*** 1) Определим отображение $\psi : \mathcal{A} \to \widetilde{A}$ по правилу

$$\psi : \theta(a_1...a_i) \mapsto a_1 \bullet ... \bullet a_i, \ i = 1, ..., n-1.$$

Ясно, что $\psi$ – сюръекция. Пусть теперь

$$\theta(a_1...a_i) \in A^{(i)}, \ \theta(a_{i+1}...a_{j+i}) \in A^{(j)}$$

произвольные элементы из $\mathcal{A}$, где $i, j \in \{1, ..., n-1\}$.
Если $i + j \leq n - 1$, то

$$\psi(\theta(a_1...a_i)\theta(a_{i+1}...a_{j+i})) = \psi(\theta(a_1...a_i a_{i+1}...a_{i+j})) =$$
$$= a_1 \bullet ... \bullet a_i \bullet a_{i+1} \bullet ... \bullet a_{i+j} = \psi(\theta(a_1...a_i))\psi(\theta(a_{i+1}...a_{i+j})).$$

Если же $i + j \geq n$, т. е. $i + j = k(n-1) + l$, где $k \geq 1$, $1 \leq l \leq n-1$, то

$$\psi(\theta(a_1...a_i)\theta(a_{i+1}...a_{j+i})) = \psi(\theta(a_1...a_i a_{i+1}...a_{i+j})) =$$
$$= \psi(\theta(a_1...a_{l-1}[a_l...a_{k(n-1)+l}]) =$$
$$= a_1 \bullet ... \bullet a_{l-1} \bullet [a_l...a_{k(n-1)+l}] = a_1 \bullet ... \bullet a_{l-1} \bullet a_l \bullet ... \bullet a_{k(n-1)+l} =$$
$$= a_1 \bullet ... \bullet a_i \bullet a_{i+1} \bullet ... \bullet a_{j+i} = \psi(\theta(a_1...a_i))\psi(\theta(a_{i+1}...a_{j+i})).$$

Мы показали, что $\psi$ гомоморфизм группы $\mathcal{A}$ на группу $<\widetilde{A},\bullet>$.

2) Ясно, что $\psi(A_o) = \widetilde{A}_o$. Кроме того, из

$$\psi(\theta(a_1...a_{n-1})) = \psi(\theta(b_1...b_{n-1})$$

следует

$$\theta(a_1...a_{n-1}) = \theta(b_1...b_{n-1}).$$

Следовательно, сужение $\psi$ на $A_o$ является изоморфизмом групп $< A_o, * >$ и $< \widetilde{A}_o, \bullet >$. Инвариантность $\widetilde{A}_o$ в $\widetilde{A}$ вытекает из инвариантности $A_o$ в $\mathcal{A}$.

3) Гомоморфизм $\widetilde{\psi}$ факторгруппы $\mathcal{A}/A_o$ на факторгруппу $\widetilde{A}/\widetilde{A}_o$ определяется по правилу

$$\widetilde{\psi} : A^{(i)} = \theta(a_1...a_i)A_o \mapsto \psi(\theta(a_1...a_i)) \bullet \widetilde{A}_o = a_1 \bullet ... \bullet a_i \bullet \widetilde{A}_o.$$

Ясно, что при этом $\widetilde{\psi}(A') = A$.

4) Вытекает из 3). ■

**1.4.10. Замечание.** Утверждение 2) теоремы 1.4.9 можно получить и непосредственно, доказав предварительно равенство

$$a_1 \bullet ... \bullet a_i \bullet \widetilde{A}_o = \widetilde{A}_o \bullet a_1 \bullet ... \bullet a_i = \widetilde{A}^{(i)}, i = 1, ..., n-1.$$

Отсюда же вытекает, что

$$\widetilde{A}/\widetilde{A}_o \subseteq \{A, \widetilde{A}'', ..., \widetilde{A}^{(n-2)}, \widetilde{A}_o\}.$$

Так как группа $<\widetilde{A}, \bullet>$ является гомоморфным образом группы $\mathcal{A}$ и содержит подгруппу $\widetilde{A}_o$ порядка $|A|$, то, учитывая следствие 1.4.5, получим

**1.4.11. Следствие.** Если $< A, [\ ] >$ – конечная n-арная группа, то $|\widetilde{A}| = |A|k$, где k делит $n - 1$.

**1.4.12. Определение.** Если $|\widetilde{A}/\widetilde{A}_o| = n - 1$, то группа $<\widetilde{A}, \bullet>$ называется *универсальной обертывающей группой* n-арной группы $< A, [\ ] >$.

**1.4.13. Определение.** Группа $A^* = <\mathcal{A}, *>$ называется *универсальной обертывающей группой Поста* n-арной группы $<A, [\ ]>$.

Из 1) теоремы 1.4.9. вытекает

**1.4.14. Следствие.** Любая универсальная обертывающая группа n-арной группы изоморфна универсальной обертывающей группе Поста этой же n-арной группы.

**1.4.15. Определение.** Подгруппа $<\widetilde{A}_o, \bullet>$ группы $<\widetilde{A}, \bullet>$ называется *соответствующей группой* n-арной группы $<A, [\ ]>$, а подгруппа $A_o$ группы $\mathcal{A}$ называется *соответствующей группой Поста* для n-арной группы $<A, [\ ]>$.

**1.4.16. Замечание.** Иногда понятие обертывающей группы расширяют, называя обертывающей группой n-арной группы $<A, [\ ]>$, любую группу, изоморфную группе из определения 1.4.7. Соответственно универсальной обертывающей группой n-арной группы $<A, [\ ]>$ называется любая группа, изоморфная группе $\mathcal{A}$.

**1.4.17. Теорема** (обратная теорема Поста). Пусть G группа, H – ее инвариантная подгруппа, факторгруппа $G/H = <gH>$ – циклическая с образующим элементом $gH$, $|G/H| = k$, $k$ делит $n - 1$. Определим на множестве $A = gH$ n-арную операцию $[\ ]$ по правилу

$$[a_1 a_2 ... a_n] = a_1 a_2 ... a_n.$$

Тогда $<A, [\ ]>$ – n-арная группа, причем G – обертывающая группа для $<A, [\ ]>$, т. е. $G = \widetilde{A}$ и, кроме того, $H = \widetilde{A}_o$.

*Доказательство.* Пусть

$$a_1 = gh_1, a_2 = gh_2, ..., a_n = gh_n \in A = gH.$$

Используя инвариантность H в G, можно показать, что

$$[a_1 a_2 ... a_n] = a_1 a_2 ... a_n = gh_1 gh_2 ... gh_n = g^n h', h' \in H,$$

т. е.

$$[a_1 a_2 ... a_n] = g^n h',\ h' \in H.$$

Так как по условию $A^k = H$ и $k$ делит $n - 1$, то $A^{n-1} = H$, откуда
$$A^{n-1} = (gH)^{n-1} = g^{n-1}H = H.$$

Следовательно, $g^{n-1} \in H$ и поэтому
$$[a_1 a_2 ... a_n] = g^n h' = g h'',$$

где $h'' \in H$. Таким образом,
$$[a_1 a_2 ... a_n] = g h'' \in gH = A,$$

т. е. множество $A$ замкнуто относительно $n$-арной операции [ ].

Ассоциативность $n$-арной операции [ ] является следствием ассоциативности бинарной операции в группе $G$.

Уравнение
$$[u a_1 ... a_{n-1}] = a \qquad (*)$$

в $A$, где $a_1 = g h_1, ..., a_{n-1} = g h_{n-1}$, $a = gh$ равносильно уравнению
$$u g h_1 ... g h_{n-1} = gh$$

в $G$, которое имеет решение $u = d$, т. е.
$$d g h_1 ... g h_{n-1} = gh.$$

Снова, используя инвариантность $H$ в $G$, получим
$$d g^{n-1} h' = gh, \; h' \in H.$$

Учитывая $g^{n-1} \in H$, имеем $g^{n-1} h' = c \in H$, т. е. $dc = gh$, откуда
$$d = gh c^{-1} = g h'' \in gH = A,$$

где $h'' = hc^{-1} \in H$. Этим доказана разрешимость в $A$ уравнения $(*)$. Аналогично доказывается разрешимость в $A$ уравнения
$$[a_1 ... a_{n-1} v] = a.$$

Так как $A = gH$ – образующий элемент факторгруппы $G/H$, то $A$ – порождающее множество группы $G$. Следовательно, $G = \widetilde{A}$. А так как $A^{n-1} = H$, то $H = \widetilde{A}_o$. ∎

Покажем, что некоторые из приведенных нами ранее примеров n-арных групп могут быть получены как следствия из теоремы 1.4.17. Заодно укажем для этих n-арных групп их обертывающие и соответствующие группы.

**1.4.18. Пример.** Пусть $T_n$ – множество всех нечетных подстановок степени n (пример 1.1.8). Так как факторгруппа $S_n/A_n$ имеет порядок 2 и порождается своим смежным классом $T_n$, то по теореме 1.4.17 $< T_n, [\ ] >$ – тернарная группа с тернарной операцией, производной от операции в группе $S_n$. Кроме того, симметрическая группа $S_n$ – универсальная обертывающая группа, а знакопеременная группа $A_n$ – соответствующая группа для тернарной группы $< T_n, [\ ] >$.

**1.4.19. Пример.** Пусть $E(2)$ – группа всех движений плоскости, $E_1(2)$ – ее нормальная подгруппа всех движений первого рода, $E_2(2)$ – множество всех движений второго рода из $E(2)$ (пример 1.1.14). Так как факторгруппа $E(2)/E_1(2)$ имеет порядок 2 и порождается своим смежным классом $E_2(2)$, то по теореме 1.4.17 $< E_2(2), [\ ] >$ – тернарная группа с тернарной операцией, производной от операции в группе $E(2)$. Универсальной обертывающей группой для $< E_2(2), [\ ] >$ является группа $E(2)$, а соответствующей группой – группа $E_1(2)$.

Аналогичный факт имеет место и в пространстве.

**1.4.20. Пример.** Множество $E_2(3)$ всех движений второго рода пространства является тернарной группой с тернарной операцией $[\ ]$ производной от операции в группе $E(3)$ – всех движений пространства. Причем, $E(3)$ – универсальная обертывающая группа для $< E_2(3), [\ ] >$, а соответствующей группой является группа $E_1(3)$ всех движений первого рода пространства.

**1.4.21. Пример.** Множество $E_2(1)$ поворотов прямой в некоторой плоскости на $180°$ (пример 1.1.15) является тернарной группой с тернарной операцией $[\ ]$, производной от операции в группе $E(1)$ самосовмещений прямой в выбранной плоскости. Сама же эта группа является универсальной обертывающей для $< E_2(1), [\ ] >$, а ее нормальная подгруппа $E_1(1)$ всех скольжений прямой по себе – соответствующей группой.

**1.4.22. Пример.** Пусть $1 + (n-1)Z$ – множество всех целых чисел, дающих при делении на $n - 1$ в остатке единицу (пример 1.1.18). Так как $Z/(n-1)Z = Z_{n-1}$ – циклическая группа порядка $n - 1$, порождаемая своим смежным классом $1 + (n-1)Z$, то по теореме 1.4.17

< 1 + (n-1)Z, [ ] > – n-арная группа с n-арной операцией, производной от операции сложения чисел. По той же теореме, группа Z всех целых чисел является универсальной обертывающей, а ее подгруппа (n – 1)Z – всех целых чисел кратных n – 1 – соответствующей группой для n-арной группы < 1 + (n – 1)Z, [ ] >.

В частности, множество всех нечетных чисел является тернарной группой, для которой универсальной обертывающей является группа Z, а соответствующей группой – ее подгруппа 2Z четных чисел.

**1.4.23. Пример.** Пусть $< B_n, [\ ] >$ – тернарная группа отражений из примера 1.2.6. Так как всякий поворот может быть представлен в виде произведения отражений, то группа $D_n$ порождается всеми отражениями. Это означает, что $D_n$ – обертывающая группа для тернарной группы $< B_n, [\ ] >$. Ясно, что $bB_n = C_n$. Поэтому группа $C_n$ является соответствующей для $< B_n, [\ ] >$. Из равенства $|D_n : C_n| = 2$ вытекает, что $D_n$ – даже универсальная обертывающая для $< B_n, [\ ] >$.

**1.4.24. Замечание.** То, что $< B_n, [\ ] >$ – тернарная группа, можно доказать и внутренними средствами теории n-арных групп. Действительно, рассмотрим факторгруппу $D_n/C_n = < B_n >$. Так как она циклическая порядка 2, делящего 3–1, то по теореме 1.4.17, образующий смежный класс $B_n$ является тернарной группой относительно тернарной операции [ ].

## §1.5. ТЕОРЕМА ГЛУСКИНА-ХОССУ

Теорема Глускина-Хоссу[13, 14] утверждает, что на всякой n-арной группе $< A, [\ ] >$ можно определить бинарную операцию $*$ и отображение $\beta$, а также выбрать элемент $d \in A$ так, что $< A, * >$ – группа, $\beta$ – ее автоморфизм, и выполняются следующие условия:

$$[x_1 x_2 ... x_n] = x_1 * x_2^{\beta} * ... * x_n^{\beta^{n-1}} * d, \ x_1, x_2, ..., x_n \in A; \quad (1)$$

$$d^{\beta} = d; \quad (2)$$

$$d * x = x^{\beta^{n-1}} * d, \ x \in A. \quad (3)$$

Верно и обратное утверждение (обратная теорема Глускина-Хоссу): если элемент d группы $< A, * >$ и ее автомор-физм β удовлетворяют условиям (2) и (3), то $< A, [\ ] >$ – n-арная группа с n-арной операцией (1).

Мы получим теорему Глускина-Хоссу в качестве следствия более общего результата [15].

**1.5.1. Лемма.** Пусть $< A, * >$ – группа, $\alpha$ – её автоморфизм, $\beta = \alpha^{-1}$; $x, d \in A$, $d^\alpha = d$. Если для фиксированных целых m и k верно

$$d * x^{\alpha^{m-k}} = x^{\beta^k} * d,$$

то последнее равенство справедливо для любого целого k.

*Доказательство.* Для произвольного целого t будем иметь

$$(d * x^{\alpha^{m-k}})^{\alpha^{k-t}} = (x^{\beta^k} * d)^{\alpha^{k-t}},$$

$$d^{\alpha^{k-t}} * (x^{\alpha^{m-k}})^{\alpha^{k-t}} = (x^{\beta^k})^{\alpha^{k-t}} * d^{\alpha^{k-t}},$$

$$d * x^{\alpha^{m-k+k-t}} = (x^{\beta^k})^{\beta^{t-k}} * d,$$

$$d * x^{\alpha^{m-t}} = x^{\beta^{k+t-k}} * d,$$

$$d * x^{\alpha^{m-t}} = x^{\beta^t} * d. \qquad \blacksquare$$

На произвольной n-арной группе $< A, [\ ] >$ определим бинарную операцию

$$x \,@\, y = [x a_1 \ldots a_{n-2} y]$$

и отображения

$$\alpha : x \mapsto x^\alpha = [a_1 \ldots a_{n-2} x a], \quad \beta : x \mapsto x^\beta = [a x a_1 \ldots a_{n-2}],$$

где a – фиксированный элемент из A, $a_1 \ldots a_{n-2}$ – обратная последовательность для элемента a. Положим также

$$d = [\underbrace{a \ldots a}_{n}].$$

Если $b_1...b_{n-2}$ – другая обратная последовательность для элемента a, то, по предложению 1.3.9, последовательности $a_1...a_{n-2}$ и $b_1...b_{n-2}$ эквивалентны, и поэтому

$$[xa_1...a_{n-2}y] = [xb_1...b_{n-2}y].$$

Следовательно, операция @ определена правильно.

**1.5.2. Предложение.** $< A, @ >$ – группа с единицей a.

*Доказательство.* Ассоциативность операции очевидна. Так $\tilde{a}a$ и $a\tilde{a}$ – нейтральные последовательности, где $\tilde{a} = a_1...a_{n-2}$ – обратная для a, то

$$x @ a = [x\tilde{a}a] = x = [a\tilde{a}x] = a @ x$$

для любого $x \in A$. Поэтому a – единица полугруппы $< A, @ >$.

Пусть теперь $\tilde{x}$ – обратная последовательность для элемента $x \in A$. Так как

$$x @ [a\tilde{x}a] = [x\tilde{a}a\tilde{x}a] = [x\tilde{x}a] = a,$$

$$[a\tilde{x}a] @ x = [a\tilde{x}a\tilde{a}x] = [a\tilde{x}x] = a,$$

то $x^{-1} = [a\tilde{x}a]$ – обратный элемент для x. ∎

В дальнейшем нам понадобится понятие автоморфизма n-арной группы.

**1.5.3. Определение.** Биекция f n-арной группы $< A, [\ ] >$ называется *автоморфизмом*, если

$$[a_1 a_2 ... a_n]^f = [a_1^f a_2^f ... a_n^f]$$

для любых $a_1, a_2, ..., a_n \in A$.

**1.5.4. Предложение.** Отображения $\alpha$ и $\beta$ являются автоморфизмами n-арной группы $< A, [\ ] >$ и бинарной группы $< A, @ >$, причем

$$\beta = \alpha^{-1}, d^\alpha = d^\beta = d.$$

***Доказательство.*** Из условия 2) определения 1.1.1 вытекает, что $\alpha$ и $\beta$ – биекции. А так как

$$[a_1a_2...a_n]^{\alpha} = [\tilde{a}\,[a_1a_2...a_n]\,a] = [\tilde{a}\,a_1a_2...a_na] =$$

$$= [\tilde{a}\,a_1a\tilde{a}\,a_2a\tilde{a}\,...a_{n-1}a\tilde{a}\,a_na] =$$

$$= [[\tilde{a}\,a_1a][\tilde{a}\,a_2a]...[\tilde{a}\,a_na]] = [a_1^{\alpha}a_2^{\alpha}...a_n^{\alpha}]$$

для любых $a_1, a_2, ..., a_n \in A$, то $\alpha$ – автоморфизм n-арной группы $< A, [\ ] >$.

Кроме того,

$$(a_1 @ a_2)^{\alpha} = [\tilde{a}\,[a_1\tilde{a}\,a_2]a] = [\tilde{a}\,a_1\tilde{a}\,a_2a] =$$

$$= [\tilde{a}\,a_1a\tilde{a}\,\tilde{a}\,a_2a] = [[\tilde{a}\,a_1a]\,\tilde{a}\,[\tilde{a}\,a_2a]] = a_1^{\alpha} @ a_2^{\alpha}$$

для любых $a_1, a_2 \in A$, т. е. $\alpha$ – автоморфизм группы $< A, @ >$.

Для $\beta$ доказательство проводится аналогично.

Так как

$$x^{\alpha\beta} = [a[\tilde{a}\,xa]\,\tilde{a}\,] = [a\tilde{a}\,xa\tilde{a}\,] = x,$$

$$x^{\beta\alpha} = [\tilde{a}\,[ax\tilde{a}\,]a] = [\tilde{a}\,ax\tilde{a}\,a] = x$$

для любого $x \in A$, то $\beta = \alpha^{-1}$.

А так как

$$d^{\alpha} = [\tilde{a}\,[\underbrace{a...a}_{n}]a] = [\tilde{a}\,a\underbrace{a...a}_{n}] = [\underbrace{a...a}_{n}] = d,$$

$$d^{\beta} = [a[\underbrace{a...a}_{n}]\tilde{a}\,] = [\underbrace{a...a}_{n}a\tilde{a}\,] = [\underbrace{a...a}_{n}] = d,$$

то $d^{\alpha} = d$, $d^{\beta} = d$. ∎

**1.5.5. Теорема.** На любой n-арной группе $< A, [\ ] >$ выполняются следующие тождества:

1) $[x_1...x_n] = x_1 @ x_2^{\beta} @ ... x_i^{\beta^{i-1}} @ d @ x_{i+1}^{\alpha^{n-1-i}} @ ... x_{n-1}^{\alpha} @ x_n$,

$i = 0, 1, ..., n$;

2) $d @ x^{\alpha^{n-1-k}} = x^{\beta^k} @ d$, $k \in Z$.

***Доказательство.*** Для доказательства обоих тождеств будем использовать нейтральность последовательностей

$$\underbrace{a...a}_{i}\underbrace{a_1^{n-2}...a_1^{n-2}}_{i},\ \underbrace{a_1^{n-2}...a_1^{n-2}}_{i}\underbrace{a...a}_{i},$$

где $i = 1, 2, \ldots$ .

1) $x_1 \text{@} x_2^{\beta} \text{@} x_3^{\beta^2} \text{@} \ldots x_i^{\beta^{i-1}} \text{@} d \text{@} x_{i+1}^{\alpha^{n-1-i}} \text{@} \ldots$

$\ldots x_{n-2}^{\alpha^2} x_{n-1}^{\alpha} \text{@} x_n =$

$= [x_1 a_1^{n-2} [a x_2 a_1^{n-2}] a_1^{n-2} [a a x_3 a_1^{n-2} a_1^{n-2}] a_1^{n-2} \ldots$

$\ldots [\underbrace{a...a}_{i-1} x_i \underbrace{a_1^{n-2}...a_1^{n-2}}_{i-1}] a_1^{n-2} [\underbrace{a...a}_{n}] a_1^{n-2}$

$[\underbrace{a_1^{n-2}...a_1^{n-2}}_{n-i-1} x_{i+1} \underbrace{a...a}_{n-i-1}] a_1^{n-2} \ldots$

$\ldots [a_1^{n-2} a_1^{n-2} x_{n-2} a a] a_1^{n-2} [a_1^{n-2} x_{n-1} a] a_1^{n-2} x_n] =$

$= [x_1 a_1^{n-2} a x_2 a_1^{n-2} a_1^{n-2} a a x_3 a_1^{n-2} a_1^{n-2} a_1^{n-2} \ldots$

$\ldots \underbrace{a...a}_{i-1} x_i \underbrace{a_1^{n-2}...a_1^{n-2}}_{i} \underbrace{a...a}_{i} \underbrace{a...a}_{n-i}$

$\underbrace{a_1^{n-2}...a_1^{n-2}}_{n-i} x_{i+1} \underbrace{a...a}_{n-i-1} a_1^{n-2} \ldots$

$\ldots a_1^{n-2} a_1^{n-2} x_{n-2} a a a_1^{n-2} a_1^{n-2} x_{n-1} a a_1^{n-2} x_n] =$

$= [x_1 x_2 x_3 \ldots x_i x_{i+1} \ldots x_{n-2} x_{n-1} x_n] = [x_1 \ldots x_n].$

2) Пусть $k = 0, 1, \ldots, n - 1$. Так как

$$d @ x^{\alpha^{n-1-k}} = [[\underbrace{a...a}_{n}]a_1^{n-2}[\underbrace{a_1^{n-2}...a_1^{n-2}}_{n-k-1}x\underbrace{a...a}_{n-k-1}]] =$$

$$= [\underbrace{a...a}_{k}\underbrace{a...a}_{n-k}\underbrace{a_1^{n-2}...a_1^{n-2}}_{n-k}x\underbrace{a...a}_{n-k-1}] = [\underbrace{a...a}_{k}x\underbrace{a...a}_{n-k-1}],$$

$$x^{\beta^k} @ d = [[\underbrace{a...a}_{k}x\underbrace{a_1^{n-2}...a_1^{n-2}}_{k}]a_1^{n-2}[\underbrace{a...a}_{n}]] =$$

$$= [\underbrace{a...a}_{k}x\underbrace{a_1^{n-2}...a_1^{n-2}}_{k+1}\underbrace{a...a}_{k+1}\underbrace{a...a}_{n-k-1}] = [\underbrace{a...a}_{k}x\underbrace{a...a}_{n-k-1}],$$

то

$$d @ x^{\alpha^{n-1-k}} = x^{\beta^k} @ d.$$

Для произвольного $k \in Z$ применяется лемма 1.5.1. ∎

Придавая n, i и k в тождествах 1) и 2) теоремы 1.5.5 конкретные значения, можно получить большое число новых тождеств, некоторые из которых приведены ниже.

**1.5.6. Следствие.** На любой n-арной группе $< A, [\ ] >$ выполнены тождества:

$$[x_1...x_n] = d @ x_{i+1}^{\alpha^{n-1}} @ ... x_{n-1}^{\alpha} @ x_n;$$

$$[x_1...x_n] = x_1 @ d @ x_{i+1}^{\alpha^{n-2}} @ ... x_{n-1}^{\alpha} @ x_n;$$

$$[x_1...x_n] = x_1 @ x_2^{\beta} @ ... x_{n-1}^{\beta^{n-2}} @ d @ x_n;$$

$$[x_1...x_n] = x_1 @ x_2^{\beta} @ ... x_n^{\beta^{n-1}} @ d;$$

$$[x_1...x_n] = x_1 @ x_2^{\beta} @ ... x_m^{\beta^{m-1}} @ d @ x_{m+1}^{\alpha^{m-1}} @ ... x_{2m-1}^{\alpha} @ x_{2m};$$

$$d @ x^{\alpha^{n-1}} = x @ d;$$

$$d @ x = x^{\beta^{n-1}} @ d.$$

Пусть $<A, *>$ – группа, $\alpha$ – её автоморфизм, $\beta = \alpha^{-1}$, $d \in A$, $d^\alpha = d$. Определим на $A$ для каждого $i = 0, 1, \ldots, n$ n-арную операцию

$$[x_1 \ldots x_n]_{(i)} = x_1 * x_2^\beta * \ldots x_i^{\beta^{i-1}} * d * x_{i+1}^{\alpha^{n-1-i}} * x_{n-1}^\alpha * x_n.$$

**1.5.7. Лемма.** Если для некоторого целого $k$ на $A$ выполнено тождество 2) теоремы 1.5.5, то

$$[x_1 \ldots x_n]_{(i)} = [x_1 \ldots x_n]_{(j)}$$

для любых $i, j \in \{0, 1, \ldots, n\}$.

*Доказательство.* По лемме 1.5.1 тождество 2) верно для любых $k \in \mathbb{Z}$. Применяя это тождество последовательно для $k = 0, 1, \ldots, n-1$, получим

$$d * x_1^{\alpha^{n-1}} * x_2^{\alpha^{n-2}} * \ldots x_{n-1}^\alpha * x_n = x_1 * d * x_2^{\alpha^{n-2}} * \ldots x_{n-1}^\alpha * x_n =$$

$$= x_1 * x_2^\beta * d * \ldots x_{n-1}^\alpha * x_n = \ldots$$

$$\ldots = x_1 * x_2^\beta * \ldots x_{i-1}^{\beta^{i-2}} * d * x_i^{\alpha^{n-1-(i-1)}} * x_{i+1}^{\alpha^{n-1-i}} * \ldots x_{n-1}^\alpha * x_n =$$

$$= x_1 * x_2^\beta * \ldots x_{i-1}^{\beta^{i-2}} * x_i^{\beta^{i-1}} * d * x_{i+1}^{\alpha^{n-1-i}} * \ldots x_{n-1}^\alpha * x_n = \ldots$$

$$\ldots = x_1 * x_2^\beta * \ldots x_{n-2}^{\beta^{n-3}} * d * x_{n-1}^\alpha * x_n =$$

$$= x_1 * x_2^\beta * \ldots x_{n-2}^{\beta^{n-3}} * x_{n-1}^{\beta^{n-2}} * d * x_n =$$

$$= x_1 * x_2^\beta * \ldots x_{n-1}^{\beta^{n-2}} * x_n^{\beta^{n-1}} * d.$$

Таким образом, мы доказали, что

$$[x_1...x_n]_{(0)} = [x_1...x_n]_{(1)} = \ldots = [x_1...x_n]_{(n)}. \qquad \blacksquare$$

Следующая теорема является обратной к теореме 1.5.5.

**1.5.8. Теорема.** Пусть $< A, * >$ – группа, $\alpha$ – её автоморфизм, $\beta = \alpha^{-1}$, $d \in A$, $d^\alpha = d$ и для некоторого целого k на A выполнено тождество

$$d * x^{\alpha^{n-1-k}} = x^{\beta^k} * d.$$

Тогда на A можно определить n-арную операцию [ ] так, что

$$[x_1...x_n] = x_1 * x_2^\beta * ... x_i^{\beta^{i-1}} * d * x_{i+1}^{\alpha^{n-1-i}} * ... x_{n-1}^\alpha * x_n$$

для любого i = 0, 1, …, n, при этом $< A, [\ ] >$ – n-арная группа, а $\alpha$ и $\beta$ – её автоморфизмы.

*Доказательство.* Положим

$$[x_1...x_n] = [x_1...x_n]_{(0)}.$$

Тогда по лемме 1.5.7

$$[x_1...x_n] = [x_1...x_n]_{(i)}$$

для любого i = 0, 1, …, n. Следовательно, требуемое тождество выполняется.

Покажем ассоциативность n-арной операции [ ]. Действительно, используя $d * x^{\alpha^{n-1}} = x * d$, будем иметь

$$[x_1...x_{n-1}[x_n...x_{2n-1}]] = [x_1...x_{n-1}[x_n...x_{2n-1}]_{(0)}]_{(0)} =$$

$$= d * x_1^{\alpha^{n-1}} * ... x_{n-2}^{\alpha^2} * (x_{n-1}^\alpha * d) * x_n^{\alpha^{n-1}} * ... x_{2n-2}^\alpha * x_{2n-1} =$$

$$= d * x_1^{\alpha^{n-1}} * ... x_{n-2}^{\alpha^2} * (d * (x_{n-1}^\alpha)^{\alpha^{n-1}}) * x_n^{\alpha^{n-1}} * ... x_{2n-2}^\alpha * x_{2n-1} =$$

$$= d * x_1^{\alpha^{n-1}} * ... x_{n-2}^{\alpha^2} * d^\alpha * x_{n-1}^{\alpha^n} * x_n^{\alpha^{n-1}} * ... x_{2n-2}^\alpha * x_{2n-1} =$$

$$= d * x_1^{\alpha^{n-1}} * \ldots x_{n-2}^{\alpha^2} * (d * x_{n-1}^{\alpha^{n-1}} * x_n^{\alpha^{n-2}} * \ldots x_{2n-2})^\alpha * x_{2n-1} =$$

$$= [x_1 \ldots x_{n-2}[x_{n-1} \ldots x_{2n-2}]_{(0)} x_{2n-1}]_{(0)} = [x_1 \ldots x_{n-2}[x_{n-1} \ldots x_{2n-2}]x_{2n-1}] =$$

$$= d * x_1^{\alpha^{n-1}} * \ldots x_{n-3}^{\alpha^3} * (x_{n-2}^{\alpha^2} * d) * x_{n-1}^{\alpha^n} * \ldots x_{2n-3}^{\alpha^2} * x_{2n-2}^{\alpha} * x_{2n-1} =$$

$$= d * x_1^{\alpha^{n-1}} * \ldots x_{n-3}^{\alpha^3} * (d * (x_{n-2}^{\alpha^2})^{\alpha^{n-1}}) * x_{n-1}^{\alpha^n} * \ldots x_{2n-3}^{\alpha^2} * x_{2n-2}^{\alpha} * x_{2n-1} =$$

$$= d * x_1^{\alpha^{n-1}} * \ldots x_{n-3}^{\alpha^3} * d^{\alpha^2} * x_{n-2}^{\alpha^{n+1}} * x_{n-1}^{\alpha^n} * \ldots x_{2n-3}^{\alpha^2} * x_{2n-2}^{\alpha} * x_{2n-1} =$$

$$= d * x_1^{\alpha^{n-1}} * \ldots x_{n-3}^{\alpha^3} * (d * x_{n-2}^{\alpha^{n-1}} * x_{n-1}^{\alpha^{n-2}} * \ldots x_{2n-3})^{\alpha^2} * x_{2n-2}^{\alpha} * x_{2n-1} =$$

$$= [x_1 \ldots x_{n-3}[x_{n-2} \ldots x_{2n-3}]_{(0)} x_{2n-2} x_{2n-1}]_{(0)} =$$

$$= [x_1 \ldots x_{n-3}[x_{n-2} \ldots x_{2n-3}] x_{2n-2} x_{2n-1}] = \ldots = [[x_1 \ldots x_n]x_{n+1} \ldots x_{2n-1}].$$

Ясно, что каждое из уравнений

$$[a_1 \ldots a_{i-1} x_i a_{i+1} \ldots a_n] = a, \, i = 1, \ldots, n$$

разрешимо в A относительно $x_i$. Таким образом, доказано, что $< A, [\,] >$ – n-арная группа.

Так как

$$[x_1 x_2 \ldots x_{n-1} x_n]^\alpha = [x_1 x_2 \ldots x_{n-1} x_n]_{(0)}^\alpha =$$

$$= (d * x_1^{\alpha^{n-1}} * x_2^{\alpha^{n-2}} * \ldots x_{n-1}^{\alpha} * x_n)^\alpha =$$

$$= d * (x_1^\alpha)^{\alpha^{n-1}} * (x_2^\alpha)^{\alpha^{n-2}} * \ldots (x_{n-1}^\alpha)^\alpha * x_n^\alpha = [x_1^\alpha x_2^\alpha \ldots x_{n-1}^\alpha x_n^\alpha]_{(0)} =$$

$$= [x_1^\alpha x_2^\alpha \ldots x_{n-1}^\alpha x_n^\alpha],$$

то $\alpha$ – автоморфизм n-арной группы $< A, [\,] >$.

Так как

$$[x_1 x_2 \ldots x_{n-1} x_n]^\beta = [x_1 x_2 \ldots x_{n-1} x_n]_{(0)}^{\alpha^{-1}} =$$

$$= (d * x_1^{\alpha^{n-1}} * x_2^{\alpha^{n-2}} * \ldots x_{n-1}^{\alpha} * x_n)^{\alpha^{-1}} =$$

$$= d^{\alpha^{-1}} * (x_1^{\alpha^{-1}})^{\alpha^{n-1}} * (x_2^{\alpha^{-1}})^{\alpha^{n-2}} * \ldots (x_{n-1}^{\alpha^{-1}})^{\alpha} * x_n^{\alpha^{-1}} =$$

$$= d * (x_1^{\beta})^{\alpha^{n-1}} * (x_2^{\beta})^{\alpha^{n-2}} * \ldots (x_{n-1}^{\beta})^{\alpha} * x_n^{\beta} =$$

$$= [x_1^{\beta} x_2^{\beta} \ldots x_{n-1}^{\beta} x_n^{\beta}]_{(0)} = [x_1^{\beta} x_2^{\beta} \ldots x_{n-1}^{\beta} x_n^{\beta}],$$

то $\beta$ – автоморфизм n-арной группы $< A, [\ ] >$. ∎

Ясно, что отображение $\beta$ из последней теоремы является автоморфизмом группы $< A, * >$ и, кроме того, $d^{\beta} = d$. Поэтому, полагая в теореме 1.5.8, $k = n - 1$, $i = n$, получим обратное утверждение к теореме Глускина-Хоссу.

Теорема 1.5.8 и, в частности, обратная теорема Глускина-Хоссу позволяют строить различные примеры n-арных групп.

**1.5.9. Пример.** Если $\varepsilon$ – тождественный автоморфизм группы A, и $a \in Z(A)$, то

$$a^{\varepsilon} = a;\ ax = xa = x^{\varepsilon^{n-1}} a, x \in A,$$

т. е. для $\varepsilon$ и $a$ выполняются условия (2) и (3). Поэтому по обратной теореме Глускина-Хоссу $< A, [\ ] >$ – n-арная группа с n-арной операцией $[a_1 a_2 \ldots a_n] = a_1 a_2 \ldots a_n a$.

Этот факт был установлен ранее (пример 1.1.6) непосредственно с помощью определения n-арной группы.

**1.5.10. Пример.** Пусть A – группа, a – произвольный элемент из A. Положим $\beta: x \to axa^{-1}$ – внутренний автоморфизм группы, $d = a^{n-1}$. Так как

$$dx = a^{n-1}x,\ x^{\beta^{n-1}} d = a^{n-1} x (a^{-1})^{n-1} a^{n-1} = a^{n-1}x,$$

то

$$dx = x^{\beta^{n-1}} d,\ x \in A,$$

т. е. для $\beta$ и $d$ выполняется условие (3). Условие (2) для них также выполняется, так как

$$d^{\beta} = a a^{n-1} a^{-1} = a^{n-1} = d.$$

Применяя обратную теорему Глускина-Хоссу, заключаем, что $<A, [\ ]>$ – n-арная группа с n-арной операцией

$$[a_1a_2...a_n] = a_1(aa_2a^{-1})(aaa_3a^{-1}a^{-1})...(\underbrace{a...a}_{n-1}a_n\underbrace{a^{-1}...a^{-1}}_{n-1})a^{n-1} =$$

$$= a_1aa_2aa_3a...aa_n,$$

т. е.

$$[a_1a_2...a_n] = a_1aa_2a...aa_n.$$

Для построения дальнейших примеров нам понадобится следующее утверждение.

**1.5.11. Предложение.** Если $a^{n-1} = 1$ для некоторого элемента $a$ тела $T$, где $n \geq 2$, то $<T, [\ ]>$ – n-арная группа с n-арной операцией

$$[x_1x_2...x_n] = x_1+ax_2+...+a^{n-2}x_{n-1}+x_n.$$

*Доказательство.* Определим преобразование $\beta : x \mapsto ax$ тела $T$. Так как $T^*$ – группа по умножению, то $\beta$ – биекция, а так как

$$(x+y)^\beta = a(x+y) = ax+ay = x^\beta+y^\beta,$$

то $\beta$ – автоморфизм группы $<T, +>$. Ясно, что $0^\beta = 0$ и, кроме того,

$$0+x = x = a^{n-1}x = x^{\beta^{n-1}} = x^{\beta^{n-1}}+0,$$

т. е.

$$0+x = x^{\beta^{n-1}}+0.$$

Так как для элемента $0$ и автоморфизма $\beta$ выполняются все условия обратной теоремы Глускина-Хоссу, то, согласно этой теореме, $<T, [\ ]>$ – n-арная группа с n-арной операцией.

$$[x_1x_2...x_n] = x_1+x_2^\beta+...x_{n-1}^{\beta^{n-2}}+x_n^{\beta^{n-1}}+0 =$$

$$= x_1+ax_2+...a^{n-2}x_{n-1}+a^{n-1}x_n = x+ax_2+...+a^{n-2}x_{n-1}+x_n. \qquad\blacksquare$$

**1.5.12. Пример** [1]**.** Пусть T = C – поле всех комплексных чисел, $\varepsilon = \cos\dfrac{2\pi}{n-1} + i\sin\dfrac{2\pi}{n-1} \in C$. Так как $\varepsilon^{n-1} = 1$, то, по предыдущему предложению < C, [ ] > – n-арная группа с n-арной операцией

$$[z_1 z_2 ... z_n] = z_1 + \varepsilon z_2 + ... + \varepsilon^{n-2} z_{n-1} + z_n.$$

**1.5.13. Пример.** Пусть снова T = C. Так как $i^4 = 1$, то, согласно предложению 1.5.11, < C, [ ] > – 5-арная группа с 5-арной операцией

$$[z_1 z_2 z_3 z_4 z_5] = z_1 + i z_2 + i^2 z_3 + i^3 z_4 + z_5,$$

$$[z_1 z_2 z_3 z_4 z_5] = z_1 + i z_2 - z_3 - i z_4 + z_5.$$

**1.5.14. Пример.** Пусть H – тело кватернионов. Так как

$$i^2 = j^2 = k^2 = -1;$$

$$i^3 = -i,\ j^3 = -j,\ k^3 = -k;$$

$$i^4 = j^4 = k^4 = 1,$$

то, согласно предложению 1.5.11, < H, [ ]$_i$ >, < H, [ ]$_j$ > и < H, [ ]$_k$ > – 5-арные группы с пятиарными операциями

$$[x_1 x_2 x_3 x_4 x_5]_i = x_1 + i x_2 - x_3 - i x_4 + x_5,$$

$$[x_1 x_2 x_3 x_4 x_5]_j = x_1 + j x_2 - x_3 - j x_4 + x_5,$$

$$[x_1 x_2 x_3 x_4 x_5]_k = x_1 + k x_2 - x_3 - k x_4 + x_5.$$

## §1.6. СВЯЗЬ МЕЖДУ ТЕОРЕМАМИ ПОСТА И ГЛУСКИНА-ХОССУ

Из формулировки теоремы Поста видно, что бинарная и n-арная операции связаны в ней самым простым образом, что является главным достоинством этой теоремы. К числу же недостатков следует отнести то, что n-арная группа и группа, в которую она вкладывается, имеют различные носители.

Теорема Глускина-Хоссу свободна от этого недостатка, так как в ней n-арная группа и группа, к которой она приводима, имеют общий носитель. Платой за такое экономное вложение служит довольно сложная зависимость между бинарной и n-арной операциями.

Приведённые соображения наводят на мысль, что теоремы Поста и Глускина-Хоссу, предлагающие, казалось бы на первый взгляд, два совершенно различных подхода к изучению n-арных групп, на самом деле должны каким-то образом вкладываться в некую общую схему. Ответу на этот вопрос посвящен этот параграф, в котором установлено, что обе интересующие нас теоремы действительно являются частными случаями общего результата.

**1.6.1. Предложение.** Группы $<A,@>$ и $<A_o,*>$ изоморфны.

*Доказательство.* Согласно 1) предложения 1.3.7, группу $A_o$ можно представить в виде
$$A_o = A^{(n-1)} = \{\theta(xa_1...a_{n-2}) \mid x \in A\}.$$
Отображение
$$\varphi: x \mapsto \theta(xa_1...a_{n-2})$$
является биекцией $A$ на $A_o$. Кроме того,
$$\varphi(x @ y) = \theta(x @ ya_1...a_{n-2}) = \theta([xa_1...a_{n-2}y]a_1...a_{n-2}) =$$
$$= \theta(xa_1...a_{n-2}ya_1...a_{n-2}) = \theta(xa_1...a_{n-2})\theta(ya_1...a_{n-2}) = \varphi(x)\varphi(y).$$
Следовательно, $\varphi$ − изоморфизм группы $<A,@>$ на группу $<A_o,*>$. ∎

**1.6.2. Следствие.** Для любых $a, c \in A$ группы $<A, @>$ и $<A, ©>$ изоморфны.

Предложение 1.6.1 устанавливает связь между теоремой Поста и теоремой Глускина-Хоссу. Ниже будет показано, что эта связь на самом деле является более тесной.

Зафиксируем в n-арной группе $<A, [\,]>$ последовательность $\alpha$ ее элементов, с помощью которой определим на $A$ новую n-арную операцию.

$$[x_1 x_2 ... x_n]_\alpha = [x_1 \alpha x_2 \alpha ... x_n].$$

Следующее утверждение является непосредственным следствием теоремы 1.3.3.

**1.6.3. Предложение.** Если $\alpha \theta \beta$, то $[\,]_\alpha = [\,]_\beta$.

Легко проверяется и справедливость следующего утверждения.

**1.6.4. Предложение.** $<A, [\,]_\alpha>$ – n-арная группа.

**1.6.5. Предложение.** Если $l(\alpha) \equiv l(\beta) \,(\mathrm{mod}\, n-1)$, то $<A, [\,]_\alpha> \simeq <A, [\,]_\beta>$.

*Доказательство.* Если $\widetilde{\beta}$ – обратная последовательность для $\beta$, то $l(\widetilde{\beta}) + l(\beta) \equiv 0 \,(\mathrm{mod}\, n-1)$, откуда, учитывая условие предложения, получаем $l(\widetilde{\beta}) + l(\alpha) \equiv 0 \,(\mathrm{mod}\, n-1)$. Поэтому можно рассмотреть отображение

$$\varphi: x \mapsto [\widetilde{\beta} x \alpha],$$

которое, как нетрудно убедиться, является биекцией множества $A$. А так как, кроме того,

$$\varphi([x_1 x_2 ... x]_\alpha) = [\widetilde{\beta}\,[x_1 \alpha x_2 \alpha ... x_n]\alpha] =$$

$$= [[\widetilde{\beta} x_1 \alpha]\beta[\widetilde{\beta} x_2 \alpha]\beta...[\widetilde{\beta} x_n \alpha]] =$$

$$= [\varphi(x_1)\beta\varphi(x_2)\beta...\varphi(x_n)] = [\varphi(x_1)\varphi(x_2)...\varphi(x_n)]_\beta,$$

то $\varphi$ – изоморфизм n-арной группы $<A, [\,]_\alpha>$ на n-арную группу $<A, [\,]_\beta>$. ∎

**1.6.6. Замечание.** Отображение

$$\psi: x \to [\alpha x \widetilde{\beta}]$$

также является изоморфизмом $<A, [\ ]_\alpha>$ на $<A, [\ ]_\beta>$.

**1.6.7. Замечание.** Ясно, что $[\ ] = [\ ]_\varnothing = [\ ]_\varepsilon$, где $\varnothing$ – пустая, $\varepsilon$ – нейтральная последовательности.

**1.6.8. Следствие.** Если $l(\alpha) \equiv 0 \pmod{n-1}$, то имеет место изоморфизм $<A,[\ ]> \simeq <A,[\ ]_\alpha>$.

На множестве $\mathcal{A}$ (см. §1.3) для любого $u \in \mathcal{A}$ определим бинарную
$$x_1*(u)x_2 = x_1*u*x_2 = x_1ux_2$$
и n-арную
$$\lfloor x_1x_2...x_n \rfloor_u = \lfloor x_1ux_2u...x_n \rfloor = x_1*(u)x_2*(u)...x_n,$$
операции, где $*$ – операция в полугруппе $\mathcal{A} = F_A/\theta$, $\lfloor\ \rfloor$ – n-арная операция из следствия 1.4.3, производная от операции $*$.

**1.6.9. Предложение.** Справедливы следующие утверждения:
1) $<\mathcal{A},*>$ и $<\mathcal{A},*(u)>$ – изоморфные группы;
2) $<\mathcal{A},\lfloor\ \rfloor>$ и $<\mathcal{A},\lfloor\ \rfloor_u>$ – изоморфные n-арные группы.

*Доказательство.* 1) легко проверяется, что $<\mathcal{A},*(u)>$ группа, а отображение $\mu: x \mapsto u^{-1}x$ – изоморфизм $<\mathcal{A},*>$ на $<\mathcal{A},*(u)>$.

2) То, что $<\mathcal{A},\lfloor\ \rfloor>$ – n-арная группа, установлено в следствии 1.4.3. Согласно предложению 1.6.4, $<\mathcal{A},\lfloor\ \rfloor_u>$ также является n-арной группой. Так как $\lfloor\ \rfloor$ и $\lfloor\ \rfloor_u$ – n-арные операции, производные от операций $*$ и $*(u)$ соответственно, то из 1) вытекает, что $\mu$ – изоморфизм $<\mathcal{A},\lfloor\ \rfloor>$ на $<\mathcal{A},\lfloor\ \rfloor_u>$. ∎

Для сокращения формулировки следующей леммы положим
$$A^{(o)} = A^{(n-1)}, \text{ т.е. } A_o = A^{(o)} = A^{(n-1)}.$$

**1.6.10. Лемма.** Если $u \in A^{(k)}$, то $\mu(A') = A^{(n-k)}$,

$\mu(A_o) = A^{(n-k-1)}$, где $k = 1, ..., n - 1$.

***Доказательство.*** Если $k = 1, ..., n - 2$, то для $u \in A^{(k)}$ имеем $u^{-1} \in A^{(n-k-1)}$. Поэтому

$$\mu(A') = u^{-1}A' = A^{(n-k)}, \quad \mu(A_o) = \mu(A^{(n-1)}) = u^{-1}A^{(n-1)} = A^{(n-k-1)}.$$

Если же $u \in A^{(n-1)}$, т.е. $k = n-1$, то $u^{-1} \in A^{n-1}$, $\mu(A') = u^{-1}A' = A'$,

$$\mu(A^{(o)}) = \mu(A^{(n-1)}) = u^{-1}A^{(n-1)} = A^{(n-1)} = A^{(o)}. \quad \blacksquare$$

Следующее предложение является следствием теоремы 1.4.2, предложения 1.6.9 и леммы 1.6.10.

**1.6.11. Предложение.** Если $u \in A^{(k)}$, где $k = 1, ..., n - 1$, то справедливы следующие утверждения:

1) группа $< \mathcal{A}, *(u) >$ порождается множеством $A^{(n-k)}$;

2) $< A^{(n-k-1)}, *(u) >$ — инвариантная подгруппа группы $< \mathcal{A}, *(u) >$;

3) $\mathcal{A}/A^{(n-k-1)}$ — циклическая группа порядка $n - 1$ с образующим элементом $A^{(n-k)}$.

Предложение 1.6.9, лемма 1.6.10 и изоморфизмы

$$< A, [\ ] > \simeq < A', \lfloor \ \rfloor > \text{ и } < A, @ > \simeq < A_o, * >$$

позволяют сформулировать еще одно утверждение.

**1.6.12. Предложение.** Если $u \in A^{(k)}$, где $k = 1, ..., n - 1$, то справедливы следующие утверждения:

1) $< A, [\ ] > \simeq < A^{(n-k)}, \lfloor \ \rfloor_u >$;

2) $< A, @ > \simeq < A^{(n-k-1)}, *(u) >$.

Из предложений 1.6.11 и 1.6.12 получаем

**1.6.13. Следствие.** Группы $< \mathcal{A}, *(u) >$ и $< A^{(n-k-1)}, *(u) >$ являются соответственно обертывающей и соответствующей для n-арной группы $< A^{(n-k)}, \lfloor \ \rfloor_u >$.

Для любого $u \in \mathcal{A}$ определим отображения $\alpha_u: \mathcal{A} \to \mathcal{A}$ и $\beta_u: \mathcal{A} \to \mathcal{A}$ по правилам

$$\alpha_u: x \mapsto uxu^{-1}, \quad \beta_u: x \mapsto u^{-1}xu$$

и зафиксируем элементы

$$d_u = u^{n-1}, \quad c_u = (u^{-1})^n.$$

Рассмотрим универсальную алгебру

$$\mathcal{A}_u = <\mathcal{A}, \{d_u, c_u, \alpha_u, \beta_u, *, *(u), \lfloor \ \rfloor, \lfloor \ \rfloor_u\}>,$$

где $d_u$ и $c_u$ – нульарные, $\alpha_u$ и $\beta_u$ – унарные, $*$ и $*(u)$ – бинарные, $\lfloor \ \rfloor$ и $\lfloor \ \rfloor_u$ – n-арные операции.

**1.6.14. Теорема.** Справедливы следующие утверждения:

1) $\lfloor x_1 x_2 ... x_n \rfloor =$

$$= x_1 *(u) x_2^{\beta_u} *(u) ... x_i^{\beta_u^{i-1}} *(u) c_u *(u) x_{i+1}^{\alpha_u^{n-1-i}} *(u) ... x_{n-1}^{\alpha_u} *(u) x_n,$$

$$\lfloor x_1 x_2 ... x_n \rfloor_u = x_1 * x_2^{\alpha_u} * ... x_i^{\alpha_u^{i-1}} * d_u * x_{i+1}^{\beta_u^{n-1-i}} * ... x_{n-1}^{\beta_u} * x_n$$

для любых $x_1, x_2, ..., x_n \in \mathcal{A}$ и любого $i = 0, 1, ..., n$;

2) $d_u * x^{\beta_u^{n-1-k}} = x^{\alpha_u^k} * d_u, \quad c_u *(u) x^{\alpha_u^{n-1-k}} = x^{\beta_u^k} *(u) c_u$

для любого $x \in \mathcal{A}$ и любого $k = 0, 1, ..., n-1$;

3) $c_u^{\alpha_u} = c_u^{\beta_u} = c_u, \quad d_u^{\alpha_u} = d_u^{\beta_u} = d_u;$

4) $\alpha_u$ и $\beta_u$ – автоморфизмы алгебры $\mathcal{A}_u$, причём $\beta_u = \alpha_u^{-1}$

*Доказательство.* 1) Так как

$$x_1 *(u) x_2^{\beta_u} *(u) ... x_i^{\beta_u^{i-1}} *(u) c_u *(u) x_{i+1}^{\alpha_u^{n-1-i}} *(u) ... x_{n-1}^{\alpha_u} *(u) x_n =$$

$$= x_1 u (u^{-1} x_2 u) u ... ((u^{-1})^{i-1} x_i u^{i-1}) u (u^{-1})^n u (u^{n-1-i} x_{i+1} (u^{-1})^{n-1-i}) u ...$$

$$... (u x_{n-1} u^{-1}) u x_n =$$

$$= x_1 u u^{-1} x_2 u^2 \ldots (u^{-1})^{i-1} x_i u^i (u^{-1})^n u^{n-i} x_{i+1} (u^{-1})^{n-1-i} \ldots u x_{n-1} u^{-1} u x_n =$$

$$= x_1 x_2 \ldots x_i x_{i+1} \ldots x_{n-1} x_n = x_1 * x_2 * \ldots x_n = \lfloor x_1 x_2 \ldots x_n \rfloor,$$

то первое тождество доказано.

Так как

$$x_1 * x_2^{\alpha_u} * \ldots x_i^{\alpha_u^{i-1}} * d_u * x_{i+1}^{\beta_u^{n-1-i}} * \ldots x_{n-1}^{\beta_u} * x_n =$$

$$= x_1 (u x_2 u^{-1}) \ldots ((u^{i-1} x_i (u^{-1})^{i-1}) u^{n-1} ((u^{-1})^{n-1-i} x_{i+1} u^{n-1-i}) \ldots (u^{-1} x_{n-1} u) x_n =$$

$$= x_1 u x_2 u^{-1} \ldots u^{i-1} x_i (u^{-1})^{i-1} u^{n-1} (u^{-1})^{n-1-i} x_{i+1} u^{n-1-i} \ldots u^{-1} x_{n-1} u x_n =$$

$$= x_1 u x_2 \ldots x_i u x_{i+1} \ldots x_{n-1} u x_n = \lfloor x_1 u x_2 u \ldots x_n \rfloor = \lfloor x_1 x_2 \ldots x_n \rfloor_u,$$

то второе тождество также доказано.

2) Так как

$$d_u * x^{\beta_u^{n-1-k}} = u^{n-1}((u^{-1})^{n-1-k} x u^{n-1-k}) = u^k x u^{n-1-k},$$

$$x^{\alpha_u^k} * d_u = (u^k x (u^{-1})^k) u^{n-1} = u^k x u^{n-1-k},$$

то первое тождество верно.

Так как

$$c_u *(u) x^{\alpha_u^{n-1-k}} = (u^{-1})^n u(n^{n-1-k} x (u^{-1})^{n-1-k}) = (u^{-1})^k x (u^{-1})^{n-1-k},$$

$$x^{\beta_u^k} *(u) c_u = ((u^{-1})^k x u^k) u (u^{-1})^n = (u^{-1})^k x (u^{-1})^{n-1-k},$$

то второе тождество также верно.

3) Очевидно.

4) Ясно, что $\alpha_u$ и $\beta_u$ – биекции. Перестановочность $\alpha_u$ и $\beta_u$ с нульарными операциями установлена в 3). Так как

$$x^{\alpha_u \beta_u} = x^{\beta_u \alpha_u} = x$$

для любого $x \in \mathcal{A}$, то операции $\alpha$ и $\beta$ перестановочны друг с другом и, кроме того, $\beta_u = \alpha_u^{-1}$. Последнее равенство позволяет ограничиться доказательством перестановочности толь-ко одной унарной операции, например α, с бинарными и n-арными операциями. Так как $\alpha_u$ – автоморфизм группы $<\mathcal{A},*>$, то

$$(x*y)^{\alpha_u} = x^{\alpha_u} * y^{\alpha_u}$$

для любых x, y ∈ $\mathcal{A}$. Так как

$$(x*(u)y)^{\alpha_u} = u(xuy)u^{-1} = uxu^{-1}uuyu^{-1} = x^{\alpha_u}uy^{\alpha_u} = x^{\alpha_u}*(u)y^{\alpha_u},$$

то

$$(x*(u)y)^{\alpha_u} = x^{\alpha_u} *(u) y^{\alpha_u}.$$

Так как n-арные операции $\lfloor \ \rfloor$ и $\lfloor \ \rfloor_u$ являются производными от операций * и *(u) соответственно, то из перестановочности α с * и *(u) вытекает

$$(\lfloor x_1 x_2 ... x_n \rfloor)^\alpha = \lfloor x_1^\alpha x_2^\alpha ... x_n^\alpha \rfloor,$$

$$(\lfloor x_1 x_2 ... x_n \rfloor_u)^\alpha = \lfloor x_1^\alpha x_2^\alpha ... x_n^\alpha \rfloor_u. \qquad \blacksquare$$

Придавая n, i и k в тождествах из 1) и 2) теоремы 1.6.14 конкретные значения, можно получить большое число новых тождеств, в том числе и аналогичных тождествам из следствия 1.5.6.

**1.6.15. Замечание.** Так как все универсальные обертывающие группы некоторой n-арной группы < A, [ ] > изоморфны универсальной обертывающей группе Поста $<\mathcal{A},*>$ этой n-арной группы, то, заменив группу $<\mathcal{A},*>$ любой другой универсальной обертывающей группой, получим утверждения, аналогичные утверждениям 1.6.9 – 1.6.14.

Теперь мы можем сформулировать основной результат, из которого вытекает, что теоремы Поста и Глускина-Хоссу являются частными случаями общего результата.

**1.6.16. Теорема.** Для всякой n-арной группы < A, [ ] > существует универсальная алгебра

$$\underline{A} = <\mathcal{A}, \Omega = \{d, c, \alpha, *, \circ, [\ ]_*, [\ ]_o\} >$$

и множества $M_*$, $N_*$, $M_o$, $N_o \subseteq \mathcal{A}$, где d, c $\in \mathcal{A}$; $<\mathcal{A},*>$ и $<\mathcal{A},\circ>$ – изоморфные группы; $<\mathcal{A}, [\ ]_*>$ и $<\mathcal{A}, [\ ]_\circ>$ – изоморфные n-арные группы; $\alpha$ – автоморфизм алгебры $\underline{A}$ такие, что выполняются следующие условия:

I. 1) $[x_1 x_2 ... x_n]_* = c \circ x_1^{\alpha^{n-1}} \circ ... x_{n-1}^{\alpha} \circ x_n,$

$[x_1 x_2 ... x_n]_\circ = x_1 * x_2^{\alpha} * ... x_n^{\alpha^{n-1}} * d$

для любых $x_1, x_2, ..., x_n \in \mathcal{A}$;

2) $x^{\alpha^{n-1}} = d*x*d^{-1} = c^{-1} \circ x \circ c$

для любого $x \in \mathcal{A}$;

II. 3) $<A, [\ ]> \simeq <M_*, [\ ]_*> \simeq <M_o, [\ ]_\circ>$;

4) $<\mathcal{A},*> = <M_*>$, $<\mathcal{A},\circ> = <M_o>$;

5) $<N_*,*> \triangleleft <\mathcal{A},*>$, $<N_o,\circ> \triangleleft <\mathcal{A},\circ>$;

6) $<\mathcal{A},*>/<N_*,*>$ и $<\mathcal{A},\circ>/<N_o,\circ>$ – циклические группы порядка n – 1.

***Доказательство.*** Для фиксированного $u \in A^{(k)} \subseteq \mathcal{A} = F_A / \theta$ положим

$d = d_u,\ c = c_u,\ \alpha = \alpha_u, \circ = *(u),\ [\ ]_* = \lfloor\ \rfloor,\ [\ ]_\circ = \lfloor\ \rfloor_u,$

$M_* = A',\ N_* = A_o,\ M_o = A^{(n-k)},\ N_o = A^{(n-k-1)}.$

1) Переписав с учетом новых обозначений первое и второе тождества из 1) теоремы 1.6.14 при i = 0 и i = n соответственно, получим требуемые тождества.

2) Снова, учитывая новые обозначения, перепишем первое и второе тождества из 2) теоремы 1.6.14 при k = n – 1 и k = 0 соответственно

$d*x = x^{\alpha^{n-1}} *d,\ c \circ x^{\alpha^{n-1}} = x \circ c,$

откуда

$$x^{\alpha^{n-1}} = d*x*d^{-1},\ x^{\alpha^{n-1}} = c^{-1} \circ x \circ c,$$

где $c^{-1}$ – обратный к c в группе $<\mathcal{A}, \circ>$, $d^{-1}$ – обратный к d в группе $<\mathcal{A}, *>$.

3) Применяется следствие 1.4.3. и 1) предложения 1.6.12.

4) – 6) Применяется теорема Поста о смежных классах и предложение 1.6.11. ∎

**1.6.17. Замечание.** Если первое и второе тождества из 1) теоремы 1.6.14 переписать при i = n и i = 0 соответственно, учтя при этом обозначения из теоремы 1.6.16, а также положив $\beta = \beta_u$, то тождества из 1) теоремы 1.6.16 можно заменить следующими тождествами:

$$[x_1 x_2 ... x_n]_* = x_1 \circ x_2^\beta \circ ... x_n^{\beta^{n-1}} \circ c,$$

$$[x_1 x_2 ... x_n]_\circ = d * x_1^{\beta^{n-1}} * ... x_{n-1}^\beta * x_n.$$

Полагая же в первом и втором тождествах из 2) теоремы 1.6.14 k = 0 и k = n – 1 соответственно, получим тождества

$$x^{\beta^{n-1}} = d^{-1} * x * d = c \circ x \circ c^{-1},$$

аналогичные тождествам из 2) теоремы 1.6.16.

**1.6.18. Замечание.** Если в теореме 1.6.16 заменить группу $<\mathcal{A}, *>$ любой универсальной обертывающей группой $\widetilde{A}$, то соответственно будем иметь

$$M_* = A,\ N_* = \widetilde{A}_\circ,\ M_\circ = \widetilde{A}^{(n-k)},\ N_\circ = \widetilde{A}^{(n-k-1)}.$$

**1.6.19. Замечание.** Если $u \in A^{(n-2)}$, то в теореме 1.6.16 будем иметь $N_\circ = A^{(n-k-1)} = A' = M_*$, т. е. $M_* = N_\circ$.

**1.6.20. Следствие.** Для всякой тернарной группы $<A, [\ ]>$ существует универсальная алгебра

$$\underline{A} = <\mathcal{A}, \Omega = \{d, c, \alpha, *, \circ, [\ ]_*, [\ ]_\circ\}>$$

и множества $A'$ и $A'' \subseteq \mathcal{A}$, где $d, c \in \mathcal{A}$; $<\mathcal{A}, *>$ и $<\mathcal{A}, \circ>$ – изоморфные группы; $<\mathcal{A}, [\ ]_*>$ и $<\mathcal{A}, [\ ]_\circ>$ – изоморфные n-арные группы; $\alpha$ – автоморфизм алгебры $\underline{A}$ такие, что выполняются следующие условия:

I. 1) $[x_1x_2x_3]_* = c \circ x_1^{\alpha^2} \circ x_2^{\alpha} \circ x_3,$

$[x_1x_2x_3]_o = x_1 * x_2^{\alpha} * x_3^{\alpha^2} * d$

для любых $x_1, x_2, x_3 \in \mathcal{A}$;

2) $x^{\alpha^2} = d*x*d^{-1} = c^{-1} \circ x \circ c$

для любого $x \in \mathcal{A}$;

II. 3) $< A, [\ ] > \simeq < A', [\ ]_* > \simeq < A'', [\ ]_o >$;

4) $< \mathcal{A}, * > = < A' >$, $< \mathcal{A}, \circ > = < A'' >$;

5) $< A'', * > \triangleleft < \mathcal{A}, * >$, $< A', \circ > \triangleleft < \mathcal{A}, \circ >$;

6) $< \mathcal{A}, * > / < A'', * >$ и $< \mathcal{A}, \circ > / < A', \circ >$ — группы второго порядка.

**1.6.21. Пример.** В примере 1.4.18 установлено, что симметрическая группа $S_n$ является обертывающей, а знакопеременная группа $A_n$ — соответствующей для тернарной группы $< T_n, [\ ] >$ нечетных подстановок.

Пусть $u = (ij)$ — произвольная транспозиция из $T_n$, операция $*$ совпадает с произведением подстановок, $s \circ t = s(ij)t$, $[str]_* = [str] = str$, $[str]_o = [s(ij)t(ij)r] = s(ij)t(ij)r$ для любых $s, t, r \in S_n$. Тогда

$$d = d_u = (ij)(ij) = e, \quad c = c_u = ((ij)^{-1})^3 = (ij)^3 = (ij),$$

$$\alpha = \beta = \alpha_u = \beta_u = x \mapsto (ij)x(ij).$$

Легко проверяется выполнимость условий 1) и 2) следствия 1.6.20, и, кроме того, условия 3) – 6) переписываются следующим образом:

$$< B_n, [\ ] > \simeq < A_n, [\ ]_o >;$$

$$< S_n, * > = < B_n >, < S_n, \circ > = < A_n >;$$

$$< A_n, * > \triangleleft < S_n, * >, < B_n, \circ > \triangleleft < S_n, \circ >;$$

$$< S_n, * > / < A_n, * > \text{ и } < S_n, \circ > / < B_n, \circ >$$

— группы второго порядка.

# ДОПОЛНЕНИЯ И КОММЕНТАРИИ

**1.** Предположим, что в n-арной (n ≥ 2) группе < A, [ ] > существует элемент e такой, что

$$[ex_1\ldots x_{n-2}x] = [xx_1\ldots x_{n-2}e] = x \qquad (*)$$

для любых $x_1, \ldots, x_{n-2}, x \in A$. Так как при n = 2 (∗) принимает вид равенств ex = xe = x, определяющих единицу группы, то может возникнуть соблазн определить с помощью (∗) еще один n-арный аналог единицы. Если n ≥ 3 и a ∈ A, то из (∗) следует

$$[e\underbrace{e\ldots e}_{n-2} e] = e, \quad [ea\underbrace{e\ldots e}_{n-3} e] = e,$$

откуда a = e, т. е. A = {e}. Таким образом, всякая n-арная группа (n ≥ 3), обладающая элементом, удовлетворяющим (∗), является одноэлементной и поэтому нет смысла вводить с помощью (∗) еще один n-арный аналог единицы.

**2.** Понятия нейтральной последовательности и единицы n-арной группы являются частными случаями следующего определения.

**Определение.** Последовательность $e_1\ldots e_{m-1}$ элементов n-арной группы < A, [ ] > называется *m-нейтральной* (n = k(m–1)+1, k≥1), если

$$[\underbrace{e_1^{m-1}\ldots e_1^{m-1}}_{i-1} x \underbrace{e_1^{m-1}\ldots e_1^{m-1}}_{k-i+1}] = x$$

для любого x ∈ A и любого i = 1, …, k+1.

Если в этом определении положить m = n, то тогда k = 1 и при i = 1 получаем $[xe_1^{n-1}] = x$, а при i = 2 получаем $[e_1^{n-1}x] = x$, т. е.

$$[e_1^{n-1}\,x] = [x\,e_1^{n-1}] = x.$$

Следовательно, n-нейтральные последовательности элементов n-арной группы – это в точности ее нейтральные последовательности.

Если в определении положить m = 2, то тогда k = n – 1 и при i = 1, …, k+1 = n получаем

$$[\underbrace{e_1\ldots e_1}_{i-1} x \underbrace{e_1\ldots e_1}_{n-i}] = x.$$

Следовательно, единицы n-арной группы – это в точности ее 2-нейтральные последовательности.

**3.** Одним определением можно объединить также понятия идемпотента и единицы n-арной группы.

**Определение.** Элемент $\varepsilon$ n-арной группы $< A, [\,] >$ называется *m-идемпотентом* (n = k(m–1)+1, k≥1), если

$$[\underbrace{\varepsilon \ldots \varepsilon}_{i-1} \, x \, \underbrace{\varepsilon \ldots \varepsilon}_{k-i+1}] = x$$

для любого $x \in A$ и любого $i = 1, \ldots, k+1$.

Легко проверяется, что n-идемпотенты n-арной группы – это в точности ее идемпотенты, а 2-идемпотенты – это в точности ее единицы. Ясно также, что элемент $\varepsilon$ n-арной группы является её m-идемпотентом тогда и только тогда, когда последовательность $\underbrace{\varepsilon \ldots \varepsilon}_{m-1}$ является m-нейтральной.

**4.** Существование бинарной операции @, отображения $\beta$ и элемента d таких, что выполнены четвертое и седьмое тождества следствия 1.5.6, было установлено независимо Л.М. Глускиным и М. Хоссу (теорема Глускина-Хоссу). Элементарное доказательство теоремы Глускина-Хоссу найдено Е.И. Соколовым [16].

**5.** Яркими примерами эффективного применения теоремы Поста являются разработанная С.А. Русаковым силовская теория n-арных групп [4], а также описание В.А. Артамоновым свободных n-арных групп [17] и шрайеровых многообразий n-арных групп [18]. В активе теоремы Глускина-Хоссу тоже немало результатов, среди которых можно отметить n-арные аналоги классических теорем Кэли и Биркгофа [19].

**6.** В связи с теоремой Глускина-Хоссу возникает естественная задача её обобщения на m-арный случай. Решение этой задачи содержится в следующих двух теоремах, в первой их которых через $\overset{m}{\theta}$ обозначено отношение эквивалентности Поста на m-арной группе $< A, [\,]_m >$.

**Теорема 1** [20, 21]. Пусть $< A, [\,]_m >$ – m-арная группа, $m \geq 2$, $n = k(m-1) + m$, $k \geq 0$. Если автоморфизм $\alpha$ m-арной группы

$< A, [\ ]_m >$ и последовательность $c_1^{t(m-1)} \in A^{t(m-1)}$ ($t \geq 1$) удовлетворяют условиям

$$c_1^\alpha c_2^\alpha \ldots c_{t(m-1)}^\alpha \overset{m}{\theta} c_1 c_2 \ldots c_{t(m-1)}, \qquad (1)$$

$$[x^{\alpha^{n-1}} c_1^{t(m-1)}]_m = [c_1^{t(m-1)} x]_m, \ x \in A, \qquad (2)$$

то $< A, [\ ]_n >$ – n-арная группа с n-арной операцией

$$[x_1 x_2 \ldots x_n]_n = [x_1 x_2^\alpha \ldots x_n^{\alpha^{n-1}} c_1^{t(m-1)}]_m, \qquad (3)$$

причем $\alpha$ становится автоморфизмом n-арной группы $< A, [\ ]_n >$.

**Теорема 2** [20, 21]. Пусть $n = k(m-1) + m$, $\nu = \mu(k+1) - 1$, $m \geq 2$, $\mu \geq 1$; $k \geq 1$ при $\mu = 1$; $k \geq 0$ при $\mu > 1$. На n-арной группе $< A, [\ ]_n >$ определим m-арную операцию

$$[x_1 x_2 \ldots x_m]_m = [x_1 a_1^\nu x_2 a_1^\nu \ldots a_1^\nu x_m]_n$$

и отображение

$$\alpha : x \mapsto x^\alpha = [b_1^l x a_1^\nu]_n,$$

где $a_1, \ldots, a_\nu \in A$, $b_1^l = b_1^{\mu(n-k-2)+1} = b_1^{\mu(k+1)(m-2)+1}$ – обратная последовательность для последовательности $a_1^\nu$. Тогда $< A, [\ ]_m >$ – m-арная группа, $\alpha$ – автоморфизм $< A, [\ ]_n >$ и $< A, [\ ]_m >$ и выполняются условия (1) – (3) для любой последовательности $c_1^{t(m-1)}$, эквивалентной в m-арной группе $< A, [\ ]_m >$ последовательности $d_1^{m-1}$, где

$$d_1 = [\underbrace{b_1^l \ldots b_1^l}_{n-2} b_1 b_1^l]_n,\ d_2 = [b_2^{\mu(k+1)+1} b_1^l]_n, \ldots,$$

$$d_{m-1} = [b_{(m-3)\mu(k+1)+2}^{(m-2)\mu(k+1)+1} b_1^l]_n.$$

Теоремы 1 и 2 охватывают теорему Глускина-Хоссу (m = 2), а также позволяют сформулировать ряд следствий, среди которых находятся результаты В. Дудека и Я. Михальского [22], а также И.И. Дериенко и О.В. Колесникова [23] о приводимости n-арных групп к m-арным группам.

**7.** Теоремам Поста и Глускина-Хоссу посвящена книга автора [24].

# ГЛАВА 2

# n-АРНЫЕ АНАЛОГИ НЕКОТОРЫХ БИНАРНЫХ ПОНЯТИЙ

Известно, что при переходе от групп к n-арным группам одно и то же бинарное понятие может иметь несколько n-арных аналогов. В данной главе определяются и изучаются различные n-арные аналоги инвариантных и сопряженных подгрупп, циклических и абелевых групп, прямых и декартовых произведений групп.

## §2.1. n-АРНЫЕ ПОДГРУППЫ. СМЕЖНЫЕ КЛАССЫ

**2.1.1. Определение.** Подалгебра $<B, [\,]>$ n-арной группы $<A, [\,]>$ называется её *n-арной подгруппой*, если она сама является n-арной группой.

n-Арная подгруппа $<B, [\,]>$ n-арной группы $<A, [\,]>$ называется *собственной*, если $B \subset A$.

Следующее утверждение является очевидным.

**2.1.2. Предложение.** Для того, чтобы подалгебра $<B, [\,]>$ n-арной группы $<A, [\,]>$ была её n-арной подгруппой, достаточно чтобы в B были разрешимы уравнения

$$[xb_2\ldots b_n] = c, \quad [b_1\ldots b_{n-1}y] = c$$

для всех $b_1\ldots b_n, c \in B$.

**2.1.3. Предложение.** Любая конечная подалгебра n-арной группы является её n-арной подгруппой.

*Доказательство.* Пусть $< B, [\ ] >$ – конечная подалгебра n-арной группы $< A, [\ ] >$. Покажем, что для любых $b_1,\ldots,b_{n-1},\ c \in B$ в B разрешимо уравнение

$$[xb_1\ldots b_{n-1}] = c.$$

Для этого рассмотрим множество

$$[Bb_1\ldots b_{n-1}] = \{[bb_1\ldots b_{n-1}] \mid b \in B\}.$$

Так как, согласно определению n-арной группы,

$$[bb_1\ldots b_{n-1}] \neq [b'b_1\ldots b_{n-1}]$$

для любых $b, b' \in B$ $(b \neq b')$, то $|[Bb_1\ldots b_{n-1}]| = |B|$, откуда, учитывая очевидное включение $[Bb_1\ldots b_{n-1}] \subseteq B$ и конечность множества B, получаем $[Bb_1\ldots b_{n-1}] = B$. Из последнего равенства вытекает разрешимость записанного выше уравнения. Аналогично доказывается разрешимость в B уравнения

$$[b_1\ldots b_{n-1}y] = c$$

для любых $b_1,\ldots,b_{n-1}, c \in A$. ∎

**2.1.4. Теорема** [1]. Для того, чтобы подалгебра $< B, [\ ] >$ n-арной группы $< A, [\ ] >$ была её n-арной подгруппой ($n \geq 3$), необходимо и достаточно, чтобы множество B вместе со всяким своим элементом b содержало и косой элемент $\bar{b}$.

**2.1.5. Пример.** Пусть $< B, [\ ] >$ – n-арная группа из примера 1.1.6, B – любая подгруппа группы A, содержащая $a \in Z(A)$. Так как $a \in Z(B)$, то $< B, [\ ] >$ – n-арная группа, а значит и n-арная подгруппа в $< A, [\ ] >$.

**2.1.6. Пример.** Пусть $< A, [\ ] >$ – n-арная группа, производная от группы A (пример 1.1.7), B – подгруппа группы A. Положив в примере 1.15 $a = 1$, получим n-арную подгруппу $< B, [\ ] >$ n-арной группы $< A, [\ ] >$.

Следующий пример показывает, что n-арные группы, производные от групп, могут содержать n-арные подгруппы, не являющиеся производными от групп.

**2.1.7. Пример.** Пусть $< S_n, [\ ] >$ – тернарная группа, производная от симметрической группы $S_n$. Ясно, что тернарная группа $< T_n, [\ ] >$ нечётных подстановок (пример 1.1.8) является тернарной подгруппой в $< S_n, [\ ] >$ и не является производной от группы. Отметим, что ввиду примера 1.1.6, тернарной подгруппой в $< S_n, [\ ] >$ будет и тернарная группа $< A_n, [\ ] >$ четных подстановок.

**2.1.8. Пример.** Если b – идемпотент n-арной группы $< A, [\ ] >$, то $< \{b\}, [\ ] >$ – её n-арная подгруппа, состоящая из одного элемента.

**2.1.9. Пример**. Пусть $< E(2), [\ ] >$ – тернарная группа, производная от группы $E(2)$ всех движений плоскости. Ясно, что тернарная группа $< E_2(2), [\ ] >$ всех движений второго рода (пример 1.1.14) является тернарной подгруппой в $< E(2), [\ ] >$. Тернарными подгруппами в $< E(2), [\ ] >$ являются и тернарная группа $< E_1(2), [\ ] >$ всех движений первого рода, а также тернарные группы из примеров 1.1.15 – 1.1.17.

Так как $A_n \cap B_n = \varnothing$, $E_1(2) \cap E_2(2) = \varnothing$, то примеры 2.1.7 и 2.1.9 показывают, что пересечение n-арных подгрупп в n-арной группе может быть пустым, что невозможно в группе.

Следующая теорема описывает строение тернарной группы $< B_n, [\ ] >$ отражений правильного n-угольника (пример 1.2.6).

**2.1.10. Теорема** [25, 26]**.** Для всякого делителя k натурального числа n существует точно m, где n = mk, тернарных подгрупп $< H_1, [\ ] >,\ldots,< H_m, [\ ] >$ порядка k тернарной группы $< B_n, [\ ] >$. Причём

$$B_n = \bigcup_{i=1}^{m} H_i, \quad H_i \cap H_j = \varnothing \ (i \neq j)$$

Следующее предложение даёт явное описание всех тернарных подгрупп тернарной группы $< B_n, [\ ] >$.

**2.1.11. Предложение** [26]**.** Тернарные подгруппы k-го порядка тернарной группы $< B_n, [\ ] >$, где $n = mk$, исчерпываются подгруппами

$$B_i^{(k)} = \{b_i, b_{m+i}, b_{2m+i}, \ldots, b_{(k-1)m+i}\},\ i = 1, 2, \ldots, m.$$

**2.1.12. Предложение.** Для любого семейства

$$\{< B_i, [\ ] > \mid i \in I\},\ B = \bigcap_{i \in I} B_i \neq \varnothing$$

n-арных подгрупп n-арной группы $< A, [\ ] >$, множество B замкнуто относительно n-арной операции $[\ ]$, а алгебра $< B, [\ ] >$ является n-арной подгруппой в $< A, [\ ] >$.

**2.1.13. Определение.** Если $< A, [\ ] >$ – n-арная группа, $M \subseteq A$, то пересечение всех n-арных подгрупп n-арной группы $< A, [\ ] >$, содержащих множество M, называется n-арной подгруппой, *порождённой множеством M*, и обозначается через $<< M >, [\ ] >$. Множество M при этом называется порождающим для n-арной группы $<< M >, [\ ] >$.

Таким образом, $< M > = \bigcap_{M \subseteq B} B$, где $< B, [\ ] >$ – n-арная подгруппа в $< A, [\ ] >$.

В теории групп группу $< M >$ можно определить с помощью элементов множества M. Аналогичный результат имеет место для n-арных групп.

Если $< A, [\ ] >$ – n-арная группа, $M \subseteq A$, $M \neq \varnothing$, то положим $\overline{M} = \{\ \overline{a}\ \mid a \in M\ \}$.

В следующей теореме для любого элемента a n-арной группы $< A, [\ ] >$ будем считать $[a] = a$.

**2.1.14. Теорема** [12]**.** Если $< A, [\ ] >$ – n-арная группа $(n \geq 3)$, $M \subseteq A$, $M \neq \varnothing$, то

$$< M > = \{[a_1 \ldots a_{k(n-1)+1}] \mid a_i \in M \cup \overline{M},\ i = 1, \ldots, k(n-1)+1,\ k = 0, 1, \ldots\}.$$

Пусть $< B, [\ ] > -$ n-арная подгруппа n-арной группы $< A, [\ ] >$, $\alpha$ и $\gamma$ – последовательности элементов множества A такие, что $l(\alpha) \geq 0$, $l(\gamma) \geq 0$ и $l(\alpha) + l(\gamma) + k = n$, где $k \in \{1, \ldots, n-1\}$. Положим

$$[\alpha \overset{k}{B} \gamma] = [\alpha \underbrace{B \ldots B}_{k} \gamma] = \{[\alpha b_1 \ldots b_k \gamma] \mid b_1, \ldots, b_k \in B\}.$$

**2.1.15. Теорема** [12]**.** Справедливы следующие утверждения:

1) $[\alpha \overset{k}{B} \gamma] \cap [\beta \overset{k}{B} \gamma] = \varnothing$ или $[\alpha \overset{k}{B} \gamma] = [\beta \overset{k}{B} \gamma]$,

$[\alpha \overset{k}{B} \gamma] \cap [\alpha \overset{k}{B} \delta] = \varnothing$ или $[\alpha \overset{k}{B} \gamma] = [\alpha \overset{k}{B} \delta]$;

2) $[\alpha \overset{k}{B} \gamma] = [\alpha b_1 \ldots b_i B\, b_{i+1} \ldots b_{k-1} \gamma]$

для фиксированных $b_1, \ldots, b_i, b_{i+1}, \ldots, b_{k-1} \in B$ и любого $i = 0, 1, \ldots, k-1$;

3) мощности множеств $[\alpha \overset{k}{B} \gamma]$ и B совпадают;

4) $A = \bigcup\limits_{u \in S} U = \bigcup\limits_{v \in T} V$,

где $S = \{[\alpha \overset{k}{B} \gamma] \mid \alpha \in F_A, l(\alpha) \neq 0, l(\alpha)$ и $\gamma -$ фиксированные$\}$,

$T = \{[\alpha \overset{k}{B} \gamma] \mid \gamma \in F_A, l(\gamma) \neq 0, l(\gamma)$ и $\alpha -$ фиксированные$\}$.

**2.1.16. Замечание.** Из теоремы 2.1.15 вытекает, что множество A можно представить в виде объединения различных попарно непересекающихся множеств вида $[\alpha \overset{k}{B} \gamma]$ для фиксированной последовательности $\gamma$. Если

$$A = [\alpha \overset{k}{B} \gamma] \cup [\beta \overset{k}{B} \gamma] \cup \ldots$$

такое объединение, то оно называется левым $(\overset{k}{B}, \gamma)$-разложением n-арной группы $< A, [\ ] >$ по $< B, [\ ] >$, а множество $\{\alpha, \beta, \ldots\}$ называется множеством представителей этого левого разложения.

Одно и то же левое разложение может иметь различные множества представителей. Однако любое множество L представителей $(\overset{k}{B}, \gamma)$-разложения n-арной группы $<A, [\ ]>$ обладает следующими свойствами:

1) для любого $a \in A$ существует $\alpha \in L$ такая, что $a \in [\alpha \overset{k}{B} \gamma]$;

2) $[\alpha \overset{k}{B} \gamma] \cap [\beta \overset{k}{B} \gamma] = \varnothing$ для любых $\alpha, \beta \in L, \alpha \neq \beta$.

Из теоремы 2.1.15 вытекает также, что множество A можно представить в виде объединения различных попарно непересекающихся множеств вида $[\alpha \overset{k}{B} \gamma]$ для фиксированной последовательности $\alpha$. Если

$$A = [\alpha \overset{k}{B} \gamma] \cup [\alpha \overset{k}{B} \delta] \cup \ldots$$

такое объединение, то оно называется правым $(\alpha, \overset{k}{B})$-разложением n-арной группы $<A, [\ ]>$ по $<B, [\ ]>$, а множество $\{\gamma, \delta, \ldots\}$ называется множеством представителей этого правого разложения. Любое множество R представителей $(\alpha, \overset{k}{B})$-разложения обладает следующими свойствами:

3) для любого $a \in A$ существует $\gamma \in R$ такая, что $a \in [\alpha \overset{k}{B} \gamma]$;

4) $[\alpha \overset{k}{B} \gamma] \cap [\alpha \overset{k}{B} \delta] = \varnothing$ для любых $\gamma, \delta \in R, \gamma \neq \delta$.

Если $k = n - 1, \alpha = a \in A, \gamma = \varnothing$, то

$$[\alpha \overset{k}{B} \gamma] = [a \overset{n-1}{B}].$$

Если же $k = n - 1, \alpha = \varnothing, \gamma = a \in A$, то

$$[\alpha \overset{k}{B} \gamma] = [\overset{n-1}{B} a].$$

**2.1.17. Определение.** Множества вида $[a \overset{n-1}{B}]$ и $[\overset{n-1}{B} a]$ называются соответственно *левыми* и *правыми смежными классами* n-арной группы $<A, [\ ]>$ по n-арной подгруппе $<B, [\ ]>$; $(\overset{n-1}{B}, \varnothing)$-разложение

$$A = [a \overset{n-1}{B}] \cup [b \overset{n-1}{B}] \cup \ldots$$

n-арной группы < A, [ ] > называется её *разложением на левые смежные классы*, ($\varnothing$, $\overset{n-1}{B}$)-разложение

$$A = [\overset{n-1}{B} a] \cup [\overset{n-1}{B} b] \cup \ldots$$

n-арной группы < A, [ ] > называется её *разложением на правые смежные классы*.

Так как

$$[b \overset{n-1}{B}] = [\overset{n-1}{B} b] = B$$

для любого $b \in B$, то любое разложение n-арной группы на левые или правые смежные классы содержит множество B.

**2.1.18. Теорема** [12]. Мощность множества представителей любого левого разложения n-арной группы совпадает с мощностью множества представителей любого правого разложения этой же n-арной группы.

Согласно теореме 2.1.18, множества представителей любых разложений n-арной группы по одной и той же n-арной подгруппе имеют одинаковую мощность. Поэтому естественно следующее

**2.1.19. Определение.** Мощность множества представителей любого разложения n-арной группы < A, [ ] > по её n-арной подгруппе < B, [ ] > называется *индексом* < B, [ ] > в < A, [ ] > и обозначается через |A : B|.

**2.1.20. Теорема Лагранжа для n-арных групп** [3]. Пусть < A, [ ] > – конечная n-арная группа, < B, [ ] > – её n-арная подгруппа. Тогда  |A| = |B| · |A : B|.

*Доказательство.* Пусть

$$A = [a_1 \overset{n-1}{B}] + \ldots + [a_k \overset{n-1}{B}]$$

– разложение < A, [ ] > на левые смежные классы, т. е. |A : B| = k. Так как согласно 3) теоремы 2.1.15,

$$| [a_1 \overset{n-1}{B} ] | = \ldots = | [a_k \overset{n-1}{B} ] | = |B|,$$

то $|A| = |B| \cdot k = |B| \cdot |A : B|$. ∎

**2.1.21. Лемма.** Если $< X, [\ ] > -$ n-арная подгруппа n-арной группы $< A, [\ ] >$; $a, b \in A$, то

$$[a \overset{n-1}{X} ] = [b \overset{n-1}{X} ]$$

тогда и только тогда, когда $b \in [a \overset{n-1}{X} ]$.

*Доказательство.* Если

$$[a \overset{n-1}{X} ] = [b \overset{n-1}{X} ],$$

то

$$b = [be_1\ldots e_{n-1}] \in [b \overset{n-1}{X} ] = [a \overset{n-1}{X} ],$$

где $e_1\ldots e_{n-1}$ – нейтральная последовательность, составленная из элементов множества X.

Обратно, пусть $b \in [a \overset{n-1}{X} ]$. Тогда $b = [ax_1\ldots x_{n-1}]$, где $x_1, \ldots, x_{n-1} \in X$, откуда

$$[b \overset{n-1}{X} ] = [[ax_1\ldots x_{n-1}] \overset{n-1}{X} ] = [a \overset{n-1}{X} ]. \quad ∎$$

**2.1.22. Предложение.** Если

$$A = \bigcup_{x \in X}[x \overset{n-1}{B} ] \text{ и } B = \bigcup_{y \in Y}[y \overset{n-1}{C} ]$$

– разложения n-арных групп $< A, [\ ] >$ и $< B, [\ ] >$ на непересекающиеся левые смежные классы по n-арным подгруппам $< B, [\ ] >$ и $< C, [\ ] >$ соответственно, и $b_1, \ldots, b_{n-2}$ – фиксированные элементы из B, то справедливы следующие утверждения:

1) $[x \overset{n-1}{B} ] = \bigcup_{y \in Y}[[xb_1\ldots b_{n-2}y] \overset{n-1}{C} ]$

для любого x ∈ A, причем любые два смежных класса в правой части не совпадают;

2) $A = \bigcup\limits_{x \in X,\ y \in Y} [[xb_1 \ldots b_{n-2} y] \overset{n-1}{C}]$

– разложение $< A, [\ ] >$ на непересекающиеся левые смежные классы по $< C, [\ ] >$.

***Доказательство.*** 1) $[x \overset{n-1}{B}] = [xb_1 \ldots b_{n-2} B] =$

$= [xb_1 \ldots b_{n-2} \bigcup\limits_{y \in Y} [y \overset{n-1}{C}]] = \bigcup\limits_{y \in Y} [xb_1 \ldots b_{n-2} [y \overset{n-1}{C}]] =$

$= \bigcup\limits_{y \in Y} [[xb_1 \ldots b_{n-2} y] \overset{n-1}{C}]$.

Предположим, что

$$[[xb_1 \ldots b_{n-2} y_1] \overset{n-1}{C}] = [[xb_1 \ldots b_{n-2} y_2] \overset{n-1}{C}],$$

где $y_1, y_2 \in Y$, $y_1 \neq y_2$. Тогда по лемме 2.1.21,

$$[xb_1 \ldots b_{n-2} y_2] \in [[xb_1 \ldots b_{n-2} y_1] \overset{n-1}{C}],$$

$$[xb_1 \ldots b_{n-2} y_2] \in [xb_1 \ldots b_{n-2} [y_1 \overset{n-1}{C}]],$$

откуда

$$y_2 \in [y_1 \overset{n-1}{C}].$$

Применяя снова лемму 2.1.21, получаем

$$[y_1 \overset{n-1}{C}] = [y_2 \overset{n-1}{C}],$$

что противоречит выбору $y_1$ и $y_2$. Следовательно,

$$[[xb_1 \ldots b_{n-2} y_1] \overset{n-1}{C}] \cap [[xb_1 \ldots b_{n-2} y_2] \overset{n-1}{C}] = \varnothing.$$

2) Применяем 1), а также разложение $<A, [\ ]>$ по $<B, [\ ])>$:

$$A = \bigcup_{x \in X}[x \overset{n-1}{B}] = \bigcup_{x \in X}\bigcup_{y \in Y}[[xb_1\ldots b_{n-2}y]\overset{n-1}{C}] =$$

$$= \bigcup_{x \in X,\ y \in Y}[[xb_1\ldots b_{n-2}y]\overset{n-1}{C}]. \qquad \blacksquare$$

## §2.2. СВЯЗЬ МЕЖДУ n-АРНЫМИ ПОДГРУППАМИ n-АРНОЙ ГРУППЫ $<A, [\ ]>$ И ПОДГРУППАМИ ГРУПП $<A, @>$, $<\mathcal{A}, *>$ И $<A_0, *>$

**2.2.1 Определение.** Пусть $<A, [\ ]>$ – n-арная группа, a – фиксированный элемент из A. Вместе со всяким подмножеством $B \subseteq A$ свяжем следующие два подмножества:

$$B_a = [\overset{n-1}{B} a] = \{[b_1\ldots b_{n-1}a] \mid b_1, \ldots, b_{n-1} \in A\}; \qquad (1)$$

$$_aB = [a \overset{n-1}{B}] = \{[ab_1\ldots b_{n-1}] \mid b_1, \ldots, b_{n-1} \in A\}. \qquad (2)$$

Ясно, что если $<B, [\ ]>$ – n-арная подгруппа в $<A, [\ ]>$, то $B_a$ и $_aB$ – смежные классы $<A, [\ ]>$ по $<B, [\ ]>$.

Из определения также вытекает, что $A_a = {}_aA = A$ для любого $a \in A$.

**2.2.2. Замечание.** Используя понятие эквивалентности последовательностей, можно показать, что если $<B, [\ ]>$ – n-арная подгруппа, то в (1) и (2) достаточно чтобы только один из элементов $b_1, \ldots, b_{n-1}$ пробегал всё множество B, а остальные могут быть фиксированными.

Через $<A, @>$ обозначается группа с операцией

$$x @ y = [x\alpha y],$$

где $\alpha$ – обратная последовательность для a (см. §1.5).

**2.2.3. Теорема** [12]**.** Если $<B,[\ ]>$ – n-арная подгруппа n-арной группы $<A,[\ ]>$, то $<B_a,@>$ и $<{}_aB,@>$ – изоморфные подгруппы группы $<A,@>$. При этом a – единица в $<A,@>$, обратными для элементов $[b_1\ldots b_{n-1}a] \in B_a$, $[ab_1\ldots b_{n-1}] \in {}_aB$ являются соответственно элементы $[\beta a] \in B_a$, $[a\beta] \in {}_aB$, где $\beta$ – обратная последовательность для последовательности $b_1\ldots b_{n-1}$, для всякого элемента x группы $<A,@>$ обратный элемент имеет вид $x^{-1} = [a\widetilde{x}a]$, где $\widetilde{x}$ – обратная последовательность для элемента x.

**2.2.4. Определение.** Если $<A,[\ ]>$ – n-арная группа, a – фиксированный элемент из A, $\alpha$ – обратная последовательность для a, $B \subseteq A$, то положим

$$\widetilde{B} = \{[ab\alpha] \mid b \in B\}, \quad \widehat{B} = \{[\alpha ba] \mid b \in B\}.$$

Легко заметить, что $\widetilde{\widehat{B}} = \widehat{\widetilde{B}} = B$.

**2.2.5. Теорема** [12]**.** Если $<A,[\ ]>$ – n-арная группа, $B \subseteq A$, то следующие утверждения равносильны:

1) $<B,@>$ – подгруппа в $<A,@>$;

2) $<\widetilde{B},@>$ – подгруппа в $<A,@>$;

3) $<\widehat{B},@>$ – подгруппа в $<A,@>$.

**2.2.6. Теорема** [12]**.** Пусть $<B,@>$ – подгруппа группы $<A,@>$ и существует элемент $x \in A$ такой, что:

1) $[\underbrace{x\ldots x}_{n-1}a] \in B$;

2) $B @ x = x @ \widetilde{B}$.

Тогда $<B @ x, [\ ]>$ – n-арная подгруппа n-арной группы $<A,[\ ]>$ такая, что $H_a = B$, где $H = B @ x = x @ \widetilde{B}$.

Полагая в теореме 2.2.6 $x = a$, получим

**2.2.7. Следствие.** Если $<B,@>$ – подгруппа группы $<A,@>$ такая, что

$$[\underbrace{a\ldots a}_{n}] \in B, \quad B = [aB\alpha],$$

то $<B, [\,]>$ – n-арная подгруппа n-арной группы $<A, [\,]>$.

**2.2.8. Следствие.** Если $<B,@>$ – характеристическая подгруппа группы $<A,@>$ такая, что $[\underbrace{a\ldots a}_{n}] \in B$, то $<B, [\,]>$ – n-арная подгруппа n-арной группы $<A, [\,]>$.

*Доказательство.* По предложению 1.5.4 отображение $x \mapsto [ax\alpha]$ является автоморфизмом группы $<A,@>$. А так как $<B,@>$ характеристична в $<A,@>$, то $B = [aB\alpha]$. Так как выполняются все условия предыдущего следствия, то $<B, [\,]>$ – n-арная подгруппа в $<A, [\,]>$. ∎

**2.2.9. Замечание.** В качестве обратной последовательности $\alpha$ для идемпотента a можно взять последовательность $\underbrace{a\ldots a}_{n-2}$.

**2.2.10. Следствие.** Если a – идемпотент n-арной группы $<A, [\,]>$, и для подгруппы $<B,@>$ группы $<A,@>$ выполняется условие

$$B = [aB\underbrace{a\ldots a}_{n-2}],$$

то $<B, [\,]>$ – n-арная подгруппа в $<A, [\,]>$.

Комбинируя следствия 2.2.8 и 2.2.10, получаем

**2.2.11. Следствие.** Если a – идемпотент n-арной группы $<A, [\,]>$, $<B,@>$ – характеристическая подгруппа в $<A,@>$, то $<B, [\,]>$ – n-арная подгруппа в $<A, [\,]>$.

Зафиксируем элемент x n-арной группы $<A, [\,]>$ и введём следующие обозначения: $L(A, [\,], x)$ – множество всех n-арных подгрупп n-арной группы $<A, [\,]>$, содержащих x;

L(A, @ , x) – множество всех подгрупп группы $<A,@>$, для которых выполняются условия 1) и 2) теоремы 2.2.6.

**2.2.12. Теорема.** Существует биекция

$$L(A, @, x) \leftrightarrow L(A, [\ ], x).$$

**2.2.13. Предложение** [22]**.** Если $<A, [\ ]> = <A, [\ ]_{\circ,\alpha,d}>$ – n-арная группа, определяемая для группы $<A, \circ>$ обратной теоремой Глускина-Хоссу, то $<A, \circ> = <A, ⓔ>$, где e – единица группы $<A, \circ>$, т. е. операции $\circ$ и ⓔ совпадают.

*Доказательство.* По условию

$$[a_1 a_2 \ldots a_n] = a_1 \circ a_2^{\beta} \circ \ldots a_n^{\beta^{n-1}} \circ d,$$

$$d^{\beta} = d, \quad a^{\beta^{n-1}} = d \circ a \circ d^{-1}$$

для любых $a_1, a_2, \ldots, a_n, a \in A$. Так как

$$[a\underbrace{e\ldots e}_{n-3} d^{-1}b] = a \circ e^{\beta} \circ \ldots e^{\beta^{n-3}} \circ (d^{-1})^{\beta^{n-2}} \circ b^{\beta^{n-1}} \circ d =$$

$$= a \underbrace{e \circ \ldots e}_{n-3} \circ d^{-1} \circ d \circ b \circ d^{-1} \circ d = a \circ b$$

для любых $a, b \in A$, то

$$[a\underbrace{e\ldots e}_{n-3} d^{-1}b] = a \circ b \qquad (1)$$

Полагая в (1) $b = e$, получим

$$[a\underbrace{e\ldots e}_{n-3} d^{-1}e] = a \circ e = a$$

для любого $a \in A$, т. е.

$$\underbrace{e\ldots e}_{n-3} d^{-1}$$

– обратная последовательность для e. Поэтому из (1) вытекает

$$a \circledcirc b = a \circ b.$$

т. е. операция $\circledcirc$ и $\circ$ совпадают. ∎

**2.2.14. Предложение.** Пусть $< A, [\ ] >$ – n-арная группа, $a, b \in A$, $\gamma$ и $\delta$ – последовательности, составленные из элементов множества A, такие, что $\gamma\delta$ – обратная последовательность для a. Тогда преобразования

$$\varphi : x \mapsto [b\delta x\gamma], \quad \psi : x \mapsto [\delta x\gamma b]$$

являются изоморфизмами группы $< A, \circledast >$ на группу $< A, \circledcirc >$.

***Доказательство.*** Ясно, что $\varphi$ – биекция. Если $\beta$ – обратная последовательность для b, то

$$\varphi(x \circledast y) = [b\delta(x \circledast y)\gamma] = [b\delta[x\gamma\delta y]\gamma] =$$

$$= [[b\delta x\gamma]\beta[b\delta y\gamma]] = \varphi(x) \circledcirc \varphi(y),$$

$$\psi(x \circledast y) = [\delta(x \circledast y)\gamma b] = [\delta[x\gamma\delta y]\gamma b] =$$

$$= [[\delta x\gamma b]\beta[\delta y\gamma b]] = \psi(x) \circledcirc \psi(y),$$

т. е.

$$\varphi(x \circledast y) = \varphi(x) \circledcirc \varphi(y), \psi(x \circledast y) = \psi(x) \circledcirc \psi(y).$$

Следовательно $\varphi$ и $\psi$ - изоморфизмы $< A, \circledast >$ на $< A, \circledcirc >$. ∎

**2.2.15. Следствие.** Если $< A, [\ ] > $ – n-арная группа, $a, b \in A$, то следующие преобразования

$$x \mapsto [b\bar{a}\underbrace{a\ldots a}_{n-3}x], \ x \mapsto [bx\bar{a}\underbrace{a\ldots a}_{n-3}],$$

$$x \mapsto [x\bar{a}\underbrace{a\ldots a}_{n-3}b], \quad x \mapsto [\bar{a}\underbrace{a\ldots a}_{n-3}xb]$$

являются изоморфизмами группы $< A, @ >$ на группу $< A, ⓑ >$.

**2.2.16. Следствие.** Если $< A, [\ ] >$ – n-арная группа, a – её идемпотент, $b \in A$, то следующие преобразования

$$x \mapsto [b\underbrace{a\ldots a}_{n-2}x], \quad x \mapsto [bx\underbrace{a\ldots a}_{n-2}],$$

$$x \mapsto [x\underbrace{a\ldots a}_{n-2}b], \quad x \mapsto [\underbrace{a\ldots a}_{n-2}x\,b]$$

являются изоморфизмами группы $< A, @ >$ на группу $< A, ⓑ >$.

**2.2.17. Предложение.** Если $< A, [\ ] >$ и $< B, [\ ] >$ – изоморфные n-арные группы, то для любых $a \in A$, $b \in B$ группы $< A, @ >$ и $< B, ⓑ >$ также изоморфны.

*Доказательство.* Так как для любых $b, c \in B$ группы $< B, © >$ и $< B, ⓑ >$ изоморфны (предложение 2.2.14), то достаточно доказать изоморфизм групп $< A, @ >$ и $< B, © >$, где $c = a^\varphi$, $\varphi$ – изоморфизм $< A, [\ ] >$ на $< B, [\ ] >$.

Так как

$$[b\bar{a}^\varphi \underbrace{a^\varphi \ldots a^\varphi}_{n-3} a^\varphi] = [(b^{\varphi^{-1}})^\varphi \bar{a}^\varphi \underbrace{a^\varphi \ldots a^\varphi}_{n-3} a^\varphi] =$$

$$= [b^{\varphi^{-1}} \bar{a} \underbrace{a\ldots a}_{n-3} a]^\varphi = (b^{\varphi^{-1}})^\varphi = b$$

для любого $b \in B$, то $\bar{a}^\varphi \underbrace{a^\varphi \ldots a^\varphi}_{n-3}$ – обратная последовательность для $a^\varphi$. Поэтому

$$(x @ y)^\varphi = [x\bar{a}\underbrace{a\ldots a}_{n-3}y]^\varphi = [x^\varphi \bar{a}^\varphi \underbrace{a^\varphi \ldots a^\varphi}_{n-3} y^\varphi] = x^\varphi ⓐ^\varphi y^\varphi = x^\varphi © y^\varphi.$$

т. е. $\varphi$ – изоморфизм $< A, @ >$ на $< B, © >$. ∎

**2.2.18. Определение.** Для всякого подмножества B n-арной группы $< A, [\ ] >$ положим

$$B^{(i)}(A) = \{\theta_A(\alpha) \in A^{(i)} \mid \exists b_1, \ldots, b_i \in B, \alpha \theta_A b_1 \ldots b_i\}, i = 1, \ldots, n-1;$$

$$B_o(A) = B^{(n-1)}(A) = \{\theta_A(\alpha) \in A_o \mid \exists b_1, \ldots, b_{n-1} \in B, \alpha \theta_A b_1 \ldots b_{n-1}\};$$

$$B^*(A) = \{\theta_A(\alpha) \in A^* \mid \exists b_1, \ldots, b_i \in B\ (i \geq 1), \alpha \theta_A b_1 \ldots b_i\}.$$

Ясно, что

$$B^{(i)}(A) \subseteq A^{(i)}, B^*(A) \subseteq A^*, B_o(A) \subseteq A_o,$$

в частности,

$$A^{(i)}(A) = A^{(i)}, A^*(A) = A^*, A_o(A) = A_o.$$

**2.2.19. Теорема.** Справедливы следующие утверждения:
1) $B^*(A)$ – подполугруппа группы $A^*$;
2) если $< B, [\ ] >$ – n-арная подгруппа n-арной группы $< A, [\ ] >$, то $B^*(A)$ – подгруппа группы $A^*$, изоморфная группе $B^*$, причём:

$$B^*(A) = \bigcup_{i=1}^{n-1} \{\theta_A(bb_1 \ldots b_{i-1}) \mid b \in B\}, \quad (1)$$

где $b_1, \ldots, b_{n-2}$ – фиксированные элементы из B;

$$B^*(A) = \bigcup_{i=1}^{n-1} B^{(i)}(A). \quad (2)$$

*Доказательство.* 1) Пусть $\theta_A(\alpha)$ и $\theta_A(\beta)$ – произвольные элементы из $B^*(A)$, т.е.

$$\theta_A(\alpha) = \theta_A(b_1 \ldots b_i), \theta_A(\beta) = \theta_A(b'_1 \ldots b'_j),$$

где $b_1, \ldots, b_i, b'_1, \ldots, b'_j \in B$. Так как $\theta_A$ конгруэнции на $F_A$, то

$$\theta_A(\alpha)\,\theta_A(\beta) = \theta_A(\alpha\beta) = \theta_A(b_1 \ldots b_i b'_1 \ldots b'_j) \in B^*(A).$$

Следовательно, $B^*(A)$ – подполугруппа группы $A^*$.

2) Так как $<B, [\ ]>$ – n-арная подгруппа в $<A, [\ ]>$, то существует нейтральная последовательность $e_1 \ldots e_{n-1}$, составленная из элементов множества B. Поэтому

$$E = \theta_A(e_1 \ldots e_{n-1}) \in B^*(A),$$

а при доказательстве теоремы 1.4.2 установлено, что E – единица группы $A^*$.

Если $\theta_A(\alpha)$ – произвольный элемент из $B^*(A)$, то $\alpha\theta_A\beta$ для некоторой последовательности $\beta$, составленной из элементов множества B. Так как $<B, [\ ]>$ – n-арная подгруппа в $<A, [\ ]>$, то для $\beta$ существует обратная последовательность $\gamma$, составленная из элементов множества B. Тогда

$$\theta_A(\alpha)\theta_A(\gamma) = \theta_A(\beta)\theta_A(\gamma) = \theta_A(\beta\gamma) = E =$$
$$= \theta_A(\gamma\beta) = \theta_A(\gamma)\theta_A(\beta) = \theta_A(\gamma)\theta_A(\alpha),$$

т.е. $\theta_A(\gamma)$ – обратный для $\theta_A(\alpha)$. Таким образом, доказано, что $B^*(A)$ – группа.

Если $b_1, \ldots, b_{n-2} \in B$, то включение

$$\bigcup_{i=1}^{n-1} \{\theta_A(bb_1 \ldots b_{i-1}) \mid b \in B\} \subseteq B^*(A)$$

очевидно. Если же $\theta_A(\alpha) \subseteq B^*(A)$, то $\theta_A(\alpha) = \theta_A(\beta)$, где

$$\beta = c_1 \ldots c_{i+k(n-1)} \ (i = 1, \ldots, n-1)$$

составлена из элементов множества B. Так как $<B, [\ ]>$ – n-арная подгруппа n-арной группы $<A, [\ ]>$, то существует $b \in B$, такой, что

$$c_1 \ldots c_{i+k(n-1)} \, \theta_A bb_1 \ldots b_{i-1},$$

откуда

$$\theta_A(\alpha) = \theta_A(\beta) = \theta_A(bb_1 \ldots b_{i-1}).$$

Таким образом, доказано включение

$$B^*(A) \subseteq \bigcup_{i=1}^{n-1} \{\theta_A(bb_1 \ldots b_{i-1}) \mid b \in B\}.$$

а значит и равенство (1).

Равенство (2) следует из (1) и очевидного равенства

$$\{\theta_A(bb_1 \ldots b_{i-1}) \mid b \in B\} = B^{(i)}(A).$$

Легко проверяется, что отображение

$$\varphi : \theta_B(bb_1 \ldots b_{i-1}) \to \theta_A(bb_1 \ldots b_{i-1})$$

является изоморфизмом $B^*$ на $B^*(A)$. ∎

**2.2.20. Замечание.** Ясно, что сужение $\varphi$ на $B_o$ является изоморфизмом $B_o$ на $B_o(A)$. Следовательно, $B_o(A)$ – инвариантная подгруппа группы $B^*(A)$.

## §2.3. ИНВАРИАНТНЫЕ И ПОЛУИНВАРИАНТНЫЕ n-АРНЫЕ ПОДГРУППЫ.

Следующие два определения, обобщающие на n-арный случай понятие инвариантности для подгрупп, были предложены Дёрнте [1].

**2.3.1. Определение.** n-Арная подгруппа $< B, [\ ] >$ n-арной группы $< A, [\ ] >$ называется *инвариантной* в ней, если

$$[x\underbrace{B\ldots B}_{n-1}] = [\underbrace{B\ldots B}_{i-1} x \underbrace{B\ldots B}_{n-i}]$$

для любого $x \in A$ и всех $i = 2, 3, \ldots, n$.

**2.3.2. Определение.** n-Арная подгруппа $< B, [\ ] >$ n-арной группы $< A, [\ ] >$ называется *полуинвариантной* в ней, если

$$[x\underbrace{B\ldots B}_{n-1}] = [\underbrace{B\ldots B}_{n-1}x]$$

для любого x ∈ A.

Ясно, что всякая инвариантная в < A, [ ] > n-арная подгруппа будет и полуинвариантной в < A, [ ] >. Каждая n-арная группа инвариантна и полуинвариантна в самой себе.

**2.3.3. Пример.** Пусть < A, [ ] > – n-арная группа, производная от группы A, B – инвариантная подгруппа группы A. Так как

$$x\underbrace{B\ldots B}_{n-1} = \underbrace{B\ldots B}_{i-1}x\underbrace{B\ldots B}_{n-i}$$

для любого x ∈ A и всех i = 2, 3, …, n, то

$$[x\underbrace{B\ldots B}_{n-1}] = [\underbrace{B\ldots B}_{i-1}x\underbrace{B\ldots B}_{n-i}],$$

т. е. < B, [ ] > – инвариантная n-арная подгруппа в < A, [ ] >.

**2.3.4. Пример.** Пусть < $S_n$, [ ] > – тернарная группа, производная от симметрической группы $S_n$. Так как $A_n$ – инвариантна в $S_n$, то < $A_n$, [ ] > – инвариантна в < $S_n$, [ ] >. Так как

$$[xB_nB_n] = [B_nxB_n] = [B_nB_nx] = A_n,\ x \in A_n,$$

$$[xB_nB_n] = [B_nxB_n] = [B_nB_nx] = B_n,\ x \in B_n,$$

то тернарная группа < $B_n$, [ ] > всех нечётных подстановок также инвариантна в < $S_n$, [ ] >.

**2.3.5. Пример.** Легко проверяется, что тернарные подгруппы < $E_1(2)$, [ ] > и < $E_2(2)$, [ ] > движений первого и второго рода плоскости являются инвариантными в тернарной группе < E(2), [ ] > всех движений плоскости. Аналогично тернарные подгруппы < $E_1(3)$, [ ] > и < $E_2(3)$, [ ] > всех движений первого и второго рода пространства являются инвариантными в тернарной группе < E(3), [ ] > всех движений пространства.

**2.3.6. Пример.** Если ε – идемпотент n-арной группы < A, [ ] >, то

$$[x\underbrace{\varepsilon\ldots\varepsilon}_{n-1}] = [\underbrace{\varepsilon\ldots\varepsilon}_{n-1}x]$$

для любого $x \in A$. Поэтому $< \{\varepsilon\}, [\ ] >$ – полуинвариантная n-арная подгруппа в $< A, [\ ] >$. Если же $\varepsilon$ – единица в $< A, [\ ] >$, то $< \{\varepsilon\}, [\ ] >$ – инвариантная в $< A, [\ ] >$ – n-арная подгруппа.

**2.3.7. Предложение.** Если индекс n-арной подгруппы $< B, [\ ] >$ в n-арной группе $< A, [\ ] >$ равен 2, то $< B, [\ ] >$ – полуинвариантна в $< A, [\ ] >$.

*Доказательство.* Так как $|A : B| = 2$, то

$$A = B \cup [x\underbrace{B\ldots B}_{n-1}] = B \cup [\underbrace{B\ldots B}_{n-1}x]$$

для любого $x \in A$, $x \notin B$. Следовательно,

$$[x\underbrace{B\ldots B}_{n-1}] = [\underbrace{B\ldots B}_{n-1}x]. \qquad \blacksquare$$

**2.3.8. Предложение** [12]**.** Любая тернарная подгруппа тернарной группы $< B_n, [\ ] >$ отражений правильного n-угольника полуинвариантна в ней.

**2.3.9. Теорема.** Если $< B, [\ ] >$ n-арная подгруппа n-арной группы $< A, [\ ] >$, то следующие утверждения эквивалентны:

1) $< B, [\ ] >$ – инвариантна в $< A, [\ ] >$;

2) $[x\underbrace{B\ldots B}_{n-1}] = [Bx\underbrace{B\ldots B}_{n-2}] = [\underbrace{B\ldots B}_{n-1}x]$ для любого $x \in A$;

3) $[\alpha B \beta] = B$ для любых взаимно обратных последовательностей $\alpha$ и $\beta$, составленных из элементов множества $A$;

4) $[xB\underbrace{x\ldots x}_{n-3}\bar{x}] = B$ для любого $x \in A$;

5) $[\bar{x}\underbrace{x\ldots x}_{n-3}Bx] = B$ для любого $x \in A$.

*Доказательство.* 1) $\Rightarrow$ 2). Очевидно.

2) $\Rightarrow$ 3). Положив
$$\beta = b_1\ldots b_k \ (k \geq 1),$$
и, используя 2), а также нейтральность последовательности $\alpha\beta$, получим

$$[\alpha B\beta] = [\alpha[\underbrace{B\ldots B}_{n}]b_1\ldots b_k] =$$

$$= [\alpha B[\underbrace{B\ldots B}_{n-1}b_1]b_2\ldots b_k] = [\alpha B[b_1\underbrace{B\ldots B}_{n-1}]b_2\ldots b_k] =$$

$$= [\alpha[Bb_1\underbrace{B\ldots B}_{n-2}]Bb_2\ldots b_k] = [\alpha[b_1\underbrace{B\ldots B}_{n-1}]B b_2\ldots b_k] =$$

$$= [\alpha b_1[\underbrace{B\ldots B}_{n}]b_2\ldots b_k] = \ldots = [\alpha b_1\ldots b_k[\underbrace{B\ldots B}_{n}]] =$$

$$= [\alpha\beta B] = B,$$

т. е. верно 3).

3) $\Rightarrow$ 4). Получается из 3) при
$$\alpha = x, \ \beta = \underbrace{x\ldots x}_{n-3}\bar{x}.$$

4) $\Rightarrow$ 5). Из 4) с учётом нейтральности последовательностей

$$\bar{x}\underbrace{x\ldots x}_{n-2}, \ \underbrace{x\ldots x}_{n-3}\bar{x}x,$$

получаем

$$[\bar{x}\underbrace{x\ldots x}_{n-3}[xB\underbrace{x\ldots x}_{n-3}\bar{x}]x] = [\bar{x}\underbrace{x\ldots x}_{n-3}Bx],$$

$$[\bar{x}\underbrace{x\ldots x}_{n-2}B\underbrace{x\ldots x}_{n-3}\bar{x}x] = [\bar{x}\underbrace{x\ldots x}_{n-3}Bx]$$

$$B = [\bar{x}\underbrace{x\ldots x}_{n-3}Bx].$$

5) $\Rightarrow$ 1). Используя равенство из 5) и учитывая нейтральность последовательности $x\bar{x}\underbrace{x\ldots x}_{n-3}$, имеем

$$[x\underbrace{B\ldots B}_{n-1}] = [x[\bar{x}\underbrace{x\ldots x}_{n-3}Bx]\underbrace{B\ldots B}_{n-2}] =$$

$$= [x\bar{x}\underbrace{x\ldots x}_{n-3}[Bx\underbrace{B\ldots B}_{n-2}]] = [Bx\underbrace{B\ldots B}_{n-2}] =$$

$$= [Bx[\bar{x}\underbrace{x\ldots x}_{n-3}Bx]\underbrace{B\ldots B}_{n-3}] =$$

$$= [Bx\bar{x}\underbrace{x\ldots x}_{n-3}Bx\underbrace{B\ldots B}_{n-3}] = [BBx\underbrace{B\ldots B}_{n-3}] = \ldots$$

$$\ldots = [\underbrace{B\ldots B}_{n-1}x],$$

т. е.

$$[x\underbrace{B\ldots B}_{n-1}] = [\underbrace{B\ldots B}_{i-1}x\underbrace{B\ldots B}_{n-i}]$$

для любого $x \in A$ и любого $i = 2, \ldots, n$. ∎

**2.3.10. Следствие.** Для тернарной подгруппы $<B, [\ ]>$ тернарной группы $<A, [\ ]>$ следующие утверждения эквивалентны:

1) $<B, [\ ]>$ – инвариантна в $<A, [\ ]>$;

2) $[xB\bar{x}] = B$ для любого $x \in A$;

3) $[\bar{x}Bx] = B$ для любого $x \in A$.

Представляет интерес следующий критерий инвариантности, обобщающий соответствующий бинарный результат [27].

**2.3.11. Теорема** [12]. n-Арная подгруппа $<B, [\ ]>$ инвариантна в n-арной группе $<A, [\ ]>$ тогда и только тогда, когда для любых $a_1, a_2, \ldots, a_n \in A$ из условия

$$[a_1 a_2 \ldots a_n] \in B$$

вытекает

$$[a_2 \ldots a_n a_1] \in B.$$

**2.3.12. Предложение** [12]. n-Арная подгруппа $< B, [\ ] >$ n-арной группы $< A, [\ ] >$ полуинвариантна в ней тогда и только тогда, когда $_aB = B_a$ и подгруппа $< {_aB},@ >$ инвариантна в группе $< A,@ >$.

***Доказательство.*** *Необходимость.* Так как $< B, [\ ] >$ — полуинвариантна в $< A, [\ ] >$, то

$$[x\underbrace{B\ldots B}_{n-1}] = [\underbrace{B\ldots B}_{n-1}x] \qquad (1)$$

для любого $x \in A$. Полагая в последнем равенстве $x = a$, получим

$$[a\underbrace{B\ldots B}_{n-1}] = [\underbrace{B\ldots B}_{n-1}a],$$

т. е. $_aB = B_a$. Тогда

$$[x\underbrace{B\ldots B}_{n-1}] = [x\alpha a\underbrace{B\ldots B}_{n-1}] = [x\alpha[a\underbrace{B\ldots B}_{n-1}]] = [x\alpha {_aB}] = x\,@\,{_aB}, \qquad (2)$$

$$[\underbrace{B\ldots B}_{n-1}x] = [\underbrace{B\ldots B}_{n-1}a\alpha x] = [[\underbrace{B\ldots B}_{n-1}a]\alpha x] =$$

$$= [B_a \alpha x] = [{_aB}\alpha x] = {_aB}\,@\,x. \qquad (3)$$

Левые части в (2) и (3) равны, поэтому равны и правые части, т. е.

$$x\,@\,{_aB} = {_aB}\,@\,x. \qquad (4)$$

для любого $x \in A$. Следовательно, подгруппа $< {_aB},@ >$ инвариантна в группе $< A,@ >$.

*Достаточность.* Пусть теперь $_aB = B_a$ и $< {_aB},@ >$ — инвариантна в $< A,@ >$, т. е. верно (4), откуда, используя (2) и (3), получаем (1). ∎

Если a ∈ B, то $_aB = B_a = B$ и имеет место

**2.3.13. Следствие.** Если a ∈ B, то n-арная подгруппа < B, [ ] > n-арной группы < A, [ ] > – полуинвариантна в ней тогда и только тогда, когда подгруппа < B,@ > инвариантна в группе < A,@ >.

Приведём ещё один полезный критерий полуинвариантности.

**2.3.14. Предложение.** Для того, чтобы n-арная подгруппа < B, [ ] > n-арной группы < A, [ ] > была полуинвариантна в ней, необходимо, чтобы для любых $b, b_1, \ldots, b_{n-1} \in B$, $x \in A$ выполнялось условие

$$[bx^{-1}b_1\ldots b_{n-1}x] \in B, \qquad (1)$$

где $x^{-1}$ – обратная последовательность для x, и достаточно, чтобы это условие выполнялось для любых $b_1, \ldots, b_{n-1} \in B$, $x \in A$ и некоторого $b \in B$.

**Доказательство.** *Необходимость.* Если < B, [ ] > – полуинвариантна в < A, [ ] >, то

$$[x\underbrace{B\ldots B}_{n-1}] = [\underbrace{B\ldots B}_{n-1}x]$$

для любого x ∈ A. Тогда

$$[bx^{-1}[x\underbrace{B\ldots B}_{n-1}]] = [bx^{-1}[\underbrace{B\ldots B}_{n-1}x]],$$

$$[bx^{-1}x\underbrace{B\ldots B}_{n-1}] = [bx^{-1}\underbrace{B\ldots B}_{n-1}x],$$

$$[b\underbrace{B\ldots B}_{n-1}] = [bx^{-1}\underbrace{B\ldots B}_{n-1}x],$$

$$B = [bx^{-1}\underbrace{B\ldots B}_{n-1}x] \qquad (2)$$

для любого b ∈ B. Из последнего равенства вытекает (1)

*Достаточность.* Если для некоторого b ∈ B и любых $b_1$, …, $b_{n-1}$ ∈ B, x ∈ A верно (1), то верно следующее включение

$$[bx^{-1}\underbrace{B\ldots B}_{n-1}x] \subseteq B, \qquad (3)$$

откуда, положив y = $[bbx^{-1}]$ и, учитывая, что $y^{-1} = xb^{-1}b^{-1}$, где $b^{-1}$ – обратная последовательность для b, составленная из элементов множества B, получим

$$[by^{-1}\underbrace{B\ldots B}_{n-1}y] \subseteq B,$$

$$[bxb^{-1}b^{-1}\underbrace{B\ldots B}_{n-1}[bbx^{-1}]] \subseteq B,$$

$$[bxb^{-1}b^{-1}[\underbrace{B\ldots B}_{n-1}b]bx^{-1}] \subseteq B,$$

$$[bxb^{-1}[b^{-1}Bb]x^{-1}] \subseteq B,$$

$$[bxb^{-1}Bx^{-1}] \subseteq B,$$

$$[bx^{-1}b^{-1}[bxb^{-1}Bx^{-1}]x] \subseteq [bx^{-1}b^{-1}Bx],$$

$$[\underbrace{bx^{-1}b^{-1}bxb^{-1}}_{\text{нейтр.}}Bx^{-1}x] \subseteq [bx^{-1}[b^{-1}Bx]],$$

$$B \subseteq [bx^{-1}[\underbrace{B\ldots B}_{n-1}x]],$$

$$B \subseteq [bx^{-1}\underbrace{B\ldots B}_{n-1}x]. \qquad (4)$$

Из (3) и (4) вытекает (2), откуда

$$[xb^{-1}B] = [xb^{-1}[bx^{-1}\underbrace{B\ldots B}_{n-1}x]],$$

$$[x\underbrace{B\ldots B}_{n-1}] = [x\,b^{-1}b\,x^{-1}\underbrace{B\ldots B}_{n-1}x],$$

$$[x\underbrace{B\ldots B}_{n-1}] = [\underbrace{B\ldots B}_{n-1}x].$$

Следовательно, $<B, [\,]>$ – полуинвариантна в $<A, [\,]>$. ∎

**2.3.15. Предложение.** Пусть

$$<B, [\,]> = <\cap B_i, [\,]>, \quad \cap B_i \neq \varnothing$$

непустое пересечение семейства $\{<B_i, [\,]> \mid i \in I\}$ n-арных подгрупп n-арной группы $<A, [\,]>$. Тогда справедливы следующие утверждения:

1) если все $<B_i, [\,]>$ – инвариантны в $<A, [\,]>$, то $<B, [\,]>$ – инвариантна в $<A, [\,]>$;

2) если все $<B_i, [\,]>$ – полуинвариантны в $<A, [\,]>$, то $<B, [\,]>$ – полуинвариантна в $<A, [\,]>$.

*Доказательство.* 1) По предложению 2.3.12, $<B, [\,]>$ – n-арная подгруппа в $<A, [\,]>$. Так как все $<B_i, [\,]>$ – инвариантны в $<A, [\,]>$, то согласно 3) теоремы 2.3.9,

$$[\alpha B_i \beta] = B_i$$

для любых взаимно обратных последовательностей $\alpha$ и $\beta$, составленных из элементов множества A, откуда

$$[\alpha B \beta] \subseteq B_i,$$

$$[\alpha B \beta] \subseteq \cap B_i = B. \qquad (1)$$

Так как $\alpha$ и $\beta$ – произвольные, то из (1) следует

$$[\beta B \alpha] \subseteq B,$$

откуда, учитывая нейтральность последовательности $\alpha\beta$, получаем

$$[\alpha[\beta B \alpha]\beta] \subseteq [\alpha B \beta],$$

$$[\alpha\beta B\alpha\beta] \subseteq [\alpha B\beta],$$

$$B \subseteq [\alpha B\beta]. \qquad (2)$$

Из (1) и (2) вытекает равенство $[\alpha B\beta] = B$. Применяя теперь 3) теоремы 2.3.9, заключаем, что $<B, [\ ]>$ – инвариантна в $<A, [\ ]>$.

2) Зафиксируем элемент $a \in B$. Так как $a \in B_i$, и все $<B_i, [\ ]>$ – полуинвариантны в $<A, [\ ]>$, то по следствию 2.3.13, все $<B_i,@>$ – инвариантны в $<A,@>$. По соответствующей бинарной теореме, подгруппа $<B,@>$ инвариантна в группе $<A,@>$. Ещё раз применяя следствие 2.3.13, заключаем, что $<B, [\ ]>$ полуинвариантна в $<A, [\ ]>$. ∎

**2.3.16. Замечание.** Утверждение 2) предыдущего предложения можно доказать и непосредственно, не используя соответствующий групповой результат.

**2.3.17. Лемма.** Пусть $<B, [\ ]>$ и $<C, [\ ]>$ – n-арные подгруппы n-арной группы $<A, [\ ]>$. Если

$$[\underbrace{B\ldots B}_{n-1}C] = [C\underbrace{B\ldots B}_{n-1}], \qquad (*)$$

то

$$<D, [\ ]> = <[\underbrace{B\ldots B}_{n-1}C], [\ ]>$$

– n-арная подгруппа в $<A, [\ ]>$.

***Доказательство.*** Если $d_1, \ldots, d_n$ – произвольные элементы из D, то используя (*), имеем

$$[d_1\ldots d_n] \in [\underbrace{D\ldots D}_{n}] = [\underbrace{[\underbrace{B\ldots B}_{n-1}C]\ldots[\underbrace{B\ldots B}_{n-1}C]}_{n}] =$$

$$= [\underbrace{\underbrace{B\ldots B}_{n-1}\ldots \underbrace{B\ldots B}_{n-1}}_{n}\underbrace{C\ldots C}_{n}] = [\underbrace{\underbrace{B\ldots B}_{n}\ldots \underbrace{B\ldots B}_{n}}_{n-1}\underbrace{C\ldots C}_{n}] =$$

$$= [\underbrace{B\ldots B}_{n-1}C] = D.$$

Следовательно, $<D, [\ ]>$ – n-арная полугруппа.

Рассмотрим в D уравнение

$$[d_1\ldots d_{n-1}x] = d,$$

где

$$d_i = [b_1^{(i)}\ldots b_{n-1}^{(i)}c_i] \in D, d = [b_1\ldots b_{n-1}c] \in D.$$

Записанное уравнение имеет в $<A, [\ ]>$ решение $x = u$. Поэтому

$$[[b_1^{(1)}\ldots b_{n-1}^{(1)}c_1]\ldots [b_1^{(n-1)}\ldots b_{n-1}^{(n-1)}c_{n-1}]u] = [b_1\ldots b_{n-1}c],$$

откуда

$$u = [c_1^{(n-1)}\ldots c_{n-2}^{(n-1)} \widetilde{b}_1^{(n-1)}\ldots \widetilde{b}_{n-1}^{(n-1)}\ldots$$

$$\ldots c_1^{(1)}\ldots c_{n-2}^{(1)} \widetilde{b}_1^{(1)}\ldots \widetilde{b}_{n-1}^{(1)}b_1\ldots b_{n-1}c],$$

где $\widetilde{b}_1^{(i)}\ldots \widetilde{b}_{n-1}^{(i)}$ – обратная последовательность для $b_1^{(i)}\ldots b_{n-1}^{(i)}$, составленная из элементов множества B, $c_1^{(i)}\ldots c_{n-2}^{(i)}$ – обратная последовательность для элемента $c_i$, составленная из элементов множества C. Из последнего равенства, используя (*), имеем

$$u \in [\underbrace{\underbrace{C\ldots C}_{n-2}\underbrace{B\ldots B}_{n-1}\ldots \underbrace{C\ldots C}_{n-2}\underbrace{B\ldots B}_{n-1}}_{n-1}\underbrace{B\ldots B}_{n-1}C] =$$

$$= [\underbrace{\underbrace{B\ldots B}_{n-1}\ldots \underbrace{B\ldots B}_{n-1}}_{n}\underbrace{C\ \ldots \ C}_{(n-1)(n-2)+1}] =$$

$$= [[\underbrace{\underbrace{B\ldots B}_{n}\ldots\underbrace{B\ldots B}_{n}}_{n-1}][\underbrace{C\ldots\phantom{C}C}_{(n-1)(n-2)+1}]] = [\underbrace{B\ldots B}_{n-1}C] = D,$$

т. е. $u \in D$.

Аналогично доказывается разрешимость в D уравнения

$$[yd_1\ldots d_{n-1}] = d. \qquad\blacksquare$$

**2.3.18. Следствие.** Если $< B, [\ ] >$ и $< C, [\ ] >$ — полуинвариантные n-арные подгруппы n-арной группы $< A, [\ ] >$, то

$$< D, [\ ] > = < [\underbrace{B\ldots B}_{n-1}C], [\ ] >$$

— полуинвариантная n-арная подгруппа в $< A, [\ ] >$.

***Доказательство.*** Так как $< B, [\ ] >$ — полуинвариантна в $< A, [\ ] >$, то верно (*) из леммы 2.3.17, согласно которой $< D, [\ ] >$ — n-арная подгруппа в $< A, [\ ] >$.

Так как $< B, [\ ] >$ и $< C, [\ ] >$ полуинвариантны в $< A, [\ ] >$, то

$$[x\underbrace{D\ldots D}_{n-1}] = [x\underbrace{[\underbrace{B\ldots B}_{n-1}C]\ldots[\underbrace{B\ldots B}_{n-1}C]}_{n-1}] =$$

$$= [x\underbrace{\underbrace{B\ldots B}_{n-1}\ldots\underbrace{B\ldots B}_{n-1}}_{n-1}\underbrace{C\ldots C}_{n-1}] = [\underbrace{\underbrace{B\ldots B}_{n-1}\ldots\underbrace{B\ldots B}_{n-1}}_{n-1}\underbrace{C\ldots C}_{n-1}x] =$$

$$= [\underbrace{[\underbrace{B\ldots B}_{n-1}C]\ldots[\underbrace{B\ldots B}_{n-1}C]}_{n-1}x] = [\underbrace{D\ldots D}_{n-1}x],$$

т. е.

$$[x\underbrace{D\ldots D}_{n-1}] = [\underbrace{D\ldots D}_{n-1}x]. \qquad\blacksquare$$

**2.3.19. Следствие.** Если $< B, [\ ] >$ и $< C, [\ ] >$ — инвариантные n-арные подгруппы n-арной группы $< A, [\ ] >$, то

$$< D, [\ ] > = < [\underbrace{B\ldots B}_{n-1}C], [\ ] >$$

– инвариантная n-арная подгруппа в $< A, [\ ] >$.

*Доказательство.* По предыдущему следствию $< D, [\ ] >$ – полуинвариантная n-арная подгруппа в $< A, [\ ] >$. Поэтому ввиду 2) теоремы 2.3.9, для инвариантности $< D, [\ ] >$ в $< A, [\ ] >$ достаточно доказать равенство

$$[x\underbrace{D\ldots D}_{n-1}] = [Dx\underbrace{D\ldots D}_{n-2}].$$

Последнее равенство верно, так как, используя инвариантность $< B, [\ ] >$ и $< C, [\ ] >$ в $< A, [\ ] >$, имеем

$$[x\underbrace{D\ldots D}_{n-1}] = [x\underbrace{[\underbrace{B\ldots B}_{n-1}C]\ldots[\underbrace{B\ldots B}_{n-1}C]}_{n-1}] =$$

$$= [[\underbrace{B\ldots B}_{n-1}C]x\underbrace{[\underbrace{B\ldots B}_{n-1}C]\ldots[\underbrace{B\ldots B}_{n-1}C]}_{n-2}] = [Dx\underbrace{D\ldots D}_{n-2}]. \qquad \blacksquare$$

**2.3.20. Лемма.** Пусть $< B, [\ ] >$ и $< C, [\ ] >$ – n-арные подгруппы n-арной группы $< A, [\ ] >$, причем $B \cap C \neq \varnothing$. Тогда

$$[\underbrace{B\ldots B}_{n-1}C] = [B\underbrace{C\ldots C}_{n-1}].$$

*Доказательство.* Зафиксировав элемент $a \in B \cap C$ и, используя утверждение 1) предложения 1.3.7, получим

$$[\underbrace{B\ldots B}_{n-1}C] = \{[b_1\ldots b_{n-1}c] \mid b_i \in B, c \in C\} =$$

$$= \{[b\underbrace{a\ldots a}_{n-2}c] \mid b \in B, c \in C\} =$$

$$= \{[bc_1\ldots c_{n-1}] \mid b \in B, c_i \in C\} = [B\underbrace{C\ldots C}_{n-1}]. \qquad \blacksquare$$

Если $< B, [\ ] >$ и $< C, [\ ] >$ – n-арные подгруппы n-арной группы $< A, [\ ] >$, то положим

$$< B, [\ ] > \wedge < C, [\ ] > = < B \cap C, [\ ] >,$$

$$< B, [\ ] > \vee < C, [\ ] > = < D, [\ ] >,$$

где $< D, [\ ] >$ – пересечение всех n-арных подгрупп, содержащих $< B, [\ ] >$ и $< C, [\ ] >$, т. е. $< D, [\ ] >$ – n-арная подгруппа, порожденная множеством $B \cup C$.

Ясно, что множество $L(A, [\ ])$ всех n-арных подгрупп n-арной группы $< A, [\ ] >$, дополненное пустым множеством, образует полную решетку относительно операций $\wedge$ и $\vee$.

**2.3.21. Теорема.** Множество всех полуинвариантных (инвариантных) n-арных подгрупп n-арной группы $< A, [\ ] >$, содержащих фиксированный элемент, образуют подрешетку решетки $L(A, [\ ])$.

***Доказательство.*** Если $< B, [\ ] >$ и $< C, [\ ] >$ – полуинвариантные (инвариантные) n-арные подгруппы из $< A, [\ ] >$, содержащие фиксированный элемент $a \in A$, то по предложению 2.3.15

$$< B, [\ ] > \wedge < C, [\ ] > = < B \cap C, [\ ] >$$

– полуинвариантная (инвариантная) n-арная подгруппа в $< A, [\ ] >$, причем $a \in B \cap C$.

Положим

$$< D, [\ ] > = < B, [\ ] > \vee < C, [\ ] >.$$

По следствию 2.3.18 (следствию 2.3.19)

$$< [\underbrace{B \ldots B}_{n-1} C], [\ ] >$$

– полуинвариантная (инвариантная) n-арная подгруппа в $< A, [\ ] >$. Так как в $< B, [\ ] >$ и в $< C, [\ ] >$ имеются нейтральные последовательности, то

$$C \subseteq [\underbrace{B \ldots B}_{n-1} C],\ B \subseteq [B \underbrace{C \ldots C}_{n-1}],$$

а так как, кроме того, $B \cap C \neq \varnothing$, то из второго включения, учитывая лемму 2.3.20, получаем

$$B \subseteq [\underbrace{B\ldots B}_{n-1}C],$$

откуда

$$B \cup C \subseteq [\underbrace{B\ldots B}_{n-1}C].$$

Следовательно,

$$D \subseteq [\underbrace{B\ldots B}_{n-1}C].$$

Включение

$$[\underbrace{B\ldots B}_{n-1}C] \subseteq D$$

является следствием теоремы 2.1.14 Таким образом,

$$D = [\underbrace{B\ldots B}_{n-1}C].$$

Ясно, что $a \in D$. ∎

Следующее следствие получается из теоремы 2.3.21 с применением леммы 2.3.20.

**2.3.22. Следствие.** Если $< B, [\ ] >$ и $< C, [\ ] >$ — полуинвариантные n-арные подгруппы n-арной группы $< A, [\ ] >$, $B \cap C \neq \varnothing$,

$$< B, [\ ] > \vee < C, [\ ] > = < D, [\ ] >,$$

то

$$D = [\underbrace{B\ldots B}_{n-1}C] = [C\underbrace{B\ldots B}_{n-1}] = [\underbrace{C\ldots C}_{n-1}B] = [B\underbrace{C\ldots C}_{n-1}].$$

**2.3.23. Лемма.** Если $< A, [\ ] >$ — n-арная группа, $B \subseteq A$, $C \subseteq A$, $B \cap C \neq \varnothing$, то

$$[x_1\ldots x_{i-1}(B \cap C)x_{i+1}\ldots x_n] =$$
$$= [x_1\ldots x_{i-1}Bx_{i+1}\ldots x_n] \cap [x_1\ldots x_{i-1}Cx_{i+1}\ldots x_n] \qquad (1)$$

для любых $x_1, \ldots, x_{i-1}, x_{i+1}, \ldots, x_n \in A$ и любого $i = 1, 2, \ldots, n$.

***Доказательство.*** Пусть
$$y \in [x_1 \ldots x_{i-1}(B \cap C)x_{i+1} \ldots x_n],$$
т. е.
$$y = [x_1 \ldots x_{i-1}dx_{i+1} \ldots x_n],$$
где $d \in (B \cap C)$. Это значит, что $d \in B$, $d \in C$, откуда
$$[x_1 \ldots x_{i-1}dx_{i+1} \ldots x_n] \in [x_1 \ldots x_{i-1}Bx_{i+1} \ldots x_n],$$
$$[x_1 \ldots x_{i-1}dx_{i+1} \ldots x_n] \in [x_1 \ldots x_{i-1}Cx_{i+1} \ldots x_n].$$
Следовательно,
$$y = [x_1 \ldots x_{i-1}dx_{i+1} \ldots x_n] \subseteq$$
$$\subseteq [x_1 \ldots x_{i-1}Bx_{i+1} \ldots x_n] \cap [x_1 \ldots x_{i-1}Cx_{i+1} \ldots x_n].$$
Так как элемент $y$ выбран произвольно, то доказано включение
$$[x_1 \ldots x_{i-1}(B \cap C)x_{i+1} \ldots x_n] \subseteq$$
$$\subseteq [x_1 \ldots x_{i-1}Bx_{i+1} \ldots x_n] \cap [x_1 \ldots x_{i-1}Cx_{i+1} \ldots x_n]. \qquad (2)$$

Пусть теперь
$$y \in [x_1 \ldots x_{i-1}Bx_{i+1} \ldots x_n] \cap [x_1 \ldots x_{i-1}Cx_{i+1} \ldots x_n],$$
т. е.
$$y \in [x_1 \ldots x_{i-1}Bx_{i+1} \ldots x_n], \ y \in [x_1 \ldots x_{i-1}Cx_{i+1} \ldots x_n],$$
откуда
$$y = [x_1 \ldots x_{i-1}bx_{i+1} \ldots x_n], \ \ b \in B,$$
$$y = [x_1 \ldots x_{i-1}cx_{i+1} \ldots x_n], \ \ c \in C.$$

Из равенства
$$[x_1 \ldots x_{i-1}bx_{i+1} \ldots x_n] = [x_1 \ldots x_{i-1}cx_{i+1} \ldots x_n]$$
и однозначной разрешимости уравнений в n-арной группе получаем $b = c$, откуда $b = c \in B \cap C$. Следовательно,

$$y \in [x_1 \ldots x_{i-1}(B \cap C)x_{i+1} \ldots x_n].$$

В силу произвольного выбора y, получаем

$$[x_1 \ldots x_{i-1}Bx_{i+1} \ldots x_n] \cap [x_1 \ldots x_{i-1}Cx_{i+1} \ldots x_n] \subseteq$$
$$\subseteq [x_1 \ldots x_{i-1}(B \cap C)x_{i+1} \ldots x_n]. \qquad (3)$$

Из включений (2) и (3) вытекает равенство (1). ∎

**2.3.24. Следствие.** Если $<B, [\ ]>$ и $<C, [\ ]>$ — n-арные подгруппы n-арной группы $<A, [\ ]>$, причём $B \cap C \neq \varnothing$, то

$$[y\underbrace{B \ldots B}_{n-1}] \cap [y\underbrace{C \ldots C}_{n-1}] = [y\underbrace{(B \cap C) \ldots (B \cap C)}_{n-1}],$$

$$[\underbrace{B \ldots B}_{n-1}y] \cap [\underbrace{C \ldots C}_{n-1}y] = [\underbrace{(B \cap C) \ldots (B \cap C)}_{n-1}y]$$

для любого $y \in A$.

*Доказательство.* Если $d \in B \cap C$, то

$$[y\underbrace{B \ldots B}_{n-1}] = [y\underbrace{d \ldots d}_{n-2}B],$$

$$[y\underbrace{C \ldots C}_{n-1}] = [y\underbrace{d \ldots d}_{n-2}C],$$

$$[y\underbrace{(B \cap C) \ldots (B \cap C)}_{n-1}] = [y\underbrace{d \ldots d}_{n-2}(B \cap C)],$$

откуда, используя лемму 2.3.23, получаем

$$[y\underbrace{B \ldots B}_{n-1}] \cap [y\underbrace{C \ldots C}_{n-1}] = [y\underbrace{d \ldots d}_{n-2}B] \cap [y\underbrace{d \ldots d}_{n-2}C] =$$

$$= [y\underbrace{d \ldots d}_{n-2}(B \cap C)] = [y\underbrace{(B \cap C) \ldots (B \cap C)}_{n-1}].$$

Второе равенство доказывается аналогично. ∎

**2.3.25. Предложение.** Если $<B, [\ ]>$ и $<C, [\ ]>$ — n-арные подгруппы n-арной группы $<A, [\ ]>$, причём

< B, [ ] > – полуинвариантна в < A, [ ] > и B ∩ C ≠ ∅, то n-арная подгруппа < B ∩ C, [ ] > – полуинвариантна в < C, [ ] >.

*Доказательство.* Если y – произвольный элемент из C, то, используя предыдущее следствие и полуинвариантность < B, [ ] > в < A, [ ] >, получим

$$[y\underbrace{(B\cap C)\ldots(B\cap C)}_{n-1}] = [y\underbrace{B\ldots B}_{n-1}] \cap [y\underbrace{C\ldots C}_{n-1}] =$$

$$= [\underbrace{B\ldots B}_{n-1}y] \cap [\underbrace{C\ldots C}_{n-1}y] = [\underbrace{(B\cap C)\ldots(B\cap C)}_{n-1}y]. \quad\blacksquare$$

**2.3.26. Предложение.** Если < B, [ ] > и < C, [ ] > – n-арные подгруппы n-арной группы < A, [ ] >, причем < B, [ ] > – инвариантна в < A, [ ] > и B ∩ C ≠ ∅, то n-арная подгруппа < B ∩ C, [ ] > – инвариантна в < C, [ ] >.

*Доказательство.* Если x – произвольный элемент из C, то, используя лемму 2.3.23, утверждение 4) следствия 2.3.9, а также то, что x, $\bar{x}$ ∈ C, получим

$$[x(B \cap C)\underbrace{x\ldots x}_{n-3}\bar{x}] = [xB\underbrace{x\ldots x}_{n-3}\bar{x}] \cap [xC\underbrace{x\ldots x}_{n-3}\bar{x}] = B \cap C,$$

что означает инвариантность < B ∩ C, [ ] > в < C, [ ] >. $\quad\blacksquare$

**2.3.27. Предложение.** Пусть < B, [ ] >, < C, [ ] > и < D, [ ] > – n-арные подгруппы n-арной группы < A, [ ] >, причем < B, [ ] > – полуинвариантна (инвариантна) в < C, [ ] > и B ∩ D ≠ ∅. Тогда n-арная подгруппа < B ∩ D, [ ] > – полуинвариантна (инвариантна) в < C ∩ D, [ ] >.

*Доказательство.* Если < B, [ ] > – полуинвариантна в < C, [ ] > и x ∈ C ∩ D, то используя следствие 2.3.24, полуинвариантность < B, [ ] > в < C, [ ] >, а также то, что x ∈ D, получим

$$[x\underbrace{(B\cap D)\ldots(B\cap D)}_{n-1}] = [x\underbrace{B\ldots B}_{n-1}] \cap [x\underbrace{D\ldots D}_{n-1}] =$$

$$= [\underbrace{B\ldots B}_{n-1}x] \cap [\underbrace{D\ldots D}_{n-1}x] = [\underbrace{(B\cap D)\ldots(B\cap D)}_{n-1}x],$$

что означает полуинвариантность $< B \cap D, [\ ] >$ в $< C \cap D, [\ ] >$.

Если $< B, [\ ] >$ – инвариантна в $< C, [\ ] >$ и $x \in C \cap D$, то используя лемму 2.3.23, утверждение 4) следствия 2.3.9, а также то, что $x, \bar{x} \in D$, получим

$$[x(B \cap D)\underbrace{x\ldots x}_{n-3}\bar{x}] = [xB\underbrace{x\ldots x}_{n-3}\bar{x}] \cap [xD\underbrace{x\ldots x}_{n-3}\bar{x}] = B \cap D,$$

что означает инвариантность $< B \cap D, [\ ] >$ в $< C \cap D, [\ ] >$. ■

Инвариантные и полуинвариантные n-арные подгруппы являются частными случаями введенных в [28] m-полуинвариантных n-арных подгрупп при $m = 2$ и $m = n$ соответственно.

**2.3.28. Определение** [28]. Полуинвариантную n-арную подгруппу $< B, [\ ] >$ n-арной группы $< A, [\ ] >$ назовем *m-полуинвариантной* ($n = k(m - 1) + 1$, $k \geq 1$), если

$$[x \overset{n-1}{B}] = [\overset{m-1}{B} x \overset{n-m}{B}], x \in A. \qquad (*)$$

Из определения вытекает, что n-полуинвариантные n-арные подгруппы n-арной группы – в точности её полуинвариантные n-арные подгруппы.

Можно показать, что n-арная подгруппа $< B, [\ ] >$ n-арной группы $< A, [\ ] >$ является инвариантной тогда и только тогда, когда она 2-полуинвариантна.

**2.3.29. Теорема** [28]. Конечная n-арная подгруппа $< B, [\ ] >$ n-арной группы $< A, [\ ] >$ является m-полу-инвариантной в $< A, [\ ] >$ тогда и только тогда, когда B удовлетворяет условию (*).

*Доказательство.* Необходимость очевидна.

*Достаточность.* Положим в лемме 6.2 [12] $s = 1$, $j = k-1$. Тогда
$$[x \overset{n-1}{B}] = [\overset{n-m}{B} x \overset{m-1}{B}].$$

Так как
$$[\overset{n-m}{B} x \overset{m-1}{B}] = [\overset{n-m}{B} x[\overset{n}{B}] \overset{m-2}{B}] = [\overset{n-m}{B} [x \overset{n-1}{B}] \overset{m-1}{B}] =$$
$$= [\overset{n-m}{B} [\overset{n-m}{B} x \overset{m-1}{B}] \overset{m-1}{B}] = [\overset{n-1}{B} x \overset{n-1}{B}],$$

то
$$[x \overset{n-1}{B}] = [\overset{n-1}{B} x \overset{n-1}{B}].$$

Если $b_1 \ldots b_{n-1}$ – нейтральная последовательность n-арной подгруппы $< B, [\ ] >$, то
$$[\overset{n-1}{B} x] = [\overset{n-1}{B} x b_1 \ldots b_{n-1}] \subseteq [\overset{n-1}{B} x \overset{n-1}{B}] = [x \overset{n-1}{B}],$$

откуда
$$[\overset{n-1}{B} x] \subseteq [x \overset{n-1}{B}].$$

Так как
$$|[\overset{n-1}{B} x]| = |[x \overset{n-1}{B}]| = |B| < \infty,$$

то
$$[x \overset{n-1}{B}] = [\overset{n-1}{B} x]. \blacksquare$$

**2.3.30. Следствие.** Конечная n-арная подгруппа $< B, [\ ] >$ n-арной группы $< A, [\ ] >$ является инвариантной в $< A, [\ ] >$ тогда и только тогда, когда
$$[x \overset{n-1}{B}] = [Bx \overset{n-2}{B}], x \in A.$$

**2.3.31. Пример.** На множестве $B_3 = \{(12), (13), (23)\}$ всех нечётных подстановок множества $\{1, 2, 3\}$ определим 5-арную операцию

$$[x_1 x_2 x_3 x_4 x_5] = x_1 x_2 x_3 x_4 x_5$$

и положим $B = \{(12)\}$. Тогда $< B, [\ ] > -$ 5-арная подгруппа 5-арной группы $< B_3, [\ ] >$. Так как

$$[x \overset{4}{B}] = \{[x(12)(12)(12)(12)]\} = \{x(12)(12)(12)(12)\} = \{x\} =$$
$$= \{(12)(12)(12)(12)x\} = \{[(12)(12)(12)(12)x]\} = [\overset{4}{B} x],$$

то $< B, [\ ] > -$ полуинвариантная 5-арная подгруппа в $< B_3, [\ ] >$.

Аналогично доказывается, что

$$[x \overset{4}{B}] = [\overset{2}{B} x \overset{2}{B}] = \{x\}.$$

Это означает, что $< B, [\ ] > -$ 3-полуинвариантная 5-арная подгруппа в $< B_3, [\ ] >$. Покажем теперь, что $< B, [\ ] >$ не является инвариантной 5-арной подгруппой в $< B_3, [\ ] >$. Действительно,

$$[B(13)\overset{3}{B}] = \{[(12)(13)(12)(12)(12)]\} = \{(12)(13)(12)(12)(12)\} =$$
$$= \{(12)(13)(12)\} = \{(23)\} \neq \{(13)\} = [(13)\overset{4}{B}].$$

Следовательно, существуют n-арные группы, обладающие m-полуинвариантными n-арными подгруппами (m > 2), которые не являются инвариантными.

**2.3.32. Теорема** [28]. Если n-арная подгруппа $< B, [\ ] >$ n-арной группы $< A, [\ ] >$ является m-полуинвариантной и k-полуинвариантной, то она является и r-полуинвариантной, где

$$r - 1 = (m - 1, k - 1).$$

**2.3.33. Следствие.** m-Полуинвариантная n-арная подгруппа $< B, [\ ] >$ n-арной группы $< A, [\ ] >$ является r-полуинвариантной, где $r - 1 = (m - 1, n - 1)$.

*Доказательство.* По определению m-полуинвариантная n-арная подгруппа является полуинвариантной, а значит, и n-полуинвариантной. Теперь применяем теорему 2.3.32. ∎

# §2.4. СОПРЯЖЁННЫЕ И ПОЛУСОПРЯЖЁННЫЕ n-АРНЫЕ ПОДГРУППЫ.

Одно из важнейших понятий теории групп – сопряженность подгрупп H и K в группе G можно определить эквивалентными равенствами

$$H = xKx^{-1}, \quad xK = Hx,$$

которые при переходе к n-арному случаю приводят к разным понятиям: сопряженности и полусопряженности.

**2.4.1. Определение** [4]. Подмножество C n-арной группы $< A, [\ ] >$ называется *сопряженным* в ней посредством последовательности $x_1^i$, где $x_i \in A$, $1 \leq i \leq n - 1$, с подмножеством B, если

$$B = [x_1^i C y_1^j], \qquad (*)$$

где $y_1^j$ – обратная последовательность для последовательности $x_1^i$. В этом случае говорят, что B и C сопряжены в $< A, [\ ] >$.

Справедливость следующего предложения устанавливается проверкой.

**2.4.2. Предложение.** Если $< C, [\ ] >$ – n-арная подгруппа n-арной группы $< A, [\ ] >$, то $< [x_1^i C y_1^j], [\ ] >$ – n-арная подгруппа в $< A, [\ ] >$, изоморфная $< C, [\ ] >$.

Легко также проверяется, что отношение сопряженности на множестве всех n-арных подгрупп n-арной группы является эквивалентностью.

**2.4.3. Лемма.** Если для n-арных подгрупп $< B, [\ ] >$ и $< C, [\ ] >$ n-арной группы $< A, [\ ] >$ верно (*), то

$$B = [xC\tilde{x}_1^t], \quad B = [\tilde{z}_i^t C z] \qquad (**)$$

для некоторых x, z ∈ A, где $\tilde{x}_1^t$ и $\tilde{z}_1^t$ – обратные последовательности соответственно для элементов x и z.

*Доказательство.* Если i = 1, то в первом равенстве доказывать нечего. Поэтому считаем $2 \leq i \leq n - 1$. Пусть $h_1, \ldots, h_{i-1}$ – произвольные элементы из B, а $h_i \ldots h_{n-1}$ – обратная последовательность для последовательности $h_1 \ldots h_{i-1}$. Так как $< B, [\ ] >$ – n-арная подгруппа в $< A, [\ ] >$, то $h_i \ldots h_{n-1}$ можно выбрать так, что $h_i, \ldots, h_{n-1} \in B$. Из (∗) имеем

$$B = [h_i^{n-1} B h_1^{i-1}] = [h_i^{n-1} [x_1^i C y_1^j] h_1^{i-1}] = [[h_i^{n-1} x_1^i] C y_1^j h_1^{i-1}],$$

т. е.
$$B = [xC\tilde{x}_1^t],$$

где
$$x = [h_i^{n-1} x_1^i], \quad \tilde{x}_1^t = y_1^j h_1^{i-1}.$$

Ясно, что $\tilde{x}_1^t = y_1^j h_1^{i-1}$ – обратная последовательность для последовательности $h_i^{n-1} x_1^i$, а значит $\tilde{x}_1^t$ – обратная последовательность и для элемента $x = [h_i^{n-1} x_1^i]$.

Для элемента z доказательство проводится аналогично. ■

**2.4.4. Лемма.** Если для подмножеств B и C n-арной группы $< A, [\ ] >$ верно (∗∗), то

$$[\overset{n-1}{x\ C}] = [\overset{i-1}{B}\ x\ \overset{n-i}{C}], \quad i = 2, \ldots, n, \qquad (\ast\ast\ast)$$

т. е.
$$[\overset{n-1}{x\ C}] = [B\overset{n-2}{x\ C}] = [\overset{2}{B}\ x\ \overset{n-3}{C}] = \ldots = [\overset{n-2}{B}\ xC] = [\overset{n-1}{B}\ x].$$

*Доказательство.* Учитывая нейтральность последовательностей $x\tilde{x}_1^t$ и $\tilde{x}_1^t x$, из (∗∗) получаем

$$[\tilde{x}_1^t Bx] = [\tilde{x}_1^t [xC\tilde{x}_1^t]x] = [\tilde{x}_1^t xC\tilde{x}_1^t x] = C,$$

т. е.
$$[\tilde{x}_1^t Bx] = C.$$

Используя снова нейтральность последовательности $x\tilde{x}_1^t$, а также последнее равенство, имеем

$$[x\overset{n-1}{C}] = [x\underbrace{CC\ldots C}_{i-1}\overset{n-i}{C}] = [x\underbrace{[\widetilde{x}_1^t Bx][\widetilde{x}_1^t Bx]\ldots[\widetilde{x}_1^t Bx]}_{i-1}\overset{n-i}{C}] =$$

$$= [\underbrace{[x\widetilde{x}_1^t B][x\widetilde{x}_1^t B]\ldots[x\widetilde{x}_1^t B]}_{i-1} x \overset{n-i}{C}] = [\underbrace{BB\ldots B}_{i-1} x \overset{n-i}{C}] = [\overset{i-1}{B} x \overset{n-i}{C}]. \quad \blacksquare$$

**2.4.5. Лемма**. Если для n-арных подгрупп $< B, [\ ] >$ и $< C, [\ ] >$ n-арной группы $< A, [\ ] >$ верно (***) для $i = 2$ и $i = n$, т. е.

$$[x \overset{n-1}{C}] = [Bx \overset{n-2}{C}] = [\overset{n-1}{B} x],$$

то для них верно (**).

**Доказательство.** Так как $< B, [\ ] >$ и $< C, [\ ] >$ — n-арные подгруппы n-арной группы $< A, [\ ] >$, то

$$[\overset{n}{B}] = B, \quad [\overset{n}{C}] = C,$$

откуда с учетом (***) сначала для $i = n$, а затем $i = 2$ будем иметь

$$[\widetilde{x}_1^t Bx] = [\widetilde{x}_1^t [\overset{n}{B}]x] = [\widetilde{x}_1^t B[\overset{n-1}{B} x]] = [\widetilde{x}_1^t B[x \overset{n-1}{C}]] = [\widetilde{x}_1^t [Bx \overset{n-2}{C}]C] =$$

$$= [\widetilde{x}_1^t [x \overset{n-1}{C}]C] = [\widetilde{x}_1^t x[\overset{n}{C}]] = [\widetilde{x}_1^t xC] = C,$$

т. е.

$$[\widetilde{x}_1^t Bx] = C,$$

где $\widetilde{x}_1^t$ — обратная последовательность для элемента x. Из последнего равенства получаем $B = [xC\widetilde{x}_1^t]$. $\blacksquare$

**2.4.6. Теорема** [29]. Для того чтобы n-арные подгруппы $< B, [\ ] >$ и $< C, [\ ] >$ n-арной группы $< A, [\ ] >$ были сопряжены в ней, необходимо, чтобы для некоторого $x \in A$ и любого $i = 2, \ldots, n$ выполнялось (***), и достаточно, чтобы (***) выполнялось для некоторого $x \in A$ при $i = 2$ и $i = n$.

*Доказательство.* Для доказательства необходимости последовательно применяются леммы 2.4.3 и 2.4.4. Для доказательства достаточности применяется лемма 2.4.5. ∎

Доказанная теорема позволяет дать новое определение сопряженности n-арных подгрупп.

**2.4.7. Определение** [30]. n-Арные подгруппы $<B, [\ ]>$ и $<C, [\ ]>$ n-арной группы $<A, [\ ]>$ называются *сопряженными* в ней, если
$$[x \overset{n-1}{C}] = [Bx \overset{n-2}{C}] = [\overset{n-1}{B} x]$$
для некоторого $x \in A$.

В этом случае говорят, что $<C, [\ ]>$ *сопряжена с* $<B, [\ ]>$ *в* $<A, [\ ]>$ *посредством элемента* $x$.

Если положить в последнем определении n = 2, то
$$xC = Bx = Bx \text{ или } C = x^{-1}Bx.$$

Следовательно, данное определение согласуется с определением сопряженности подгрупп в группе.

Следующее определение принадлежит Воробьёву [30].

**2.4.8. Определение** [30]. n-Арная подгруппа $<C, [\ ]>$ n-арной группы $<A, [\ ]>$ называется *полусопряженной в ней посредством элемента* $x \in A$ *с n-арной подгруппой* $<B, [\ ]>$, если
$$[x \overset{n-1}{C}] = [\overset{n-1}{B} x].$$

В этом случае говорят, что $<B, [\ ]>$ и $<C, [\ ]>$ *полусопряжены* в $<A, [\ ]>$.

Ясно, что при n = 2 понятия сопряженности и полусопряженности подгрупп совпадают.

Из утверждения 3) теоремы 2.1.12 сразу же вытекает

**2.4.9. Следствие.** Полусопряженные n-арные подгруппы имеют одинаковую мощность.

Из теоремы 2.4.6 получаем

**2.4.10. Следствие.** Сопряженные n-арные подгруппы являются полусопряженными.

Следующий пример показывает, что понятия сопряженности и полусопряженности не тождественны, т. е. существуют полусопряженные n-арные подгруппы, не являющиеся сопряженными.

**2.4.11. Пример.** Тернарные подгруппы $<A_n, [\,]>$ и $<T_n, [\,]>$ четных и нечетных подстановок тернарной группы $<S_n, [\,]>$ всех подстановок степени n не являются сопряженными, однако являются полусопряженными.

**2.4.12. Теорема** [30]. Отношение полусопряженности на множестве всех n-арных подгрупп n-арной группы является отношением эквивалентности.

**2/4.13. Предложение.** n-Арные подгруппы $<B, [\,]>$ и $<C, [\,]>$ n-арной группы $<A, [\,]>$ полусопряжены в ней тогда и только тогда, когда группы $<B_a, @>$ и $<{}_aC, @>$ сопряжены в группе $<A, @>$.

*Доказательство.* Пусть n-арные подгруппы $<B, [\,]>$ и $<C, [\,]>$ полусопряжены в n-арной группе $<A, [\,]>$, т. е. существует $x \in A$ такой, что

$$[x \overset{n-1}{C}] = [\overset{n-1}{B} x].$$

Тогда последовательно будем иметь

$$[x\alpha a \overset{n-1}{C}] = [\overset{n-1}{B} a\alpha x], \quad [x\alpha {}_aC] = [B_a \alpha x],$$

$$x @ {}_aC = B_a @ x,$$

где $\alpha$ – обратная последовательность для a. Следовательно, $<B_a, @>$ и $<{}_aC, @>$ сопряжены в $<A, @>$.

Обратно, если $<B_a, @>$ и $<{}_aC, @>$ сопряжены в $<A, @>$, то верно последнее равенство для некоторого $x \in A$. Рассуж-

дая в обратном порядке, придем к первому равенству, что означает полусопряженность $<B, [\,]>$ и $<C, [\,]>$ в $<A, [\,]>$. ∎

**2.4.14. Следствие.** n-Арные подгруппы $<B, [\,]>$ и $<C, [\,]>$ n-арной группы $<A, [\,]>$, имеющие общий элемент a, полусопряжены в ней тогда и только тогда, когда группы $<B, @>$ и $<C, @>$ сопряжены в группе $<A, @>$.

## § 2.5 ЦИКЛИЧЕСКИЕ И ПОЛУЦИКЛИЧЕСКИЕ n-АРНЫЕ ГРУППЫ

**2.5.1. Определение** [3, 4]. Для любого элемента a n-арной ($n \geq 3$) группы $<A, [\,]>$ и любого целого s определим *s-ую n-адическую степень* следующим образом

$$a^{[s]} = \begin{cases} a, & s = 0, \\ [\overset{s(n-1)+1}{a}], & s > 0, \\ [\overset{-2s}{\overline{a}} \overset{-s(n-3)+1}{a}], & s < 0. \end{cases}$$

Так как

$$[[\overset{-2s}{\overline{a}} \overset{-s(n-3)+1}{a}] \overset{-s(n-1)}{a}] = [\overset{-2s}{\overline{a}} \overset{-2s(n-2)}{a}] = a,$$

то s-ую n-адическую степень элемента a при $s < 0$ можно определить [4] как решение уравнения

$$[x \overset{-s(n-1)}{a}] = a.$$

Полагая в определении 2.5.1 $s = -1$, получим

$$a^{[-1]} = [\overline{a}\,\overline{a}\,\overset{n-2}{a}] = \overline{a},$$

т. е. $a^{[-1]} = \overline{a}$.

**2.5.2. Лемма.** Для любого элемента a n-арной группы $<A, [\,]>$ и любого целого s справедливо равенство

$$\theta(a^{[s]}) = \theta^{s(n-1)+1}(a).$$

***Доказательство.*** Если s = 0, то
$$\theta(a^{[s]}) = \theta(a^{[0]}) = \theta(a) = \theta^{0(n-1)+1}(a) = \theta^{s(n-1)+1}(a).$$

Если же s > 0, то
$$\theta(a^{[s]}) = \theta([\overset{s(n-1)+1}{a}]) = \theta^{s(n-1)+1}(a).$$

Пусть теперь s < 0. Тогда из
$$[a^{[s]}\overset{-s(n-1)}{a}] = a$$
следует
$$\theta([a^{[s]}\overset{-s(n-1)}{a}]) = \theta(a),$$
откуда
$$\theta(a^{[s]})\theta^{-s(n-1)}(a) = \theta(a),$$
$$\theta(a^{[s]}) = \theta(a)\theta^{s(n-1)}(a) = \theta^{s(n-1)+1}(a). \qquad \blacksquare$$

**2.5.3. Предложение** [3, 4]. Пусть $< A, [\ ] >$ – n-арная группа, $a \in A$, $k_1, k_2, \ldots, k_n$ – целые. Тогда:

1) $[a^{[k_1]} a^{[k_2]} \ldots a^{[k_n]}] = a^{[k_1+k_2\ldots k_n+1]}$;

2) $(a^{[k_1]})^{[k_2]} = a^{[k_1 k_2 (n-1) + k_1 + k_2]}$.

***Доказательство.*** 1) Используя лемму 2.5.2, получаем
$$\theta(a^{[k_1]} a^{[k_2]} \ldots a^{[k_n]}]) = \theta(a^{[k_1]})\theta(a^{[k_2]}) \ldots \theta(a^{[k_n]}) =$$
$$= \theta^{k_1(n-1)+1}(a)\theta^{k_2(n-1)+1}(a) \ldots \theta^{k_n(n-1)+1}(a) =$$
$$= \theta^{(k_1+k_2+\ldots+k_n)(n-1)+n}(a) = \theta^{(k_1+k_2+\ldots+k_n+1)(n-1)+1}(a) =$$
$$= \theta(a^{[k_1+k_2+\ldots+k_n+1]}),$$
т. е.
$$\theta([a^{[k_1]} a^{[k_2]} \ldots a^{[k_n]}]) = \theta(a^{[k_1+k_2+\ldots+k_n+1]}),$$
откуда и из предложения 1.3.5 следует 1).

3) Снова используя лемму 2.5.2, получим

$$\theta((a^{[k_1]})^{[k_2]}) = \theta^{k_2(n-1)+1}(a^{[k_1]}) = \underbrace{\theta(a^{[k_1]})\ldots\theta(a^{[k_1]})}_{k_2(n-1)+1} =$$

$$= \underbrace{\theta^{k_1(n-1)+1}(a)\ldots\theta^{k_1(n-1)+1}(a)}_{k_2(n-1)+1} = \theta^{(k_1(n-1)+1)(k_2(n-1)+1)}(a) =$$

$$= \theta^{(k_1k_2(n-1)+k_1+k_2)(n-1)+1}(a) = \theta(a^{[k_1k_2(n-1)+k_1+k_2]}),$$

т. е.

$$\theta((a^{[k_1]})^{[k_2]}) = \theta(a^{[k_1k_2(n-1)+k_1+k_2]}),$$

откуда и из предложения 1.3.5 следует 2). ∎

**2.5.4. Следствие.** Для любых целых $k_1$, $k_2$, …, $k_{m(n-1)+1}$, $m > 0$ и любого элемента a n-арной группы $<A, [\ ]>$ справедливо равенство

$$[a^{[k_1]}a^{[k_2]}\ldots a^{[k_{m(n-1)+1}]}] = a^{[k_1+k_2+\ldots+k_{m(n-1)+1}+m]}.$$

Теорема 2.1.10 обобщает на n-арный случай бинарный результат о совпадении подгруппы, порожденной множеством M, со множеством всех конечных произведений элементов из M и обратных к ним. С.А. Русаков, используя понятие n-арной степени элемента, получил [4] другое описание n-арной подгруппы, порожденной некоторым множеством.

**2.5.5. Теорема** [4]**.** Если M – непустое подмножество n-арной группы $<A, [\ ]>$, то

$$\langle M \rangle = \{(a_1^{[k_1]}a_2^{[k_2]}\ldots a_{m(n-1)+1}^{[k_{m(n-1)+1}]}) \mid a_j \in M, k_j \gtrless 0, m \in \mathbb{N}\}.$$

**2.5.6. Предложение.** Пусть $<A, [\ ]>$ – n-арная группа, $a \in A$, $s_1$ и $s_2$ – целые. Тогда следующие утверждения равносильны:

1) $a^{[s_1]} = a^{[s_2]}$;

2) $a^{[s_1-s_2]} = a$;

3) $a^{[s_2 - s_1]} = a$.

**Доказательство.** 1) $\Rightarrow$ 2). Из 1) и леммы 2.5.2 следует

$$\theta^{s_1(n-1)+1}(a) = \theta^{s_2(n-1)+1}(a),$$

откуда

$$\theta^{s_1(n-1)+1}(a)\theta^{-s_2(n-1)}(a) = \theta(a),$$

$$\theta^{(s_1-s_2)(n-1)+1}(a) = \theta(a). \qquad (*)$$

Снова применяя лемму 2.5.2, получаем

$$\theta(a^{[s_1-s_2]}) = \theta(a),$$

откуда и из предложения 1.3.5. следует 2).

2) $\Rightarrow$ 3). Из 2) следует (*), откуда

$$\theta^{-(s_1-s_2)(n-1)-1}(a) = \theta^{-1}(a),$$

$$\theta^{(s_2-s_1)(n-1)-1}(a)\theta^2(a) = \theta^{-1}(a)\theta^2(a);$$

$$\theta^{(s_2-s_1)(n-1)+1}(a) = \theta(a). \qquad (**)$$

Применяя к последнему равенству лемму 2.5.2, получаем

$$\theta(a^{[s_2-s_1]}a) = \theta(a),$$

откуда и из предложения 1.3.5 следует 3).

3) $\Rightarrow$ 1). Из 3) следует (**), откуда

$$\theta^{s_2(n-1)+1}(a)\theta^{-s_1(n-1)}(a) = \theta(a),$$

Из последнего равенства получаем

$$\theta^{s_2(n-1)+1}(a) = \theta^{s_1(n-1)+1}(a),$$

откуда

$$\theta(a^{[s_2]}) = \theta(a^{[s_1]}).$$

Применяя предложение 1.3.5, получаем 1). ∎

**2.5.7. Следствие.** Пусть a – элемент n-арной группы

< A, [ ] >. Если выполняется равенство 1) предыдущего предложения при $s_1 \neq s_2$, то существует такое целое положительное число r, что $a^{[r]} = a$.

**2.5.8. Следствие.** Пусть a – элемент n-арной группы < A, [ ] >. Если $a^{[1]} = a$, то и $a^{[-1]} = a$.

**2.5.9. Определение** [3, 4]. *Конечным n-адическим порядком элемента a* n-арной группы < A, [ ] > называется наименьшее целое положительное число m, для которого выполняется равенство $a^{[m]} = a$. Если же все n-адические степени элемента a различны, то a называется *элементом бесконечного n-адического порядка*. n-Адический порядок элемента a обозначают через |a|.

**2.5.10. Предложение.** Пусть < A, [ ] > – n-арная группа, $a \in A$. Тогда:

1) a имеет в < A, [ ] > – конечный n-адический порядок, равный m тогда и только тогда, когда элемент $\theta(a)$ имеет в $A^*$ конечный порядок, равный m(n - 1);

2) a имеет в < A, [ ] > бесконечный n-адический порядок тогда и только тогда, когда $\theta(a)$ имеет в $A^*$ бесконечный порядок.

*Доказательство.* 1) Так как |a| = m, то $a^{[m]} = a$ и все элементы a, $a^{[1]}$, …, $a^{[m-1]}$ – различны. Положим

$$\mathscr{B} = \mathscr{B}_1 \cup \mathscr{B}_2 \cup \ldots \cup \mathscr{B}_{m-1},$$

где

$$\mathscr{B}_1 = \{\theta(a), \theta(a^{[1]}), \ldots, \theta(a^{[m-1]})\},$$

$$\mathscr{B}_2 = \{\theta(aa), \theta(a^{[1]}a), \ldots, \theta(a^{[m-1]}a)\},$$

………………………………………

$$\mathscr{B}_{n-1} = \{\theta(a\underbrace{a\ldots a}_{n-2}), \theta(a^{[1]}\underbrace{a\ldots a}_{n-2}), \ldots, \theta(a^{[m-1]}\underbrace{a\ldots a}_{n-2}a)\}.$$

Ясно, что в каждом $\mathscr{B}_i$ (i = 1, 2, …, n – 1) все элементы различны, а так как, кроме того, $\mathscr{B}_i \cap \mathscr{B}_j = \varnothing$ для любых i, j ∈ {1, 2, …, n – 1} то все элементы в $\mathscr{B}$ различны.

Легко проверяется, что

$$\mathscr{B}_1 = \{\theta(a), \theta^n(a), \ldots, \theta^{(m-1)n+1}(a)\},$$

$$\mathscr{B}_2 = \{\theta^2(a), \theta^{n+1}(a), \ldots, \theta^{(m-1)n+2}(a)\},$$

$$\ldots\ldots\ldots\ldots\ldots\ldots\ldots\ldots\ldots\ldots\ldots\ldots\ldots$$

$$\mathscr{B}_{n-1} = \{\theta^{n-1}(a), \theta^{2(n-1)}(a), \ldots, \theta^{m(n-1)}(a)\}$$

и $\theta^{m(n-1)}(a) = E$ – единица группы $A^*$. Таким образом, установлено, что

$$\mathscr{B} = \{\theta(a), \theta^2(a), \ldots, \theta^{m(n-1)}(a) = E\}$$

и все степени, входящие в $\mathscr{B}$ различны. Следовательно, порядок $\theta(a)$ в $A^*$ равен $m(n-1)$.

Если теперь порядок $\theta(a)$ в $A^*$ равен $m(n-1)$, то что все элементы

$$\theta(a), \theta^n(a) = \theta(a^{[1]}), \ldots, \theta^{(m-1)(n-1)}(a) = \theta(a^{[m-1]})$$

различны и $\theta^{m(n-1)}(a) = E$ – единица в $A^*$. Отсюда следует, что все элементы $a, a^{[1]}, \ldots, a^{[m-1]}$ – различны и $a^{[m]} = a$, т. е. n-адический порядок $a$ в $< A, [\ ] >$ равен $m$.

3) Доказывается аналогично 1). ∎

**2.5.11. Лемма.** Для всякого элемента $a$ n-арной группы $< A, [\ ] >$ и любого целого $s$ справедливо равенство

$$\theta^s(\underbrace{a \ldots a}_{n-1}) = \theta(a^{[s-1]} \underbrace{a \ldots a}_{n-2}).$$

*Доказательство.* Согласно лемме 2.5.2,

$$\theta(a^{[s-2]}) = \theta^{(s-2)(n-1)-1}(a),$$

откуда, учитывая $a^{[0]} = a$ и применяя предложение 2.5.3, получим

$$\theta(a^{[s-2]})\underbrace{\theta(a)\ldots\theta(a)}_{2n-3} = \theta^{(s-2)(n-1)+1}(a)\underbrace{\theta(a)\ldots\theta(a)}_{2n-3},$$

$$\theta([a^{[s-2]}\underbrace{a^{[0]}\ldots a^{[0]}}_{n-1}]\underbrace{a\ldots a}_{n-2} = \theta^{s(n-1)}(a),$$

$$\theta(a^{[s-2+\underbrace{0+\ldots+0}_{n-1}+1]}\underbrace{a\ldots a}_{n-2} = \theta^s(\underbrace{a\ldots a}_{n-1}),$$

$$\theta(a^{[s-1]}\underbrace{a\ldots a}_{n-2}) = \theta^s(\underbrace{a\ldots a}_{n-1}). \blacksquare$$

Так как

$$\theta(a^{[s]}\underbrace{a\ldots a}_{n-2}) = \theta(a^{[t]}\underbrace{a\ldots a}_{n-2})$$

тогда и только тогда, когда $a^{[s]} = a^{[t]}$, то из леммы 2.5.11 вытекает

**2.5.12. Предложение.** Пусть $< A, [\ ] >$ – n-арная группа, $a \in A$. Тогда n-адический порядок элемента $a$ в $< A, [\ ] >$ и порядок элемента $\theta(\underbrace{a\ldots a}_{n-1})$ в $A^*$ совпадают.

**2.5.13. Лемма.** Пусть $< A, [\ ] >$ – n-арная группа, $a \in A$,

$$d = [\underbrace{a\ldots a}_{n}],\ d^k = \underbrace{d \text{\scriptsize @} \ldots \text{\scriptsize @} d}_{k},\quad k \in N.$$

Тогда $d^k = (a^{[1]})^{[k]} = a^{[k]}$.

***Доказательство.*** Равенство $d^1 = a^{[1]}$ следует из условия леммы и определения n-адической степени.

Предположим, что $d^{k-1} = a^{[k-1]}$. Тогда

$$d^k = d^{k-1} \text{\scriptsize @}\ d = [a^{[k-1]}\overline{a}\underbrace{a\ldots a}_{n-3}[\underbrace{a\ldots a}_{n}]] = [a^{[k-1]}\underbrace{a\ldots a}_{n-1}] = a^{[k]},$$

т. е. $d^k = a^{[k]}$. $\blacksquare$

Из леммы 2.5.13 вытекает

**2.5.14. Предложение.** Пусть $< A, [\ ] >$ – n-арная группа, $a \in A$, $d = [\underbrace{a \ldots a}_{n}]$. Тогда n-адический порядок элемента a в $< A, [\ ] >$ совпадает с порядком элемента d в $< A, @ >$.

**2.5.15. Предложение.** Если элемент a n-арной группы $< A, [\ ] >$ имеет конечный n-адический порядок m, то $a^{[s]} = a$ тогда и только тогда, когда s кратно m.

*Доказательство. Необходимость.* По предложению 2.5.14 порядок элемента $d = [\underbrace{a \ldots a}_{n}]$ в группе $< A, @ >$ равен m, а по лемме 2.5.13 $d^s = a^{[s]} = a$, где a – единица группы $< A, @ >$. Поэтому применяя соответствующий групповой результат, заключаем, что m делит s.

*Достаточность.* Так как порядок d равен m и делит s, то из соответствующего бинарного результата следует $d^s = a$, откуда и из леммы 2.5.13 вытекает $a^{[s]} = a$. ∎

Из предложений 2.5.6 и 2.5.15 вытекает

**2.5.16. Следствие.** Если элемент a n-арной группы $< A, [\ ] >$ имеет конечный n-адический порядок m, то $a^{[s]} = a^{[l]}$ тогда и только тогда, когда s – l кратно m.

Аналогично предложению 2.5.15 с использованием предложения 2.5.14, леммы 2.5.13 и соответствующего бинарного результата, доказывается следующее

**2.5.17. Предложение.** Если a – элемент конечного n-адического порядка m n-арной группы $< A, [\ ] >$ и s – любое целое число, то существует такое целое число $0 \leq r < m$, что $s = mq + r$ и $a^{[s]} = a^{[r]}$.

**2.5.18. Замечание.** При доказательстве предложений 2.5.15 и 2.5.17 можно вместо предложения 2.5.14 и леммы 2.5.13 использовать соответственно предложение 2.5.12 и лемму 2.5.11.

**2.5.19. Определение** [3, 4]. Пусть a – элемент n-арной группы < A, [ ] >. n-Арная подгруппа, порожденная одноэлементным множеством {a}, называется *циклической n-арной подгруппой n-арной группы* < A, [ ] >, *порожденной элементом* a, и обозначается через < ⟨a⟩, [ ] >; сам же элемент a называется *порождающим элементом* этой n-арной подгруппы. Если < A, [ ] > = < ⟨a⟩, [ ] >, то < A, [ ] > называется *циклической n-арной группой, порожденной элементом* a.

Из утверждения 1) предложения 2.5.3 и теоремы 2.5.5 следует

**2.5.20. Следствие.** Если a – элемент n-арной группы < A, [ ] >, то
$$< ⟨a⟩, [\ ] > = \{a^{[s]} \mid s \in \mathbf{Z}\}.$$

**2.5.21. Предложение.** Если a – элемент конечного n-адического порядка m n-арной группы < A, [ ] >, то $|⟨a⟩| = m$ и
$$⟨a⟩ = \{\, a^{[0]} = a, a^{[1]}, \ldots, a^{[m-1]}\,\}.$$

*Доказательство.* Из предложения 2.5.18 и следствия 2.5.20 следует
$$⟨a⟩ \subseteq \{\, a^{[0]} = a, a^{[1]}, \ldots, a^{[m-1]}\,\},$$

откуда, учитывая очевидное обратное включение, получаем требуемое равенство.

Согласно следствию 2.5.10, все степени в правой части доказанного равенства различны. Поэтому $|⟨a⟩| = m$. ∎

**2.5.22. Теорема.** n-Арная группа < A, [ ] > – является конечной (бесконечной) циклической, порожденной элементом a тогда и только тогда, когда группа $A^*$ – конечная (бесконечная) циклическая, порожденная элементом θ(a).

*Доказательство.* Сразу же заметим, что n-арная группа < A, [ ] > и группа $A^*$ конечны или бесконечны одновременно.

*Необходимость.* Так как
$$A = \{a^{[s]} \mid s \in \mathbf{Z}\},$$
то ввиду предложения 1.3.7,
$$A^* = \bigcup_{i=1}^{n-1}\{\theta(a^{[s]}\underbrace{a\ldots a}_{i-1}) \mid s \in \mathbf{Z}\},$$
откуда, учитывая лемму 2.5.2, имеем
$$A^* = \bigcup_{i=1}^{n-1}\{\theta(a^{s(n-1)+i}) \mid s \in \mathbf{Z}\}.$$

Следовательно, $A^*$ – циклическая порождённая элементом $\theta(a)$.

*Достаточность.* Пусть $A^*$ – циклическая группа, порождённая $\theta(a)$, и $b$ – произвольный элемент из $A$. Тогда, если $A^*$ конечная, то $\theta(b) = \theta^k(a)$ для некоторого целого $k > 0$, откуда
$$\theta(b) = \theta(\underbrace{a\ldots a}_{k}).$$

Последнее равенство возможно тогда и только тогда, когда $k = s(n-1)+1$ для некоторого целого $s \geq 0$. Тогда по лемме 2.5.2,
$$\theta(b) = \theta^k(a) = \theta^{s(n-1)+1}(a) = \theta(a^{[s]}),$$
откуда $b = a^{[s]}$. Следовательно, $<A, [\ ]>$ – циклическая, порождённая $a$.

Если $A^*$ – бесконечная, то $\theta(b) = \theta^k(a)$ для некоторого целого $k \neq 0$. Если $k > 0$, то аналогично конечному случаю доказывается равенство $b = a^{[s]}$. Если же $k = -t < 0$, где $t > 0$, то из $\theta(b) = \theta^{-t}(a)$ следует
$$\theta(b)\theta^{t+1}(a) = \theta^{-t}(a)\theta^{t+1}(a),$$
откуда

$$\theta(b\underbrace{a\ldots a}_{t+1}) = \theta(a)$$

Последнее равенство возможно тогда и только тогда, когда $t + 2 = r(n - 1) + 1$ для некоторого целого $r \geq 1$, откуда

$$-k + 2 = r(n - 1) + 1,$$
$$k = -r(n - 1) + 1,$$
$$k = s(n - 1) + 1, \ s = -r.$$

Тогда по лемме 2.5.2

$$\theta(b) = \theta^k(a) = \theta^{s(n-1)+1}(a) = \theta(a^{[s]}),$$

откуда $b = a^{[s]}$. Следовательно, $< A, [\ ] >$ – циклическая, порождённая $a$. ∎

**2.5.23. Теорема.** n-Арная группа $< A, [\ ] >$ – является конечной (бесконечной) циклической, порождённой элементом $a$, тогда и только тогда, когда группа $A_o$ – конечная (бесконечная) циклическая, порожденная элементом $\theta(\underbrace{a\ldots a}_{n-1})$.

*Доказательство.* Сразу же заметим, что множества $A$ и $A_o$ равномощны.

*Необходимость.* Так как $< A, [\ ] >$ – циклическая, то согласно утверждению 1) предложения 1.3.7, любой элемент группы $A_o$ может быть представлен в виде

$$u = \theta(a^{[s]}\underbrace{a\ldots a}_{n-2}), \ s \in \mathbf{Z},$$

откуда, ввиду леммы 2.5.2, имеем

$$u = \theta(a^{[s]})\theta^{n-2}(a) = \theta^{s(n-1)+1}(a)\theta^{n-2}(a) =$$
$$= \theta^{(s+1)(n-1)}(a) = (\theta^{n-1}(a))^{s+1} = (\theta(\underbrace{a\ldots a}_{n-1}))^{s+1}.$$

Следовательно, $A_o$ – циклическая группа, порождённая $\theta(\underbrace{a \ldots a}_{n-1})$.

*Достаточность*. Так как $A_o$ – циклическая, порождённая $\theta(\underbrace{a \ldots a}_{n-1})$, то согласно утверждению 4) теорема 1.4.2, для любого $b \in B$ элемент $\theta(b)$ может быть представлен в виде

$$\theta(b) = (\theta(\underbrace{a \ldots a}_{n-1}))^s \theta(a), \quad s \in \mathbf{Z},$$

откуда, снова применяя лемму 2.5.2, получаем

$$\theta(b) = (\theta^{n-1}(a))^s \theta(a) = \theta^{s(n-1)+1}(a) = \theta(a^{[s]}).$$

Это значит, что $b = a^{[s]}$, и $< A, [\ ] >$ – циклическая n-арная группа, порождённая a. ∎

Следствием леммы 2.5.13 является

**2.5.24. Теорема.** n-Арная группа $< A, [\ ] >$ является циклической, порождённой элементом a, тогда и только тогда, когда группа $< A, @ >$ – циклическая, порождённая элементом $[\underbrace{a \ldots a}_{n}]$.

**2.5.25. Лемма**. Определим на циклической группе $A = < a >$ n-арную операцию

$$[b_1 \ldots b_n] = b_1 \ldots b_n a.$$

Тогда справедливы следующие утверждения:

1) $< A, [\ ] >$ – циклическая n-арная группа, порождённая единицей e группы A;

2) если $|A| = n - 1$, то в $< A, [\ ] >$ нет n-арных подгрупп, отличных от $< A, [\ ] >$.

*Доказательство*. 1) Если $k > 0$, то для любого

$$a^k \in A = \{a, a^2, \ldots, a^k, \ldots\}$$

имеем
$$a^k = \underbrace{e\ldots e}_{n}a\underbrace{e\ldots e}_{n-1}a\ldots\underbrace{e\ldots e}_{n-1}a =$$
$$\underbrace{\phantom{a\underbrace{e\ldots e}_{n-1}a\ldots\underbrace{e\ldots e}_{n-1}a}}_{k-1}$$

$$= [\ldots[[\underbrace{e\ldots e}_{n}]\underbrace{e\ldots e}_{n-1}]\ldots\underbrace{e\ldots e}_{n-1}] = [\underbrace{e\ \ldots\ e}_{k(n-1)+1}] = e^{[k]},$$

т. е. $a^k = e^{[k]}$.

Если $k = 0$, то $a^k = e^{[k]} = e$.

Если $k < 0$, то
$$e = a^k a^{-k} = a^k \underbrace{\underbrace{e\ldots e}_{n-1}a\ldots\underbrace{e\ldots e}_{n-1}a}_{-k} =$$

$$= [a^k \underbrace{\underbrace{e\ldots e}_{n-1}\ldots\underbrace{e\ldots e}_{n-1}}_{-k}] = [a^k \underbrace{e\ldots e}_{-k(n-1)}],$$

т. е.
$$[a^k \underbrace{e\ldots e}_{k(n-1)}] = e,$$

и значит $a^k = e^{[k]}$. Таким образом, $a^k = e^{[k]}$ для любого целого $k$. Следовательно, $<A, [\ ]>$ – циклическая n-арная группа, порождённая элементом $e$.

2) Если $|A| = n - 1$, то для любого
$$a^k \in A = \{a, a^2, \ldots, a^{n-2}, a^{n-1} = e\}$$
имеем
$$[\underbrace{a^k \ldots a^k}_{n}] = \underbrace{a^k \ldots a^k}_{n} a = (a^k)^{n-1} a^k a = e a^{k+1} = a^{k+1},$$
т. е.
$$[\underbrace{a^k \ldots a^k}_{n}] = a^{k+1}.$$

Предположим, что $<B, [\ ]>$ – собственная n-арная подгруппа n-арной группы $<A, [\ ]>$, и $b$ – произвольный элемент из B.

Так как b = $a^k$ для некоторого k ∈ {1, 2, ..., n −1}, то согласно доказанному,

$$[\underbrace{b \ldots b}_{n}] = [\underbrace{a^k \ldots a^k}_{n}] = a^{k+1},$$

откуда, учитывая

$$[\underbrace{b \ldots b}_{n}] \in B,$$

получаем $a^{k+1}$ ∈ B. Аналогично из $a^{k+1}$ ∈ B следует $a^{k+2}$ ∈ B. Продолжая процесс, получим B = A. ∎

**2.5.26. Замечание**. Утверждение 2) предыдущего предложения является частным случаем следующего результата Поста, который мы докажем позже вместе с обратным к нему утверждением.

**2.5.27. Теорема** ([3], с.285). Всякая циклическая n-арная группа < A, [ ] > порядка |A|, где π(|A|) ⊆ π(n − 1), не содержит n-арных подгрупп, отличных от < A, [ ] >.

**2.5.28. Теорема** [3, 4]. Любые две циклические n-арные группы одной и той же мощности изоморфны.

*Доказательство.* Пусть < C = < c >, ⌊ ⌋ > – произвольная конечная порядка m (бесконечная) циклическая n-арная группа и A = < a > – конечная порядка m (бесконечная) циклическая группа. Определив на A n-арную операцию [ ] как в лемме 2.5.25, получим циклическую n-арную группу < A = < e >, [ ] >, где e – единица группы A.

Ясно, что отображение φ: $c^{[k]} \mapsto e^{[k]}$ является биекцией C на A. Если < C, ⌊ ⌋ > бесконечная, то

$$\varphi(\lfloor c^{[k_1]} \ldots c^{[k_n]} \rfloor) = \varphi(c^{[k_1 + \ldots + k_n + 1]}) = e^{[k_1 + \ldots + k_n + 1]} =$$

$$= [e^{[k_1]} \ldots e^{[k_n]}] = [\varphi(c^{[k_1]}) \ldots \varphi(c^{[k_n]})]$$

для любых $c^{[k_1]}, \ldots, c^{[k_n]}$ ∈ C, т.е. φ – изоморфизм.

Если же |A| = m, то $k_1 + \ldots + k_n + 1 = mq + r$, где 0 ≤ r < m. Тогда, применяя следствие 2.5.16, получим

$$\varphi(\lfloor c^{[k_1]} \ldots c^{[k_n]} \rfloor) = \varphi(c^{[k_1 + \ldots + k_n + 1]}) = \varphi(c^{[mq+r]}) = \varphi(c^{[r]}) = e^{[r]} =$$

$$= e^{[mq+r]} = e^{[k_1 + \ldots + k_n + 1]} = [e^{[k_1]} \ldots e^{[k_n]}] = [\varphi(c^{[k_1]}) \ldots \varphi(c^{[k_n]})],$$

т. е. $\varphi$ – изоморфизм. ∎

Согласно теореме 2.5.23, необходимым условием цикличности n-арной группы $< A, [\,] >$ является цикличность её соответствующей группы $A_o$. Однако, как показывает следующий пример, существуют нециклические n-арные группы, обладающие циклической соответствующей группой, т. е. цикличность группы $A_o$ не являются достаточным условием цикличности n-арной группы $< A, [\,] >$.

**2.5.29. Пример.** Пусть $< A = T_3 = \{(12), (13), (23)\}, [\,] > -$ тернарная группа нечетных подстановок на трех символах с тернарной операцией $[\,]$, производной от операции в симметрической группе $S_3$ (пример 1.1.8). Так как в $T_3$ все элементы идемпотенты, то $< A = T_3, [\,] > -$ нециклическая тернарная группа. А так как $|A_o| = |A| = 3$, то $A_o$ – циклическая группа.

Определим еще один n-арный аналог циклических групп.

**2.5.30. Определение.** n-Арная группа $< A, [\,] >$ называется полуциклической, если для любого $a \in A$ группа $< A, \circledast >$ – циклическая.

**2.5.31. Теорема.** Для n-арной группы $< A, [\,] >$ следующие утверждения равносильны:
1) $< A, [\,] >$ – полуциклическая;
2) группа $< A, \circledast >$ – циклическая для некоторого $a \in A$;
3) группа $A_o$ – циклическая;
4) некоторая соответствующая группа $\widetilde{A}_o$ – циклическая.

*Доказательство.* 1) $\Leftrightarrow$ 2) Согласно следствию 2.2.15, для любых $a, b \in A$ группы $< A, \circledast >$ и $< A, \circledcirc >$ изоморфны.

2) $\Leftrightarrow$ 3) Согласно предложению 1.6.1, группы $< A, \circledast >$ и $A_o$ изоморфны.

3) $\Leftrightarrow$ 4) Согласно утверждению 2) теоремы 1.4.9, группа

$A_o$ изоморфна любой соответствующей группе $\widetilde{A}_o$. ∎

Из теоремы 2.5.24 и утверждения 2) предыдущего предложения вытекает

**2.5.32. Следствие.** Всякая циклическая n-арная группа является полуциклической.

**2.5.33. Пример.** В тернарной группе $<B_n, [\,]>$ всех отражений правильного n-угольника, ввиду предложения 1.2.7, все элементы являются идемпотентами. Поэтому $<B_n, [\,]>$ – нециклическая тернарная группа. А так как циклическая группа $C_n$ всех поворотов правильного n-угольника является соответствующей для $<B_n, [\,]>$ (см. пример 1.4.23), то $<B_n, [\,]>$ – полуциклическая тернарная группа.

Из примеров 2.5.29 и 2.5.33 следует, что существуют полуциклические n-арные группы, не являющиеся циклическими.

Из этих примеров также вытекает, что существуют нециклические n-арные группы простого порядка, что невозможно в группах. Однако, имеет место очевидное

**2.5.34. Предложение.** Всякая n-арная группа простого порядка является полуциклической.

**2.5.35. Лемма.** Если $<A, [\,]>$ – n-арная группа порядка g, то $a^{[g]} = a$ для любого $a \in A$.

*Доказательство.* Пусть m – n-адический порядок элемента a, т.е. $a^{[m]} = a$, $1 \leq m \leq g$, откуда

$$[\underbrace{a \ldots a\, a}_{m(n-1)}] = a.$$

Следовательно, $\underbrace{a \ldots a}_{m(n-1)}$ – нейтральная последовательность. Так как m делит g, то $g = mk$, и поэтому

$$a^{[g]} = a^{[mk]} = [\underbrace{a \ldots a}_{mk(n-1)} a] = [\underbrace{a \ldots a}_{m(n-1)} \ldots \underbrace{a \ldots a}_{m(n-1)} a]_k = a,$$

т. е. $a^{[g]} = a$. ∎

**2.5.36. Лемма**. Если $< B, [\ ] >$ – собственная n-арная подгруппа конечной циклической n-арной группы $< A, [\ ] >$, то

$$(\frac{|A|}{|B|}, n-1) = 1.$$

*Доказательство.* Пусть $|B| = \gamma$, $|A| = g$, $1 < \gamma < g$, $A = \langle a \rangle$ и пусть b – произвольный элемент из B. Тогда по предыдущей лемме $b^{[\gamma]} = b$, а из цикличности $< A, [\ ] >$ следует $b = a^{[k]}$ для некоторого $k = 1, \ldots, g - 1$. Из последних двух равенств следует $(a^{[k]})^{[\gamma]} = a^{[k]}$, откуда, применяя предложение 2.5.3, получаем $a^{[k\gamma(n-1)+k+\gamma]} = a^{[k]}$. Так как порядок элемента a равен g, то по следствию 2.5.16

$$k\gamma(n-1) + k + \gamma - k = sg$$

для некоторого целого s, откуда

$$\gamma(k(n-1)+1) = s\gamma\frac{g}{\gamma}, \quad k(n-1) + 1 = s\frac{g}{\gamma}, \quad s\frac{g}{\gamma} - k(n-1) = 1.$$

Следовательно, $(\frac{g}{\gamma}, n-1) = 1$. ∎

Доказательство следующего известного результата из теории чисел заимствовано из книги В. Серпинского [31].

**2.5.37. Лемма.** Пусть $n - 1$, s и $\gamma$ – натуральные числа, причем $n - 1$ и s – взаимно простые. Тогда существует натуральное число t такое, что числа $s + t(n - 1)$ и $\gamma$ – взаимно простые.

*Доказательство.* Определим три числа: P – произведе-

ние всех простых делителей числа γ, являющихся делителями n – 1, причем P = 1, если таких делителей нет; Q – произведение всех простых делителей числа γ, являющихся делителями s, причем Q = 1, если таких делителей нет; R – произведение всех простых делителей числа γ, которые не являются делителями ни числа s, ни числа n – 1, причем R = 1, если таких делителей нет.

Так как (n – 1, s) = 1, то (P, Q) = 1. Ясно также, что (P, R) = 1 и (Q, R) = 1.

Допустим, что существует простое p такое, что

$$p \mid \gamma \text{ и } p \mid s + PR(n-1).$$

Если p | P, то из p | s + PR(n – 1) следует p | s, откуда p | Q, что противоречит (P, Q) = 1. Если же p | Q, то p | s, откуда, учитывая p | s + PR(n – 1), получаем p | PR(n – 1), что невозможно, так как

$$(s, P) = 1, (s, R) = 1 \text{ и } (s, n-1) = 1.$$

Если, наконец, p | R, то из p | s + PR(n – 1) следует p | s, откуда p | Q, что противоречит (R, Q) = 1.

Так как любой простой делитель p числа γ является делителем одного из чисел P, Q, R, то мы показали, что общих простых делителей числа γ и s + PR(n – 1) не имеют, т. е. (γ, s + t(n – 1)) =1, где t = PR. ∎

Теперь мы можем доказать цикличность n-арных подгрупп конечной циклической n-арной группы.

**2.5.38. Теорема** [3]. Любая n-арная подгруппа конечной циклической n-арной группы является циклической.

***Доказательство.*** Пусть < B, [ ] > – n-арная подгруппа циклической n-арной группы < A, [ ] >, причем

$$A = \langle a \rangle, |A| = g, |B| = \gamma, \ 1 \leq \gamma < g.$$

Если γ = 1, то < B, [ ] > – циклическая первого порядка. Поэтому считаем 1 < γ < g. По теореме 2.5.22 A* – циклическая группа, порожденная элементом $\theta(a)$. Тогда ее циклическая подгруппа $B^*(A)$ порождается элементом

$$\theta_A^{\frac{|A^*|}{|B^*(A)|}}(a) = \theta_A^{\frac{g(n-1)}{\gamma(n-1)}}(a) = \theta_A^{\frac{g}{\gamma}}(a) = u \in B^*(A).$$

При доказательстве леммы 2.5.36 установлено существование целых положительных $k$ и $s$ таких, что

$$s\frac{g}{\gamma} - k(n-1) = 1, \qquad (1)$$

откуда $(s, n-1) = 1$. Если $(s, \gamma) = 1$, то

$$(s, \gamma(n-1) = |B^*(A)|) = 1.$$

Следовательно, элемент

$$u^s = \theta_A^{s\frac{g}{\gamma}}(a) = \theta_A^{k(n-1)+1}(a) = \theta_A(a^{k(n-1)+1}) = \theta_A(a^{[k]})$$

порождает циклическую группу $B^*(A)$. По теореме 2.2.19 существует изоморфизм $\varphi$ группы $B^*$ на группу $B^*(A)$, по которому $\varphi(\theta_B(b)) = \theta_A(b)$. Поэтому $B^*$ – циклическая группа, порожденная элементом $\theta_B(a^{[k]})$.

Применяя теперь теорему 2.5.22, заключаем, что n-арная подгруппа $< B, [\ ] >$ – циклическая, порожденная элементом $a^{[k]}$.

Если $(s, \gamma) \neq 1$, то положим

$$s' = s + PR(n-1), \quad k' = k + PR\frac{g}{\gamma},$$

где $P$ и $R$ такие же, как в лемме 2.5.37. Так как согласно (1), $(s, n-1) = 1$, то по лемме 2.5.37 $(s', \gamma) = 1$. Из (1) легко также получается

$$s'\frac{g}{\gamma} - k'(n-1) = 1,$$

откуда $(s', n-1) = 1$. Таким образом,

$$(s', \gamma(n-1) = |B^*(A)|) = 1.$$

Следовательно, элемент

$$u^{s'} = \theta_A^{s'\frac{g}{\gamma}}(a) = \theta_A^{k'(n-1)+1}(a) = \theta_A(a^{[k']})$$

порождает циклическую группу $B^*(A)$. Далее как и выше показывается, что $<B, [\ ]>$ – циклическая, порожденная элементом $a^{[k']}$. ∎

**2.5.39. Теорема.** n-Арная подгруппа $<B, [\ ]>$ полуциклической n-арной группы $<A, [\ ]>$ является полуциклической.

*Доказательство.* По теореме 2.5.31 группа $A_o$ циклическая, а согласно замечанию 2.2.20 группа $B_o$ изоморфна подгруппе $B_o(A)$ группы $A_o$. Поэтому группа $B_o$ циклическая. Снова, применяя теорему 2.5.31, заключаем, что n-арная подгруппа $<B, [\ ]>$ – полуциклическая. ∎

Выясним теперь какими своими степенями может порождаться конечная циклическая n-арная группа.

**2.5.40. Предложение** [3]. Пусть $<A, [\ ]>$ – конечная циклическая n-арная группа, порождаемая элементом a. Тогда $<A, [\ ]>$ порождается элементом $a^{[k]}$ тогда и только тогда, когда

$$(k(n-1)+1, |A|) = 1. \qquad (1)$$

*Доказательство.* Так как $A = \langle a \rangle$, то по теореме 2.5.22

$$A^* = \langle \theta(a^{[k]}) = \theta(a^{k(n-1)+1}) = \theta^{k(n-1)+1}(a) \rangle.$$

Следовательно,

$$(k(n-1)+1, |A^*| = |A|(n-1)) = 1, \qquad (2)$$

откуда следует (1).

Обратно, если верно (1), то из $(k(n-1)+1, (n-1)) = 1$ следует (2). Следовательно, $A^* = \langle \theta(a^{[k]}) \rangle$, откуда, применяя теорему 2.5.22, получаем $A = \langle a^{[k]} \rangle$. ∎

**2.5.41. Следствие.** Множество всех порождающих конечной циклической n-арной группы $< A = \langle a \rangle, [\ ] >$ имеет вид

$$\{ a^{[k]} \mid 0 \leq k < |A|, (k(n-1)+1, |A|) = 1 \}.$$

**2.5.42. Лемма.** Пусть $< A = \langle a \rangle, [\ ] >$ – конечная циклическая n-арная группа порядка $g$, $k \geq 0$, $d$ – общий делитель чисел $k(n-1)+1$ и $g$, $g = d\gamma$. Тогда

$$(a^{[k]})^{[\gamma]} = a^{[k]}. \qquad (1)$$

*Доказательство.* Так как $d \mid k(n-1)+1$, то

$$k(n-1)+1 = ds \qquad (2)$$

для некоторого целого $s \geq 1$, откуда, учитывая $d = \dfrac{g}{\gamma}$, получаем

$$\gamma(k(n-1)+1) = gs. \qquad (3)$$

Из последнего равенства следует

$$\gamma k(n-1) + \gamma + k - k = gs. \qquad (4)$$

Тогда, применяя следствие 2.5.16, получаем

$$a^{[\gamma k(n-1) + \gamma + k]} = a^{[k]}, \qquad (5)$$

откуда следует (1). ∎

**2.5.43. Лемма.** Пусть $< A = \langle a \rangle, [\ ] >$ – конечная циклическая n-арная группа порядка $g$, $k \geq 0$, $\gamma$ – порядок элемента $a^{[k]}$. Тогда $\gamma = \dfrac{g}{d}$, где $d$ – общий делитель чисел $k(n-1)+1$ и $g$.

*Доказательство.* Так как γ – порядок элемента $a^{[k]}$, то верно (1) из предыдущей леммы, откуда следует справедливость равенств (5) – (3) из той же леммы. А так как порядок элемента совпадает с порядком циклической n-арной подгруппы, порождаемой этим элементом, то по теореме Лагранжа γ | g, т.е. g = dγ для некоторого d ≥ 1. Поэтому из (3) вытекает (2). Следовательно, $\gamma = \frac{g}{d}$, где d – общий делитель чисел k(n – 1) + 1 и g. ∎

**2.5.44. Теорема** [3]. Если $< A = \langle a \rangle, [\ ] >$ – конечная циклическая n-арная группа порядка g, k ≥ 0, то порядок элемента $a^{[k]}$ равен γ тогда и только тогда, когда

$$\gamma = \frac{g}{d},\ d = (k(n-1)+1, g). \qquad (*)$$

*Доказательство. Необходимость.* Так как порядок элемента $a^{[k]}$ равен γ, то по предыдущей лемме $\gamma = \frac{g}{d}$, где d – общий делитель чисел k(n – 1) + 1 и g. Предположим, что d < (k(n – 1) + 1, g) = d'. Тогда $\gamma' = \frac{g}{d'} < \gamma$ и по лемме 2.5.42

$$(a^{[k]})^{[\gamma']} = a^{[k]},$$

что противоречит определению порядка элемента. Следовательно, d = (k(n – 1) + 1, g).

*Достаточность.* Если имеет место (*), то по лемме 2.5.42

$$(a^{[k]})^{[\gamma]} = a^{[k]},$$

откуда следует, что порядок элемента $a^{[k]}$, равен γ' ≤ γ. Если γ' < γ, то

$$d' = \frac{g}{\gamma'} > \frac{g}{\gamma} = d = (k(n-1)+1, \gamma),$$

что невозможно, так как по лемме 2.5.43 d' – общий делитель чисел k(n – 1) + 1 и g. Следовательно, γ' = γ – порядок элемента $a^{[k]}$. ∎

Ясно, что предложение 2.5.40 является следствием теоремы 2.5.44.

**2.5.45. Теорема** [3]. Конечная циклическая n-арная группа < A = ⟨ a ⟩, [ ] > порядка g имеет элемент порядка γ тогда и только тогда, когда

$$\gamma \mid g, \ (\frac{g}{\gamma}, n-1) = 1. \tag{1}$$

***Доказательство.*** *Необходимость*. Если элемент $a^{[k]} \in A$ имеет порядок γ, то по теореме 2.5.44

$$\gamma \mid g, \ \frac{g}{\gamma} = (k(n-1)+1, g).$$

Из второго равенства следует $(\frac{g}{\gamma}, n-1) = 1$ и значит верно (1).

*Достаточность*. Так $d = \frac{g}{\gamma}$ взаимно просто с n – 1, то существуют целые $k_o$ и $s_o$ такие, что

$$k_o(n-1)+1 = s_o d, \tag{2}$$

т. е. ($k_o$, $s_o$) является частным решением диофантова уравнения

$$k(n-1)+1 = sd, \tag{3}$$

общее решение которого имеет вид

$$(k = k_o + td, \ s = s_o + t(n-1)), \ t - \text{целое}.$$

Согласно (2), ($s_o$, n – 1) = 1. Поэтому в арифметической прогрессии

$$s_o + (n-1), \ s_o + 2(n-1), \ \ldots, \ s_o + t(n-1), \ \ldots$$

бесконечно много простых чисел, среди которых можно выбрать простое число s = s₀ + t(n – 1) такое, что (s, γ) = 1, откуда и из (3) следует

$$(k(n-1)+1, g) = (sd, \gamma d) = d.$$

Таким образом,

$$\gamma = \frac{g}{d}, \quad d = (k(n-1)+1, g)$$

и по теореме 2.5.44 элемент $a^{[k]} \in A$ имеет порядок γ. ∎

**2.5.46. Следствие** [3]. Конечная циклическая n-арная группа порядка g имеет идемпотентный элемент тогда и только тогда, когда (g, n – 1) = 1.

Ясно, что идемпотент циклической n-арной группы является её единицей. А так как n-арная группа является производной от группы тогда и только тогда, когда она обладает единицей (пример 1.2.4, предложение 1.2.5), то из следствия 2.5.46 вытекает

**2.5.47. Следствие** [3]. Конечная циклическая n-арная группа порядка g является производной от группы тогда и только тогда, когда (g, n – 1) = 1.

Напомним, что через π(g) обозначается множество всех простых делителей числа g.

**2.5.48. Теорема** [3]. В циклической n-арной группе < A, [ ] > порядка g для любого ее делителя γ такого, что $(\frac{g}{\gamma}, n-1) = 1$ существует единственная n-арная подгруппа порядка γ. Других n-арных подгрупп < A, [ ] > не имеет.

*Доказательство.* Если γ делит g и $(\frac{g}{\gamma}, n-1) = 1$, то по теореме 2.5.45 в < A, [ ] > существует элемент $a^{[k]}$ порядка γ, который порождает циклическую n-арную подгруппу < B = ⟨ $a^k$ ⟩, [ ] > порядка γ.

Пусть теперь $< C, [\ ] >$ – произвольная n-арная подгруппа из $< A, [\ ] >$ того же порядка γ. Тогда по теореме 2.5.38 $< C = \langle a^m \rangle, [\ ] >$ – циклическая, порождаемая некоторым элементом $a^{[m]} \in A$ порядка γ, а по предложению 2.5.10 элемент $\theta(a^{[m]}) \in A^*$ имеет порядок $\gamma(n-1)$, который совпадает с порядком подгруппы $B^*(A)$ группы $A^*$. А так как по теореме 2.5.22 $A^*$ – циклическая группа, то $B^*(A)$ – единственная подгруппа в $A^*$ порядка $\gamma(n-1)$. Поэтому $\theta(a^{[m]}) \in B^*(A)$, откуда $a^{[m]} \in B$. Следовательно, $C \subseteq B$, а так как $|C| = |B|$, то $C = B$. Таким образом, $< B, [\ ] >$ – единственная n-арная подгруппа в $< A, [\ ] >$, имеющая порядок γ.

Если $< B', [\ ] >$ – n-арная подгруппа из $< A, [\ ] >$ порядка γ', то по теореме Лагранжа $\gamma' \mid g$, а по лемме 2.5.36

$$(\frac{g}{\gamma'}, n-1) = 1. \qquad \blacksquare$$

Если $\pi(g) \subseteq \pi(n-1)$, то для любого делителя γ числа g верно $(\frac{g}{\gamma}, n-1) \neq 1$. Поэтому теорема 2.5.27 является следствием теоремы 2.5.48.

Если $g = rs$, $(r, n-1) = 1$, $\pi(s) \subseteq \pi(n-1)$, то следующие условия равносильны:

1) $(\frac{g}{\gamma}, n-1) = 1$;

2) $\gamma = r's$, $r' \mid r$.

Поэтому из теоремы 2.5.48 вытекает

**2.5.49. Следствие** [3]. Если $< A, [\ ] >$ – конечная циклическая n-арная группа порядка $g = rs$, где $(r, n-1) = 1$, $\pi(s) \subseteq \pi(n-1)$, то для любого делителя r' числа r в $< A, [\ ] >$ существует единственная n-арная подгруппа порядка r's. Других n-арных подгрупп $< A, [\ ] >$ не имеет.

Делителю $r' = r$ в следствии 2.5.49 соответствует сама n-арная группа $< A, [\ ] >$. Поэтому имеет место

**2.5.50. Следствие** [3]. Если $<A, [\ ]>$ – конечная циклическая n-арная группа порядка $g = rs$, где $(r, n - 1) = 1$, $\pi(s) \subseteq \pi(n - 1)$, то число всех её собственных n-арных подгрупп равно $\tau(r) - 1$, где $\tau(r)$ – число всех делителей числа r. В частности, если $r = p_1^{\alpha_1} \ldots p_k^{\alpha_k}$ – разложение r на простые множители, то число всех собственных n-арных подгрупп в $<A, [\ ]>$ равно $(\alpha_1 + 1) \ldots (\alpha_k + 1) - 1$.

Если в следствии 2.5.50 положить $r = 1$, то

$$\pi(g) = \pi(s) \subseteq \pi(n - 1)$$

и число всех n-арных подгрупп в $<A, [\ ]>$ равно $\tau(1) - 1 = 0$. Следовательно, следствие 2.5.50 обобщает теорему 2.5.27.

Если в следствии 2.5.49 положить $s = 1$, то получим

**2.5.51. Следствие** [3]. Если $<A, [\ ]>$ – конечная циклическая n-арная группа порядка g, где $(g, n - 1) = 1$, то для любого делителя $\gamma$ её порядка g в $<A, [\ ]>$ существует единственная n-арная подгруппа порядка $\gamma$. Других n-арных подгрупп $<A, [\ ]>$ не имеет.

**2.5.52. Замечание**. Ясно, что условие $g = rs$, $(r, n - 1) = 1$, $\pi(s) \subseteq \pi(n - 1)$ в следствиях 2.5.49 и 2.5.50 можно заменить равносильным ему условием: r – наибольший взаимно простой с $n - 1$ делитель числа g.

**2.5.53. Следствие** [3]. Если r – наибольший взаимно простой с $n - 1$ делитель числа g, то в n-арной группе $<A, [\ ]>$ порядка g существует единственная n-арная подгруппа порядка $\frac{g}{r}$, которая содержится в любой другой n-арной подгруппе. В частности, если $(g, n - 1) = 1$, то в $<A, [\ ]>$ имеется единственная одноэлементная n-арная подгруппа, которая содержится в любой другой n-арной подгруппе.

**2.5.54. Теорема.** Конечная n-арная группа $<A, [\ ]>$ не имеет собственных n-арных подгрупп тогда и только тогда, когда она циклическая и $\pi(A) \subseteq \pi(n - 1)$.

*Доказательство. Необходимость.* Пусть a ∈ A и < B, [ ] > – циклическая n-арная подгруппа в < A, [ ] >, порожденная элементом a. Так как в < A, [ ] > нет собственных n-арных подгрупп, то B = A. Следовательно, < A, [ ] > – циклическая n-арная группа, порожденная любым своим элементом.

Пусть p ∈ π(A) и p ∉ π(n – 1), т. е. (p, n – 1) = 1. По теореме 2.5.48 в < A, [ ] > существует n-арная подгруппа порядка $\frac{|A|}{p}$, что невозможно. Таким образом, p ∈ π(n – 1) для любого p ∈ π(A), т. е. π(A) ⊆ π(n – 1).

*Достаточность.* Если в < A, [ ] > существует собственная n-арная подгруппа < B, [ ] > порядка γ, то по теореме 2.5.38 < B, [ ] > – циклическая. Следовательно элемент b, порождающий < B, [ ] >, имеет порядок γ. Согласно теореме 2.5.45, $(\frac{|A|}{\gamma}, n – 1) = 1$, что невозможно, так как $\pi(\frac{|A|}{\gamma}) \subseteq \pi(A) \subseteq \pi(n – 1)$. Таким образом, в < A, [ ] > нет собственных n-арных подгрупп. Теорема доказана.

**2.5.55. Теорема.** Если полуциклическая n-арная группа < A, [ ] > порядка g содержит идемпотент a, то для любого делителя γ числа g в < A, [ ] > существует единственная n-арная подгруппа < B, [ ] > порядка γ, содержащая идемпотент a.

*Доказательство.* Так как < A, [ ] > – полуциклическая n-арная группа, то < A, ⓐ > – циклическая группа порядка g, в которой для любого делителя γ числа g имеется единственная подгруппа < B, ⓐ > порядка γ. Ясно, что a ∈ B и < B, ⓐ > – характеристична в < A, ⓐ >. Тогда по следствию 2.2.11 < B, [ ] > – n-арная подгруппа в < A, [ ] >.

Предположим, что < C, [ ] > ещё одна подгруппа в < A, [ ] >, отличная от < B, [ ] >, имеющая порядок γ и содер-жащая a. Тогда в < A, @ > существуют две подгруппы < B, @ > и < C, @ > порядка γ, что невозможно. ∎

**2.5.56. Теорема.** Пусть < A, [ ] > – полуциклическая n-арная группа, содержащая идемпотент a. Тогда множество всех ее n-арных подгрупп, содержащих a, образует решетку, изоморфную решетке всех подгрупп циклической группы < A, @ >.

*Доказательство.* Ясно, что множество L(A, [ ], a) всех n-арных подгрупп n-арной группы < A, [ ] >, содержащих фиксированный элемент a, образует решетку с операциями пересечения и порождения n-арных подгрупп. При этом наименьший элемент этой решетки совпадает с < {a}, [ ] >, а наибольший – с самой n-арной группой < A, [ ] >.

Обозначим через L(A, @) решетку всех подгрупп группы < A, @ > и определим отображение $f : L(A, [\ ], a) \to L(A, @)$ по правилу

$$f : < B, [\ ] > \to < B, @ >.$$

Ясно, что f – инъекция.

Если < B, @ > – подгруппа группы < A, @ >, то по следствию 2.2.11 < B, [ ] > – n-арная подгруппа в < A, [ ] >, причем a ∈ B. Следовательно, f-сюръекция, а значит и биекция. Если < B, [ ] > и < C, [ ] > – n-арные подгруппы в < A, [ ] >, то B ⊆ C тогда и только тогда, когда B = f(B) ⊆ f(C) = C. Таким образом, f – изоморфизм решетки L(A, [ ], a) на решетку L(A, @). ∎

**2.5.57. Следствие.** Пусть < A, [ ] > – конечная циклическая n-арная группа, порядок g которой взаимно прост с n – 1. Тогда решетка всех n-арных подгрупп n-арной группы < A, [ ] > изоморфна решетке всех подгрупп циклической группы порядка g.

*Доказательство.* Согласно следствию 2.5.51, < A, [ ] > обладает единственным идемпотентом. Если < B, [ ] > – n-арная подгруппа в < A, [ ] >, то из |B| | g и (g, n – 1) = 1 следует

(|B|, n – 1) = 1. Кроме того, по теореме 2.5.38 < B, [ ] > – циклическая. Снова, применяя следствие 2.5.51, заключаем, что < B, [ ] > обладает единственным идемпотентом, который в таком случае должен совпадать с a. Таким образом, a содержится во всех n-арных подгруппах n-арной группы < A, [ ] >. Следовательно, решетка всех n-арных подгрупп n-арной группы < A, [ ] > совпадает с решеткой всех n-арных подгрупп n-арной группы < A, [ ] >, содержащих a. Так как всякая циклическая n-арная группа является полуциклической, то для завершения доказательства следствия достаточно применить теорему 2.5.55 к полуциклической n-арной группе < A, [ ] >. ∎

## §2.6. АБЕЛЕВЫ n-АРНЫЕ ГРУППЫ И ИХ ОБОБЩЕНИЯ

**2.6.1. Определение** [1]. n-Арная группа < A, [ ] > называется *абелевой*, если

$$[a_1 a_2 \ldots a_n] = [a_{\sigma(1)} a_{\sigma(2)} \ldots a_{\sigma(n)}]$$

для всех $a_1, \ldots, a_n \in A$ и любой подстановки $\sigma$ множества $\{1, \ldots, n\}$.

**2.6.2. Определение** [3]. n-Арная группа < A, [ ] > называется *m-полуабелевой*, если m – 1 делит n – 1 и последовательности

$$a a_1 \ldots a_{m-2} b, \qquad b a_1 \ldots a_{m-2} a$$

эквивалентны в < A, [ ] > для всех $a, a_1, \ldots, a_{m-2}, b \in A$.

Используя понятие эквивалентности последовательностей в n-арной группе, можно показать, что 2-полуабелевы n-арные группы это в точности абелевы n-арные группы [3, с. 217].

n-Полуабелевые n-арные группы называются также *полуабелевыми*, т. е. имеет место

**2.6.3. Определение** [1]. n-Арная группа $< A, [\ ] >$ называется полуабелевой, если

$$[aa_1 \ldots a_{n-2}b] = [ba_1 \ldots a_{n-2}a]$$

для всех $a, a_1, \ldots, a_{n-2}, b \in A$.

Ясно, что тождество определения 2.6.3 получается из тождества определения 2.6.1 при $\sigma = (1\ n)$. Поэтому абелевы n-арные группы являются полуабелевыми.

При $n = 2$ понятия абелевости, полуабелевости и m-полуабелевости совпадают. Если же $n > 2$, то все три указанные понятия различны. Например, легко проверяется [26], что тернарная группа $< B_n, [\ ] >$ всех отражений правильного n-угольника является полуабелевой, но не является абелевой.

**2.6.4. Теорема.** Для n-арной группы $< A, [\ ] >$ следующие утверждения равносильны:
1) $< A, [\ ] >$ – абелева;
2) универсальная обёртывающая группа Поста $A^* = < \mathcal{A}, * > $ – абелева;
3) любая обёртывающая группа $< \widetilde{A}, \bullet > $ – абелева;
4) некоторая обёртывающая группа $< \widetilde{A}, \bullet > $ – абелева.

*Доказательство.* 1) $\Rightarrow$ 2) Пусть

$$u = \theta(a_1 \ldots a_t), v = \theta(b_1 \ldots b_s),$$

где $s, t \in \{1, \ldots, n-1\}$, произвольные элементы из $A^*$. Если $s + t = i \leq n$, то, зафиксировав элементы $c_{i+1}, \ldots, c_n \in A$, и, используя абелевость $< A, [\ ] >$, получим

$$u * v * \theta(c_{i+1} \ldots c_n) = \theta(a_1 \ldots a_t)\theta(b_1 \ldots b_s)\theta(c_{s+t+1} \ldots c_n) =$$

$$= \theta(a_1 \ldots a_t b_1 \ldots b_s c_{s+t+1} \ldots c_n) = \theta([a_1 \ldots a_t b_1 \ldots b_s c_{s+t+1} \ldots c_n]) =$$

$$= \theta([b_1 \ldots b_s a_1 \ldots a_t c_{s+t+1} \ldots c_n]) = \theta(b_1 \ldots b_s a_1 \ldots a_t c_{s+t+1} \ldots c_n) =$$

$$= \theta(b_1 \ldots b_s)\theta(a_1 \ldots a_t)\theta(c_{s+t+1} \ldots c_n) = v * u * \theta(c_{i+1} \ldots c_n),$$

т. е
$$u * v * \theta(c_{i+1}\ldots c_n) = v * u * \theta(c_{i+1}\ldots c_n),$$

откуда $u * v = v * u$.

Если $n < s + t = i \leq 2(n - 1)$, то снова, используя абелевость $<A, [\ ]>$, получим

$$u * v = \theta(a_1\ldots a_t)\theta(b_1\ldots b_s) = \theta(a_1\ldots a_t b_1\ldots b_s) =$$

$$= \theta(\underbrace{[a_1\ldots a_t b_1 \ldots b_k}_{n}]b_{k+1}\ldots b_s) = \theta([b_1\ldots b_k a_1\ldots a_t]b_{k+1}\ldots b_s) =$$

$$= \theta(b_1[b_2\ldots b_k a_1\ldots a_t b_{k+1}]b_{k+2}\ldots b_s) =$$

$$= \theta(b_1[b_2\ldots b_k b_{k+1} a_1\ldots a_t]b_{k+2}\ldots b_s) =$$

$$= \theta(b_1\ldots b_{k+1} a_1\ldots a_t b_{k+2}\ldots b_s) = \ldots = \theta(b_1\ldots b_s a_1\ldots a_t) =$$

$$= \theta(b_1\ldots b_s)\theta(a_1\ldots a_t) = v * u,$$

т. е. $u * v = v * u$. Абелевость группы $A^*$ доказана.

2) $\Rightarrow$ 3) Применяется утверждение 1) теоремы 1.4.9, согласно которому любая обёртывающая группа $<\widetilde{A}, \bullet>$ является гомоморфным образом группы $A^*$.

3) $\Rightarrow$ 4) Очевидно.

4) $\Rightarrow$ 1) Из абелевости группы $<\widetilde{A}, \bullet>$ вытекает

$$a_1 \bullet a_2 \bullet \ldots \bullet a_u = a_{\sigma(1)} \bullet a_{\sigma(2)} \bullet \ldots \bullet a_{\sigma(u)}$$

для любых $a_1, a_2, \ldots, a_u \in A$ и любой подстановки $\sigma$ множества $\{1, 2, \ldots, n\}$, откуда

$$[a_1 a_2 \ldots a_n] = [a_{\sigma(1)}\ a_{\sigma(2)}\ \ldots\ a_{\sigma(n)}].$$

Следовательно, n-арная группа $<A, [\ ]>$ – абелева. ∎

**2.6.5. Теорема.** Для n-арной группы $<A, [\ ]>$ следующие утверждения равносильны:

1) $<A, [\ ]>$ – полуабелева;

2) соответствующая группа Поста $A_o$ – абелева;

3) любая соответствующая группа $\widetilde{A}_o$ – абелева;

4) некоторая соответствующая группа $A_o$ – абелева.

***Доказательство.*** 1)$\Rightarrow$2) Пусть $u$ и $v$ – произвольные элементы из $A_o$, которые согласно 1) предложения 1.3.7, можно представить в виде

$$u = \theta(ac_1 \ldots c_{n-2}),\ v = \theta(bc_1 \ldots c_{n-2}),$$

где $a, b, c_1, \ldots, c_{n-2} \in A$. Используя полуабелевость $< A, [\ ] >$, получаем

$$u * v = \theta(ac_1 \ldots c_{n-2})\theta(bc_1 \ldots c_{n-2}) = \theta(ac_1 \ldots c_{n-2}bc_1 \ldots c_{n-2}) =$$

$$= \theta([ac_1 \ldots c_{n-2}b]c_1 \ldots c_{n-2}) = \theta([bc_1 \ldots c_{n-2}a]c_1 \ldots c_{n-2}) =$$

$$= \theta([bc_1 \ldots c_{n-2}ac_1 \ldots c_{n-2}) = \theta(bc_1 \ldots c_{n-2})\theta(ac_1 \ldots c_{n-2}) = v * u,$$

т. е $u * v = v * u$. Следовательно, $A_o$ – абелева.

2)$\Rightarrow$3) Следует из утверждения 2) теоремы 1.4.9, согласно которому, группа $A_o$ изоморфна любой соответствующей группе $A_o$ n-арной группы $< A, [\ ] >$.

3)$\Rightarrow$4) Очевидно.

4)$\Rightarrow$1) Из абелевости группы $\widetilde{A}_o$ вытекает

$$(a\bullet c_1 \bullet \ldots \bullet c_{n-2})\bullet(b\bullet c_1 \bullet \ldots \bullet c_{n-2}) = (b\bullet c_1 \bullet \ldots \bullet c_{n-2})\bullet(a\bullet c_1 \bullet \ldots \bullet c_{n-2})$$

для любых $a, b, c_1, \ldots, c_{n-2} \in A$, откуда

$$(a\bullet c_1 \bullet \ldots \bullet c_{n-2} \bullet b)\bullet c_1 \bullet \ldots \bullet c_{n-2} = (b\bullet c_1 \bullet \ldots \bullet c_{n-2} \bullet a)\bullet c_1 \bullet \ldots \bullet c_{n-2},$$

$$[ac_1 \ldots c_{n-2}b] = [bc_1 \ldots c_{n-2}a],$$

т. е. n-арная группа $< A, [\ ] >$ – полуабелева. ∎

**2.6.6. Теорема.** n-Арная группа $< A, [\ ] >$ является m-полуабелевой тогда и только тогда, когда для любых $a, b \in A$ и некоторых $c_1, \ldots, c_{n-2} \in A$ последовательности.

$$ac_1^{m-2}b,\ \ bc_1^{m-2}a \tag{1}$$

эквивалентны в $< A, [\ ] >$.

*Доказательство.* Необходимость очевидна.

*Достаточность.* Пусть для любых a, b ∈ A и некоторых $c_1, \ldots, c_{m-2} \in A$ последовательности (1) эквивалентны в < A, [ ] >, откуда при $b = c_i$ следует эквивалентность последовательностей

$$ac_1^{m-2} c_i, \ c_i c_1^{m-2} a \ (i = 1, \ldots, m-2). \tag{2}$$

Из эквивалентности последовательностей (2) следует эквивалентность последовательностей

$$a\alpha_i, \beta_i a \ (i = 1, \ldots, m-2), \tag{3}$$

где

$$\alpha_i = \underbrace{c_1^{m-2} c_i \ldots c_1^{m-2} c_i}_{k-1}, \ \beta_i = c_i c_1^{m-2} \ldots c_i c_1^{m-2}.$$

Для любых $a, b, a_1, \ldots, a_{m-2} \in A$ и любого $i = 1, \ldots, m-2$ обозначим через $a_i^*$ решение уравнения

$$a_i = [c_i^{m-2} u_i \alpha_i c_1^i], \ \text{т. е.} \ a_i = [c_i^{m-2} a_i^* \alpha_i c_1^i].$$

Тогда для любого $c \in A$, используя эквивалентность последовательностей (1), а также эквивалентность последовательностей (3), будем иметь

$$[a a_1^{m-2} b \overset{n-m}{c}] = [a[c_1^{m-2} a_1^* \alpha_1 c_1] c_2^{m-2} a_2^* \alpha_2 c_1^2][c_3^{m-2} a_3^* \alpha_3 c_1^3] \ldots$$

$$\ldots [a_{m-3}^{m-2} a_{m-3}^* \alpha_{m-3} c_1^{m-3}][c_{m-2} a_{m-2}^* \alpha_{m-2} c_1^{m-2}] b \overset{n-m}{c}] =$$

$$= [a c_1^{m-2} a_1^* \alpha_1 c_1^{m-2} a_2^* \alpha_2 c_1^{m-2} a_3^* \alpha_3 \ldots$$

$$\ldots a_{m-3}^* \alpha_{m-3} c_1^{m-2} a_{m-2}^* \alpha_{m-2} c_1^{m-2} b \overset{n-m}{c}] =$$

$$= [a_1^* c_1^{m-2} a \alpha_1 c_1^{m-2} a_2^* \alpha_2 \ldots \underbrace{c_1^{m-2} a_{m-2}^* \alpha_{m-2} c_1^{m-2} b \overset{n-m}{c}}_{\alpha}] =$$

$$= [a_1^* c_1^{m-2} \beta_1 a c_1^{m-2} a_2^* \alpha_2 \ldots \alpha] = [a_1^* c_1^{m-2} \beta_1 a_2^* c_1^{m-2} a \alpha_2 \ldots \alpha] =$$

$$= [\,a_1^* c_1^{m-2} \beta_1 a_2^* c_1^{m-2} \beta_2 a \ldots \alpha\,] = \ldots$$

$$\ldots [\underbrace{a_1^* c_1^{m-2} \beta_1 a_2^* c_1^{m-2} \beta_2 a_3^*}_{\beta} \ldots \beta_{m-3} a c_1^{m-2} a_{m-2}^* \alpha_{m-2} c_1^{m-2} b \overset{n-m}{c}\,] =$$

$$= [\,\beta \ldots \beta_{m-3} a_{m-2}^* c_1^{m-2} a \alpha_{m-2} c_1^{m-2} b \overset{n-m}{c}\,] =$$

$$= [\,\beta \ldots \beta_{m-3} a_{m-2}^* c_1^{m-2} \beta_{m-2} a c_1^{m-2} b \overset{n-m}{c}\,] =$$

$$= [\,\beta \ldots \beta_{m-3} a_{m-2}^* c_1^{m-2} \beta_{m-2} b c_1^{m-2} a \overset{n-m}{c}\,] =$$

$$= [\,\beta \ldots \beta_{m-3} a_{m-2}^* c_1^{m-2} b \alpha_{m-2} c_1^{m-2} a \overset{n-m}{c}\,] =$$

$$= [\,\beta \ldots \beta_{m-3} b c_1^{m-2} a_{m-2}^* \alpha_{m-2} c_1^{m-2} a \overset{n-m}{c}\,] =$$

$$= [\,\beta \ldots \underbrace{b \alpha_{m-3} c_1^{m-2} a_{m-2}^* \alpha_{m-2} c_1^{m-2} a \overset{n-m}{c}}_{\gamma}\,] = \ldots$$

$$\ldots = [\,a_1^* c_1^{m-2} \beta_1 a_2^* c_1^{m-2} b \alpha_2 \ldots \gamma\,] =$$

$$= [\,a_1^* c_1^{m-2} \beta_1 b c_1^{m-2} a_2^* \alpha_2 \ldots \gamma\,] = [\,a_1^* c_1^{m-2} b \alpha_1 c_1^{m-2} a_2^* \alpha_2 \ldots \gamma\,] =$$

$$= [\,b c_1^{m-2} a_1^* \alpha_1 c_1^{m-2} a_2^* \alpha_2 c_1^{m-2} a_3^* \alpha_3 \ldots$$

$$\ldots c_1^{m-2} a_{m-2}^* \alpha_{m-2} c_1^{m-2} a \overset{n-m}{c}\,] =$$

$$= [\,b[c_1^{m-2} a_1^* \alpha_1 c_1][c_2^{m-2} a_2^* \alpha_2 c_1^2] \ldots$$

$$\ldots [c_{m-2} a_{m-2}^* \alpha_{m-2} c_1^{m-2}] a \overset{n-m}{c}\,] = [\,b a_1^{m-2} a \overset{n-m}{c}\,],$$

т. е.

$$[\,a a_1^{m-2} b \overset{n-m}{c}\,] = [\,b a_1^{m-2} a \overset{n-m}{c}\,].$$

Следовательно, последовательности $a a_1^{m-2} b$, $b a_1^{m-2} a$ эквивалентны, а n-арная группа $< A, [\ ] > -$ m-полуабелева. ∎

Полагая в теореме 2.6.6 $c_1 = \ldots = c_{m-2} = c$, получим

**2.6.7. Следствие** [19, 32]. n-Арная группа $< A, [\ ] >$ является m-полуабелевой тогда и только тогда, когда для любых $a, b \in A$ и некоторого $c \in A$ последовательности

$$a\underbrace{c\ldots c}_{m-2}b,\ b\underbrace{c\ldots c}_{m-2}a$$

эквивалентны в $< A, [\ ] >$.

Полагая в теореме 2.6.6. $m = n$, получим

**2.6.8. Следствие**. n-Арная группа $< A, [\ ] >$ является полуабелевой тогда и только тогда, когда

$$[ac_1 \ldots c_{n-2}b] = [bc_1 \ldots c_{n-2}a]$$

для любых $a, b \in A$ и некоторых $c_1, \ldots, c_{n-2} \in A$.

Полагая в следствии 2.6.7 $m = n$ или в следствии 2.6.8 $c_1 = \ldots = c_{n-2} = c$, получим

**2.6.9. Следствие** [33]. n-Арная группа $< A, [\ ] >$ является полуабелевой тогда и только тогда, когда

$$[a\underbrace{c\ldots c}_{n-2}b] = [b\underbrace{c\ldots c}_{n-2}a]$$

для любых $a, b \in A$ и некоторого $c \in A$.

**2.6.10. Следствие** [34]. n-Арная группа $< A, [\ ] >$ является полуабелевой тогда и только тогда, когда она удовлетворяет одному из следующих тождеств:

$$[a_1 \ldots a_{n-1}b_1 \ldots b_{n-1}c] = [b_1 \ldots b_{n-1}a_1 \ldots a_{n-1}c]; \qquad (1)$$

$$[ca_1 \ldots a_{n-1}b_1 \ldots b_{n-1}] = [cb_1 \ldots b_{n-1}a_1 \ldots a_{n-1}]. \qquad (2)$$

*Доказательство.* Для фиксированного $d \in A$, согласно 1) предложения 1.3.7, существуют $a, b \in A$ такие, что

$$a_1 \ldots a_{n-1} \theta a \underbrace{d \ldots d}_{n-2}, \quad b_1 \ldots b_{n-1} \theta b \underbrace{d \ldots d}_{n-2}.$$

*Необходимость.* Применяя для полуабелевой n-арной группы $< A, [\ ] >$ следствие 2.6.9, получаем

$$[a_1 \ldots a_{n-1} b_1 \ldots b_{n-1} c] = [a \underbrace{d \ldots d}_{n-2} b \underbrace{d \ldots d}_{n-2} c] =$$

$$= [b \underbrace{d \ldots d}_{n-2} a \underbrace{d \ldots d}_{n-2} c] = [b_1 \ldots b_{n-1} a_1 \ldots a_{n-1} c],$$

т. е. верно (1).

*Достаточность.* Если на n-арной группе $< A, [\ ] >$ выполнено тождество (1), то на $< A, [\ ] >$ выполнено также тождество

$$[a \underbrace{d \ldots d}_{n-2} b \underbrace{d \ldots d}_{n-2} c] = [b \underbrace{d \ldots d}_{n-2} a \underbrace{d \ldots d}_{n-2} c],$$

а значит, и тождество

$$[a \underbrace{d \ldots d}_{n-2} b] = [b \underbrace{d \ldots d}_{n-2} a].$$

Тогда по следствию 2.6.9, $< A, [\ ] >$ – полуабелева. Для второго тождества доказательство проводится аналогично. ∎

**2.6.11. Теорема.** Для n-арной группы $< A, [\ ] >$ следующие утверждения равносильны:
1) $< A, [\ ] >$ – полуабелева;
2) для любого $a \in A$ группа $< A, @ >$ – абелева;
3) для некоторого $a \in A$ группа $< A, @ >$ – абелева.

*Доказательство.* 1)$\Rightarrow$2) Для любого $a \in A$ зафиксируем обратную последовательность $a_1 \ldots a_{n-2}$. Из полуабелевости $< A, [\ ] >$ следует

$$[x a_1 \ldots a_{n-2} y] = [y a_1 \ldots a_{n-2} x] \tag{1}$$

для любых $x, y \in A$, откуда

$$x @ y = y @ x. \tag{2}$$

Следовательно, < A, @ > – абелева.

2)⇒3) Очевидно.

3)⇒1) Так как < A, @ > – абелева, то верно (2), откуда следует (1). Тогда по следствию 2.6.8 < A, [ ] > – полуабелева. ∎

Следующее определение является n-арной версией соответствующего определения для произвольных универсальных алгебр [35, с.32].

**2.6.12. Определение.** n-Арная группа < A, [ ] > называется *коммутативной*, если она удовлетворяет тождеству

$$[[a_{11}a_{12} \ldots a_{1n}][a_{21}a_{22} \ldots a_{2n}] \ldots [a_{n1}a_{n2} \ldots a_{nn}]] =$$
$$= [[a_{11}a_{21} \ldots a_{n1}][a_{12}a_{22} \ldots a_{n2}] \ldots [a_{1n}a_{2n} \ldots a_{nn}]].$$

**2.6.13. Теорема** [34, 36]. n-Арная группа < A, [ ] > является полуабелевой тогда и только тогда, когда она коммутативна.

*Доказательство. Необходимость.* Заметим, что правая часть тождества в определении 2.6.12 получается из левой части, если в ней поменять местами элементы $x_{ij}$ и $x_{ji}$, т. е. в левой части тождества в определении 2.6.12 внешняя n-арная операция применяется к строкам матрицы $(a_{ij})_{n \times n}$, а в правой части тождества определения 2.6.12 внешняя n-арная операция применяется к столбцам этой же матрицы.

Найдём длину последовательности $a_{ij} \ldots a_{ji}$ в последовательности

$$a_{11}a_{12} \ldots a_{1n}a_{21}a_{22} \ldots a_{2n} \ldots a_{n1}a_{n2} \ldots a_{nn}.$$

Для этого рассмотрим расположение элементов $a_{ij}$ и $a_{ji}$ в матрице $(a_{ij})_{n \times n}$, полагая для определенности $i < j$. В i-ой строке в столбцах $j, j + 1, \ldots, n$ расположено $n - j + 1$ элементов. Между i-ой и j-ой строками расположено $j - i - 1$ строк, содержащих по n элементов. В j-ой строке в столбцах $1, 2, \ldots, i$ расположено i элементов. Поэтому длина последовательно-сти $x_{ij} \ldots x_{ji}$ равна
$$n - j + 1 + (j - i - 1)n + i = n + (j - i - 1)n - (j - i - 1) =$$
$$= n + (j - i - 1)(n - 1) = (j - 1)(n - 1) + 1.$$

В частности, при $i = j$ последовательность $x_{ij} \ldots x_{ji} = x_{ij}$ состоит из одного элемента.

Так как $< A, [\ ] >$ полуабелевая n-арная группа и длина последовательности $a_{ij} \ldots a_{ji}$ равна $(j - i)(n - 1) + 1$, то, используя ассоциативность операции $[\ ]$, можно так переставить элементы $a_{ij}$ и $a_{ji}$, что будет выполнятся тождество из определения 2.6.12. Следовательно, $< A, [\ ] >$ — коммутативная n-арная группа.

*Достаточность.* Пусть $a_1, \ldots, a_n, a \in A$. Так как

$$[[aa_1\ldots a_{n-1}][a_n\underbrace{a\ldots a}_{n-1}][a_2\underbrace{a\ldots a}_{n-1}]\ldots[a_{n-1}\underbrace{a\ldots a}_{n-1}]] =$$

$$= [[aa_na_2\ldots a_{n-1}][a_1\underbrace{a\ldots a}_{n-1}][a_2\underbrace{a\ldots a}_{n-1}]\ldots[a_{n-1}\underbrace{a\ldots a}_{n-1}]],$$

то, используя ассоциативность n-арной операции, получаем

$$[a[a_1a_2\ldots a_{n-1}a_n][\underbrace{a\ldots a}_{n-1}a_2\underbrace{a\ldots a}_{n-1}\ldots a_{n-1}aa]\underbrace{a\ldots a}_{n-3}] =$$

$$= ]a]a_na_2\ldots a_{n-1}a_1][\underbrace{a\ldots a}_{n-1}a_2\underbrace{a\ldots a}_{n-1}\ldots a_{n-1}aa]\underbrace{a\ldots a}_{n-3}],$$

откуда
$$[a_1a_2\ldots a_{n-1}a_n] = [a_n\, a_2\, \ldots a_{n-1}\, a_1],$$

т. е. $< A, [\ ] >$ полуабелева. ∎

**2.6.14. Предложение** [32]. Если $n = k(m - 1) + m$, $k \geq 0$, то n-арная группа $< A, [\ ] >$ является (k+2)-полуабелевой тогда и только тогда, когда $< A, [\ ]_{m, a_1^k} >$ — абелева m-арная группа для всех $a_1, \ldots, a_к \in A$, где

$$[x_1x_2\ldots x_m]_{m, a_1^k} = [x_1 a_1^k x_2 a_1^k \ldots a_1^k x_m].$$

*Доказательство*. Пусть $< A, [\ ] > -$ (k+2)-полуабелева n-арная группа. Тогда
$$aa_1^k b \ \theta \ ba_1^k a \qquad (1)$$
для всех $a_1, \ldots, a_к \in A$, откуда
$$[aa_1^k b\ a_1^k b_1 \ldots a_1^k b_{m-2}] = [ba_1^k aa_1^k b_1 \ldots a_1^k b_{m-2}], \qquad (2)$$
где $b_1, \ldots, b_{m-2}$ – произвольные элементы из A. Из (2) получаем
$$[abb_1^{m-2}]_{m,a_1^k} = [bab_1^{m-2}]_{m,a_1^k}, \qquad (3)$$
откуда
$$ab \ \theta \ ba. \qquad (4)$$
Следовательно, $< A, [\ ]_{m,a_1^k} > -$ 2-абелева, а значит и абелева m-арная группа.

Обратно, если $< A, [\ ]_{m,a_1^k} > -$ абелева m-арная группа, то последовательно выполняются (4), (3), (2) и (1). ∎

Из теорем 2.5.20 и 2.6.4 вытекает

**2.6.15. Следствие.** Всякая циклическая n-арная группа является абелевой.

Это следствие может быть также выведено из определений циклической и абелевой n-арной группы.

Из определения 2.5.30 и теоремы 2.6.5 вытекает

**2.6.16. Следствие.** Всякая полуциклическая n-арная группа является полуабелевой.

**2.6.17. Определение** [32]. n-Арная группа $< A, [\ ] >$ называется *слабо m-полуабелевой* (m − 1 делит n − 1], если последовательности
$$a\underbrace{b\ldots b}_{m-1}, \ \underbrace{b\ldots b}_{m-1}a$$
эквивалентны в $< A, [\ ] >$ для любых $a, b \in A$.

Слабо 2-полуабелевы n-арные группы – это в точности абелевы n-арные группы. Слабо n-полуабелевы n-арные группы называются также слабо полуабелевыми, т. е. имеет место

**2.6.18. Определение** [32]**.** n-Арная группа $< A, [\ ] >$ называется *слабо полуабелевой*, если в ней выполняется тождество

$$[a\underbrace{b\ldots b}_{n-1}] = [\underbrace{b\ldots b}_{n-1}a].$$

Из определения 2.6.18 вытекает, что любая m-полуабелевая n-арная группа является слабо m-полуабелевой.

Следующие примеры показывают, что в общем случае класс всех слабо m-полуабелевых n-арных групп шире класса всех m-полуабелевых n-арных групп.

**2.6.19. Пример**. Пусть $< A, [\ ] >$ – n-арная группа кватернионов (пример 1.2.8). Так как

$$x^{4k} = x^{n-1} = 1$$

для любого $x \in A$, то

$$x^{n-1}ya^2 = yx^{n-1}a^2$$

для всех $x, y \in A$, откуда

$$[y\underbrace{x\ldots x}_{n-1}] = [\underbrace{x\ldots x}_{n-1}y],$$

и n-арная группа $< A, [\ ] >$ является слабо полуабелевой. А так как

$$[a\underbrace{1\ldots 1}_{4k-1}b] = ba \neq ba \cdot a^2 = [b\underbrace{1\ldots 1}_{4k-1}a],$$

то n-арная группа $< A, [\ ] >$ не является полуабелевой.

**2.6.20. Пример**. Определим на n-арной группе $< A, [\ ] >$ кватернионов $\nu$-арную операцию $\lfloor\ \rfloor$, $\nu = t(n-1) + 1$, $t > 1$, производную от n-арной операции $[\ ]$. Из примера 2.6.19 вытекает, что $\nu$-арная группа $< A, \lfloor\ \rfloor >$ не является n-полуабелевой, но является слабо n-полуабелевой.

Справедливость следующей леммы устанавливается простой проверкой.

**2.6.21. Лемма**. Если для любых элементов a и b n-арной группы $\langle A, [\ ] \rangle$ и любого натурального $m \geq 1$ последовательности

$$a \overset{m-1}{b}, \quad \overset{m-1}{b} a$$

эквивалентны, то для любого $k \geq 1$ эквивалентны также последовательности

$$a \overset{k(m-1)}{b}, \quad \overset{k(m-1)}{b} a.$$

**2.6.22. Лемма.** Пусть $\langle A, [\ ] \rangle$ – n-арная группа, $k > m \geq 1$, $a, b \in A$, и выполняются следующие условия:

$$a \overset{m-1}{b} \theta \overset{m-1}{b} a, \quad a \overset{k-1}{b} \theta \overset{k-1}{b} a. \qquad (*)$$

Тогда

$$a \overset{k-m}{b} \theta \overset{k-m}{b} a.$$

*Доказательство.* Пусть $t \geq 1$, $t(n-1) + 1 \geq k$. Так как

$$[a \overset{k-m}{b} \overset{t(n-1)-k+m}{b}] = [a \overset{k-1}{b} \overset{t(n-1)-k+1}{b}] =$$

$$= [\overset{k-1}{b} a \overset{t(n-1)-k+1}{b}] = [\overset{k-m}{b} \overset{m-1}{b} a \overset{t(n-1)-k+1}{b}] =$$

$$= [\overset{k-m}{b} a \overset{m-1}{b} \overset{t(n-1)-k+1}{b}] = [\overset{k-m}{b} a \overset{t(n-1)-k+m}{b}],$$

то последовательности

$$a \overset{k-m}{b}, \quad \overset{k-m}{b} a$$

эквивалентны в $\langle A, [\ ] \rangle$. ∎

**2.6.23. Лемма.** Если в n-арной группе $\langle A, [\ ] \rangle$ для элементов $a, b \in A$ выполняются условия (*), то

$$a \overset{t-1}{b} \theta \overset{t-1}{b} a, \qquad (**)$$

где t − 1 = (m − 1, k − 1).

***Доказательство***. Пусть t − 1 = (m − 1, k − 1). Тогда существуют целые числа α и β такие, что
$$\alpha(m-1) + \beta(k-1) = t - 1.$$

Пусть для определенности α > 0, β < 0, т. е.
$$\alpha(m-1) > -\beta(k-1).$$

Тогда по лемме 2.6.21 для a, b ∈ A эквивалентны последовательности

$$\overset{\alpha(m-1)}{a}\overset{\alpha(m-1)}{b}, \quad \overset{\alpha(m-1)}{b}\overset{}{a},$$

а также последовательности

$$\overset{-\beta(k-1)}{a}\overset{-\beta(k-1)}{b}, \quad \overset{-\beta(k-1)}{b}\overset{}{a}.$$

Применяя теперь лемму 2.6.22, получаем эквивалентность последовательностей

$$\overset{\alpha(m-1)-(-\beta(k-1))}{a}\overset{\alpha(m-1)-(-\beta(k-1))}{b}, \quad \overset{}{b}\overset{}{a},$$

$$\overset{\alpha(m-1)+\beta(k-1)}{a}\overset{\alpha(m-1)+\beta(k-1)}{b}, \quad \overset{}{b}\overset{}{a},$$

откуда следует эквивалентность последовательностей (∗∗). ■

Из лемм 2.6.21 − 2.6.23 вытекает

**2.6.24. Теорема** [32]**.** Если n-арная группа является слабо m-полуабелевой и слабо k-полуабелевой, то она является слабо t-полуабелевой, где t − 1 = (m − 1, k − 1).

**2.6.25. Следствие**. Любая слабо m-полуабелева n-арная группа < A, [ ] > является слабо k-полуабелевой, где
$$k - 1 = (m - 1, k - 1).$$

***Доказательство***. Так как m – 1 делит n – 1, то

$$n - 1 = t(m - 1).$$

По лемме 2.6.21 последовательности $a\overset{n-1}{b}$, $\overset{n-1}{b}a$ эквивалентны в $< A, [\ ] >$ для всех $a, b \in A$. Тогда по теореме 2.6.24 эквивалентны последовательности $a\overset{k-1}{b}$, $\overset{k-1}{b}a$. ∎

## § 2.7. ПРОИЗВЕДЕНИЯ n-АРНЫХ ГРУПП

Декартово произведение n-арных групп является частным случаем соответствующего определения для произвольных универсальных алгебр.

Пусть $\mathcal{A}_1 = < A_1, [\ ]_1 >, \ldots, \mathcal{A}_k = < A_k, [\ ]_k >$ – n-арные группы и

$$A = A_1 \times \ldots \times A_k = \{(a_1, \ldots, a_k) \mid a_1 \in A_1, \ldots, a_k \in A_k\}$$

– декартово произведение их носителей. Определим на $A$ n-арную операцию по формуле

$$[(a_1^{(1)}, \ldots, a_k^{(1)}) \ldots (a_1^{(n)}, \ldots, a_k^{(n)})] = ([a_1^{(1)} \ldots a_1^{(n)}], \ldots, [a_k^{(1)} \ldots a_k^{(n)}]).$$

Легко проверяется, что алгебра $\mathcal{A} = < A, [\ ] >$ является n-арной группой с косым элементом

$$\overline{(a_1, \ldots, a_k)} = (\overset{1}{\overline{a}}_1, \ldots, \overset{k}{\overline{a}}_k),$$

где $\overset{i}{\overline{\phantom{a}}}$ – унарная операция взятия косого элемента в n-арной группе $\mathcal{A}_i$ $(i = 1, \ldots, k)$.

Построенная таким образом n-арная группа $\mathcal{A} = < A, [\ ] >$ называется *декартовым* или *внешним прямым* произведением n-арных групп $\mathcal{A}_1, \ldots, \mathcal{A}_k$ и обозначается через $\mathcal{A}_1 \times \ldots \times \mathcal{A}_k$.

Определение декартового произведения конечного числа n-арных групп можно распространить на случай произвольного семейства n-арных групп.

Пусть $\{\mathcal{A}_i = \langle A_i, [\ ]_i \rangle \mid i \in I\}$ – непустое семейство n-арных групп,
$$A = \prod_{i \in I} A_i = \{a: I \to \bigcup_{i \in I} A_i \mid a(i) \in A_i\}$$
– декартово произведение их носителей, $[\ ]$ – n-арная операция, определенная на $A$ покомпонентно:

$$[a_1 a_2 \ldots a_n](i) = [a_1(i) a_2(i) \ldots a_n(i)]_i,\ i \in I$$

для всех $a_1, a_2, \ldots, a_n \in A$. Тогда алгебра $\mathcal{A} = \langle A, [\ ] \rangle = \prod_{i \in I} \mathcal{A}_i$ называется *декартовым* или *полным прямым* произведением n-арных групп $\mathcal{A}_i$ ($i \in I$), и, как несложно проверить (см., например, предложение 5.1 из [4]), является n-арной группой с косым элементом $\bar{a}$, где $\bar{a}(i) = \overset{i}{\bar{a}_i}$.

Если $I = \{1, \ldots, k\}$, то положив $a(i) = a_i$ и отождествив элементы $a$ и $(a_1, \ldots, a_k)$, видим, что последнее определение включает в себя случай конечного числа сомножителей.

**2.7.1. Теорема.** Если $\mathcal{A}_i = \langle A_i, [\ ]_i \rangle$ – n-арные группы, производные от групп $A_i$ ($i \in I$), то n-арная группа $\mathcal{A} = \prod \mathcal{A}_i = \langle A, [\ ] \rangle$ является производной от группы $\prod A_i$.

*Доказательство.* Обозначим через $\odot_i$ операцию в группе $A_i$ ($i \in I$), а через $\circ$ – операцию в группе $\prod A_i$. Так как

$$[a_1 a_2 \ldots a_n](i) = [a_1(i) a_2(i) \ldots a_n(i)]_i =$$
$$= a_1(i) \odot_i a_2(i) \odot_i \ldots \odot_i a_n(i) = (a_1 \circ a_2 \circ \ldots \circ a_n)(i),$$

т. е.
$$[a_1 a_2 \ldots a_n](i) = (a_1 \circ a_2 \circ \ldots \circ a_n)(i),$$
то
$$[a_1 a_2 \ldots a_n] = a_1 \circ a_2 \circ \ldots \circ a_n$$

для любых $a_1, a_2, \ldots, a_n \in A$. Следовательно, n-арная группа $\mathcal{A}$ является производной от группы $\prod A_i$. ∎

**2.7.2. Лемма.** Пусть $< A_i, [\ ]_i >$ – n-арные группы ($i \in I$), $a_i$ – фиксированные элементы из $A_i$, $c_1^{(i)} \ldots c_{n-2}^{(i)}$ – обратная последовательность для $a_i$ ($c_1^{(i)}, \ldots, c_{n-2}^{(i)} \in A_i$). Пусть также $a, c_j \in \prod A_i$, где

$$a(i) = a_i \in A_i, \quad c_j(i) = c_j^{(i)} \in A_i \ (j = 1, \ldots, n-2).$$

Тогда $c_1 \ldots c_{n-2}$ – обратная последовательность для элемента $a$ в n-арной группе

$$\prod < A_i, [\ ]_i > = < \prod A_i, [\ ] >.$$

*Доказательство.* Так как

$$[ac_1 \ldots c_{n-2}a](i) = [a(i)c_1(i) \ldots c_{n-2}(i)a(i)]_i =$$
$$= [\underbrace{a_i c_1^{(i)} \ldots c_{n-2}^{(i)}}_{\text{нейтр}} a_i]_i = a_i = a(i),$$

то

$$[ac_1 \ldots c_{n-2}a] = a.$$

Следовательно, $ac_1 \ldots c_{n-2}$ – нейтральная последовательность, а значит, $c_1 \ldots c_{n-2}$ – обратная последовательность для элемента $a$. ∎

**2.7.3. Теорема** [37]. Пусть $< A_i, [\ ]_i >$ – n-арные группы ($i \in I$), $a_i$ – фиксированные элементы из $A_i$, $a \in \prod A_i$, где $a(i) = a_i \in A_i$ Тогда

$$< \prod A_i, \textcircled{a} > = \prod < A_i, \textcircled{a_i} >.$$

*Доказательство.* Если положить $A = \prod A_i$, то доказываемое равенство примет вид

$$< A, \textcircled{a} > = < A, \circ >,$$

при этом операция $\textcircled{a}$ определяется по правилу

$$x \textcircled{a} y = [xc_1 \ldots c_{n-2}y], \ x, y \in A,$$

где $<A, [\ ]> = \prod <A_i, [\ ]_i>$, $c_1 \ldots c_{n-2}$ — обратная последовательность, определённая в лемме 2.7.2, а операция $\circ$ определяется покомпонентно следующим образом

$$(x \circ y)(i) = x(i) \circledast_i y(i).$$

Так как для любых $x, y \in A$ и любого $i \in I$ верно

$$(x \circledast y)(i) = [x\, c_1 \ldots c_{n-2}\, y](i) = [x(i) c_1(i) \ldots c_{n-2}(i) y(i)]_i =$$
$$= [x(i) c_1^{(i)} \ldots c_{n-2}^{(i)} y(i)]_i = x(i) \circledast_i y(i) = (x \circledast_i y)(i),$$

т. е.

$$(x \circledast y)(i) = (x \circledast_i y)(i),$$

то

$$x \circledast y = x \circ y. \qquad \blacksquare$$

**2.7.4. Следствие.** Пусть $<A_1, [\ ]_1>, \ldots, <A_k, [\ ]_k>$ — n-арные группы, $a = (a_1, \ldots, a_k)$, где $a_i$ — фиксированные элементы из $A_i$ ($i = 1, \ldots, k$). Тогда

$$<A_1 \times \ldots \times A_k, \circledast> = <A_1, \circledast_1> \times \ldots \times <A_k, \circledast_k>.$$

**2.7.5. Теорема.** Декартово произведение $\prod(A_i)_\circ$ соответствующих групп $(A_i)_\circ$, $i \in I$, n-арных групп $<A_i, [\ ]_i>$ изоморфно соответствующей группе $(\prod A_i)_\circ$ декартова произведения $<\prod A_i, [\ ]>$ — n-арных групп $<A_i, [\ ]_i>$:

$$\prod(A_i)_\circ \simeq (\prod A_i)_\circ.$$

*Доказательство.* По определению

$$\prod(A_i)_\circ = \{u : I \to U(A_i)_\circ \mid u(i) = \theta_{A_i}(a_{i1} \ldots a_{i(n-1)}) \in (A_i)_\circ\},$$

$$(\prod A_i)_\circ = \{\theta_{\prod A_i}(v_1 \ldots v_{n-1}) \mid v_j \in \prod A_i, v_j : I \to UA_i,$$

$$v_j(i) = b_{ij} \in A_i,\ j = 1, \ldots, n-1\}.$$

Определим отображение $f : \prod(A_i)_\circ \to (\prod A_i)_\circ$ по правилу

$$f : u \to \theta_{\Pi A_i}(v_1 \ldots v_{n-1}),$$

где

$$u(i) = \theta_{A_i}(a_{i1} \ldots a_{i(n-1)}), \ v_j(i) = a_{ij}.$$

Ясно, что f – сюръекция.

Пусть теперь $u'$, $u'' \in \Pi(A_i)_o$, причем

$$u'(i) = \theta_{A_i}(a'_{i1} \ldots a'_{i(n-1)}), \ u''(i) = \theta_{A_i}(a''_{i1} \ldots a''_{i(n-1)})$$

и предположим, что $f(u') = f(u'')$, где

$$f(u') = \theta_{\Pi A_i}(v'_1 \ldots v'_{n-1}), \ f(u'') = \theta_{\Pi A_i}(v''_1 \ldots v''_{n-1}),$$

$$v'_j(i) = a'_{ij}, \ v''_j(i) = a''_{ij}.$$

Если v – произвольный элемент из $\Pi A_i$, причем $v(i) = a_i \in A_i$, то из $f(u') = f(u'')$ следует

$$\theta_{\Pi A_i}(v'_1 \ldots v'_{n-1})\theta_{\Pi A_i}(v) = \theta_{\Pi A_i}(v''_1 \ldots v''_{n-1})\theta_{\Pi A_i}(v),$$

$$\theta_{\Pi A_i}([v'_1 \ldots v'_{n-1}v]) = \theta_{\Pi A_i}([v''_1 \ldots v''_{n-1}v]),$$

$$[v'_1 \ldots v'_{n-1}v] = [v''_1 \ldots v''_{n-1}v],$$

$$[v'_1 \ldots v'_{n-1}v](i) = [v''_1 \ldots v''_{n-1}v](i),$$

$$[v'_1(i) \ldots v'_{n-1}(i)v(i)]_i = [v''_1(i) \ldots v''_{n-1}(i)v(i)]_i,$$

$$[a'_{i1} \ldots a'_{i(n-1)}a_i]_i = [a''_{i1} \ldots a''_{i(n-1)}a_i]_i,$$

$$\theta_{A_i}(a'_{i1} \ldots a'_{i(n-1)}) = \theta_{A_i}(a''_{i1} \ldots a''_{i(n-1)}),$$

$$u'(i) = u''(i), \ i \in I,$$

т. е. $u' = u''$. Следовательно, f – инъекция, а значит и биекция.

Пусть по-прежнему $u'$ и $u''$ произвольные элементы из $\Pi(A_i)_o$, $u'(i)$ и $u''(i)$ определены так же как и выше.

Так как

$$u'\,u''(i) = u'(i)\,u''(i) = \theta_{A_i}(a'_{i1}\ \ldots\ a'_{i(n-1)})\,\theta_{A_i}(a''_{i1}\ \ldots\ a''_{i(n-1)}) =$$

$$= \theta_{A_i}(a'_{i1}\ \ldots\ a'_{i(n-2)}[\,a'_{i(n-1)}\,a''_{i1}\ \ldots\ a''_{i\,(n-1)}]_i) =$$

$$= \theta_{A_i}(v'_1(i)\ \ldots\ v'_{n-2}(i)[\,v'_{n-1}(i)\,v''_1(i)\ \ldots\ v''_{n-1}(i)]_i) =$$

$$= \theta_{A_i}(v'_1(i)\ \ldots\ v'_{n-2}(i)[\,v'_{n-1}\,v''_1\ \ldots\ v''_{n-1}](i)),$$

т.е.

$$u'\,u''(i) = \theta_{A_i}(v'_1(i)\ \ldots\ v'_{n-2}(i)\ [\,v'_{n-1}\,v''_1\ \ldots\ v''_{n-1}](i)),$$

то

$$f(uv) = \theta_{\Pi A_i}(v'_1\ \ldots\ v'_{n-2}[\,v'_{n-1}\,v''_1\ \ldots\ v''_{n-1}]) =$$

$$= \theta_{\Pi A_i}(v'_1\ \ldots\ v'_{n-1})\,\theta_{\Pi A_i}(v''_1\ \ldots\ v''_{n-1}) = f(u')f(u'').$$

Следовательно, $f$ – изоморфизм. ∎

**2.7.6. Следствие.** Пусть $< A_1,\ [\ ]_1 >,\ \ldots,\ < A_k,\ [\ ]_k >$ – n-арные группы. Тогда

$$(A_1)_o \times \ldots \times (A_k)_o \simeq (< A_1,\ [\ ]_1 > \times \ldots \times < A_k,\ [\ ]_k >)_o.$$

**2.7.7. Замечание.** Так как по предложению 1.6.1

$$(A_i)_o \simeq < A_i,\ \textcircled{a}_i >,\ < \Pi A_i,\ \textcircled{a} > \simeq (\Pi A_i)_o,$$

то теоремы 2.7.3 и 2.7.5 являются следствиями друг друга.

Пусть по-прежнему $\{\mathscr{A}_i = < A_i,\ [\ ]_i > |\ i \in I\}$ – семейство n-арных групп, $\mathscr{A} = < A = \Pi A_i,\ [\ ] >$ – их декартово произведение. Обозначим через $\mathscr{A}^* = (\Pi A_i)^*$ обёртывающую группу Поста декартового произведения $\mathscr{A}$, а через $\Pi(A_i)^*$ декартово произведение обёртывающих групп Поста $(A_i)^*$ n-арных групп $\mathscr{A}_i$.

Группу $(\Pi A_i)^*$ можно представить (см., например, предложение 1.3.7.) в виде

$$(\Pi A_i)^* = \bigcup_{j=1}^{n-1} \{\theta_{\Pi A_i}(v\underbrace{s\ldots s}_{j-1}) \mid v \in \Pi A_i\},$$

где s – фиксированный из $\Pi A_i$, при этом

$$v: I \to UA_i,\ v(i) \in A_i;\ s \in I \to UA_i,\ s(i) \in A_i.$$

Аналогично группу $\Pi(A_i)^*$ можно представить в виде

$$\Pi(A_i)^* = \{u: I \to U(A_i)^* \mid u(i) \in (A_i)^*\}.$$

Заметим, что в этом представлении группы $\Pi(A_i)^*$ в записи

$$u'(i) = \theta_{A_i}(a_i' \underbrace{b_i \ldots b_i}_{j-1}),\ u''(i) = \theta_{A_i}(a_i'' \underbrace{b_i \ldots b_i}_{k})$$

различных элементов $u', u'' \in \Pi(A_i)^*$ элемент $b_i$ – фиксированный.

**2.7.8. Лемма.** Тогда и только тогда

$$\theta_{\Pi A_i}(v_1' \ldots v_j') = \theta_{\Pi A_i}(v_1'' \ldots v_k''),\ j \geq 1,\ k \geq 1, \quad (1)$$

когда

$$\theta_{A_i}(v_1'(i) \ldots v_j'(i)) = \theta_{A_i}(v_1''(i) \ldots v_k''(i)) \quad (2)$$

для любого $i \in I$.

*Доказательство.* Из (1) следует, что существуют

$$v_{j+1}, \ldots, v_{m(n-1)+1} \in \Pi A_i$$

такие, что

$$[v_1' \ldots v_j'\ v_{j+1} \ldots v_{m(n-1)+1}] = [v_1'' \ldots v_k''\ v_{j+1} \ldots v_{m(n-1)+1}], \quad (3)$$

откуда

$$[v_1' \ldots v_j'\ v_{j+1} \ldots v_{m(n-1)+1}](i) = [v_1'' \ldots v_k''\ v_{j+1} \ldots v_{m(n-1)+1}](i), \quad (4)$$

$$[v_1'(i) \ldots v_j'(i) v_{j+1}(i) \ldots v_{m(n-1)+1}(i)]_i =$$

$$= [\,v''_1(i) \ldots v''_k(i)v_{j+1}(i) \ldots v_{m(n-1)+1}(i)]_i, \qquad (5)$$

$$\theta_{A_i}(v'_1(i) \ldots v'_j(i)) = \theta_{A_i}(v''_1(i) \ldots v''_k(i)), \, i \in I.$$

Следовательно, для любого $i \in I$ верно (2).

Если теперь для любого $i \in I$ верно (2), то существуют $a_{i(j+1)}, \ldots, a_{i(m(n-1)+1)}$ такие, что

$$[\,v'_1(i) \ldots v'_j(i)a_{i(j+1)} \ldots a_{i(m(n-1)+1)}]_i =$$

$$= [\,v''_1(i) \ldots v''_k(i)\, a_{i(j+1)} \ldots a_{i(m(n-1)+1)}]_i.$$

Полагая, $a_{i(j+1)} = v_{j+1}(i), \ldots, a_{i(m(n-1)+1)} = v_{m(n-1)+1}(i)$, определим элементы $v_{j+1}, \ldots, v_{m(n-1)+1} \in \prod A_i$, а последнее равенство примет вид (5), откуда последовательно получаются (4), (3) и (1). ∎

В связи с теоремой 2.7.5 возникает вопрос: будут ли изоморфными группы $\prod(A_i)^*$ и $(\prod A_i)^*$?

Если $<A_1, [\,]_1>$ и $<A_2, [\,]_2>$ — конечные n-арные группы, то

$$|(A_1 \times A_2)^*| = |A_1 \times A_2|(n-1) = |A_1||A_2|(n-1),$$

$$|A_1^* \times A_2^*| = |A_1^*||A_2^*| = |A_1|(n-1)|A_2|(n-1) = |A_1||A_2|(n-1)^2.$$

Поэтому группы $(A_1 \times A_2)^*$ и $A_1^* \times A_2^*$ не могут быть изоморфными. Таким образом, получен отрицательный ответ на поставленный выше вопрос. Однако, имеет место

**2.7.9. Теорема.** Группа $(\prod A_i)^*$ изоморфно вкладывается в группу $\prod(A_i)^*$: $(\prod A_i)^* \simeq H \subseteq \prod(A_i)^*$.

*Доказательство.* Определим отображение

$$\varphi: (\prod A_i)^* \to \prod(A_i)^*$$

по правилу

$$\varphi: \theta_{\Pi A_i}(v\underbrace{s\ldots s}_{j-1}) \to u, \ u(i) = \theta_{A_i}(v(i)\underbrace{s(i)\ldots s(i)}_{j-1}), \ i \in I.$$

Пусть

$$t' = \theta_{\Pi A_i}(v'\underbrace{s\ldots s}_{j-1}), \ t'' = \theta_{\Pi A_i}(v''\underbrace{s\ldots s}_{k-1}); \ j, k \in \{1, \ldots, n\text{-}1\}$$

– произвольные элементы из $(\Pi A_i)^*$ и предположим, что $\varphi(t')=\varphi(t'')$. Тогда, если $\varphi(t') = u'$, $\varphi(t'') = u''$, то $u'(i) = u''(i)$ для любого $i \in I$, т. е.

$$\theta_{A_i}(v'(i)\underbrace{s(i)\ldots s(i)}_{j-1}) = \theta_{A_i}(v''(i)\underbrace{s(i)\ldots s(i)}_{k-1}).$$

Из последнего равенства следует $j = k$. Поэтому

$$t'' = \theta_{\Pi A_i}(v''\underbrace{s\ldots s}_{j-1})$$

и имеет место равенство

$$\theta_{A_i}(v'(i)\underbrace{s(i)\ldots s(i)}_{j-1}) = \theta_{A_i}(v''(i)\underbrace{s(i)\ldots s(i)}_{j-1}),$$

откуда

$$\theta_{A_i}(v'(i)\underbrace{s(i)\ldots s(i)}_{j-1})\theta_{A_i}(\underbrace{s(i)\ldots s(i)}_{n-j}) =$$

$$= \theta_{A_i}(v''(i)\underbrace{s(i)\ldots s(i)}_{j-1})\theta_{A_i}(\underbrace{s(i)\ldots s(i)}_{n-j}),$$

$$\theta_{A_i}(v'(i)\underbrace{s(i)\ldots s(i)}_{n-1}) = \theta_{A_i}(v''(i)\underbrace{s(i)\ldots s(i)}_{n-1}).$$

$$[v'(i)\underbrace{s(i)\ldots s(i)}_{n-1}]_i = [v''(i)\underbrace{s(i)\ldots s(i)}_{n-1}]_i,$$

$$v'(i) = v''(i), \ i \in I.$$

Следовательно, $v' = v''$ и значит $t' = t''$, а отображение $\varphi$ является инъекцией.

Пусть снова $t'$ и $t''$ – произвольные элементы из $(\Pi A_i)^*$, определенные выше. Так как

$$\varphi(t't'')(i) = \varphi(\theta_{\Pi A_i}(v'\underbrace{s\ldots s}_{j-1})\theta_{\Pi A_i}(v''\underbrace{s\ldots s}_{k-1}))(i) =$$

$$= \varphi(\theta_{\Pi A_i}(v'\underbrace{s\ldots s}_{j-1}v''\underbrace{s\ldots s}_{k-1}))(i) = \varphi(\theta_{\Pi A_i}(v\underbrace{s\ldots s}_{m-1}))(i)$$

для некоторых $v \in \Pi A_i$ и $m \in \{1, \ldots, n-1\}$, то

$$\varphi(t't'')(i) = \theta_{A_i}(v(i)\underbrace{s(i)\ldots s(i)}_{m-1}).$$

Применяя к правой части последнего равенства лемму 2.7.8, получим

$$\varphi(t't'')(i) = \theta_{A_i}(v'(i)\underbrace{s(i)\ldots s(i)}_{j-1}v''(i)\underbrace{s(i)\ldots s(i)}_{k-1}) =$$

$$= \theta_{A_i}(v'(i)\underbrace{s(i)\ldots s(i)}_{j-1})\theta_{A_i}(v''(i)\underbrace{s(i)\ldots s(i)}_{k-1}) =$$

$$= \varphi(t')(i)\varphi(t'')(i) = \varphi(t')\varphi(t'')(i),$$

т. е. $\varphi(t't'') = \varphi(t')\varphi(t'')$. Следовательно, $\varphi$-изоморфизм группы $(\Pi A_i)^*$ на подгруппу $H = \{\varphi(t)| t \in (\Pi A_i)^*\}$ группы $\Pi(A_i)^*$. ∎

**2.7.10. Следствие.** Пусть $< A_1, [\ ]_1 >, \ldots, < A_k, [\ ]_k >$ – n-арные группы и $(b_1, \ldots, b_k)$ – фиксированный элемент из $A_1 \times \ldots \times A_k$. Тогда отображение

$$\psi: \theta_{A_1 \times \ldots \times A_k}((a_1, \ldots, a_k)\underbrace{(b_1,\ldots,b_k)\ldots(b_1,\ldots,b_k)}_{j-1}) \to$$

$$\to (\theta_{A_1}(a_1\underbrace{b_1\ldots b_1}_{j-1}), \ldots, \theta_{A_k}(a_k\underbrace{b_k\ldots b_k}_{j-1}))$$

является изоморфизмом группы

$$(A_1 \times \ldots \times A_k)^* = \bigcup_{j=1}^{n-1} \{\theta_{A_1 \times \ldots \times A_k}((a_1, \ldots, a_k)\underbrace{(b_1, \ldots, b_k)\ldots(b_1, \ldots, b_k)}_{j-1}) \mid (a_1, \ldots, a_k) \in A_1 \times \ldots \times A_k\},$$

на подгруппу

$$H = \bigcup_{j=1}^{n-1} \{(\theta_{A_1}(a_1\underbrace{b_1\ldots b_1}_{j-1}), \ldots, \theta_{A_k}(a_k\underbrace{b_k\ldots b_k}_{j-1})) \mid (a_1, \ldots, a_k) \in A_1 \times \ldots \times A_k\}$$

группы

$$A_1^* \times \ldots \times A_k^* = \bigcup_{j_1=1,\ldots,j_k=1}^{n-1} \{(\theta_{A_1}(a_1\underbrace{b_1\ldots b_1}_{j_1-1}), \ldots, \theta_{A_k}(a_k\underbrace{b_k\ldots b_k}_{j_k-1})) \mid$$

$$(a_1, \ldots, a_k) \in A_1 \times \ldots \times A_k.$$

**2.7.11. Замечание.** Образ $H = \varphi((\prod A_i)^*)$ группы $(\prod A_i)^*$ из предыдущей теоремы можно представить следующим образом

$$H = \{u \in \prod(A_i)^* \mid \exists j \in \{1, \ldots, n-1\}, \forall i \in I, u(i) \in A^{*(j)}\},$$

где

$$A^{*(j)} = \{\theta_A(a\underbrace{b\ldots b}_{j-1}) \mid a \in A\},$$

b – фиксированный из A.

Декартово произведение n-арных групп называют еще *внешним прямым произведением*. В теории n-арных групп ($n \geq 3$), так же как и в теории групп, рассматривают и внутренние прямые произведения n-арных подгрупп. Мы ограничимся здесь случаем конечного числа сомножителей.

**2.7.12. Лемма.** Пусть

$$\mathscr{B}_1 = <B_1, [\ ]>, \mathscr{B}_2 = <B_2, [\ ]>, \ldots, \mathscr{B}_k = <B_k, [\ ]>$$

– полуинвариантные n-арные подгруппы n-арной группы $<A, [\ ]>$, пересечение которых непусто. Тогда

$$<[B_1 \overset{n-1}{B_2} \ldots \overset{n-1}{B_k}], [\ ]>$$

– полуинвариантная n-арная подгруппа в $<A, [\ ]>$, порождаемая n-арными подгруппами $\mathscr{B}_1, \mathscr{B}_2, \ldots, \mathscr{B}_k$, т.е.

$$<B_1, B_2, \ldots, B_k> = [B_1 \overset{n-1}{B_2} \ldots \overset{n-1}{B_k}].$$

***Доказательство.*** Из следствия 2.3.18 следует, что

$$<[B_1 \overset{n-1}{B_2} \ldots \overset{n-1}{B_k}], [\ ]>$$

– полуинвариантная n-арная подгруппа в $<A, [\ ]>$, а по лемме 2.3.20

$$[B_1 \overset{n-1}{B_2} \ldots \overset{n-1}{B_k}] = [\overset{n-1}{B_1} \ldots \overset{n-1}{B_{i-1}} \overset{n-1}{B_i} \overset{n-1}{B_{i+1}} \ldots \overset{n-1}{B_k}], i = 2, \ldots, k.$$

Поэтому

$$B_i \subseteq [B_1 \overset{n-1}{B_2} \ldots \overset{n-1}{B_k}], \quad i = 1, 2, \ldots, k,$$

и значит

$$<B_1, B_2, \ldots, B_k> \subseteq [B_1 \overset{n-1}{B_2} \ldots \overset{n-1}{B_k}].$$

С другой стороны, по теореме 2.1.14

$$[B_1 \overset{n-1}{B_2} \ldots \overset{n-1}{B_k}] \subseteq <B_1, B_2, \ldots, B_k>.$$

Из доказанных включений следует требуемое равенство. ∎

Пусть $\mathscr{A} = <A, [\ ]> = \mathscr{B}_1 \times \ldots \times \mathscr{B}_m$ – внешнее прямое произведение n-арных групп $\mathscr{B}_i = <B_i, [\ ]_i>$, каждая из которых содержит идемпотент $a_i$ ($i = 1, \ldots, a_m$). Положим $a = (a_1, \ldots, a_m)$,

$$A_i = \{(a_1, \ldots, a_{i-1}, b_i, a_{i+1}, \ldots, a_m) | b_i \in B_i\} =$$

$$= (a_1, \ldots, a_{i-1}, B_i, a_{i+1}, \ldots, a_m).$$

В следующей лемме доказано, что все $<A_i, [\ ]>$ – n-арные подгруппы в $<A, [\ ]>$. Некоторые важные свойства этих n-арных подгрупп описаны в теореме 2.7.14.

**2.7.13. Лемма.** Справедливы следующие утверждения:
1) $a$ – идемпотент n-арной группы $\mathscr{A}$, причем $\bigcap A_i = \{a\}$;
2) $\mathscr{A}_i = <A_i, [\ ]>$ – n-арная подгруппа в $\mathscr{A}$;
3) $A = [A_1 \underbrace{a \ldots a}_{n-2} A_2 \underbrace{a \ldots a}_{n-2} \ldots A_{m-1} \underbrace{a \ldots a}_{n-2} A_m]$;
4) отображение $\varepsilon_i : b_j \to (a_1, \ldots, a_{i-1}, b_j, a_{i+1}, \ldots, a_m)$ является изоморфизмом $\mathscr{B}_i$ на $<A_i, [\ ]>$, $i = 1, \ldots, m$.

***Доказательство.*** 1) Из определения n-арной операции $[\ ]$ в $\mathscr{A}$ и того, что $a_i$ – идемпотент в $\mathscr{B}_i$ следует, что $a$ – идемпотент $\mathscr{A}$.

Если $b = \{b_1, \ldots, b_m\} \in \bigcap A_i$, то $b \in A_1$, откуда $b = \{b_1, a_2, \ldots, a_m\}$. А так как $b \in A_2$, то $b_1 = a_1$. Следовательно, $b = a$.

2) Следует из определения n-арной операции $[\ ]$ в $\mathscr{A}$ и того, что все $\mathscr{B}_i$ – n-арные группы.

3) Ясно, что

$$[A_1 \underbrace{a \ldots a}_{n-2} A_2 \underbrace{a \ldots a}_{n-2} \ldots A_{m-1} \underbrace{a \ldots a}_{n-2} A_m] \subseteq A.$$

Если теперь $b = (b_1, b_2, \ldots, b_m)$ – произвольный элемент из $A$, то

$$b = (b_1, b_2, \ldots, b_m) = ([b_1 \underbrace{a_1 \ldots a_1}_{(m-1)(n-1)}]_1, [\underbrace{a_2 \ldots a_2}_{n-1} b_2 \underbrace{a_2 \ldots a_2}_{(m-2)(n-1)}]_2,$$

$$\ldots, [\underbrace{a_i \ldots a_i}_{(i-1)(n-1)} b_i \underbrace{a_i \ldots a_i}_{(m-i)(n-1)}]_i, \ldots, [\underbrace{a_m \ldots a_m}_{(m-1)(n-1)} b_m]_m) =$$

$$= [(b_1, a_2, \ldots, a_m) \underbrace{(a_1, a_2, \ldots, a_m) \ldots (a_1, a_2, \ldots, a_m)}_{n-2} (a_1, b_2, a_3, \ldots, a_m)$$

$$\underbrace{(a_1, a_2, \ldots a_m) \ldots (a_1, a_2, \ldots, a_m)}_{n-2} \ldots (a_1, \ldots, a_{m-2}, b_{m-1}, a_m)$$

$$\underbrace{(a_1, a_2, \ldots a_m) \ldots (a_1, a_2, \ldots, a_m)}_{n-2} (a_1, a_2, \ldots, a_{m-1}, b_m)] \subseteq$$

$$\subseteq [A_1 \underbrace{a \ldots a}_{n-2} A_2 \underbrace{a \ldots a}_{n-2} \ldots A_{m-1} \underbrace{a \ldots a}_{n-2} A_m],$$

т. е.

$$A \subseteq [A_1 \underbrace{a \ldots a}_{n-2} A_2 \underbrace{a \ldots a}_{n-2} \ldots A_{m-1} \underbrace{a \ldots a}_{n-2} A_m].$$

Из доказанных включений следует равенство 3).

4) Следует из определения n-арной операции [ ] в $\mathscr{A}$ и того, что $a_i$ — идемпотент в $\mathscr{B}_i$. ∎

Далее для сокращения записей иногда будем употреблять символ $\overset{k}{A}$ вместо $\underbrace{A \ldots A}_{k}$ (см. с. 69).

**2.7.14. Теорема.** Справедливы следующие утверждения:

1) $\mathscr{A}_i = <A_i, [\ ]>$ – полуинвариантная n-арная подгруппа n-арной группы $\mathscr{A} = <A, [\ ]>$, $i = 1, \ldots, m$;

2) $A = [A_1 \overset{n-1}{A_2} \ldots \overset{n-1}{A_m}]$;

3) $<A_1, \ldots, A_{i-1}, A_{i+1}, \ldots, A_m> \cap A_i = \{a\}$, $i = 1, \ldots, m$.

***Доказательство.*** 1) Если $x = (x_1, \ldots, x_m)$ – произвольный элемент из A, то, учитывая $x_i \in B_i$, получаем

$$[x \overset{n-1}{A_i}] = [(x_1, \ldots, x_{i-1}, x_i, x_{i+1}, \ldots, x_m)$$

$$\underbrace{(a_1, \ldots, a_{i-1}, B_i, a_{i+1}, \ldots, a_m) \ldots (a_1, \ldots, a_{i-1}, B_i, a_{i+1}, \ldots, a_m)}_{n-1}] =$$

$$= ([x_1 \underbrace{a_1 \ldots a_1}_{n-1}]_1, \ldots, [x_{i-1} \underbrace{a_{i-1} \ldots a_{i-1}}_{n-1}]_{i-1}, [x_i \overset{n-1}{B_i}]_i,$$

$$[x_{i+1} \underbrace{a_{i+1} \ldots a_{i+1}}_{n-1}]_{i+1}, \ldots, [x_m \underbrace{a_m \ldots a_m}_{n-1}]_m) = (x_1, \ldots, x_{i-1}, B_i, x_{i+1}, \ldots, x_m),$$

т. е.
$$[x \overset{n-1}{A_i}] = (x_1, \ldots, x_{i-1}, B_i, x_{i+1}, \ldots, x_m).$$

Аналогично доказывается равенство

$$[\overset{n-1}{A_i} x] = (x_1, \ldots, x_{i-1}, B_i, x_{i+1}, \ldots, x_m).$$

Следовательно,
$$[x \overset{n-1}{A_i}] = [\overset{n-1}{A_i} x],$$

и n-арная группа $\mathscr{A}_i$ – полуинвариантна в $\mathscr{A}$.

2) Следует из утверждения 3) леммы 2.7.13 и того, что $a \in A_i$ для любого $i = 1, \ldots, m$.

3) Докажем вначале справедливость следующих равенств

$$[A_2 \overset{n-1}{A_3} \ldots \overset{n-1}{A_m}] \cap A_1 = \{a\}, \qquad (1)$$

$$[\overset{n-1}{A_1} \overset{n-1}{A_2} \ldots \overset{n-1}{A_{i-1}} \overset{n-1}{A_{i+1}} \ldots \overset{n-1}{A_m}] \cap A_i = \{a\}, \; i = 2, \ldots, m. \qquad (2)$$

Если
$$b \in [A_2 \overset{n-1}{A_3} \ldots \overset{n-1}{A_m}] \cap A_1,$$

то $b \in A_1$, откуда $b = (b_1, a_2, \ldots, a_m)$ для некоторого $b_1 \in B_1$. А так как

$$b \in [A_2 \overset{n-1}{A_3} \ldots \overset{n-1}{A_m}],$$

то
$$b_1 = [\underbrace{a_1 \ldots a_1}_{(m-2)(m-1)+1}]_1 = a_1,$$

Таким образом, $b = (a_1, \ldots, a_m) = a$ и верно (1).

Равенство (2) доказывается аналогично.

Из (1), (2) и леммы 2.7.12 следует

$$< A_2, A_3, \ldots, A_m > \cap A_1 = \{a\},$$

$$< A_1, A_2, \ldots, A_{i-1}, A_{i+1}, \ldots, A_m > \cap A_i = \{a\}, \ i = 2, \ldots, m.$$

Следовательно, утверждение 3) верно. ∎

**2.7.15. Предложение** [4]. Пусть

$$\mathscr{A}_1 = < A_1, [\ ] >, \ldots, \mathscr{A}_m = < A_m, [\ ] >$$

– n-арные подгруппы n-арной группы $\mathscr{A} = < A, [\ ] >$, содержащей идемпотент a, и пусть выполняются утверждения 1) - 3) теоремы 2.7.14. Тогда справедливы следующие утверждения:

1) $\bigcap A_i = \{a\}$;

2) $[\overset{n-1}{A_1} \overset{n-1}{A_2} \ldots \overset{n-1}{A_{r-1}}] \cap A_r = \{a\}, r = 2, \ldots, m$;

3) $A = [\overset{n-1}{A_1} \ldots \overset{n-1}{A_{i-1}} \overset{n-1}{A_i} \overset{n-1}{A_{i+1}} \ldots A_m], i = 1, 2, \ldots, m$;

4) для любых $r = 2, \ldots, m$ и $j = 1, 2, \ldots, r-1$ верно

$$[\overset{n-1}{A_1} \ldots \overset{n-1}{A_{j-1}} \overset{n-1}{A_j} \overset{n-1}{A_{j+1}} \ldots \overset{n-1}{A_{r-1}}] \cap A_r = \{a\};$$

5) для любых $x \in A_k, y \in A_s$, где $k, s \in \{1, 2, \ldots, m\}, k \neq s$ верно

$$[x \underbrace{a \ldots a}_{n-2} y] = [y \underbrace{a \ldots a}_{n-2} x];$$

6) любой элемент $b \in A$ единственным образом представим в виде

$$[b_1 \underbrace{a \ldots a}_{n-2} b_2 \underbrace{a \ldots a}_{n-2} \ldots b_{m-1} \underbrace{a \ldots a}_{n-2} b_m], b_i \in A_i (i = 1, \ldots, m).$$

*Доказательство.* 1) Следует из 3) теоремы 2.7.14.

2) Предположим, что утверждение 2) не выполняется, т. е. существуют $r \in \{2, \ldots, m\}$ и $d \neq a$ такие, что

$$d \in [\overset{n-1}{A_1} \overset{n-1}{A_2} \ldots A_{r-1}] \cap A_r,$$

откуда

$$d \in [\overset{n-1}{A_1} \overset{n-1}{A_2} \ldots \overset{n-1}{A_{r-1}} \overset{n-1}{A_{r+1}} \ldots A_m] \cap A_r.$$

Тогда согласно 3) теоремы 2.7.14, $d = a$, что противоречит предположению $d \neq a$.

3) По лемме 2.3.20

$$[\overset{n-1}{A_1} A_i] = [A_1 \overset{n-1}{A_i}].$$

Поэтому из 2) теоремы 2.7.14, используя полуинвариантность n-арных подгрупп $\mathcal{A}_1, \ldots, \mathcal{A}_{i-1}$, получаем

$$A = [\overset{n-1}{A_1} \overset{n-1}{A_2} \ldots \overset{n-1}{A_{i-1}} \overset{n-1}{A_i} \overset{n-1}{A_{i+1}} \ldots A_m] =$$

$$= [\overset{n-1}{A_2} \ldots \overset{n-1}{A_{i-1}} [\overset{n-1}{A_1} A_i] \overset{n-1}{A_{i+1}} \ldots A_m] =$$

$$= [\overset{n-1}{A_2} \ldots \overset{n-1}{A_{i-1}} [A_1 \overset{n-1}{A_i}] \overset{n-1}{A_{i+1}} \ldots A_m] =$$

$$= [\overset{n-1}{A_1} \ldots \overset{n-1}{A_{i-1}} A_i \overset{n-1}{A_{i+1}} \ldots A_m].$$

4) Вытекает из 2) и доказывается аналогично предыдущему пункту 3).

5) Ясно, что в A существует c такой, что

$$[x\underbrace{a\ldots a}_{n-2}y] = [[y\underbrace{a\ldots a}_{n-2}x]\underbrace{a\ldots a}_{n-2}c], \qquad (1)$$

откуда

$$c = [a\bar{x}\underbrace{x\ldots x}_{n-3}a\bar{y}\underbrace{y\ldots y}_{n-3}x\underbrace{a\ldots a}_{n-2}y].$$

Так как $A_k$ и $A_s$ – полуинвариантны в $\mathcal{A}$, то существуют

$$u_1, \ldots, u_{n-1} \in A_k, \qquad v_1, \ldots, v_{n-1} \in A_s$$

такие, что

$$[x\underbrace{a\ldots a}_{n-2}y] = [yu_1\ldots u_{n-1}], \quad [a\bar{y}\underbrace{y\ldots y}_{n-3}x] = [xv_1\ldots v_{n-1}].$$

Поэтому

$$c = [a\bar{x}\underbrace{x\ldots x}_{n-3}a\bar{y}\underbrace{y\ldots y}_{n-3}yu_1\ldots u_{n-1}],$$

$$c = [a\bar{x}\underbrace{x\ldots x}_{n-3}xv_1\ldots v_{n-1}\underbrace{a\ldots a}_{n-2}y],$$

откуда

$$c = [a\bar{x}\underbrace{x\ldots x}_{n-3}au_1\ldots u_{n-1}] \in A_k,$$

$$c = [a v_1 \ldots v_{n-1}\underbrace{a\ldots a}_{n-2}y] \in A_s.$$

Следовательно, $c \in A_k \cap A_s$.

Пусть для определенности $k < s$. Тогда, согласно утверждению 4) этого предложения,

$$T = [\overset{n-1}{A_1}\ldots \overset{n-1}{A_{k-1}}\overset{n-1}{A_k}\overset{n-1}{A_{k+1}}\ldots \overset{n-1}{A_{s-1}}] \cap A_s = \{a\}.$$

Так $c \in A_k$, то $c \in T$, а так как $c \in A_s$, то $c \in T \cap A_s$. Поэтому из последнего равенства следует $c = a$, а из равенства (1) следует доказываемое равенство 5).

6) Из 2) теоремы 2.7.14 вытекает, что любой элемент $b \in A$ может быть представлен в виде

$$b = [b_1 b_{21}\ldots b_{2(n-1)}\ldots b_{m1}\ldots b_{m(n-1)}],$$

где $b_1 \in A_1$, $b_{ij} \in A_i$ ($i = 2, \ldots, m$). Так как $a \in A_i$, то последовательность $b_{i1}\ldots b_{i(n-1)}$ эквивалентна последовательности $\underbrace{a\ldots a}_{n-2}b_i$, где $b_i$ – некоторый элемент из $A_i$, который определяется единственным образом. Следовательно,

$$b = [b_1\underbrace{a\ldots a}_{n-2}b_2\underbrace{a\ldots a}_{n-2}\ldots b_{m-1}\underbrace{a\ldots a}_{n-2}b_m].$$

Пусть теперь

$$b = [c_1 \underbrace{a \ldots a}_{n-2} c_2 \underbrace{a \ldots a}_{n-2} \ldots c_{m-1} \underbrace{a \ldots a}_{n-2} c_m]$$

еще одно представление элемента b, где $c_i \in A_i$ (i = 1, ..., m), т. е.

$$[b_1 \overset{n-2}{a} \ldots b_{i-1} \overset{n-2}{a} b_i \overset{n-2}{a} b_{i+1} \overset{n-2}{a} \ldots \overset{n-2}{a} b_m] =$$

$$= [c_1 \overset{n-2}{a} \ldots c_{i-1} \overset{n-2}{a} c_i \overset{n-2}{a} c_{i+1} \overset{n-2}{a} \ldots \overset{n-2}{a} c_m],$$

откуда

$$c_i = [a \overline{c_{i-1}} \overset{n-3}{c_{i-1}} \ldots a \overline{c_1} \overset{n-3}{c_1}$$

$$[b_1 \overset{n-2}{a} \ldots b_{i-1} \overset{n-2}{a} b_i \overset{n-2}{a} b_{i+1} \overset{n-2}{a} \ldots \overset{n-2}{a} b_m]$$

$$\overline{c_m} \overset{n-3}{c_m} a \ldots \overline{c_{i+1}} \overset{n-3}{c_{i+1}} a].$$

Из последнего равенства, используя полуинвариантность n-арных групп $\mathcal{A}_i$, для i = 2, ..., m, получаем

$$c_i \in [\overset{n-1}{A_{i-1}} \ldots \overset{n-1}{A_1} \overset{n-1}{A_1} \ldots \overset{n-1}{A_{i-1}} b_i \overset{n-1}{A_{i+1}} \ldots \overset{n-1}{A_m} \overset{n-1}{A_m} \ldots \overset{n-1}{A_{i+1}}] =$$

$$= [b_i \overset{n-1}{A_1} \overset{n-1}{A_2} \ldots \overset{n-1}{A_{i-1}} \overset{n-1}{A_{i+1}} \ldots \overset{n-1}{A_m}] =$$

$$= [b_i \underbrace{a \ldots a}_{n-2} [\overset{n-1}{A_1} \overset{n-1}{A_2} \ldots \overset{n-1}{A_{i-1}} \overset{n-1}{A_{i+1}} \ldots \overset{n-1}{A_m}]] = [b_i \underbrace{a \ldots a}_{n-2} F_i],$$

т. е.

$$c_i = [b_i \underbrace{a \ldots a}_{n-2} d_i], d_i \in F_i, \ F_i = [\overset{n-1}{A_1} \overset{n-1}{A_2} \ldots \overset{n-1}{A_{i-1}} \overset{n-1}{A_{i+1}} \ldots \overset{n-1}{A_m}].$$

Если же i = 1, то

$$c_1 \in [b_1 \overset{n-1}{A_2} \ldots \overset{n-1}{A_m}], \ F_1 = [\overset{n-1}{A_2} \overset{n-1}{A_3} \ldots \overset{n-1}{A_m}].$$

А так как $c_i, b_i, a \in A_i$, то $d_i \in A_i$. Следовательно, $d_i \in F_i \cap A_i$. По лемме 2.7.12, $F_i = \langle A_1, \ldots, A_{i-1}, A_{i+1}, \ldots, A_m \rangle$. Тогда со-гласно 3) теоремы 2.7.14, $F_i \cap A_i = \{a\}$, поэтому $d_i = a$ и $c_i = [b_i \underbrace{a \ldots a}_{n-2} a] = b_i$, т.е. $c_i = b_i$ ($i = 1, \ldots, m$). ∎

Следующая лемма вытекает из утверждения 5) предложения 2.7.15 и может быть доказана индукцией по $m$.

**2.7.16. Лемма** [4]. Пусть

$$\mathscr{A}_1 = <A_1, [\ ]>, \ldots, \mathscr{A}_m = <A_m, [\ ]>$$

– n-арные подгруппы n-арной группы $<A, [\ ]>$, содержащие общий идемпотент $a$, и удовлетворяющие 1) – 3) теоремы 2.7.14 и пусть $b_{ij}$ – произвольные элементы из $A_j$ ($i = 1, \ldots, n$; $j = 1, \ldots, m$). Тогда

$$[[b_{11} \ldots b_{n1}] \underbrace{a \ldots a}_{n-2} [b_{12} \ldots b_{n2}] \underbrace{a \ldots a}_{n-2} \ldots [b_{1m} \ldots b_{nm}]] =$$

$$= [[b_{11} \underbrace{a \ldots a}_{n-2} b_{12} \underbrace{a \ldots a}_{n-2} \ldots b_{1m}] \ldots [b_{n1} \underbrace{a \ldots a}_{n-2} b_{n2} \underbrace{a \ldots a}_{n-2} \ldots b_{nm}]].$$

**2.7.17. Теорема** [4]. Пусть

$$\mathscr{A}_1 = <A_1, [\ ]>, \ldots, \mathscr{A}_m = <A_m, [\ ]>$$

– n-арные подгруппы n-арной группы $\mathscr{A} = <A, [\ ]>$, содержащий идемпотент $a$, и удовлетворяющие 1) – 3) теоремы 2.7.14. Тогда отображение

$$\alpha : <A_1 \times \ldots \times A_m, [\ ]> \to <A, [\ ]>$$

по правилу

$$\alpha : (b_1, \ldots, b_m) \to [b_1 \underbrace{a \ldots a}_{n-2} b_2 \underbrace{a \ldots a}_{n-2} \ldots b_{m-1} \underbrace{a \ldots a}_{n-2} b_m]$$

является изоморфизмом.

*Доказательство.* Так как, согласно 1) предложения 2.7.15, $a \in A_i$ ($i_1, \ldots, m$), то

$$[A_1 \overset{n-1}{A_2} \ldots \overset{n-1}{A_{m-1}}] = [A_1 \underbrace{a \ldots a}_{n-2} A_2 \underbrace{a \ldots a}_{n-2} \ldots A_{m-1} \underbrace{a \ldots a}_{n-2} A_m].$$

Поэтому из 2) теоремы 2.7.14 следует, что $\alpha$ – сюръекция.

Предположим, что

$$\alpha((b_1, \ldots, b_m)) = \alpha((b'_1, \ldots, b'_m)),$$

откуда

$$[b_1 \underbrace{a \ldots a}_{n-2} b_2 \underbrace{a \ldots a}_{n-2} \ldots b_{m-1} \underbrace{a \ldots a}_{n-2} b_m] =$$

$$= [b'_1 \underbrace{a \ldots a}_{n-2} b'_2 \underbrace{a \ldots a}_{n-2} \ldots b'_{m-1} \underbrace{a \ldots a}_{n-2} b'_m].$$

Тогда, согласно утверждению 6) предложения 2.7.15, $b_i = b'_i$ ($i = 1, \ldots, m$) и значит

$$(b_1, \ldots, b_m) = (b'_1, \ldots, b'_m).$$

Следовательно, $\alpha$ – инъекция и значит и биекция.

Применяя теперь лемму 2.7.16, получим

$$\alpha([(b_{11},b_{12}, \ldots, b_{1m})(b_{21},b_{22}, \ldots, b_{2m}) \ldots (b_{n1},b_{n2}, \ldots, b_{nm})]) =$$

$$= \alpha(([b_{11}\, b_{21} \ldots b_{n1}], [b_{12}\, b_{22} \ldots b_{n2}], \ldots, [b_{1m}\, b_{2m} \ldots b_{nm}])) =$$

$$= [[b_{11}\, b_{21} \ldots b_{n1}] \underbrace{a \ldots a}_{n-2} [b_{12}\, b_{22} \ldots b_{n2}] \underbrace{a \ldots a}_{n-2} \ldots [b_{1m}\, b_{2m} \ldots b_{nm}]] =$$

$$= [[b_{11}\underbrace{a \ldots a}_{n-2} b_{12}\underbrace{a \ldots a}_{n-2} \ldots b_{1m}] \ldots [b_{n1}\underbrace{a \ldots a}_{n-2} b_{n2}\underbrace{a \ldots a}_{n-2} \ldots b_{nm}]] =$$

$$= [\alpha(b_{11},b_{12}, \ldots, b_{1m})\, \alpha(b_{21},b_{22}, \ldots, b_{2m}) \ldots \alpha(b_{n1},b_{n2}, \ldots, b_{nm})].$$

Следовательно, $\alpha$ – изоморфизм. ∎

Ясно, что теоремы 2.7.14 и 2.7.17 являются взаимообратными.

Теперь мы можем привести определение внутреннего прямого произведения n-арных подгрупп n-арной группы.

**2.7.18. Определение** [4]. n-Арная группа $< A, [\ ] >$, содержащая идемпотент $a$, называется *a-прямым произведением* своих n-арных подгрупп $< B_1, [\ ] >, ..., < B_m, [\ ] >$, если при $m = 1$ $A = B_1$, а при $m \geq 2$:

1) $< B_i, [\ ] >$ – полуинвариантна в $< A, [\ ] >$ ($i = 1, ..., m$);

2) $A = [B_1 \underbrace{B_2 ... B_2}_{n-1} ... \underbrace{B_m ... B_m}_{n-1}]$;

3) $< B_1, ..., B_{i-1}, B_{i+1}, ..., B_m > \bigcap B_i = \{a\}$, $i = 2, ..., m$.

Если n-арная группа $< A, [\ ] >$ является a-прямым произведением своих n-арных подгрупп $< B_1, [\ ] >, ..., < B_m, [\ ] >$, то будем употреблять обозначение

$$< A, [\ ] > = < B_1, [\ ] > \overset{a}{\times} ... \overset{a}{\times} < B_m, [\ ] >.$$

**2.7.19. Лемма.** Если $< B_1, [\ ] >, ..., < B_m, [\ ] >$ – n-арные подгруппы n-арной группы $< A, [\ ] >$, содержащие элемент $a \in A$, то

$$[B_1 \underbrace{B_2 ... B_2}_{n-1} ... \underbrace{B_m ... B_m}_{n-1}] = B_1 \text{\textcircled{a}} B_2 \text{\textcircled{a}} ... \text{\textcircled{a}} B_m,$$

где $< B_1, \text{\textcircled{a}} >, ..., < B_m, \text{\textcircled{a}} >$ – подгруппы группы $< A, \text{\textcircled{a}} >$.

*Доказательство.* Так как $a \in B_i$, то $< B_i, \text{\textcircled{a}} >$ – подгруппа в $< A, \text{\textcircled{a}} >$ для любого $i = 1, ..., m$. Используя нейтральность последовательности $\alpha a$ и учитывая $a \in A$, а также то, что $< B_i, [\ ] >$ – n-арная подгруппа в $< A, [\ ] >$, получим

$$[B_1 \underbrace{B_2 ... B_2}_{n-1} ... \underbrace{B_m ... B_m}_{n-1}] =$$

$$= [B_1 \alpha a \underbrace{B_2 ... B_2}_{n-1} \alpha a \underbrace{B_3 ... B_3}_{n-1} ... \alpha a \underbrace{B_m ... B_m}_{n-1}] =$$

$$= [B_1 \alpha [a \underbrace{B_2 \ldots B_2}_{n-1}] \alpha [a \underbrace{B_3 \ldots B_3}_{n-1}] \ldots \alpha [a \underbrace{B_m \ldots B_m}_{n-1}]] =$$

$$= [B_1 \alpha B_2 \alpha B_3 \ldots \alpha B_m] = B_1 \, @ \, B_2 \, @ \, \ldots \, @ \, B_m,$$

т. е.

$$[B_1 \underbrace{B_2 \ldots B_2}_{n-1} \ldots \underbrace{B_m \ldots B_m}_{n-1}] = B_1 \, @ \, B_2 \, @ \, \ldots \, @ \, B_m. \qquad \blacksquare$$

**2.7.20. Следствие.** Пусть $< B_1, [\,] >, \ldots, < B_m, [\,] >$ — полуинвариантные n-арные подгруппы n-арной группы $< A, [\,] >$, содержащие общий элемент a, $< F, [\,] >$ — n-арная подгруппа, порожденная ими, $< H, @ >$ — подгруппа группы $< A, @ >$, порожденная подгруппами $< B_1, @ >, \ldots, < B_m, @ >$. Тогда $F = H$.

*Доказательство.* По следствию 2.3.13 все $< B_1, @ >, \ldots,$ $< B_m, @ >$ инвариантны в $< A, @ >$. Поэтому

$$H = B_1 \, @ \, B_2 \, @ \, \ldots \, @ \, B_m.$$

С другой стороны, по лемме 2.7.13

$$F = [B_1 \overset{n-1}{B_2} \ldots \overset{n-1}{B_m}].$$

Применяя к доказанным равенствам лемму 2.7.19, получаем $F = H$. $\blacksquare$

Следующая теорема вытекает из следствия 2.3.13, леммы 2.7.19 и следствия 2.7.20.

**2.7.21. Теорема** [37]**.** n-Арная группа $< A, [\,] >$ является a-прямым произведением своих n-арных подгрупп $< B_1, [\,] >, \ldots, < B_m, [\,] >$ тогда и только тогда, когда группа $< A, @ >$ является прямым произведением своих подгрупп $< B_1, @ >, \ldots, < B_m, @ >$, т. е.

$$< A, [\,] > = < B_1, [\,] > \overset{a}{\times} \ldots \overset{a}{\times} < B_m, [\,] >$$

тогда и только тогда, когда

$$< A, @ > = < B_1, @ > \times ... \times < B_m, @ >.$$

Из теоремы 2.7.17 вытекает

**2.7.22 Следствие.** n-Арная группа $< A, [\ ] >$, являющаяся a-прямым произведением своих n-арных подгрупп $< B_1, [\ ] >, ..., < B_m, [\ ] >$, изоморфна их прямому произведению:

$$< A, [\ ] > = < B_1, [\ ] > \overset{a}{\times} ... \overset{a}{\times} < B_m, [\ ] > \simeq < B_1, [\ ] > \times ... \times < B_m, [\ ] >.$$

## § 2.8. ПОЛУАБЕЛЕВЫ n-АРНЫЕ ГРУППЫ С ИДЕМПОТЕНТАМИ

Известно, что при переходе от групп к n-арным группам возможны различные обобщения абелевости, среди которых одним из самых широких является полуабелевость. Поэтому, распространяя результаты об абелевых группах на произвольные n-арные группы, естественно пытаться это делать сразу для полуабелевых n-арных групп. Такой подход реализован в данном параграфе при получении его основного результата – n-арного аналога теоремы о разложении абелевой группы в прямое произведение своих силовских подгрупп. При этом, в качестве n-арного аналога внутреннего прямого произведения подгрупп, используется введенное С.А. Русаковым a-прямое произведение n-арных подгрупп n-арной группы (определение 2.7.18), предполагающее наличие в последней идемпотентного элемента. По этой причине в параграфе рассматриваются только полуабелевы n-арные группы с идемпотентами.

Холловские n-арные подгруппы определяются [4] аналогично холловским подгруппам: если $\pi$ – множество простых чисел, то n-арная подгруппа $< B, [\ ] > –$ n-арной группы $< A, [\ ] >$ называется $\pi$ – холловской, если её порядок |B| яв-

ляется наибольшим π – делителем порядка $|A|$, т.е. $|B| = |A|_\pi$. Если $\pi = \{p\}$, то $\{p\}$ – холловская n-арная подгруппа называется p – силовской.

**2.8.1. Теорема** [37]**.** Конечная полуабелева n-арная группа $< A, [\ ] >$ порядка $|A| = p_1^{\alpha_1} \ldots p_m^{\alpha_m}$ ($p_1, \ldots, p_m$ – простые), содержащая идемпотент a, единственным образом разлагается в a-прямое произведение

$$< A, [\ ] > = < A(p_1), [\ ] > \overset{a}{\times} \ldots \overset{a}{\times} < A(p_m), [\ ] > \qquad (*)$$

своих $p_i$-силовских ($i = 1, \ldots, m$) n-арных подгрупп $< A(p_i), [\ ] >$.

*Доказательство.* Из полуабелевости n-арной группы $< A, [\ ] >$ вытекает абелевость группы $< A, @ >$, которая по соответствующей бинарной теореме разлагается в прямое произведение

$$< A, @ > = < A(p_i), @ > \times \ldots \times < A(p_m), @ > \qquad (**)$$

своих $p_i$-силовских ($i = 1, \ldots, m$) n-арных подгрупп $< A(p_i), @ >$.

Легко проверяется, что преобразование

$$x \to [ax\underbrace{a \ldots a}_{n-2}]$$

является автоморфизмом группы $< A, @ >$. Отсюда, с учетом характеристичности силовских подгрупп в абелевой группе, получаем

$$[aA(p_i)\underbrace{a \ldots a}_{n-2}] = A(p_i)$$

для любого $i = 1, \ldots, m$. Применяя теперь следствие 2.2.10, заключаем, что $< A(p_i), [\ ] >$ – n-арная подгруппа n-арной группы $< A, [\ ] >$. Ясно, что $< A(p_i), [\ ] >$ – $p_i$-силовская в $< A, [\ ] >$.

Применяя к разложению (**) теорему 2.7.20, получаем разложение (*) n-арной группы $< A, [\ ] >$ в a-прямое произведение своих $p_i$-силовских n-арных подгрупп.

Предположим, что

$$< A, [\,] > = < A'(p_1), [\,] > \overset{a}{\times} \ldots \overset{a}{\times} < A'(p_m), [\,] >$$

еще одно разложение n-арной группы $< A, [\,] >$ в a-прямое произведение своих $p_i$-силовских n-арных подгрупп $< A'(p_i), [\,] >$, содержащих идемпотент a, отличное от разложения (∗). Это означает, что существует, по крайней мере, один индекс j такой, что $A(p_j) \neq A'(p_j)$, откуда, учитывая $a \in A(p_j)$, $a \in A'(p_j)$, заключаем, что $< A(p_j), @ >$ и $< A'(p_j), @ > - $ различные p-силовские подгруппы группы $< A, @ >$, что невозможно в силу единственности $p_j$-силовской подгруппы в абелевой группе $< A, @ >$. Следовательно, предположение о существовании двух различных разложений неверно. ∎

Согласно следствию 2.7.22, n-арная группа $< A, [\,] >$, являющаяся a-прямым произведением своих n-арных подгрупп $< B_1, [\,] >, \ldots, < B_m, [\,] >$, изоморфна прямому произведению

$$< B_1, [\,] > \times \ldots \times < B_m, [\,] >.$$

Поэтому справедливо

**2.8.2. Следствие.** Конечная полуабелева n-арная группа $< A, [\,] >$ порядка $|A| = p_1^{\alpha_1} \ldots p_m^{\alpha_m}$ ($p_1, \ldots, p_m$ – простые), содержащая идемпотент a, изоморфна прямому произведению $< A(p_1), [\,] > \times \ldots \times < A(p_m), [\,] >$ своих $p_i$-силовских ($i = 1, \ldots, m$) n-арных подгрупп $< A(p_i), [\,] >$, содержащих идемпотент a.

**2.8.3. Лемма** [4]. Конечная n-арная группа $< A, [\,] >$, порядок g которой взаимно прост с $n - 1$, обладает, по крайней мере, одним идемпотентом.

*Доказательство.* Пусть b – произвольный элемент из A, $< B, [\,] > -$ циклическая n-арная подгруппа, порождённая b. Так как $|B|$ делит g и $(g, n-1) = 1$, то $(|B|, n-1) = 1$. Поэтому, согласно следствию 2.5.46, $< B, [\,] >$, а значит и $< A, [\,] >$ обладают идемпотентом. ∎

Лемма 2.8.3 позволяет сформулировать ещё два следствия.

**2.8.4. Следствие.** Конечная полуабелева n-арная группа $< A, [\,] >$ порядка $|A| = p_1^{\alpha_1}...p_m^{\alpha_m}$ ($p_1, ..., p_m$ – простые), взаимно простого с $n-1$, единственным образом разлагается в a-прямое произведение своих $p_i$-силовских ($i = 1, ..., m$) n-арных подгрупп для любого идемпотента $a \in A$.

**2.8.5. Следствие.** Конечная полуабелева n-арная группа $< A, [\,] >$ порядка $|A| = p_1^{\alpha_1}...p_m^{\alpha_m}$ ($p_1, ..., p_m$ – простые), взаимно простого с $n-1$, изоморфна прямому произведению своих $p_i$-силовских ($i = 1, ..., m$) n-арных подгрупп, содержащих идемпотент $a \in A$.

**2.8.6. Лемма.** Если конечная n-арная группа $< A, [\,] >$ порядка $g$, где $(n-1, g) = 1$, обладает единицей, то эта единица является единственным идемпотентом в $< A, [\,] >$.

*Доказательство.* Если $e$ – единица n-арной группы $< A, [\,] >$, то $e$ – единица группы $< A, * >$, для которой $< A, [\,] >$ – производная. Так как

$$a * e = [a \underbrace{e...e}_{n-2} e] = a, \qquad e * a = [\underbrace{e e ... e}_{n-2} a] = a$$

для любого $a \in A$, то $e$ – единица группы $< A, * >$. Если $e_1$ – идемпотент из $< A, [\,] >$, то

$$[\underbrace{e_1 ... e_1}_{n-1} e_1] = e_1, \qquad [\underbrace{e_1 ... e_1}_{n}] = e_1^n,$$

откуда $e_1^n = e_1$, $e_1^{n-1} = e$. Предположив, что $e_1 \neq e$, получим $n - 1 = mt$ и $m$ делит $g$, где $m > 1$ порядок циклической подгруппы, порождённой элементом $e_1$, что противоречит условию $(n - 1, g) = 1$. Следовательно, $e_1 = e$.

Так как в абелевой n-арной группе любой её идемпотент является единицей, то леммы 2.8.3 и 2.8.6 позволяют сформулировать ещё одну лемму.

**2.8.7. Лемма** [4]. Конечная абелева n-арная группа, порядок которой взаимно прост с $n - 1$, обладает единственным идемпотентом.

Теперь можно сформулировать

**2.8.8. Следствие.** Конечная абелева n-арная группа $< A, [\ ] >$ порядка $|A| = p_1^{\alpha_1}...p_m^{\alpha_m}$ ($p_1, ..., p_m$ – простые), взаимно простого с $n - 1$, единственным образом разлагается в прямое произведение (внутреннее) своих $p_i$-силовских ($i = 1, ..., m$) n-арных подгрупп.

В [26] установлено, что тернарная группа $< B_n, [\ ] >$ отражений правильного n-угольника является полуабелевой, и все её элементы – идемпотенты. Поэтому имеет место

**2.8.9. Следствие.** Если $n = p_1^{\alpha_1}...p_m^{\alpha_m}$ ($p_1, ..., p_m$ – простые), $b_j$ – фиксированный элемент из $B_n$ ($j = 1, ..., n$), то тернарная группа $< B_n, [\ ] >$ единственным образом разлагается в $b_j$-прямое произведение

$$< B_n, [\ ] > = < B(p_1), [\ ] ) \overset{b_j}{\times} ... \overset{b_j}{\times} < B(p_m), [\ ] >$$

своих $p_i$-силовских ($i = 1, ..., m$) тернарных подгрупп $< B(p_i), [\ ] >$.

**2.8.10. Пример.** Пусть $< B_6, [\ ] >$ – тернарная группа отражений правильного шестиугольника, все тернарные подгруппы которой исчерпываются тремя тернарными подгруппами.

$$< K_1 = \{b_1, b_4\}, [\ ] >, \ < K_2 = \{b_2, b_5\}, [\ ] >, \ < K_3 = \{b_3, b_6\}, [\ ] >$$

второго порядка и двумя подгруппами

$$< H_1 = \{b_1, b_3, b_5\}, [\ ] >, \ < H_2 = \{b_2, b_4, b_6\}, [\ ] >$$

третьего порядка. Так как $|B_6| = 6$, то все перечисленные тернарные подгруппы являются силовскими. Все элементы в $< B_6, [\ ] >$ являются идемпотентами, а сама она – полуабелева. Выпишем для каждого $b_j \in B_6$ (j = 1, ..., 6) соответствующее прямое разложение.

$$B_6 = \{b_1, b_4\} \overset{b_1}{\times} \{b_1, b_3, b_5\} = \{b_2, b_5\} \overset{b_2}{\times} \{b_2, b_4, b_6\} =$$

$$= \{b_3, b_6\} \overset{b_3}{\times} \{b_1, b_3, b_5\} = \{b_1, b_4\} \overset{b_4}{\times} \{b_2, b_4, b_6\} =$$

$$= \{b_2, b_5\} \overset{b_5}{\times} \{b_1, b_3, b_5\} = \{b_3, b_6\} \overset{b_6}{\times} \{b_2, b_4, b_6\}.$$

Отметим, что если в полуабелевой n-арной группе отсутствуют идемпотенты, то она не только не разлагается в a-прямое произведение своих силовских n-арных подгрупп, но в ней вообще могут отсутствовать силовские n-арные подгруппы. Это вытекает из теоремы 2.5.27, согласно которой существуют конечные циклические n-арные группы, являющиеся очевидно полуабелевыми, не обладающие ни одной n-арной подгруппой, в том числе и идемпотентами.

В следующей теореме, обобщающей теорему 2.8.1, через $\pi(A)$, как обычно, обозначается множество всех простых делителей порядка $|A|$ n-арной группы $< A, [\,] >$.

**2.8.11. Теорема.** Пусть $< A, [\,] >$ – конечная полуабелева n-арная группа, содержащая идемпотент a,

$$\pi(A) = \pi_1 \cup ... \cup \pi_m, \quad \pi_i \cap \pi_j = \varnothing \ (i \neq j).$$

Тогда $< A, [\,] >$ единственным образом разлагается в a-прямое произведение

$$< A, [\,] > \; = \; < A(\pi_1), [\,] > \overset{a}{\times} ... \overset{a}{\times} < A(\pi_m), [\,] >$$

своих $\pi_i$-холловских ( i = 1, ..., m ) n-арных подгрупп $< A(\pi_i), [\,] >$.

Доказательство теоремы 2.8.11 дословно повторяет доказательство теоремы 2.8.1.

# ДОПОЛНЕНИЯ И КОММЕНТАРИИ

**1.** Смежные классы n-арной группы по ее n-арной подгруппе впервые появились у Дёрнте [1]. Там же имеется доказательство теоремы 2.1.15. Теорему Лагранжа для n-арных групп получил Пост [3, с. 222].

**2.** Теорема 2.2.6 устанавливает соответствие между n-арными подгруппами n-арной группы $< A, [\ ] >$ и подгруппами группы $< A, @ >$, к которой согласно теореме Глускина-Хоссу приводима $< A, [\ ] >$, и является аналогом теоремы Поста [3, с. 222] о соответствии между n-арными подгруппами n-арной группы $< A, [\ ] >$ и подгруппами соответствующей группы Поста $A_0$.

**3.** Теорема 2.2.12 является аналогом результата В.А. Артамонова [17] о существовании биекции множества всех n-арных подгрупп n-арной группы $< A, [\ ] >$, содержащих x, на множество всех подгрупп соответствующей группы Поста $A_0$, содержащих $x^{n-1}$ и инвариантных относительно сопряжения с помощью x.

**4.** Утверждения 2.3.17 – 2.3.20 имеются в [4] и получены там, как следствия более общих результатов.

**5.** В [26] помещен созданный с помощью ЭВМ атлас, включающий в себя 34 рисунка, на которых для $n \leq 30$ содержится информация о всех тернарных подгруппах тернарной группы $< B_n, [\ ] >$ и их нормализаторах.

**6.** m-Полуинвариантные n-арные подгруппы n-арной группы впервые были определены и изучались в [28] по аналогии с m-полуабелевыми n-арными подгруппами Поста [3].

**7.** В [38] Г.Н. Воробьевым по аналогии с m-полуабелевыми и m-полуинвариантными n-арными подгруппами определяются и изучаются m-полусопряженные n-арные подгруппы n-арной группы. При этом сопряженные и полусопряженные n-арные подгруппы n-арной группы являются частными случаями этого нового понятия.

**8.** Степени элементов n-арной группы впервые появились у Дернте [1, 6]. Конечные циклические n-арные группы определил и подробно изучил Пост [3], а бесконечные циклические n-арные группы изучал С.А. Русаков [4], доказавший, в частности, изоморфизм бесконечных циклических n-арных групп.

**9.** Из утверждения 1) леммы 2.5.25 следует, что для любого n ≥ 3 существует циклическая n-арная группа любого порядка.

**10.** Доказанная в теореме 2.5.47 единственность n-арных подгрупп одного и того же порядка в циклической n-арной группе фактически является следствием следующего общего утверждения.

**Предложение.** Если $<B_1, [\ ]>$ и $<B_2, [\ ]>$ – n-арные подгруппы n-арной группы $<A, [\ ]>$, то $B_1 = B_2$ тогда и только тогда, когда $B_1^*(A) = B_2^*(A)$.

**11.** Полуциклические n-арные группы можно определить, используя понятие степени последовательности [39].

**12.** Приведем еще один критерий m-полуабелевости n-арной группы.

**Теорема** [32]. Если $n = k(m-1) + m$, $n > 2$, $k \geq 0$, то n-арная группа $<A, [\ ]>$ является (k+2)-полуабелевой тогда и только тогда, когда

$$[x_1\ldots x_n]_n = [x_1 x_2^{\alpha} \ldots x_{k+1}^{\alpha^k} x_{k+2} x_{k+3}^{\alpha^{k+2}} \ldots x_n^{\alpha^{n-1}} c^{(n-k-1)(m-1)}]_m,$$

где $<A, [\ ]_m>$ – абелева m-арная группа, а элемент $c \in A$ и автоморфизм $\alpha$ m-арной группы $<A, [\ ]_m>$ удовлетворяют условиям

$$c^{\alpha} = c,\ x_k^{\alpha^{k+1}} = x,\ x \in A.$$

**Следствие** [23]. n-Арная группа $<A, [\ ]>$ является полуабелевой тогда и только тогда, когда

$$[x_1 x_2 \ldots x_n] = x_1 + x_2^{\alpha} + \ldots + x_{n-1}^{\alpha^{n-2}} + x_n + c,$$

где $<A, +>$ – абелева группа, $c \in A$, $\alpha$ – автоморфизм $<A, +>$, причем $c^{\alpha} = c$, $x^{\alpha^{n-1}} = x$ для всех $x \in A$.

Свойства m-полуабелевых n-арных групп изучал также В.А. Дудек [40]. Им, в частности, доказано существование первообразных, т. е. неприводимых m-полуабелевых n-арных групп для любого n ≥ 3.

**13.** Слабо m-полуабелевы n-арные группы впервые были определены и изучались А.М. Гальмаком [32]. Дальнейшее развитие теория слабо m-полуабелевых n-арных групп получила в работах

В.А. Дудека, среди которых отметим [41]. В этой работе найдены необходимые и достаточные условия слабой m-полуабелевости n-арной группы, а также доказано существование первообразных слабо m-полуабелевых n-арных групп для любого $n \geq 3$.

**14.** Пусть $< A, [\,] > -$ n-арная группа ($n \geq 2$). Преобразования множества A, имеющие вид

$$x \mapsto [a_1 \ldots a_{i-1} x a_{i+1} \ldots a_n], \qquad (1)$$

где $a_1, \ldots, a_{i-1}, a_{i+1}, \ldots, a_n \in A$, $i \in \{1, \ldots, n\}$, называются [42] главными трансляциями n-арной группы $< A, [\,] >$ и обозначаются через $t(a_1^{i-1}, a_{i+1}^n)$. Композиция конечного числа главных трансляций называется [42] элементарной трансляцией. Выделим также преобразования

$$t(a_1^{n-1}, a_{n+1}^{2n-1}) : x \mapsto [a_1 \ldots a_{n-1} x a_{n+1} \ldots a_{2n-1}], \qquad (2)$$

где $a_1, \ldots, a_{i-1}, a_{i+1}, \ldots, a_n \in A$. Введем следующие обозначения: $T_i(A)$ – множество всех преобразований вида (1); $T_0(A)$ – множество всех преобразований вида (2); $T(A)$ – множество всех элементарных трансляций n-арной группы $< A, [\,] >$.

Известно [43], что $T(A)$ – транзитивная группа. Строение этой группы описывают следующие две теоремы.

**Теорема** [19]. Пусть $< A, [\,] > -$ n-арная группа. Тогда справедливы следующие утверждения:

1) $T_1(A)$ и $T_n(A)$ – транзитивные нормальные подгруппы группы $T(A)$;

2) подгруппа $T_1(A) \cap T_n(A)$ абелева;

3) подгруппы $T_1(A)$ и $T_n(A)$ изоморфны;

4) $T_i(A) T_{n+1-i}(A) = T_{n+1-i}(A) T_i(A) = T_0(A)$ – нормальная подгруппа группы $T(A)$, $i \in \{1, 2, \ldots, n\}$;

5) если $n \geq 3$, то $T_i(A) = \underbrace{T_2(A) \ldots T_2(A)}_{i-1} = \quad$, $i \in \{2, \ldots, n-1\}$,

$T_0(A) = \underbrace{T_2(A) \ldots T_2(A)}_{n-1} = \underbrace{T_{n-1}(A) \ldots T_{n-1}(A)}_{n-1}$;

6) $T(A)/T_0(A)$ – циклическая группа порядка, делящего $n - 1$, причем

$T(A)/T_0(A) = \{T_0(A), T_2(A), \ldots, T_{n-1}(A)\} = < T_2(A) > = < T_{n-1}(A) >$

при $n \geq 3$, $T(A)/T_0(A) = \{T(A)\}$ при $n = 2$.

**Теорема** [19]. Для всякой n-арной группы $< A, [\ ] >$ множества $T_2(A)$ и $T_{n-1}(A)$ являются антиизоморфными n-арными группами относительно n-арной операции

$$[u_1 u_2 \ldots u_n] = u_1 u_2 \ldots u_n,$$

при $n \geq 3$ обертывающей для каждой из них является группа $T(A)$, а соответствующей – её подгруппа $T_0(A)$.

Теперь мы можем сформулировать критерии различных типов абелевости на языке трансляций.

**Теорема** [19]. Пусть $< A, [\ ] >$ – m-полуабелева n-арная группа. Тогда $T_1(A) = T_m(A) = T_{2(m-1)+1}(A) = \ldots = T_n(A) = T_0(A)$.

**Теорема** [19]. Для n-арной группы $< A, [\ ] >$ следующие утверждения равносильны:
1) n-арная группа $< A, [\ ] >$ – полуабелева;
2) $T_1(A) = T_n(A)$;
3) группа $T_1(A)$ – абелева;
4) группа $T_n(A)$ – абелева.

**Теорема** [19]. Для n-арной группы $< A, [\ ] >$ следующие утверждения равносильны:
1) n-арная группа $< A, [\ ] >$ – абелева;
2) группа $T(A)$ – абелева;
3) n-арная группа $< T_2(A), [\ ] >$ – абелева;
3) n-арная группа $< T_{n-1}(A), [\ ] >$ – абелева.

**15.** Ю.И. Кулаженко показал [44], что полуабелевы n-арные группы можно определить как n-арные группы, для которых хотя бы один из четырехугольников специального вида является параллелограммом. Он же определяет полуабелевость n-арной группы с помощью векторных равенств на ней.

**16.** Свободные абелевы n-арные группы и тензорное произведение абелевых n-арных групп изучал Сиосон [45 – 47].

**17.** Если в декартовом произведении $\prod_{i \in I} \mathcal{A}_i$ каждая n-арная группа $\mathcal{A}_i$ обладает по крайней мере одним идемпотентом, то для любой функции $a \in \prod A_i$, расширяя соответствующее бинарное понятие, можно ввести понятие носителя. Далее, как и в теории групп, рассматривают функции с конечными носителями и определяют [4] еще один тип произведений n-арных групп – прямое произведение . Ос-новные свойства прямых произведений n-арных групп установлены С.А. Русаковым [4].

В книге С. А. Русакова [4] для обозначения декартова произведения n-арных групп, по аналогии с теорией групп, используется символ $\overline{\prod} \mathscr{A}_i$, а символом $\prod \mathscr{A}_i$, обозначается подгруппа из $\overline{\prod} \mathscr{A}_i$, являющаяся прямым произведением n-арных групп $\mathscr{A}_i$. Как и в теории групп, декартово и прямое произведения n-арных групп совпадают, если множество I конечное.

Для всякого целого m ≥ 2 определим функцию $\varphi_m(n)$ натурального аргумента n, полагая значение $\varphi_m(n)$ равным числу чисел вида k(m – 1) + 1, где k = 0, 1, …, n – 1, взаимно простых с n. Иначе говоря $\varphi_m(n)$ – это количество чисел из множества

$$\{1, m - 1, 2(m - 1) + 1, …, (n - 1)(m - 1) + 1\},$$

взаимно простых с n. В частности $\varphi_2(n)$ совпадает с количеством чисел из множества {1, 2, …, n}, взаимно простых с n. Следовательно, $\varphi_2(n)$ совпадает с функцией Эйлера $\varphi(n)$.

Теперь мы можем сформулировать n-арный аналог бинарного результата о числе порождающей конечной циклической n-арной группы.

**Предложение.** Число порождающих циклической n-арной группы порядка g равно $\varphi_n(g)$.

# ЛИТЕРАТУРА

1. **Dörnte W.** Untersuchungen über einen verallgemeinerten Gruppenbegrieff // Math. Z. 1928. Bd. 29. S. 1-19.

2. **Prüfer H.** Theorie der abelshen Gruppen. I. Grundeigenschaften // Math. Z. 1924. Bd. 20. S. 165-187.

3. **Post E.L.** Polyadic groups // Trans. Amer. Math. Soc. 1940. Vol. 48, N2. P.208-350.

4. **Русаков С.А.** Алгебраические n-арные системы. Мн.: Навука i тэхнiка, 1992. 245 с.

5. **Русаков С.А.** Некоторые приложения теории n-арных групп. Мн.: Беларуская навука, 1998. 167 с.

6. **Сушкевич А.К.** Теория обобщенных групп. Харьков; Киев, 1937.

7. **Курош А.Г.** Общая алгебра: Лекции 1969/70 учебного года. М.: Наука, 1974. 160 с.

8. **Bruck R.H.** A survey of binary systems. Berlin-Heldelberg-New York: Springer-Verlad, 1966. 185 p.

9. **Бурбаки Н.** Алгебра. Алгебраические структуры, линейная и полилинейная алгебра. М.: Физматгиз, 1962. 516 с.

10. **Артамонов В.А.** Универсальные алгебры //Итоги науки и техники. Сер. Алгебра. Топология. Геометрия. 1976. С. 191-248.

11. **Glazek K.** Bibliographi of n-groups (poliadic groups) and same group like n-ary sistems// Proc. of the sympos. n-ary structures. Skopje, 1982. P. 259-289.

12. **Гальмак А. М.** Конгруэнции полиадических групп. Мн.: «Беларуская навука», 1999. 182 с.

13. **Глускин Л.М.** Позиционные оперативы// Мат.сборник. 1965. Т.68(110), №3. С.444-472.

14. **Hosszu M.** On the explicit form of n-group operations //Publ. Math. 1963. V.10, №1- 4. P.88-92.

15. **Гальмак А.М.** О приводимости n-арных групп // Вопр. алгебры. 1996. Вып. 10. С. 164-169.

16. **Соколов Е.И.** О теореме Глускина-Хоссу для n-групп Дёрнте //Мат. исследования. Вып.39. С.187-189.

17. **Артамонов В.А.** Свободные n-арные группы //Мат. заметки. 1970. Т.8, №4. С. 499-507.

18. **Артамонов В.А.** О шрайеровых многообразиях n-групп и n-полугрупп // Труды семинара им. И.Г. Петровского. 1979. Вып.5. С. 193-202.

19. **Гальмак А.М.** Трансляции n-арных групп //Докл. АН БССР. 1986. Т.30, №8. С. 677-680.

20. **Гальмак А.М.** О приводимости n-арных групп // Препринты ИМ АН БССР, 1976. 6(242). 36 с.

21. **Гальмак А.М.** Приводимость полиадических групп //Докл. АН БССР. 1985. Т.29, №10. С. 874–877.

22. **Dudek W.A., Michalski J.** On a generalisation of Hosszu theorem // Denconstratio Math. 1982. Vol. 15, №3. P. 783 – 805.

23. **Дириенко И.И., Колесников О.В.** К теореме Глускина-Хоссу об n-группах // Деп. в ВИНИТИ №374 – 80. Харьков, 1980. 10 с.

24. **Гальмак А.М.** Теоремы Поста и Глускина-Хоссу. Гомель, 1997. 85 с.

25. **Гальмак А.М.** Тернарные группы отражений. //Междунар. Мат. Конф. Тез. докл. Гомель, 1994. С. 33.

26. **Гальмак А. М., Воробьёв Г. Н.** Тернарные группы отражений. Мн.: «Беларуская навука», 1998. 128с.

27. **Masat Fransis E.** A useful characterisation of a normal subgroups // Math. Mag.1979. Vol 52, № 3. P. 171-173.

28. **Гальмак А.М.** Инвариантные подгруппы n-арных групп и их обобщения // Вопросы алгебры. Мн.: Университетское, 1990. Вып. 5. С. 91-94.

29. **Воробьев Г.Н.** О сопряженности n-арных подгрупп // Весці Акадэміі навук Беларусі. Сер. фіз.-мат. навук. 1996. №1. С. 121.

30. **Воробьев Г.Н.** О полусопряженности n-арных подгрупп. // Вопросы алгебры. Гомель, 1997. Вып. 10. С.157 – 163.

31. **Серпинский В.** 250 задач по элементарной теории чисел. М.: Просвещение, 1968. 160 с.

32. **Гальмак А.М.** Абелевы n-арные группы и их обобщения // Вопросы алгебры. Минск: Университетское, 1987. Вып. 3. С. 86-93.

33. **Dudek W.A.** Remarcs on n-groups // Demonstratio Math. Vol. 13, №1. 1980. P.165–181.

34. **Колесников О.В.** Разложение n-групп // Мат. исслед. Вып. 51. Квазигруппы и лупы. Кишинёв: Штиинца, 1979. С. 88–92.

35. **Плоткин Б.И.** Группы автоморфизмов алгебраических систем. М. Наука, 1966. 603 с.

36. **Glazek K., Gleichgewicht B.** Abelian n-groups // Proc. Congr. Math. Soc. J. Bolyai. Esztergom. 1977. P. 321 – 329.

37. **Гальмак А.М.** Полуабелевые n-арные группы с идемпотентами // Веснік ВДУ ім. П.М. Машэрава. 1999. № 2(12). С. 56–60.

38. **Воробьев Г.Н.** Сопряженные n-арные подгруппы и их обобщения. // Веснік ВДУ ім. П.М. Машэрава. 1997. № 2(4). С.59–64.

39. **Гаврилов В.В.** О полуциклических n-арных группах // Конф. математиков Беларуси. Тез. докл. Гродно, 1992. С. 15.

40. **Дудек В.А.** m-Полуабелевые n-арные группы // Изв. АН ССР Молдова. Математика. 1990. №2. С. 66 – 70.

41. **Dudek W.A.** On the class of weakly semiabelian polyadic groups // Discrete Math. Appl. Vol. 6, №5. P. 427 – 433.

42. **Мальцев А.И.** К общей теории алгебраических систем // Мат. сб. 1954. Т.35, №1. С.3-20.

43. **Monk J.D., Sioson F.M.** On the general theory of m-groups // Fund. Math: 1971. №72. С. 233-244.

44. **Кулаженко Ю.И.** Критерии полуабелевости n-арной группы // Веснік ВДУ ім. П.М. Машэрава. 1997. №3(5). С. 61–64.

45. **Sioson F.M.** On Free Abelian n-Groups I // Proc. Japan Acad. 1967. Vol. 43. С. 876–879.

46. **Sioson F.M.** On Free Abelian n-Groups II // Proc. Japan Acad. 1967. Vol. 43. С. 880–883.

47. **Sioson F.M.** On Free Abelian n-Groups III // Proc. Japan Acad. 1967. Vol. 43. С. 884–888.

# ПРЕДМЕТНЫЙ УКАЗАТЕЛЬ

# УСЛОВНЫЕ ОБОЗНАЧЕНИЯ

$$a_m^k = \begin{cases} a_m a_{m+1} \ldots a_k, & \text{если } m \leq k, \\ \varnothing, & \text{если } m > k; \end{cases}$$

$$\overset{k}{a} = \begin{cases} \underbrace{a \ldots a}_{k}, & k > 0, \\ \varnothing, & k = 0. \end{cases}$$

$$a^{[s]} = \begin{cases} a, & s = 0, \\ [\overset{s(n-1)+1}{a}], & s > 0, \\ [\overset{-2s}{\bar{a}}\; \overset{-s(n-3)+1}{a}], & s < 0. \end{cases}$$

$$\overset{m}{B} = \begin{cases} \underbrace{B \ldots B}_{m}, & \text{если } m \geq 1; \\ \varnothing, & \text{если } m \leq 0; \end{cases}$$

$\bar{a}$ – косой элемент для элемента $a$;
$\alpha^{-1}$ – обратная последовательность для последовательности $\alpha$;
$l(\alpha)$ – длина последовательности $\alpha$;
$F_A$ – свободная полугруппа над алфавитом $A$;
$\sim$ или $\theta$ – отношение эквивалентности Поста на $F_A$;
$\mathscr{A} = F_A / \theta$;
$A^{(i)} = \{\theta(\alpha a) \mid a \in A\} = \{\theta(a\alpha) \mid a \in A\}$;
$A_n$ – знакопеременная группа степени $n$;
$< B_n, [\;] > $ – тернарная группа всех отражений правильного n-угольника;
$C_n$ – циклическая группа порядка $n$;
$D_n$ – диэдральная группа порядка $2n$;
$S_n$ – симметрическая группа степени $n$;
$[\;]$ – n-арная операция;
$< A, [\;] >$ – n-арная группа;
$x@y = [xa^{-1}y]$, где $a$ – элемент n-арной группы $< A, [\;] >$;
$<< M >, [\;] >$ – n-арная подгруппа, порождённая множеством $M$;

$< \langle a \rangle, [\ ] > $ – циклическая n-арная группа, порожденная элементом a;

$\overline{M} = \{\overline{a} \mid a \in M\}$;

$B_a = \{[b_1 \ldots b_{n-1}a] \mid b_1, \ldots, b_{n-1} \in B\}$;

$_aB = \{[ab_1 \ldots b_{n-1}] \mid b_1, \ldots, b_{n-1} \in B\}$;

$\widetilde{B} = \{[aba^{-1}] \mid b \in B\}$;

$\mathring{B} = \{[a^{-1}ba] \mid b \in B\}$;

$B^{(i)}(A) = \{\theta_A(\alpha) \in A^{(i)} \mid \exists b_1, \ldots, b_i \in B, \alpha\theta_A b_1 \ldots b_i\}$, $i = 1, \ldots, n-1$;

$B_o(A) = B^{(n-1)}(A) = \{\theta_A(\alpha) \in A_o \mid \exists b_1, \ldots, b_{n-1} \in B, \alpha\theta_A b_1 \ldots b_{n-1}\}$;

$B^*(A) = \{\theta_A(\alpha) \in A^* \mid \exists b_1, \ldots, b_i \in B \ (i \geq 1), \alpha\theta_A b_1 \ldots b_i\}$;

$[a \overset{n-1}{B}]$ – левый смежный класс n-арной группы $< A, [\ ] >$ по её n-арной подгруппе $< B, [\ ] >$;

$[\overset{n-1}{B} a]$ – правый смежный класс n-арной группы $< A, [\ ] >$ по её n-арной подгруппе $< B, [\ ] >$;

$|A : B|$ – индекс n-арной подгруппы $< B, [\ ] >$ в n-арной группе $< A, [\ ] >$;

$< B, [\ ] > \vee < C, [\ ] >$ – n-арная подгруппа, порождённая n-арными подгруппами $< B, [\ ] >$ и $< C, [\ ] >$;

$L(A, [\ ])$ – множество всех n-арных подгрупп n-арной группы $< A, [\ ] >$;

$L(A, [\ ], x)$ – множество всех n-арных подгрупп n-арной группы $< A, [\ ] >$, содержащих x;

$L(A, \text{@}, x)$ – множество всех подгрупп $< B, \text{@} >$ группы $< A, \text{@} >$, удовлетворяющих условиям

$$[\underbrace{x \ldots x}_{n-1} a] \in B, \ B \text{@} x = x \text{@} \widetilde{B};$$

$< B_1, [\ ] > \overset{a}{\times} \ldots \overset{a}{\times} < B_m, [\ ] >$ – a-прямое произведение n-арных подгрупп $< B_1, [\ ] >, \ldots, < B_m, [\ ] >$;

$A(\pi)$ – $\pi$-холловская n-арная подгруппа в $< A, [\ ] >$;

$A(p)$ – p-силовская n-арная подгруппа в $< A, [\ ] >$.

# СОДЕРЖАНИЕ

А.М. Гальмак

# n-АРНЫЕ ГРУППЫ

## ЧАСТЬ 2

Приведены новые результаты о разрешимости в n-арной полугруппе уравнений с числом неизвестных большим единицы, n-арных подстановках и n-арных морфизмах, n-арных подгруппах, смежных классах, n-арных аналогах нормализатора подмножества в группе и центра группы. Как и в первой части, много внимания уделено изучению связи между полиадическими аналогами бинарных понятий и результатов и их прототипами в группах, к которым приводима n-арная группа согласно теоремам Поста и Глускина-Хоссу. Рассмотрены вопросы, связанные со строением идемпотентных n-арных групп, в том числе n-арных групп, допускающих регулярный автоморфизм.

Библиогр.: 138 назв.

# ВВЕДЕНИЕ

Данная книга является продолжением изданной в 2003 году книги "n-Арные группы. Часть 1" и составляет с ней единое целое. Поэтому во второй части, включающей семь глав с третьей по девятую, нумерация глав и параграфов продолжает соответствующую нумерацию глав и параграфов в первой части.

В главе 3 изучается разрешимость в n-арной полугруппе уравнений с числом неизвестных большим единицы, в частности, равным $n - 1$. Показано, что n-арную группу можно определить как n-арную полугруппу, в которой разрешимы такие уравнения. Установлено, что почти все известные определения n-арной группы являются непосредственными следствиями этого результата. Приведено большое число новых определений n-арной группы.

Глава 4 посвящена n-арным подстановкам, значительный вклад в изучение которых внесли Э. Пост и С.А. Русаков, а также n-арным морфизмам. Из результатов этой главы можно отметить n-арный аналог теоремы Биркгофа, утверждающий, что всякая n-арная группа изоморфна n-арной группе автоморфизмов некоторой последовательности универсальных алгебр.

В главе 5 приводятся различные критерии существования n-арных подгрупп в n-арной группе. Определяются и изучаются новые n-арные аналоги нормальных подгрупп. Исследуется связь между сопряженностью и полусопряженностью n-арных подгрупп в n-арной группе и сопряженностью подгрупп в универсальной обертывающей группе Поста.

Установлению связей между разложениями n-арной группы по ее n-арной подгруппе и соответствующими разложениями в универсальной обертывающей группе Поста посвящена глава 6.

В седьмой и восьмой главах определяются различные n-арные аналоги нормализатора подмножества в группе и центра группы и изучается их связь со своими бинарными прототипами в группах, к которым приводима n-арная группа согласно теоремам Поста и Глускина-Хоссу.

Предметом изучения заключительной главы 9 являются n-арные группы с идемпотентами, в том числе n-арные группы, все элементы которых являются идемпотентами, и n-арные группы, допускающие автоморфизм с единственным неподвижным элементом.

Во введении к первой части отмечалось, что дальнейшему прогрессу в изучении n-арных групп будет способствовать издание новых книг по n-арным группам. В связи с этим можно указать опубликованную в 2003 году монографию Я. Ушана "n-Groups in the light of the neutral operations", в которой ее автор для изучения n-арных групп и близких к ним алгебраических систем использует введенные им нейтральные и обратные операции. В 2005 году появилась электронная версия этой книги.

# Г Л А В А  3

# ОПРЕДЕЛЕНИЯ n-АРНОЙ ГРУППЫ

Основным результатом данной главы является теорема, характеризующаяся тем, что почти все известные до сих пор определения n-арной группы, являются непосредственными ее следствиями. Кроме того, эта теорема позволяет дать большое число новых определений n-арной группы.

## §3.1. ОСНОВНАЯ ТЕОРЕМА

Согласно Дёрнте [1], n-арная полугруппа $< A, [\ ] >$ называется n-арной группой, если каждое из уравнений

$$[a_1 \ldots a_{i-1} x_i a_{i+1} \ldots a_n] = b,\ i = 1, 2, \ldots, n$$

однозначно разрешимо в ней относительно $x_i$ для всех $a_1, \ldots, a_{i-1}, a_{i+1}, \ldots, a_n, b \in A$ (определение 1.1.1).

Пост заметил [2], что n-арную группу $< A, [\ ] >$ можно определить как n-арную полугруппу, в которой разрешимы уравнения

$$[x a_2 \ldots a_n] = b, \qquad (i)$$

$$[a_1 \ldots a_{n-1} y] = b \qquad (ii)$$

для всех $a_1, \ldots, a_n, b \in A$ (определение 1.1.2).

В [2] Пост также заметил, что n-арную группу $< A, [\ ] >$ можно определить как n-арную полугруппу, в которой разрешимо уравнение

$$[a_1 \ldots a_{i-1} x a_{i+1} \ldots a_n] = b$$

для всех $a_1, \ldots, a_{i-1}, a_{i+1}, \ldots, a_n, b \in A$ и некоторого $i \in \{2, \ldots, n-1\}$ (определение 1.1.3).

Понятно, что уравнения в определении 1.1.1 Дернте, а также заменяющие его уравнения (i) и (ii) в определении 1.1.2 Поста суть n-арные аналоги соответствующей групповой аксиомы о разрешимости уравнений

$$xa = b, \quad ay = b.$$

Нетрудно заметить, что n-арными аналогами последних являются также уравнения

$$[x_1 \ldots x_{n-1}a] = b, \qquad (j)$$

$$[ay_1 \ldots y_{n-1}] = b \qquad (jj)$$

с $n-1$ неизвестными.

**3.1.1. Теорема** [48]. n-Арная полугруппа $<A, [\ ]>$ является n-арной группой тогда и только тогда, когда для любых $a, b \in A$ в $A$ разрешимы уравнения (j) и (jj).

*Доказательство. Необходимость.* Пусть $<A, [\ ]>$ – n-арная группа, $d$ – произвольный элемент из $A$. Определение 1.1.2 Поста гарантирует существование решений $d_1$ и $d_{n-1}$ уравнений

$$[x_1 \underbrace{d \ldots d}_{n-2} a] = b,$$

$$[a \underbrace{d \ldots d}_{n-2} y_{n-1}] = b$$

соответственно, то есть $x_1 = d_1$, $y_{n-1} = d_{n-1}$. Тогда

$$x_1 = d_1, \ x_2 = \ldots = x_{n-1} = d$$

– решение уравнения (j),

$$y_1 = \ldots = y_{n-2} = d, \ y_{n-1} = d_{n-1}$$

– решение уравнения (jj).

*Достаточность.* Пусть теперь $< A, [\ ] >$ – n-арная полугруппа, в которой для любых $a, b \in A$ разрешимы уравнения (j) и (jj). Для доказательства теоремы достаточно показать, что в $A$ разрешимы оба уравнения (i) и (ii) из определения 1.1.2 Поста.

Покажем разрешимость уравнения (i). Из (jj) следует, что существуют элементы $c_1, \ldots, c_{n-1} \in A$ такие, что

$$[bc_1 \ldots c_{n-1}] = b. \qquad (1)$$

Аналогично из (j) вытекает существование элементов $c'_1, \ldots, c'_{n-1} \in A$ таких, что

$$[c'_1 \ldots c'_{n-1} a_n] = c_{n-1},$$

откуда и из (1) получаем

$$[bc_1 \ldots c_{n-2}[c'_1 \ldots c'_{n-1} a_n]] = b.$$

Теперь, используя ассоциативность n-арной операции $[\ ]$, получим

$$[bc_1 \ldots c_{n-3}[c_{n-2} c'_1 \ldots c'_{n-1}] a_n] = b. \qquad (2)$$

Снова применяя (j), устанавливаем существование элементов $c''_1, \ldots, c''_{n-1} \in A$ таких, что

$$[c''_1 \ldots c''_{n-1} a_{n-1}] = [c_{n-2} c'_1 \ldots c'_{n-1}],$$

откуда, используя (2), а также ассоциативность n-арной операции $[\ ]$, последовательно получаем

$$[bc_1 \ldots c_{n-3}[c''_1 \ldots c''_{n-1} a_{n-1}] a_n] = b,$$

$$[bc_1 \ldots c_{n-4}[c_{n-3} c''_1 \ldots c''_{n-1}] a_{n-1} a_n] = b.$$

Рассуждая аналогично, на $(n-2)$-ом шаге устанавливаем существование элементов $c_1^{(n-2)}, \ldots, c_{n-1}^{(n-2)} \in A$ таких, что

$$[b[c_1 c_1^{(n-2)} \ldots c_{n-1}^{(n-2)}] a_3 \ldots a_n] = b. \qquad (3)$$

Еще раз применяя (j), находим элементы $c_1^{(n-1)}, \ldots, c_{n-1}^{(n-1)} \in A$ такие, что

$$[c_1^{(n-1)} \ldots c_{n-1}^{(n-1)} a_2] = [c_1 c_1^{(n-2)} \ldots c_{n-1}^{(n-2)}],$$

откуда, используя (3) и ассоциативность n-арной операции [ ], получаем

$$[b[c_1^{(n-1)} \ldots c_{n-1}^{(n-1)} a_2] a_3 \ldots a_n] = b,$$

$$[[b c_1^{(n-1)} \ldots c_{n-1}^{(n-1)}] a_2 \ldots a_n] = b. \qquad (4)$$

Из (4) заключаем, что

$$x = [b c_1^{(n-1)} \ldots c_{n-1}^{(n-1)}]$$

– решение уравнения (i) определения 1.1.2 Поста.

Покажем теперь разрешимость уравнения (ii) из определения 1.1.2 Поста. Из (j) вытекает существование элементов $e_1, \ldots, e_{n-1} \in A$ таких, что

$$[e_1 \ldots e_{n-1} b] = b, \qquad (5)$$

а из (jj) вытекает существование элементов $e_1', \ldots, e_{n-1}' \in A$, таких, что

$$[a_1 e_1' \ldots e_{n-1}'] = e_1.$$

Подставляя последнее выражение для $e_1$ в (5) и, используя ассоциативность n-арной операции [ ], получаем

$$[[a_1 e_1' \ldots e_{n-1}'] e_2 \ldots e_{n-1} b] = b,$$

$$[a_1 [e_1' \ldots e_{n-1}' e_2] e_3 \ldots e_{n-1} b] = b. \qquad (6)$$

Снова, применяя (jj), находим элементы $e_1'', \ldots, e_{n-1}'' \in A$ такие, что

$$[a_2 e_1'' \ldots e_{n-1}''] = [e_1' \ldots e_{n-1}' e_2].$$

Тогда из (6) вытекает

$$[a_1 [a_2 e_1'' \ldots e_{n-1}''] e_3 \ldots e_{n-1} b] = b,$$

$$[a_1 a_2 [e_1'' \ldots e_{n-1}'' e_3] e_4 \ldots e_{n-1} b] = b.$$

На (n – 2)-ом шаге устанавливаем существование элементов $e_1^{(n-2)}, \ldots, e_{n-1}^{(n-2)} \in A$ таких, что

$$[a_1 \ldots a_{n-2} [e_1^{(n-2)} \ldots e_{n-1}^{(n-2)} e_{n-1}] b] = b. \qquad (7)$$

Еще раз применив (jj), находим элементы $e_1^{(n-1)}, \ldots, e_{n-1}^{(n-1)} \in A$ такие, что

$$[a_{n-1} e_1^{(n-1)} \ldots e_{n-1}^{(n-1)}] = [e_1^{(n-2)} \ldots e_{n-1}^{(n-2)} e_{n-1}],$$

откуда, используя (7) и ассоциативность n-арной операции [ ], получаем

$$[a_1 \ldots a_{n-2} [a_{n-1} e_1^{(n-1)} \ldots e_{n-1}^{(n-1)}] b] = b,$$

$$[a_1 \ldots a_{n-1} [e_1^{(n-1)} \ldots e_{n-1}^{(n-1)} b]] = b. \qquad (8)$$

Из (8) заключаем, что

$$y = [e_1^{(n-1)} \ldots e_{n-1}^{(n-1)} b]$$

есть решение уравнения (ii) определения 1.1.2 Поста. ∎

Число неизвестных в уравнениях (j) и (jj) можно уменьшить до двух.

**3.1.2. Предложение** [48]. n-Арная полугруппа $< A, [\,] >$ является n-арной группой тогда и только тогда, когда для любых $a, b \in A$ в $A$ разрешимы уравнения

$$[x \underbrace{u \ldots u}_{n-2} a] = b, \quad [a \underbrace{v \ldots v}_{n-2} y] = b.$$

Справедливо также

**3.1.3. Предложение.** n-арная полугруппа $< A, [\,] >$ является n-арной группой тогда и только тогда, когда для любых $a, b \in A$ в $A$ разрешима система

$$\begin{cases} [x\underbrace{z\ldots z}_{n-2}a] = b, \\ [a\underbrace{z\ldots z}_{n-2}y] = b. \end{cases}$$

Заметим, что в предложениях 3.1.2 и 3.1.3 переменная x может стоять на любом месте слева от a, переменная y может стоять на любом месте справа от a.

Рассмотрим n-арные полугруппы < A, [ ] >, в которых для любых a, b ∈ A разрешимы уравнения

$$[\underbrace{x\ldots x}_{n-1}a] = b, \qquad (*)$$

$$[a\underbrace{y\ldots y}_{n-1}] = b. \qquad (**)$$

Возникает естественный

**3.1.4. Вопрос**. Можно ли в теореме 3.1.1 заменить уравнения (j) и (jj) внешне более простыми уравнениями (*) и (**)?

Так как из разрешимости уравнений (*) и (**) следует разрешимость уравнений (j) и (jj), то по теореме 3.1.1, n-арные полугруппы, в которых разрешимы уравнения (*) и (**), являются n-арными группами. Это утверждение может быть доказано и независимо от теоремы 3.1.1.

Следующее предложение утверждает, что класс всех n-арных групп, в которых разрешимы уравнения (*) и (**), уже класса всех n-арных групп, то есть ответ на вопрос 3.1.4 – отрицательный.

**3.1.5. Предложение** [48]. Для любого n ≥ 3 существует n-арная группа, в которой неразрешимы уравнения (*) и (**).

*Доказательство.* Пусть A – циклическая группа порядка $n - 1$ ($n \geq 3$), то есть $x^{n-1} = 1$ для любого $x \in A$. Определим на A n-арную операцию

$$[x_1 x_2 \ldots x_n] = x_1 x_2 \ldots x_n,$$

производную от операции в группе A. Тогда $< A, [\ ] >$ – n-арная группа, причем

$$[\underbrace{c \ldots c}_{n-1} a] = c^{n-1} a = a, \qquad [a \underbrace{c \ldots c}_{n-1}] = a c^{n-1} = a$$

для любых $c, a \in A$, откуда вытекает, что в A неразрешимы уравнения

$$[\underbrace{x \ldots x}_{n-1} a] = b, \qquad [a \underbrace{y \ldots y}_{n-1}] = b.$$

где $a \neq b$. Такой выбор элементов $a \neq b$ возможен, так как $|A| = n - 1 \geq 2$. ∎

Представляется естественным также следующий

**3.1.6. Вопрос.** Можно ли в теореме 3.1.1 вместо разрешимости уравнений (j) и (jj) потребовать их однозначную разрешимость? Известно, что в бинарном случае ($n = 2$) это действительно так. Однако, как показывает следующая теорема, при $n > 2$ ответ на вопрос 3.1.6 является отрицательным.

**3.1.7. Теорема** [48]. Тогда и только тогда в n-арной полугруппе $< A, [\ ] >$ ($n > 2$) однозначно разрешимы уравнения (j) и (jj), когда A – одноэлементное множество.

*Доказательство. Необходимость.* По теореме 3.1.1, n-арная полугруппа $< A, [\ ] >$, в которой разрешимы уравнения (j) и (jj) (необязательно однозначно) является n-арной группой. Предположим, что множество A содержит более одного элемента, и пусть

$$x_1 = a_1, \quad x_2 = a_2, \ldots, x_{n-1} = a_{n-1}$$

решение уравнения (j). Выберем элемент $b_1 \in A$, отличным от элемента $a_1$, а решение уравнения

$$[b_1 \underbrace{a \ldots a}_{n-3} xa] = [a_1 \ldots a_{n-1}a].$$

обозначим через d, то есть x = d. Так как $a_1, \ldots, a_{n-1}$ – решение уравнения (j), то $[a_1 \ldots a_{n-1}a] = b$, откуда, с учетом выбора d, имеем

$$[b_1 \underbrace{a \ldots a}_{n-3} da] = b.$$

Это означает, что

$$(b_1, \underbrace{a, \ldots, a}_{n-3}, d)$$

решение уравнения (j), отличное от решения $(a_1, \ldots, a_{n-1})$, что противоречит однозначной разрешимости уравнения (j). Следовательно, множество A – одноэлементное.

*Достаточность.* В одноэлементной n-арной полугруппе уравнения (j) и (jj) разрешимы однозначно. ∎

В связи с определением 1.1.3 Поста, по аналогии с теоремой 3.1.1 возникает

**3.1.8. Вопрос.** Будет ли n-арной группой n-арная полугруппа, в которой для некоторого $2 \leq i \leq n - 1$ и любых a, b разрешимо уравнение

$$[x_1 \ldots x_{i-1}ax_{i+1} \ldots x_n] = b?$$

Следующий пример показывает, что ответ на вопрос 3.1.8 – отрицательный.

**3.1.9. Пример.** Пусть A – произвольное множество, содержащее более одного элемента, Определим на нем n-арную операцию ($n \geq 3$)

$$[a_1 \ldots a_{n-1}a_n] = a_n.$$

Так как

$$[[a_1 \ldots a_{n-1}a_n]a_{n+1} \ldots a_{2n-2}a_{2n-1}] =$$

$$= [a_1 \ldots a_i[a_{i+1} \ldots a_{i+n-1}a_{i+n}]a_{i+n+1} \ldots a_{2n-2}a_{2n-1}] = a_{2n-1},$$

то $<A, [\ ]>$ – n-арная полугруппа.

Для любых $a, b \in A$ и любого $i = 1, \ldots, n-1$ в A разрешимо уравнение

$$[x_1 \ldots x_{i-1} a x_{i+1} \ldots x_n] = b.$$

Например,

$$x_1 = \ldots = x_{i-1} = x_{i+1} = \ldots = x_n = b$$

одно из решений. Однако, $<A, [\ ]>$ не является n-арной группой, так как в A неразрешимо уравнение

$$[x_1 \ldots x_{n-1} a] = b,$$

где $a \neq b$.

Интересным представляется также следующий

**3.1.10. Вопрос.** Будет ли n-арной группой n-арная полугруппа $<A, [\ ]>$, в которой для любых $a, b \in A$ разрешимы уравнения

$$[x_1 \ldots x_{i-1} a x_{i+1} \ldots x_n] = b, \quad [y_1 \ldots y_{j-1} a y_{j+1} \ldots y_n] = b, \qquad (***)$$

где $i \neq j$, $i, j = 1, 2, \ldots, n$?

Ясно, что при $i = j$ получается рассмотренный выше случай одного уравнения. Случай $i = n$, $j = 1$ совпадает с теоремой 3.1.1.

Если обе последовательности $x_{i+1} \ldots x_n$ и $y_{j+1} \ldots y_n$ непустые, то множество всех n-арных полугрупп с разрешимостью уравнений (***) содержит n-арную полугруппу из примера 3.1.6 (если положить $x_n = y_n = b$), которая, как показано, не является n-арной группой. Поэтому, по крайней мере, одна из последовательностей должна быть пустой. Пусть для определенности $x_{i+1} \ldots x_n$ – пустая последовательность. Тогда уравнения (***) трансформируются в уравнения

$$[x_1 \ldots x_{n-1} a] = b, \quad [y_1 \ldots y_{j-1} a y_{j+1} \ldots y_n] = b, \qquad (****)$$

На произвольном множестве A содержащем более одного элемента определим еще одну n-арную операцию ($n \geq 3$)

$$(a_1 a_2 \ldots a_n) = a_1.$$

Также, как к в примере 3.1.6, показывается, что $<A, [\ ]>$ – n-арная полугруппа, не являющаяся n-арной группой. Если теперь считать последовательность $y_1 \ldots y_{j-1}$ во втором уравнении (****) непустой, то множество всех n-арных полугрупп с разрешимостью уравнений (****) содержит построенную n-арную полугруппу $<A, [\ ]>$ (если положить $x_1 = y_1 = b$), которая не является n-арной группой. Поэтому $y_1 \ldots y_{j-1}$ – пустая последовательность, а уравнения (****) трансформируются в уравнения (j) и (jj).

Таким образом, мы показали, что ответ на вопрос 3.1.10 будет положительным только в случае $i = n$, $j = 1$.

Следующая теорема обобщает результат Поста (определение 1.1.2) и основную теорему 3.1.1.

**3.1.11. Теорема.** n-Арная полугруппа $<A, [\ ]>$ является n-арной группой тогда и только тогда, когда в ней для некоторых $i, j \in \{1, \ldots, n-1\}$ и любых $a_{i+1}, \ldots, a_n, b_1, \ldots, b_{n-j}, b \in A$ разрешимы уравнения

$$[x_1 \ldots x_i a_{i+1} \ldots a_n] = b, \qquad (1)$$

$$[b_1 \ldots b_{n-j} y_1 \ldots y_j] = b. \qquad (2)$$

*Доказательство. Необходимость.* Пусть $<A, [\ ]>$ – n-арная группа, $d$ – произвольный элемент из A. Определение 1.1.1 Дёрнте гарантирует существование решений $x_1 = d_1$ и $y_j = d_j$ уравнений

$$[x_1 \underbrace{d \ldots d}_{i-1} a_{i+1} \ldots a_n] = b,$$

$$[b_1 \ldots b_{n-j} \underbrace{d \ldots d}_{j-1} y_j] = b.$$

Тогда
$$x_1 = d_1\ x_2 = \ldots = x_i = d$$

– решение уравнения (1), а

$$y_1 = \ldots = y_{j-1} = d,\ y_j = d_j$$

– решение уравнения (2).

*Достаточность.* Пусть теперь $< A, [\ ] >$ – n-арная полугруппа, в которой для некоторых $i, j \in \{1, \ldots, n-1\}$ и любых $a_{i+1}, \ldots, a_n, b_1, \ldots, b_{n-j}, b \in A$ разрешимы уравнения (1) и (2).

Если $c_1, \ldots, c_i$ и $d_1, \ldots, d_j$ – решения уравнений (1) и (2) соответственно, то

$$x_1 = c_1, \ldots, x_i = c_i, x_{i+1} = a_{i+1}, \ldots, x_{n-1} = a_{n-1}$$

– решение уравнения

$$[x_1 \ldots x_{n-1} a_n] = b,$$

для любых $a_n, b \in A$, а

$$z_1 = b_2, \ldots, z_{n-j-1} = b_{n-j}, z_{n-j} = d_1, \ldots, z_{n-1} = d_j$$

– решение уравнения

$$[b_1 z_1 \ldots z_{n-1}] = b$$

для любых $b_1, b \in A$. Поэтому, согласно основной теореме 3.1.1, $< A, [\ ] >$ – n-арная группа. ∎

Аналогично предыдущей доказывается следующая

**3.1.12. Теорема.** n-Арная полугруппа $< A, [\ ] >$ является n-арной группой тогда и только тогда, когда в ней для некоторых $i, j \in \{1, \ldots, n-1\}$ и любых $a, b \in A$ разрешимы уравнения

$$[x_1 \ldots x_i \underbrace{a \ldots a}_{n-i}] = b,$$

$$[\underbrace{a \ldots a}_{n-j} y_1 \ldots y_j] = b.$$

Полагая в теоремах 3.1.11 и 3.1.12 $i = j = 1$, получаем соответственно результат Поста (определение 1.1.2) и результат А.Н. Скибы и В.И. Тютина [49], который будет приведен позже (определение 3.3.21). При i = j = n – 1 теоремы 3.1.11 и 3.1.12 включают основную теорему 3.1.1.

Теоремы 3.1.11 и 3.1.12 являются формальным обобщением основной теоремы 3.1.1, так как она использовалась при их доказательстве.

## §3.2. ДРУГОЙ ПОДХОД К ДОКАЗАТЕЛЬСТВУ ОСНОВНОЙ ТЕОРЕМЫ

В доказательстве основной теоремы разрешимость каждого из уравнений (i) и (ii) определения 1.1.2 Поста является следствием разрешимости уравнений (j) и (jj). А можно ли доказать разрешимость каждого из уравнений (i) и (ii) определения 1.1.2 Поста, используя разрешимость только одного из уравнений (j) или (jj) и не используя разрешимости второго из них? Ниже приводится такое доказательство (предложения 3.2.1 и 3.2.2).

**3.2.1. Предложение.** В n-арной полугруппе $< A, [\ ] >$ уравнение (i) разрешимо для любых $a_2, \ldots, a_n, b \in A$ тогда и только тогда, когда в ней для любых $a, b \in A$ разрешимо уравнение (j).

*Доказательство*. *Необходимость*. Для фиксированных $d_2, \ldots, d_{n-1} \in A$ существует решение $x_1 = d_1$ уравнения

$$[x_1 d_2 \ldots d_{n-1} a] = b.$$

Тогда
$$x_1 = d_1, \ldots, x_{n-1} = d_{n-1}$$

– решение уравнения (j).

*Достаточность*. Обозначим через $(c_1, \ldots, c_{n-1})$ решение уравнения

$$[x_1 \ldots x_{n-1} a_n] = b,$$

а через $(c_1^{(k)}, \ldots, c_{n-1}^{(k)})$ решение уравнения

$$[x_1 \ldots x_{n-1}a_{n-k}] = c_{n-1}^{(k-1)}, \ k = 1, \ldots, n-2,$$

считая при этом $c_{n-1}^{(0)} = c_{n-1}$. Положим также

$$g = [c_1 \ldots c_{n-2}c_1' \ldots c_{n-2}' \ldots c_1^{(n-3)} \ldots c_{n-2}^{(n-3)} c_1^{(n-2)} \ldots c_{n-1}^{(n-2)}].$$

Так как

$$[ga_2 \ldots a_n] = [c_1 \ldots c_{n-2}c_1' \ldots c_{n-2}' \ldots$$

$$\ldots c_1^{(n-3)} \ldots c_{n-2}^{(n-3)} [c_1^{(n-2)} \ldots c_{n-1}^{(n-2)}a_2]a_3 \ldots a_n] =$$

$$= [c_1 \ldots c_{n-2}c_1' \ldots c_{n-2}' \ldots c_1^{(n-3)} \ldots c_{n-2}^{(n-3)} c_{n-1}^{(n-3)}a_3 \ldots a_n] =$$

$$= [c_1 \ldots c_{n-2}c_1' \ldots c_{n-2}' \ldots [c_1^{(n-3)} \ldots c_{n-2}^{(n-3)} c_{n-1}^{(n-3)}a_3]a_4 \ldots a_n] =$$

$$= [c_1 \ldots c_{n-2}c_1' \ldots c_{n-2}' \ldots c_{n-1}^{(n-4)}a_4 \ldots a_n] = \ldots$$

$$\ldots [c_1 \ldots c_{n-2}[c_1' \ldots c_{n-2}' c_{n-1}'a_{n-1}]a_n] = [c_1 \ldots c_{n-2} c_{n-1}a_n] = b,$$

то $g$ – решение уравнения (i). ∎

**3.2.2. Предложение.** В n-арной полугруппе $<A, [\ ]>$ уравнение (ii) разрешимо для любых $a_1, \ldots, a_{n-1}, b \in A$ тогда и только тогда, когда в ней для любых $a, b \in A$ разрешимо уравнение (jj).

*Доказательство*. *Необходимость*. Для фиксированных $d_1, \ldots, d_{n-2} \in A$ существует решение $y_{n-1} = d_{n-1}$ уравнения

$$[ad_1 \ldots d_{n-2}y_{n-1}] = b.$$

Тогда

$$y_1 = d_1, \ldots, y_{n-1} = d_{n-1}$$

– решение уравнения (jj).

*Достаточность*. Обозначим через $(d_1, \ldots, d_{n-1})$ решение уравнения

$$[a_1y_1 \ldots y_{n-1}] = b,$$

а через $(d_1^{(k)}, \ldots, d_{n-1}^{(k)})$ решение уравнения

$$[a_{k+1}y_1 \ldots y_{n-1}] = d_1^{(k-1)}, k = 1, \ldots, n-2,$$

считая при этом $d_1^{(0)} = d_1$. Положим также

$$h = [d_1^{(n-2)} \ldots d_{n-1}^{(n-2)} d_2^{(n-3)} \ldots d_{n-1}^{(n-3)} \ldots d_2' \ldots d_{n-1}' d_2 \ldots d_{n-1}].$$

Так как

$$[a_1 \ldots a_{n-1}h] = [a_1 \ldots a_{n-2}[a_{n-1}d_1^{(n-2)} \ldots d_{n-1}^{(n-2)}]$$

$$d_2^{(n-3)} \ldots d_{n-1}^{(n-3)} \ldots d_2' \ldots d_{n-1}' d_2 \ldots d_{n-1}] =$$

$$= [a_1 \ldots a_{n-2}d_1^{(n-3)} d_2^{(n-3)} \ldots d_{n-1}^{(n-3)} \ldots d_2' \ldots d_{n-1}' d_2 \ldots d_{n-1}] =$$

$$= [a_1 \ldots a_{n-3}[a_{n-2}d_1^{(n-3)} d_2^{(n-3)} \ldots d_{n-1}^{(n-3)}] \ldots d_2' \ldots d_{n-1}' d_2 \ldots d_{n-1}] =$$

$$= [a_1 \ldots a_{n-3}d_1^{(n-4)} \ldots d_2' \ldots d_{n-1}' d_2 \ldots d_{n-1}] = \ldots$$

$$\ldots = [a_1[a_2d_1' d_2' \ldots d_{n-1}']d_2 \ldots d_{n-1}] = [a_1d_1d_2 \ldots d_{n-1}] = b,$$

то h – решение уравнения (ii). ∎

Достаточные утверждения в предложениях 3.2.1 и 3.2.2 можно обобщить.

**3.2.3. Предложение.** Если в n-арной полугруппе $< A, [\ ] >$ для любых $a, b \in A$ разрешимо уравнение (j), то в ней для любого $i \in \{1, \ldots, n-2\}$ и любых $a_{i+1}, \ldots a_n, b \in A$ разрешимо уравнение

$$[x_1 \ldots x_i a_{i+1} \ldots a_n] = b. \tag{v}$$

*Доказательство*. Случай $i = 1$ доказан в предложении 3.2.1.

Также как и в доказательстве предложения 3.2.1 обозначим через $(c_1, \ldots, c_{n-1})$ решение уравнения

$$[x_1 \ldots x_{n-1}a_n] = b,$$

а через $(c_1^{(k)}, \ldots, c_{n-1}^{(k)})$ решение уравнения

$$[x_1 \ldots x_{n-1}a_{n-k}] = c_{n-1}^{(k-1)}, \ k = 1, \ldots, n-i-1,$$

считая при этом $c_{n-1}^{(0)} = c_{n-1}$. Так как

$$[c_1 \ldots c_{i-1}[c_i \ldots c_{n-2}c_1' \ldots c_{n-2}' \ldots$$

$$\ldots c_1^{(n-i-2)} \ldots c_{n-2}^{(n-i-2)} c_1^{(n-i-1)} \ldots c_{n-1}^{(n-i-1)}]a_{i+1} \ldots a_n] =$$

$$= [c_1 \ldots c_{i-1}c_i \ldots c_{n-2}c_1' \ldots c_{n-2}' \ldots$$

$$\ldots c_1^{(n-i-2)} \ldots c_{n-2}^{(n-i-2)} [c_1^{(n-i-1)} \ldots c_{n-1}^{(n-i-1)}a_{i+1}]a_{i+2} \ldots a_n] =$$

$$= [c_1 \ldots c_{n-2}c_1' \ldots c_{n-2}' \ldots$$

$$\ldots [c_1^{(n-i-2)} \ldots c_{n-2}^{(n-i-2)} c_{n-1}^{(n-i-2)}a_{i+2}]a_{i+3} \ldots a_n] =$$

$$= [c_1 \ldots c_{n-2}c_1' \ldots c_{n-2}' \ldots c_1^{(n-i-3)} \ldots c_{n-2}^{(n-i-3)} c_{n-1}^{(n-i-3)}a_{i+3} \ldots a_n] = \ldots$$

$$\ldots = [c_1 \ldots c_{n-2}[c_1' \ldots c_{n-1}'a_{n-1}]a_n] = [c_1 \ldots c_{n-2} c_{n-1}a_n] = b,$$

то

$$x_1 = c_1, \ldots, x_{i-1} = c_{i-1},$$

$$x_i = [c_i \ldots c_{n-2}c_1' \ldots c_{n-2}' \ldots c_1^{(n-i-2)} \ldots c_{n-2}^{(n-i-2)} c_1^{(n-i-1)} \ldots c_{n-1}^{(n-i-1)}]$$

– решение уравнения (ν). ∎

**3.2.4. Предложение.** Если в n-арной полугруппе $< A, [\ ] >$ для любых $a, b \in A$ разрешимо уравнение (jj), то в ней для любого $j \in \{1, \ldots, n-2\}$ и любых $a_1, \ldots a_{n-j}, b \in A$ разрешимо уравнение

$$[a_1 \ldots a_{n-j}y_1 \ldots y_j] = b. \qquad (\nu\nu)$$

*Доказательство*. Случай $j = 1$ доказан в предложении 3.2.2.

Также как и в доказательстве предложения 3.2.2 обозначим через $(d_1, \ldots, d_{n-1})$ решение уравнения

$$[a_1y_1 \ldots y_{n-1}] = b,$$

а через $(d_1^{(k)}, \ldots, d_{n-1}^{(k)})$ решение уравнения

$$[a_{k+1} y_1 \ldots y_{n-1}] = d_1^{(k-1)}, \ k = 1, \ldots, n-j-1,$$

считая при этом $d_1^{(0)} = d_1$. Так как

$$[a_1 \ldots a_{n-j}[d_1^{(n-j-1)} \ldots d_{n-1}^{(n-j-1)} d_2^{(n-j-2)} \ldots d_{n-1}^{(n-j-2)} \ldots$$

$$\ldots d_2' \ldots d_{n-1}' d_2 \ldots d_{n-j}] d_{n-j+1} \ldots d_{n-1}] =$$

$$= [a_1 \ldots a_{n-j-1}[a_{n-j} d_1^{(n-j-1)} \ldots d_{n-1}^{(n-j-1)}] d_2^{(n-j-2)} \ldots d_{n-1}^{(n-j-2)} \ldots$$

$$\ldots d_2' \ldots d_{n-1}' d_2 \ldots d_{n-j}] d_{n-j+1} \ldots d_{n-1}] =$$

$$= [a_1 \ldots a_{n-j-2}[a_{n-j-1} d_1^{(n-j-2)} \ldots d_{n-1}^{(n-j-2)}] \ldots d_2' \ldots d_{n-1}' d_2 \ldots d_{n-1}] =$$

$$= [a_1 \ldots a_{n-j-2} d_1^{(n-j-3)} d_2^{(n-j-3)} \ldots d_{n-1}^{(n-j-3)} \ldots d_2' \ldots d_{n-1}' d_2 \ldots d_{n-1}] = \ldots$$

$$\ldots = [a_1[a_2 d_1' \ldots d_{n-1}'] d_2 \ldots d_{n-1}] = [a_1 d_1 d_2 \ldots d_{n-1}] = b,$$

то

$$y_1 = [d_1^{(n-j-1)} \ldots d_{n-1}^{(n-j-1)} d_2^{(n-j-2)} \ldots d_{n-1}^{(n-j-2)} \ldots$$

$$\ldots d_2' \ldots d_{n-1}' d_2 \ldots d_{n-j}], \ y_2 = d_{n-j+1}, \ldots, y_j = d_{n-1}$$

– решение уравнения (νν). ∎

**3.2.5. Пример.** Полагая в предложении 3.2.3 $i = 2$ и $i = n - 2$, а в предложении 3.2.4 $j = 2$ и $j = n - 2$, выпишем в явном виде решения соответствующих уравнений.

Решение уравнения

$$[x_1 x_2 a_3 \ldots a_n] = b$$

имеет вид

$$x_1 = c_1, \ [c_2 \ldots c_{n-2} c_1' \ldots c_{n-2}' \ldots c_1^{(n-4)} \ldots c_{n-2}^{(n-4)} c_1^{(n-3)} \ldots c_{n-1}^{(n-3)}],$$

а решение уравнения

$$[x_1 \ldots x_{n-2} a_{n-1} a_n] = b$$

имеет вид

$$x_1 = c_1, \ldots, x_{n-3} = c_{n-3}, \ x_{n-2} = [c_{n-2} c_1' \ldots c_{n-1}'].$$

Решение уравнения
$$[a_1 \ldots a_{n-2} y_1 y_2] = b$$
имеет вид
$$y_1 = [\,d_1^{(n-3)} \ldots d_{n-1}^{(n-3)} d_2^{(n-4)} \ldots d_{n-1}^{(n-4)} \ldots d_2' \ldots d_{n-1}' d_2 \ldots d_{n-2}], \; y_2 = d_{n-1},$$
а решение уравнения
$$[a_1 a_2 y_1 \ldots y_{n-2}] = b$$
имеет вид
$$y_1 = [\,d_1' \ldots d_{n-1}' d_2], \; y_2 = d_3, \ldots, y_{n-2} = d_{n-1}.$$

**3.2.6. Замечание.** При $i \neq 1$ и $j \neq 1$ решения уравнений ($\nu$) и ($\nu\nu$), найденные при доказательстве предложений 3.2.3 и 3.2.4, в том числе и приведенные в примере 3.2.5, не являются единственными.

Можно получать различные решения уравнения ($\nu$), разбивая последовательность
$$\alpha = c_1 \ldots c_{n-2} c_1' \ldots c_{n-2}' \ldots c_1^{(n-i-2)} \ldots c_{n-2}^{(n-i-2)} c_1^{(n-i-1)} \ldots c_{n-1}^{(n-i-1)}$$
на $i$ подпоследовательностей $\alpha_1, \ldots, \alpha_i$ так, что $\alpha = \alpha_1 \ldots \alpha_i$ и либо $\alpha_s \in A$, либо $[\alpha_s] \in A$, $s = 1, \ldots, i$.

Аналогично различные решения уравнения ($\nu\nu$) можно получать разбивая последовательность
$$\beta = d_1^{(n-j-1)} \ldots d_{n-1}^{(n-j-1)} d_2^{(n-j-2)} \ldots d_{n-1}^{(n-j-2)} \ldots d_2' \ldots d_{n-1}' d_2 \ldots d_{n-1}$$
на $j$ подпоследовательностей $\beta_1, \ldots, \beta_j$ так, что $\beta = \beta_1 \ldots \beta_j$ и либо $\beta_s \in A$, либо $[\beta_s] \in A$, $s = 1, \ldots, j$.

**3.2.7. Пример.** Методом, описанном в замечании 3.2.6, из решений, приведенных в примере 3.2.5, можно получить решение
$$x_1 = [c_1 \ldots c_{n-2} c_1' \ldots c_{n-2}' \ldots c_1^{(n-3)} \ldots c_{n-2}^{(n-3)}], \; x_2 = c_{n-1}^{(n-3)}$$
уравнения
$$[x_1 x_2 a_3 \ldots a_n] = b$$
и решение
$$y_1 = d_1^{(n-3)}, \; y_2 = [\,d_2^{(n-3)} \ldots d_{n-1}^{(n-3)} \ldots d_2' \ldots d_{n-1}' d_2 \ldots d_{n-1}]$$

уравнения
$$[a_1 \ldots a_{n-2}y_1y_2] = b.$$

**3.2.8. Пример.** Запишем все n – 2 решения уравнения
$$[x_1 \ldots x_{n-2}a_{n-1}a_n] = b,$$
которые можно получить методом, описанном в замечании 3.2.6:

$$x_1 = c_1, \ldots, x_{n-3} = c_{n-3}, x_{n-2} = [c_{n-2}c'_1 \ldots c'_{n-1}];$$

$$x_1 = c_1, \ldots, x_{n-4} = c_{n-4}, x_{n-3} = [c_{n-3}c_{n-2}c'_1 \ldots c'_{n-2}], x_{n-2} = c'_{n-1};$$

$$\ldots\ldots\ldots\ldots\ldots\ldots\ldots\ldots\ldots\ldots\ldots\ldots\ldots\ldots\ldots\ldots\ldots\ldots\ldots\ldots\ldots\ldots\ldots\ldots$$

$$x_1 = c_1, x_2 = [c_2 \ldots c_{n-2}c'_1 c'_2 c'_3], x_3 = c'_4, \ldots, x_{n-2} = c'_{n-1};$$

$$x_1 = [c_1 \ldots c_{n-2}c'_1 c'_2], x_2 = c'_3, \ldots, x_{n-2} = c'_{n-1}.$$

Запишем все n – 2 решения уравнения
$$[a_1a_2y_1 \ldots y_{n-2}] = b,$$
которые можно получить методом, описанном в замечании 3.2.6:

$$y_1 = [d'_1 \ldots d'_{n-1}d_2], y_2 = d_3, \ldots, y_{n-2} = d_{n-1};$$

$$y_1 = d'_1, y_2 = [d'_2 \ldots d'_{n-1}d_2d_3], y_3 = d_4, \ldots, y_{n-2} = d_{n-1};$$

$$\ldots\ldots\ldots\ldots\ldots\ldots\ldots\ldots\ldots\ldots\ldots\ldots\ldots\ldots\ldots\ldots\ldots\ldots\ldots\ldots\ldots\ldots\ldots\ldots$$

$$y_1 = d'_1, \ldots, y_{n-4} = d'_{n-4}, y_{n-3} = [d'_{n-3} d'_{n-2} d'_{n-1}d_2 \ldots d_{n-2}], y_{n-2} = d_{n-1};$$

$$y_1 = d'_1, \ldots, y_{n-3} = d'_{n-3}, y_{n-2} = [d'_{n-2} d'_{n-1}d_2 \ldots d_{n-1}].$$

Ясно, что если в n-арной полугруппе $< A, [\ ] >$ для некоторого $i \in \{1, \ldots, n-2\}$ и любых $a_{i+1}, \ldots, a_n, b \in A$ разрешимо уравнение (v), то в ней для любых $a, b \in A$ разрешимо уравнение (j). Поэтому предложение 3.2.3 позволяет сформулировать следующую теорему.

**3.2.9. Теорема.** Если в n-арной полугруппе $< A, [\ ] >$ для некоторого $i \in \{1, \ldots, n-1\}$ и любых $a_{i+1}, \ldots, a_n, b \in A$ разрешимо уравнение (v), то в ней разрешимо каждое из уравнений

$$[x_1a_2 \ldots a_n] = b, [x_1x_2a_3 \ldots a_n] = b, \ldots, [x_1 \ldots x_{n-1}a_n] = b \quad (vvv)$$

для всех $a_2, \ldots, a_n, b \in A$.

Ясно, что если в n-арной полугруппе < A, [ ] > для некоторого j ∈ {1, …, n – 2} и любых $a_1, …, a_{n-j}, b \in A$ разрешимо уравнение (vv), то в ней для любых a, b ∈ A разрешимо уравнение (jj). Поэтому предложение 3.2.4 позволяет сформулировать следующую теорему.

**3.2.10. Теорема.** Если в n-арной полугруппе < A, [ ] > для некоторого j ∈ {1, …, n – 1} и любых $a_1, …, a_{n-j}, b \in A$ разрешимо уравнение (vv), то в ней разрешимо каждое из уравнений

$$[a_1 … a_{n-1} y_1] = b, [a_1 … a_{n-2} y_1 y_2] = b, …, [a_1 y_1 … y_{n-1}] = b \quad (vvvv)$$

для всех $a_1, …, a_{n-1}, b \in A$.

Условия теорем 3.2.9 и 3.2.10 можно ослабить, сохранив их утверждения.

**3.2.11. Теорема.** Если в n-арной полугруппе < A, [ ] > для некоторого i ∈ {1, …, n – 1} и любых a, b ∈ A разрешимо уравнение

$$[x_1 … x_i \underbrace{a … a}_{n-i}] = b,$$

то в ней разрешимо каждое из уравнений (vvv) для любых $a_2, …, a_n, b \in A$.

**3.2.12. Теорема.** Если в n-арной полугруппе < A, [ ] > для некоторого j ∈ {1, …, n – 1} и любых a, b ∈ A разрешимо уравнение

$$[\underbrace{a … a}_{n-j} y_1 … y_j] = b,$$

то в ней разрешимо каждое из уравнений (vvvv) для любых $a_1, …, a_{n-1}, b \in A$.

**3.2.13. Лемма.** Пусть в n-арной полугруппе < A, [ ] > для некоторого i ∈ {2, …, n – 1} и любых $a_{i+1}, …, a_n, b \in A$ разрешимо уравнение (v). Для любых $\alpha_2, …, \alpha_i \in A$ обозначим че-рез $X(\alpha_2, \ldots, \alpha_i)$ множество всех решений уравнения (v) вида $(\alpha, \alpha_2, \ldots, \alpha_i)$, $\alpha \in A$. Тогда:

1) $X(\alpha_2, \ldots, \alpha_i) \neq \varnothing$;

2) $X(\alpha_2, \ldots, \alpha_i) \cap X(\gamma_2, \ldots, \gamma_i) = \varnothing$ для несовпадающих последовательностей $\alpha_2 \ldots \alpha_i$ и $\gamma_2 \ldots \gamma_i$;

3) множество всех решений уравнения (v) совпадает с объединением $\bigcup\limits_{\alpha_2, \ldots, \alpha_i \in A} X(\alpha_2, \ldots, \alpha_i)$.

*Доказательство.* 1) По теореме 3.2.9 в $< A, [\ ] >$ разрешимо уравнение

$$[x\alpha_2 \ldots \alpha_i a_{i+1} \ldots a_n] = b,$$

то есть существует такой $\alpha \in A$, что

$$[\alpha\alpha_2 \ldots \alpha_i a_{i+1} \ldots a_n] = b.$$

Следовательно, $(\alpha, \alpha_2, \ldots, \alpha_i)$ – решение уравнения (v), а значит $X(\alpha_2, \ldots, \alpha_i) \neq \varnothing$.

Утверждения 2) и 3) очевидны. ∎

Аналогично лемме 3.2.13 доказывается следующая

**3.2.14. Лемма.** Пусть в n-арной полугруппе $< A, [\ ] >$ для некоторого $j \in \{2, \ldots, n-1\}$ и любых $a_1, \ldots, a_{n-j}, b \in A$ разрешимо уравнение (vv). Для любых $\beta_1, \ldots, \beta_{j-1} \in A$ обозначим через $Y(\beta_1, \ldots, \beta_{j-1})$ множество всех решений уравнения (vv) вида $(\beta_1, \ldots, \beta_{j-1}, \beta)$, $\beta \in A$. Тогда:

1) $Y(\beta_1, \ldots, \beta_{j-1}) \neq \varnothing$;

2) $Y(\beta_1, \ldots, \beta_{j-1}) \cap Y(\gamma_1, \ldots, \gamma_{j-1}) = \varnothing$ для несовпадающих последовательностей $\beta_1 \ldots \beta_{j-1}$ и $\gamma_1 \ldots \gamma_{j-1}$;

3) множество всех решений уравнения (vv) совпадает с объединением $\bigcup\limits_{\beta_1, \ldots, \beta_{j-1} \in A} Y(\beta_1, \ldots, \beta_{j-1})$.

**3.2.15. Предложение.** Пусть в конечной n-арной полугруппе $< A, [\ ] >$ порядка r для некоторого $i \in \{2, \ldots, n-1\}$ и любых $a_{i+1}, \ldots, a_n, b \in A$ разрешимо уравнение (ν). Тогда:

1) число различных решений уравнения (ν) ограничено снизу числом $r^{i-1}$;

2) если для любых $c_2, \ldots, c_n \in A$ из

$$[xc_2 \ldots c_n] = [yc_2 \ldots c_n]$$

всегда следует $x = y$, то число различных решений уравнения (ν) равно $r^{i-1}$.

***Доказательство.*** 1) Так как $|A| = r$, то число различных последовательностей вида $\alpha_2 \ldots \alpha_i$, где все $\alpha_2, \ldots, \alpha_i$ пробегают множество A, равно $r^{i-1}$. Тогда, ввиду леммы 3.2.13, множество всех решений уравнения (ν) разбивается на $r^{i-1}$ непересекающихся подмножеств, каждое из которых не пусто. Поэтому число различных решений уравнения (ν) не меньше $r^{i-1}$.

2) Предположим, что $(\alpha, \alpha_2, \ldots, \alpha_i)$ и $(\delta, \alpha_2, \ldots, \alpha_i)$ – два различных решения уравнения (ν) из множества $X(\alpha_2, \ldots, \alpha_i)$, то есть

$$[\alpha\alpha_2 \ldots \alpha_i a_{i+1} \ldots a_n] = [\delta\alpha_2 \ldots \alpha_i a_{i+1} \ldots a_n] = b.$$

Из последнего равенства, учитывая условие утверждения 2), получаем $\alpha = \delta$. Следовательно, для любых $\alpha_2, \ldots, \alpha_i \in A$ все множества $X(\alpha_2, \ldots, \alpha_i)$ являются одноэлементными. Применяя лемму 3.2.13, видим, что число различных решений уравнения (ν) равно $r^{i-1}$. ∎

Аналогично предложению 3.2.15 при помощи леммы 3.2.14 доказывается следующее

**3.2.16. Предложение.** Пусть в конечной n-арной полугруппе $< A, [\,] >$ порядка r для некоторого $j \in \{2, \ldots, n-1\}$ и любых $a_1, \ldots, a_{n-j}, b \in A$ разрешимо уравнение (νν). Тогда:

1) число различных решений уравнения (νν) ограничено снизу числом $r^{j-1}$;

2) если для любых $c_1, \ldots, c_{n-1} \in A$ из

$$[c_1 \ldots c_{n-1}x] = [c_1 \ldots c_{n-1}y]$$

всегда следут $x = y$, то число различных решений уравнения (νν) равно $r^{j-1}$.

**3.2.17. Замечание.** Для несовпадающих последовательностей $\alpha_2 \ldots \alpha_i$ и $\gamma_2 \ldots \gamma_i$ и некоторых

$$(\alpha, \alpha_2, \ldots, \alpha_i) \in X(\alpha_2, \ldots, \alpha_i), (\gamma, \gamma_2, \ldots, \gamma_i) \in X(\gamma_2, \ldots, \gamma_i)$$

возможно равенство $\alpha = \gamma$. Аналогично для несовпадающих последовательностей $\beta_1 \ldots \beta_{j-1}$ и $\gamma_1 \ldots \gamma_{j-1}$ и некоторых

$$(\beta_1, \ldots, \beta_{j-1}, \beta) \in Y(\beta_1, \ldots, \beta_{j-1}), (\gamma_1 \ldots \gamma_{j-1}, \gamma) \in Y(\gamma_1 \ldots \gamma_{j-1})$$

возможно равенство $\beta = \gamma$.

Следующее предложение доказывается аналогично утверждениям 2) предложений 3.2.15 и 3.2.16.

**3.2.18. Предложение.** Для n-арной группы $< A, [\ ] >$ мощности множеств всех решений каждого из уравнений

$$[x_1 x_2 a_3 \ldots a_n] = b, \quad [a_1 \ldots a_{n-2} y_1 y_2] = b$$

совпадают с мощностью множества $A$.

Из утверждений 2) предложений 3.2.15 и 3.2.16 вытекает

**3.2.19. Следствие.** В конечной n-арной группе $< A, [\ ] >$ число различных решений уравнений (ν) и (νν) равно соответственно $|A|^{i-1}$ и $|A|^{j-1}$.

**3.2.20. Замечание.** Каждое из уравнений

$$[xa_2 \ldots a_n] = b \text{ и } [a_1 \ldots a_{n-1} y] = b$$

имеет в n-арной группе единственное решение. Поэтому следствие 3.2.19 формально включает и случаи $i = 1$ и $j = 1$.

## §3.3. ОПРЕДЕЛЕНИЯ, СВЯЗАННЫЕ С РАЗРЕШИМОСТЬЮ ДВУХ УРАВНЕНИЙ

Теорема 3.1.1 позволяет дать следующее

**3.3.1. Определение.** n-Арная полугруппа $< A, [\ ] >$ называется n-арной группой, если для любых $a, b \in A$ в $A$ разрешимы уравнения

$$[x_1 \ldots x_{n-1}a] = b, \quad [ay_1 \ldots y_{n-1}] = b.$$

Предложение 3.1.2 позволяет дать следующее

**3.3.2. Определение.** n-Арная полугруппа $< A, [\ ] >$ называется n-арной группой, если для любых $a, b \in A$ в $A$ разрешимы уравнения

$$[x\underbrace{u \ldots u}_{n-2}a] = b, \quad [a\underbrace{v \ldots v}_{n-2}y] = b.$$

Теорема 3.1.1 позволяет дать большое число определений n-арной группы, как уже известных, так и новых, некоторые из которых приведены ниже.

**3.3.3. Определение.** n-Арная полугруппа $< A, [\ ] >$ называется n-арной группой, если для любых $a, b \in A$ существуют последовательности $\alpha(a, b)$ и $\beta(a, b)$ длины $n - 2$ каждая, составленные из элементов множества $A$, такие, что в $A$ разрешимы уравнения

$$[x\alpha(a, b)a] = b, \quad [a\beta(a, b)y] = b.$$

**3.3.4. Определение.** n-Арная полугруппа $< A, [\ ] >$ называется n-арной группой, если для любого $a, \in A$ существуют последовательности $c_1(a) \ldots c_{n-2}(a)$ и $d_1(a) \ldots d_{n-2}(a)$, составленные из элементов множества $A$, такие, что для любого $b \in A$ в $A$ разрешимы уравнения

$$[xc_1(a) \ldots c_{n-2}(a)a] = b, \quad [ad_1(a) \ldots d_{n-2}(a)y] = b.$$

**3.3.5. Определение.** n-Арная полугруппа $< A, [\,] >$ называется n-арной группой, если для любого $b \in A$ существуют последовательности $c_1(b) \ldots c_{n-2}(b)$ и $d_1(b) \ldots d_{n-2}(b)$, составленные из элементов множества A, такие, что для любого $a \in A$ в A разрешимы уравнения

$$[xc_1(b) \ldots c_{n-2}(b)a] = b, \quad [ad_1(b) \ldots d_{n-2}(b)y] = b.$$

**3.3.6. Определение.** n-Арная полугруппа $< A, [\,] >$ называется n-арной группой, если для любых $a, b \in A$ существуют $c(a, b), d(a, b) \in A$ такие, что в A разрешимы уравнения

$$[x\underbrace{c(a,b) \ldots c(a,b)}_{n-2}a] = b, \quad [a\underbrace{d(a,b) \ldots d(a,b)}_{n-2}y] = b.$$

**3.3.7. Определение.** n-Арная полугруппа $< A, [\,] >$ называется n-арной группой, если для любого $a \in A$ существуют $c(a), d(a) \in A$ такие, что для любого $b \in A$ в A разрешимы уравнения

$$[x\underbrace{c(a) \ldots c(a)}_{n-2}a] = b, \quad [a\underbrace{d(a) \ldots d(a)}_{n-2}y] = b.$$

**3.3.8. Определение.** n-Арная полугруппа $< A, [\,] >$ называется n-арной группой, если для любого $b \in A$ существуют $c(b), d(b) \in A$ такие, что для любого $a \in A$ в A разрешимы уравнения

$$[x\underbrace{c(b) \ldots c(b)}_{n-2}a] = b, \quad [a\underbrace{d(b) \ldots d(b)}_{n-2}y] = b.$$

**3.3.9. Определение.** n-Арная полугруппа $< A, [\,] >$ называется n-арной группой, если существуют $c_1, \ldots, c_{n-2}, d_1, \ldots, d_{n-2} \in A$ такие, что для любых $a, b \in A$ в A разрешимы уравнения

$$[xc_1 \ldots c_{n-2}a] = b, \quad [ad_1 \ldots d_{n-2}y] = b.$$

**3.3.10. Определение.** n-Арная полугруппа $<A, [\,]>$ называется n-арной группой, если существуют $c, d \in A$ такие, что для любых $a, b \in A$ в $A$ разрешимы уравнения

$$[x\underbrace{c \ldots c}_{n-2}a] = b, \quad [a\underbrace{d \ldots d}_{n-2}y] = b.$$

**3.3.11. Определение.** n-Арная полугруппа $<A, [\,]>$ называется n-арной группой, если для любых $a, b \in A$ существует последовательность $\alpha(a, b)$ длины $n - 2$, составленная из элементов множества $A$, такая, что в $A$ разрешимы уравнения

$$[x\alpha(a, b)a] = b, \quad [a\alpha(a, b)y] = b.$$

**3.3.12. Определение.** n-Арная полугруппа $<A, [\,]>$ называется n-арной группой, если для любого $a \in A$ существует последовательность $c_1(a) \ldots c_{n-2}(a)$, составленная из элементов множества $A$, такая, что для любого $b \in A$ в $A$ разрешимы уравнения

$$[xc_1(a) \ldots c_{n-2}(a)a] = b, \quad [ac_1(a) \ldots c_{n-2}(a)y] = b.$$

**3.3.13. Определение.** n-Арная полугруппа $<A, [\,]>$ называется n-арной группой, если для любого $b \in A$ существует последовательность $c_1(b) \ldots c_{n-2}(b)$, составленная из элементов множества $A$, такая, что для любого $a \in A$ в $A$ разрешимы уравнения

$$[xc_1(b) \ldots c_{n-2}(b)a] = b, \quad [ac_1(b) \ldots c_{n-2}(b)y] = b.$$

**3.3.14. Определение.** n-Арная полугруппа $<A, [\,]>$ называется n-арной группой, если для любых $a, b \in A$ существует $c(a, b) \in A$ такой, что в $A$ разрешимы уравнения

$$[x\underbrace{c(a,b)\ldots c(a,b)}_{n-2}a] = b, \quad [a\underbrace{c(a,b)\ldots c(a,b)}_{n-2}y] = b.$$

**3.3.15. Определение.** n-Арная полугруппа $< A, [\ ] >$ называется n-арной группой, если для любого $a \in A$ существует $c(a) \in A$ такой, что для любого $b \in A$ в $A$ разрешимы уравнения

$$[x\underbrace{c(a)\ldots c(a)}_{n-2}a] = b, \quad [a\underbrace{c(a)\ldots c(a)}_{n-2}y] = b.$$

**3.3.16. Определение.** n-Арная полугруппа $< A, [\ ] >$ называется n-арной группой, если для любого $b \in A$ существует $c(b) \in A$ такой, что для любого $a \in A$ в $A$ разрешимы уравнения

$$[x\underbrace{c(b)\ldots c(b)}_{n-2}a] = b, \quad [a\underbrace{c(b)\ldots c(b)}_{n-2}y] = b.$$

**3.3.17. Определение.** n-Арная полугруппа $< A, [\ ] >$ называется n-арной группой, если существуют $c_1, \ldots, c_{n-2} \in A$ такие, что для любых $a, b \in A$ в $A$ разрешимы уравнения

$$[xc_1 \ldots c_{n-2}a] = b, \quad [ac_1 \ldots c_{n-2}y] = b.$$

**3.3.18. Определение** [50]. n-Арная полугруппа $< A, [\ ] >$ называется n-арной группой, если существует $c \in A$ такой, что для любых $a, b \in A$ в $A$ разрешимы уравнения

$$[x\underbrace{c \ldots c}_{n-2}a] = b, \quad [a\underbrace{c \ldots c}_{n-2}y] = b.$$

**3.3.19. Определение** [51]. n-Арная полугруппа $< A, [\ ] >$ называется n-арной группой, если для фиксированных $i, j \in \{1, \ldots, n-1\}$ и любых $a, b \in A$ и в $A$ разрешимы уравнения

$$[x\underbrace{b\ldots b}_{n-1-i}\underbrace{a\ldots a}_{i}] = b, \quad [\underbrace{a\ldots a}_{j}\underbrace{b\ldots b}_{n-1-j}y] = b.$$

**3.3.20. Определение** [51]. n-Арная полугруппа $< A, [\,] >$ называется n-арной группой, если для любых $a, b \in A$ в $A$ разрешимы уравнения

$$[x\underbrace{b\ldots b}_{n-2}a] = b, \quad [a\underbrace{b\ldots b}_{n-2}y] = b.$$

**3.3.21. Определение** [49, Скиба А.Н., Тютин В.И.]. n-Арная полугруппа $< A, [\,] >$ называется n-арной группой, если для любых $a, b \in A$ в $A$ разрешимы уравнения

$$[x\underbrace{a\ldots a}_{n-1}] = b, \quad [\underbrace{a\ldots a}_{n-1}y] = b.$$

Большое число новых определений можно получить, комбинируя уже имеющиеся определения. Например, с помощью определений 3.3.20 и 3.3..21 можно получить следующие два определения.

**3.3.22. Определение.** n-Арная полугруппа $< A, [\,] >$ называется n-арной группой, если для любых $a, b \in A$ в $A$ разрешимы уравнения

$$[x\underbrace{b\ldots b}_{n-2}a] = b, \quad [\underbrace{a\ldots a}_{n-1}y] = b.$$

**3.3.23. Определение.** n-Арная полугруппа $< A, [\,] >$ называется n-арной группой, если для любых $a, b \in A$ в $A$ разрешимы уравнения

$$[x\underbrace{a\ldots a}_{n-1}] = b, \quad [a\underbrace{b\ldots b}_{n-2}y] = b.$$

Ниже в определениях этого параграфа будем считать $k$ и $m$ фиксированными из множества $\{1, \ldots, n-1\}$.

**3.3.24. Определение.** n-Арная полугруппа $< A, [\ ] >$ называется n-арной группой, если для любых

$$a_1, \ldots, a_k, b_1, \ldots, b_m, b \in A$$

существуют последовательности $c_1 \ldots c_{n-k-1}$, $d_1 \ldots d_{n-m-1}$, составленные из элементов множества $A$, такие, что в $A$ разрешимы уравнения

$$[xc_1 \ldots c_{n-k-1}a_1 \ldots a_k] = b, \quad [b_1 \ldots b_m d_1 \ldots d_{n-m-1}y] = b.$$

Взяв за основу определение 3.3.24, можно получить серию определений, аналогичных определениям 3.3.4 – 3.3.21. Приведем только несколько определений этой серии.

**3.3.25. Определение.** n-Арная полугруппа $< A, [\ ] >$ называется n-арной группой, если существуют $c_1, \ldots, c_{n-k-1}$, $d_1, \ldots, d_{n-m-1} \in A$ такие, что для любых $a_1, \ldots, a_k, b_1, \ldots, b_m, b \in A$ в $A$ разрешимы уравнения

$$[xc_1 \ldots c_{n-k-1}a_1 \ldots a_k] = b, \quad [b_1 \ldots b_m d_1 \ldots d_{n-m-1}y] = b.$$

**3.3.26. Определение.** n-Арная полугруппа $< A, [\ ] >$ называется n-арной группой, если существуют такие элементы $c_1, \ldots, c_{n-k-1} \in A$, что для любых $a_1, \ldots, a_k, b \in A$ в $A$ разрешимы уравнения

$$[xc_1 \ldots c_{n-k-1}a_1 \ldots a_k] = b, \quad [a_1 \ldots a_k c_1 \ldots c_{n-k-1}y] = b.$$

**3.3.27. Определение.** n-Арная полугруппа $< A, [\ ] >$ называется n-арной группой, если существует такой $c \in A$, что для любых $a, b \in A$ в ней разрешимы уравнения

$$[x \underbrace{c \ldots c}_{n-k-1} \underbrace{a \ldots a}_{k}] = b, \quad [\underbrace{a \ldots a}_{m} \underbrace{c \ldots c}_{n-m-1} y] = b.$$

**3.3.28. Определение.** n-Арная полугруппа $< A, [\ ] >$ называется n-арной группой, если существует такой $c \in A$, что для любых $a, b \in A$ в ней разрешимы уравнения

$$[\underbrace{xc\ldots c}_{n-k-1}\underbrace{a\ldots a}_{k}] = b, \quad [\underbrace{a\ldots a}_{k}\underbrace{c\ldots c}_{n-k-1}y] = b.$$

Теореме 3.1.11 соответствует

**3.3.29. Определение.** n-Арная полугруппа $< A, [\ ] >$ называется n-арной группой, если в ней для некоторых $i, j \in \{1, \ldots, n-1\}$ и любых $a_{i+1}, \ldots, a_n, b_1, \ldots, b_{n-j}, b \in A$ разрешимы уравнения

$$[x_1 \ldots x_i a_{i+1} \ldots a_n] = b, \quad [b_1 \ldots b_{n-j} y_1 \ldots y_j] = b.$$

**3.3.30. Определение.** n-Арная полугруппа $< A, [\ ] >$ называется n-арной группой, если в ней для некоторого $i \in \{1, \ldots, n-1\}$ и любых $a_{i+1}, \ldots, a_n, b_1, \ldots, b_{n-i}, b \in A$ разрешимы уравнения

$$[x_1 \ldots x_i a_{i+1} \ldots a_n] = b, \quad [b_1 \ldots b_{n-i} y_1 \ldots y_i] = b.$$

Теореме 3.1.12 соответствует

**3.3.31. Определение.** n-Арная полугруппа $< A, [\ ] >$ называется n-арной группой, если в ней для некоторых $i, j \in \{1, \ldots, n-1\}$ и любых $a, b \in A$ разрешимы уравнения

$$[x_1 \ldots x_i \underbrace{a \ldots a}_{n-i}] = b, \quad [\underbrace{a \ldots a}_{n-j} y_1 \ldots y_j] = b.$$

**3.3.32. Определение.** n-Арная полугруппа $< A, [\ ] >$ называется n-арной группой, если в ней для некоторого $i \in \{1, \ldots, n-1\}$ и любых $a, b \in A$ разрешимы уравнения

$$[x_1 \ldots x_i \underbrace{a \ldots a}_{n-i}] = b, \quad [\underbrace{a \ldots a}_{n-i} y_1 \ldots y_i] = b.$$

Комбинируя уже имеющиеся определения, можно получать новые определения, например, такие

**3.3.33. Определение.** n-Арная полугруппа $< A, [\ ] >$ называется n-арной группой, если в ней для любых $a, b \in A$ разрешимы уравнения

$$[x\underbrace{a \ldots a}_{n-1}] = b, \quad [ay_1 \ldots y_{n-1}] = b.$$

**3.3.34. Определение.** n-Арная полугруппа $< A, [\ ] >$ называется n-арной группой, если в ней для любых $a, b \in A$ разрешимы уравнения

$$[x_1 \ldots x_{n-1}a] = b, \quad [\underbrace{a \ldots a}_{n-1}y] = b.$$

## §3.4. ОПРЕДЕЛЕНИЯ, СВЯЗАННЫЕ С РАЗРЕШИМОСТЬЮ ОДНОГО УРАВНЕНИЯ

Всюду в этом параграфе $n \geq 3$.

**3.4.1. Теорема.** n-Арная полугруппа $< A, [\ ] >$ является n-арной группой тогда и только тогда, когда в ней для некоторых $k, i \in \{1, \ldots, n-2\}$, удовлетворяющих неравенству $k + i \leq n - 1$, и любых $a_1, \ldots, a_k, a_{k+i+1}, \ldots, a_n, b \in A$ разрешимо уравнение

$$[a_1 \ldots a_k x_1 \ldots x_i a_{k+i+1} \ldots a_n] = b. \qquad (1)$$

*Доказательство. Необходимость.* Пусть $< A, [\ ] >$ – n-арная группа, $d$ – фиксированный элемент из $A$. Определение 1.1.3 Поста гарантирует существование решения $x_1 = d_1$ уравнения

$$[a_1 \ldots a_k x_1 \underbrace{d \ldots d}_{i-1} a_{k+i+1} \ldots a_n] = b.$$

Тогда
$$x_1 = d_1, x_2 = \ldots = x_i = d$$

– решение уравнения (1).

*Достаточность.* Пусть теперь $< A, [\ ] >$ – n-арная полугруппа, в которой для некоторых $k, i \in \{1, \ldots, n-2\}$ таких, что $k + i \leq n - 1$, и любых $a_1, \ldots, a_k, a_{k+i+1}, \ldots, a_n, b \in A$ разрешимо уравнение (1).

Если $x_1 = d_1, \ldots, x_i = d_i$ – решение уравнения (1), то

$$u_1 = a_1, \ldots, u_k = a_k, u_{k+1} = d_1, \ldots,$$

$$u_{k+i} = d_i, u_{k+i+1} = a_{k+i+1} \ldots, u_{n-1} = a_{n-1}$$

– решение уравнения

$$[u_1 \ldots u_{n-1}a_n] = b,$$

для любых $a_n, b \in A$, а

$$v_1 = a_2, \ldots, v_{k-1} = a_k, v_k = d_1, \ldots,$$

$$v_{k+i-1} = d_i, v_{k+i} = a_{k+i+1}, \ldots, v_{n-1} = a_n$$

– решение уравнения

$$[a_1 v_1 \ldots v_{n-1}] = b$$

для любых $a_1, b \in A$. Поэтому, согласно основной теореме 3.1.1, $< A, [\ ] >$ – n-арная группа. ∎

Полагая в теореме 3.4.1 $i = 1$, получаем результат Поста (определение 1.1.3), а полагая в ней $i = n - 2$, $k = 1$ получаем следующий результат.

**3.4.2. Теорема.** n-Арная полугруппа $< A, [\ ] >$ является n-арной группой тогда и только тогда, когда для любых $a, b, c \in A$ в ней разрешимо уравнение

$$[ax_1 \ldots x_{n-2}c] = b$$

с $n - 2$ неизвестными.

Аналогично теореме 3.4.1 доказывается следующая

**3.4.3. Теорема.** n-Арная полугруппа $< A, [\ ] >$ является n-арной группой тогда и только тогда, когда для некоторых k, i ∈ {1, …, n – 2} таких, что k + i ≤ n – 1 и любых a, b ∈ A в ней разрешимо уравнение

$$[\underbrace{a \ldots a}_{k} x_1 \ldots x_i \underbrace{a \ldots a}_{n-k-i}] = b.$$

Полагая в теореме 3.4.3 i = 1, получаем результат А.Н. Скибы и В.И. Тютина [50] (определение 3.4.17), а, полагая в ней i = n – 2, k = 1 получаем следующий результат

**3.4.4. Теорема** [50]. n-Арная полугруппа < A, [ ] > является n-арной группой тогда и только тогда, когда для любых a, b ∈ A в ней разрешимо уравнение

$$[ax_1 \ldots x_{n-2}a] = b.$$

Последнюю теорему можно сформулировать иначе.

**3.4.5. Теорема** [50]. n-Арная полугруппа < A, [ ] > является n-арной группой тогда и только тогда, когда для любых a, b ∈ A существуют последовательности γ(a, b) и δ(a, b), составленные из элементов множества A и сумма длин которых равна n – 3, такие, что в A разрешимо уравнение

$$[a\gamma(a, b)x\delta(a, b)a] = b.$$

**3.4.6. Замечание.** Для уравнения (1) из теоремы 3.4.1 имеют место утверждения, аналогичные утверждениям из предложения 3.2.15. В частности, число различных решений этого уравнения в конечной n-арной группе < A, [ ] > равно $|A|^{i-1}$.

Теоремам 3.4.1 – 3.4.5 соответствуют следующие определения

**3.4.7. Определение.** n-Арная полугруппа < A, [ ] > называется n-арной группой если для некоторых k, i ∈ {1,…, n–2}, удовлетворяющих неравенству k + i ≤ n – 1, и любых $a_1, \ldots, a_k, a_{k+i+1}, \ldots, a_n, b \in A$ в ней разрешимо уравнение

$$[a_1 \ldots a_k x_1 \ldots x_i a_{k+i+1} \ldots a_n] = b.$$

**3.4.8. Определение.** n-Арная полугруппа $< A, [\,] >$ называется n-арной группой, если для любых a, b, c ∈ A в в ней разрешимо уравнение

$$[ax_1 \ldots x_{n-2}c] = b$$

с n − 2 неизвестными.

**3.4.9. Определение.** n-Арная полугруппа $< A, [\,] >$ называется n-арной группой, если существуют k, i ∈ {1, …, n − 2}, удовлетворяющие неравенству k + i ≤ n − 1, такие, что для любых a, b ∈ A в ней разрешимо уравнение

$$[\underbrace{a \ldots a}_{k} x_1 \ldots x_i \underbrace{a \ldots a}_{n-k-i}] = b.$$

**3.4.10. Определение.** n-Арная полугруппа $< A, [\,] >$ называется n-арной группой, если для любых a, b ∈ A в ней разрешимо уравнение

$$[ax_1 \ldots x_{n-2}a] = b.$$

**3.4.11. Определение.** n-Арная полугруппа $< A, [\,] >$ называется n-арной группой, если для любых a, b ∈ A существуют последовательности γ(a, b) и δ(a, b), составленные из элементов множества A и сумма длин которых равна n − 3, такие, что в ней разрешимо уравнение

$$[a\gamma(a, b)x\delta(a, b)a] = b.$$

Взяв за основу определение 3.4.11, можно по схеме, использованной в §3.3 (определение 3.3.4 – 3.3.21), получить серию других определений n-арной группы (см., например, [58]). Приведем только несколько определений из этой серии.

Ниже во всех определениях этого параграфа зафиксируем i ∈ {2, …, n − 1}.

**3.4.12. Определение.** n-Арная полугруппа $< A, [\ ] >$ называется n-арной группой, если существуют $c_1, \ldots, c_{i-2}, d_1, \ldots, d_{n-i-1} \in A$ такие, что для любых $a, b \in A$ в ней разрешимо уравнение

$$[ac_1 \ldots c_{i-2} x d_1 \ldots d_{n-i-1} a] = b.$$

**3.4.13. Определение.** n-Арная полугруппа $< A, [\ ] >$ называется n-арной группой, если существуют $c, d \in A$ такие, что для любых $a, b \in A$ в ней разрешимо уравнение

$$[a\underbrace{c \ldots c}_{i-2} x \underbrace{d \ldots d}_{n-i-1} a] = b.$$

**3.4.14. Определение** [50]. n-Арная полугруппа $< A, [\ ] >$ называется n-арной группой, если существует $c \in A$ такой, что для любых $a, b \in A$ в ней разрешимо уравнение

$$[a\underbrace{c \ldots c}_{i-2} x \underbrace{c \ldots c}_{n-i-1} a] = b.$$

**3.4.15. Определение.** [51]. n-Арная полугруппа $< A, [\ ] >$ называется n-арной группой, если для фиксированных

$$i \in \{2, \ldots, n-1\}, k \in \{1, \ldots, i-1\}, m \in \{1, \ldots, n-1\}$$

и любых $a, b \in A$ в ней разрешимо уравнение

$$[\underbrace{a \ldots a}_{k} \underbrace{b \ldots b}_{i-k-1} x \underbrace{b \ldots b}_{n-i-m} \underbrace{a \ldots a}_{m}] = b.$$

**3.4.16. Определение.** n-Арная полугруппа $< A, [\ ] >$ называется n-арной группой, если для любых $a, b \in A$ в ней разрешимо уравнение

$$[a\underbrace{b \ldots b}_{i-2} x \underbrace{b \ldots b}_{n-i-1} a] = b.$$

**3.4.17. Определение.** [49, Скиба А.Н., Тютин В.И.]. n-Арная полугруппа $<A, [\,]>$ называется n-арной группой, если для любых $a, b \in A$ в ней разрешимо уравнение

$$[\underbrace{a \ldots a}_{i-1} x \underbrace{a \ldots a}_{n-i}] = b.$$

## §3.5. ОПРЕДЕЛЕНИЯ, НЕ СВЯЗАННЫЕ С РАЗРЕШИМОСТЬЮ УРАВНЕНИЙ

В этом параграфе по-прежнему $n \geq 3$, если не указаны другие значения n.

**3.5.1. Теорема.** n-Арная полугруппа $<A, [\,]>$ является n-арной группой тогда и только тогда, когда для любого $a \in A$ существуют такие последовательности $\alpha(a)$ и $\beta(a)$ длины $n-2$ каждая, составленные из элементов множества A, что для любого $b \in A$ верно

$$[b\alpha(a)a] = b = [a\beta(a)b]. \qquad (1)$$

*Доказательство. Необходимость.* Если $<A, [\,]>$ – n-арная группа, то для любых $a, b \in A$ верно

$$[b\bar{a}\underbrace{a \ldots a}_{n-3}a] = b, \quad [a\underbrace{a \ldots a}_{n-3}\bar{a}b] = b,$$

где $\bar{a}$ – косой элемент для a, то есть существуют последовательности

$$\alpha(a) = \bar{a}\underbrace{a \ldots a}_{n-3}, \quad \beta(a) = \underbrace{a \ldots a}_{n-3}\bar{a}$$

такие, что верно (1).

*Достаточность.* Пусть теперь $<A, [\,]>$ – n-арная полугруппа, в которой для любого $a \in A$ существуют последовательности $\alpha(a)$ и $\beta(a)$ такие, что для любого $b \in A$ верно (1).

Это означает, что для любых a, b ∈ A в A разрешимы уравнения

$$[x_1 \ldots x_{n-1}a] = b, \quad [ay_1 \ldots y_{n-1}] = b.$$

Тогда по основной теореме 3.1.1 < A, [ ] > является n-арной группой. ∎

Теореме 3.5.1 соответствует

**3.5.2. Определение.** n-Арная полугруппа < A, [ ] > называется n-арной группой, если для любого a ∈ A существуют последовательности α(a) и β(a) длины n − 2 каждая, составленные из элементов множества A, такие, что для любого b ∈ A верно

$$[b\alpha(a)a] = b = [a\beta(a)b].$$

Следствием основной теоремы 3.1.1 является и следующее

**3.5.3. Определение** [52, Celakoski N.]. n-Арная полугруппа < A, [ ] > называется n-арной группой, если для фиксированного k ∈ {1, …, n − 2} и любых $a_1, \ldots, a_k \in A$ существуют $a'_1, \ldots, a'_{n-k-1} \in A$ такие, что для любого b ∈ A верно

$$[b a'_1 \ldots a'_{n-k-1} a_1 \ldots a_k] = b = [a_1 \ldots a_k a'_1 \ldots a'_{n-k-1} b].$$

Следующее определение может быть получено и как следствие основной теоремы 3.1.1 и как следствие теоремы 3.5.1 при

$$\alpha(a) = \underbrace{a \ldots a}_{n-3} \bar{a}, \quad \beta(a) = \bar{a} \underbrace{a \ldots a}_{n-3}.$$

**3.5.4. Определение** [53, Dudek W., Glazek K., Gleichgewicht B.]. n-Арная полугруппа < A, [ ] > называется n-арной группой, если для любого a ∈ A существует $\bar{a} \in A$ такой, что для любого b ∈ A верно

$$[b\underbrace{a \ldots a}_{n-3} \bar{a} a] = b = [a \bar{a} \underbrace{a \ldots a}_{n-3} b].$$

В дальнейшем нам понадобится

**3.5.5. Лемма.** Пусть $k \geq 1$ при $n \geq 3$, $k \geq 2$ при $n = 2$. Тогда для всякой n-арной полугруппы $< A, [\ ] >$ справедливы следующие утверждения:

1) если для любого $a \in A$ существует $\bar{a} \in A$ такой, что для любого $b \in A$ верно

$$[\bar{a}\underbrace{a \ldots a}_{k(n-1)-1} b] = b, \qquad (1)$$

то

$$[\bar{a}\underbrace{a \ldots a}_{k(n-1)-1} b] = [\underbrace{a \ldots a}_{k(n-1)-1} \bar{a} b];$$

2) если для любого $a \in A$ существует $\bar{a} \in A$ такой, что для любого $b \in A$ верно

$$[b\underbrace{a \ldots a}_{k(n-1)-1} \bar{a}] = b,$$

то

$$[b\underbrace{a \ldots a}_{k(n-1)-1} \bar{a}] = [b \bar{a}\underbrace{a \ldots a}_{k(n-1)-1}].$$

*Доказательство.* 1) Положим

$$[\underbrace{a \ldots a}_{k(n-1)-1} \bar{a} b] = c. \qquad (2)$$

По условию для $\bar{a} \in A$ существует $\bar{\bar{a}} \in A$ такой, что

$$[\bar{\bar{a}} \underbrace{\bar{a} \ldots \bar{a}}_{k(n-1)-1} d] = d$$

для любого $d \in A$, в частности,

$$[\underbrace{\bar{\bar{a}}\,\bar{a}\,\ldots\,\bar{a}}_{k(n-1)-1}b] = b \qquad (3)$$

$$[\underbrace{\bar{\bar{a}}\,\bar{a}\,\ldots\,\bar{a}}_{k(n-1)-1}c] = c. \qquad (4)$$

Подставляя в (4) вместо c левую часть равенства (2), а затем, используя (1) и (3), получаем:

$$c = [\underbrace{\bar{\bar{a}}\,\bar{a}\,\ldots\,\bar{a}}_{k(n-1)-1}c] = [\underbrace{\bar{\bar{a}}\,\bar{a}\,\ldots\,\bar{a}}_{k(n-1)-1}[\underbrace{a\,\ldots\,a\,\bar{a}}_{k(n-1)-1}b]] =$$

$$= [\underbrace{\bar{\bar{a}}\,\bar{a}\,\ldots\,\bar{a}}_{k(n-1)-2}[\underbrace{\bar{a}\,a\,\ldots\,a\,\bar{a}}_{k(n-1)-1}]b] = [\underbrace{\bar{\bar{a}}\,\bar{a}\,\ldots\,\bar{a}}_{k(n-1)-1}b] = b.$$

Таким образом, мы показали, что b = c, то есть

$$[\underbrace{\bar{a}\,a\,\ldots\,a}_{k(n-1)-1}b] = [\underbrace{a\,\ldots\,a\,\bar{a}}_{k(n-1)-1}b].$$

Утверждение 2) доказывается аналогично. ∎

Полагая в теореме 3.5.1

$$\alpha(a) = \bar{a}\underbrace{a\ldots a}_{n-3}, \quad \beta(a) = \underbrace{a\ldots a}_{n-3}\bar{a}$$

и применяя лемму 3.5.5, получим

**3.5.6. Определение** [53, Dudek W., Glazek K., Cleichgewicht B.5]. n-Арная полугруппа < A, [ ] > называется n-арной группой, если для любого a ∈ A существует $\bar{a}$ ∈ A такой, что для любого b ∈ A верно

$$[b\underbrace{a\ldots a}_{n-2}\bar{a}] = b = [\bar{a}\underbrace{a\,a\,\ldots\,a}_{n-2}b].$$

Для формулировки следующего определения полагаем в теореме 3.5.1:

1) если $s \neq 0$, то $\alpha(a) = \underbrace{a \ldots a}_{n-s-2} \bar{a} \underbrace{a \ldots a}_{s-1}$;

2) если $r \neq 0$, то $\beta(a) = \underbrace{a \ldots a}_{r-1} \bar{a} \underbrace{a \ldots a}_{n-r-2}$;

3) если $s = 0$, то $\alpha(a) = \bar{a} \underbrace{a \ldots a}_{n-3}$ и применяем лемму 3.5.5;

4) если $r = 0$, то $\beta(a) = \underbrace{a \ldots a}_{n-3} \bar{a}$ и применяем лемму 3.5.5.

**3.5.7. Определение** [53, Dudek W., Glazek K., Gleichgewicht B.; 52, Celakoski N.]. n-Арная полугруппа $< A, [\ ] >$ называется n-арной группой, если для фиксированных $r, s$ из множества $\{0, 1, \ldots, n-2\}$ и любого $a \in A$ существует $\bar{a} \in A$ такой, что для любого $b \in A$ верно

$$[b\underbrace{a \ldots a}_{n-s-2} \bar{a} \underbrace{a \ldots a}_{s}] = b = [\underbrace{a \ldots a}_{r} \bar{a} \underbrace{a \ldots a}_{n-r-2} b].$$

## §3.6. ОПРЕДЕЛЕНИЯ, СВЯЗАННЫЕ С ДЛИННЫМИ ОПЕРАЦИЯМИ

Во всех предыдущих определениях n-арной группы внутри квадратных скобок $[\ ]$, обозначающих n-арную операцию, стояли n элементов. В этом параграфе мы рассмотрим определения n-арной группы с длинными операциями, то есть когда внутри скобок $[\ ]$ стоит $k(n-1) + 1$ элементов, где $k \geq 1$.

**3.6.1. Теорема.** n-Арная полугруппа $< A, [\ ] >$ является n-арной группой тогда и только тогда, когда для любых $a, b \in A$ в $A$ разрешимы уравнения

$$[u_1 \ldots u_{k(n-1)}a] = b, \quad [av_1 \ldots v_{m(n-1)}] = b, \qquad (1)$$

где $k \geq 1, m \geq 1$.

***Доказательство.*** При k = 1, m = 1 доказывать нечего. Поэтому k > 1, m > 1 и пусть < A, [ ] > – n-арная полугруппа, в которой разрешимы уравнения (1), то есть существуют

$$c_i = u_i \, (i = 1, \ldots, k(n-1)), \quad d_j = v_j \, (j = 1, \ldots, m(n-1))$$

такие, что

$$[c_1 \ldots c_{n-2} c_{n-1} \ldots c_{k(n-1)} a] = b, \quad [a d_1 \ldots d_{n-2} d_{n-1} \ldots d_{m(n-1)}] = b,$$

откуда

$$[c_1 \ldots c_{n-2}[c_{n-1} \ldots c_{k(n-1)}]a] = b, \; [a d_1 \ldots d_{n-2}[d_{n-1} \ldots d_{m(n-1)}]] = b.$$

Это означает, что

$$x_1 = c_1, \ldots, x_{n-2} = c_{n-2}, \quad x_{n-1} = [c_{n-1} \ldots c_{k(n-1)}]$$

являются решениями уравнения

$$[x_1 \ldots x_{n-1} a] = b,$$

$$y_1 = d_1, \ldots, y_{n-2} = d_{n-2}, \quad y_{n-1} = [d_{n-1} \ldots d_{m(n-1)}]$$

являются решениями уравнения

$$[a y_1 \ldots y_{n-1}] = b.$$

Тогда по основной теореме 3.1.1 < A, [ ] > – n-арная группа.

Если теперь < A, [ ] > – n-арная группа, то по основной теореме 3.1.1 для любых

$$c_n, \ldots, c_{k(n-1)}, \; d_1, \ldots, d_{(m-1)(n-1)}, \; a, b \in A$$

в A разрешимы уравнения

$$[x_1 \ldots x_{n-1}[c_n \ldots c_{k(n-1)} a]] = b, \quad [[a d_1 \ldots d_{(m-1)(n-1)}] y_1 \ldots y_{n-1}] = b,$$

то есть существуют $c_1, \ldots, c_{n-1}, d_{(m-1)(n-1)+1}, \ldots, d_{m(n-1)} \in A$ такие, что

$$[c_1 \ldots c_{n-1}[c_n \ldots c_{k(n-1)} a]] = b,$$

$$[[ad_1 \ldots d_{(m-1)(n-1)}]d_{(m-1)(n-1)+1}, \ldots, d_{m(n-1)}] = b,$$

откуда

$$[c_1 \ldots c_{k(n-1)}a] = b, \quad [ad_1 \ldots d_{m(n-1)}] = b.$$

Это означает, что

$$u_1 = c_1, \ldots, u_{k(n-1)} = c_{k(n-1)}$$

– решение первого уравнения (1),

$$v_1 = d_1, \ldots, v_{m(n-1)} = d_{m(n-1)}$$

– решение второго уравнения (1). ∎

Заметим, что при $k = m = 1$ теорема 3.6.1 формально включает в себя основную теорему 3.1.1.

Ясно, что теорема 3.6.1 является «длинным» аналогом теоремы 3.1.1. Сформулируем «длинные» аналоги теорем 3.4.4 и 3.5.1.

**3.6.2. Теорема.** n-Арная полугруппа $< A, [\ ] >$ является n-арной группой тогда и только тогда, когда для любых $a, b \in A$ в ней разрешимо уравнение

$$[ax_1 \ldots x_{k(n-1)-1}a] = b.$$

**3.6.3. Теорема.** n-Арная полугруппа $< A, [\ ] >$ является n-арной группой тогда и только тогда, когда для любого $a \in A$ существуют последовательности

$$\alpha(a) = c_1(a) \ldots c_{k(n-1)-1}(a) \quad (k \geq 1, \ k(n-1) \geq 2),$$

$$\beta(a) = d_1(a) \ldots d_{m(n-1)-1}(a) \quad (m \geq 1, \ m(n-1) \geq 2),$$

составленные из элементов множества A, такие, что для любого $b \in A$ верно

$$[b\alpha(a)a] = b = [a\beta(a)b].$$

Заметим, что при $k > 1$ теоремы 3.6.2 и 3.6.3 справедливы для всех $n \geq 2$, то есть включают в себя бинарный случай ($n = 2$).

Теореме 3.6.1 соответствует

**3.6.4. Определение.** n-Арная полугруппа $< A, [\ ] >$ называется n-арной группой, если для любых $a, b \in A$ в ней разрешимы уравнения

$$[x_1 \ldots x_{k(n-1)}a] = b, \quad [ay_1 \ldots y_{m(n-1)}] = b,$$

где $k \geq 1$, $m \geq 1$.

Теореме 3.6.2 соответствует

**3.6.5. Определение.** n-Арная полугруппа $< A, [\ ] >$ называется n-арной группой, если для любых $a, b \in A$ в ней разрешимо уравнение

$$[ax_1 \ldots x_{k(n-1)-1}a] = b.$$

Теореме 3.6.3 соответствует

**3.6.6. Определение.** n-Арная полугруппа $< A, [\ ] >$ называется n-арной группой, если для любого $a \in A$ существуют элементы

$$c_1(a), \ldots, c_{k(n-1)-1}(a) \in A, \quad (k \geq 1,\ k(n-1) \geq 2),$$

$$d_1(a), \ldots, d_{m(n-1)-1}(a) \in A, \quad (m \geq 1,\ m(n-1) \geq 2),$$

такие, что верно

$$[bc_1(a) \ldots c_{k(n-1)-1}(a)a] = b = [ad_1(a) \ldots d_{m(n-1)-1}(a)b].$$

С помощью теорем 3.6.1 – 3.6.3 можно получать "длинные" аналоги многих приведенных в этой главе определений n-арной группы. Мы не будем перечислять все эти определения. Приведем только некоторые из них, полученные с помощью теоремы 3.6.3.

Для формулировки следующего определения полагаем в теореме 3.6.3:

1) если r ≠ 0, то α(a) = $\underbrace{a \ldots a}_{k(n-1)-r-1} \bar{a} \underbrace{a \ldots a}_{r-1}$;

2) если s ≠ 0, то β(a) = $\underbrace{a \ldots a}_{s-1} \bar{a} \underbrace{a \ldots a}_{m(n-1)-s-1}$;

3) если r = 0, то α(a) = $\bar{a} \underbrace{a \ldots a}_{k(n-1)-2}$ и применяем лемму 3.5.5;

4) если s = 0, то β(a) = $\underbrace{a \ldots a}_{m(n-1)-2} \bar{a}$ и применяем лемму 3.5.5.

**3.6.7. Определение.** n-Арная полугруппа $<A, [\,]>$ называется n-арной группой, если для фиксированных

$$r \in \{0, 1, \ldots, k(n-1)-1\}, s \in \{0, 1, \ldots, m(n-1)-1\},$$

где $k(n-1) \geq 2$, $m(n-1) \geq 2$ и любого $a \in A$ существует $\bar{a} \in A$ такой, что для любого $b \in A$ верно

$$[b \underbrace{a \ldots a}_{k(n-1)-r-1} \bar{a} \underbrace{a \ldots a}_{r}] = b = [\underbrace{a \ldots a}_{s} \bar{a} \underbrace{a \ldots a}_{m(n-1)-s-1} b].$$

Если в определении 3.6.7 положить $k = m$, то получим

**3.6.8. Определение** [54, Русаков С.А.]. n-Арная полугруппа $<A, [\,]>$ называется n-арной группой, если для фиксированных $r, s \in \{0, 1, \ldots, k(n-1)-1\}$, где $k(n-1) \geq 2$ и любого $a \in A$ существует $\bar{a} \in A$ такой, что для любого $b \in A$ верно

$$[b \underbrace{a \ldots a}_{k(n-1)-r-1} \bar{a} \underbrace{a \ldots a}_{r}] = b = [\underbrace{a \ldots a}_{s} \bar{a} \underbrace{a \ldots a}_{k(n-1)-s-1} b].$$

Полагая в теореме 3.6.3

$$\alpha(a) = \underbrace{a \ldots a}_{2(n-2)} \bar{a}, \quad \beta(a) = \bar{a} \underbrace{a \ldots a}_{2(n-2)}$$

и учитывая ассоциативность n-арной операции в n-арной полугруппе, получим

**3.6.9. Определение** [54, Русаков С.А.]. n-Арная полугруппа < A, [ ] > называется n-арной группой, если для любого a ∈ A существует $\bar{a}$ ∈ A такой, что для любого b ∈ A верно

$$[b\underbrace{a \ldots a}_{n-2}[\underbrace{a \ldots a}_{n-2}\bar{a}a]] = b = [[a\bar{a}\underbrace{a \ldots a}_{n-2}]\underbrace{a \ldots a}_{n-2}b].$$

Полагая в теореме 3.6.3

$$\alpha(a, b) = \bar{a}\underbrace{a \ldots a}_{2(n-2)}, \quad \beta(a, b) = \underbrace{a \ldots a}_{2(n-2)}\bar{a},$$

а также учитывая лемму 3.5.5 и ассоциативность n-арной операции в n-арной полугруппе, получим

**3.6.10. Определение** [54, Русаков С.А.]. n-Арная полугруппа < A, [ ] > называется n-арной группой, если для любого a ∈ A существует $\bar{a}$ ∈ A такой, что для любого b ∈ A верно

$$[[b\underbrace{a \ldots a}_{n-1}]\underbrace{a \ldots a}_{n-2}\bar{a}] = b = [\bar{a}\underbrace{a \ldots a}_{n-2}[\underbrace{a \ldots a}_{n-1}b]].$$

## §3.7. ДРУГИЕ ОПРЕДЕЛЕНИЯ n-АРНОЙ ГРУППЫ

Полученные в данном параграфе теоремы и следствия из них обобщают соответствующий групповой результат из [55]. При этом существенно используется следующая теорема, являющаяся следствием основной теоремы 3.1.1 (ср. с определением 3.3.19).

**3.7.1. Теорема.** n-Арная полугруппа < A, [ ] > является n-арной группой тогда и только тогда, когда для любых a, b ∈ A и фиксированных i, j ∈ {1, ..., n – 1} в A разрешимы уравнения

$$[\underbrace{a \ldots a}_{i} \underbrace{b \ldots b}_{n-1-i} x] = b, \quad [y \underbrace{b \ldots b}_{n-1-j} \underbrace{a \ldots a}_{j}] = b.$$

**3.7.2. Теорема** [56]. n-Арная полугруппа < A, [ ] > является n-арной группой (n ≥ 3) тогда и только тогда, когда существует такой элемент d ∈ A, что для любых a, b ∈ A и фиксированных i, j ∈ {1, ..., n – 1} в A разрешимы уравнения

$$[\underbrace{a \ldots a}_{i} \underbrace{b \ldots b}_{n-1-i} x] = d, \tag{1}$$

$$[y \underbrace{b \ldots b}_{n-1-j} \underbrace{a \ldots a}_{j}] = b. \tag{2}$$

*Доказательство.* Необходимость очевидна.

*Достаточность.* Если положить a = b = d, то из (1) вытекает разрешимость в A уравнения

$$[\underbrace{d \ldots d}_{n-1} x] = d,$$

т. е. существует элемент u ∈ A такой, что

$$[\underbrace{d \ldots d}_{n-1} u] = d.$$

Для любого b ∈ A с помощью уравнения (2) определим v ∈ A такой, что

$$[v \underbrace{b \ldots b}_{n-1-j} \underbrace{d \ldots d}_{j}] = b.$$

Так как

$$[b \underbrace{d \ldots d}_{n-2} u] = [[v \underbrace{b \ldots b}_{n-1-j} \underbrace{d \ldots d}_{j}] \underbrace{d \ldots d}_{n-2} u] =$$

$$= [v\underbrace{b\ldots b}_{n-1-j}\underbrace{d\ldots d}_{j-1}[\underbrace{d\ldots d}_{n-1}u]] =$$

$$= [v\underbrace{b\ldots b}_{n-1-j}\underbrace{d\ldots d}_{j-1}dd] = [v\underbrace{b\ldots b}_{n-1-j}\underbrace{d\ldots d}_{j}] = b,$$

то

$$[b\underbrace{d\ldots d}_{n-2}u] = b \qquad (3)$$

для любого $b \in A$.

Положив в (2) $a = b$, получим разрешимость в A уравнения вида

$$[y\underbrace{a\ldots a}_{n-1}] = a$$

для любого $a \in A$. Следовательно, для любого $a \in A$ существует $w \in A$ такой, что

$$[w\underbrace{a\ldots a}_{n-1}] = a.$$

Из (1) вытекает, что для любых $a, b \in A$ существует $z \in A$ такой, что

$$[\underbrace{a\ldots a}_{i}\underbrace{b\ldots b}_{n-1-i}z] = d. \qquad (4)$$

Так как

$$[w\underbrace{a\ldots a}_{n-2}d] = [w\underbrace{a\ldots a}_{n-2}[\underbrace{a\ldots a}_{i}\underbrace{b\ldots b}_{n-1-i}z]] =$$

$$= [[w\underbrace{a\ldots a}_{n-1}]\underbrace{a\ldots a}_{i-1}\underbrace{b\ldots b}_{n-1-i}z] =$$

$$= [a\underbrace{a\ldots a}_{i-1}\underbrace{b\ldots b}_{n-1-i}z] = [\underbrace{a\ldots a}_{i}\underbrace{b\ldots b}_{n-1-i}z] = d,$$

то
$$[w\underbrace{a\ldots a}_{n-2}d] = d$$

для любого $a \in A$.

Используя последнее равенство, а также (3) при $b = a$, получим

$$[\underbrace{d\ldots d}_{n-2}ua] = [\underbrace{dd\ldots d}_{n-3}ua] = [[w\underbrace{a\ldots a}_{n-2}d]\underbrace{d\ldots d}_{n-3}ua] =$$

$$= [w\underbrace{a\ldots a}_{n-3}[a\underbrace{d\ldots d}_{n-2}u]a] = [w\underbrace{a\ldots a}_{n-3}aa] = [w\underbrace{a\ldots a}_{n-1}] = a,$$

т. е.

$$[\underbrace{d\ldots d}_{n-2}ua] = a \qquad (5)$$

для любого $a \in A$.

Из (4) вытекает, что

$$[[\underbrace{a\ldots a}_{i}\underbrace{b\ldots b}_{n-1-i}z]\underbrace{d\ldots d}_{n-3}ub] = [\underbrace{dd\ldots d}_{n-3}ub],$$

откуда

$$[\underbrace{a\ldots a}_{i}\underbrace{b\ldots b}_{n-1-i}[z\underbrace{d\ldots d}_{n-3}ub]] = [\underbrace{d\ldots d}_{n-2}ub].$$

Из последнего равенства и из (5) при $a = b$ следует

$$[\underbrace{a\ldots a}_{i}\underbrace{b\ldots b}_{n-1-i}[z\underbrace{d\ldots d}_{n-3}ub]] = b,$$

а это означает, что элемент

$$[z\underbrace{d\ldots d}_{n-3}ub]$$

является решением уравнения

$$[\underbrace{a \ldots a}_{i} \underbrace{b \ldots b}_{n-1-i} x] = b. \qquad (6)$$

Так как для любых $a, b \in A$ и фиксированных $i, j \in \{1, \ldots, n-1\}$ в A разрешимы уравнения (2) и (6), то согласно теореме 3.1.1, $< A, [\ ] > -$ n-арная группа. ∎

Аналогично теореме 3.7.2 доказывается двойственная к ней

**3.7.3. Теорема.** n-Арная полугруппа $< A, [\ ] >$ является n-арной группой ($n \geq 3$) тогда и только тогда, когда существует такой элемент $d \in A$, что для любых $a, b \in A$ и фиксированных $i, j \in \{1, \ldots, n-1\}$ в A разрешимы уравнения

$$[\underbrace{a \ldots a}_{i} \underbrace{b \ldots b}_{n-1-i} x] = b, \quad [y \underbrace{b \ldots b}_{n-1-j} \underbrace{a \ldots a}_{j}] = d.$$

Придавая i и j в теоремах 3.7.2 и 3.7.3 различные значения, можно получить большое число новых определений n-арной группы, некоторые из которых приведены ниже.

Полагая в теореме 3.7.2 $i = j = 1$, получим

**3.7.4. Определение.** n-Арная полугруппа $< A, [\ ] >$ является n-арной группой тогда и только тогда, когда существует такой элемент $d \in A$, что для любых $a, b \in A$ в A разрешимы уравнения

$$[a \underbrace{b \ldots b}_{n-2} x] = d, \quad [y \underbrace{b \ldots b}_{n-2} a] = b,$$

Полагая в теореме 3.7.2 $i = j = n - 1$, получим

**3.7.5. Определение**. n-Арная полугруппа $< A, [\ ] >$ является n-арной группой тогда и только тогда, когда существует такой элемент $d \in A$, что для любых $a, b \in A$ в A разрешимы уравнения

$$[\underbrace{a \ldots a}_{n-1} x] = d, \quad [y \underbrace{a \ldots a}_{n-1}] = b.$$

Полагая в теореме 3.7.2 i = 1, j = n – 1, получим

**3.7.6. Определение.** n-Арная полугруппа $<A, [\,]>$ является n-арной группой тогда и только тогда, когда существует такой элемент $d \in A$, что для любых $a, b \in A$ в $A$ разрешимы уравнения

$$[a\underbrace{b \ldots b}_{n-2} x] = d, \quad [y\underbrace{a \ldots a}_{n-1}] = b.$$

Полагая в теореме 3.7.2 i = n – 1, j = 1, получим

**3.7.7. Определение.** n-Арная полугруппа $<A, [\,]>$ является n-арной группой тогда и только тогда, когда существует такой элемент $d \in A$, что для любых $a, b \in A$ в $A$ разрешимы уравнения

$$[\underbrace{a \ldots a}_{n-1} x] = d, \quad [y\underbrace{b \ldots b}_{n-2} a] = b.$$

Полагая в теореме 3.7.3 i = j = 1, получим

**3.7.8. Определение.** n-Арная полугруппа $<A, [\,]>$ является n-арной группой тогда и только тогда, когда существует такой элемент $d \in A$, что для любых $a, b \in A$ в $A$ разрешимы уравнения

$$[a\underbrace{b \ldots b}_{n-2} x] = b, \quad [y\underbrace{b \ldots b}_{n-2} a] = d,$$

Полагая в теореме 3.7.3 i = j = n – 1, получим

**3.7.9. Определение** [57, Тютин В.И.]. n-Арная полугруппа $<A, [\,]>$ является n-арной группой тогда и только тогда, когда существует такой элемент $d \in A$, что для любых $a, b \in A$ в $A$ разрешимы уравнения

$$[\underbrace{a \ldots a}_{n-1} x] = b, \quad [y\underbrace{a \ldots a}_{n-1}] = d.$$

Полагая в теореме 3.7.3 i = 1,  j = n – 1, получим

**3.7.10. Определение.** n-Арная полугруппа < A, [ ] > является n-арной группой тогда и только тогда, когда существует такой элемент d ∈ A, что для любых a, b ∈ A в A разрешимы уравнения

$$[a\underbrace{b \ldots b}_{n-2}x] = b, \quad [y\underbrace{a \ldots a}_{n-1}] = d.$$

Полагая в теореме 3.7.3 i = n – 1,  j = 1, получим

**3.7.11. Определение.** n-Арная полугруппа < A, [ ] > является n-арной группой тогда и только тогда, когда существует такой элемент d ∈ A, что для любых a, b ∈ A в A разрешимы уравнения

$$[\underbrace{a \ldots a}_{n-1} x] = b, \quad [y\underbrace{b \ldots b}_{n-2}a] = d.$$

## ДОПОЛНЕНИЯ И КОММЕНТАРИИ

**1.** Большинство результатов главы 3 были опубликованы в [58].

**2.** Аксиоматикой n-арных групп, помимо указанных выше авторов, занимались также Твермоес X. [59, 60], Робинсон Д. [61], Слипенко А.К. [62], Монк Д. и Сиосон Ф. [43], Ушан Я. [63]. Информация о работах по аксиоматике n-арных групп имеется в обзоре Глазека К. [11].

**3.** Б. Гляйхгевихт и К. Глазек первыми установили [64], что n-арную группу можно определить с помощью n-арной и унарной операций.

**4.** Представляет интерес следующее

**Определение** [43, Monk J.D., Sioson F.M.]. n-арная полугруппа < A, [ ] > называется n-арной группой, если для любых $a_2, \ldots, a_{n-2} \in A$ существует единственный элемент $(a_1, \ldots, a_{n-2})^{-1} \in A$ такой, что для любого b ∈ A верно

$$[(a_1, \ldots, a_{n-2})^{-1} a_1 \ldots a_{n-2} b] = b,$$

$$[b(a_1, \ldots, a_{n-2})^{-1}a_1 \ldots a_{n-2}] = b,$$

$$[a_1 \ldots a_{n-2}(a_1, \ldots, a_{n-2})^{-1}b] = b,$$

$$[ba_1 \ldots a_{n-2}(a_1, \ldots, a_{n-2})^{-1}] = b.$$

Фактически в определении Монка и Сиосона на множестве A определена $(n-2)$-арная операция. Поэтому естественным выглядит следующее

**Определение** [52, Celanoski N.]. n-Арная полугруппа $< A, [\ ] >$ называется n-арной группой, если на A существует такая $(n-2)$-арная операция $^{-1}$, что для любых $a_1, \ldots, a_{n-2}, b \in A$ верно

$$[(a_1, \ldots, a_{n-2})^{-1}a_1 \ldots a_{n-2}b] = b,$$

$$[ba_1 \ldots a_{n-2}(a_1, \ldots, a_{n-2})^{-1}] = b.$$

**5.** Еще дальше пошел Ушан Я., определив [63] на множестве A с одной n-арной операцией $[\ ]$ еще две операции:

1) $(n-2)$-арную операцию **e**, удовлетворяющую условию

$$[\mathbf{e}(a_1, \ldots, a_{n-2})a_1 \ldots a_{n-2}b] = b = [ba_1 \ldots a_{n-2}\mathbf{e}(a_1, \ldots, a_{n-2})]$$

для любого $b \in A$;

2) $(n-1)$-арную операцию $^{-1}$, удовлетворяющую условию

$$[(a_1, \ldots, a_{n-2}, a)^{-1}a_1 \ldots a_{n-2}a] = \mathbf{e}(a_1, \ldots, a_{n-2}),$$

$$[aa_1 \ldots a_{n-2}(a_1, \ldots, a_{n-2}, a)^{-1}] = \mathbf{e}(a_1, \ldots, a_{n-2}).$$

Операции **e** и $^{-1}$ позволили Ушану получить ряд определений n-арной группы в терминах операций **e** и $^{-1}$.

Ясно, что если $< A, [\ ], ^- >$ – n-арная группа с n-арной операцией $[\ ]$ и унарной операцией $^-$, то определение операции **e** равносильно определению обратного элемента для любой последовательности $a_1 \ldots a_{n-2}$, где $a_1, \ldots, a_{n-2} \in A$:

$$\mathbf{e}(a_1, \ldots, a_{n-2}) = [\,\overline{a}_{n-2}\underbrace{a_{n-2} \ldots a_{n-2}}_{n-3} \ldots \overline{a}_1 \underbrace{a_1 \ldots a_1}_{n-3}\,]. \qquad (*)$$

Поэтому **e** является $(n-2)$-арной операцией, производной от основных операций $[\ ]$ и $^-$ n-арной группы $< A, [\ ], ^- >$.

Ясно также, что в n-арной группе $< A, [\ ], ^- >$ верно

$$(a_1, \ldots, a_{n-2}, a)^{-1} = [\mathbf{e}(a_1, \ldots, a_{n-2})\overline{a}\underbrace{a \ldots a}_{n-3}\mathbf{e}(a_1, \ldots, a_{n-2})].$$

А так как правая часть последнего равенства является обратной для последовательности

$$a_1 \ldots a_{n-2}aa_1 \ldots a_{n-2},$$

то определение операции $^{-1}$ равносильно определению обратного элемента для любой такой последовательности. Учитывая (*), получим

$$(a_1, \ldots, a_{n-2}, a)^{-1} = [\overline{a}_{n-2}\underbrace{a_{n-2} \ldots a_{n-2}}_{n-3} \ldots \overline{a}_1\underbrace{a_1 \ldots a_1}_{n-3}$$

$$\overline{a}\underbrace{a \ldots a}_{n-3}\overline{a}_{n-2}\underbrace{a_{n-2} \ldots a_{n-2}}_{n-3} \ldots \overline{a}_1\underbrace{a_1 \ldots a_1}_{n-3}],$$

т.е. операция $^{-1}$ также является производной от основных операций [ ] и $^-$ n-арной группы $<A, [\ ], ^->$.

Так как обе операции **e** и $^{-1}$ являются производными от основных операций [ ] и $^-$ в n-арной группе $<A, [\ ], ^->$, то все свойства n-арных групп, в формулировках которых присутствуют нейтральные и обратные последовательности, можно переформулировать, используя операции **e** и $^{-1}$. В качестве примера рассмотрим следующее известное

**Предложение 1**[3, 4]. Если $b_1 \ldots b_{n-2}b$ – нейтральная последовательность n-арной группы $<A, [\ ]>$, то для любого $i = 2, \ldots, n-1$ последовательность $b_i \ldots b_{n-2}bb_1 \ldots b_{i-1}$ также является нейтральной в $<A, [\ ]>$.

Это предложение с помощью операции **e** формулируется следующим образом

**Предложение 2** [63, предложение 1.1 на с. 26]. Если $<A, [\ ]>$ – n-арная группа ($n \geq 3$), то для любых $a_1, \ldots, a_{n-2}, b_1, \ldots, b_{n-2}, x \in A$ и всех $i = 1, \ldots, n-1$ верно

$$[xb_i \ldots b_{n-2}\mathbf{e}(b_1, \ldots, b_{n-2})b_1 \ldots b_{i-1}] = [\mathbf{e}(a_1, \ldots, a_{n-2})a_1 \ldots a_{n-2}x],$$

$$[b_i \ldots b_{n-2}\mathbf{e}(b_1, \ldots, b_{n-2})b_1 \ldots b_{i-1}x] = [xa_1 \ldots a_{n-2}\mathbf{e}(a_1, \ldots, a_{n-2})].$$

**6.** Помимо ассоциативности n-арной операции, используемой в определениях n-арной группы, приведенных в данной главе, возможны и другие виды ассоциативности. Соответственно, возможны и другие n-арные обобщения понятия группы, отличные от определения Дёрнте. Одно из таких обобщений принадлежит Ф.Н. Сохацкому, ко-торый ввел понятие полиагруппы [65].

n-Арная квазигруппа < A, [ ] > называется полиагруппой сорта (s, n), где s делит n – 1, если для всех i, j таких, что i ≡ j(mod s) в < A, [ ] > выполняется тождество

$$[x_1 \ldots x_i[x_{i+1} \ldots x_{i+n}]x_{i+n+1} \ldots x_{2n-1}] = [x_1 \ldots x_j[x_{j+1} \ldots x_{j+n}]x_{j+n+1} \ldots x_{2n-1}].$$

Ясно, что n-арные группы – это в точности полиагруппы сорта (1, n).

Многие результаты для n-арных групп, например, по аксиоматике, обобщаются на случай полиагрупп [66 – 70].

**7.** Обобщением понятия n-арной группы являются (i, j)-ассоциативные n-арные квазигруппы, где (i, j) ∈ {1, …, n}, изучавшиеся В.Д. Белоусовым [71], Е.И. Соколовым [72] и другими.

**8.** Еще одно обобщение понятия n-арной группы предложил Г. Чупона, определив (n, m)-группы [73]. При этом в определении Дёрнте видоизменяется не только условие 1), но и условие 2).

**9.** В данной книге мы не рассматриваем топологические n-арные группы [4, 74 – 79] и упорядоченные n-арные группы [80, 81].

# Г Л А В А  4

# n-АРНЫЕ ПОДСТАНОВКИ И МОРФИЗМЫ

Многие математические понятия, зависящие от параметра n, первоначально были определены и изучались для фиксированного n, например, n = 1 или n = 2. Так, развитие теории групп стимулировало появление n-арных аналогов группы. n-Арными аналогами соответствующих бинарных понятий являются также n-арные подстановки и n-арные морфизмы, рассматриваемые в данной главе.

## §4.1. n-АРНЫЕ ПОДСТАНОВКИ

Изучавшиеся в работах Поста [3], С.А. Русакова [4] и Сиосона [45] последовательности $\{f_1, f_2, \ldots, f_{n-1}\}$ взаимно однозначных отображений

$$A_1 \xrightarrow{f_1} A_2 \xrightarrow{f_2} \ldots \xrightarrow{f_{n-2}} A_{n-1} \xrightarrow{f_{n-1}} A_1$$

являются естественным обобщением понятия обычной (бинарной) подстановки, и, как нетрудно заметить, определяются циклической подстановкой $(1\ 2\ \ldots\ n-1) \in S_{n-1}$. В данном параграфе изучаются последовательности взаимно однозначных отображений, определяемые произвольной подстановкой $\sigma \in S_{n-1}$.

Пусть $A_1, \ldots, A_{n-1}$ ($n \geq 2$) произвольные множества одинаковой мощности. Для всякой подстановки $\sigma \in S_{n-1}$ определим множество $S_{A_1, \ldots, A_{n-1}}(\sigma)$ всех последовательностей

$$f(\sigma) = \{\sigma, f_1, \ldots, f_{n-1}\}, \text{ где } f_j : A_j \to A_{\sigma(j)}$$

– взаимно однозначные отображения, $j = 1, \ldots, n-1$.

Иногда подстановку σ в записи f(σ) = {σ, f₁, ..., f_{n−1}} указывать не будем, то есть будем писать f = {f₁, ..., f_{n−1}}. Множество всех таких последовательностей будем обозначать символом $S_{A_1, ..., A_{n-1}}$.

**4.1.1. Пример.** Пусть n = 3, A₁ = A, A₂ = A, $\sigma = \begin{pmatrix} 1 & 2 \\ 1 & 2 \end{pmatrix} \in S_2$, $\delta = \begin{pmatrix} 1 & 2 \\ 2 & 1 \end{pmatrix} \in S_2$, f₁ и f₂ – биекции A на A. Тогда

f(σ) = {σ, f₁, f₂}, где f₁: A₁ → A₁, f₂: A₂ → A₂,

f(δ) = {δ, f₁, f₂}, где f₁: A₁ → A₂, f₂: A₂ → A₁.

Следовательно, для различных подстановок σ и δ последовательности f(σ) и f(δ) могут различаться только первыми элементами.

Ясно, что $S_{A_1, ..., A_{n-1}}(\sigma) \cap S_{A_1, ..., A_{n-1}}(\delta) = \varnothing$ для σ ≠ δ.

Если T ⊆ S_{n−1}, то положим

$$S_{A_1, ..., A_{n-1}}(T) = \bigcup_{\sigma \in T} S_{A_1, ..., A_{n-1}}(\sigma).$$

**4.1.2. Определение** [82]. *n-Арными подстановками последовательности $A_1, ..., A_{n-1}$* называются элементы множества $S_{A_1, ..., A_{n-1}}(S_{n-1})$.

При n = 2 определение 4.1.2 превращается в определение обычной (бинарной) подстановки. В этом случае S₁ состоит из единственной тождественной подстановки ε,

f(ε) = {ε, f}, f: A → A, $S_A(\varepsilon) = S_A(S_1) = S_A$.

Для любых m (m ≥ 1) последовательностей

$f_k = f_k(\sigma_k) = \{\sigma_k, f_{k1}, ..., f_{k(n-1)}\}$ из $S_{A_1, ..., A_{n-1}}(S_{n-1})$,

где k = 1, ..., m определим на $S_{A_1, ..., A_{n-1}}(S_{n-1})$ m-арную операцию

$$(f_1 \ldots f_m)_{m,\sigma_1,\ldots,\sigma_m} = \{\sigma, g_1, \ldots, g_{n-1}\} = g(\sigma),$$

где

$$\sigma = \sigma_1 \ldots \sigma_m,$$

$$g_j = f_{1j} f_{2\sigma_1(j)} f_{3\sigma_1\sigma_2(j)} \ldots f_{m\sigma_1\ldots\sigma_{m-1}(j)} \colon A_j \to A_{\sigma(j)} = A_{\sigma_1\ldots\sigma_m(j)}, \quad (*)$$

$j = 1, \ldots, n - 1$.

Как обычно, полагаем

$$\sigma_s(\ldots(\sigma_2(\sigma_1(j)))\ldots) = \sigma_1 \ldots \sigma_s(j).$$

Заметим, что в определении m-арной операции $(\ )_m$ подстановки $\sigma_1, \ldots, \sigma_m$ не обязательно все различные. В частном случае при $m = 1$ имеем одну последовательность

$$f(\sigma_1) = \{\sigma_1, f_1, \ldots, f_{n-1}\} \in S_{A_1,\ldots,A_{n-1}}(\sigma_1)$$

и по определению $\sigma = \sigma_1$, $g_j = f_j$ ($j = 1, \ldots, n-1$). Следовательно, $g(\sigma) = g(\sigma_1)$ и, таким образом, $(f)_1 = f$.

**4.1.3. Теорема** [82]. Для всех i и k таких, что

$$1 \leq i + 1 \leq i + k \leq m$$

и любых

$$f_1 = f_1(\sigma_1), \ldots, f_m = f_m(\sigma_m) \in S_{A_1,\ldots,A_{n-1}}(S_{n-1})$$

имеет место равенство

$$(f_1 \ldots f_m)_{m,\sigma_1,\ldots,\sigma_m} =$$
$$= (f_1 \ldots f_i (f_{i+1} \ldots f_{i+k})_{k,\sigma_{i+1},\ldots,\sigma_{i+k}} f_{i+k+1} \ldots$$
$$\ldots f_m)_{m-k+1,\sigma_1,\ldots,\sigma_i,\mu,\sigma_{i+k+1},\ldots,\sigma_m}, \qquad (1)$$

где $\mu = \sigma_{i+1} \ldots \sigma_{i+k}$.

***Доказательство.*** Случай k = 1 очевиден. Поэтому считаем k ≥ 2.

Положив

$$(f_{i+1} \ldots f_{i+k})_{k, \sigma_{i+1}, \ldots, \sigma_{i+k}} = \{\mu, h_1, \ldots, h_{n-1}\} = h(\mu),$$

имеем по определению

$$h_j = f_{(i+1)j} f_{(i+2)\sigma_{i+1}(j)} \ldots f_{(i+k)\sigma_{i+1} \ldots \sigma_{i+k-1}(j)}, \; j = 1, \ldots, n-1. \quad (2)$$

Положим также

$$(f_1 \ldots f_i (f_{i+1} \ldots f_{i+k})_{k, \, \sigma_{i+1}, \, \ldots, \, \sigma_{i+k}} f_{i+k+1} \ldots$$
$$\ldots f_m)_{m-k+1, \, \sigma_1, \, \ldots, \, \sigma_i, \, \mu, \, \sigma_{i+k+1}, \, \ldots, \, \sigma_m} = \{\delta, l_1, \ldots, l_{n-1}\} = l(\delta), \quad (3)$$

где по определению

$$\delta = \sigma_1 \ldots \sigma_i \mu \sigma_{i+k+1} \ldots \sigma_m = \sigma_1 \ldots \sigma_i \sigma_{i+1} \ldots \sigma_{i+k} \sigma_{i+k+1} \ldots \sigma_m = \sigma.$$

Таким образом,

$$\delta = \sigma. \quad (4)$$

1) Если i = 0, то $\mu = \sigma_1 \ldots \sigma_k$ и, согласно (2) и (3),

$$h_j = f_{1j} f_{2\sigma_1(j)} \ldots f_{k\sigma_1 \ldots \sigma_{k-1}(j)}, \; j = 1, \ldots, n-1,$$

$$((f_1 \ldots f_k)_{k, \, \sigma_1, \, \ldots, \, \sigma_k} f_{k+1} \ldots f_m)_{m-k+1, \, \mu, \, \sigma_{k+1}, \, \ldots, \, \sigma_m} = \{\delta, l_1, \ldots, l_{n-1}\},$$

где ввиду (*),

$$l_j = h_j f_{(k+1)\mu(j)} f_{(k+2)\mu\sigma_{k+1}(j)} \ldots f_{m\mu\sigma_{k+1} \ldots \sigma_{m-1}(j)} =$$

$$= \underbrace{f_{1j} f_{2\sigma_1(j)} \ldots f_{k\sigma_1 \ldots \sigma_{k-1}(j)}}_{h_j} f_{(k+1)\sigma_1 \ldots \sigma_k(j)} f_{(k+2)\sigma_1 \ldots \sigma_k \sigma_{k+1}(j)} \ldots$$

$$\ldots f_{m\sigma_1 \ldots \sigma_k \sigma_{k+1} \ldots \sigma_{m-1}(j)} = f_{1j} f_{2\sigma_1(j)} f_{3\sigma_1\sigma_2(j)} \ldots f_{m\sigma_1 \ldots \sigma_{m-1}(j)} = g_j,$$

то есть $l_j = g_j$ для всех $j = 1, \ldots, n - 1$. А так как, кроме того, ввиду (4), $\delta = \sigma$, то

$$(f_1 \ldots f_m)_{m, \sigma_1, \ldots, \sigma_m} = ((f_1 \ldots f_k)_{k, \sigma_1, \ldots, \sigma_k} f_{k+1} \ldots$$

$$\ldots f_m)_{m-k+1, \mu, \sigma_{k+1}, \ldots, \sigma_m},$$

то есть при $i = 0$ равенство (1) верно.

2) Если $i = m - k$, то $\mu = \sigma_{m-k+1} \ldots \sigma_m$ и, согласно (2) и (3),

$$h_{\sigma_1 \ldots \sigma_{m-k}(j)} = f_{(m-k+1)\sigma_1 \ldots \sigma_{m-k}(j)} f_{(m-k+2)\sigma_{m-k+1}(\sigma_1 \ldots \sigma_{m-k}(j))} \cdots$$

$$\ldots f_{m\sigma_{m-k+1} \ldots \sigma_{m-1}(\sigma_1 \ldots \sigma_{m-k}(j))} =$$

$$= f_{(m-k+1)\sigma_1 \ldots \sigma_{m-k}(j)} f_{(m-k+2)\sigma_1 \ldots \sigma_{m-k} \sigma_{m-k+1}(j)} \cdots$$

$$\ldots f_{m\sigma_1 \ldots \sigma_{m-k} \sigma_{m-k+1} \ldots \sigma_{m-1}(j)}, \; j = 1, \ldots, n-1,$$

$$(f_1 \ldots f_{m-k}(f_{m-k+1} \ldots f_m)_{k, \sigma_{m-k+1}, \ldots, \sigma_m})_{m-k+1, \sigma_1, \ldots, \sigma_{m-k}, \mu} =$$

$$= \{\delta, l_1, \ldots, l_{n-1}\},$$

где ввиду (*),

$$l_j = f_{1j} f_{2\sigma_1(j)} f_{3\sigma_1 \sigma_2(j)} \ldots f_{(m-k)\sigma_1 \ldots \sigma_{m-k-1}(j)} h_{\sigma_1 \ldots \sigma_{m-k}(j)} =$$

$$= f_{1j} f_{2\sigma_1(j)} f_{3\sigma_1 \sigma_2(j)} \ldots$$

$$\ldots f_{(m-k)\sigma_1 \ldots \sigma_{m-k-1}(j)} \underbrace{f_{(m-k+1)\sigma_1 \ldots \sigma_{m-k}(j)} \ldots f_{m\sigma_1 \ldots \sigma_{m-1}(j)}}_{h_{\sigma_1 \ldots \sigma_{m-k}(j)}} =$$

$$= f_{1j} f_{2\sigma_1(j)} f_{3\sigma_1 \sigma_2(j)} \ldots f_{m\sigma_1, \ldots, \sigma_{m-1}(j)} = g_j,$$

то есть $l_j = g_j$ для всех $j = 1, \ldots n - 1$. Таким образом, установлено, что

$$(f_1 \ldots f_m)_{m, \sigma_1, \ldots, \sigma_m} =$$

$$= (f_1 \ldots f_{m-k}(f_{m-k+1} \ldots f_m)_{k, \sigma_{m-k+1}, \ldots, \sigma_m})_{m-k+1, \sigma_1, \ldots, \sigma_{m-k}, \mu},$$

то есть при i = m – k равенство (1) также верно.

3) Если $0 < i < m - k$, то ввиду (*)

$$l_j = f_{1j}f_{2\sigma_1(j)} \ldots f_{i\sigma_1\ldots\sigma_{i-1}(j)} h_{\sigma_1\ldots\sigma_{i-1}\sigma_i(j)} f_{(i+k+1)\sigma_1\ldots\sigma_{i-1}\sigma_i\mu(j)} \ldots$$

$$\ldots f_{m\sigma_1\ldots\sigma_{i-1}\sigma_i\mu\ldots\sigma_{m-1}(j)} = f_{1j}f_{2\sigma_1(j)} \ldots f_{i\sigma_1\ldots\sigma_{i-1}(j)}$$

$$\underbrace{f_{(i+1)\sigma_1\ldots\sigma_{i-1}\sigma_i(j)}f_{(i+2)\sigma_{i+1}(\sigma_1\ldots\sigma_{i-1}\sigma_i(j))} \ldots f_{(i+k)\sigma_{i+1}\ldots\sigma_{i+k-1}(\sigma_1\ldots\sigma_{i-1}\sigma_i(j))}}_{h_{\sigma_1\ldots\sigma_{i-1}\sigma_i(j)}}$$

$$f_{(i+k+1)\sigma_1\ldots\sigma_i\underbrace{\sigma_{i+1}\ldots\sigma_{i+k}}_{\mu}(j)} \ldots f_{m\sigma_1\ldots\sigma_i\underbrace{\sigma_{i+1}\ldots\sigma_{i+k}}_{\mu}\ldots\sigma_{m-1}(j)} =$$

$$= f_{1j}f_{2\sigma_1(j)} \ldots f_{i\sigma_1\ldots\sigma_{i-1}(j)} \underbrace{f_{(i+1)\sigma_1\ldots\sigma_i(j)} \ldots f_{(i+k)\sigma_1\ldots\sigma_i\sigma_{i+1}\ldots\sigma_{i+k-1}(j)}}_{h_{\sigma_1\ldots\sigma_{i-1}\sigma_i(j)}}$$

$$f_{(i+k+1)\sigma_1\ldots\sigma_{i+k}(j)} \ldots f_{m\sigma_1\ldots\sigma_{m-1}(j)} = f_{1j}f_{2\sigma_1(j)} \ldots f_{m\sigma_1\ldots\sigma_{m-1}(j)} = g_j,$$

то есть $l_j = g_j$ для всех $j = 1, \ldots n - 1$. откуда $l(\delta) = g(\sigma)$. Следовательно, верно равенство (1). ∎

Если в теореме 4.1.3 положить

$$m = 2k - 1, \sigma_1 = \ldots = \sigma_{2k-1} = \sigma, \sigma^k = \sigma,$$

то равенство (1) перепишется следующим образом

$$(f_1 \ldots f_{2k-1})_{2k-1,\underbrace{\sigma,\ldots,\sigma}_{2k-1}} =$$

$$= (f_1 \ldots f_i(f_{i+1} \ldots f_{i+k})_{k,\underbrace{\sigma,\ldots,\sigma}_{k}} f_{i+k+1} \ldots$$

$$\ldots f_{2k-1})_{k,\underbrace{\sigma,\ldots,\sigma}_{i},\sigma^k = \sigma,\underbrace{\sigma,\ldots,\sigma}_{k-i-1}} =$$

$$= (f_1 \ldots f_i(f_{i+1} \ldots f_{i+k})_{k,\underbrace{\sigma,\ldots,\sigma}_{k}} f_{i+k+1} \ldots f_{2k-1})_{k,\underbrace{\sigma,\ldots,\sigma}_{k}}.$$

Так как в правой части последнего равенства i принимает любое значение из множества $\{0, \ldots, k-1\}$, то

$$(f_1 \ldots f_i (f_{i+1} \ldots f_{i+k})_{k,\underbrace{\sigma,\ldots,\sigma}_{k}} f_{i+k+1} \ldots f_{2k-1})_{k,\underbrace{\sigma,\ldots,\sigma}_{k}} =$$

$$= (f_1 \ldots f_j (f_{j+1} \ldots f_{j+k})_{k,\underbrace{\sigma,\ldots,\sigma}_{k}} f_{j+k+1} \ldots f_{2k-1})_{k,\underbrace{\sigma,\ldots,\sigma}_{k}}$$

для любых $i, j \in \{0, \ldots, k-1\}$. Это означает ассоциативность k-арной операции $(\;)_{k,\underbrace{\sigma,\ldots,\sigma}_{k}}$. Таким образом, имеет место

**4.1.4. Следствие.** Если $\sigma^k = \sigma \in S_{n-1}$, где $k \geq 2$, то алгебра $<S_{A_1,\ldots,A_{n-1}}(\sigma),(\;)_{k,\underbrace{\sigma,\ldots,\sigma}_{k}}>$ является k-арной полугруппой.

В действительности имеет место более сильное утверждение

**4.1.5. Теорема.** Если $\sigma^k = \sigma \in S_{n-1}$, где $k \geq 2$, то алгебра $<S_{A_1,\ldots,A_{n-1}}(\sigma),(\;)_{k,\underbrace{\sigma,\ldots,\sigma}_{k}}>$ является k-арной группой.

*Доказательство.* Ввиду следствия 4.1.4, алгебра $<S_{A_1,\ldots,A_{n-1}}(\sigma),(\;)_{k,\underbrace{\sigma,\ldots,\sigma}_{k}}>$ является k-арной полугруппой. Поэтому достаточно доказать, что в $S_{A_1,\ldots,A_{n-1}}(\sigma)$ разрешимы уравнения

$$(uf_2 \ldots f_k)_{k,\underbrace{\sigma,\ldots,\sigma}_{k}} = g, \; (f_1 \ldots f_{k-1} v)_{k,\underbrace{\sigma,\ldots,\sigma}_{k}} = g, \qquad (1)$$

где

$$f_i = f_i(\sigma) = \{\sigma, f_{i1}, \ldots, f_{i(n-1)}\}, \; i = 1, \ldots, k;$$

$$g = g(\sigma) = \{\sigma, g_1, \ldots, g_{n-1}\}.$$

Покажем, что

$$u = u(\sigma) = \{\sigma, u_1, \ldots, u_{n-1}\},$$

где

$$u_j = g_j f^{-1}_{k\sigma^{k-1}(j)} \cdots f^{-1}_{3\sigma^2(j)} f^{-1}_{2\sigma(j)}, j = 1, \ldots n-1$$

является решением первого уравнения из (1). Действительно, если положить

$$(uf_2 \ldots f_k)_{k,\underbrace{\sigma,\ldots,\sigma}_{k}} = \{h_1, \ldots, h_{n-1}\},$$

то

$$h_j = u_j f_{2\sigma(j)} f_{3\sigma^2(j)} \cdots f_{k\sigma^{k-1}(j)} =$$

$$= g_j f^{-1}_{k\sigma^{k-1}(j)} \cdots f^{-1}_{3\sigma^2(j)} f^{-1}_{2\sigma(j)} f_{2\sigma(j)} f_{3\sigma^2(j)} \cdots f_{k\sigma^{k-1}(j)} =$$

$$= g_j f^{-1}_{k\sigma^{k-1}(j)} \cdots f^{-1}_{3\sigma^2(j)} f_{3\sigma^2(j)} \cdots f_{k\sigma^{k-1}(j)} = \ldots$$

$$\ldots = g_j f^{-1}_{k\sigma^{k-1}(j)} f_{k\sigma^{k-1}(j)} = g_j,$$

то есть $h_j = g_j$ для любого $j = 1, \ldots n-1$.

Докажем теперь разрешимость второго уравнения из (1). Для этого положим, что

$$v = v(\sigma) = \{\sigma, v_1, \ldots, v_{n-1}\},$$

где все $v_j$ определяются следующим образом

$$v_{\sigma^{k-1}(j)} = f^{-1}_{(k-1)\sigma^{k-2}(j)} \cdots f^{-1}_{2\sigma(j)} f^{-1}_{1j} g_j, j = 1, \ldots n-1$$

является решением этого уравнения. Действительно, если положить

$$(f_1 \ldots f_{k-1} v)_{k,\underbrace{\sigma,\ldots,\sigma}_{k}} = \{l_1, \ldots, l_{n-1}\},$$

то

$$l_j = f_{1j} f_{2\sigma(j)} \ldots f_{(k-1)\sigma^{k-2}(j)} v_{\sigma^{k-1}(j)} =$$

$$= f_{1j} f_{2\sigma(j)} \ldots f_{(k-1)\sigma^{k-2}(j)} f^{-1}_{(k-1)\sigma^{k-2}(j)} \ldots f^{-1}_{2\sigma(j)} f^{-1}_{1j} g_j = \ldots$$

$$\ldots = f_{1j} f_{2\sigma(j)} f^{-1}_{2\sigma(j)} f^{-1}_{1j} g_j = f_{1j} f^{-1}_{1j} g_j = g_j,$$

то есть $l_j = g_j$ для любого j = 1, …, n. ∎

**4.1.6. Следствие.** Если $\sigma^k = \sigma \in S_{n-1}$, где $k \geq 2$, то алгебра $< S_{A_1, \ldots, A_{n-1}}(\sigma^{-1}), (\ )_{k, \underbrace{\sigma^{-1}, \ldots, \sigma^{-1}}_{k}} >$ является k-арной группой.

Так как для циклической подстаноки $\alpha = (1\ 2\ \ldots\ n-1)$ верно равенство $\alpha^n = \alpha$, то справедливо

**4.1.7. Следствие** [3, Post; 4, Русаков; 45, Sioson]. Алгебра $< S_{A_1, \ldots, A_{n-1}}(\alpha), (\ )_{n, \underbrace{\alpha, \ldots, \alpha}_{n}} >$ является n-арной группой.

**4.1.8. Замечание.** Если $\sigma^k = \sigma \in S_{n-1}$, $k \geq 2$, то для обозначения k-арной операции $(\ )_{k, \underbrace{\sigma, \ldots, \sigma}_{k}}$ будем употреблять более экономный символ $(\ )_{k, \sigma}$. Таким образом, если

$$f_1, \ldots, f_k \in S_{A_1, \ldots, A_{n-1}}(\sigma),$$

$$(f_1 \ldots f_k)_{k, \sigma} = \{g_1, \ldots, g_{n-1}\},$$

то, согласно (*),

$$g_j = f_{1j} f_{2\sigma(j)} \ldots f_{k\sigma^{k-1}}(j).$$

**4.1.9. Замечание.** Так как тождественная подстановка $\varepsilon \in S_{n-1}$ удовлетворяет условию $\varepsilon^{k-1} = \varepsilon$, то по теореме 4.1.5 $< S_{A_1, \ldots, A_{n-1}}(\varepsilon), (\ )_{k, \varepsilon} >$ – k-арная группа с k-арной операцией

$$(f_1 \ldots f_k)_{k, \varepsilon} = \{g_1, \ldots, g_{n-1}\},$$

где

$$g_j = f_{1j}f_{2j} \ldots f_{kj}, \, j = 1, \ldots, n - 1.$$

Удостовериться в том, что $< S_{A_1, \ldots, A_{n-1}}(\varepsilon), (\,)_{k,\varepsilon} >$ k-арная группа можно и непосредственно, не используя теорему 4.1.5.

**4.1.10. Предложение.** Если $\alpha = (1\ 2\ \ldots\ n-1)$,

$$(f_1 \ldots f_n)_{n,\alpha} = \{g_1, \ldots, g_{n-1}\},$$

то

$$g_j = f_{1j}f_{2(j+1)} \ldots f_{(n-j)(n-1)}f_{(n-j+1)1} \ldots f_{(n-1)(j-1)}f_{nj}, \, j = 1, \ldots, n - 1.$$

*Доказательство.* Так как

$$\alpha(1) = 2, \alpha^2(1) = 3, \ldots \alpha^{n-2}(1) = n - 1, \alpha^{n-1}(1) = 1;$$

$$\alpha(2) = 3, \ldots, \alpha^{n-3}(2) = n - 1, \alpha^{n-2}(2) = 1, \alpha^{n-1}(2) = 2;$$

$$\ldots\ldots\ldots\ldots\ldots\ldots\ldots\ldots\ldots\ldots\ldots\ldots\ldots\ldots\ldots\ldots\ldots$$

$$\alpha(n-2) = n - 1, \alpha^2(n-2) = 1, \ldots, \alpha^{n-1}(n-2) = n - 2;$$

$$\alpha(n-1) = 1, \alpha^2(n-1) = 2, \ldots, \alpha^{n-1}(n-1) = n - 1,$$

то, полагая в (*) $m = n$, $\sigma_1 = \ldots = \sigma_n = \alpha$, получим

$$g_j = f_{1j}\, f_{2\alpha(j)}\, f_{3\alpha^2(j)} \ldots f_{(n-j)\alpha^{n-j-1}(j)}\, f_{(n-j+1)\alpha^{n-j}(j)} \ldots f_{n\alpha^{n-1}(j)},$$

откуда и из записанных выше равенств для степеней подстановки $\alpha$ следует требуемое равенство. ∎

**4.1.11. Замечание.** Так как в записи

$$g_j = f_{1j}f_{2(j+1)} \ldots f_{(n-j)(n-1)}f_{(n-j+1)1} \ldots f_{nj}, \, j = 1, \ldots, n - 1 \qquad (**)$$

подстановка $\alpha$ явно не присутствует, то положим

$$(f_1 \ldots f_n)_{n,\alpha} = (\,)_n.$$

Именно с помощью равенств (**) Пост определил n-арную операцию $(\,)_n$ для n-арных подстановок. Хотя, как мы видели выше (замечание 4.1.8), n-арная операция $(\,)_n$ может быть определена также равенством

$$g_j = f_{1j} f_{2\alpha(j)} f_{3\alpha^2(j)} \ldots f_{n\alpha^{n-1}(j)}, j = 1, \ldots, n-1 \qquad (***)$$

**4.1.12. Пример.** В $S_2$ помимо тождественной подстановки имеется еще подстановка $\alpha = (1\ 2)$, которая удовлетворяет условию $\alpha^3 = \alpha$. Поэтому по следствию 4.1.7 $< S_{A_1, A_2}(\alpha), (\ )_{3, \alpha} >$ – тернарная группа с тернарной операцией

$$(fgh)_{3, \alpha} = \{f_1 g_{\alpha(1)} h_{\alpha^2(1)}, f_2 g_{\alpha(2)} h_{\alpha^2(2)}\} = \{f_1 g_2 h_1, f_2 g_1 h_2\},$$

где

$$f = (f_1, f_2),\ g = (g_1, g_2),\ h = (h_1, h_2).$$

Если $n \geq 4$, то в $S_{n-1}$ помимо тождественной подстановки и траспозиции $\alpha = (1\ 2)$ есть и другие подстановки $\sigma$, удовлетворяющие условию $\sigma^3 = \sigma$, например, все транспозиции $\sigma_{ij} = (i\ j)$, где $i, j \in \{1, \ldots, n-1\}$.

По теореме 4.1.5 в этом случае для любой подстановки $\sigma_{ij} = (i\ j)$ алгебра $< S_{A_1, \ldots, A_{n-1}}(\sigma), (\ )_{3, (i\ j)} >$ является тернарной группой с тернарной операцией

$$(fgh)_{3, (i, j)} = \{g_1, \ldots, g_{n-1}\} = \{f_1 g_{\sigma_{ij}(1)} h_{\sigma_{ij}^2(1)}, \ldots, f_{n-1} g_{\sigma_{ij}(n-1)} h_{\sigma_{ij}^2(n-1)}\},$$

где

$$f = \{f_1, \ldots, f_{n-1}\},\ g = \{g_1, \ldots, g_{n-1}\},\ h = \{h_1, \ldots, h_{n-1}\}.$$

Запишем тернарные операции $(\ )_{3, (1\ 2)}, (\ )_{3, (1\ 3)}, (\ )_{3, (2\ 3)}$:

$$(fgh)_{3, (1, 2)} = \{f_1 g_2 h_1, f_2 g_1 h_2, f_3 g_3 h_3, \ldots, f_{n-1} g_{n-1} h_{n-1}\};$$

$$(fgh)_{3, (1, 3)} = \{f_1 g_3 h_1, f_2 g_2 h_2, f_3 g_1 h_3, f_4 g_4 h_4, \ldots, f_{n-1} g_{n-1} h_{n-1}\};$$

$$(fgh)_{3, (2, 3)} = \{f_1 g_1 h_1, f_2 g_3 h_2, f_3 g_2 h_3, f_4 g_4 h_4, \ldots, f_{n-1} g_{n-1} h_{n-1}\}.$$

В частности, при $n = 4$ последние три операции примут следующий вид:

$$(fgh)_{3, (1, 2)} = \{f_1 g_2 h_1, f_2 g_1 h_2, f_3 g_3 h_3\};$$

$$(fgh)_{3, (1, 3)} = \{f_1 g_3 h_1, f_2 g_2 h_2, f_3 g_1 h_3\};$$

$$(fgh)_{3, (2, 3)} = \{f_1 g_1 h_1, f_2 g_3 h_2, f_3 g_2 h_3\},$$

где

$$f = \{f_1, f_2, f_3\}, g = \{g_1, g_2, g_3\}, h = \{h_1, h_2, h_3\}.$$

Заметим, что в ассоциативности указанных тернарных операций можно убедиться непосредственно, проделав соответствующие вычисления.

**4.1.13. Предложение.** Если $\sigma = (n - 1 \ n - 2 \ \ldots \ 2 \ 1) \in S_{n-1}$, то алгебра $< S_{A_1, \ldots, A_{n-1}}(\sigma), (\ )_{n, \sigma} >$ – n-арная группа с n-арной операцией

$$(f_1 \ \ldots \ f_n)_{n, \sigma} = (f_{11}f_{2(n-1)} \ \ldots \ f_{(n-1)2}f_{n1},$$

$$f_{12}f_{21}f_{3(n-1)} \ \ldots \ f_{(n-1)3}f_{n2},$$

$$\ldots\ldots\ldots\ldots\ldots\ldots\ldots\ldots\ldots\ldots\ldots\ldots\ldots$$

$$f_{1(n-2)} \ \ldots \ f_{(n-2)1}f_{(n-1)(n-1)}f_{n(n-2)},$$

$$f_{1(n-1)}f_{2(n-2)} \ \ldots \ f_{(n-1)1}f_{n(n-1)})$$

или более кратко

$$(f_1 \ \ldots \ f_n)_{n, \sigma} = (g_1, \ldots, g_{n-1}),$$

где

$$g_j = f_{1j}f_{2(j-1)} \ \ldots \ f_{j1}f_{(j+1)(n-1)} \ \ldots \ f_{nj}, j = 1, \ldots, n - 1.$$

***Доказательство.*** Так $\sigma = \alpha^{-1}$, где $\alpha = (1 \ 2 \ \ldots \ n - 1)$, $\alpha^n = \alpha$, то по следствию 4.1.6 $< S_{A_1, \ldots, A_{n-1}}(\sigma), (\ )_{n, \sigma} >$ – n-арная группа.

Так как

$$\sigma(1) = n - 1, \sigma^2(1) = n - 2, \ldots \sigma^{n-2}(1) = 2, \sigma^{n-1}(1) = 1;$$

$$\sigma(2) = 1, \sigma^2(2) = n - 1, \sigma^3(2) = n - 2, \ldots \sigma^{n-2}(2) = 3, \sigma^{n-1}(2) = 2;$$

$$\ldots\ldots\ldots\ldots\ldots\ldots\ldots\ldots\ldots\ldots\ldots\ldots\ldots\ldots\ldots\ldots$$

$$\sigma(n - 2) = n - 3, \ldots, \sigma^{n-3}(n - 2) = 1, \sigma^{n-2}(n - 2) = n - 1,$$

$$\sigma^{n-1}(n-2) = n-2;$$

$$\sigma(n-1) = n-2, \sigma^2(n-1) = n-3 \ldots, \sigma^{n-2}(n-1) = 1,$$

$$\sigma^{n-1}(n-1) = n-1,$$

то, полагая в (*) $m = n$, $\sigma_1 = \ldots = \sigma_n = \sigma$, получим

$$g_j = f_{1j}f_{2(j-1)} \ldots f_{j1}f_{(j+1)(n-1)} \ldots f_{nj}, j = 1, \ldots, n-1. \quad \blacksquare$$

**4.1.14. Пример.** Если в предыдущем предложении положить $n = 4$, то $\sigma = (3\ 2\ 1)$, а операция $(\ )_{4,\sigma}$ определяется следующим образом

$$(fghu)_{4,\sigma} = (f_1g_3h_2u_1, f_2g_1h_3u_2, f_3g_2h_1u_3),$$

где

$$f = (f_1, f_2, f_3), g = (g_1, g_2, g_3), h = (h_1, h_2, h_3), u = (u_1, u_2, u_3).$$

Следующий пример показывает, что в n-арной группе $< S_{A_1, \ldots, A_{n-1}}(\alpha), (\ )_{n,\underbrace{\alpha, \ldots, \alpha}_{n}} >$ могут быть собственные n-арные подгруппы

**4.1.15. Пример.** Пусть $S_{z,w}((12)) = S_{z,w}$ – множество всех последовательностей

$$\{\alpha = (12), f_1, f_2\} = \{f_1, f_2\},$$

где $f_1$ – биекция комплексной плоскости $z$ на комплексную плоскость $w$, $f_2$ – биекция $w$ на $z$. По следствию 4.1.7 $< S_{z,w}, (\ )_{3,\alpha^2} >$ – тернарная группа. Выделим в $S_{z,w}$ подмножество $L_{z,w}$ последовательностей $\{\alpha, f_1, f_2\}$, где $f_1$ и $f_2$ – дробно-линейные преобразования. Так как существуют биекции одной комплексной плоскости на другую комплексную плоскость, отличные от дробно-линейных, то $L_{z,w}$ является собственным подмножеством тернарной группы $< S_{z,w}, (\ )_{3,\alpha^2} >$.

Пусть

$$f_1 = \{\alpha, f_{11}, f_{12}\}\ f_2 = \{\alpha, f_{21}, f_{22}\}\ f_3 = \{\alpha, f_{31}, f_{32}\}$$

– произвольные элементы из $L_{z,w}$. По определению тернарной операции $(\ )_{3,\alpha^2}$ имеем

$$(f_1 f_2 f_3)_3 = \{\alpha^3, f_{11} f_{2\alpha(1)} f_{3\alpha(\alpha(1))}, f_{12} f_{2\alpha(2)} f_{3\alpha(\alpha(2))}\} =$$

$$= \{\alpha, f_{11} f_{22} f_{3\alpha(2)}, f_{12} f_{21} f_{3\alpha(1)}\} = \{\alpha, f_{11} f_{22} f_{31}, f_{12} f_{21} f_{32}\}.$$

Так как последовательное выполнение дробно-линейных преобразований является дробно-линейным преобразованием, то $f_{11} f_{22} f_{31}$ – дробно-линейное преобразование плоскости z на плоскость w, $f_{12} f_{21} f_{32}$ – дробно-линейное преобразование плоскости w на плоскость z, и поэтому

$$(f_1 f_2 f_3)_{3,\alpha^2} \in L_{z,w},$$

то есть множество $L_{z,w}$ замкнуто относительно тернарной операции $(\ )_{3,\alpha^2}$.

Рассмотрим теперь в $L_{z,w}$ уравнение

$$(u f_2 f_3)_{3,\alpha^2} = f,$$

где

$$f_2 = \{\alpha, f_{21}, f_{22}\},\ f_3 = \{\alpha, f_{31}, f_{32}\},\ f = \{\alpha, g_1, g_2\}.$$

Покажем, что

$$u = \{\alpha,\ g_1 f_{31}^{-1} f_{22}^{-1},\ g_2 f_{32}^{-1} f_{21}^{-1}\}$$

является решением последнего уравнения. Действительно,

$$(u f_2 f_3)_3 = \{\alpha^3,\ (g_1 f_{31}^{-1} f_{22}^{-1})\ f_{22} f_{31},\ (g_2 f_{32}^{-1} f_{21}^{-1})\ f_{21} f_{32}\} = \{\alpha, g_1, g_2\} = f.$$

Аналогично показывается, что

$$v = \{\alpha,\ f_{22}^{-1} f_{11}^{-1} g_1,\ f_{21}^{-1} f_{12}^{-1} g_2\}$$

является решением уравнения

$$(f_1 f_2 v)_{3,\alpha^2} = f,\ \text{где}\ f_1 = \{\alpha, f_{12}, f_{11}\},\ f_2 = \{\alpha, f_{21}, f_{22}\},\ f = \{\alpha, g_1, g_2\}.$$

Мы показали, что $<L_{z,w},\ (\ )_{3,\alpha^2}>$ – собственная тернарная подгруппа в $<S_{z,w},\ (\ )_{3,\alpha^2}>$.

## § 4.2. n-АРНЫЕ МОРФИЗМЫ АЛГЕБРАИЧЕСКИХ СИСТЕМ

Изучению n-арных морфизмов алгебраических систем посвящен настоящий параграф. На этом пути естественно возникают n-арные полугруппы, n-арные группы, (m, n)-кольца и другие алгебраические системы.

**4.2.1**. **Определение** [83, 84]. Пусть

$$\{A_1, A_2, \ldots, A_{n-1}, A_n\} \qquad (1)$$

– последовательность однотипных универсальных алгебр. Назовем *n-арным гомомоморфизмом* последовательности (1) последовательность

$$f = \{f_1, f_2, \ldots, f_{n-1}\} \qquad (2)$$

гомоморфизмов

$$A_1 \xrightarrow{f_1} A_2 \xrightarrow{f_2} \ldots \xrightarrow{f_{n-2}} A_{n-1} \xrightarrow{f_{n-1}} A_n.$$

При n = 2 получаем понятие гомоморфизма универсальных алгебр. n-Арный гомоморфизм последовательности

$$\{A_1, A_2, \ldots, A_{n-1}, A_1\} \qquad (3)$$

называется ее *n-арным эндоморфизмом*, а n-арный эндоморфизм последовательности

$$\{\underbrace{A, A, \ldots, A, A}_{n}\} \qquad (4)$$

– *n-арным эндоморфизмом алгебры A*.

Если все $f_i$ в последовательности (2) изоморфизмы, то f – *n-арный изоморфизм*.

n-Арный изоморфизм последовательности (3) называется ее *n-арным автоморфизмом*, а n-арный автоморфизм последовательности (4) – *n-арным автоморфизмом алгебры A*.

Пусть

$$f_i = \{f_{i1}, f_{i2}, \ldots, f_{i(n-1)}\}$$

– n-арные гомоморфизмы соответственно

$$f_i : B_i \xrightarrow{f_{i1}} A_1 \xrightarrow{f_{i2}} \ldots \xrightarrow{f_{i(n-2)}} A_{n-2} \xrightarrow{f_{i(n-1)}} B_{i+1},$$

$$i = 1, \ldots, n.$$

Определим "n-арное произведение" n-арных гомоморфизмов

$$g = [f_1 f_2 \ldots f_{n-1} f_n] = \{f_{11} f_{22} \ldots f_{(n-2)(n-2)} f_{(n-1)(n-1)} f_{n1},$$

$$\ldots\ldots\ldots\ldots\ldots\ldots\ldots\ldots\ldots\ldots\ldots\ldots\ldots\ldots$$

$$f_{1k}\, f_{2(k+1)} \ldots f_{(n-k)(n-1)} f_{(n-k+1)1} \ldots f_{(n-1)(k-1)} f_{nk},$$

$$\ldots\ldots\ldots\ldots\ldots\ldots\ldots\ldots\ldots\ldots\ldots\ldots$$

$$f_{1(n-1)}\, f_{21} \ldots f_{(n-2)(n-3)}\, f_{(n-1)(n-2)}\, f_{n(n-1)}\} = \{g_1, \ldots, g_k, \ldots, g_{n-1}\}.$$

Так как произведение гомоморфизмов есть гомоморфизм, то $g_1, \ldots, g_{n-1}$ – гомоморфизмы соответственно

$$B_1 \xrightarrow{g_1} A_1 \xrightarrow{g_2} A_2 \xrightarrow{g_3} \ldots \xrightarrow{g_{n-2}} A_{n-2} \xrightarrow{g_{n-1}} B_{n+1},$$

поэтому g – n-арный гомоморфизм последовательности

$$\{B_1, A_1, \ldots, A_{n-2}, B_{n+1}\}.$$

Следовательно, имеет место

**4.2.2. Предложение** [83, 84]. n-Арное произведение n-арных гомоморфизмов является n-арным гомоморфизмом. Аналогичные утверждения имеют место для n-арных эндоморфизмов, n-арных изоморфизмов, n-арных автоморфизмов.

Обозначим через $\text{End}(A_1, A_2, \ldots, A_{n-1})$ множество всех n-арных эндоморфизмов последовательности (3), множество всех n-арных эндоморфизмов алгебры A обозначим End(n, A), то есть

$$\text{End}(n,A) = \text{End}(\underbrace{A, \ldots, A}_{n-1}).$$

В частности, End(2, A) = EndA.

**4.2.3. Теорема** [83, 84]. < End($A_1, A_2, \ldots, A_{n-1}$), [ ] > – n-арная полугруппа.

*Доказательство.* Пусть

$f_j = \{f_{ji}, f_{j2}, \ldots, f_{j(n-1)}\} \in \text{End}(A_1, A_2, \ldots, A_{n-1}), j = 1, 2, \ldots, n-1.$

Покажем, что

$$[[f_1 \ldots f_n]f_{n+1} \ldots f_{2n-1}] = [f_1^i \, [f_{i+1}^{i+n}] f_{i+n+1}^{2n-1}]$$

для любого $i = 1, \ldots, n-1$.

Положим

$$[[f_1 \ldots f_n] f_{n+1} \ldots f_{2n-1}] = \{g_1, \ldots, g_{n-1}\},$$

$$[f_1^i \, [f_{i+1}^{i+n}] f_{i+n+1}^{2n-1}] = \{h_1, \ldots, h_{n-1}\},$$

и покажем, что $g_k = h_k$, $k = 1, \ldots, n-1$.

1) $k = 1$. В этом случае

$$g_1 = (f_{11}f_{22} \ldots f_{(n-1)(n-1)}f_{n1})f_{(n+1)2} \ldots f_{(2n-2)(n-1)}f_{(2n-1)1} =$$

$$= f_{11}f_{22} \ldots f_{(n-1)(n-1)}f_{n1}f_{(n+1)2} \ldots f_{(2n-2)(n-1)}f_{(2n-1)1},$$

$$h_1 = f_{11} \ldots f_{ii}(f_{(i+1)(i+1)} \ldots f_{(n-1)(n-1)}f_{n1} \ldots f_{(i+n)(i+1)})$$

$$f_{(i+n+1)(i+2)} \ldots f_{(2n-2)(n-1)}f_{(2n-1)1} =$$

$$= f_{11}f_{22} \ldots f_{(n-1)(n-1)} f_{n1}f_{(n+1)2} \ldots f_{(2n-2)(n-1)}f_{(2n-1)1}.$$

Следовательно, $g_1 = h_1$.

2) $k = n - 1$. В этом случае

$$g_{n-1} = (f_{1(n-1)}f_{21} \ldots f_{n(n-1)})f_{(n+1)1} \ldots f_{(2n-2)(n-2)}f_{(2n-1)(n-1)} =$$

$$= f_{1(n-1)} f_{21} \ldots f_{n(n-1)}f_{(n+1)1} \ldots f_{(2n-2)(n-2)} f_{(2n-1)(n-1)},$$

$$h_{n-1} = f_{1(n-1)}f_{21} \ldots f_{i(i-1)}(f_{(i+1)i} \ldots f_{n(n-1)}f_{(n+1)1} \ldots$$

$$f_{(i+n)i})f_{(i+n+1)(i+1)} \ldots f_{(2n-2)(n-2)} f_{(2n-1)(n-1)} =$$

$$= f_{1(n-1)} f_{21} \ldots f_{n(n-1)}f_{(n+1)1}\ldots f_{(2n-2)(n-2)} f_{(2n-1)(n-1)}.$$

Следовательно, $g_{n-1} = h_{n-1}$.

3) $2 \leq k \leq n-2$. Рассмотрим три возможности.

a) $1 \leq i \leq n-k-1$. Отсюда вытекает

$$2 \leq i + 1 \leq n - k, n + 1 < n + i \leq 2n - k - 1.$$

Тогда

$$h_k = f_{1k} \ldots f_{i(i+k-1)}(f_{(i+1)(i+k)} \ldots f_{(n-k)(n-1)}f_{(n-k+1)1} \ldots$$

$$\ldots f_{nk}f_{(n+1)(k+1)} \ldots f_{(i+n)(i+k)})f_{(i+n+1)(i+k+1)} \ldots$$

$$\ldots f_{(2n-k-1)(n-1)}f_{(2n-k)1} \ldots f_{(2n-1)k} =$$

$$= f_{1k} \ldots f_{(n-k)(n-1)}f_{(n-k+1)1} \ldots f_{nk}f_{(n+1)(k+1)} \ldots$$

$$\ldots f_{(2n-k-1)(n-1)}f_{(2n-k)1} \ldots f_{(2n-1)k} =$$

$$= (f_{1k} \ldots f_{(n-k)(n-1)}f_{(n-k+1)1} \ldots f_{nk})f_{(n+1)(k+1)} \ldots$$

$$\ldots f_{(2n-k-1)(n-1)}f_{(2n-k)1} \ldots f_{(2n-1)k} = g_k.$$

Следовательно, $g_k = h_k$.

b) $n - k + 1 \leq i \leq n - 1$. Так как в этом случае

$$n - k + 2 \leq i + 1 \leq n, 2n - k + 1 \leq n + i \leq 2n - 1,$$

то

$$h_k = f_{1k} \ldots f_{(n-k)(n-1)}f_{(n-k+1)1} \ldots f_{i(i-n+k)}(f_{(i+1)(i-n+k+1)} \ldots$$

$$\ldots f_{nk}f_{(n+1)(k+1)} \ldots f_{(2n-k-1)(n-1)}f_{(2n-k)1} \ldots f_{(i+n)(i-n+k+1)})$$

$$f_{(i+n+1)(i-n+k+2)} \ldots f_{(2n-1)k} =$$

$$= (f_{1k} \ldots f_{(n-k)(n-1)} f_{(n-k+1)1} \ldots f_{nk}) f_{(n+1)(k+1)} \ldots$$

$$\ldots f_{(2n-k-1)(n-1)} f_{(2n-k)1} \ldots f_{(2n-1)k} = g_k.$$

Следовательно, $g_k = h_k$.

c) $i = n - k$. В этом случае $i + 1 = n - k + 1$, $n + i = 2n - k$,

$$h_k = f_{1k} \ldots f_{i(n-1)}(f_{(i+1)1} \ldots f_{nk} f_{(n+1)(k+1)} \ldots$$

$$\ldots f_{(2n-k-1)(n-1)} f_{(i+n)1}) f_{(i+n+1)2} \ldots f_{(2n-1)k} =$$

$$= f_{1k} \ldots f_{(n-k)(n-1)} f_{(n-k+1)1} \ldots f_{nk} f_{(n+1)(k+1)} \ldots$$

$$\ldots f_{(2n-k-1)(n-1)} f_{(2n-k)1} \ldots f_{(2n-1)k} = g_k.$$

Следовательно, $g_k = h_k$. ∎

**4.2.4. Следствие** [84]. $< \text{End}(n, A), [\,] > -$ n-арная полугруппа.

Обозначим через $\text{Aut}(A_1, A_2, \ldots, A_{n-1})$ множество всех n-арных автоморфизмов последовательности (3), а через $\text{Aut}(n, A) -$ множество всех n-арных автоморфизмов алгебры $A$, то есть

$$\text{Aut}(n, A) = \text{Aut}(\underbrace{A, \ldots, A}_{n-1}).$$

В частности, $\text{Aut}(2, A) = \text{Aut}A$.

Исторически n-арные группы возникли как обобщение понятия группы. Однако, это не единственный путь, на котором они естественно возникают, о чем свидетельствует следующая

**4.2.5. Теорема** [83, 84]. $< \text{Aut}(A_1, A_2, \ldots, A_{n-1}), [\,] > -$ n-арная группа.

***Доказательство.*** Из теоремы 4.2.3 следует, что
$< \text{Aut}(A_1, A_2, \ldots, A_{n-1}), [\,] >$ – n-арная полугруппа. Покажем, что в ней разрешимы уравнения

$$[f_1 \ldots f_{n-1} u] = g, \quad (5)$$

$$[v f_2 \ldots f_n] = g, \quad (6)$$

где

$$f_j = \{f_{j1}, f_{j2}, \ldots, f_{j(n-1)}\} \in \text{Aut}(A_1, A_2, \ldots, A_{n-1}), j = 1, 2, \ldots, n,$$

$$g = \{g_1, g_2, \ldots, g_{(n-1)}\} \in \text{Aut}(A_1, A_2, \ldots, A_{n-1}).$$

Положим

$$h_1 = f^{-1}_{(n-1)\,(n-1)} f^{-1}_{(n-2)\,(n-2)} \ldots f^{-1}_{22} f^{-1}_{11} g_1,$$

$$\ldots\ldots\ldots\ldots\ldots\ldots\ldots\ldots\ldots\ldots\ldots\ldots$$

$$h_k = f^{-1}_{(n-1)\,(k-1)} \ldots f^{-1}_{(n-k+1)1} f^{-1}_{(n-k)\,(n-1)} \ldots f^{-1}_{2(k+1)} f^{-1}_{1k} g_k,$$

$$\ldots\ldots\ldots\ldots\ldots\ldots\ldots\ldots\ldots\ldots\ldots$$

$$h_{n-1} = f^{-1}_{(n-1)\,(n-2)} f^{-1}_{(n-2)\,(n-3)} \ldots f^{-1}_{21} f^{-1}_{1\,(n-1)} g_{n-1}.$$

Ясно, что $h = \{h_1, \ldots, h_{n-1}\} \in \text{Aut}(A_1, A_2, \ldots, A_{n-1})$.
Так как

$$[f_1 \ldots f_{(n-h)} h] = \{\, f_{11}\, f_{22} \ldots f_{(n-2)\,(n-2)}\, f_{(n-1)\,(n-1)}$$

$$f^{-1}_{(n-1)\,(n-1)} f^{-1}_{(n-2)\,(n-2)} \ldots f^{-1}_{21} f^{-1}_{11}\, g_1, \ldots$$

$$\ldots, f_{1k} f_{2(k+1)} \ldots f_{(n-k)(n-1)} f_{(n-k+1)1} \ldots f_{(n-k)(k-1)}$$

$$f^{-1}_{(n-k)\,(k-1)} \ldots f^{-1}_{(n-k+1)1} f^{-1}_{(n-k)\,(n-1)} \ldots f^{-1}_{2\,(k+1)} f^{-1}_{1k}\, g_k, \ldots$$

$$\ldots, f_{1(n-1)} f_{21} \ldots f_{(n-2)(n-3)} f_{(n-1)(n-2)}$$

$$f^{-1}_{(n-1)(n-2)} f^{-1}_{(n-2)(n-3)} \ldots f^{-1}_{21} f^{-1}_{1(n-1)} g_{n-1}\} =$$

$$= \{g_1, \ldots, g_k, \ldots, g_{n-1}\} = g,$$

то $u = h$ – решение уравнения (5).

Аналогично показывается, что $v = d$ – решение уравнения (6), где

$$d = (d_1, \ldots, d_k, \ldots, d_{n-1}),$$

$$d_1 = g_1 f^{-1}_{n1} f^{-1}_{(n-1)(n-1)} f^{-1}_{(n-2)(n-2)} \ldots f^{-1}_{22},$$

$$\ldots\ldots\ldots\ldots\ldots\ldots\ldots\ldots\ldots\ldots\ldots\ldots$$

$$d_k = g_k f^{-1}_{nk} f^{-1}_{(n-1)(k-1)} \ldots f^{-1}_{(n-k+1)1} f^{-1}_{(n-k)(n-1)} \ldots f^{-1}_{2(k+1)},$$

$$\ldots\ldots\ldots\ldots\ldots\ldots\ldots\ldots\ldots\ldots\ldots\ldots$$

$$d_{n-1} = g_{n-1} f^{-1}_{n(n-1)} f^{-1}_{(n-1)(n-2)} f^{-1}_{(n-2)(n-3)} \ldots f^{-1}_{21}. \quad\blacksquare$$

**4.2.6. Следствие** [84]. $< \text{Aut}(n, A), [\ ] > -$ n-арная группа.

Напомним, что универсальная алгебра $< A, \Omega >$ называется коммутативной (определение 2.6.12), если любые две операции из $\Omega$, в том числе и совпадающие, перестановочны на A, то есть в $< A, \Omega >$ выполняется тождество

$$[(a_{11}a_{12} \ldots a_{1n})(a_{21}a_{22} \ldots a_{2n}) \ldots (a_{m1}a_{m2} \ldots a_{mn})] =$$

$$= ([a_{11}a_{21} \ldots a_{m1}][a_{12}a_{22} \ldots a_{m2}] \ldots [a_{1n}a_{2n} \ldots a_{mn}]),$$

где $(\ )$ – произвольная n-арная ($n \geq 1$) операция из $\Omega$, $[\ ]$ – произвольная m-арная ($m \geq 1$) операция из $\Omega$. Если a – нуль-арная операция, $(\ )$ – n-арная операция, то перестановочность этих операций означает выполнимость условия: $(\underbrace{a \ldots a}_{n}) = a$.

Пусть $<A_1, \Omega>, ..., <A_{n-1}, \Omega>$ – коммутативные универсальные алгебры, $(\ )_m$ – произвольная m-арная операция из $\Omega$. Определим отображение $\varphi_j$ по правилу

$$a^{\varphi_j} = (a^{f_{1j}} \ ... \ a^{f_{mj}})_m, a \in A_j, j = 1, ..., n-1,$$

где

$$\{f_{i1}, ..., f_{i(n-1)}\} = f_i \in \text{End}(A_1, ..., A_{n-1}), \quad i = 1, ..., m.$$

Если $\nu$ – нульарная операция из $\Omega$ и $\nu(A_i)$ – выделенный ею элемент из $A_i$, то n-арный эндоморфизм $\theta = \{\theta_1, ..., \theta_{n-1}\}$ определяется следующим образом:

$$a^{\theta_j} = \begin{cases} \nu(A_{j+1}), j = 1, ..., n-2 \\ \nu(A_1), \quad j = n-1, \end{cases}$$

где $a \in A_i$.

Так как универсальные алгебры $<A_1, \Omega>, ..., <A_{n-1}, \Omega>$ коммутативны, то для любой k-арной операции $(\ )_k$ и любых $a_1, ..., a_k \in A_j$ получаем

$$(a_1 \ ... \ a_k)_k^{\varphi_j} = ((a_1 \ ... \ a_k)_k^{f_{1j}} \ ... \ (a_1 \ ... \ a_k)_k^{f_{mj}})_m =$$

$$= ((a_1^{f_{1j}} \ ... \ a_k^{f_{1j}})_k \ ... \ (a_1^{f_{mj}} \ ... \ a_k^{f_{mj}})_k)_m =$$

$$= ((a_1^{f_{1j}} \ ... \ a_1^{f_{mj}})_m \ ... \ (a_k^{f_{1j}} \ ... \ a_k^{f_{mj}})_m)_k = (a_1^{\varphi_j} \ ... \ a_k^{\varphi_j})_k.$$

Если $\nu$ – нульарная операция из $\Omega$, то для любого $a \in A_i$ получаем

$$\nu(a)^{\varphi_j} = (\nu(a)^{f_{1j}} \ ... \ \nu(a)^{f_{mj}})_m = (\nu(b) \ ... \ \nu(b))_m = \nu(b),$$

где $b \in A_{i+1}$ при $j = 1, ..., n-2$; $b \in A_1$ при $j = n-1$. Мы показали, что $\varphi_j$ – гомоморфизм и, следовательно,

$$\{\varphi_1, ..., \varphi_{n-1}\} = \varphi \in \text{End}(A_1, ..., A_{n-1}).$$

Определим теперь на $\text{End}(A_1, ..., A_{n-1})$ $m$-арную операцию $(\ )_m \in \Omega$ по правилу

$$(f_1 f_2 ... f_m)_m = \varphi,$$

в частности,

$$\nu(\text{End}(A_1, ..., A_{n-1})) = \theta.$$

Из сказанного следует, что на множество $\text{End}(A_1, ..., A_{n-1})$ естественно переносятся все операции из $\Omega$, и можно рассматривать универсальную алгебру $< \text{End}(A_1, ..., A_{n-1}), \Omega >$.

**4.2.7. Теорема** [83, 84]. Если $< A_1, \Omega >, ..., < A_{n-1}, \Omega >$ – коммутативные универсальные алгебры, то коммутативной будет и универсальная алгебра $< \text{End}(A_1, ..., A_{n-1}, \Omega >$.

*Доказательство.* Пусть $(\ )_k$ и $(\ )_m$ – произвольные $k$-арная и $m$-арная операции из $\Omega$. Положим

$$(f_{i1} ... f_{ik})_k = \{g_{i1}, ..., g_{i(n-1)}\}, i = 1, ..., m,$$

$$(f_{1j} ... f_{mj})_m = \{h_{j1}, ..., h_{j(n-1)}\}, j = 1, ..., k,$$

$$((f_{11} ... f_{1k})_k ... (f_{m1} ... f_{mk})_k)_m = \{s_1, ... s_{n-1}\},$$

$$((f_{11} ... f_{m1})_m ... (f_{1k} ... f_{mk})_m)_k = \{r_1, ... r_{n-1}\},$$

где

$$f_{ij} = \{f_{ij1}, ..., f_{ij(n-1)}\}, i = 1, ..., m; j = 1, ..., k.$$

Используя коммутативность универсальных алгебр

$$< A_1, \Omega >, ..., < A_{n-1}, \Omega >,$$

получим

$$a^{s_t} = (a^{g_{1t}} ... a^{g_{mt}})_m = ((a^{f_{11t}} ... a^{f_{1kt}})_k ... (a^{f_{m1t}} ... a^{f_{mkt}})_k)_m =$$

$$= ((a^{f_{11t}} ... a^{f_{m1t}})_m ... (a^{f_{1kt}} ... a^{f_{mkt}})_m)_k = (a^{h_{1t}} ... a^{h_{kt}})_k = a^{r_t},$$

откуда $a^{s_t} = a^{r_t}$, $s_t = r_t$ для любого $t = 1, \ldots, n - 1$. Следовательно,

$$((f_{11} \ldots f_{1k})_k \ldots (f_{m1} \ldots f_{mk})_k)_m = ((f_{11} \ldots f_{m1})_m \ldots (f_{1k} \ldots f_{mk})_m)_k.$$

Пусть

$$\nu(A_j) \in A_j \ (j = 1,\ldots,n - 1), \ \theta = \{\theta_1,\ldots,\theta_{n-1}\} \in \mathrm{End}(A_1,\ldots,A_{n-1})$$

элементы, выделенные нульарной операцией $\nu$ из $\Omega$. Из коммутативности универсальных алгебр $<A_j, \Omega>$ вытекает $(\nu(A_j)\ldots\nu(A_j))_m = \nu(A_j)$ для любой m-арной операции $(\ )_m \in \Omega$, откуда $(\theta \ldots \theta)_m = \theta$. ∎

**4.2.8. Теорема** [83, 84]**.** Если универсальные алгебры $<A_1, \Omega>,\ldots, <A_{n-1}, \Omega>$ – коммутативны, то универсальная алгебра $<\mathrm{End}(A_1, \ldots, A_{n-1},) \ \Omega \cup \{[\ ]\}>$ также коммутативна и в ней выполняются тождества:

$$[[f_1 \ldots f_n]f_{n+1} \ldots f_{2n-1}] = [f_1 \ldots f_i[f_{i+1} \ldots f_{i+n}]f_{i+n+1} \ldots f_{2n-1}] \quad (7)$$

для любого $i = 1, \ldots, n - 1$;

$$[f_1 \ldots f_{i-1}(g_1 \ldots g_m)_m f_{i+1} \ldots f_n] =$$
$$= ([f_1 \ldots f_{i-1}g_1 f_{i+1} \ldots f_n] \ldots [f_1 \ldots f_{i-1}g_m f_{i+1} \ldots f_n])_m \quad (8)$$

для любого $i = 1, \ldots, n$ и любой m-арной операции $(\ )_m \in \Omega$.

*Доказательство.* Тождество (7) доказано в теореме 4.2.3. Для доказательства тождества (8) введем обозначения

$$[f_1 \ldots f_{i-1}(g_1 \ldots g_m)_m f_{i+1} \ldots f_n] = \{s_1, \ldots, s_{n-1}\},$$

$$([f_1 \ldots f_{i-1}g_1 f_{i+1} \ldots f_n] \ldots [f_1 \ldots f_{i-1}g_m f_{i+1} \ldots f_n])_m =$$
$$= \{r_1, \ldots, r_{n-1}\},$$

$$f_j = \{f_{j1}, \ldots, f_{j(n-1)}\}, \ j = 1, \ldots, n,$$

$$g_j = \{g_{j1}, \ldots, g_{j(n-1)}\} \ j = 1, \ldots, m,$$

$$[f_1 \ldots f_{i-1} g_j f_{i+1} \ldots f_n] = \{t_{j1}, \ldots, t_{j(n-1)}\} = t_j, j = 1, \ldots, m,$$

$$(g_1 \ldots g_m)_m = \{\varphi_1, \ldots, \varphi_{n-1}\}.$$

Покажем, что $s_k = r_k$, $k = 1, \ldots, n-1$. Для этого рассмотрим два случая.

1) $1 \leq i \leq n - k$. Для сокращения записей положим

$$\alpha = f_{1k} \ldots f_{(i-1)(i+k-2)}$$

$$\beta = f_{(i+1)(i+k)} \ldots f_{(n-k)(n-1)} f_{(n-k+1)1} \ldots f_{(n-1)(k-1)} f_{nk}.$$

Заметим, что при $i = n - k$

$$f_{(i+1)(i+k)} \ldots f_{(n-k)(n-1)} = \varnothing \text{ и } \beta = f_{(n-k+1)1} \ldots f_{(n-1)(k-1)} f_{nk}.$$

Так как

$$a^{s_k} = a^{\alpha \varphi_i \beta} = (a^{\alpha g_{1i}} \ldots a^{\alpha g_{mi}})_m^\beta = (a^{\alpha g_{1i} \beta} \ldots a^{\alpha g_{mi} \beta})_m,$$

$$a^{r_k} = (a^{t_{1k}} \ldots a^{t_{mk}})_m = (a^{\alpha g_{1i} \beta} \ldots a^{\alpha g_{mi} \beta})_m,$$

то $a^{s_k} = a^{r_k}$, откуда $s_k = r_k$.

2) $n - k + 1 \leq i \leq n$. Положим

$$\gamma = f_{1k} \ldots f_{(n-k)(n-1)} f_{(n-k+1)1} \ldots f_{(i-1)(i+k-n-1)},$$

$$\delta = f_{(i+1)(i+k-n+1)} \ldots f_{(n-1)(k-1)} f_{nk}.$$

Заметим, что при $i = n - k + 1$

$$f_{(n-k+1)1} \ldots f_{(i-1)(i+k-n-1)} = \varnothing \text{ и } \gamma = f_{1k} \ldots f_{(n-k)(n-1)}.$$

Так как

$$a^{s_k} = a^{\gamma \varphi_i \delta} = (a^{\gamma g_{1i}} \ldots a^{\gamma g_{mi}})_m^\delta = (a^{\gamma g_{1i} \delta} \ldots a^{\gamma g_{mi} \delta})_m,$$

$$a^{r_k} = (a^{t_{1k}} \ldots a^{t_{nk}})_m = (a^{\gamma g_{1i} \delta} \ldots a^{\gamma g_{mi} \delta})_m,$$

то $a^{s_k} = a^{r_k}$, откуда $s_k = r_k$. ∎

При n = 2 теорема 4.2.8 включает в себя результат Б.И. Плоткина о том, что множество всех эндоморфизмов коммутативной $\Omega$-алгебры является дистрибутивной $\Omega$-полугруппой ([85], теорема 6.2).

### § 4.3. ПОЛИАДИЧЕСКИЕ АНАЛОГИ ТЕОРЕМ КЭЛИ И БИРКГОФА.

Известно, что один и тот же бинарный результат может иметь несколько различных n-арных аналогов. В этом можно убедиться на примере теоремы Кэли, n-арные аналоги для которой из [45] и [19] не совпадают, хотя оба подхода основаны на использовании обычных, то есть бинарных подстановок. Между тем, при изучении n-арных групп наряду с бинарными подстановками рассматриваются (§ 4.1) и их n-арные аналоги – конечные последовательности обычных подстановок. Поэтому вполне естественной является следующая задача: получить n-арный аналог теоремы Кэли, в которой роль симметрической группы играла бы n-арная группа n-арных подстановок.

Доказанная в данном параграфе теорема решает поставленную задачу в самом общем виде. С ее помощью получен также n-арный аналог известной теоремы Биркгофа о представлении группы автоморфизмами подходящей универсальной алгебры.

Пусть $A$ – n-арная группа с n-арной операцией $f$. На множестве $A^{n-1}$ определим n-арную операцию $g$ аналогично n-арной операции, введенной Постом для n-арных подстановок:

$$g((a_1',...,a_{n-1}')\ (a_1'',...,a_{n-1}'')\ ...\ (a_1^{(n)},...,a_{n-1}^{(n)})) =$$

$$= (f(a_1'a_2''...a_{n-1}^{(n-1)}a_1^{(n)}),\ f(a_2'...a_{n-1}^{(n-2)}a_1^{(n-1)}a_2^{(n)}),\ ...,$$

$$f(a'_{n-1}a''_1 \ldots a^{(n)}_{n-1})).$$

Ясно, что декартова степень $A^{n-1}$ вместе с n-арной операцией g является n-арной группой.

В n-арной группе $A^{n-1}$ выделим подмножество

$$A_0 = \{\underbrace{(a,\ldots,a)}_{n-1} | a \in A\}.$$

Так как

$$g((\underbrace{a_1,\ldots,a_1}_{n-1}) \ldots (\underbrace{a_n,\ldots,a_n}_{n-1})) = (f(a_1 \ldots a_n), \ldots, f(a_1 \ldots a_n)) \in A_0,$$

то множество $A_0$ замкнуто относительно операции g.

**4.3.1. Предложение** [86]. $A_0$ вместе с n-арной операцией g является n-арной группой, изоморфной n-арной группе A.

*Доказательство*. Определим отображение $\alpha: A \to A_0$ по правилу $\alpha: a \mapsto \underbrace{(a,\ldots,a)}_{n-1}$. Ясно, что $\alpha$ – биекция. Так как

$$\alpha(f(a_1 \ldots a_n)) = \underbrace{(f(a_1 \ldots a_n), \ldots, f(a_1 \ldots a_n))}_{n-1} =$$

$$= g((\underbrace{a_1,\ldots,a_1}_{n-1}) \ldots (\underbrace{a_n,\ldots,a_n}_{n-1})) = g(\alpha(a_1) \ldots \alpha(a_n)),$$

то $\alpha$ является изоморфизмом n-арной группы A на n-арную группу $A_0$. ∎

Для n-арной группы A с n-арной операцией f, на множестве $A^i$ (i = 1, …, n − 1) определим отношение $\Theta_i$ по правилу:

$$(a_1,\ldots,a_i) \Theta_i (b_1,\ldots,b_i)$$

тогда и только тогда, когда

$$f(a_1 \ldots a_i x_{i+1} \ldots x_n) = f(b_1 \ldots b_i x_{i+1} \ldots x_n)$$

для любых $x_{i+1}, \ldots, x_n \in A$. Отношение $\Theta_i$ является эквивалентностью на $A^i$. Класс эквивалентности, определяемый элементом $(a_1,\ldots,a_i)$, обозначим через $\Theta_i[a_1,\ldots,a_i]$, положим также $A_i = A^i/\Theta_i$. Можно показать, что

$$A_i = \{\Theta_i[a_1,\ldots,a_{i-1},c] \mid c \in A\} = \{\Theta_i[c,a_1,\ldots,a_{i-1}] \mid c \in A\}, \quad (1)$$

где $a_1, \ldots, a_{i-1}$ – фиксированные элементы из $A$.

Пусть $\sigma = (1\ 2\ \ldots\ n-1)$ – циклическая подстановка. Для произвольного элемента $c \in A$ и любого $i = 1, \ldots, n-1$ определим отображения

$$r_{ci}: A_i \to A_{\sigma(i)}, \quad l_{ci}: A_i \to A_{\sigma(i)}$$

по формулам

$$r_{ci}(\Theta_i[a_1,\ldots,a_i]) = \Theta_{i+1}[a_1,\ldots,a_i,c], \quad i = 1, \ldots, n-2,$$

$$r_{c(n-1)}(\Theta_{n-1}[a_1,\ldots,a_{n-1}]) = \Theta_1[f(a_1\ldots a_{n-1}c)],$$

$$l_{ci}(\Theta_i[a_1,\ldots,a_i]) = \Theta_{i+1}[c,a_1,\ldots,a_i], \quad i = 1, \ldots, n-2,$$

$$l_{c(n-1)}(\Theta_{n-1}[a_1,\ldots,a_{n-1}]) = \Theta_1[f(ca_1\ldots a_{n-1})].$$

Используя (1), можно показать, что $r_{ci}$ и $l_{ci}$ – биекции.

Положим

$$r_{c_1,\ldots,c_{n-1}} = (r_{c_1 1},\ldots,r_{c_{n-1}(n-1)}), R(A) = \{r_{c_1,\ldots,c_{n-1}} \mid c_1,\ldots,c_{n-1} \in A\},$$

$$l_{c_1,\ldots,c_{n-1}} = (l_{c_1 1},\ldots,l_{c_{n-1}(n-1)}), L(A) = \{l_{c_1,\ldots,c_{n-1}} \mid c_1,\ldots,c_{n-1} \in A\}.$$

Ясно, что $r_{c_1,\ldots,c_{n-1}}$ и $l_{c_1,\ldots,c_{n-1}}$ для любых $c_1, \ldots, c_{n-1} \in A$ являются n-арными подстановками. Поэтому $R(A)$ и $L(A)$ являются подмножествами множества $S_{A_1,\ldots,A_{n-1}}$ всех n-арных подстановок последовательности $(A_1, \ldots, A_{n-1})$.

На множестве R(A) определим n-арную операцию h по правилу

$$h(r_{c_1^I, \ldots, c_{n-1}^I}, r_{c_1^{II}, \ldots, c_{n-1}^{II}} \ldots r_{c_1^{(n)}, \ldots, c_{n-1}^{(n)}}) =$$

$$= r_{f(c_1^I c_2^{II} \ldots c_{n-1}^{(n-1)}, c_1^{(n)}),\ f(c_2^I \ldots c_{n-1}^{(n-2)} c_1^{(n-1)} c_2^{(n)}),\ \ldots,\ f(c_{n-1}^I c_1^{II} \ldots c_{n-1}^{(n)})},$$

где

$$r_{c_1^{(j)}, \ldots, c_{n-1}^{(j)}} = (r_{c_1^{(j)} 1}, \ldots, r_{c_{n-1}^{(j)}(n-1)}),\ j = 1, \ldots, n.$$

**4.3.2. Предложение** [86]. Множество R(A) вместе n-арной операцией h является n-арной группой.

*Доказательство.* Легко заметить, что n-арная операция h совпадает с ассоциативной n-арной операцией, определённой Постом на множестве $S_{A_1, \ldots, A_{n-1}}$ всех n-арных подстановок последовательности $(A_1, \ldots, A_{n-1})$ произвольных множеств одинаковой мощности.

Полагая $c_1' = x_1$ – решение уравнения

$$f(x_1 c_2'' \ldots c_{n-1}^{(n-1)} c_1^{(n)}) = c_1,$$

………………………………………………………………

$c_{n-1}' = x_{n-1}$ – решение уравнения

$$f(x_{n-1} c_1'' \ldots c_{n-1}^{(n)}) = c_{n-1},$$

убеждаемся, что $u = r_{c_1^I, \ldots, c_{n-1}^I}$ является решением уравнения

$$h(u r_{c_1^{II}, \ldots, c_{n-1}^{II}} \ldots r_{c_1^{(n)}, \ldots, c_{n-1}^{(n)}}) = r_{c_1, \ldots, c_{n-1}}.$$

Полагая $c_1^{(n)} = y_1$ – решение уравнения

$$f(c_1^{'} c_2^{''} \ldots c_{n-1}^{(n-1)} y_1) = c_1,$$

..............................................................................

$c_{n-1}^{(n)} = y_{n-1}$ — решение уравнения

$$f(c_{n-1}^{'} c_1^{''} \ldots c_{n-2}^{(n-1)} y_{n-1}) = c_{n-1},$$

убеждаемся, что $v = r_{c_1^{(n)},\ldots,c_{n-1}^{(n)}}$ является решением уравнения

$$h(r_{c_1^{I},\ldots,c_{n-1}^{I}} \ldots r_{c_1^{(n-1)},\ldots,c_{n-1}^{(n-1)}} v) = r_{c_1,\ldots,c_{n-1}}. \quad\blacksquare$$

На $L(A)$ n-арная операция $h$ определяется также как на $R(A)$ и аналогично доказывается, что множество $L(A)$ вместе с n-арной операцией $h$ является n-арной группой.

**4.3.3. Определение** [86]. Элементы множества $R(A)$ назовем *правыми n-арными сдвигами* n-арной группы $A$. Элементы множества $L(A)$ назовем *левыми n-арными сдвигами* n-арной группы $A$.

Среди всех n-арных сдвигов выделим правые n-арные сдвиги вида

$$(r_{c1},\ldots,r_{c(n-1)}) = r_{\underbrace{c,\ldots,c}_{n-1}} = r_c$$

и левые n-арные сдвиги вида

$$(l_{c1},\ldots,l_{c(n-1)}) = l_{\underbrace{c,\ldots,c}_{n-1}} = l_c.$$

Положим $R_0(A) = \{r_c \mid c \in A\}$, $L_0(A) = \{l_c \mid c \in A\}$.

**4.3.4. Предложение** [86]. Множества $R_0(A)$ и $L_0(A)$ являются n-арными подгруппами соответственно n-арных групп $R(A)$ и $L(A)$.

***Доказательство***. Согласно определению n-арной операции h для любых $r_{c_1}, \ldots, r_{c_n} \in R_0(A)$ будем иметь

$$h(r_{c_1} \ldots r_{c_n}) = h(\underbrace{r_{c_1, \ldots, c_1}}_{n-1} \ldots \underbrace{r_{c_n, \ldots, c_n}}_{n-1}) =$$

$$= \underbrace{r_{f(c_1 \ldots c_n), \ldots, f(c_1 \ldots c_n)}}_{n-1} = r_{f(c_1 \ldots c_n)},$$

то есть

$$h(r_{c_1} \ldots r_{c_n}) = r_{f(c_1 \ldots c_n)}.$$

Следовательно, множество $R_0(A)$ замкнуто относительно ассоциативной n-арной операции h.

Разрешимость уравнений

$$h(ur_{c_2} \ldots r_{c_n}) = r_c, \quad h(r_{c_1} \ldots r_{c_{n-1}}v) = r_c$$

является следствием разрешимости в n-арной группе A уравнений

$$f(xc_1 \ldots c_n) = c, \quad f(c_1 \ldots c_n y) = c.$$

Мы показали, что $R_0(A)$ – n-арная подгруппа n-арной группы $R(A)$. Для $L_0(A)$ доказательство проводится аналогично. ∎

Ясно, что при n = 2 группы $R(A)$ и $R_0(A)$, соответственно группы $L(A)$ и $L_0(A)$ совпадают.

Пост заметил ([3], с.248 ), что если равномощные множества $A_1, \ldots, A_{n-1}$ попарно не пересекаются, то всякой n-арной подстановке $t = (t_1, \ldots, t_{n-1})$ последовательности $(A_1, \ldots, A_{n-1})$ можно поставить в соответствие бинарную подстановку $\tau$ множества $\bigcup_{i=1}^{n-1} A_i$, которая действует на $A_i$ также как $t_i$. Обозначив через $\widetilde{S}_{A_1, \ldots, A_{n-1}}$ множество всех таких подстановок, и определив n-арную операцию

$$\widetilde{h}(\tau_1 \tau_2 \ldots \tau_n) = \tau_1 \tau_2 \ldots \tau_n$$

для любых

$$\tau_1, \tau_2, \ldots, \tau_n \in \widetilde{S}_{A_1,\ldots,A_{n-1}},$$

можно показать, что $\widetilde{S}_{A_1, \ldots, A_{n-1}}$ вместе с n-арной операцией $\widetilde{h}$ является n-арной группой, изоморфной n-арной группе $S_{A_1, \ldots, A_{n-1}}$.

В n-арной группе $\widetilde{S}_{A_1, \ldots, A_{n-1}}$ естественно выделяются n-арные подгруппы $\widetilde{R}(A)$, $\widetilde{L}(A)$, $\widetilde{R}_0(A)$ и $\widetilde{L}_0(A)$, которые изоморфны соответственно n-арным группам $R(A)$, $L(A)$, $R_0(A)$ и $L_0(A)$.

**4.3.5. Теорема** [86]. Для всякой n-арной группы $< A, f >$ существует изоморфизм n-арной группы $< A^{n-1}, g >$ на n-арную группу $< R(A), h >$.

*Доказательство.* Определим отображение

$$\beta: A^{n-1} \to R(A)$$

по правилу

$$\beta: (c_1, \ldots, c_{n-1}) \mapsto r_{c_1, \ldots, c_{n-1}}.$$

Ясно, что $\beta$ – сюръекция. Предположим, что

$$r_{c_1, \ldots, c_{n-1}} = r_{b_1, \ldots, b_{n-1}}.$$

Это означает $r_{c_i} = r_{b_i}$ для любого $i = 1, \ldots, n-1$, то есть

$$r_{c_i}(\Theta_i[a_1,\ldots,a_i]) = r_{b_i}(\Theta_i[a_1,\ldots,a_i])$$

для любого класса эквивалентности $\Theta_i[a_1, \ldots, a_i] \in A_i$. Из последнего равенства имеем

$\Theta_{i+1}[a_1, \ldots, a_i, c_i] = \Theta_{i+1}[a_1, \ldots, a_i, b_i]$ при $i = 1, \ldots, n-2$;

$\Theta_1[f(a_1 \ldots a_{n-1} c_{n-1})] = \Theta_1[f(a_1 \ldots a_{n-1} b_{n-1})]$ при $i = n-1$;

откуда вытекает

$(a_1, \ldots, a_i, c_i)\Theta_{i+1}(a_1, \ldots, a_i, b_i)$ при $i = 1, \ldots, n-2$;

$f(a_1 \ldots a_{n-1} c_{n-1}) = f(a_1 \ldots a_{n-1} b_{n-1})$ при $i = n-1$.

Учитывая определение отношения $\Theta_i$, окончательно получаем $c_i = b_i$ для любого $i = 1, \ldots, n-1$, то есть

$$(c_1, \ldots, c_{n-1}) = (b_1, \ldots, b_{n-1}).$$

Мы показали, что $\beta$ – инъекция, а значит и биекция.

Так как

$$\beta(g((c_1', \ldots, c_{n-1}') \ldots (c_1^{(n)}, \ldots, c_{n-1}^{(n)}))) =$$

$$= \beta(f(c_1' \ldots c_{n-1}^{(n-1)} c_1^{(n)}), \ldots, f(c_{n-1}' c_1'' \ldots c_{n-1}^{(n)})) =$$

$$= r_{f(c_1' \ldots c_{n-1}^{(n-1)} c_1^{(n)}), \ldots, f(c_{n-1}' c_1'' \ldots c_{n-1}^{(n)})} =$$

$$= h(r_{c_1', \ldots, c_{n-1}'} \ldots r_{c_1^{(n)}, \ldots, c_{n-1}^{(n)}}) =$$

$$= h(\beta(c_1', \ldots, c_{n-1}') \ldots \beta(c_1^{(n)}, \ldots, c_{n-1}^{(n)})),$$

то $\beta$ – изоморфизм. ∎

**4.3.6 Замечание.** Теорема 4.3.5, как и другие n-арные аналоги групповых результатов, может быть доказана при помощи теоремы Поста о смежных классах с использованием соответствующего бинарного прототипа, в данном случае – теоремы Кэли для групп.

Ясно, что $\beta(A_0) = R_0(A)$, то есть сужение $\beta$ на $A_0$ является изоморфизмом n-арных групп $<A_0, g>$ и $<R_0(A), h>$.

**4.3.7. Следствие** [3, с.312–313, Post E.]. Всякая n-арная группа $<A, f>$ изоморфна n-арной группе $<R_0(A), h>$.

Искомый изоморфизм определяется произведением $\alpha\beta$, где $\alpha$ – изоморфизм предложения 4.3.1, $\beta$ – изоморфизм из теоремы 4.3.5.

**4.3.8. Следствие.** Для всякой n-арной группы $<A, f>$ существует изоморфизм n-арной группы $<A^{n-1}, g>$ на n-арную группу $<\widetilde{R}(A), \widetilde{h}>$.

Искомый изоморфизм равен произведению $\beta\gamma$, где $\gamma$ – изоморфизм n-арной группы $<R(A), h>$ на n-арную группу $<\widetilde{R}(A), \widetilde{h}>$.

**4.3.9. Следствие** [45, Sioson F.]. Всякая n-арная группа $<A, f>$ изоморфна n-арной группе $<\widetilde{R}_0(A), \widetilde{h}>$.

Искомый изоморфизм равен произведению $\alpha\beta\gamma$.

Теорема 4.3.5 и каждое из следствий 4.3.7 – 4.3.9 являются аналогами теоремы Кэли для полиадических групп. Нетрудно заметить, что при n = 2 эта теорема и все её следствия совпадают с теоремой Кэли.

Напомним (§4.2), что n-арным автоморфизмом последовательности $\{A_1, \ldots, A_{n-1}\}$ однотипных универсальных алгебр называется последовательность $\alpha = \{\alpha_1, \ldots, \alpha_{n-1}\}$ изоморфизмов

$$A_1 \xrightarrow{\alpha_1} A_2 \xrightarrow{\alpha_2} \ldots \xrightarrow{\alpha_{n-2}} A_{n-1} \xrightarrow{\alpha_{n-1}} A_1.$$

Множество всех n-арных автоморфизмов последовательности $\{A_1, \ldots, A_{n-1}\}$ обозначается через $\mathrm{Aut}(A_1, \ldots, A_{n-1})$.

Следующая теорема и следствия из нее являются n-арными аналогами известной теоремы Биркгофа для групп [87].

**4.3.10. Теорема** [86, 88]. Для всякой n-арной группы $< A, f >$ существует изоморфизм n-арной группы $< A^{n-1}, g >$ на n-арную группу всех n-арных автоморфизмов некоторой последовательности универсальных алгебр.

*Доказательство.* Для любого $b \in A$ определим преобразование

$$\varphi_b : \bigcup_{i=1}^{n-1} A_i \to \bigcup_{i=1}^{n-1} A_i$$

по правилу

$$\varphi_b : \Theta_i[a_1, \ldots, a_i] \mapsto \Theta_i[f(b\, b_1 \ldots b_{n-2}\, a_1), a_2, \ldots, a_i],$$

где $b_1, \ldots, b_{n-2}$ – фиксированные элементы из $A$. Положим $\Omega = \{\varphi_b \mid b \in A\}$ и рассмотрим последовательность универсальных алгебр

$$< A_1, \Omega >, \ldots, < A_{n-1}, \Omega >. \qquad (2)$$

Так как

$$r_{c_i i}(\varphi_b(\Theta_i[a_1, \ldots, a_i])) = r_{c_i i}(\Theta_i[f(bb_1 \ldots b_{n-2}a_1), a_2, \ldots, a_i]) =$$

$$= \Theta_{i+1}[f(bb_1 \ldots b_{n-2}a_1), a_2, \ldots, a_i, c_i] =$$

$$= \varphi_b(\Theta_{i+1}[a_1, \ldots, a_i, c_i]) = \varphi_b(r_{c_i i}(\Theta_i[a_1, \ldots, a_i]))$$

для любой операции $\varphi_b \in \Omega$, то $r_{c_i i}$ – изоморфизм алгебры $< A_i, \Omega >$ на алгебру $< A_{i+1}, \Omega >$ ($i = 1, \ldots, n-2$).

Аналогично доказывается, что $r_{c_{n-1}(n-1)}$ – изоморфизм алгебры $< A_{n-1}, \Omega >$ на алгебру $< A_1, \Omega >$. Следовательно, $r_{c_1, \ldots, c_{n-1}}$ – n-арный автоморфизм последовательности (2). Мы показали включение

$$R(A) \subseteq \mathrm{Aut}(A_1, \ldots, A_{n-1}). \qquad (3)$$

Пусть теперь $\delta = (\delta_1, \ldots, \delta_{n-1}) \in \mathrm{Aut}(A_1, \ldots, A_{n-1})$, то есть все $\delta_i$ – биекции и верно

$$\delta_i(\varphi_b(\Theta_i[a_1, \ldots, a_i])) = \varphi_b(\delta_i(\Theta_i[a_1, \ldots, a_i])), \ i = 1, \ldots, n-1 \quad (4)$$

для любой операции $\varphi_b \in \Omega$ и любого $\Theta_i[a_1, \ldots, a_i] \in A_i$. Последнее равенство можно переписать следующим образом:

$$\delta_i(\Theta_i[f(bb_1 \ldots b_{n-2}a_1), a_2, \ldots, a_i]) =$$
$$= \Theta_{i+1}[f(bb_1 \ldots b_{n-2}d_1), d_2, \ldots, d_{i+1}] \quad (5)$$

при $i = 1, \ldots, n-2$, где $\delta_i(\Theta_i[a_1, \ldots, a_i]) = \Theta_{i+1}[d_1, \ldots, d_{i+1}]$.

$$\delta_{n-1}(\Theta_{n-1}[f(bb_1 \ldots b_{n-2}a_1), a_2, \ldots, a_{n-1}]) = \Theta_1[f(bb_1 \ldots b_{n-2}d)] \quad (6)$$

при $i = n-1$, где $\delta_{n-1}(\Theta_{n-1}[a_1, \ldots, a_{n-1}]) = \Theta_1[d]$.

Так как (4) справедливо для любого $\Theta_i[a_1, \ldots, a_i] \in A_i$, то элемент $a_1$ можно выбрать так, что $(b_1, \ldots, b_{n-2}, a_1)$ – нейтральная последовательность. С учётом этого (5) и (6) примут соответственно вид

$$\delta_i(\Theta_i[b, a_2, \ldots, a_i]) = \Theta_{i+1}[f(bb_1 \ldots b_{n-2}d_1), d_2, \ldots, d_{i+1}], \quad (7)$$

$$\delta_{n-1}(\Theta_{n-1}[b, a_2, \ldots, a_{n-1}]) = \Theta_1[f(bb_1 \ldots b_{n-2}d)]. \quad (8)$$

В n-арной группе всегда существует элемент $c_i$ такой, что

$$(d_1, \ldots, d_{i+1})\Theta_i(a_1, \ldots, a_i, c_i),$$

и элемент $c_{n-1}$ такой, что $d = f(a_1 \ldots a_{n-1}c_{n-1})$. Теперь (7) и (8) перепишутся в виде

$$\delta_i(\Theta_i[b, a_2, \ldots, a_i]) = \Theta_{i+1}[f(bb_1 \ldots b_{n-2}a_1), a_2, \ldots, a_i, c_i],$$

$$\delta_{n-1}(\Theta_{n-1}[b, a_2, \ldots, a_{n-1}]) = \Theta_1[f(bb_1 \ldots b_{n-2}a_1), a_2, \ldots, a_{n-1}, c_{n-1}].$$

Учитывая нейтральность последовательности $(b_1, \ldots, b_{n-2}, a_1)$, из полученных равенств получаем

$$\delta_i(\Theta_i[b, a_2, \ldots, a_i]) = \Theta_{i+1}[b, a_2, \ldots, a_i, c_i], \ i = 1, \ldots, n-2,$$

$$\delta_{n-1}(\Theta_{n-1}[b, a_2, \ldots, a_{n-1}]) = \Theta_1[f(ba_2 \ldots a_{n-1}c_{n-1})].$$

Первое из полученных равенств справедливо для всех элементов множества $A_i$ ($i = 1, \ldots, n-2$), а второе справедливо для всех элементов множества $A_{n-1}$. Следовательно, отображение $\delta_i$ совпадает с $r_{c_i i}$ ($i = 1, \ldots, n-1$), а $\delta = (\delta_1, \ldots, \delta_{n-1})$ является правым n-арным сдвигом. В силу произвольного выбора $\delta \in \mathrm{Aut}(A_1, \ldots, A_{n-1})$ доказано включение

$$\mathrm{Aut}(A_1, \ldots, A_{n-1}) \subseteq R(A),$$

откуда и из (3) получаем

$$R(A) = \mathrm{Aut}(A_1, \ldots, A_{n-1}).$$

Применяя теорему 4.3.5, получаем изоморфизм $\beta$ n-арной группы $<A^{n-1}, g>$ на n-арную группу $<\mathrm{Aut}(A_1, \ldots, A_{n-1}), h>$. ∎

Так как $R_0(A) \subseteq R(A) = \mathrm{Aut}(A_1, \ldots, A_{n-1})$, и согласно следствию 4.3.7 существует изоморфизм n-арных групп $<A, f>$ и $<R_0, h>$, то справедливо

**4.3.11. Следствие.** Всякая n-арная группа изоморфна n-арной группе n-арных автоморфизмов некоторой последовательности универсальных алгебр.

**4.3.12. Следствие.** Для всякой n-арной группы $<A, f>$ существует изоморфизм n-арной группы $<A^{n-1}, g>$ на n-арную группу автоморфизмов некоторой универсальной алгебры.

*Доказательство.* Согласно следствию 4.3.8, существует изоморфизм $\mu = \beta\gamma$ n-арной группы $<A^{n-1}, g>$ на n-арную группу $<\widetilde{R}(A), \widetilde{h}>$, где

$$\beta : A^{n-1} \to R(A), \quad \beta : (c_1, \ldots c_{n-1}) \mapsto r_{c_1, \ldots, c_{n-1}};$$

$$\gamma : R(A) \to \widetilde{R}(A), \quad \gamma : r_{c_1 \ldots c_{n-1}} \mapsto r,$$

причём r – биекция множества $\bigcup_{i=1}^{n-1} A_i$, которая действует на $A_i$ также, как $r_{c_i i}$ (i = 1, …, n – 1).

Рассмотрим универсальную алгебру $< \bigcup_{i=1}^{n-1} A_i, \Omega >$. При доказательстве теоремы 4.3.10 было установлено, что все $r_{c_i i}$ – изоморфизмы. Поэтому

$$r(\varphi_b(\Theta_i[a_1, \ldots, a_i])) = r_{c_i i}(\varphi_b(\Theta_i[a_1, \ldots, a_i])) =$$
$$= \varphi_b(r_{c_i i}(\Theta_i[a_1, \ldots, a_i])) = \varphi_b(r(\Theta_i[a_1, \ldots, a_i]))$$

для любого $\Theta_i[a_1, \ldots, a_i] \in \bigcup_{i=1}^{n-1} A_i$ и любой операции $\varphi_b \in \Omega$, то есть r – автоморфизм алгебры $< \bigcup_{i=1}^{n-1} A_i, \Omega >$. ∎

Следствию 4.3.9 из теоремы 4.3.5 соответствует следующее

**4.3.13. Следствие.** Всякая n-арная группа изоморфна n-арной группе автоморфизмов некоторой универсальной алгебры.

## § 4.4. ОБОБЩЁННЫЕ МОРФИЗМЫ АБЕЛЕВЫХ m-АРНЫХ ГРУПП

n-Арные морфизмы, в том числе и n-арные эндоморфизмы последовательности

$$\{A_1, A_2, \ldots, A_{n-1}\} \qquad (1)$$

однотипных универсальных алгебр, были определены в §4.2. Там же определено n-арное произведение

$$[f_1 f_2 \ldots f_n] \qquad (2)$$

n-арных гомоморфизмов а также доказано, что множество End($A_1$, $A_2$, …, $A_{n-1}$) всех n-арных эндоморфизмов последовательности (1) образует n-арную полугруппу относительно n-арной операции (2).

В данном параграфе продолжается изучение n-арных эндоморфизмов. В нем доказывается n-арный аналог теоремы, утверждающей, что множество всех эндоморфизмов абелевой группы является кольцом.

Напомним определение (m, n)-кольца [89].

Универсальная алгебра < A, ( ), [ ] > с двумя, m-арной и n-арной операциями

$$( ) : A^m \to A, \quad [ ] : A^n \to A$$

называется *(m, n)-кольцом*, если выполняются следующие условия:

1) < A, ( ) > – абелева m-арная группа;

2) < A, [ ] > – n-арная полугруппа;

3) в < A, ( ), [ ] > для i = 1, …, n выполняется тождество

$$[a_1^{i-1}(b_1^m)a_{i+1}^n] = ([a_1^{i-1}b_1 a_{i+1}^n] \ldots [a_1^{i-1}b_m a_{i+1}^n]).$$

При m = 2 условие 3) примет вид

$$[a_1^{i-1}(b_1 + b_2)a_{i+1}^n] = [a_1^{i-1}b_1 a_{i+1}^n] + [a_1^{i-1}b_2 a_{i+1}^n],$$

где "+" – групповая операция. Поэтому (2, n) – кольцо является мультиоператорным кольцом [90], обратное не всегда верно.

Нам понадобятся следующие две легко проверяемые леммы.

**4.4.1. Лемма.** Если $a_1^i$ и $b_1^{i+k(m-1)}$ – эквивалентные последовательности m-арной группы < A, ( ) > и φ – гомоморфизм m-арной группы < A, ( ) > в m-арную группу < B, ( ) >, то последовательности

$$a_1^\varphi \ldots a_i^\varphi \text{ и } b_1^\varphi \ldots b_{i+k(m-1)}^\varphi$$

– эквивалентны в m-арной группе < B, ( ) >.

**4.4.2. Лемма.** Пусть $\varphi$ – гомоморфизм m-арной группы < A, ( ) > в m-арную группу < B, ( ) >, $a_1 \ldots a_k$ – обратная последовательность для элемента $a \in A$. Тогда

$$a^\varphi a_1^\varphi \ldots a_k^\varphi \text{ и } a_1^\varphi \ldots a_k^\varphi a^\varphi$$

– нейтральные последовательности m-арной группы < B, ( ) >.

Пусть $< A_1, (\ ) >, \ldots, < A_{n-1}, (\ ) >$ – абелевы m–арные группы. Определим отображение $\varphi_j$ по правилу

$$a^{\varphi_j} = (a^{f_{1j}} \ldots a^{f_{mj}}), \ a \in A_j, \ j = 1, \ldots, n-1,$$

где

$$\{f_{i1}, \ldots, f_{i(n-1)}\} = f_i \in \text{End}(A_1, \ldots, A_{n-1}), \ i = 1, \ldots, m.$$

В § 4.2. показано, что $\varphi_j$ – гомоморфизм и, следовательно,

$$\{\varphi_1, \ldots, \varphi_{n-1}\} = \varphi \in \text{End}(A_1, \ldots, A_{n-1}).$$

Определив на $\text{End}(A_1, \ldots, A_{n-1})$ m-арную операцию ( ) по правилу

$$(f_1 f_2 \ldots f_m) = \varphi,$$

видим, что на множество $\text{End}(A_1, \ldots, A_{n-1})$ естественно переносится m-арная операция ( ).

**4.4.3. Теорема** [83, 91]. Алгебра

$$< End(A_1, \ldots, A_{n-1}), (\ ), [\ ] >$$

является (m, n) – кольцом.

*Доказательство.* Для сокращения записей введём обозначения:

$$((f_1 f_2 \ldots f_m) f_{m+1} \ldots f_{2m-1}) = g = \{g_1, g_2, \ldots, g_{n-1}\};$$

$$(f_1 \ldots f_i (f_{i+1} \ldots f_{i+m}) f_{i+m+1} \ldots f_{2m-1}) =$$
$$= h = \{h_1, h_2, \ldots, h_{n-1}\}, \ i = 1, 2, \ldots, m-1;$$

$$(f_1 f_2 \ldots f_m) = \varphi = \{\varphi_1, \varphi_2, \ldots, \varphi_{n-1}\};$$

$$(f_{i+1} \ldots f_{i+m}) = \psi = \{\psi_1, \psi_2, \ldots, \psi_{n-1}\},$$

где

$$f_j = \{f_{j1}, f_{j2}, \ldots, f_{j(n-1)}\}, \ j = 1, 2, \ldots, 2m-1.$$

Покажем ассоциативность m-арной операции ( ). Так как

$$a^{\varphi_j} = (a^{f_{1j}} a^{f_{2j}} \ldots a^{f_{mj}}),$$

то

$$a^{g_j} = (a^{\varphi_j} a^{f_{(m+1)j}} \ldots a^{f_{(2m-1)j}}) =$$
$$= ((a^{f_{1j}} a^{f_{2j}} \ldots a^{f_{mj}}) a^{f_{(m+1)j}} \ldots a^{f_{(2m-1)j}}) =$$
$$= (a^{f_{1j}} a^{f_{2j}} \ldots a^{f_{(2m-1)j}}),$$

то есть

$$a^{g_j} = (a^{f_{1j}} a^{f_{2j}} \ldots a^{f_{(2m-1)j}}). \qquad (3)$$

где $a \in A_j, \ j = 1, \ldots, n-1.$

Так как

$$a^{\psi_j} = (a^{f_{(i+1)j}} \ldots a^{f_{(i+m)j}}),$$

то

$$a^{h_j} = (a^{f_{1j}} \ldots a^{f_{ij}} a^{\psi_j} a^{f_{(i+m+1)j}} \ldots a^{f_{(2m-1)j}}) =$$
$$= (a^{f_{1j}} \ldots a^{f_{ij}} (a^{f_{(i+1)j}} \ldots a^{f_{(i+m)j}}) a^{f_{(i+m+1)j}} \ldots a^{f_{(2m-1)j}}) =$$
$$= (a^{f_{1j}} a^{f_{2j}} \ldots a^{f_{(2m-1)j}}),$$

то есть

$$a^{h_j} = (a^{f_{1j}} a^{f_{2j}} \ldots a^{f_{(2m-1)j}}), \tag{4}$$

где $a \in A_j$, $j = 1, \ldots, n-1$.

Из (3) и (4) получаем $g_j = h_j$ ($j = 1, \ldots, n-1$), откуда $q = h$, то есть

$$((f_1^m) f_{m+1}^{2m-1}) = (f_1^i (f_{i+1}^{i+m}) f_{i+m+1}^{2m-1}), \ i = 1, \ldots, m-1.$$

Следовательно, $<\text{End}(A_1, A_2, \ldots, A_{n-1}), (\ )> - $ m-арная полугруппа.

Рассмотрим теперь уравнения

$$(f_1 f_2 \ldots f_{m-1} u) = \varphi, \tag{5}$$

$$(v f_1 f_2 \ldots f_{m-1}) = \varphi, \tag{6}$$

где

$$f_1, f_2, \ldots, f_{m-1}, \varphi \in \text{End}(A_1, A_2, \ldots, A_{n-1}),$$
$$f_i = \{f_{i1}, f_{i2}, \ldots, f_{i(n-1)}\}, \ i = 1, 2, \ldots, m-1,$$
$$\varphi = \{\varphi_1, \varphi_2, \ldots, \varphi_{n-1}\}.$$

Пусть $a_1 \ldots a_k$ – обратная последовательность для элемента $a \in A_j$. Определим отображение

$$u_j : a \to (a_1^{f_{(m-1)j}} \ldots a_k^{f_{(m-1)j}} \ldots a_1^{f_{1j}} \ldots a_k^{f_{1j}} a^{\varphi_j}).$$

Покажем, что $u_j$ – гомоморфизм. Пусть $b_i$ – произвольные элементы из $A_j$, $b_{i1} \ldots b_{ik}$ – обратная последовательность для элемента $b_i$ ($i = 1, 2, \ldots, m$), $d_1 \ldots d_k$ – обратная последовательность для элемента $(b_1 b_2 \ldots b_m)$.

Ясно, что

$$b_{m1} \ldots b_{mk} \ldots b_{11} \ldots b_{1k}$$

– обратная последовательность для последовательности $b_1 b_2 \ldots b_m$, а значит и для элемента $(b_1 b_2 \ldots b_m)$. Используя эквивалентность последовательностей (лемма 4.4.1)

$$b_{m1}^{f_{ij}} \ldots b_{mk}^{f_{ij}} \ldots b_{11}^{f_{ij}} \ldots b_{1k}^{f_{ij}}, \; d_1^{f_{ij}} \ldots d_k^{f_{ij}},$$

а также абелевость m-арных групп $A_1, A_2, \ldots, A_{n-1}$, получим

$$(b_1 b_2 \ldots b_m)^{u_j} =$$

$$= (d_1^{f_{(m-1)j}} \ldots d_k^{f_{(m-1)j}} \ldots d_1^{f_{1j}} \ldots d_k^{f_{1j}} (b_1 b_2 \ldots b_m)^{\varphi_j}) =$$

$$= (b_{m1}^{f_{(m-1)j}} \ldots b_{mk}^{f_{(m-1)j}} \ldots b_{11}^{f_{(m-1)j}} \ldots b_{1k}^{f_{(m-1)j}} \ldots$$

$$\ldots b_{m1}^{f_{1j}} \ldots b_{mk}^{f_{1j}} \ldots b_{11}^{f_{1j}} \ldots b_{1k}^{f_{1j}} b_1^{\varphi_j} b_2^{\varphi_j} \ldots b_m^{\varphi_j}) =$$

$$= ((b_{11}^{f_{(m-1)j}} \ldots b_{1k}^{f_{(m-1)j}} \ldots b_{11}^{f_{1j}} \ldots b_{1k}^{f_{1j}} b_1^{\varphi_j}) \ldots$$

$$\ldots (b_{m1}^{f_{(m-1)j}} \ldots b_{mk}^{f_{(m-1)j}} \ldots b_{m1}^{f_{1j}} \ldots b_{mk}^{f_{1j}} b_m^{\varphi_j})) = (b_1^{u_j} \ldots b_m^{u_j}),$$

откуда

$$(b_1 \ldots b_m)^{u_j} = (b_1^{u_j} \ldots b_m^{u_j}).$$

Следовательно, $u_j$ – гомоморфизм ($j = 1, 2, \ldots, n-1$), и $u = \{u_1, \ldots, u_{n-1}\} \in \mathrm{End}\,(A_1, \ldots, A_{n-1})$. Теперь применяя лемму 4.4.2, получим

$$(a^{f_{1j}} \ldots a^{f_{(m-1)j}} a^{u_j}) =$$

$$= (a^{f_{1j}} \ldots a^{f_{(m-1)j}} (a_1^{f_{(m-1)j}} \ldots a_k^{f_{(m-1)j}} \ldots a_1^{f_{1j}} \ldots a_k^{f_{1j}} a^{\varphi_j})) =$$

$$= (\underbrace{a^{f_{1j}} a_1^{f_{1j}} \ldots a_k^{f_{1j}}}_{\text{нейтр.}} \ldots \underbrace{a^{f_{(m-1)j}} a_1^{f_{(m-1)j}} \ldots a_k^{f_{(m-)j}}}_{\text{нейтр.}} a^{\varphi_j}) = a^{\varphi_j},$$

откуда

$$(a^{f_{1j}} \ldots a^{f_{(m-1)j}} a^{u_j}) = a^{\varphi_j}.$$

Следовательно, $u = \{u_1, u_2, \ldots, u_{n-1}\}$ – решение уравнения (5). Аналогично показывается, что $v = \{v_1, v_2, \ldots, v_{n-1}\}$ – решение уравнения (6), где гомоморфизмы $v_j$ ($j = 1, 2, \ldots, n-1$} определены по правилу

$$v_j : a \to (a^{\varphi_j} a_1^{f_{(m-1)j}} \ldots a_k^{f_{(m-1)j}} \ldots a_1^{f_{1j}} \ldots a_k^{f_{1j}}).$$

Этим показано, что $<\text{End}(A_1, A_2, \ldots, A_{n-1}), (\ )>$, – m-арная группа. Абелевость этой m-арной группы вытекает из абелевости m-арных групп $A_1, A_2, \ldots, A_{n-1}$.

Осталось показать

$$[f_1^{i-1}(g_1^m)f_{i+1}^n] = ([f_1^{i-1}g_1 f_{i+1}^n] \ldots [f_1^{i-1}g_m f_{i+1}^n]), \ i = 1, \ldots, n.$$

Введём следующие обозначения

$$[f_1^{i-1}(g_1^m)f_{i+1}^n] = \{s_1, s_2, \ldots, s_{n-1}\},$$

$$([f_1^{i-1}g_1 f_{i+1}^n] \ldots [f_1^{i-1}g_m f_{i+1}^n]) = \{r_1, r_2, \ldots, r_{n-1}\},$$

$$f_j = \{f_{j1}, f_{j2}, \ldots, f_{j(n-1)}\}, \ j = 1, \ldots, n,$$

$$g_j = \{g_{j1}, g_{j2}, \ldots, g_{j(n-1)}\}, \ j = 1, \ldots, m,$$

$$[f_1^{i-1}g_j f_{i+1}^n] = \{t_{j1}, t_{j2}, \ldots, t_{j(n-1)}\}, \ j = 1, \ldots, m,$$

$$(g_1^m) = \{\varphi_1, \varphi_2, \ldots, \varphi_{n-1}\}.$$

Покажем, что $s_k = r_k$, $k = 1, \ldots, n-1$.
Рассмотрим возможные случаи.
1) $1 \le i \le n - k$. Так как

$$a^{s_k} = a^{f_{1k}\ldots f_{(i-1)(i+k-2)}\varphi_i f_{(i+1)(i+k)}\ldots f_{(n-k)(n-1)}f_{(n-k+1)1}\ldots f_{(n-1)(k-1)}f_{nk}} =$$

$$= (a^{f_{1k}\ldots f_{(i-1)(i+k-2)}g_{1i}} \ldots a^{f_{1k}\ldots f_{(i-1)(i+k-2)}g_{mi}})^{f_{(i+1)(i+k)}\ldots}$$

$$\ldots f_{(n-k)(n-1)}f_{(n-k+1)1}\ldots f_{(n-k)(k-1)}f_{nk} =$$

$$= (a^{f_{1k}\ldots f_{(i-1)(i+k-2)}g_{1i}f_{(i+1)(i+k)}\ldots f_{(n-k)(n-1)}f_{(n-k+1)1}\ldots f_{(n-1)(k-1)}f_{nk}} \ldots$$

$$\ldots a^{f_{1k}\ldots f_{(i-1)(i+k-2)}g_{mi}f_{(i+1)(i+k)}\ldots f_{(n-k)(n-1)}f_{(n-k+1)1}\ldots f_{(n-1)(k-1)}f_{nk}}).$$

$$a^{r_k} = (a^{t_{1k}} \ldots a^{t_{mk}}) =$$

$$= (a^{f_{1k}\ldots f_{(i-1)(i+k-2)}g_{1i}f_{(i+1)(i+k)}\ldots f_{(n-k)(n-1)}f_{(n-k+1)1}\ldots f_{(n-1)(k-1)}f_{nk}} \ldots$$

$$\ldots a^{f_{1k}\ldots f_{(i-1)(i+k-2)}g_{mi}f_{(i+1)(i+k)}\ldots f_{(n-k)(n-1)}f_{(n-k+1)1}\ldots f_{(n-1)(k-1)}f_{nk}}),$$

то $a^{s_k} = a^{r_k}$, откуда $s_k = r_k$.
2) $n - k < i \le n$. Так как

$$a^{s_k} = a^{f_{1k}\ldots f_{(n-k)(n-1)}f_{(n-k+1)1}\ldots f_{(i-1)(i+k-n-1)}\varphi_i f_{(i+1)(i+k-n+1)}\ldots f_{(n-1)(k-1)}f_{nk}} =$$

$$= (a^{f_{1k}\ldots f_{(n-k)(n-1)}f_{(n-k+1)1}\ldots f_{(i-1)(i+k-n-1)}g_{1i}} \ldots$$

$$\ldots a^{f_{1k}\ldots f_{(n-1)(k-1)}f_{(n-k+1)1}\ldots f_{(i-1)(i+k-n-1)}g_{mi}})^{f_{(i+1)(i+k-n+1)}\ldots f_{(n-1)(k-1)}f_{nk}} =$$

$$= (a^{f_{1k}\ldots f_{(n-k)(n-1)}f_{(n-k+1)1}\ldots f_{(i-1)(i+k-n-1)}g_{1i}f_{(i+1)(i+k-n+1)}\ldots f_{(n-1)(k-1)}f_{nk}} \ldots$$

$$\ldots a^{f_{1k}\ldots f_{(n-k)(n-1)}f_{(n-k+1)1}\ldots f_{(i-1)(i+k-n-1)}g_{mi}f_{(i+1)(i+k-n+1)}\ldots f_{(n-1)(k-1)}f_{nk}}),$$

$$a^{r_k} = (a^{t_{1k}} \ldots a^{t_{nk}}) =$$

$$= (a^{f_{1k}\ldots f_{(n-k)(n-1)}f_{(n-k+1)1}\ldots f_{(i-1)(i+k-n-1)}g_{1i}f_{(i+1)(i+k-n+1)}\ldots f_{(n-1)(k-1)}f_{nk}} \ldots$$

$$\ldots a^{f_{1k}\ldots f_{(n-k)(n-1)}f_{(n-k+1)1}\ldots f_{(i-1)(i+k-n-1)}g_{mi}f_{(i+1)(i+k-n+1)}\ldots f_{(n-1)(k-1)}f_{nk}}),$$

то $a^{s_k} = a^{r_k}$, откуда $s_k = r_k$. ∎

**4.4.4. Следствие.** Если

$$< A_1, \{+, -, 0\} >, \ldots, < A_{n-1}, \{+, -, 0\} >$$

– абелевы группы, то

$$< \text{End}(A_1, \ldots, A_{n-1}), \{+, -, \Theta, [\,] \} >$$

– мультиоператорное кольцо, где $\Theta = \{\underbrace{0,\ldots,0}_{n-1}\}$.

**4.4.5. Следствие** [36, Dudek W.]. Множество всех эндоморфизмов абелевой m-арной группы является (m, 2)-кольцом.

## ДОПОЛНЕНИЯ И КОММЕНТАРИИ

**1.** Последовательности $\{f_1, f_2, \ldots, f_{n-1}\}$ биекций

$$A_1 \xrightarrow{f_1} A_2 \xrightarrow{f_2} \ldots \xrightarrow{f_{n-2}} A_{n-1} \xrightarrow{f_{n-1}} A_1,$$

где $A_1, \ldots, A_{n-1}$ – конечные множества, первым начал изучать Пост [3], называя такие последовательности n-арными подстановками. На множестве всех определенных таким образом n-арных подстановок Пост определил n-арную операцию ( ) и показал, что вместе с этой n-арной операцией множество всех введенных им n-арных подстановок является n-арной группой, которую он назвал n-арной симметрической группой степени k, где k – мощность множеств $A_1, \ldots, A_{n-1}$. В наших обозначениях (§ 4.2) такая n-арная группа обозначается символом

$$< S_{A_1, \ldots, A_{n-1}}(\alpha), (\,) >,$$

где $\alpha = (12 \ldots n - 1)$ – циклическая подстановка. Пост показал, что n-арная группа $< S_{A_1, \ldots, A_{n-1}}(\alpha), (\ ) >$ имеет порядок $(k!)^{n-1}$ и содержит $(k!)^{n-2}$ идемпотентов. Ясно, что при $n = 2$ получается порядок $k!$ симметрической группы степени $k$ и число 1 ее идемпотентов. Отметим, что в [3] Пост получил n-арные аналоги многих результатов о группах подстановок, известных на момент написания его работы.

**2.** Более общим, чем n-арная подстановка является понятие последовательности отображений множеств, введенное С.А. Русаковым [4] следующим образом.

**Определение** Пусть

$$X = \{ X_1, X_2, \ldots, X_k \}, Y = \{ Y_1, Y_2, \ldots, Y_k \}$$

– последовательности длины $k \geq 1$, составленные из произвольных непустых множеств, и пусть $\delta$ некоторая подстановка из $S_k$. Если для каждого $i = 1, 2, \ldots, k$ определено отображение

$$f_i : X_i \to Y_{\delta(i)},$$

то последовательность

$$f = \{ f_1, f_2, \ldots, f_k \}$$

называется последовательностью отображений из $X$ в $Y$, определенной перестановкой $\delta$. Если при этом для каждого $i = 1, 2, \ldots, k$ отображение $f_i$ является биекцией $X_i$ на $Y_{\delta(i)}$, то $f$ называется последовательностью биективных отображений $X$ на $Y$, определенной подстановкой $\delta$.

Полагая в приведенном определении $X = Y = \{ A_1, \ldots, A_{n-1}, A_1\}$, получим определение n-арной подстановки из $S_{A_1, \ldots, A_{n-1}}(\delta)$.

**3.** n-Арные морфизмы алгебраических систем впервые были определены в [83] по аналогии с n-арными подстановками Поста. Изучению n-арных морфизмов посвящена работа автора [84].

**4.** Понятие n-арного эндоморфизма можно расширить, если воспользоваться конструкцией из §4.1, с помощью которой определялись n-арные подстановки. А именно, если

$$\{ A_1, A_2, \ldots, A_{n-1} \}$$

последовательность однотипных универсальных алгебр, то для всякой подстановки $\sigma \in S_{n-1}$ определим множество $\text{End}(\sigma, A_1, \ldots, A_{n-1})$ всех последовательностей

$$\{\sigma, f_1, \ldots, f_{n-1}\}, \text{ где } f_j : A_j \to A_{\sigma(j)}$$

гомоморфизмы, $j = 1, \ldots, n - 1$. Такие последовательности будем называть n-арными эндоморфизмами. Если $\sigma = (1 2 \ldots n - 1)$ – циклическая подстановка, то определенные таким образом n-арные эндоморфизмы совпадают с n-арными эндоморфизмами из §4.2, и, кроме того,

$$\text{End}(\sigma, A_1, \ldots, A_{n-1}) = \text{End}(A_1, \ldots, A_{n-1}).$$

Аналогично для всякой подстановки $\sigma \in S_{n-1}$ определяется множество $\text{Aut}(\sigma, A_1, \ldots, A_{n-1})$ n-арных автоморфизмов, совпадающее при $\sigma = (1\, 2 \ldots n - 1)$ с множеством $\text{Aut}(A_1, \ldots, A_{n-1})$.

Для всякого подмножества $T \subseteq S_{n-1}$ полагаем

$$\text{End}(T, A_1, \ldots, A_{n-1}) = \bigcup_{\sigma \in T} \text{End}(\sigma, A_1, \ldots, A_{n-1}),$$

$$\text{Aut}(T, A_1, \ldots, A_{n-1}) = \bigcup_{\sigma \in T} \text{Aut}(\sigma, A_1, \ldots, A_{n-1}).$$

**5.** Первый n-арный аналог теоремы Кэли был получен Постом [3] (следствие 4.3.7.).

**6.** Следующая теорема из [19] является n-арным аналогом теоремы Кэли, при получении которого использовались трансляции n-арных групп. Используемые в ее доказательстве обозначения имеются в дополнениях и комментариях к главе 2.

**Теорема 1** [19]. Для всякой n-арной группы $< A, [\ ] >$ существует её гомоморфизм на факторалгебру некоторой n-арной группы подстановок на $A$.

*Доказательство*. Если $\bigcup$ – n-арная операция, производная от операции в группе $T(A)$, то легко проверяется, что $< T_{n-1}(A), \bigcup > -$ n-арная группа. Определим на $T_{n-1}(A)$ отношение $\pi$ следующим образом: $(u, v) \in \pi$ тогда и только тогда, когда

$$u = t(a_1^{n-2}, a), v = t(b_1^{n-2}, a)$$

для некоторых $a_1, \ldots, a_{n-2}, b_1, \ldots, b_{n-2}, a \in A$. Если $n = 2$, то $\pi$ – тривиальная нулевая конгруэнция. Можно показать, что при $n \geq 3$ $\pi$ – конгруэнция на $<T_{n-1}(A), \cup>$. Тогда отображение $\gamma\colon A \to T_{n-1}(A)/\pi$ по правилу $\gamma\colon a \to T_{n-1}(A, a)$, где

$$T_{n-1}(A, a) = \{t(b_1^{n-2}, a) \mid b_1, \ldots, b_{n-2} \in A\}$$

является гомоморфизмом $<A, [\ ]>$ на фактор алгебру $<T_{n-1}(A)/\pi, [\ ]_\pi>$ с n-арной операцией

$$[T_{n-1}(A, a_1) \ldots T_{n-1}(A, a_n)]_\pi = T_{n-1}(A, [a_1 \ldots a_n])$$

(при $n = 2$ полагаем $T_{n-1}(A, a) = \{t(, a)\}$). ■

Если в последней теореме $n = 2$, то $A$ – группа, $T_{n-1}(A) = R(A)$ – множество всех правых сдвигов группы $A$, $\pi$ –тривиальная конгруэнция, $\gamma$ – изоморфизм групп $A$ и $R(A)$. Следовательно, доказанная теорема является аналогом теоремы Кэли для групп.

В [19] показано, что для любой n-арной группы $<A, [\ ]>$ группа автоморфизмов алгебры $<A, T_n(A)>$ совпадает с группой $T_1(A)$, то есть

$$\mathrm{Aut}<A, T_n(A)> = T_1(A).$$

Там же, с использованием теоремы Глускина-Хоссу, показано, что на $T_1(A)$ можно определить n-арную операцию $(\ )$, так, что n-арные группы $<A, [\ ]>$ и $<T_1(A), (\ )>$ изоморфны, то есть имеет место

**Теорема 2** [19]. Каждая n-арная группа изоморфна n-арной группе всех автоморфизмов некоторой универсальной алгебры.

**Следствие** [87]. Каждая группа изоморфна группе всех автоморфизмов некоторой универсальной алгебры.

**7.** Изучением (m, n)-колец занимались Чупона Г. [92], Crombez G. [89, 93], Timm J. [93], Dudek W.[94]. (2, n)-Кольца изучали Никитин А.Н. [95 – 97] и Артамонов В.А. [97]. Celakoski N. изучал (F, G)-кольца [98], являющиеся обобщением (m, n)-колец. (M, N)-Кольца и (M, N)-полукольца изучала Кондратова-Суворова А.Д. [99 – 102]. Полиадическим мультикольцам посвящены работы Кравченво Ю.В. [103 – 105] и Новикова С.П. [105].

# Г Л А В А  5

# n-АРНЫЕ ПОДГРУППЫ

В данной главе продолжено изучение n-арных аналогов нормальных и сопряженных подгрупп в группе. В частности, определяются и изучаются $\Sigma$-нормальные n-арные подгруппы. Приведено большое число критериев существования n-арных подгрупп в n-арной группе.

## §5.1. КРИТЕРИИ СУЩЕСТВОВАНИЯ n-АРНЫХ ПОДГРУПП

**5.1.1. Теорема** [12]. Для того, чтобы подалгебра $< B, [\ ] >$ n-арной ($n \geq 3$) группы $< A, [\ ] >$ была её n-арной подгруппой необходимо и достаточно, чтобы для любого $b \in B$ существовала последовательность $\alpha(b)$ длины $n - 2$, составленная из элементов множества B, такая, что

$$[b\alpha(b)b] = b. \qquad (1)$$

*Доказательство*. *Необходимость*. Если $< B, [\ ] >$ – n-арная ($n \geq 3$) подгруппа в $< A, [\ ] >$, то для любого $b \in B$ существует обратная последовательность $\alpha(b)$ длины $n - 2$, составленная из элементов множества B. Тогда последовательность $b\alpha(b)$ нейтральная, и поэтому верно (1).

*Достаточность*. Из (1) вытекает, что последовательности $b\alpha(b)$ и $\alpha(b)b$ – нейтральные в $< A, [\ ] >$. Поэтому

$$[c\alpha(b)b] = c = (b\alpha(b)c]$$

для любого $c \in A$. Применяя теорему 3.5.1, заключаем, что $< B, [\ ] >$ n-арная группа. ∎

В дальнейшем, чтобы каждый раз не повторять, что некоторая последовательность α составлена из элементов множества B, будем писать α ∈ $F_B$, где $F_B$ – свободная полугруппа над алфавитом B.

**5.1.2. Следствие** [12]. Для того, чтобы подалгебра < B, [ ] > n-арной ( n ≥ 3) группы < A, [ ] > была её n-арной подгруппой, необходимо и достаточно, чтобы для любого b ∈ B существовала обратная последовательность α(b) ∈ $F_B$ длины n − 2.

**5.1.3. Следствие**. Для того, чтобы подалгебра < B, [ ] > тернарной группы < A, [ ] > была её тернарной подгруппой, необходимо и достаточно, чтобы для любого b ∈ B существовал элемент $b^{-1}$ ∈ B такой, что

$$[bb^{-1}b] = b.$$

**5.1.4. Следствие**. Подалгебра < B, [ ] > тернарной группы < A, [ ] > является её тернарной подгруппой тогда и только тогда, когда она вместе со всяким элементом содержит и его обратный.

**5.1.5. Следствие. [**1, Дёрнте]. Для того, чтобы подалгебра < B, [ ] > n-арной группы < A, [ ] > была её n-арной подгруппой (n ≥ 3), необходимо и достаточно, чтобы множество B вместе со всяким своим элементом b содержало и косой элемент $\bar{b}$.

*Доказательство*. Необходимость является следствием определения косого элемента.

*Достаточность*. Если для любого b ∈ B верно $\bar{b}$ ∈ B, то последовательность

$$\alpha(b) = \bar{b}\underbrace{b\ldots b}_{n-3}$$

составлена из элементов множества B и ввиду предложения 1.2.22, является обратной к элементу b. Поэтому

$$[c\,\overline{b}\,\underbrace{b\ldots b}_{n-3}\,b] = c = [b\,\overline{b}\,\underbrace{b\ldots b}_{n-3}\,c],$$

то есть

$$[c\alpha(b)b] = c = [b\alpha(b)c]$$

для любого $c \in B$. Применяя теперь теорему 5.1.1, заключаем, что $< B, [\ ] >$ – n-арная подгруппа в $< A, [\ ] >$. ∎

Для непустого множества A положим

$$F_{A,n-1} = \{\alpha \in F_A \mid 1 \leq \ell(\alpha) \leq n-1\}$$

В частности,

$$F_{A,1} = \{\alpha \in F_A \mid \ell(\alpha) = 1\} = A.$$

**5.1.6. Теорема** [12]. Пусть B – непустое подмножество n-арной ($n \geq 3$) группы $< A, [\ ] >$. Если для любых $\alpha, \beta \in F_{B,n-1}$ и некоторой обратной к $\alpha$ последовательности $\alpha^{-1} \in F_{A,n-1}$ последовательность $\alpha^{-1}\beta(\beta\alpha^{-1})$ эквивалентна некоторой последовательности из $F_{B,n-1}$, то $< B, [\ ] >$ – n-арная подгруппа n-арной группы $< A, [\ ] >$.

***Доказательство***. Пусть $\alpha = \beta \in F_{B,n-1}$, $\alpha^{-1} \in F_{A,n-1}$. Так как $\alpha^{-1}\beta = \alpha^{-1}\alpha$ – нейтральная и, кроме того, по условию эквивалентна некоторой последовательности из $F_{B,n-1}$, то существует нейтральная последовательность $\varepsilon \in F_{B,n-1}$. Полагая $\alpha \in F_{B,n-1}$, $\beta = \varepsilon$, заключаем, что $\alpha^{-1}\beta = \alpha^{-1}\varepsilon \sim \alpha^{-1}$ – эквивалентна некоторой последовательности из $F_{B,n-1}$, которая также является обратной к $\alpha$. Таким образом, мы показали, что для любой $\alpha \in F_{B,n-1}$ существует обратная последовательность из $F_{B,n-1}$.

Покажем теперь замкнутость n-арной операции $[\ ]$ на множестве B. Если $b_1, b_2, \ldots, b_n \in B$, то по доказанному существует обратная к $b_1$ последовательность $b_1^{-1} \in F_{B,n-1}$. Поэтому, полагая $\beta = b_2\ldots b_n$, и, учитывая, что любая обратная к $b_1^{-1}$ последовательность из $F_{A,n-1}$ совпадает с $b_1$, то есть $(b_1^{-1})^{-1} = b_1$, получим эквивалентность последовательности

$$(b_1^{-1})^{-1}b_2\ldots b_n = b_1 b_2 \ldots b_n$$

некоторой последовательности $\gamma \in F_{B,n-1}$. А так как $\gamma \sim b_1 b_2 \ldots b_n$ и $1 \leq \ell(\gamma) \leq n-1$, то $\ell(\gamma) = 1$, то есть $\gamma = c$ для некоторого $c \in B$. Тогда из $b_1 b_2 \ldots b_n \sim c$ следует

$$[b_1 b_2 \ldots b_n] = c \in B.$$

Если теперь $b$ – произвольный элемент из $B$, то по доказанному существует обратная к $b$ последовательность $\alpha(b) = b^{-1} \in F_{B,n-1}$. Так как $b\alpha(b)$ и $\alpha(b)b$ – нейтральные последовательности, то

$$[b\alpha(b)b] = b$$

для любого $b \in b$. Применяя теорему 5.1.1, заключаем, что $< B, [\ ] >$ – n-арная подгруппа n-арной группы $< A, [\ ] >$.

Для последовательности $\beta\alpha^{-1}$ доказательство проводится аналогично. ∎

**5.1.7. Замечание**. Так как при n = 2, $F_{A,n-1} = A$, $F_{B,n-1} = B$, то теорема 5.1.6 формально включает в себя соответствующий бинарный результат.

Теорему 5.1.6 можно усилить.

**5.1.8. Предложение** [12]. Пусть B – непустое подмножество n-арной (n ≥ 3) группы $< A, [\ ] >$. Если для любых $a \in B$, $\beta \in F_{B,n-1}$ и некоторой обратной к $a$ последовательности $a_1 \ldots a_{n-2}$, последовательность $a_1 \ldots a_{n-2}\beta$ ($\beta a_1 \ldots a_{n-2}$) эквивалентна некоторой последовательности из $F_{B,n-1}$, то $< B, [\ ] >$ – n-арная подгруппа n-арной группы $< A, [\ ] >$.

*Доказательство*. Полагаем $\alpha = a$ и дословно повторяем доказательства теоремы 5.1.6. ∎

**5.1.9 Следствие** [12]. Пусть B – непустое подмножество n-арной ($n \geq 3$) группы $<A, [\ ]>$. Если для любых $a \in B$, $\beta \in F_{B,n-1}$ и некоторого $i \in \{0, 1, \ldots, n-3\}$ последовательность

$$\underbrace{a\ldots a}_{i}\,\overline{a}\,\underbrace{a\ldots a}_{n-i-3}\beta \quad (\beta\underbrace{a\ldots a}_{i}\,\overline{a}\,\underbrace{a\ldots a}_{n-i-3})$$

эквивалентна некоторой последовательности из $F_{B,n-1}$, то $<B, [\ ]>$ – n-арная подгруппа n-арной группы $<A, [\ ]>$.

**5.1.10. Теорема** [12]. Пусть $<B, [\ ]>$ – подалгебра n-арной ($n \geq 3$) группы $<A, [\ ]>$. Если существует такая $\beta \in F_{B,n-1}$, что для любой $\alpha \in F_{B,n-1}$ и некоторой обратной к $\alpha$ последовательности $\alpha^{-1} \in F_{A,n-1}$, последовательность $\alpha^{-1}\beta$ ($\beta\alpha^{-1}$) эквивалентна некоторой последовательности из $F_{B,n-1}$, то $<B, [\ ]>$ – n-арная подгруппа n-арной группы $<A, [\ ]>$.

*Доказательство*. Так как $<B, [\ ]>$ – подалгебра и $\beta$ составлена из элементов множества B, то последовательность $\beta\beta$ эквивалентна некоторой последовательности $\gamma \in F_{B,n-1}$. В частности, если верно неравенство $\ell(\beta\beta) \leq n-1$, то $\beta\beta = \gamma$. По условию, $\gamma^{-1}\beta \sim \widetilde{\beta} \in F_{B,n-1}$ для некоторой обратной к $\gamma$ последовательности $\gamma^{-1} \in F_{A,n-1}$, которая всегда существует. С другой стороны,

$$\gamma^{-1}\beta \sim (\beta\beta)^{-1}\beta \sim \beta^{-1}\beta^{-1}\beta \sim \beta^{-1}$$

для любой обратной $\beta^{-1}$ к последовательности $\beta$, откуда, учитывая транзитивность отношения эквивалентности, получаем $\beta^{-1} \sim \widetilde{\beta}$. Следовательно, существует обратная к $\beta$ последовательность $\widetilde{\beta} \in F_{B,n-1}$.

Если теперь $\alpha^{-1}\beta \sim \delta \in F_{B,n-1}$, то $\alpha^{-1}\beta\widetilde{\beta} \sim \delta\widetilde{\beta}$, откуда, учитывая $\alpha^{-1}\beta\widetilde{\beta} \sim \alpha^{-1}$, получаем $\alpha^{-1} \sim \delta\widetilde{\beta}$. А так как $<B, [\ ]>$ – подалгебра и последовательности $\delta$ и $\widetilde{\beta}$ составлены из элементов множества B, то последовательность $\delta\widetilde{\beta}$ эквивалентна некоторой последовательности $\rho \in F_{B,n-1}$, откуда $\alpha^{-1} \sim \rho \in F_{B,n-1}$, и ρ – обратная для α. Следовательно, для любой α ∈ $F_{B,n-1}$ существует обратная из $F_{B,n-1}$. В частности, это верно для любого b ∈ B, т. е. для любого b ∈ B существует обратная последовательность α(b) длины n – 2, составленная из элементов множества B. Поэтому по следствию 5.1.2 < B, [ ] > – n-арная подгруппа n-арной группы < A, [ ] >.

Для последовательности βα$^{-1}$ доказательство проводится аналогично. ∎

**5.1.11. Следствие** [12]. Пусть < B, [ ] > – подалгебра n-арной (n ≥ 3) группы < A, [ ] >. Если существует такая β ∈ $F_{B,n-1}$, что для любой α ∈ $F_{B,n-1}$ и некоторой обратной к α последовательности α$^{-1}$ ∈ $F_{A,n-1}$, последовательность βα$^{-1}$β эквивалентна некоторой последовательности из $F_{B,n-1}$, то < B, [ ] > – n-арная подгруппа n-арной группы < A, [ ] >.

*Доказательство.* Так как < B, [ ] > – подалгебра, α и β составлены из элементов множества B, то αβ ~ γ ∈ $F_{B,n-1}$. Если γ$^{-1}$ – некоторая обратная к γ из $F_{A,n-1}$, то по условию имеем βγ$^{-1}$β ~ δ ∈ $F_{B,n-1}$. С другой стороны,

$$βγ^{-1}β ~ β(αβ)^{-1}β ~ ββ^{-1}α^{-1}β ~ α^{-1}β,$$

откуда, учитывая транзитивность отношения эквивалентности, получаем α$^{-1}$β ~ δ ∈ $F_{B,n-1}$. Применяя теперь теорему 5.1.10, заключаем, что < B, [ ] > – n-арная подгруппа n-арной группы < A, [ ] >. ∎

Приведем еще один, принадлежащий Посту, критерий существования n-арной подгруппы в n-арной группе.

**5.1.12. Теорема** [3]. Пусть < A, [ ] > n-арная группа, H – подгруппа соответствующей группы $A_o$, и существует такой x ∈ A, что:
1) $θ_A^{n-1}(x) ∈ H$;
2) $θ_A(x)H = Hθ_A(x)$.

Тогда < B, [ ] > является n-арной подгруппой в < A, [ ] >, где

$$B = [xF] = \{[xα] \mid α ∈ F\}, F = \{α = a_1 … a_{n-1} \mid θ_A(α) ∈ H\},$$

причем $B_o(A) = H$.

*Доказательство.* Пусть

$$b_1 = [x\alpha_1], \ldots, b_n = [x\alpha_n]$$

– произвольные элементы из В. Так как $\theta_A(\alpha_i) \in H$, то, используя 2), а также замкнутость H относительно бинарной операции в группе $A_o$, получим

$$\theta_A([b_1 \ldots b_n]) = \theta_A([x\alpha_1] \ldots [x\alpha_n]) =$$

$$\theta_A(x)\theta_A(\alpha_1)\theta_A(x)\theta_A(\alpha_2) \ldots \theta_A(x)\theta_A(\alpha_n) = \theta_A^n(x)\theta_A(\alpha),$$

то есть

$$\theta_A([b_1 \ldots b_n]) = \theta_A^n(x)\theta_A(\alpha)$$

для некоторого $\theta_A(\alpha) \in H$. Из последнего равенства, учитывая 1), получаем

$$\theta_A([b_1 \ldots b_n]) = \theta_A(x)\theta_A^{n-1}(x)\theta_A(\alpha) = \theta_A(x)\theta_A(\alpha'),$$

где $\theta_A(\alpha') \in H$, то есть

$$\theta_A([b_1 \ldots b_n]) = \theta_A([x\alpha']).$$

Следовательно,

$$[b_1 \ldots b_n] = [x\alpha'] \in B.$$

Пусть теперь $[x\alpha]$ – произвольный элемент из В. Так как

$$[\overline{[x\alpha]}\underbrace{[x\alpha] \ldots [x\alpha]}_{n-1}] = [x\alpha],$$

то

$$\theta_A(\overline{[x\alpha]}\underbrace{[x\alpha] \ldots [x\alpha]}_{n-1}) = \theta_A([x\alpha]),$$

откуда

$$\theta_A(\overline{[x\alpha]})\underbrace{\theta_A(x)\theta_A(\alpha)\ldots\theta_A(x)\theta_A(\alpha)}_{n-1} = \theta_A(x)\theta_A(\alpha).$$

Из последнего равенства, учитывая 1), 2), а также то, что H – подгруппа в $A_o$, получим

$$\theta_A(\overline{[x\alpha]})\theta_A^{n-1}(x)\theta_A(\beta) = \theta_A(x)\theta_A(\alpha), \quad \theta_A(\beta) \in H,$$

$$\theta_A(\overline{[x\alpha]})\theta_A(\gamma) = \theta_A(x)\theta_A(\alpha), \quad \theta_A(\gamma) \in H,$$

$$\theta_A(\overline{[x\alpha]}) = \theta_A(x)\theta_A(\alpha)\theta_A^{-1}(\gamma),$$

$$\theta_A(\overline{[x\alpha]}) = \theta_A(x)\theta_A(\delta), \quad \theta_A(\delta) \in H,$$

$$\theta_A(\overline{[x\alpha]}) = \theta_A([x\delta]).$$

Следовательно,

$$\overline{[x\alpha]} = [x\delta] \in B.$$

Согласно критерию Дернте, $< B, [\ ] >$ – n-арная подгруппа в $< A, [\ ] >$.

Если снова

$$b_1 = [x\alpha_1], \ldots, b_{n-1} = [x\alpha_{n-1}] \in B,$$

то снова, используя 1), 2), а также то, что H – подгруппа в $A_o$, получим

$$\theta_A(b_1 \ldots b_{n-1}) = \theta_A(x\alpha_1) \ldots \theta_A(x\alpha_{n-1}) =$$

$$= \theta_A(x)\theta_A(\alpha_1) \ldots \theta_A(x)\theta_A(\alpha_{n-1}) = \theta_A^{n-1}(x)\theta_A(\alpha) = \theta_A(\beta),$$

где $\theta_A(\alpha), \theta_A(\beta) \in H$, т.е.

$$B_o(A) \subseteq H. \qquad (1)$$

Любой элемент $\theta_A(\beta) \in H$ можно представить в виде

$$\theta_A(\beta) = \theta_A^{n-1}(x)\theta_A(\alpha) = \theta_A^{n-1}(x)\theta_A(\alpha_1') \ldots \theta_A(\alpha_{n-1}'), \qquad (2)$$

где

$$\theta_A(\alpha), \theta_A(\alpha'_1), \ldots, \theta_A(\alpha'_{n-1}) \in H.$$

Например,

$$\alpha'_1 = \alpha,\ \alpha'_2 = \ldots = \alpha'_{n-1} = \varepsilon,$$

где $\varepsilon$ – нейтральная последовательность.

Из (2), ввиду 2), следует

$$\theta_A(\beta) = \theta_A(x)\theta_A(\beta_1) \ldots \theta_A(x)\theta_A(\beta_{n-1}),$$

где $\theta_A(\beta_1), \ldots, \theta_A(\beta_{n-1}) \in H$. Тогда

$$\theta_A(\beta) = \theta_A([x\beta_1] \ldots [x\beta_{n-1}]) = \theta_A(b_1 \ldots b_{n-1}),$$

где $b_i = [x\beta_i] \in B$, $i = 1, \ldots, n-1$. Следовательно,

$$H \subseteq B_o(A). \tag{3}$$

Из (1) и (3) следует $B_o(A) = H$. ∎

## §5.2. n-АРНЫЕ АНАЛОГИ НОРМАЛЬНЫХ ПОДГРУПП

Согласно Дёрнте, n-арная подгруппа $<B, [\ ]>$ n-арной группы $<A, [\ ]>$ называется инвариантной в ней, если

$$[x\underbrace{B\ldots B}_{n-1}] = [\underbrace{B\ldots B}_{i-1}\, x\, \underbrace{B\ldots B}_{n-i}]$$

для любого $x \in A$ и всех $i = 2, 3, \ldots, n$ (определение 2.3.1).

n-Арную подгруппу $<B, [\ ]>$ n-арной группы $<A, [\ ]>$, удовлетворяющую условию

$$[x\underbrace{B\ldots B}_{n-1}] = [\underbrace{B\ldots B}_{n-1}x]$$

для любого х ∈ A, Дёрнте назвал полуинвариантной в ней (определение 2.3.1).

В данном параграфе определяются и изучаются некоторые новые n-арные аналоги нормальных подгрупп группы, отличные от указанных выше и подробно изученных в главе 2 (§2.3).

**5.2.1. Определение** [106]. n-Арная подгруппа $<B, [\ ]>$ n-арной группы $<A, [\ ]>$ называется *нормальной* в ней, если

$$[xx_1\ldots x_{n-2}B] = [Bx_1\ldots x_{n-2}x]$$

для любых $x_1, \ldots, x_{n-2}, x \in A$.

При n = 2 все три понятия: инвариантности, полуинвариантности и нормальности совпадают с понятием нормальности для подгрупп.

Ясно, что каждая n-арная группа инвариантна, нормальна и полуинвариантна в самой себе.

**5.2.2. Предложение.** Если $<B, [\ ]>$ – нормальная n-арная подгруппа n-арной группы $<A, [\ ]>$, то она и полуинвариантна в $<A, [\ ]>$.

*Доказательство.* Из определения 5.1.1 следует справедливость следующего равенства

$$[x\underbrace{b\ldots b}_{n-2}B] = [B\underbrace{b\ldots b}_{n-2}x]$$

для любого $b \in B$, откуда

$$[x\underbrace{B\ldots B}_{n-1}] = [\underbrace{B\ldots B}_{n-1}x]. \qquad \blacksquare$$

Пусть $<A, [\ ]>$ полуабелева n-арная группа, то есть n-арная группа, удовлетворяющая тождеству

$$[xx_1 \ldots x_{n-2}y] = [yx_1 \ldots x_{n-2}x].$$

Если < B, [ ] > – n-арная подгруппа в < A, [ ] >, то

$$[xx_1 \ldots x_{n-2}B] = \{[xx_1 \ldots x_{n-2}b] \mid b \in B\} =$$

$$= \{[bx_1 \ldots x_{n-2}x] \mid b \in B\} = [Bx_1 \ldots x_{n-2}x],$$

то есть < B, [ ] > нормальна в < A, [ ] >. Таким образом, имеет место

**5.2.3. Предложение.** В полуабелевой n-арной группе все её n-арные подгруппы являются нормальными.

Из предложений 5.2.2 и 5.2.3 следует известное утверждение о полуинвариантности всех n-арных подгрупп в полуабелевой n-арной группе.

В качестве n-арной группы из предложения 5.2.3 можно взять, например, полуабелеву n-арную группу < $B_n$, [ ] > всех отражений правильного n-угольника [26].

Так как в < $B_3$, [ ] > все три тернарные подгруппы первого порядка не являются инвариантными, то из сказанного выше следует, что при n ≥ 3 существуют n-арные группы, обладающие нормальными n-арными подгруппами, не являющимися инвариантными.

А существуют ли при n ≥ 3 n-арные группы, обладающие инвариантными n-арными подгруппами, которые не являются нормальными? Положительный ответ на этот вопрос даёт следующий

**5.2.4. Пример.** Рассмотрим тернарную группу < $S_4$, [ ] >, производную от симметрической группы $S_4$ на четырёх символах. Так как четверная подгруппа Клейна

$$V_4 = \{e, (12)(34), (13)(24), (14)(23)\}$$

нормальна в группе $S_4$, то тернарная группа < $V_4$, [ ] > инвариантна в < $S_4$, [ ] > (см. пример 2.3.3). Так как

$$[V_4(12)(13)] = \left\{ e\begin{pmatrix}1234\\2314\end{pmatrix}, \begin{pmatrix}1234\\2143\end{pmatrix}\begin{pmatrix}1234\\2314\end{pmatrix}, \begin{pmatrix}1234\\3412\end{pmatrix}\begin{pmatrix}1234\\2314\end{pmatrix}, \begin{pmatrix}1234\\4321\end{pmatrix}\begin{pmatrix}1234\\2314\end{pmatrix} \right\} =$$

$$= \left\{ \begin{pmatrix}1234\\2314\end{pmatrix}, \begin{pmatrix}1234\\3241\end{pmatrix}, \begin{pmatrix}1234\\1423\end{pmatrix}, \begin{pmatrix}1234\\4132\end{pmatrix} \right\} = \{(123),(134),(243),(142)\},$$

то

$$(13)(12) = (132) \notin [V_4(12)(13)].$$

А так как

$$(13)(12) \in [(13)(12)V_4],$$

то

$$[(13)(12)V_4] \neq [V_4(12)(13)].$$

Следовательно, инвариантная в $<S_4, [\ ]>$ тернарная подгруппа $<V_4, [\ ]>$ не является нормальной в $<S_4, [\ ]>$.

Отметим, что даже одноэлементная n-арная подгруппа, инвариантная в содержащей её n-арной группе может не быть нормальной в этой же n-арной группе. Об этом свидетельствует следующий

**5.2.5. Пример.** Пусть $<A, [\ ]>$ – n-арная группа, производная от группы A с центром $Z(A) \neq A$. Если e – единица группы A, то $<E = \{e\}, [\ ]>$ – инвариантная n-арная подгруппа в $<A, [\ ]>$. Если $x \notin Z(A)$, то существует $y \in A$ такой, что $xy \neq yx$. Тогда

$$[xyE] \neq [Eyx].$$

Следовательно, инвариантная в $<A, [\ ]>$ – n-арная подгруппа $<E, [\ ]>$ не является нормальной в $<A, [\ ]>$.

В качестве группы A можно взять, например, симметрическую группу $S_m(m \geq 3)$, знакопеременную группу $A_m(m \geq 4)$ или диэдральную группу $D_m(m \geq 2)$.

Покажем, что в n-арной группе n-арная подгруппа может быть одновременно и нормальной и инвариантной.

**5.2.6. Пример.** Пусть $< S_n, [\ ] >$ – тернарная группа, производная от симметрической группы $S_n$. Так как $A_n$ нормальна в $S_n$, то $< A_n, [\ ] >$ инвариантна в $< S_n, [\ ] >$. Так как

$$[xB_nB_n] = [B_nxB_n] = [B_nB_nx] = A_n, x \in A_n,$$

$$[xB_nB_n] = [B_nxB_n] = [B_nB_nx] = B_n, x \in B_n,$$

то тернарная группа $< B_n, [\ ] >$ всех нечётных подстановок также инвариантна в $< S_n, [\ ] >$.

Если x и y – произвольные подстановки из $S_n$ разной чётности, то

$$[xyA_n] = [A_nyx] = B_n, [xyB_n] = [B_nyx] = A_n.$$

Если x и y имеют одинаковую чётность, то

$$[xyA_n] = [Ayx] = A_n, [xyB_n] = [B_nyx] = B_n.$$

Таким образом,

$$[xyA_n] = [A_nyx], [xyB_n] = [B_nyx]$$

для любых $x, y \in S_n$. Следовательно, $< A_n, [\ ] >$ и $< B_n, [\ ] >$ нормальные в $< S_n, [\ ] >$ тернарные подгруппы.

Так как инвариантные и нормальные n-арные подгруппы n-арной группы являются и полуинвариантными в ней, то возникает естественный вопрос: исчерпываются ли полуинвариантные n-арные подгруппы любой n-арной группы нормальными n-арными подгруппами и инвариантными n-арными подгруппами? Иначе говоря, может ли в n-арной группе быть полуинвариантная n-арная подгруппа, которая не является нормальной и не является инвариантной?

Положительный ответ на последний вопрос следует из следующего примера.

**5.2.7. Пример.** Пусть $< D_6, [\ ] >$ – тернарная группа, производная от диэдральной группы $D_6 = C_6 \cup B_6$, где

$$C_6 = \{e, c, c^2, c^3, c^4, c^5\}$$

– циклическая подгруппа порядка 6, порождённая элементом c,

$$B_6 = \{b, bc, bc^2, bc^3, bc^4, bc^5\}$$

множество всех отражений правильного шестиугольника.

Если $B = \{b, bc^3\}$, то по предложению 2.1.11 (см. также [26]) $<B, [\ ]>$ – тернарная подгруппа в $<D_6, [\ ]>$.

Так как

$$[ccB] = c^2B = \{c^2b, c^2bc^3\} = \{bc^4, bc^4c^3\} = \{bc^4, bc\},$$

$$[Bcc] = Bc^2 = \{bc^2, bc^3c^2\} = \{bc^2, bc^5\},$$

то

$$[ccB] \neq [Bcc].$$

Следовательно, $<B, [\ ]>$ не является нормальной в $<D_6, [\ ]>$.

Так как

$$[BBc] = BBc = \{b^2c, bc^3bc, bbc^3c, bc^3bc^3c\} =$$
$$= \{c, bbc^3c, c^4, bbc^3c^3c\} = \{c, c^4, c^4, c\} = \{c, c^4\},$$

$$[BcB] = BcB = \{bcb, bc^3cb, bcbc^3, bc^3cbc^3\} =$$
$$= \{bbc^5, bbc^2, bbc^5c^3, bbc^2c^3\} = \{c^5, c^2, c^2, c^5\} = \{c^2, c^5\},$$

то

$$[BBc] \neq [BcB].$$

Следовательно, $<B, [\ ]>$ не является инвариантной в $<D_6, [\ ]>$.

Так как

$$[cBB] = cBB = \{cb^2, cbc^3b, cbbc^3, cbc^3bc^3\} =$$
$$= \{c, cbbc^3, c^4, cbbc^3c^3\} = \{c, c^4, c^4, c^7\} = \{c, c^4\},$$

то

$$[BBc] = [cBB] = \{c, c^4\},$$

откуда следует

$$[BBx] = [xBB]$$

для любого $x \in C_6$.

Так как $<B_6, [\ ]>$ полуабелева [26], то в $<B_6, [\ ]>$ все тернарные подгруппы полуинвариантны. Поэтому

$$[BBx] = [xBB]$$

для любого $x \in B_6$.

Таким образом, $<B, [\ ]>$ – полуинвариантная в $<D_6, [\ ]>$ тернарная подгруппа, которая не является нормальной в $<D_6, [\ ]>$ и не является инвариантной в $<D_6, [\ ]>$.

Если при $n \geq 3$ зафиксировать элементы $a_1, \ldots, a_{n-3} \in A$, то для любых $x_1, \ldots, x_{n-2} \in A$ существуют $y, z \in A$, такие, что все три последовательности

$$x_1 \ldots x_{n-2}, \quad ya_1 \ldots a_{n-3}, \quad a_1 \ldots a_{n-3}z$$

эквивалентны. Поэтому справедливо

**5.2.8. Предложение** [106]. Если $<B, [\ ]>$ – n-арная подгруппа n-арной группы $<A, [\ ]>$, $n \geq 3$, то следующие утверждения эквивалентны:

1) $<B, [\ ]>$ – нормальна в $<A, [\ ]>$;

2) $[xya_1 \ldots a_{n-3}B] = [Bya_1 \ldots a_{n-3}x]$ для любых $x, y \in A$ и некоторых $a_1, \ldots, a_{n-3} \in A$;

3) $[xa_1 \ldots a_{n-3}zB] = [Ba_1 \ldots a_{n-3}zx]$ для любых $x, z \in A$ и некоторых $a_1, \ldots, a_{n-3} \in A$.

Полагая в предложении 5.2.8

$$a_1 = \ldots = a_{n-3} = b \in B,$$

получим

**5.2.9. Предложение** [106]. Если $<B, [\ ]>$ – n-арная подгруппа n-арной группы $<A, [\ ]>$, $n \geq 3$, то следующие утверждения эквивалентны:

1) $<B, [\ ]>$ – нормальна в $<A, [\ ]>$;

2) $[xy\underbrace{b \ldots b}_{n-3}B] = [By\underbrace{b \ldots b}_{n-3}x]$ для любых $x, y \in A$ и некоторого $b \in B$;

3) $[x\underbrace{b \ldots b}_{n-3}zB] = [B\underbrace{b \ldots b}_{n-3}zx]$ для любых $x, z \in A$ и некоторого $b \in B$.

Так как

$$[xy\underbrace{b\ldots b}_{n-3}B] = [xy\underbrace{B\ldots B}_{n-2}], \quad [B\underbrace{b\ldots b}_{n-3}zx] = [\underbrace{B\ldots B}_{n-2}zx]$$

для любого $b \in B$, то из предложения 5.2.9 вытекает

**5.2.10. Следствие** [106]. Для n-арной подгруппы $< B, [\ ] >$ n-арной группы $< A, [\ ] >$, где $n \geq 3$, следующие утверждения эквивалентны:

1) $< B, [\ ] >$ — нормальна в $< A, [\ ] >$;

2) $[xy\underbrace{B\ldots B}_{n-2}] = [By\underbrace{b\ldots b}_{n-3}x]$ для любых $x, y \in A$ и некоторого $b \in B$;

3) $[x\underbrace{b\ldots b}_{n-3}zB] = [\underbrace{B\ldots B}_{n-2}zx]$ для любых $x, z \in A$ и некоторого $b \in B$.

Используя определение 5.2.1, предложение 5.2.9 и следствие 5.2.10, получим

**5.2.11. Предложение** [106]. Для n-арной подгруппы $< B, [\ ] >$ n-арной группы $< A, [\ ] >$, где $n \geq 3$, следующие утверждения эквивалентны:

1) $< B, [\ ] >$ — нормальна в $< A, [\ ] >$;

2) $[xx_1 \ldots x_{n-2}By_1 \ldots y_{n-1}] = B$ для любых $x, x_1, \ldots, x_{n-2} \in A$ и любой последовательности $y_1\ldots y_{n-1}$, обратной для $x_1\ldots x_{n-2}x$;

3) $[z_1 \ldots z_{n-1}Bx_1 \ldots x_{n-2}x] = B$ для любых $x_1, \ldots, x_{n-2}, x \in A$ и любой последовательности $z_1\ldots z_{n-1}$ обратной для $xx_1\ldots x_{n-2}$;

4) $[xy\underbrace{b\ldots b}_{n-3}B\overline{x}\underbrace{x\ldots x}_{n-3}\overline{b}b\overline{y}\underbrace{y\ldots y}_{n-3}] = B$ для любых $x, y \in A$ и некоторого $b \in B$;

5) $[\overline{z}\underbrace{z\ldots z}_{n-3}\overline{b}b\overline{x}\underbrace{x\ldots x}_{n-3}B\underbrace{b\ldots b}_{n-3}zx] = B$ для любых $x, z \in A$ и некоторого $b \in B$;

6) $[xy\underbrace{B\ldots B}_{n-2}\overline{x}\underbrace{x\ldots x}_{n-3}\overline{b}b\overline{y}\underbrace{y\ldots y}_{n-3}] = B$ для любых $x, y \in A$ и некоторого $b \in B$;

7) $[\bar{z}\underbrace{z\ldots z}_{n-3}\bar{b}bb\bar{x}\underbrace{x\ldots x}_{n-3}\underbrace{B\ldots B}_{n-2}zx] = B$ для любых $x, z \in A$ и некоторого $b \in B$.

Полагая в предложении 5.2.11 n = 3, получим

**5.2.12. Предложение** [106]. Для тернарной подгруппы $< B, [\,]>$ тернарной группы $< A, [\,]>$ следующие утверждения эквивалентны:
1) $< B, [\,]>$ – нормальна в $< A, [\,]>$;
2) $[xyB\bar{x}\bar{y}] = B$ для любых $x, y \in A$;
3) $[\bar{z}\bar{x}Bzx] = B$ для любых $x, z \in A$.

**5.2.13. Лемма** [106]. Если

$$< B, [\,]> = <\cap B_i, [\,]>, \cap B_i \neq \varnothing$$

– непустое пересечение семейства $\{< B_i, [\,]> | i \in I\}$ нормальных n-арных подгрупп n-арной группы $< A, [\,]>$, то n-арная подгруппа $< B, [\,]>$ нормальна в $< A, [\,]>$.

*Доказательство.* Так как все $< B_i, [\,]>$ нормальны в $< A, [\,]>$, то согласно 2) и 3) предложения 5.2.11

$$[\alpha B_i \beta^{-1}] = B_i, [\alpha^{-1} B_i \beta] = B_i,$$

где

$$\alpha = xx_1\ldots x_{n-2}, \ \beta = x_1\ldots x_{n-2}x,$$

$\alpha^{-1}$ и $\beta^{-1}$ – обратные последовательности для $\alpha$ и $\beta$ соответственно. Тогда

$$[\alpha B \beta^{-1}] = [\alpha(\cap B_i)\beta^{-1}] \subseteq [\alpha B_i \beta^{-1}] = B_i, i \in I,$$

$$[\alpha^{-1} B \beta] = [\alpha^{-1}(\cap B_i)\beta] \subseteq [\alpha^{-1} B_i \beta] = B_i, i \in I,$$

откуда

$$[\alpha B \beta^{-1}] \subseteq \cap B_i = B, [\alpha^{-1} B \beta] \subseteq \cap B_i = B,$$

то есть

$$[\alpha B \beta^{-1}] \subseteq B, \tag{1}$$

$$[\alpha^{-1} B \beta] \subseteq B. \tag{2}$$

Из (2) следует

$$[\alpha[\alpha^{-1}B\beta]\beta^{-1}] \subseteq [\alpha B \beta^{-1}],$$

$$B \subseteq [\alpha B \beta^{-1}]. \tag{3}$$

Из (1) и (3) следует

$$[\alpha B \beta^{-1}] = B.$$

Тогда согласно 2) предложения 5.2.11, $< B, [\ ] >$ нормальна в $< A, [\ ] >$. ∎

**5.2.14. Лемма** [106]. Если $< B, [\ ] >$ – полуинвариантная, $< C, [\ ] >$ – нормальная n-арные подгруппы n-арной группы $< A, [\ ] >$, то

$$< D, [\ ] > = < [\underbrace{B \ldots B}_{n-1} C], [\ ] >$$

– нормальная n-арная подгруппа в $< A, [\ ] >$.

*Доказательство.* Ввиду предложения 5.2.2, $< C, [\ ] >$ – полуинвариантна в $< A, [\ ] >$. Поэтому по следствию 2.3.18 $< D, [\ ] >$ – полуинвариантная n-арная подгруппа в $< A, [\ ] >$.

Используя полуинвариантность $< B, [\ ] >$, а также нормальность $< C, [\ ] >$, получаем

$$[xx_1 \ldots x_{n-2} D] = [xx_1 \ldots x_{n-2}[\underbrace{B \ldots B}_{n-1} C]] =$$

$$= [\underbrace{B \ldots B}_{n-1}[xx_1 \ldots x_{n-2} C]] = [\underbrace{B \ldots B}_{n-1}[Cx_1 \ldots x_{n-2} x]] =$$

$$= [[\underbrace{B \ldots B}_{n-1} C] x_1 \ldots x_{n-2} x] = [D x_1 \ldots x_{n-2} x],$$

то есть

$$[xx_1\ldots x_{n-2}D] = [D\ x_1\ldots x_{n-2}x].$$

Следовательно, $<D, [\ ]>$ – нормальна в $<A, [\ ]>$. ∎

**5.2.15. Лемма** [106]. Пусть $<B, [\ ]>$ и $<C, [\ ]>$ – n-арные подгруппы n-арной группы $<A, [\ ]>$, причем $<B, [\ ]>$ – полуинвариантна в $<A, [\ ]>$ и $B\cap C \neq \varnothing$. Тогда $<[\underbrace{B\ldots B}_{n-1}C], [\ ]>$ – n-арная подгруппа, порождённая множеством $B\cup C$.

*Доказательство.* Пусть

$$<D, [\ ]> = <B, [\ ]> \vee <C, [\ ]>$$

– n-арная подгруппа, порождённая множеством $B\cup C$. Так как в $<B, [\ ]>$ и в $<C, [\ ]>$ имеются нейтральные последовательности, то

$$C \subseteq [\underbrace{B\ldots B}_{n-1}C],\ B \subseteq [B\underbrace{C\ldots C}_{n-1}],$$

а так как, кроме того, $B\cap C \neq \varnothing$, то из второго включения, учитывая лемму 2.3.20, получаем

$$B \subseteq [\underbrace{B\ldots B}_{n-1}C],$$

откуда

$$B\cup C \subseteq [\underbrace{B\ldots B}_{n-1}C]. \qquad (1)$$

Так как $<B, [\ ]>$ – полуинвариантна в $<A, [\ ]>$, то

$$[\underbrace{B\ldots B}_{n-1}C] = [C\underbrace{B\ldots B}_{n-1}],$$

откуда и из леммы 2.3.17 следует, что $<[\underbrace{B\ldots B}_{n-1}C], [\ ]>$ – n-арная подгруппа в $<A, [\ ]>$. Тогда из (1) вытекает

$$D \subseteq [\underbrace{B\ldots B}_{n-1}C].$$

Включение

$$[\underbrace{B\ldots B}_{n-1}C] \subseteq D$$

является следствием теоремы 2.1.14. Таким образом,

$$D = [\underbrace{B\ldots B}_{n-1}C]. \qquad \blacksquare$$

**5.2.16. Теорема** [106]. Множество всех нормальных n-арных подгрупп n-арной группы $<A, [\,]>$, содержащих фиксированный элемент, образуют подрешетку решетки $L(A, [\,])$ всех n-арных подгрупп n-арной группы $<A, [\,]>$.

*Доказательство.* Если $<B, [\,]>$ и $<C, [\,]>$ – нормальные n-арные подгруппы из $<A, [\,]>$, содержащие фиксированный элемент $a \in A$, то по лемме 5.2.13

$$<B, [\,]> \wedge <C, [\,]> = <B \cap C, [\,]>$$

является нормальной n-арной подгруппой в $<A, [\,]>$, причем $a \in B \cap C$.

Положим

$$<D, [\,]> = <B, [\,]> \vee <C, [\,]>.$$

Тогда ввиду предложения 5.2.2 и леммы 5.2.15, имеем

$$<D, [\,]> = <[\underbrace{B\ldots B}_{n-1}C], [\,]>,$$

а по лемме 5.2.14 $<D, [\,]>$ – нормальная n-арная подгруппа в $<A, [\,]>$. Ясно, что $a \in D$. $\blacksquare$

## §5.3 Σ-НОРМАЛЬНЫЕ n-АРНЫЕ ПОДГРУППЫ

Инвариантные n-арные подгруппы, а также нормальные n-арные подгруппы являются частными случаями более общего понятия – σ-нормальной n-арной подгруппы.

**5.3.1. Определение.** Пусть $\Sigma$ – подмножество множества $S_{n-1}$ всех подстановок на множестве $\{1, 2, \ldots, n-1\}$. n-Арная подгруппа $< B, [\,] >$ n-арной группы $< A, [\,] >$ называется *Σ-нормальной* в $< A, [\,] >$, если

$$[x_1 x_2 \ldots x_{n-1} B] = [B x_{\sigma(1)} x_{\sigma(2)} \ldots x_{\sigma(n-1)}] \qquad (*)$$

для всех $x, x_2, \ldots, x_{n-1} \in A$ и любой подстановки $\sigma$ из $\Sigma$.

Если $\Sigma = \{\sigma\}$, то Σ-нормальную n-арную подгруппу будем называть *σ-нормальной*.

**5.3.2. Предложение.** n-Арная подгруппа $< B, [\,] >$ n-арной группы $< A, [\,] >$ является инвариантной в $< A, [\,] >$ тогда и только тогда, когда она ε-нормальна в $< A, [\,] >$, где ε – тождественная подстановка из $S_{n-1}$, то есть тогда и только тогда, когда

$$[x_1 x_2 \ldots x_{n-1} B] = [B x_1 x_2 \ldots x_{n-1}] \qquad (1)$$

для всех $x, x_2, \ldots, x_{n-1} \in A$.

*Доказательство*. *Необходимость*. Так как $< B, [\,] >$ инвариантна в $< A, [\,] >$, то, согласно 3) теоремы 2.3.9, $[\alpha B \beta] = B$, где $\alpha = x_1 x_2 \ldots x_{n-1}$, $\beta$ – обратная последовательность для $\alpha$. Тогда, учитывая нейтральность последовательности $\beta\alpha$, получим

$$[[\alpha B \beta]\alpha] = [B\alpha], [\alpha B] = [B\alpha],$$

то есть верно (1).

*Достаточность*. Полагая в (1)

$$x_1 = x \in A, \ x_2 = \ldots = x_{n-1} = b \in B,$$

получаем

$$[x \underbrace{b \ldots b}_{n-2} B] = [Bx \underbrace{b \ldots b}_{n-2}],$$

откуда, учитывая нейтральность последовательностей

$$\overline{b} \underbrace{b \ldots b}_{n-2}, \ \underbrace{b \ldots b}_{n-2} \overline{b}, \ \underbrace{x \ldots x}_{n-2} \overline{x},$$

а также тот факт, что $\overline{b} \in B$, последовательно будем иметь

$$[\overline{b} \underbrace{b \ldots b}_{n-2} [x \underbrace{b \ldots b}_{n-2} B] \overline{b} \underbrace{x \ldots x}_{n-3} \overline{x}] = [\overline{b} \underbrace{b \ldots b}_{n-2} [Bx \underbrace{b \ldots b}_{n-2}] \overline{b} \underbrace{x \ldots x}_{n-3} \overline{x}],$$

$$[x[\underbrace{b \ldots b}_{n-2} B \overline{b}] \underbrace{x \ldots x}_{n-3} \overline{x}] = [\overline{b} \underbrace{b \ldots b}_{n-2} B], \ [xB \underbrace{x \ldots x}_{n-3} \overline{x}] = B$$

для любого $x \in A$. Тогда, согласно 4) теоремы 2.3.9, $< B, [\ ] >$ инвариантна в $< A, [\ ] >$. ∎

Если $\sigma = (1 \ 2 \ \ldots \ n-1)$ – циклическая подстановка из $S_{n-1}$, то

$$x_{\sigma(1)} x_{\sigma(2)} \ldots x_{\sigma(n-1)} = x_2 x_3 \ldots x_{n-1} x_1,$$

а равенство (*) примет вид

$$[x_1 x_2 \ldots x_{n-1} B] = [B x_2 x_3 \ldots x_{n-1} x_1].$$

Поэтому имеет место

**5.3.3. Предложение.** n-Арная подгруппа $< B, [\ ] >$ n-арной группы $< A, [\ ] >$ является нормальной в $< A, [\ ] >$ тогда и только тогда, когда она σ-нормальна в $< A, [\ ] >$, где $\sigma = (1 \ 2 \ \ldots \ n-1)$.

Таким образом, инвариантные n-арные подгруппы – это в точности ε-нормальные n-арные подгруппы для тождественной подстановки ε, а нормальные n-арные подгруппы – это в точности σ-нормальные n-арные подгруппы для подстановки σ = (1 2 … n – 1).

Понятно, что n-арная подгруппа < B, [ ] > n-арной группы < A, [ ] > является одновременно и инвариантной и нормальной в < A, [ ] > тогда и только тогда, когда она Σ-нормальна в < A, [ ] >, где Σ = {ε, σ = (1 2 … n – 1)}.

**5.3.4. Лемма.** Если < B, [ ] > – σ-нормальная n-арная подгруппа n-арной группы < A, [ ] > и σ(j) = i для некоторых i, j ∈ {1, …, n – 1), то

$$[\underbrace{B\ldots B}_{i-1} x \underbrace{B\ldots B}_{n-i}] = [\underbrace{B\ldots B}_{j} x \underbrace{B\ldots B}_{n-j-1}]$$

для любого x ∈ A.

*Доказательство*. Сразу же заметим, что если j = i – 1, то равенство из условия леммы принимает вид

$$[\underbrace{B\ldots B}_{j} x \underbrace{B\ldots B}_{n-j-1}] = [\underbrace{B\ldots B}_{j} x \underbrace{B\ldots B}_{n-j-1}],$$

не требующий доказательства. Поэтому случай j = i – 1 можно не рассматривать.

Положим $x_i = x \in A$. Так как σ(j) = i, то равенство (*) примет вид

$$[x_1 \ldots x_{i-1} x x_{i+1} \ldots x_{n-1} B] = [B x_{\sigma(1)} \ldots x_{\sigma(j-1)} x x_{\sigma(j+1)} \ldots x_{\sigma(n-1)}]. \quad (1)$$

Если $u = [b_1 \ldots b_{i-1} x b_{i+1} \ldots b_{n-1} b]$ – произвольный элемент из

$$[\underbrace{B\ldots B}_{i-1} x \underbrace{B\ldots B}_{n-i}],$$

то ввиду (1),

$$[b_1 \ldots b_{i-1} x b_{i+1} \ldots b_{n-1} b] = [c b_{\sigma(1)} \ldots b_{\sigma(j-1)} x b_{\sigma(j+1)} \ldots b_{\sigma(n-1)}]$$

для некоторого c ∈ B, где, как легко заметить, из условия σ(j) = i вытекает

$$b_{\sigma(1)}, \ldots, b_{\sigma(j-1)}, b_{\sigma(j+1)}, \ldots, b_{\sigma(n-1)} \in B.$$

Следовательно,

$$u = [b_1 \ldots b_{i-1} x b_{i+1} \ldots b_{n-1} b] \in [\underbrace{B \ldots B}_{j} x \underbrace{B \ldots B}_{n-j-1}],$$

откуда, ввиду произвольного выбора u, следует включение

$$[\underbrace{B \ldots B}_{i-1} x \underbrace{B \ldots B}_{n-i}] \subseteq [\underbrace{B \ldots B}_{j} x \underbrace{B \ldots B}_{n-j-1}]. \qquad (2)$$

Если $v = [bb_1 \ldots b_{j-1} x b_{j+1} \ldots b_{n-1}]$ – произвольный элемент из

$$[\underbrace{B \ldots B}_{j} x \underbrace{B \ldots B}_{n-j-1}],$$

то, положив

$$c_{\sigma(1)} = b_1, \ldots, c_{\sigma(j-1)} = b_{j-1}, c_{\sigma(j+1)} = b_{j+1}, \ldots, c_{\sigma(n-1)} = b_{n-1},$$

и, используя (1), получим

$$[c_1 \ldots c_{i-1} x c_{i+1} \ldots c_{n-1} c] = [b c_{\sigma(1)} \ldots c_{\sigma(j-1)} x c_{\sigma(j+1)} \ldots c_{\sigma(n-1)}]$$

для некоторого c ∈ B, где, как легко заметить, из условия σ(j) = i вытекает

$$c_1, \ldots, c_{i-1}, c_{i+1}, \ldots, c_{n-1} \in B.$$

Следовательно,

$$v = [bb_1 \ldots b_{j-1} x b_{j+1} \ldots b_{n-1}] \in [\underbrace{B \ldots B}_{i-1} x \underbrace{B \ldots B}_{n-i}],$$

откуда, ввиду произвольного выбора v, следует включение

$$[\underbrace{B\ldots B}_{j}x\underbrace{B\ldots B}_{n-j-1}] \subseteq [\underbrace{B\ldots B}_{i-1}x\underbrace{B\ldots B}_{n-i}]. \qquad (3)$$

Из (2) и (3) следует требуемое равенство. ∎

**5.3.5. Лемма.** n-Арная подгруппа $<B, [\ ]>$ n-арной группы $<A, [\ ]>$ является инвариантной в $<A, [\ ]>$, если для любого $x \in A$ и некоторого $i \in \{1, \ldots, n-1\}$ выполняется по крайней мере одно из следующих условий:

1) $[x\underbrace{B\ldots B}_{n-1}] = [Bx\underbrace{B\ldots B}_{n-2}]$, $[\underbrace{B\ldots B}_{i-1}x\underbrace{B\ldots B}_{n-i}] = [\underbrace{B\ldots B}_{n-1}x]$;

2) $[x\underbrace{B\ldots B}_{n-1}] = [\underbrace{B\ldots B}_{n-2}xB]$, $[\underbrace{B\ldots B}_{i-1}x\underbrace{B\ldots B}_{n-i}] = [\underbrace{B\ldots B}_{n-1}x]$;

3) $[\underbrace{B\ldots B}_{n-1}x] = [\underbrace{B\ldots B}_{n-2}xB]$, $[\underbrace{B\ldots B}_{n-i}x\underbrace{B\ldots B}_{i-1}] = [x\underbrace{B\ldots B}_{n-1}]$;

4) $[\underbrace{B\ldots B}_{n-1}x] = [Bx\underbrace{B\ldots B}_{n-2}]$, $[\underbrace{B\ldots B}_{n-i}x\underbrace{B\ldots B}_{i-1}] = [x\underbrace{B\ldots B}_{n-1}]$.

*Доказательство*. 1) Используя первое равенство, получим

$$[x\underbrace{B\ldots B}_{n-1}] = [Bx\underbrace{B\ldots B}_{n-2}] = [B[x\underbrace{B\ldots B}_{n-1}]\underbrace{B\ldots B}_{n-2}] =$$

$$= [B[Bx\underbrace{B\ldots B}_{n-2}]\underbrace{B\ldots B}_{n-2}] = [BBx\underbrace{B\ldots B}_{n-3}] = \ldots = [\underbrace{B\ldots B}_{n-2}xB],$$

то есть

$$[x\underbrace{B\ldots B}_{n-1}] = [Bx\underbrace{B\ldots B}_{n-2}] = [BBx\underbrace{B\ldots B}_{n-3}] = \ldots = [\underbrace{B\ldots B}_{n-2}xB], \quad (1)$$

откуда и из второго равенства условия 1) следует инвариантность $<B, [\ ]>$ в $<A, [\ ]>$.

2) Используя первое равенство, получим

$$[x\underbrace{B\ldots B}_{n-1}] = [\underbrace{B\ldots B}_{n-2}xB] = [\underbrace{B\ldots B}_{n-2}[x\underbrace{B\ldots B}_{n-1}]B] =$$

$$= [\underbrace{B\ldots B}_{n-2}[\underbrace{B\ldots B}_{n-2}xB]B] = [\underbrace{B\ldots B}_{n-3}xBB] = \ldots = [Bx\underbrace{B\ldots B}_{n-2}],$$

то есть верно (1), откуда и из второго равенства условия 2) следует инвариантность $< B, [\ ] >$ в $< A, [\ ] >$.

3) Используя первое равенство, получим

$$[\underbrace{B\ldots B}_{n-1}x] = [\underbrace{B\ldots B}_{n-2}xB] = [\underbrace{B\ldots B}_{n-2}[\underbrace{B\ldots B}_{n-1}x]B] =$$

$$= [\underbrace{B\ldots B}_{n-2}[\underbrace{B\ldots B}_{n-2}xB]B] = [\underbrace{B\ldots B}_{n-3}xBB] = \ldots = [Bx\underbrace{B\ldots B}_{n-2}],$$

то есть

$$[\underbrace{B\ldots B}_{n-1}x] = [\underbrace{B\ldots B}_{n-2}xB] = [\underbrace{B\ldots B}_{n-3}xBB] = \ldots = [Bx\underbrace{B\ldots B}_{n-2}], \quad (2)$$

откуда и из второго равенства условия 3) следует инвариантность $< B, [\ ] >$ в $< A, [\ ] >$.

4) Используя первое равенство, получим

$$[\underbrace{B\ldots B}_{n-1}x] = [Bx\underbrace{B\ldots B}_{n-2}] = [B[\underbrace{B\ldots B}_{n-1}x]\underbrace{B\ldots B}_{n-2}] =$$

$$= [B[Bx\underbrace{B\ldots B}_{n-2}]\underbrace{B\ldots B}_{n-2}] = [BBx\underbrace{B\ldots B}_{n-3}] = \ldots = [\underbrace{B\ldots B}_{n-2}xB],$$

то есть верно (2), откуда и из второго равенства условия 4) следует инвариантность $< B, [\ ] >$ в $< A, [\ ] >$. ∎

**5.3.6. Теорема.** σ-Нормальная n-арная подгруппа $< B, [\ ] >$ n-арной группы $< A, [\ ] >$ является инвариантной в $< A, [\ ] >$, если выполняется по крайней мере одно из следующих условий:

1) $\sigma(i) = i$ для некоторого $i \in \{1, \ldots, n-1\}$;
2) $\sigma(j) = j + 2$ для некоторого $j \in \{1, \ldots, n-3\}$;

3) $\sigma(n-2) = 1$;
4) $\sigma(n-1) = 2$.

*Доказательство*. 1) Полагая в (*)

$$x_1 = \ldots = x_{i-1} = x_{i+1} = \ldots = x_{n-1} = b \in B, \, x_i = x \in A,$$

получаем

$$[\underbrace{b \ldots b}_{i-1} x \underbrace{b \ldots b}_{n-i-1} B] = [B \underbrace{b \ldots b}_{i-1} x \underbrace{b \ldots b}_{n-i-1}],$$

откуда, учитывая нейтральность последовательностей

$$\overline{b} \underbrace{b \ldots b}_{n-2}, \; \underbrace{b \ldots b}_{n-2} \overline{b}, \; \underbrace{x \ldots x}_{n-2} \overline{x},$$

а также тот факт, что $\overline{b} \in B$, последовательно будем иметь

$$[\overline{b} \underbrace{b \ldots b}_{n-i-1} [\underbrace{b \ldots b}_{i-1} x \underbrace{b \ldots b}_{n-i-1} B] \underbrace{b \ldots b}_{i-1} \overline{b} \underbrace{x \ldots x}_{n-3} \overline{x}] =$$

$$= [\overline{b} \underbrace{b \ldots b}_{n-i-1} [B \underbrace{b \ldots b}_{i-1} x \underbrace{b \ldots b}_{n-i-1}] \underbrace{b \ldots b}_{i-1} \overline{b} \underbrace{x \ldots x}_{n-3} \overline{x}],$$

$$[\overline{b} \underbrace{b \ldots b}_{n-2} x [\underbrace{b \ldots b}_{n-i-1} B \underbrace{b \ldots b}_{i-1} \overline{b}] \underbrace{x \ldots x}_{n-3} \overline{x}] =$$

$$= [[\overline{b} \underbrace{b \ldots b}_{n-i-1} B \underbrace{b \ldots b}_{i-1}] x \underbrace{b \ldots b}_{n-2} \overline{b} \underbrace{x \ldots x}_{n-3} \overline{x}],$$

$$[xB \underbrace{x \ldots x}_{n-3} \overline{x}] = B$$

для любого $x \in A$. Тогда, согласно 4) теоремы 2.3.9, $< B, [\,] >$ инвариантна в $< A, [\,] >$.

2) Полагая в (*)

$$x_1 = \ldots = x_{j+1} = x_{j+3} = \ldots = x_{n-1} = b \in B, \, x_{j+2} = x \in A,$$

получим

$$[\underbrace{b\ldots b}_{j+1}x\underbrace{b\ldots b}_{n-j-3}B] = [B\underbrace{b\ldots b}_{j-1}x\underbrace{b\ldots b}_{n-j-1}],$$

откуда

$$[\overline{b}\underbrace{b\ldots b}_{n-j-3}[\underbrace{b\ldots b}_{j+1}x\underbrace{b\ldots b}_{n-j-3}B]\underbrace{B\ldots B}_{j+1}] =$$

$$= [\overline{b}\underbrace{b\ldots b}_{n-j-3}[B\underbrace{b\ldots b}_{j-1}x\underbrace{b\ldots b}_{n-j-1}]\underbrace{B\ldots B}_{j+1}],$$

$$[\overline{b}\underbrace{b\ldots b}_{n-2}[x\underbrace{b\ldots b}_{n-j-3}\underbrace{B\ldots B}_{j+2}]] = [\overline{b}\underbrace{b\ldots b}_{n-j-3}B\underbrace{b\ldots b}_{j-1}x[\underbrace{b\ldots b}_{n-j-1}\underbrace{B\ldots B}_{j+1}]],$$

$$[x\underbrace{b\ldots b}_{n-j-3}\underbrace{B\ldots B}_{j+2}] = [\overline{b}\underbrace{b\ldots b}_{n-j-3}B\underbrace{b\ldots b}_{j-1}xB],$$

$$[x\underbrace{B\ldots B}_{n-1}] = [\underbrace{B\ldots B}_{n-2}xB].$$

Так как по условию $j \neq n-1$, то $\sigma(n-1) = i$ для некоторого $i \in \{1, \ldots, n-1\}$, $i \neq j+2$. Применяя лемму 5.3.4, получим

$$[\underbrace{B\ldots B}_{i-1}x\underbrace{B\ldots B}_{n-i}] = [\underbrace{B\ldots B}_{n-1}x].$$

Тогда, согласно утверждению 2) леммы 5.3.5, $<B, [\ ]>$ инвариантна в $<A, [\ ]>$.

3) По лемме 5.3.4

$$[x\underbrace{B\ldots B}_{n-1}] = [\underbrace{B\ldots B}_{n-2}xB]$$

для любого $x \in A$. Так как $\sigma(n-2) = 1$, то $\sigma(n-1) = i$ для некоторого $i \in \{2, \ldots, n-1\}$. Снова, применяя лемму 5.3.4, получаем

$$[\underbrace{B\ldots B}_{i-1} x \underbrace{B\ldots B}_{n-i}] = [\underbrace{B\ldots B}_{n-1} x]$$

для любого x ∈ A. Тогда, согласно утверждению 2) леммы 5.3.5, < B, [ ] > инвариантна в < A, [ ] >.

4) По лемме 5.3.4

$$[Bx\underbrace{B\ldots B}_{n-2}] = [\underbrace{B\ldots B}_{n-1} x]$$

для любого x ∈ A. Так как σ(n – 1) = 2, то σ(j) = 1 для некоторого j ∈ {1, …, n – 2}. Если положить j = n – i, то i ∈ {2, …, n – 1}. Снова, применяя лемму 5.3.4, получаем

$$[x\underbrace{B\ldots B}_{n-1}] = [\underbrace{B\ldots B}_{j} x \underbrace{B\ldots B}_{n-j-1}] = [\underbrace{B\ldots B}_{n-i} x \underbrace{B\ldots B}_{i-1}]$$

для любого x ∈ A. Тогда, согласно утверждению 4) леммы 5.3.5, < B, [ ] > инвариантна в < A, [ ] >. ∎

**5.3.7. Следствие.** Если множество Σ содержит подстановку σ, удовлетворяющую по крайней мере одному из условий 1) – 4) теоремы 5.3.6, то Σ-нормальная n-арная подгруппа < B, [ ] > n-арной группы < A, [ ] > является инвариантной в < A, [ ] >.

**5.3.8. Предложение.** σ-Нормальная n-арная подгруппа < B, [ ] > n-арной группы < A, [ ] > является полуинвариантной в < A, [ ] >, если σ(n – 1) = 1.

*Доказательство*. Полагая в (*)

$$x_1 = x \in A, x_2 = \ldots = x_{n-1} = b \in B,$$

и, учитывая равенство σ(n – 1) = 1, получим

$$[x\underbrace{b\ldots b}_{n-2} B] = [B\underbrace{b\ldots b}_{n-2} x],$$

откуда

$$[x\underbrace{B\ldots B}_{n-1}] = [\underbrace{B\ldots B}_{n-2}x],$$

то есть $<B, [\ ]>$ полуинвариантна в $<A, [\ ]>$. ∎

Так как нормальные n-арные подгруппы – это в точности ее σ-нормальные n-арные подгруппы для циклической подстановки $\sigma = (1\ 2\ n-1)$, то из предыдущего предложения вытекает уже отмечавшаяся (предложение 5.2.2) полуинвариантность нормальных n-арных подгрупп.

**5.3.9. Следствие.** Если множество Σ содержит подстановку σ такую, что $\sigma(n-1) = 1$, то Σ-нормальная n-арная подгруппа $<B, [\ ]>$ n-арной группы $<A, [\ ]>$ является полуинвариантной в $<A, [\ ]>$.

Так как $S_2 = \{\sigma_1 = \varepsilon$ – тождественная подстановка, $\sigma_2 = (1\ 2)\}$, то имеет место

**5.3.10. Следствие.** Для любой подстановки $\sigma \in S_2$ всякая σ-нормальная тернарная подгруппа $<B, [\ ]>$ тернарной группы $<A, [\ ]>$ является полуинвариантной в $<A, [\ ]>$, причем:
1) если $\sigma = \varepsilon$, то $<B, [\ ]>$ инвариантна в $<A, [\ ]>$;
2) если $\sigma = (1\ 2)$, то $<B, [\ ]>$ нормальна в $<A, [\ ]>$.

Пусть $<B, [\ ]>$ – σ-нормальная 4-арная подгруппа 4-арной группы $<A, [\ ]>$. Так как $n = 4$, то

$$S_{n-1} = S_3 = \{\sigma_1 = \varepsilon, \sigma_2 = (1\ 2), \sigma_3 = (1\ 3),$$

$$\sigma_4 = (2\ 3), \sigma_5 = (1\ 2\ 3), \sigma_6 = (1\ 3\ 2)\}.$$

Если $\sigma = \sigma_1$, то по предложению 5.3.2 $<B, [\ ]>$ инвариантна в $<A, [\ ]>$.

Если $\sigma = \sigma_2$, то $\sigma(3) = 3$, $\sigma(2) = \sigma(n-2) = 1$. Поэтому, либо, ввиду утверждения 1) теоремы 5.3.6, либо, ввиду утверждения 3) этой же теоремы, $<B, [\ ]>$ инвариантна в $<A, [\ ]>$.

Если $\sigma = \sigma_3$, то $\sigma(2) = 2$, $\sigma(1) = \sigma(3) = \sigma(1 + 2)$ и, согласно

любому из утверждений 1) или 2) теоремы 5.3.6, < B, [ ] > инвариантна в < A, [ ] >.

Если $\sigma = \sigma_4$, то $\sigma(1) = 1$, $\sigma(3) = \sigma(n-1) = 2$. Поэтому, либо, ввиду утверждения 1) теоремы 5.3.6, либо, ввиду утверждения 4) этой же теоремы, < B, [ ] > инвариантна в < A, [ ] >.

Если $\sigma = \sigma_5$, то по предложению 5.3.3 < B, [ ] > инвариантна в < A, [ ] >.

Если $\sigma = \sigma_6$, то $\sigma(1) = \sigma(3) = \sigma(1+2)$, $\sigma(2) = \sigma(n-2) = 1$, $\sigma(3) = \sigma(n-1) = 2$. Поэтому, ввиду любого из утверждений 2), 3) или 4) теоремы 5.3.6, < B, [ ] > инвариантна в < A, [ ] >.

Таким образом, верно

**5.3.11. Следствие.** Для любой подстановки $\sigma \in S_3$ всякая $\sigma$-нормальная 4-арная подгруппа < B, [ ] > 4-арной группы < A, [ ] > является полуинвариантной в < A, [ ] >, причем:

1) если $\sigma \in \{\varepsilon, (1\ 2), (1\ 3), (2\ 3), (1\ 3\ 2)\}$, то < B, [ ] > инвариантна в < A, [ ] >;

2) если $\sigma = (1\ 2\ 3)$, то < B, [ ] > нормальна в < A, [ ] >.

Согласно следствиям 5.3.10 и 5.3.11, если n = 3 или n = 4, то в n-арной группе любая $\sigma$-нормальная n-арная подгруппа является полуинвариантной, и все они являются либо инвариантными, либо нормальными.

Наличие в n-арной группе $\sigma$-нормальных n-арных подгрупп для некоторых подстановок $\sigma$ существенно связано с тождествами, которые могут выполняться в этой n-арной группе. Например, в абелевой n-арной группе любая ее n-арная подгруппа является $\sigma$-нормальной для любой подстановки $\sigma \in S_{n-1}$, а в полуабелевой n-арной группе все ее n-арные подгруппы являются $\sigma$-нормальными для подстановки $\sigma = (1\ 2\ \ldots\ n-1)$, то есть нормальными.

Если в n-арной группе < A, [ ] > выполняется тождество

$$[x_1 x_2 \ldots x_{n-1} x_n] = [x_{\tau(1)} x_{\tau(2)} \ldots x_{\tau(n-1)} x_{\tau(n)}], \qquad (**)$$

где $\tau \in S_n$, и $\tau(1) = n$, то

$$[x_1 x_2 \ldots x_{n-1} B] = [B x_{\tau(2)} \ldots x_{\tau(n-1)} x_{\tau(n)}]$$

для любой n-арной подгруппы $< B, [\ ] >$ из $< A, [\ ] >$. Полагая в последнем равенстве

$$\sigma(1) = \tau(2), \ldots, \sigma(n-2) = \tau(n-1),\ \sigma(n-1) = \tau(n), \quad (***)$$

получим

$$[x_1 x_2 \ldots x_{n-1} B] = [B x_{\sigma(1)} x_{\sigma(2)} \ldots x_{\sigma(n-1)}],$$

то есть $< B, [\ ] >$ $\sigma$-нормальна $< A, [\ ] >$. Таким образом, имеет место

**5.3.12. Предложение.** Если в n-арной группе $< A, [\ ] >$ выполняется тождество (**) для подстановки $\tau \in S_n$, такой, что $\tau(1) = n$, то в $< A, [\ ] >$ все n-арные подгруппы являются $\sigma$-нормальными для подстановки $\sigma \in S_{n-1}$, определяемой равенствами (***).

## §5.4 СОПРЯЖЕННОСТЬ И ПОЛУСОПРЯЖЕННОСТЬ n-АРНЫХ ПОДГРУПП В n-АРНОЙ ГРУППЕ

Так как по теореме Поста о смежных классах всякая n-арная группа изоморфно вкладывается в n-арную группу, производную от группы, то естественно предположить наличие связи между сопряженностью и полусопряженностью n-арных подгрупп в n-арной группе и сопряженностью подгрупп в группе, к которой приводима n-арная группа согласно теореме Поста о смежных классах. Изучению указанной связи и посвящен данный параграф.

Напомним, что подмножество C n-арной группы $< A, [\ ] >$ называется [4] сопряженным в ней посредством последовательности $x_1 \ldots x_i$, $x_i \in A$ с подмножеством B, если

$$B = [x_1 \ldots x_i C y_1 \ldots y_j],$$

где $y_1 \ldots y_j$ – обратная последовательность для последовательности $x_1 \ldots x_i$ (определение 2.4.1). Если же

$$[x\underbrace{C \ldots C}_{n-1}] = [\underbrace{B \ldots B}_{n-1}x],$$

то подмножество $< C, [\ ] >$ n-арной группы $< A, [\ ] >$ называется [30] полусопряженным в ней посредством элемента $x \in A$ с подмножеством $< B, [\ ] >$ (определение 2.4.8). Ясно, что при n = 2 понятия сопряженности и полусопряженности совпадают.

В §2.2 для подмножества B n-арной группы $< A, [\ ] >$ определены подмножества $B^{(i)}(A)$, где $i \in \{1, 2, \ldots n-1\}$, $B_o(A)$ и $B^*(A)$ (определение 2.2.18) и доказано, что: 1) $B^*(A)$ – подполугруппа группы $A^*$; 2) если $< B, [\ ] >$ – n-арная подгруппа n-арной группы $< A, [\ ] >$, то $B^*(A)$ – подгруппа группы $A^*$, изоморфная группе $B^*$, а $B_o(A)$ – инвариантная подгруппа группы $B^*(A)$, изоморфная группе $B_o$ (теорема 2.2.19).

Заметим, что определение множества $B^{(i)}(A)$ можно расширить, считая i = 1, 2, … .

Справедливость следующего предложения устанавливается простой проверкой.

**5.4.1. Предложение**. Справедливы следующие утверждения:

1) $\bigcup\limits_{j=1}^{\infty} B^{(j(n-1))}(A)$ – подполугруппа группы $A_o$;

2) если $< B, [\ ] >$ – n-арная подполугруппа в $< A, [\ ] >$, то

$$B_o(A) = \bigcup\limits_{j=1}^{\infty} B^{(j(n-1))}(A).$$

**5.4.2. Лемма** [107, Воробьев Г.Н.]. Если подмножество C n-арной группы $< A, [\ ] >$ сопряжено в ней посредством по-следовательности $\alpha = x_1 \ldots x_i$ с подмножеством B, то подполугруппа C*(A) сопряжена в группе A* посредством элемента $\theta_A(\alpha)$ с подполугруппой B*(A).

*Доказательство.* По условию леммы

$$B = [x_1 \ldots x_i C\, y_1 \ldots y_j], \qquad (1)$$

где $\beta = y_1 \ldots y_j$ – обратная последовательность для последовательности $\alpha$.

Пусть

$$u = \theta_A(\alpha)\theta_A(c_1 \ldots c_k)\theta_A^{-1}(\alpha)$$

– произвольный элемент из множества $\theta_A(\alpha)C^*(A)\theta_A^{-1}(\alpha)$, где $c_1, \ldots, c_k \in C$, $k \geq 1$. Так как $\theta_A^{-1}(\alpha) = \theta_A(\beta)$, то, используя (1), получим

$$u = \theta_A(\alpha)\theta_A(c_1 \ldots c_k)\theta_A(\beta) = \theta_A(\alpha c_1 \ldots c_k \beta) =$$
$$= \theta_A(\alpha c_1 \beta \alpha c_2 \beta \alpha \ldots c_k \beta) =$$
$$= \theta_A([\alpha c_1 \beta][\alpha c_2 \beta] \ldots [\alpha c_k \beta]) = \theta_A(b_1 b_2 \ldots b_k)$$

для некоторых $b_1, \ldots, b_k \in B$. Следовательно, $u \in B^*(A)$ и верно включение

$$\theta_A(\alpha)C^*(A)\theta_A^{-1}(\alpha) \subseteq B^*(A). \qquad (2)$$

Если теперь $v = \theta_A(b_1 \ldots b_k)$ – произвольный элемент из B*(A), где $b_1, \ldots, b_k \in B$, $k \geq 1$, то, используя (1), получим

$$v = \theta_A(b_1 \ldots b_k) = \theta_A([\alpha c_1 \beta][\alpha c_2 \beta] \ldots [\alpha c_k \beta]) =$$
$$= \theta_A(\alpha c_1 \beta \alpha c_2 \beta \ldots \alpha c_k \beta) = \theta_A(\alpha c_1 \ldots c_k \beta) =$$
$$= \theta_A(\alpha)\theta_A(c_1 \ldots c_k)\theta_A(\beta) = \theta_A(\alpha)\theta_A(c_1 \ldots c_k)\theta_A^{-1}(\alpha)$$

для некоторых $c_1, \ldots, c_k \in C$. Поэтому $v \in \theta_A(\alpha)C^*(A)\theta_A^{-1}(\alpha)$ и верно включение

$$B^*(A) \subseteq \theta_A(\alpha)C^*(A)\theta_A^{-1}(\alpha). \qquad (3)$$

Из (2) и (3) следует равенство

$$B^*(A) = \theta_A(\alpha)C^*(A)\theta_A^{-1}(\alpha). \qquad (4)$$

Таким образом, подполугруппа $C^*(A)$ сопряжена в группе $A^*$ посредством элемента $\theta_A(\alpha)$ с подполугруппой $B^*(A)$. ∎

**5.4.3. Теорема** [107, Воробьев Г.Н.]. n-Арная подполугруппа $< C, [\ ] >$ n-арной группы $< A, [\ ] >$ сопряжена в ней посредством последовательности $\alpha$ с n-арной подполугруппой $< B, [\ ] >$ тогда и только тогда, когда подполугруппа $C^*(A)$ сопряжена в группе $A^*$ посредством элемента $\theta_A(\alpha)$ с подполугруппой $B^*(A)$.

*Доказательство.* *Необходимость.* Лемма 5.4.2.

*Достаточность.* По условию верно (4) из леммы 5.4.2. Если $u = [\alpha c \beta]$ – произвольный элемент из $[\alpha C \beta]$, где $c \in C$, $\beta$ – обратная последовательность для $\alpha$, то согласно (4) из леммы 5.4.2,

$$\theta_A(\alpha)\theta_A(c)\theta_A^{-1}(\alpha) \in B^*(A),$$

откуда

$$\theta_A(\alpha)\theta_A(c)\theta_A(\beta) = \theta_A(\alpha c \beta) \in B^*(A).$$

Из последнего соотношения следует

$$\theta_A(\alpha c \beta) = \theta_A(b_1 \ldots b_k)$$

для некоторых $b_1, \ldots, b_k \in B$, $k \geq 1$. Так как $< B, [\ ] >$ – n-арная подполугруппа, то $[b_1 \ldots b_k] = b \in B$, откуда и из предыдущего равенства получаем $[\alpha c \beta] = b \in B$ и верно включение

$$[\alpha C \beta] \subseteq B. \qquad (1)$$

Если теперь $b$ – произвольный элемент из $B$, то согласно (4) из леммы 5.4.2,

$$\theta_A(b) \in \theta_A(\alpha) C^*(A) \theta_A^{-1}(\alpha),$$

то есть

$$\theta_A(b) = \theta_A(\alpha)\theta_A(c_1 \ldots c_k)\theta_A(\beta)$$

для некоторых $c_1, \ldots, c_k \in C$, $k \geq 1$, откуда

$$\theta_A(b) = \theta_A(\alpha c_1 \ldots c_k \beta),\ b = [\alpha c_1 \ldots c_k \beta].$$

Так как длина последовательности $\alpha\beta$ кратна $n-1$, то из последнего равенства следует $k \equiv 1 \bmod(n-1)$, а так как $< C, [\ ] >$ n-арная подполугруппа, то

$$[c_1 \ldots c_k] = c \in C.$$

Следовательно,

$$b = [\alpha c \beta] \in [\alpha C \beta]$$

и верно включение

$$B \subseteq [\alpha C \beta]. \tag{2}$$

Из (1) и (2) следует равенство $B = [\alpha C \beta]$. Таким образом, n-арная подполугруппа $< C, [\ ] >$ сопряжена в n-арной группе $< A, [\ ] >$ посредством последовательности $\alpha$ с n-арной подполугруппой $< B, [\ ] >$. ∎

**5.4.4. Следствие** [107, Воробьев Г.Н.]. n-Арная подгруппа $< C, [\ ] >$ n-арной группы $< A, [\ ] >$ сопряжена в ней посредством последовательности $\alpha$ с n-арной подгруппой $< B, [\ ] >$ тогда и только тогда, когда подгруппа $C^*(A)$ сопряжена в группе $A^*$ посредством элемента $\theta_A(\alpha)$ с подгруппой $B^*(A)$.

**5.4.5. Предложение** [107, Воробьев Г.Н.]. Если подмножество $C$ n-арной группы $< A, [\ ] >$ полусопряжено в ней посредством элемента $x$ с подмножеством $B$, то подмножество $C^{j(n-1)}$ сопряжено в группе $A^*$ посредством элемента $\theta_A(x)$ с подмножеством $B^{j(n-1)}$ для любого $j \geq 1$.

**Доказательство.** По условию предложения

$$[x\underbrace{C \ldots C}_{n-1}] = [\underbrace{B \ldots B}_{n-1}x]. \qquad (1)$$

Пусть $u = \theta_A(x)\theta_A(c_1 \ldots c_{j(n-1)})$ – произвольный элемент из $\theta_A(x)C^{j(n-1)}$, где $c_1, \ldots, c_{j(n-1)} \in C$, $j \geq 1$. Используя $j$ раз (1), получим

$$u = \theta_A(x)\theta_A(c_1 \ldots c_{j(n-1)}) =$$
$$= \theta_A([xc_1 \ldots c_{n-1}])\theta_A(c_n \ldots c_{j(n-1)}) =$$
$$= \theta_A([b_1 \ldots b_{n-1}x])\theta_A(c_n \ldots c_{j(n-1)}) =$$
$$= \theta_A(b_1 \ldots b_{n-1})\theta_A([xc_n \ldots c_{2(n-1)}])\theta_A(c_{2(n-1)+1} \ldots c_{j(n-1)}) =$$
$$= \theta_A(b_1 \ldots b_{n-1})\theta_A([b_n \ldots b_{2(n-1)}x])\theta_A(c_{2(n-1)+1} \ldots c_{j(n-1)}) =$$
$$= \theta_A(b_1 \ldots b_{2(n-1)})\theta_A(xc_{2(n-1)+1} \ldots c_{j(n-1)}) = \ldots$$
$$\ldots = \theta_A(b_1 \ldots b_{j(n-1)}x) = \theta_A(b_1 \ldots b_{j(n-1)})\theta_A(x)$$

для некоторых $b_1, \ldots, b_{j(n-1)} \in B$. Поэтому $u \in B^{j(n-1)}\theta_A(x)$ и верно включение

$$\theta_A(x)C^{j(n-1)} \subseteq B^{j(n-1)}\theta_A(x). \qquad (2)$$

Аналогично доказывается включение

$$B^{j(n-1)}\theta_A(x) \subseteq \theta_A(x)C^{j(n-1)}. \qquad (3)$$

Из (2) и (3) следует равенство

$$\theta_A(x)C^{j(n-1)} = B^{j(n-1)}\theta_A(x).$$

Таким образом, подмножество $C^{j(n-1)}$ сопряжено в группе $A^*$ посредством элемента $\theta_A(x)$ с подмножеством $B^{j(n-1)}$. ∎

**5.4.6. Следствие** [107, Воробьев Г.Н.]. Если подмножест-

во C n-арной группы $< A, [\ ] >$ полусопряжено в ней посредством элемента x с подмножеством B, то подполугруппа $\bigcup_{j=1}^{\infty} C^{(j(n-1))}(A)$ сопряжена в группе $A^*$ посредством элемента $\theta_A(x)$ с подполугруппой $\bigcup_{j=1}^{\infty} B^{(j(n-1))}(A)$.

***Доказательство.*** Так как

$$\theta_A(x)(\bigcup_{j=1}^{\infty} C^{(j(n-1))}(A)) = \bigcup_{j=1}^{\infty} \theta_A(x) C^{j(n-1)}(A) =$$

$$= \bigcup_{j=1}^{\infty} B^{j(n-1)}(A)\theta_A(x) = (\bigcup_{j=1}^{\infty} B^{(j(n-1))}(A))\theta_A(x),$$

то подполугруппа $\bigcup_{j=1}^{\infty} C^{(j(n-1))}(A)$ сопряжена в группе $A^*$ посредством элемента $\theta_A(x)$ с подполугруппой $\bigcup_{j=1}^{\infty} B^{(j(n-1))}(A)$. ∎

В §2.4 установлена связь между полусопряженностью n-арных подгрупп в n-арной группе $< A, [\ ] >$ и сопряженностью подгрупп в группе $< A, \text{\textcircled{a}} >$, к которой приводима n арная группа $< A, [\ ] >$ согласно теореме Глускина-Хоссу. Так как группы $< A, \text{\textcircled{a}} >$ и $A_o$ изоморфны, то должна существовать подобная связь между полусопряженностью n-арных подгрупп в n-арной группе $< A, [\ ] >$ и сопряженностью подгрупп в группе $A_o$. Следующая теорема устанавливает подобную связь.

**5.4.7. Теорема** [107, Воробьев Г.Н.]. Пусть $< A, [\ ] >$ – n-арная группа, B и C ее n-арные подгруппы, $B \cap C \neq \varnothing$, и пусть $a_1, \ldots, a_{n-2} \in B \cap C$. Тогда $< C, [\ ] >$ полусопряжена в $< A, [\ ] >$ посредством элемента x с $< B, [\ ] >$ тогда и только тогда, когда подгруппа $C_o(A)$ сопряжена в группе $A_o$ посредством элемента $\theta_A(xa_1 \ldots a_{n-2})$ с подгруппой $B_o(A)$.

***Доказательство.*** *Необходимость*. По условию теоремы

верно

$$[x\underbrace{C \ldots C}_{n-1}] = [\underbrace{B \ldots B}_{n-1}x].$$

Тогда, используя предложение 5.4.5 при $j = n - 1$, а также легко проверяемое равенство

$$C_o(A) = C^{(1)}(A)\theta_A(a_1 \ldots a_{n-2}) = \theta_A(a_1 \ldots a_{n-2})C^{(1)}(A),$$

получим

$$\theta_A(xa_1 \ldots a_{n-2})C_o(A) = \theta_A(x)\theta_A(a_1 \ldots a_{n-2})C^{(1)}(A)\theta_A(a_1 \ldots a_{n-2}) =$$

$$= \theta_A(x)C_o(A)\theta_A(a_1 \ldots a_{n-2}) = B_o(A)\theta_A(x)\theta_A(a_1 \ldots a_{n-2}) =$$

$$= B_o(A)\theta_A(xa_1 \ldots a_{n-2}),$$

то есть

$$\theta_A(xa_1 \ldots a_{n-2})C_o(A) = B_o(A)\theta_A(xa_1 \ldots a_{n-2}). \qquad (1)$$

Следовательно подгруппа $C_o(A)$ сопряжена в группе $A_o$ посредством элемента $\theta_A(xa_1 \ldots a_{n-2})$ с подгруппой $B_o(A)$.

*Достаточность.* По условию теоремы верно (1).

Пусть $u = [xc_1 \ldots c_{n-1}]$ – произвольный элемент из

$$[x\underbrace{C \ldots C}_{n-1}], c_1, \ldots, c_{n-1} \in C.$$

Так как $< C, [\ ] >$ – n-арная подгруппа в $< A, [\ ] >$, то элемент $u$ можно представить в виде

$$u = [xa_1 \ldots a_{n-2}c], c \in C,$$

а элемент с в свою очередь может быть представлен в виде

$$c = [d_1 \ldots d_{n-1}a], d_1, \ldots, d_{n-1} \in C,$$

a – обратный элемент для последовательности $a_1 \ldots a_{n-2}$. Поэтому, используя (1) и нейтральность последовательности $a_1 \ldots a_{n-2}a$, получим

$$\theta_A(u) = \theta_A([xc_1 \ldots c_{n-1}]) = \theta_A([xa_1 \ldots a_{n-2}c]) =$$

$$= \theta_A(xa_1 \ldots a_{n-2}[d_1 \ldots d_{n-1}a]) =$$

$$= \theta_A(xa_1 \ldots a_{n-2})\theta_A(d_1 \ldots d_{n-1})\theta_A(a) =$$

$$= \theta_A(b_1 \ldots b_{n-1})\theta_A(xa_1 \ldots a_{n-2})\theta_A(a) =$$

$$= \theta_A([b_1 \ldots b_{n-1}x])\theta_A(a_1 \ldots a_{n-2}a) = \theta_A([b_1 \ldots b_{n-1}x])$$

для некоторых $b_1, \ldots, b_{n-1} \in B$. Таким образом,

$$\theta_A(u) = \theta_A([b_1 \ldots b_{n-1}x]),$$

откуда

$$u = [b_1 \ldots b_{n-1}x] \in [\underbrace{B \ldots B}_{n-1}x]$$

и доказано включение

$$[x\underbrace{C \ldots C}_{n-1}] \subseteq [\underbrace{B \ldots B}_{n-1}x]. \qquad (2)$$

Пусть теперь $v = [b_1 \ldots b_{n-1}x]$ – произвольный элемент из

$$[\underbrace{B \ldots B}_{n-1}x], b_1, \ldots, b_{n-1} \in B.$$

Тогда

$$v = [b_1 \ldots b_{n-1}xa_1 \ldots a_{n-2}a],$$

где a – обратный элемент для последовательности $a_1 \ldots a_{n-2}$. Снова используя (1), получим

$$\theta_A(v) = \theta_A([b_1 \ldots b_{n-1}xa_1 \ldots a_{n-2}a]) =$$

$$= \theta_A(b_1 \ldots b_{n-1})\theta_A(xa_1 \ldots a_{n-2})\theta_A(a) =$$

$$= \theta_A(xa_1 \ldots a_{n-2})\theta_A(c_1 \ldots c_{n-1})\theta_A(a) =$$

$$= \theta_A([xa_1 \ldots a_{n-2}[c_1 \ldots c_{n-1}a]])$$

для некоторых $c_1, \ldots, c_{n-1} \in C$. Таким образом,

$$\theta_A(v) = \theta_A([xa_1 \ldots a_{n-2}[c_1 \ldots c_{n-1}a]]),$$

откуда, учитывая $a \in C$, получаем

$$v = [xa_1 \ldots a_{n-2}[c_1 \ldots c_{n-1}a]] \in [x\underbrace{C \ldots C}_{n-1}]$$

и доказано включение

$$[\underbrace{B \ldots B}_{n-1} x] \subseteq [x \underbrace{C \ldots C}_{n-1}] \qquad (3)$$

Из (2) и (3) следует равенство

$$[x\underbrace{C \ldots C}_{n-1}] = [\underbrace{B \ldots B}_{n-1} x].$$

Таким образом, n-арная подгруппа $< C, [\ ] >$ полусопряжена в $< A, [\ ] >$ посредством элемента x с n-арной подгруппой $< B, [\ ] >$. ∎

**5.4.8. Следствие** [107, Воробьев Г.Н.]. Пусть $< A, [\ ] >$ – n-арная группа, B и C ее n-арные подгруппы, $B \cap C \neq \varnothing$, $a \in B \cap C$. Тогда $< C, [\ ] >$ полусопряжена в $< A, [\ ] >$ посредством элемента x с $< B, [\ ] >$ тогда и только тогда, когда подгруппа $C_o(A)$ сопряжена в группе $A_o$ посредством элемента $\theta_A(x\underbrace{a \ldots a}_{n-2})$ с подгруппой $B_o(A)$.

Аналогично теореме 5.4.7 доказывается следующая

**5.4.9. Теорема** [107, Воробьев Г.Н.]. Пусть $< A, [\ ] >$ – n-арная группа, B и C ее n-арные подгруппы, $B \cap C \neq \varnothing$, и пусть $a_1, \ldots, a_{n-2} \in B \cap C$. Тогда $< C, [\ ] >$ полусопряжена в $< A, [\ ] >$ посредством элемента x с $< B, [\ ] >$ тогда и только тогда, когда подгруппа $C_o(A)$ сопряжена в группе $A_o$ посред-ством элемента $\theta_A(a_1 \ldots a_{n-2}x)$ с подгруппой $B_o(A)$.

**5.4.10. Следствие** [107, Воробьев Г.Н.]. Пусть $< A, [\ ] > -$ n-арная группа, B и C ее n-арные подгруппы, $B \cap C \neq \varnothing$, $a \in B \cap C$. Тогда $< C, [\ ] >$ полусопряжена в $< A, [\ ] >$ посредством элемента x с $< B, [\ ] >$ тогда и только тогда, когда подгруппа $C_o(A)$ сопряжена в группе $A_o$ посредством элемента $\theta_A(\underbrace{a \ldots a}_{n-2} x)$ с подгруппой $B_o(A)$.

## ДОПОЛНЕНИЯ И КОММЕНТАРИИ

**1.** Введём ещё один n-арный аналог нормальных подгрупп, расширяющий понятие нормальности для n-арных подгрупп из §5.2.

**Определение.** n-Арная подгруппа $< B, [\ ] >$ n-арной группы $< A, [\ ] >$ называется *слабо нормальной* в ней, если

$$[\underbrace{x \ldots x}_{n-1} B] = [B \underbrace{x \ldots x}_{n-1}]$$

для любого $x \in A$.

Ясно, что в идемпотентной n-арной группе, т. е. в n-арной группе, все элементы которой являются идемпотентами, любая её n-арная подгруппа является слабо нормальной. Слабо нормальными будут и все n-арные подгруппы в слабо полуабелевой n-арной группе, в том числе и в полуабелевой n-арной группе.

**Предложение 1.** Если n-арная подгруппа $< B, [\ ] >$ нормальна в n-арной группе $< A, [\ ] >$, то она и слабо нормальна в $< A, [\ ] >$.

*Доказательство.* Так как $< B, [\ ] >$ нормальна в $< A, [\ ] >$, то, полагая в определении 5.2.1 $x = x_1 = \ldots = x_{n-2}$, получаем

$$[x \underbrace{x \ldots x}_{n-2} B] = [B \underbrace{x \ldots x}_{n-2} x]$$

Следовательно, $< B, [\ ] >$ слабо нормальна в $< A, [\ ] >$. ∎

**Предложение 2.** Определим на группе A n-арную операцию

$$[x_1 \ldots x_n] = x_1 \ldots x_n a,$$

где $a \in Z(A)$. Тогда:

1) если $n = \kappa m + 1 \geq 3$, где m – экспонента группы A, то в $< A, [\ ] >$ все n-арные подгруппы являются слабо нормальными;

2) если подгруппа B группы A содержит a и не является нормальной в группе A, то n-арная подгруппа $< B, [\ ] >$ не является полуинвариантной в $< A, [\ ] >$.

***Доказательство.*** 1) Так как $a \in Z(A)$, то $< A, [\ ] >$ – n-арная группа. Если $< B, [\ ] >$ – n-арная подгруппа в $< A, [\ ] >$, $x \in A$, то

$$[\underbrace{x \ldots x}_{n-1} B] = \underbrace{x \ldots x}_{n-1} Ba = x^{\kappa m} Ba = Ba,$$

$$[B \underbrace{x \ldots x}_{n-1}] = B \underbrace{x \ldots x}_{n-1} a = B x^{\kappa m} a = Ba,$$

откуда

$$[\underbrace{x \ldots x}_{n-1} B] = [B \underbrace{x \ldots x}_{n-1}].$$

Следовательно, $< B, [\ ] >$ слабо нормальна в $< A, [\ ] >$.

2) Так как $a \in B$, $a \in Z(A)$, то $b \in Z(A)$. Поэтому $< B, [\ ] >$ – n-арная подгруппа в $< A, [\ ] >$. Так как B не является нормальной в A, то существует такой $x \in A$, что $xB \neq Bx$.

Если теперь $< B, [\ ] >$ полуинвариантна в $< A, [\ ] >$, то

$$[(xa^{-1}) \underbrace{B \ldots B}_{n-1}] = [\underbrace{B \ldots B}_{n-1} (xa^{-1})],$$

$$xa^{-1} \underbrace{B \ldots B}_{n-1} a = \underbrace{B \ldots B}_{n-1} xa^{-1} a,$$

$$xa^{-1} Ba = Bx.$$

Так как a, $a^{-1} \in B$, то из последнего равенства следует $xB = Bx$, что противоречит выбору x. Следовательно, $< B, [\ ] >$ – не является полуинвариантной в $< A, [\ ] >$. ∎

Покажем, что в n-арной группе могут быть слабо нормальные n-арные подгруппы, которые не являются полуинвариантными.

**Пример 1.** Пусть $< S_3, [\,] >$ – 7-арная группа, производная от симметрической группы $S_3$. Так как в $S_3$ подгруппы $B = \{e, (12)\}$, $C = \{e, (13)\}$, $D = \{e, (23)\}$ не являются нормальными, то, согласно, утверждения 2) предыдущего предложения, в идемпотентной 7-арной группе $< S_3, [\,] >$ 7-арные подгруппы второго порядка $< B, [\,] >$, $< C, [\,] >$ и $< D, [\,] >$ не являются полуинвариантными. Согласно же 1) того же предложения, 7-арные подгруппы $< B, [\,] >$, $< C, [\,] >$ и $< D, [\,] >$ являются слабо нормальными.

Существование n-арных групп, обладающих полуинвариантными n-арными подгруппами, не являющимися полунормальными, следует из примера 5.2.7. Тернарная подгруппа B из этого примера является полуинвариантной в $< D_6, [\,] >$, но не является слабо нормальной в $< D_6, [\,] >$, так как $[ccB] \neq [Bcc]$.

Следующий пример показывает, что в n-арной группе n-арные подгруппы могут быть одновременно и слабо нормальными и полуинвариантными.

**Пример 2.** Пусть $< R, [\,] >$ – 5-арная группа С. А. Русакова (пример 1.2.8). Так как экспонента группы R равна 4 и $Z(R) = \{1, a^2\}$, то согласно 1) предложения 2, в $< R, [\,] >$ все 5-арные подгруппы слабо нормальны. Учитывая, что в $< R, [\,] >$ все 5-арные подгруппы являются и полувариантными [4], видим, что в 5-арной группе $< R, [\,] >$ все 5-арные подгруппы являются одновременно и слабо нормальными и полуинвариантными.

Отметим, что $< R, [\,] >$ не является ни идемпотентной ни полуабелевой.

**2.** Приведем еще одно доказательство достаточности следствия 5.3.6.

Так как для тождественной подстановки ε верно $\varepsilon(i) = i$ для любого $i \in \{1, \ldots, n-1\}$, то, применяя лемму 5.3.3 последовательно для $i = 1, \ldots, n-1$, получим

$$[x\underbrace{B \ldots B}_{n-1}] = [Bx\underbrace{B \ldots B}_{n-2}],$$

$$[Bx\underbrace{B \ldots B}_{n-2}] = [BBx\underbrace{B \ldots B}_{n-3}],$$

$$\ldots\ldots\ldots\ldots\ldots\ldots\ldots\ldots\ldots\ldots\ldots\ldots$$

$$[\underbrace{B\ldots B}_{n-3} xBB] = [\underbrace{B\ldots B}_{n-2} xB],$$

$$[\underbrace{B\ldots B}_{n-2} xB] = [\underbrace{B\ldots B}_{n-1} x],$$

то есть

$$[x\underbrace{B\ldots B}_{n-1}] = [Bx\underbrace{B\ldots B}_{n-2}] = \ldots = [\underbrace{B\ldots B}_{n-2} xB] = [\underbrace{B\ldots B}_{n-1} x]$$

для любого x ∈ A. Следовательно, < B, [ ] > инвариантна в < A, [ ] >.

**3.** Известно, что множество всех инвариантных (нормальных) n-арных подгрупп n-арной группы, содержащих фиксированный элемент, образует подрешетку решетки всех n-арных подгрупп этой n-арной группы (теорема 2.3.21, теорема 5.2.16). Поэтому, если Σ = {ε}, Σ = {(1 2 … n – 1)} или Σ = {ε, (1 2 … n – 1)}, то множество всех Σ-нормальных n-арных подгрупп n-арной группы, содержащих фиксированный элемент, образует подрешетку решетки всех n-арных подгрупп этой n-арной группы.

**4.** Результаты §5.3 могут быть распространены на m-полуинвариантные n-арные подгруппы.

**Вопрос.** Для каких еще множеств Σ, отличных от указанных выше, множество всех Σ-нормальных n-арных подгрупп n-арной группы, содержащих фиксированный элемент, образует подрешетку решетки всех n-арных подгрупп этой n-арной группы?

**5.** Сопряженные n-арные подгруппы впервые появились у Поста [3], показавшего, что силовские n-арные подгруппы n-арной группы сопряжены в ней при довольно сильном ограничении на n. Обобщая результаты Поста, С.А. Русаков при аналогичном ограничении на n изучал [4] сопряженность холловых n-арных подгрупп в n-арной группе. Для устранения отмеченных выше ограничений Г.Н. Воробьев в работе [30] определил в n-арной группе полусопряженные n-арные подгруппы и получил n-арные аналоги теоремы Силова и теоремы Холла-Чунихина без каких либо ограничений на арность операции.

# Г Л А В А  6

# СМЕЖНЫЕ КЛАССЫ

В данной главе строится представление n-арной группы подстановками на смежных классах и изучается связь между разложениями n-арной группы по ее n-арной подгруппе и соответствующими разложениями в универсальной обертывающей группе Поста.

## §6.1. ПРЕДСТАВЛЕНИЕ n-АРНОЙ ГРУППЫ ПОДСТАНОВКАМИ НА СМЕЖНЫХ КЛАССАХ

Пусть $<A, [\ ]>$ – n-арная группа, $<B, [\ ]>$ – её n-арная подгруппа. Для любых элементов $a_1, \ldots, a_{n-1} \in A$ определим преобразование

$$\delta_{a_1 \ldots a_{n-1}} : [\underbrace{B \ldots B}_{n-1} x] \to [[\underbrace{B \ldots B}_{n-1} x] a_1 \ldots a_{n-1}]$$

множества

$$\Omega = \{[\underbrace{B \ldots B}_{n-1} x] \mid x \in A\}$$

всех смежных классов $<A, [\ ]>$ по $<B, [\ ]>$ вида $[\underbrace{B \ldots B}_{n-1} x]$ и положим

$$\Delta = \{\delta_{a_1 \ldots a_{n-1}} \mid a_1, \ldots, a_{n-1} \in A\}.$$

Если $e_1 \ldots e_{n-1}$ – нейтральная последовательность, то

$$\delta_{e_1 \ldots e_{n-1}} : [\underbrace{B \ldots B}_{n-1} x] \to [\underbrace{B \ldots B}_{n-1} x].$$

Поэтому Δ содержит тождественное преобразование.

Зафиксируем элемент b ∈ B, и пусть $\widetilde{b}$ – обратная последовательность для b. Так как последовательность $a_1 \ldots a_{n-1}$ эквивалентна в n-арной группе < A, [ ] > последовательности $\widetilde{b}\,a$ для некоторого a ∈ A, то

$$\delta_{a_1 \ldots a_{n-1}} = \delta_{\widetilde{b}a} : [\underbrace{B\ldots B}_{n-1}x] \to [[\underbrace{B\ldots B}_{n-1}x]\widetilde{b}\,a]$$

$$\Delta = \{\delta_{\widetilde{b}a} \mid a \in A\}.$$

Для сокращения записей положим $\delta_{\widetilde{b}a} = \delta_a$. Тогда

$$\Delta = \{\delta_a \mid a \in A\}.$$

Ясно, что $\delta_b$ – тождественное преобразование.

Определим отображение γ: A → Δ по правилу

$$\gamma: a \to a^\gamma = \delta_a.$$

**6.1.1. Теорема** [108]. Пусть < A, [ ] > n-арная группа, < B, [ ] > – её n-арная подгруппа конечного индекса, b ∈ B. Тогда справедливы следующие утверждения:

1) Δ – транзитивная группа подстановок на Ω, при этом γ – гомоморфизм группы < A, ⓑ > на группу Δ;

2) на множестве Δ можно так определить n-арную операцию ⟨ ⟩, что <Δ, ⟨ ⟩ > n-арная группа, при этом γ – гомоморфизм < A, [ ] > на <Δ, ⟨ ⟩ >.

*Доказательство.* Напомним (§1.5), что операция ⓑ на множестве A определяется по правилу

$$x \;ⓑ\; y = [x\,\widetilde{b}\,y],$$

где $\widetilde{b}$, как отмечалось выше, обратная последовательность для элемента b, с помощью которой определяется отображение

$$\alpha : x \to x^\alpha = [bx\widetilde{b}\,].$$

По теореме Глускина-Хоссу $< A, \textcircled{b} >$ – группа, $< B, \textcircled{b} >$ – ее подгруппа, $\alpha$ – автоморфизм $< A, \textcircled{b} >$, сужение которого на B является автоморфизмом $< B, \textcircled{b} >$, и выполняются следующие условия:

$$[x_1 x_2 \ldots x_n] = x_1 \textcircled{b} x_2^\alpha \textcircled{b} \ldots \textcircled{b} x_n^{\alpha^{n-1}} \textcircled{b} c, \quad x_1, x_2, \ldots, x_n \in A; \quad (1)$$

$$c^\alpha = c; \quad (2)$$

$$x_n^{\alpha^{n-1}} \textcircled{b} c = c \textcircled{b} x, \quad x \in A, \quad (3)$$

где

$$c = [\underbrace{b \ldots b}_{n}].$$

1) Так как

$$[\underbrace{B \ldots B}_{n-1} x] = \{[s\widetilde{b} x] \mid s \in B\} = \{s \textcircled{b} x \mid s \in B\} = B \textcircled{b} x,$$

то

$$[\underbrace{B \ldots B}_{n-1} x] = B \textcircled{b} x. \quad (4)$$

Следовательно, множество $\Omega$ совпадает с множеством всех правых смежных классов $< A, \textcircled{b} >$ по $< B, \textcircled{b} >$, то есть

$$\Omega = \{[\underbrace{B \ldots B}_{n-1} x] \mid x \in A\} = \{B \textcircled{b} x \mid x \in A\}.$$

Так как

$$[[\underbrace{B \ldots B}_{n-1} x]\widetilde{b} a] = [(B \textcircled{b} x)\widetilde{b} a] = (B \textcircled{b} x) \textcircled{b} a,$$

то

$$[[\underbrace{B \ldots B}_{n-1} x]\widetilde{b} a] = (B \textcircled{b} x) \textcircled{b} a. \quad (5)$$

Из (4), (5) и определения преобразования $\delta_a$ получаем

$$\delta_a : B \circledb x \to (B \circledb x) \circledb a,$$

$$\Delta = \{\delta_a : B \circledb x \to (B \circledb x) \circledb a \mid a \in A\}.$$

По известной теореме для бинарных групп (см., например, [109]) $\Delta$ – транзитивная группа подстановок на $\Omega$, при этом $\gamma$ – гомоморфизм группы $< A, \circledb >$ на группу $\Delta$, ядром которого является подгруппа $< N, \circledb >$, где

$$N = \bigcap_{x \in A} x^{-1} \circledb B \circledb x,$$

причем $N \subseteq B$.

2) Пусть x – произвольный элемент из A. По теореме 2.2.3 обратный элемент в группе $< A, \circledb >$ для элемента x имеет вид

$$x^{-1} = [b\,\widetilde{x}\,b], \qquad (6)$$

где $\widetilde{x}$ – обратная последовательность для x.

Так как

$$x^{-1} \circledb B \circledb x = [[b\,\widetilde{x}\,b]\,\widetilde{b}\,B\,\widetilde{b}\,x] = [b\,\widetilde{x}\,B\,\widetilde{b}\,x],$$

то

$$x^{-1} \circledb B \circledb x = [b\,\widetilde{x}\,\underbrace{B\ldots B}_{n-1}\,x]. \qquad (7)$$

Для элемента $y = [bx\,\widetilde{b}]$ обратный элемент $y^{-1}$ в группе $< A, \circledb >$, ввиду (6), может быть представлен в виде

$$y^{-1} = [b\,\widetilde{y}\,b] = [bb\,\widetilde{x}\,\widetilde{\widetilde{b}}\,b] = [bb\,\widetilde{x}\,],$$

то есть

$$y^{-1} = [bb\,\widetilde{x}\,]. \qquad (8)$$

Так как, с одной стороны, ввиду (7),

$$(x^{-1} \text{\textcircled{b}} B \text{\textcircled{b}} x)^\alpha = [b[b\tilde{x} \underbrace{B\ldots B}_{n-1} x]\tilde{b}] = [bb\tilde{x} \underbrace{B\ldots B}_{n-1} x\tilde{b}],$$

а, с другой стороны, ввиду (8),

$$y^{-1} \text{\textcircled{b}} B \text{\textcircled{b}} y = [[bb\tilde{x}]\tilde{b} B \tilde{b}[bx\tilde{b}]] =$$
$$= [bb\tilde{x}\,\tilde{b} Bx\tilde{b}] = [bb\tilde{x} \underbrace{B\ldots B}_{n-1} x\tilde{b}],$$

то

$$(x^{-1} \text{\textcircled{b}} B \text{\textcircled{b}} x)^\alpha = y^{-1} \text{\textcircled{b}} B \text{\textcircled{b}} y,$$

где $y = [bx\tilde{b}]$. Кроме того, из разрешимости в $< A, [\ ] >$ уравнения $[bu\tilde{b}] = y$ относительно $u$ следует, что любой элемент $y \in A$ может быть представлен в виде $y = [bx\tilde{b}]$, $x \in A$. Таким образом, автоморфизм $\alpha$ переставляет множества $x^{-1} \text{\textcircled{b}} B \text{\textcircled{b}} x$, и поэтому $N^\alpha = N$. Из последнего равенства следует, что отображение $\hat{\alpha}$, определенное на $< A/N, \text{\textcircled{b}} >$ по правилу

$$\hat{\alpha} : x \text{\textcircled{b}} N \to x^\alpha \text{\textcircled{b}} N \qquad (9)$$

является автоморфизмом группы $< A/N, \text{\textcircled{b}} >$.

Группа $< A/N, \text{\textcircled{b}} >$ как факторгруппа по ядру гомоморфизма $\gamma$ изоморфна образу $\Delta$ этого гомоморфизма. Указанный изоморфизм определяется следующим образом (см., например, [109]):

$$\tau : x \text{\textcircled{b}} N \to x^\gamma = \delta_x. \qquad (10)$$

Ясно, что тогда

$$\beta = \tau^{-1}\hat{\alpha}\tau \qquad (11)$$

– автоморфизм группы $\Delta$.

Положив

$$c^{\gamma} = d, \qquad (12)$$

получим

$$d^{\beta} \stackrel{(11)}{=} d^{\tau^{-1}\hat{\alpha}\tau} = (d^{\tau^{-1}})^{\hat{\alpha}\tau} \stackrel{(12)}{=} ((c^{\gamma})^{\tau^{-1}})^{\hat{\alpha}\tau} \stackrel{(10)}{=}$$

$$\stackrel{(10)}{=} ((c \ⓑ N)^{\hat{\alpha}})^{\tau} \stackrel{(9)}{=} (c^{\alpha} \ⓑ N)^{\tau} \stackrel{(2)}{=} (c \ⓑ N)^{\tau} \stackrel{(10)}{=} c^{\gamma} \stackrel{(12)}{=} d,$$

откуда

$$d^{\beta} = d. \qquad (13)$$

Если теперь $t \in \Delta$,

$$x^{\gamma} = t \qquad (14)$$

для некоторого $x \in A$, то

$$t^{\beta^{n-1}} d \stackrel{(11)}{=} t^{\underbrace{\tau^{-1}\hat{\alpha}\tau \ldots \tau^{-1}\hat{\alpha}\tau}_{n-1}} d = t^{\tau^{-1}\hat{\alpha}^{n-1}\tau} d \stackrel{(12)}{=}$$

$$\stackrel{(12)}{=} t^{\tau^{-1}\hat{\alpha}^{n-1}\tau} c^{\gamma} \stackrel{(10)}{=} (t^{\tau^{-1}\hat{\alpha}^{n-1}})^{\tau}(c \ⓑ N)^{\tau} =$$

$$= (t^{\tau^{-1}\hat{\alpha}^{n-1}} \ⓑ (c \ⓑ N))^{\tau} \stackrel{(14)}{=} ((x^{\gamma})^{\tau^{-1}})^{\hat{\alpha}^{n-1}} \ⓑ c \ⓑ N)^{\tau} \stackrel{(10)}{=}$$

$$\stackrel{(10)}{=} ((x \ⓑ N)^{\hat{\alpha}^{n-1}} \ⓑ c \ⓑ N)^{\tau} \stackrel{(9)}{=} (x^{\alpha^{n-1}} \ⓑ N \ⓑ c \ⓑ N)^{\tau} =$$

$$= (x^{\alpha^{n-1}} \ⓑ c \ⓑ N)^{\tau} \stackrel{(3)}{=} (c \ⓑ x \ⓑ N)^{\tau} = ((c \ⓑ N) \ⓑ (x \ⓑ N))^{\tau} =$$

$$= (c \ⓑ N)^{\tau}(x \ⓑ N)^{\tau} \stackrel{(10)}{=} c^{\gamma} x^{\gamma} \stackrel{(12),(14)}{=} dt.$$

Следовательно,

$$t^{\beta^{n-1}} d = dt. \qquad (15)$$

Так как для автоморфизма $\beta$ группы $\Delta$ и элемента $d \in \Delta$ выполняются условия (13) и (15), то по обратной теореме Глу-скина-Хоссу $<\Delta, \langle\ \rangle>$ – n-арная группа с n-арной операцией $\langle\ \rangle$, определяемой по правилу

$$\langle t_1 t_2 \ldots t_n \rangle = t_1 t_2^{\beta} \ldots t_n^{\beta^{n-1}} d, \ t_1, t_2, \ldots, t_n \in \Delta. \qquad (16)$$

Покажем теперь, что $\gamma$ является не только гомоморфизмом группы $<A, \textcircled{b}>$ на группу $\Delta$, но также и гомоморфизмом n-арной группы $<A, [\ ]>$ на n-арную группу $<\Delta, \langle\ \rangle>$. Действительно,

$$[x_1 x_2 \ldots x_n]^{\gamma} \stackrel{(1)}{=} (x_1 \textcircled{b} x_2^{\alpha} \textcircled{b} \ldots \textcircled{b} x_n^{\alpha^{n-1}} \textcircled{b} c)^{\gamma} =$$

$$= x_1^{\gamma} x_2^{\alpha\gamma} \ldots x_n^{\alpha^{n-1}\gamma} c^{\gamma} \stackrel{(10)}{=}$$

$$\stackrel{(10)}{=} (x_1 \textcircled{b} N)^{\tau} (x_2^{\alpha} \textcircled{b} N)^{\tau} \ldots (x_n^{\alpha^{n-1}} \textcircled{b} N)^{\tau} (c \textcircled{b} N)^{\tau} \stackrel{(9)}{=}$$

$$\stackrel{(9)}{=} (x_1 \textcircled{b} N)^{\tau} (x_2 \textcircled{b} N)^{\hat{\alpha}\tau} \ldots (x_n \textcircled{b} N)^{\hat{\alpha}^{n-1}\tau} (c \textcircled{b} N)^{\tau} =$$

$$= (x_1 \textcircled{b} N)^{\tau} (x_2 \textcircled{b} N)^{\tau\tau^{-1}\hat{\alpha}\tau} \ldots (x_n \textcircled{b} N)^{\tau\tau^{-1}\hat{\alpha}^{n-1}\tau} (c \textcircled{b} N)^{\tau} \stackrel{(10)}{=}$$

$$\stackrel{(10)}{=} x_1^{\gamma} (x_2^{\gamma})^{\tau^{-1}\hat{\alpha}\tau} \ldots (x_n^{\gamma})^{\tau^{-1}\hat{\alpha}^{n-1}\tau} c^{\gamma} \stackrel{(11)}{=}$$

$$\stackrel{(11)}{=} x_1^{\gamma} (x_2^{\gamma})^{\beta} \ldots (x_n^{\gamma})^{\beta^{n-1}} d \stackrel{(16)}{=} \langle x_1^{\gamma} x_2^{\gamma} \ldots x_n^{\gamma} \rangle,$$

то есть

$$[x_1 x_2 \ldots x_n]^{\gamma} = \langle x_1^{\gamma} x_2^{\gamma} \ldots x_n^{\gamma} \rangle. \qquad \blacksquare$$

В следующих предложениях укажем явный вид автоморфизма $\beta$.

**6.1.2. Предложение.** Если $\alpha$ и $\beta$ – автоморфизмы из доказательства теоремы 6.1.1, то образ элемента $\delta_a$ при отображе-нии β совпадает с $\delta_{a^\alpha}$. В частности, автоморфизм β переводит $\delta_c$ в $\delta_c$.

*Доказательство.* Так как

$$\delta_a^\beta = \delta_a^{\tau^{-1}\hat\alpha\hat\tau} = (a \circledb N)^{\hat\alpha\tau} = (a^\alpha \circledb N)^\tau = \delta_{a^\alpha},$$

то $\delta_a^\beta = \delta_{a^\alpha}$. Из последнего равенства и условия $c^\alpha = c$ следует $\delta_c^\beta = \delta_c$. ∎

**Предложение 6.1.3.** Если в условии теоремы 6.1.1 b – идемпотентный элемент, то

$$\beta : t \to \langle \varepsilon t \underbrace{\varepsilon \ldots \varepsilon}_{n-2} \rangle$$

для любого $t \in \Delta$, где ε – тождественная подстановка, а отображение β определено так же, как в доказательстве теоремы 6.1.1.

*Доказательство.* Так как b – идемпотент, и по определению

$$c = [\underbrace{b \ldots b}_{n}],$$

то $c = b$, а так как, согласно (12), $d = c^\gamma$, то $d = b^\gamma$. Из последнего равенства, учитывая, что при гомоморфизме γ единица b группы $\langle A, \circledb \rangle$ переходит в единицу группы Δ, получаем, что $d = \varepsilon$ – тождественная подстановка.

Так как любая степень автоморфизма β переводит единицу ε группы Δ в себя, то, ввиду (16),

$$t^\beta = \varepsilon t^\beta \underbrace{\varepsilon \ldots \varepsilon}_{n-2} \varepsilon = \varepsilon t^\beta \underbrace{\varepsilon^{\beta^2} \ldots \varepsilon^{\beta^{n-1}}}_{n-2} d = \langle \varepsilon t \underbrace{\varepsilon \ldots \varepsilon}_{n-2} \rangle,$$

то есть $\beta : t \to \langle \varepsilon t \underbrace{\varepsilon \ldots \varepsilon}_{n-2} \rangle$. ■

**6.1.4. Теорема**. Пусть $<A, [\,]>$, $<B, [\,]>$ и $b$ те же, что и в формулировке теоремы 6.1.1, и пусть

$$c = [\underbrace{b \ldots b}_{n}] \in N,$$

где $N$ – ядро гомоморфизма $\gamma$ группы $<A, \circledb>$ на группу $\Delta$. Тогда:

1) $<N, [\,]>$ – полуинвариантная n-арная подгруппа n-арной группы $<A, [\,]>$, максимальная среди содержащихся в $<B, [\,]>$;

2) если $\gamma$ – изоморфизм, то $N = \{b\}$, $b$ – идемпотент, и $<B, [\,]>$ не содержит полуинвариантных в $<A, [\,]>$ n-арных подгрупп, отличных от $<\{b\}, [\,]>$.

*Доказательство.* 1) При доказательстве теоремы 6.1.1 установлено, что $N^\alpha = N$. Поэтому сужение автоморфизма $\alpha$ группы $<A, \circledb>$ на $N$ является автоморфизмом группы $<N, \circledb>$. Кроме того, так как по условию $c \in N$, то, согласно теореме Глускина-Хоссу, $<N, [\,]>$ – n-арная подгруппа n-арной группы $<A, [\,]>$.

Так как $<N, \circledb>$ инвариантна в $<A, \circledb>$, и $b \in N$, то по следствию 2.3.13 $<N, [\,]>$ полуинвариантна в $<A, [\,]>$.

Включение $N \subseteq B$ отмечалось при доказательстве 1) теоремы 6.1.1.

Пусть теперь $<K, [\,]>$ – полуинвариантная n-арная подгруппа n-арной группы $<A, [\,]>$ такая, что $N \subseteq K \subseteq B$. Ясно, что $b \in K$ и $<K, \circledb>$ – подгруппа в $<A, \circledb>$. Согласно следствию 2.3.13, $<K, \circledb>$ инвариантна в $<A, \circledb>$. Так как $N \subseteq K \subseteq B$, то по соответствующей теореме для бинарных групп $N = K$ (см., например, [109]).

2) Так как $\gamma$ – взаимно однозначное отображение, то $|N| = 1$, а так как $b \in N$, то $N = \{b\}$. Тогда из условия

$$[\underbrace{b \ldots b}_{n}] \in N$$

следует

$$[\underbrace{b \ldots b}_{n}] = b.$$

Следовательно, b – идемпотент.

Согласно 1), $< \{b\}, [\,] >$ – полуинвариантная n-арная подгруппа n-арной группы $< A, [\,] >$, максимальная среди содержащихся в $< B, [\,] >$. Поэтому $< B, [\,] >$ не содержит полуинвариантных в $< A, [\,] >$ n-арных подгрупп, отличных от $< \{b\}, [\,] >$. ■

**6.1.5. Следствие.** Пусть $< B, [\,] >$ – n-арная подгруппа n-арной группы $< A, [\,] >$, $b \in B$, b – идемпотент, и пусть $< B, [\,] >$ не содержит полуинвариантных в $< A, [\,] >$ n-арных подгрупп, отличных от $< \{b\}, [\,] >$. Тогда $\gamma$ – изоморфизм.

*Доказательство.* Предположим, что $\gamma$ не является взаимно однозначным отображением. Тогда $\gamma$, как гомоморфизм группы $< A, ⓑ >$ на группу $\Delta$ имеет ядро N такое, что $\{b\} \subset N \subseteq B$. А так как b – идемпотент, то

$$[\underbrace{b \ldots b}_{n}] = b \in N.$$

Тогда из 1) теоремы 6.1.4 вытекает, что $< N, [\,] >$ – полуинвариантная n-арная подгруппа в $< A, [\,] >$, отличная от $< \{b\}, [\,] >$, что противоречит условию. ■

**6.1.6. Предложение.** Если индекс n-арной подгруппы $< B, [\,] >$ в n-арной группе $< A, [\,] >$ равен m, то m делит порядок $|\Delta|$, который делит m!.

*Доказательство.* Так как индекс $< B, [\,] >$ в $< A, [\,] >$ равен m, то мощность множества $\Omega$ равна m. Тогда число $|\Delta|$,

равное порядку подгруппы $\Delta$ симметрической группы $S_\Omega$ делит ее порядок, равный $m!$.

При доказательстве теоремы 6.1.1 установлено, что

$$\ker \gamma = N \subseteq B,$$

а группы $< A/N, \circledR >$ и $\Delta$ изоморфны. Последнее означает равенство порядков $|A/N|$ и $|\Delta|$, откуда, учитывая конечность порядка $|\Delta|$, получаем конечность индекса $< N, \circledR >$ в $< A, \circledR >$. Тогда, согласно соответствующему бинарному результату,

$$|A : N| = |A : B| \cdot |B : N|,$$

откуда и из $|A : N| = |\Delta|$, $|A : B| = m$ следует

$$|\Delta| = m \cdot |B : N|,$$

то есть $m$ делит порядок $|\Delta|$. ∎

Если $b$ – идемпотентный элемент n-арной группы $< A [\,] >$, $b \in B$, то $[\underbrace{b \ldots b}_{n}] = b \in N$. Поэтому из 1) теоремы 6.1.4 и предложения 6.1.6 вытекает

**6.1.7. Следствие.** Пусть $< A, [\,] >$ – n-арная группа, $< B, [\,] >$ – ее n-арная подгруппа, $b \in B$, $b$ – идемпотент. Тогда:

1) $< N, [\,] >$ – полуинвариантная n-арная подгруппа n-арной группы $< A, [\,] >$, максимальная среди содержащихся в $< B, [\,] >$;

2) если $< B, [\,] >$ имеет в $< A, [\,] >$ конечный индекс $m$, то $m$ делит порядок n-арной факторгруппы $< A/N, [\,] >$, который делит $m!$.

## §6.2. СВЯЗЬ МЕЖДУ РАЗЛОЖЕНИЯМИ B < A, [ ] > И РАЗЛОЖЕНИЯМИ В $A^{(K)}$

В данном параграфе устанавливается связь между разложением n-арной группы < A, [ ] > по ее n-арной подгруппе < B, [ ] > и разложениями множеств

$$A^{(k)} = \{\theta_A(a_1 \ldots a_k) \,|\, a_1, \ldots, a_k \in A\},\ k = 1, \ldots, n-1,$$

которые были определены в §1.3, на непересекающиеся подмножества. Как следствия, получены соответствия между разложением < A, [ ] > по < B, [ ] > и разложениями $A_o$ по $B_o(A)$ и $A^*$ по $B^*(A)$.

**6.2.1. Теорема** [110]**.** Пусть < B, [ ] > – n-арная подгруппа n-арной группы < A, [ ] >, $k \in \{1, \ldots, n-1\}$, $b_1, \ldots, b_{k-1}$ – фиксированные элементы из B. Тогда:

1) если

$$A = \bigcup_{i \in I}\ [x_i \underbrace{B \ldots B}_{n-1}] \qquad (1.1)$$

– разложение < A, [ ] > на непересекающиеся левые смежные классы по < B, [ ] >, то

$$A^{(k)} = \bigcup_{i \in I}\ \theta_A(x_i b_1 \ldots b_{k-1})B_o(A) \qquad (1.2)$$

– разложение $A^{(k)}$ на непересекающиеся подмножества, а отображение

$$[x_i \underbrace{B \ldots B}_{n-1}] \to \theta_A(x_i b_1 \ldots b_{k-1})B_o(A) \qquad (1.3)$$

является биекцией множества всех левых смежных классов < A, [ ] > по < B, [ ] > на множество

$$\{\theta_A(x_i b_1 \ldots b_{k-1})B_o(A) \,|\, i \in I\}; \qquad (1.4)$$

2) если (1.2) – разложение $A^{(k)}$ на непересекающиеся подмножества, то (1.1) – разложение $<A, [\,]>$ на непересекающиеся левые смежные классы по $<B, [\,]>$, а отображение

$$\theta_A(x_i b_1 \ldots b_{k-1}) B_o(A) \to [x_i \underbrace{B \ldots B}_{n-1}] \qquad (1.5)$$

является биекцией множества (1.4) на множество всех левых смежных классов $<A, [\,]>$ по $<B, [\,]>$.

***Доказательство.*** 1) Пусть $\theta_A(a_1 \ldots a_k)$ – произвольный элемент из $A^{(k)}$. Для фиксированных $b_1, \ldots, b_{k-1} \in B$ найдется $y \in A$ такой, что

$$\theta_A(a_1 \ldots a_k) = \theta_A(y b_1 \ldots b_{k-1}). \qquad (1.6)$$

Если $b_k, \ldots, b_{n-1} \in B$, то по условию

$$[y b_1 \ldots b_{k-1} b_k \ldots b_{n-1}] \in [x_i \underbrace{B \ldots B}_{n-1}]$$

для некоторого $i \in I$, откуда

$$[y b_1 \ldots b_{k-1} b_k \ldots b_{n-1}] = [x_i b_1 \ldots b_{k-1} b_k \ldots b_{n-2} b]$$

для некоторого $b \in B$. Тогда

$$\theta_A(y b_1 \ldots b_{k-1} b_k \ldots b_{n-1}) = \theta_A(x_i b_1 \ldots b_{n-2} b),$$

$$\theta_A(y b_1 \ldots b_{k-1}) \theta_A(b_k \ldots b_{n-1}) = \theta_A(x_i b_1 \ldots b_{k-1}) \theta_A(b_k \ldots b_{n-2} b),$$

$$\theta_A(y b_1 \ldots b_{k-1}) = \theta_A(x_i b_1 \ldots b_{k-1}) \theta_A(b_k \ldots b_{n-2} b) \theta_A^{-1}(b_k \ldots b_{n-1}),$$

$$\theta_A(y b_1 \ldots b_{k-1}) \in \theta_A(x_i b_1 \ldots b_{k-1}) B_o(A),$$

откуда и из (1.6) следует $\theta_A(a_1 \ldots a_k) \in \theta_A(x_i b_1 \ldots b_{k-1}) B_o(A)$. Следовательно,

$$A^{(k)} \subseteq \bigcup_{i \in I} \theta_A(x_i b_1 \ldots b_{k-1}) B_o(A).$$

Обратное включение

$$\bigcup_{i \in I} \theta_A(x_i b_1 \ldots b_{k-1}) B_o(A) \subseteq A^{(k)}$$

очевидно. Таким образом, доказано равенство (1.2).

Предположим, что

$$\theta_A(x_i b_1 \ldots b_{k-1}) B_o(A) \cap \theta_A(x_j b_1 \ldots b_{k-1}) B_o(A) \neq \varnothing, \ i \neq j,$$

то есть

$$\theta_A(x_i b_1 \ldots b_{k-1}) \theta_A(c_1 \ldots c_{n-1}) = \theta_A(x_j b_1 \ldots b_{k-1}) \theta_A(d_1 \ldots d_{n-1})$$

для $c_1, \ldots, c_{n-1}, d_1, \ldots, d_{n-1} \in B$. Тогда

$$\theta_A(x_i b_1 \ldots b_{k-1}) \theta_A(c_1 \ldots c_{n-1}) \theta_A(b_k \ldots b_{n-1}) =$$
$$= \theta_A(x_j b_1 \ldots b_{k-1}) \theta_A(d_1 \ldots d_{n-1}) \theta_A(b_k \ldots b_{n-1})$$

для любых $b_k, \ldots, b_{n-1} \in B$. Так как

$$b_1 \ldots b_{k-1} c_1 \ldots c_{n-1} b_k \ldots b_{n-1} \theta_A c'_1 \ldots c'_{n-1},$$
$$b_1 \ldots b_{k-1} d_1 \ldots d_{n-1} b_k \ldots b_{n-1} \theta_A d'_1 \ldots d'_{n-1}$$

для некоторых $c'_1, \ldots, c'_{n-1}, d'_1, \ldots, d'_{n-1} \in B$, то из последнего равенства следует

$$\theta_A(x_i) \theta_A(c'_1 \ldots c'_{n-1}) = \theta_A(x_j) \theta_A(d'_1 \ldots d'_{n-1}),$$

откуда

$$\theta_A(x_i c'_1 \ldots c'_{n-1}) = \theta_A(x_j d'_1 \ldots d'_{n-1}),$$
$$[x_i c'_1 \ldots c'_{n-1}] = [x_j d'_1 \ldots d'_{n-1}].$$

Последнее равенство противоречит тому, что (1.1) – разложение < A, [ ] > на непересекающиеся левые смежные классы по < B, [ ] >. Следовательно, (1.2) является разложением $A^{(k)}$ на непересекающиеся подмножества.

Из доказанного следует, что (1.3) – биекция.

2) Пусть a – произвольный элемент из A. Тогда, если зафиксировать $b_k, \ldots, b_{n-1} \in B$, то найдется $y \in A$ такой, что

$$a = [yb_1 \ldots b_{k-1}b_k \ldots b_{n-1}] \qquad (1.7)$$

В силу условия,

$$\theta_A(yb_1 \ldots b_{k-1}) \in \theta_A(x_ib_1 \ldots b_{k-1})B_o(A)$$

для некоторого $i \in I$, откуда

$$\theta_A(yb_1 \ldots b_{k-1}) = \theta_A(x_ib_1 \ldots b_{k-1})\theta_A(b_1 \ldots b_{n-2}b)$$

для некоторого $b \in B$. Тогда

$$\theta_A(yb_1 \ldots b_{k-1})\theta_A(b_k \ldots b_{n-1}) =$$
$$= \theta_A(x_ib_1 \ldots b_{k-1})\theta_A(b_1 \ldots b_{n-2}b)\theta_A(b_k \ldots b_{n-1}),$$
$$[yb_1 \ldots b_{n-1}] = [x_ib_1 \ldots b_{k-1}b_1 \ldots b_{n-2}bb_k \ldots b_{n-1}],$$

откуда и из (1.7) следует $a \in [x_i \underbrace{B \ldots B}_{n-1}]$. Следовательно,

$$A \subseteq \bigcup_{i \in I} [x_i \underbrace{B \ldots B}_{n-1}].$$

Обратное включение

$$\bigcup_{i \in I} [x_i \underbrace{B \ldots B}_{n-1}] \subseteq A$$

очевидно. Таким образом, доказано равенство (1.1).

Предположим, что

$$[x_i \underbrace{B \ldots B}_{n-1}] \cap [x_j \underbrace{B \ldots B}_{n-1}] \neq \varnothing, i \neq j,$$

то есть

$$[x_ic_1 \ldots c_{n-1}] = [x_jd_1 \ldots d_{n-1}]$$

для $c_1, \ldots, c_{n-1}, d_1, \ldots, d_{n-1} \in B$. Так как

$$c_1 \ldots c_{n-1}\theta_A b_1 \ldots b_{k-1}c'_k \ldots c'_{n-1},$$

$$d_1 \ldots d_{n-1}\theta_A b_1 \ldots b_{k-1}d'_k \ldots d'_{n-1}$$

для некоторых $c'_k, \ldots, c'_{n-1}, d'_k, \ldots, d'_{n-1} \in B$, то из последнего равенства следует

$$\theta_A(x_ib_1 \ldots b_{k-1})\theta_A(c'_k \ldots c'_{n-1}) = \theta_A(x_jb_1 \ldots b_{k-1})\theta_A(d'_k \ldots d'_{n-1}),$$

$$\theta_A(x_ib_1 \ldots b_{k-1})\theta_A(c'_k \ldots c'_{n-1})\theta_A(b_1 \ldots b_{k-1}) =$$

$$= \theta_A(x_jb_1 \ldots b_{k-1})\theta_A(d'_k \ldots d'_{n-1})\theta_A(b_1 \ldots b_{k-1}),$$

$$\theta_A(x_ib_1 \ldots b_{k-1})B_o(A) \cap \theta_A(x_jb_1 \ldots b_{k-1})B_o(A) \neq \varnothing,$$

что противоречит условию.

Из доказанного следует, что (1.5) – биекция. ∎

Аналогично теореме 6.2.1 доказывается двойственная к ней

**6.2.2. Теорема** [110]**.** Пусть $<B, [\,]>$ – n-арная подгруппа n-арной группы $<A, [\,]>$, $k \in \{1, \ldots, n-1\}$, $b_1, \ldots, b_{k-1}$ – фиксированые элементы из B. Тогда:

1) если

$$A = \bigcup_{i \in I} [\underbrace{B \ldots B}_{n-1} x_i] \tag{2.1}$$

– разложение $<A, [\,]>$ на непересекающиеся правые смежные классы по $<B, [\,]>$, то

$$A^{(k)} = \bigcup_{i \in I} B_o(A)\theta_A(b_1 \ldots b_{k-1}x_i) \tag{2.2}$$

– разложение $A^{(k)}$ на непересекающиеся подмножества, а отображение

$$[\underbrace{B \ldots B}_{n-1} x_i] \to B_o(A)\theta_A(b_1 \ldots b_{k-1}x_i) \qquad (2.3)$$

является биекцией множества всех правых смежных классов $<A, [\ ]>$ по $<B, [\ ]>$ на множество

$$\{B_o(A)\theta_A(b_1 \ldots b_{k-1}x_i) \mid i \in I\}; \qquad (2.4)$$

2) если (2.2) – разложение $A^{(k)}$ на непересекающиеся подмножества, то (2.1) – разложение $<A, [\ ]>$ на непересекающиеся правые смежные классы по $<B, [\ ]>$, а отображение

$$B_o(A)\theta_A(b_1 \ldots b_{k-1}x_i) \to [\underbrace{B \ldots B}_{n-1} x_i] \qquad (2.5)$$

является биекцией множества (2.4) на множество всех правых смежных классов $<A, [\ ]>$ по $<B, [\ ]>$.

Ясно, что отображения (1.3) и (1.5) являются взаимно обратными. То же самое можно сказать об отображениях (2.3) и (2.5).

Полагая в теореме 6.2.1 $k = n - 1$, получим

**6.2.3. Следствие** [110]. Пусть $<B, [\ ]>$ – n-арная подгруппа n-арной группы $<A, [\ ]>$, $b_1, \ldots, b_{n-2}$ – фиксированные элементы из B. Тогда:

1) если (1.1) – разложение $<A, [\ ]>$ на непересекающиеся левые смежные классы по $<B, [\ ]>$, то

$$A_o = \bigcup_{i \in I} \theta_A(x_i b_1 \ldots b_{n-2}) B_o(A)$$

– разложение $A_o$ на непересекающиеся левые смежные классы по $B_o(A)$, а отображение

$$[x_i \underbrace{B \ldots B}_{n-1}] \to \theta_A(x_i b_1 \ldots b_{n-2}) B_o(A)$$

является биекцией множества всех левых смежных классов $<A, [\ ]>$ по $<B, [\ ]>$ на множество всех левых смежных классов $A_o$ по $B_o(A)$;

2) если равенство из 1) является разложением $A_o$ на непересекающиеся левые смежные классы по $B_o(A)$, то (1.1) – разложение $<A, [\ ]>$ на непересекающиеся левые смежные классы по $<B, [\ ]>$, а отображение

$$\theta_A(x_i b_1 \ldots b_{n-2}) B_o(A) \to [x_i \underbrace{B \ldots B}_{n-1}]$$

является биекцией множества всех левых смежных классов $A_o$ по $B_o(A)$ на множество всех левых смежных классов $<A, [\ ]>$ по $<B, [\ ]>$.

Полагая в теореме 6.2.2 $k = n - 1$, получим

**6.2.4. Следствие** [110]. Пусть $<B, [\ ]>$ – n-арная подгруппа n-арной группы $<A, [\ ]>$, $b_1, \ldots, b_{n-2}$ – фиксированные элементы из B. Тогда:

1) если (2.1) – разложение $<A, [\ ]>$ на непересекающиеся правые смежные классы по $<B, [\ ]>$, то

$$A_o = \bigcup_{i \in I} B_o(A) \theta_A(b_1 \ldots b_{n-2} x_i)$$

– разложение $A_o$ на непересекающиеся правые смежные классы по $B_o(A)$, а отображение

$$[\underbrace{B \ldots B}_{n-1} x_i] \to B_o(A) \theta_A(b_1 \ldots b_{n-2} x_i)$$

является биекцией множества всех правых смежных классов $<A, [\ ]>$ по $<B, [\ ]>$ на множество всех правых смежных классов $A_o$ по $B_o(A)$;

2) если равенство из 1) является разложением $A_o$ на непересекающиеся правые смежные классы по $B_o(A)$, то (2.1) – разложение $<A, [\ ]>$ на непересекающиеся правые смежные классы по $<B, [\ ]>$, а отображение

$$B_o(A)\theta_A(b_1 \ldots b_{n-2}x_i) \to [\underbrace{B \ldots B}_{n-1} x_i]$$

является биекцией множества всех правых смежных классов $A_o$ по $B_o(A)$ на множество всех правых смежных классов $<A, [\ ]>$ по $<B, [\ ]>$.

Из теорем 6.2.1 и 6.2.2 вытекает

**6.2.5. Следствие.** Индекс n-арной подгруппы $<B, [\ ]>$ в n-арной группе $<A, [\ ]>$ совпадает с мощностями множеств (1.4) и (2.4).

В частности, из следствия 6.2.3 (также из следствия 6.2.4) вытекает

**6.2.6. Следствие.** Индекс n-арной подгруппы $<B, [\ ]>$ в n-арной группе $<A, [\ ]>$ совпадает с индексом подгруппы $B_o(A)$ в группе $A_o$: $|A : B| = |A_o : B_o(A)|$.

Полагая в теореме 6.1.1 $k = 1$, получим

**6.2.7. Следствие** [110]**.** Пусть $<B, [\ ]>$ – n-арная подгруппа n-арной группы $<A, [\ ]>$. Тогда:
1) если (1.1) – разложение $<A, [\ ]>$ на непересекающиеся левые смежные классы по $<B, [\ ]>$, то

$$A^{(1)} = \bigcup_{i \in I} \theta_A(x_i)B_o(A)$$

– разложение $A^{(1)}$ на непересекающиеся подмножества, а отображение

$$[x_i \underbrace{B \ldots B}_{n-1}] \to \theta_A(x_i)B_o(A)$$

является биекцией множества всех левых смежных классов $<A, [\ ]>$ по $<B, [\ ]>$ на множество $\{\theta_A(x_i)B_o(A) \mid i \in I\}$;

2) если равенство из 1) является разложением $A^{(1)}$ на непересекающиеся подмножества, то (1.1) – разложение $<A, [\ ]>$ на непересекающиеся левые смежные классы по $<B, [\ ]>$, а отображение

$$\theta_A(x_i)B_o(A) \to [x_i \underbrace{B \ldots B}_{n-1}]$$

является биекцией множества из 1) на множество всех левых смежных классов $<A, [\ ]>$ по $<B, [\ ]>$.

Полагая в теореме 6.2.2 $k = 1$, получим

**6.2.8. Следствие** [110]. Пусть $<B, [\ ]>$ – n-арная подгруппа n-арной группы $<A, [\ ]>$. Тогда:

1) если (2.1) – разложение $<A, [\ ]>$ на непересекающиеся правые смежные классы по $<B, [\ ]>$, то

$$A^{(1)} = \bigcup_{i \in I} B_o(A)\theta_A(x_i)$$

– разложение $A^{(1)}$ на непересекающиеся подмножества, а отображение

$$[\underbrace{B \ldots B}_{n-1} x_i] \to B_o(A)\theta_A(x_i)$$

является биекцией множества всех правых смежных классов $<A, [\ ]>$ по $<B, [\ ]>$ на множество $\{B_o(A)\theta_A(x_i) \mid i \in I\}$;

2) если равенство из 1) является разложением $A^{(1)}$ на непересекающиеся подмножества, то (2.1) является разложением $<A, [\ ]>$ на непересекающиеся правые смежные классы по $<B, [\ ]>$, а отображение

$$B_o(A)\theta_A(x_i) \to [\underbrace{B \ldots B}_{n-1} x_i]$$

является биекцией множества из 1) на множество всех правых смежных классов $< A, [\ ] >$ по $< B, [\ ] >$.

**6.2.9. Замечание.** Замена в утверждении 2) теоремы 6.1.1 разложения (1.2) на разложение

$$A^{(k)} = \bigcup_{i \in I} \theta_A(x_{i1} \ldots x_{ik})B_o(A),$$

$$\theta_A(x_{i1} \ldots x_{ik})B_o(A) \cap \theta_A(x_{j1} \ldots x_{jk})B_o(A) = \varnothing, i \neq j$$

не приводит к более общей ситуации, так как в этом случае для фиксированных $b_1, \ldots b_{k-1} \in B$ найдутся $x_i \in A$ ($i \in I$) такие, что $\theta_A(x_{i1} \ldots x_{ik}) = \theta_A(x_i b_1 \ldots b_{k-1})$. Поэтому записанное выше разложение совпадает с разложением (1.2).

То же самое можно сказать о разложении (2.2) из утверждения 2) теоремы 6.2.2.

**6.2.10. Замечание.** В приведенных выше разложениях только при $k = n - 1$ речь идет о разложении группы по подгруппе, так как $B_o(A)$ – подгруппа группы $A_o = A^{(n-1)}$. Во всех остальных случаях нельзя говорить даже о разложении множества по подмножеству, так как $B_o(A)$ не является подмножеством в $A^{(k)}$ для любого $k = 1, \ldots, n - 2$.

Так как $A^* = \bigcup_{k=1}^{n-1} A^{(k)}$, то теоремы 6.2.1 и 6.2.2 позволяют по разложению $< A, [\ ] >$ по $< B, [\ ] >$ получать разложения $A^*$ по $B_o(A)$.

**6.2.11. Теорема** [110]. Пусть $< B, [\ ] >$ – n-арная подгруппа n-арной группы $< A, [\ ] >$, $b_1, \ldots, b_{n-2}$ – фиксированные элементы из В. Тогда:

1) если (1.1) – разложение $< A, [\ ] >$ на непересекающиеся левые смежные классы по $< B, [\ ] >$, то

$$A^* = \bigcup_{k=1}^{n-1} \left( \bigcup_{i \in I} \theta_A(x_i b_1 \ldots b_{k-1})B_o(A) \right)$$

– разложение A* на непересекающиеся левые смежные классы по $B_o(A)$;

2) если равенство из 1) является разложением A* на непересекающиеся левые смежные классы по $B_o(A)$, то (1.1) – разложение $<A, [\ ]>$ на непересекающиеся левые смежные классы по $<B, [\ ]>$.

***Доказательство.*** 1) Следует из 1) теоремы 1 и равенств $A^{(k)} \cap A^{(m)} = \varnothing$, где $m \neq k$.

2) Если предположить, что $\bigcup_{i \in I} \theta_A(x_i) B_o(A) \neq A^{(1)}$, то есть

$$\bigcup_{i \in I} \theta_A(x_i) B_o(A) \subset A^{(1)},$$

то из

$$A^* = \bigcup_{k=1}^{n-1} A^{(k)}, \ A^{(k)} \cap A^{(m)} = \varnothing, \ k \neq m$$

вытекает

$$\bigcup_{k=1}^{n-1} (\bigcup_{i \in I} \theta_A(x_i b_1 \ldots b_{k-1}) B_o(A)) \subset A^*,$$

что противоречит равенству из 1). Следовательно,

$$\bigcup_{i \in I} \theta_A(x_i) B_o(A) = A^{(1)}.$$

Точно так же

$$\bigcup_{i \in I} \theta_A(x_i b_1 \ldots b_{k-1}) B_o(A) = A^{(k)}.$$

Применяя к любому из полученных равенств утверждения 2) теоремы 6.2.1, видим, что (1.1) – разложение $<A, [\ ]>$ на непересекающиеся левые смежные классы по $<B, [\ ]>$. ∎

Двойственной к теореме 6.2.11 является

**6.2.12. Теорема** [110]**.** Пусть $< B, [\ ] >$ – n-арная подгруппа n-арной группы $< A, [\ ] >$, $b_1, \ldots, b_{n-2}$ – фиксированные элементы из B. Тогда:

1) если (2.1) – разложение $< A, [\ ] >$ на непересекающиеся правые смежные классы по $< B, [\ ] >$, то

$$A^* = \bigcup_{k=1}^{n-1} (\bigcup_{i \in I} B_0(A) \theta_A(b_1 \ldots b_{k-1} x_i))$$

– разложение $A^*$ на непересекающиеся правые смежные классы по $B_0(A)$;

2) если равенство из 1) является разложением $A^*$ на непересекающиеся правые смежные классы по $B_0(A)$, то (2.1) – разложение $< A, [\ ] >$ на непересекающиеся правые смежные классы по $< B, [\ ] >$.

**6.2.13. Предложение** [110]**.** Пусть $< B, [\ ] >$ – n-арная подгруппа n-арной группы $< A, [\ ] >$, $x \in A$, $i \in \{1, 2, \ldots, n-1\}$. Тогда:

$$\theta_A(xb_1 \ldots b_i)B_o(A) = \theta_A(xc_1 \ldots c_i)B_o(A),$$

$$B_o(A)\theta_A(b_1 \ldots b_i x) = B_o(A)\theta_A(c_1 \ldots c_i x)$$

для любых $b_1, \ldots, b_i, c_1, \ldots, c_i \in B$.

*Доказательство.* В B существуют элементы $b_{i+1}, \ldots, b_n$ такие, что

$$b_1 = [c_1 \ldots c_i b_{i+1} \ldots b_n].$$

Тогда

$$\theta_A(xb_1 \ldots b_i)B_o(A) = \theta_A(x[c_1 \ldots c_i b_{i+1} \ldots b_n]b_2 \ldots b_i)B_o(A) =$$

$$= \theta_A(xc_1 \ldots c_i)\theta_A(b_{i+1} \ldots b_n b_2 \ldots b_i)B_o(A) = \theta_A(xc_1 \ldots c_i)B_o(A),$$

то есть верно первое равенство.

Второе равенство доказывается аналогично. ∎

**6.2.14. Замечание.** Предложение 6.2.13 показывает, что в формулировке теоремы 6.2.1 элементы $b_1, \ldots, b_{k-1} \in B$ можно не фиксировать, а разложение (1.2) можно записать в виде

$$A^{(k)} = \bigcup_{i \in I} \theta_A(x_i b_{i1} \ldots b_{i(k-1)}) B_o(A),$$

где $b_{i1}, \ldots, b_{i(k-1)}$ – произвольные из $B$.

Сказанное справедливо также для теоремы 6.2.2, следствий 6.2.3 и 6.2.4 и теорем 6.2.11 и 6.2.12.

## §6.3. СВЯЗЬ МЕЖДУ РАЗЛОЖЕНИЯМИ В < A, [ ] > И РАЗЛОЖЕНИЯМИ В A*

В данном параграфе устанавливается связь между разложением n-арной группы $< A, [\ ] >$ по ее n-арной подгруппе $< B, [\ ] >$ и разложением универсальной обертывающей группы Поста $A^*$ по её подгруппе, изоморфной универсальной обертывающей группе Поста $B^*$.

**6.3.1. Теорема** [111]. Пусть $< B, [\ ] >$ – n-арная подгруппа n-арной группы $< A, [\ ] >$. Тогда:

1) если

$$A = \bigcup_{i \in I} [x_i \underbrace{B \ldots B}_{n-1}] \tag{1.1}$$

– разложение $< A, [\ ] >$ на непересекающиеся левые смежные классы по $< B, [\ ] >$, то

$$A^* = \bigcup_{i \in I} \theta_A(x_i) B^*(A) \tag{1.2}$$

– разложение $A^*$ на непересекающиеся левые смежные классы по $B^*(A)$, а отображение

$$[x_i \underbrace{B \ldots B}_{n-1}] \to \theta_A(x_i) B^*(A) \tag{1.3}$$

является биекцией множества всех левых смежных классов $<A, [\ ]>$ по $<B, [\ ]>$ на множество всех левых смежных классов $A^*$ по $B^*(A)$;

2) если (1.2) – разложение $A^*$ на непересекающиеся левые смежные классы по $B^*(A)$, то (1.1) – разложение $<A, [\ ]>$ на непересекающиеся левые смежные классы по $<B, [\ ]>$, а отображение

$$\theta_A(x_i)B^*(A) \to [x_i \underbrace{B \ldots B}_{n-1}] \tag{1.4}$$

является биекцией множества всех левых смежных классов $A^*$ по $B^*(A)$ на множество всех левых смежных классов $<A, [\ ]>$ по $<B, [\ ]>$.

***Доказательство.*** 1) Пусть $\theta_A(a_1 \ldots a_k)$ – произвольный элемент из группы $A^*$, $k = 1, \ldots, n-1$. Если зафиксировать $b_1, \ldots, b_{k-1} \in B$, то найдется $y \in A$ такой, что

$$\theta_A(a_1 \ldots a_k) = \theta_A(yb_1 \ldots b_{k-1}). \tag{1.5}$$

Если $b_k, \ldots, b_{n-1} \in B$, то

$$[yb_1 \ldots b_{k-1}b_k \ldots b_{n-1}] \in [x_i \underbrace{B \ldots B}_{n-1}]$$

для некоторого $i \in I$, откуда

$$[yb_1 \ldots b_{k-1}b_k \ldots b_{n-1}] = [x_i b_1 \ldots b_{k-1}b_k \ldots b_{n-2}b]$$

для некоторого $b \in B$. Тогда

$$\theta_A(yb_1 \ldots b_{k-1}b_k \ldots b_{n-1}) = \theta_A(x_i b_1 \ldots b_{n-2}b),$$

$$\theta_A(yb_1 \ldots b_{k-1})\theta_A(b_k \ldots b_{n-1}) = \theta_A(x_i)\theta_A(b_1 \ldots b_{n-2}b),$$

$$\theta_A(yb_1 \ldots b_{k-1}) = \theta_A(x_i)\,\theta_A(b_1 \ldots b_{n-2}b)\theta_A^{-1}(b_k \ldots b_{n-1}),$$

$$\theta_A(yb_1 \ldots b_{k-1}) \in \theta_A(x_i)B^*(A),$$

откуда и из (1.5) следует $\theta_A(a_1 \ldots a_k) \in \theta_A(x_i)B^*(A)$. Следовательно,

$$A^* \subseteq \bigcup_{i \in I} \theta_A(x_i)B^*(A).$$

Обратное включение

$$\bigcup_{i \in I} \theta_A(x_i)B^*(A) \subseteq A^*$$

очевидно. Таким образом, доказано равенство (1.2).

Предположим, что

$$\theta_A(x_i)B^*(A) \cap \theta_A(x_j)B^*(A) \neq \varnothing, \; i \neq j,$$

то есть

$$\theta_A(x_i)\theta_A(c_1 \ldots c_k) = \theta_A(x_j)\theta_A(d_1 \ldots d_m)$$

для $c_1, \ldots, c_k, d_1, \ldots, d_m \in B$, где $k, m \in \{1, \ldots, n-1\}$. Ясно, что $k = m$.

Если $k = n - 1$, то из последнего равенства следует

$$[x_i c_1 \ldots c_{n-1}] = [x_j d_1 \ldots d_{n-1}].$$

Полученное равенство противоречит тому, что (1.1) – разложение $< A, [\,] >$ на непересекающиеся левые смежные классы по $< B, [\,] >$.

Если $k < n - 1$, то

$$\theta_A(x_i)\theta_A(c_1 \ldots c_k)\theta_A(c_{k+1} \ldots c_{n-1}) =$$
$$= \theta_A(x_j)\theta_A(d_1 \ldots d_k)\theta_A(c_{k+1} \ldots c_{n-1})$$

для любых $c_{k+1}, \ldots, c_{n-1} \in B$, откуда

$$[x_i c_1 \ldots c_{n-1}] = [x_j d_1 \ldots d_k c_{k+1} \ldots c_{n-1}].$$

Последнее равенство противоречит тому, что (1.1) – разложение $< A, [\,] >$ на непересекающиеся левые смежные классы по

< B, [ ] >. Следовательно, равенство (1.2) является разложением A* на непересекающиеся левые смежные классы по B*(A).

Из доказанного следует, что (1.3) – биекция.

2) Пусть a – произвольный элемент из A. Тогда, согласно условию 2) теоремы, $\theta_A(a) \in \theta_A(x_i)B^*(A)$ для некоторого $i \in I$, откуда $a = [x_i b_1 \ldots b_{n-1}]$ для некоторых $b_1, \ldots, b_{n-1} \in B$, то есть в действительности $\theta_A(a) \in \theta_A(x_i)B_o(A)$ Следовательно,

$$A \subseteq \bigcup_{i \in I} [x_i \underbrace{B \ldots B}_{n-1}].$$

Обратное включение

$$\bigcup_{i \in I} [x_i \underbrace{B \ldots B}_{n-1}] \subseteq A$$

очевидно. Таким, образом, доказано равенство (1.1).

Если

$$[x_i \underbrace{B \ldots B}_{n-1}] \cap [x_j \underbrace{B \ldots B}_{n-1}] \neq \varnothing, i \neq j,$$

то

$$\theta_A(x_i)B_o(A) \cap \theta_A(x_j)B_o(A) \neq \varnothing,$$

откуда, учитывая

$$\theta_A(x_i)B_o(A) \subseteq \theta_A(x_i)B^*(A), \theta_A(x_j)B_o(A) \subseteq \theta_A(x_j)B^*(A),$$

получаем

$$\theta_A(x_i)B^*(A) \cap_A(x_j)B^*(A) \neq \varnothing,$$

что противоречит тому, что (1.2) – разложение A* на непересекающиеся левые смежные классы по B*(A).

Ясно, что (1.4) – биекция. ∎

Аналогично теореме 6.3.1 доказывается "правая" теорема.

**6.3.2. Теорема** [111]. Пусть $< B, [\ ] >$ – n-арная подгруппа n-арной группы $< A, [\ ] >$. Тогда:

1) если

$$A = \bigcup_{i \in I} [\underbrace{B \ldots B}_{n-1} x_i] \qquad (2.1)$$

– разложение $< A, [\ ] >$ на непересекающиеся правые смежные классы по $< B, [\ ] >$, то

$$A^* = \bigcup_{i \in I} B^*(A)\theta_A(x_i) \qquad (2.2)$$

– разложение $A^*$ на непересекающиеся правые смежные классы по $B^*(A)$, а отображение

$$[\underbrace{B \ldots B}_{n-1} x_i] \to B^*(A)\theta_A(x_i) \qquad (2.3)$$

является биекцией множества всех правых смежных классов $< A, [\ ] >$ по $< B, [\ ] >$ на множество всех правых смежных классов $A^*$ по $B^*(A)$;

2) если (2.2) – разложение $A^*$ на непересекающиеся правые смежные классы по $B^*(A)$, то (2.1) – разложение $< A, [\ ] >$ на непересекающиеся правые смежные классы по $< B, [\ ] >$, а отображение

$$B^*(A)\theta_A(x_i) \to [\underbrace{B \ldots B}_{n-1} x_i] \qquad (2.4)$$

является биекцией множества всех правых смежных классов $A^*$ по $B^*(A)$ на множество всех правых смежных классов $< A, [\ ] >$ по $< B, [\ ] >$.

**6.3.3. Замечание.** Ясно, что отображения (1.3) и (1.4) являются взаимно обратными. То же самое можно сказать об отображениях (2.3) и (2.4).

Из любой из теорем 6.3.1 или 6.3.2 вытекает

**6.3.4. Следствие**. Индекс n-арной подгруппы $<B,[\ ]>$ в n-арной группе $<A,[\ ]>$ совпадает с индексом подгруппы $B^*(A)$ в группе $A^*$: $|A:B| = |A^*:B^*(A)|$.

## §6.4. ИЗОМОРФИЗМ ГРУПП $A^*/B^*(A)$ и $A_o/B_o(A)$

Если $<B,[\ ]>$ – инвариантная n-арная подгруппа n-арной группы $<A,[\ ]>$, то, как будет показано ниже, подгруппа $B^*(A)$ инвариантна в группе $A^*$ (предложение 7.3.14), а подгруппа $B_o(A)$ инвариантна в группе $A_o$ (предложение 7.3.15). Таким образом, в случае инвариантности $<B,[\ ]>$ в $<A,[\ ]>$ можно рассматривать факторгруппы $A^*/B^*(A)$ и $A_o/B_o(A)$.

Из следствий 6.2.6 и 6.3.4 вытекает равенство мощностей n-арной факторгруппы $<A/B,[\ ]>$ и факторгрупп $A^*/B^*(A)$ и $A_o/B_o(A)$. Совпадение мощностей факторгрупп $A^*/B^*(A)$ и $A_o/B_o(A)$ наводит на мысль об их возможном изоморфизме.

Покажем, что это действительно так.

**6.4.1. Лемма.** Если $<B,[\ ]>$ – n-арная подгруппа n-арной группы $<A,[\ ]>$, то $A^* = B^*(A)A_o$.

*Доказательство.* Зафиксируем $b \in B$. Тогда

$$\theta_A(b), \theta_A(bb) = \theta_A^2(b), \ldots, \theta_A(\underbrace{b \ldots b}_{n-2}) = \theta_A^{n-2}(b) \in B^*(A),$$

откуда

$$A_o \cup \theta_A(b)A_o \cup \ldots \cup \theta_A^{n-2}(b)A_o \in B^*(A)A_o.$$

А так как, согласно предложению 1.4.6,

$$A_o \cup \theta_A(b)A_o \cup \ldots \cup \theta_A^{n-2}(b)A_o = A^*,$$

то $A^* \subseteq B^*(A)A_o$. Обратное включение очевидно. ∎

Следующая лемма является следствием определений.

**6.4.2. Лемма.** Если $< B, [\ ] > $ – n-арная подгруппа n-арной группы $< A, [\ ] >$, то

$$B_o(A) = B^*(A) \cap A_o.$$

**6.4.3. Теорема** [111, 112]**.** Если $< B, [\ ] >$ – инвариантная n-арная подгруппа n-арной группы $< A, [\ ] >$, то факторгруппы $A^*/B^*(A)$ и $A_o/B_o(A)$ изоморфны.

*Доказательство.* Применяя леммы 6.4.1 и 6.4.2, а также первую теорему об изоморфизмах для групп, получим

$$A^*/B^*(A) = B^*(A)A_o/B^*(A) \simeq A_o/B^*(A) \cap A_o = A_o/B_o(A). \quad \blacksquare$$

**6.4.4. Замечание.** Если отождествить группу $B^*(A)$ с группой $B^*$, а группу $B_o(A)$ с группой $B_o$, то изоморфизм из теоремы 6.4.3 может быть записан более компактно

$$A^*/B^* \simeq A_o/B_o.$$

В доказательстве теоремы 6.4.3 явный вид изоморфизма факторгрупп $A^*/B^*(A)$ и $A_o/B_o(A)$ не указан. Для нахождения явного вида этого изоморфизма нам понадобиться несколько вспомогательных результатов.

**6.4.5. Лемма** [111]**.** Пусть $< B, [\ ] >$ – n-арная подгруппа n-арной группы $< A, [\ ] >$, $b_1, \ldots, b_{n-2}$ – фиксированные элементы из B, $x, y \in A$. Тогда следующие равенства равносильны:

1) $[x \underbrace{B \ldots B}_{n-1}] = [y \underbrace{B \ldots B}_{n-1}]$;

2) $\theta_A(x)B^*(A) = \theta_A(y)B^*(A)$;

3) $\theta_A(xb_1 \ldots b_{n-2})B_o(A) = \theta_A(yb_1 \ldots b_{n-2})B_o(A)$.

*Доказательство.* 1) $\Rightarrow$ 2) Так как

$$x = [x\overline{b\underbrace{b \dots b}_{n-2}}] \in [x\underbrace{B \dots B}_{n-1}]$$

для любого $b \in B$, то из 1) следует $x = [yc_1 \dots c_{n-1}]$ для некоторых $c_1, \dots, c_{n-1} \in B$. Тогда

$$\theta_A(x)B^*(A) = \theta([yc_1 \dots c_{n-1}])B^*(A) =$$

$$= \theta_A(y)\theta_A(c_1 \dots c_{n-1})B^*(A) = \theta_A(y)B^*(A).$$

2)$\Rightarrow$3) Так как $\theta_A(x) \in \theta_A(x)B^*(A)$, то из 2) следует

$$\theta_A(x) = \theta_A(y)\theta_A(d_1 \dots d_i)$$

для некоторых $d_1, \dots, d_i \in B$. А так как $x, y \in A$, то в последнем равенстве можно считать $i = n - 1$, то есть

$$x = [yd_1 \dots d_{n-1}].$$

В B всегда найдется элемент $d$ такой, что

$$d_1 \dots d_{n-1}\theta_A b_1 \dots b_{n-2}d,$$

откуда

$$x = [yd_1 \dots d_{n-1}] = [yb_1 \dots b_{n-2}d].$$

Тогда

$$\theta_A(xb_1 \dots b_{n-2})B_o(A) = \theta_A([yb_1 \dots b_{n-2}d]b_1 \dots b_{n-2})B_o(A) =$$

$$= \theta_A(yb_1 \dots b_{n-2})\theta_A(db_1 \dots b_{n-2})B_o(A) = \theta_A(yb_1 \dots b_{n-2})B_o(A).$$

3)$\Rightarrow$1) Так как

$$\theta_A(xb_1 \dots b_{n-2}) \in \theta_A(xb_1 \dots b_{n-2})B_o(A),$$

то из 3) следует

$$\theta_A(xb_1 \dots b_{n-2}) = \theta_A(yb_1 \dots b_{n-2})\theta_A(g_1 \dots g_{n-1})$$

для некоторых $g_1, \dots, g_{n-1} \in B$, откуда

$$\theta_A(x) = \theta_A(yb_1 \ldots b_{n-2})\theta_A(g_1 \ldots g_{n-1})\theta_A^{-1}(b_1 \ldots b_{n-2}),$$

$$x = [yb_1 \ldots b_{n-2}g_1 \ldots g_{n-1}b],$$

где b – обратный элемент для последовательности $b_1 \ldots b_{n-2}$. Ясно, что $b \in B$. Тогда

$$[x\underbrace{B \ldots B}_{n-1}] = [[yb_1 \ldots b_{n-2}g_1 \ldots g_{n-1}b]\underbrace{B \ldots B}_{n-1}] = [y\underbrace{B \ldots B}_{n-1}]. \quad \blacksquare$$

Аналогично лемме 6.4.5 доказывается "правая" лемма.

**6.4.6. Лемма** [111]. Пусть $< B, [\ ] >$ – n-арная подгруппа n-арной группы $< A, [\ ] >$, $b_1, \ldots, b_{n-2}$ – фиксированные элементы из B, $x, y \in A$. Тогда следующие равенства равносильны:

1) $[\underbrace{B \ldots B}_{n-1} x] = [\underbrace{B \ldots B}_{n-1} y]$;
2) $B^*(A)\theta_A(x) = B^*(A)\theta_A(y)$;
3) $B_o(A)\theta_A(b_1 \ldots b_{n-2}x) = B_o(A)\theta_A(b_1 \ldots b_{n-2}y)$.

**6.4.7. Предложение** [111]. Пусть $< B, [\ ] >$ – инвариантная n-арная подгруппа n-арной группы $< A, [\ ] >$, $b_1, \ldots, b_{n-2} \in B$,

$$A/B = \{[x_i\underbrace{B \ldots B}_{n-1}] \mid i \in I\}.$$

Тогда отображение

$$\varphi : \theta_A(x_i)B^*(A) \to \theta_A(x_ib_1 \ldots b_{n-2})B_o(A)$$

является изоморфизмом группы $A^*/B^*(A)$ на группу $A_o/B_o(A)$.

*Доказательство.* По теореме 6.3.1

$$A^*/B^* = \{\theta_A(x_i)B^*(A) \mid i \in I\},$$

и

$$\varphi_1 : \theta_A(x_i)B^*(A) \to [x_i \underbrace{B \ldots B}_{n-1}]$$

– биекция $A^*/B^*$ на $A/B$, а ввиду 1) следствия 6.2.3,

$$A_o/B_o = \{\theta_A(x_i b_1 \ldots b_{n-2})B_o(A) \mid i \in I\},$$

и

$$\varphi_2 : [x_i \underbrace{B \ldots B}_{n-1}] \to \theta_A(x_i b_1 \ldots b_{n-2})B_o(A)$$

– биекция $A/B$ на $A_o/B_o$.

Так как $\varphi = \varphi_1 \varphi_2$, то $\varphi$ является биекцией. Кроме того, если $\theta_A(x_i)B^*(A)$ и $\theta_A(x_j)B^*(A)$ – произвольные элементы из $A^*/B^*(A)$, то, учитывая инвариантность $B^*(A)$ в $A^*$, а также то, что $b_1, \ldots, b_{n-2} \in B$, получим

$$\varphi(\theta_A(x_i)B^*(A)\theta_A(x_j)B^*(A)) = \varphi(\theta_A(x_i)\theta_A(x_j)B^*(A)B^*(A)) =$$
$$= \varphi(\theta_A(x_i)\theta_A(x_j)B^*(A)) =$$
$$= \varphi(\theta_A(x_i)\theta_A(x_j)\theta_A(b_1 \ldots b_{n-2})B^*(A)) =$$
$$= \varphi(\theta_A([x_i x_j b_1 \ldots b_{n-2}])B^*(A)),$$

то есть

$$\varphi(\theta_A(x_i)B^*(A)\theta_A(x_j)B^*(A)) = \varphi(\theta_A([x_i x_j b_1 \ldots b_{n-2}])B^*(A)).$$

Так как

$$\theta_A([x_i x_j b_1 \ldots b_{n-2}])B^*(A)) = \theta_A(x_k)B^*(A)$$

для некоторого $k \in I$, то, используя лемму 6.4.5, инвариантность $<B, [\,]>$ в $<A, [\,]>$ и инвариантность $B_o(A)$ в $A_o$, получим

$$\varphi(\theta_A(x_i)B^*(A)\theta_A(x_j)B^*(A)) = \varphi(\theta_A(x_k)B^*(A)) =$$
$$= \theta_A(x_k b_1 \ldots b_{n-2})B_o(A) = \theta_A([x_i x_j b_1 \ldots b_{n-2}] b_1 \ldots b_{n-2})B_o(A) =$$

$$= \theta_A(x_i b'_1 \ldots b'_{n-2} x_j b_1 \ldots b_{n-2}) B_o(A) =$$

$$= \theta_A(x_i b'_1 \ldots b'_{n-2}) \theta_A(x_j b_1 \ldots b_{n-2}) B_o(A) B_o(A) =$$

$$= \theta_A(x_i b'_1 \ldots b'_{n-2}) B_o(A) \theta_A(x_j b_1 \ldots b_{n-2}) B_o(A) =$$

$$= \theta_A(x_i) \theta_A(b_1 \ldots b_{n-2}) \theta_A(b b'_1 \ldots b'_{n-2}) B_o(A) \theta_A(x_j b_1 \ldots b_{n-2}) B_o(A) =$$

$$= \theta_A(x_i b_1 \ldots b_{n-2}) B_o(A) \theta_A(x_j b_1 \ldots b_{n-2}) B_o(A) =$$

$$= \varphi(\theta_A(x_i) B^*(A)) \varphi(\theta_A(x_j) B^*(A)),$$

где $b'_1, \ldots, b'_{n-2} \in B$, $b$ – обратный элемент для последовательности $b_1 \ldots b_{n-2}$. Следовательно, $\varphi$ – изоморфизм группы $A^*/B^*(A)$ на группу $A_o/B_o(A)$. ∎

Следующее предложение получается с использованием теоремы 6.3.2, леммы 6.4.6 и утверждения 1) следствия 6.2.4.

**6.4.8. Предложение** [111]. Пусть $< B, [\ ] >$ – инвариантная n-арная подгруппа n-арной группы $< A, [\ ] >$, $b_1, \ldots, b_{n-2} \in B$,

$$A/B = \{[\underbrace{B \ldots B}_{n-1} x_i] \mid i \in I\}.$$

Тогда отображение

$$\psi : B^*(A) \theta_A(x_i) \to B_o(A) \theta_A(b_1 \ldots b_{n-2} x_i)$$

является изоморфизмом группы $A^*/B^*(A)$ на группу $A_o/B_o(A)$.

## ДОПОЛНЕНИЯ И КОММЕНТАРИИ

**1.** Теорема 6.1.1 остается верной, если так же, как и в бинарном случае, не требовать конечности индекса $< B, [\ ] >$ в $< A, [\ ] >$ (см., например, [113]). В этом случае под подстановкой понимается взаимно однозначное отображение.

**2.** Так как множество $\Delta = \{\delta_a \mid a \in A\}$ содержит тождественное преобразование $\delta_b$, то можно определить множество $T$ всех элементов из $A$, которые отображение $\gamma: a \to a^\gamma = \delta_a$ переводит в тождественное преобразование. Так как согласно 1) теоремы 6.1.1, $\gamma$ – гомоморфизм группы $<A,[\ ]>$ на группу $\Delta$, то $T$ совпадает с ядром этого гомоморфизма.

**3.** Отмечавшееся при доказательстве теоремы 6.1.1 включение $\ker \gamma = N \subseteq B$ можно доказать непосредственно, используя определения множества $T$ и отображения $\gamma$.

Если $u \in T$, то согласно определениям множества $T$ и отображения $\gamma$, $u^\gamma = \delta_u = \varepsilon$ – тождественное преобразование. Поэтому, учитывая нейтральность последовательности $b\widetilde{b}$, получим

$$\delta_u([\underbrace{B\ldots B}_{n-1}b]) = [\underbrace{B\ldots B}_{n-1}b],$$

$$[[\underbrace{B\ldots B}_{n-1}b]\widetilde{b}\,u] = B,$$

$$[\underbrace{B\ldots B}_{n-1}u] = B,$$

откуда $u \in B$. Следовательно, $T \subseteq B$. Из последнего включения и равенства $T = \ker \gamma$, полученного в предыдущем пункте, следует $\ker \gamma \subseteq B$.

**4.** О связи между разложением $A^*$ по $B_o(A)$ и разложением $<A,[\ ]>$ по $<B,[\ ]>$ писал Пост [3, с. 223], отождествляя при этом $B_o(A)$ и $B_o$.

**5.** Утверждения 1) следствия 6.2.3 и теоремы 6.3.1 использовались в [112] для получения n-арных аналогов теоремы Шура о конечности коммутанта группы, центр которой имеет в ней конечный индекс.

**6.** Для конечных n-арных групп равенства индексов из следствий 6.2.6 и 6.3.4 могут быть получены без использования теорем 6.2.1 и 6.3.1 соответственно:

$$|A^* : B^*(A)| = |A^*| : |B^*(A)| = |A^*| : |B^*| =$$

$$= |A|(n-1) : |B|(n-1) = |A| : |B| = |A : B|,$$

$$|A_o : B_o(A)| = |A_o| : |B_o(A)| = |A_o| : |B_o| = |A| : |B| = |A : B|.$$

**7.** Так как из инвариантности n-арной подгруппы $< B, [\ ] >$ в n-арной группе $< A, [\ ] >$ следует инвариантность $B^*(A)$ в $A^*$ и $B_o(A)$ в $A_o$, то изоморфизмы $\varphi$ и $\psi$ из предложений 6.4.7 и 6.4.8 совпадают.

# ГЛАВА 7

## n-АРНЫЕ АНАЛОГИ НОРМАЛИЗАТОРА ПОДМНОЖЕСТВА В ГРУППЕ

В данной главе определяются различные n-арные аналоги нормализатора подмножества в группе и изучается их связь со своими бинарными прототипами. Получены также новые критерии инвариантности и полуинвариантности n-арной подгруппы в n-арной группе.

### §7.1. ОПРЕДЕЛЕНИЯ

В [12] для всякой n-арной подгруппы $< B, [\ ] >$ n-арной группы $< A, [\ ] >$, где $n = k(m - 1) + 1$, $k \geq 1$, введено понятие m-полунормализатора

$$N_A(B, m) = \{x \in A \mid [x\underbrace{B\ldots B}_{n-1}] = [\underbrace{B\ldots B}_{m-1} x \underbrace{B\ldots B}_{n-m}] = [\underbrace{B\ldots B}_{n-1} x]\}.$$

Если $N_A(B, m) = A$, то n-арная подгруппа $< B, [\ ] >$ называется [12, 28] m-полуинвариантной в $< A, [\ ] >$. Таким образом, n-арная подгруппа $< B, [\ ] >$ n-арной группы $< A, [\ ] >$ называется m-полуинвариантной. в ней, если

$$[x\underbrace{B\ldots B}_{n-1}] = [\underbrace{B\ldots B}_{m-1} x \underbrace{B\ldots B}_{n-m}] = [\underbrace{B\ldots B}_{n-1} x] \qquad (*)$$

для любого $x \in A$, $n = k(m - 1) + 1$, $k \geq 1$.

**7.1.1 Лемма.** Пусть n-арная подгруппа $< B, [\ ] >$ n-арной группы $< A, [\ ] >$, где $n = k(m - 1) + 1$, $k \geq 1$, удовлетворяет условию (*) для некоторого $x \in A$. Тогда

$$[x\underbrace{B\ldots B}_{n-1}] = [\underbrace{B\ldots B}_{m-1}x\underbrace{B \ldots B}_{(k-1)(m-1)}] = [\underbrace{B\ldots B}_{2(m-1)}x\underbrace{B \ldots B}_{n-1-2(m-1)}] = \ldots$$

$$\ldots = [\underbrace{B \ldots B}_{(k-1)(m-1)}x\underbrace{B\ldots B}_{m-1}] = [\underbrace{B\ldots B}_{n-1}x]. \qquad (**)$$

*Доказательство.* Так как $n - m = (k - 1)(n - 1)$, то используя равенства

$$[\underbrace{B\ldots B}_{n}] = B, \ [x\underbrace{B\ldots B}_{n-1}] = [\underbrace{B\ldots B}_{m-1}x\underbrace{B\ldots B}_{n-m}],$$

получим

$$[x\underbrace{B\ldots B}_{n-1}] = [\underbrace{B\ldots B}_{m-1}x\underbrace{B \ldots B}_{(k-1)(m-1)}] = [\underbrace{B\ldots B}_{m-1}x[\underbrace{B\ldots B}_{n}]\underbrace{B\ldots B}_{n-m-1}] =$$

$$= [\underbrace{B\ldots B}_{m-1}[x\underbrace{B\ldots B}_{n-1}]\underbrace{B\ldots B}_{n-m}] = [\underbrace{B\ldots B}_{m-1}[\underbrace{B\ldots B}_{m-1}x\underbrace{B\ldots B}_{n-m}]\underbrace{B\ldots B}_{n-m}] =$$

$$= [\underbrace{B\ldots B}_{2(m-1)}x\underbrace{B\ldots B}_{2n-2m}] = [\underbrace{B\ldots B}_{2(m-1)}x[\underbrace{B\ldots B}_{n}]\underbrace{B \ldots B}_{n-2-2(m-1)}] =$$

$$= [\underbrace{B\ldots B}_{2(m-1)}x\underbrace{B \ldots B}_{n-1-2(m-1)}] = \ldots = [\underbrace{B \ldots B}_{(k-1)(m-1)}x\underbrace{B\ldots B}_{m-1}].$$

Равенство

$$[x\underbrace{B\ldots B}_{n-1}] = [\underbrace{B\ldots B}_{n-1}x]$$

следует из (*). ∎

**7.1.2. Следствие.** Пусть n-арная подгруппа $< B, [\ ] >$ n-арной группы $< A, [\ ] >$, где $n = k(m - 1) + 1$, $k \geq 1$, удовлетворяет условию (*) для некоторого $x \in A$ и пусть $s \geq 1$. Тогда

$$[x\underbrace{B\ldots B}_{n-1}] = [x\underbrace{B\ldots B}_{s(n-1)}] = [\underbrace{B\ldots B}_{i(m-1)}x\underbrace{B \ldots B}_{s(n-1)-i(m-1)}]$$

для любого i = 1, …, sk – 1.

**7.1.3. Следствие.** Если n-арная подгруппа $< B, [\,] >$ n-арной группы $< A, [\,] >$, где $n = k(m - 1) + 1$, $k \geq 1$, m-полуинвариантна в $< A, [\,] >$, то она удовлетворяет условию (**) для любого $x \in A$.

Таким образом, имеет место

**7.1.4. Предложение.** n-Арная подгруппа $< B, [\,] >$ n-арной группы $< A, [\,] >$, где $n = k(m - 1) + 1$, $k \geq 1$, m-полуинвариантна в $< A, [\,] >$ тогда и только тогда, когда она удовлетворяет условию (**) для любого $x \in A$. В частности, n-арная подгруппа $< B, [\,] >$ n-арной группы $< A, [\,] >$ является инвариантной в ней тогда и только тогда, когда она 2-полуинвариантна в ней.

Теперь мы можем расширить понятия m-полуинвариантности и m-полунормализатора, рассматривая подмножества n-арной группы.

**7.1.5. Определение.** Подмножество B n-арной группы $< A, [\,] >$, где $n = k(m - 1) + 1$, $k \geq 1$, называется *m-полуинвариантным* в ней, если

$$[x\underbrace{B\ldots B}_{n-1}] = [\underbrace{B\ldots B}_{i(m-1)} x \underbrace{B \ldots B}_{n-1-i(m-1)}]$$

для любого $x \in A$ и всех $i = 1, \ldots, k$.

**7.1.6. Определение.** *m-Полунормализатором* подмножества B в n-арной группе $< A, [\,] >$, где $n = k(m - 1) + 1$, $k \geq 1$, называется множество

$$N_A(B, m) = \{x \in A \mid [x\underbrace{B\ldots B}_{n-1}] = [\underbrace{B\ldots B}_{i(m-1)} x \underbrace{B \ldots B}_{n-1-i(m-1)}], \forall i = 1, \ldots, k\}.$$

Для 2-полунормализатора $N_A(B, 2)$ будем употреблять обозначение $N_A(B)$ и называть его *нормализатором* подмно-жества B в n-арной группе $< A, [\ ] >$. Таким образом,

$$N_A(B) = \{x \in A \mid [x\underbrace{B\ldots B}_{n-1}] = [\underbrace{B\ldots B}_{i} x \underbrace{B\ldots B}_{n-1-i}], \forall i = 1, \ldots, n-1\}.$$

n-Полунормализатор $N_A(B, n)$ будем обозначать через $HN_A(B)$ и называть *полунормализатором* [4]. Таким образом,

$$HN_A(B) = \{x \in A \mid [x\underbrace{B\ldots B}_{n-1}] = [\underbrace{B\ldots B}_{n-1} x]\}.$$

Ясно, что для n-арной подгруппы $< B, [\ ] >$ n-арной группы $< A, [\ ] >$ верно включение $B \subseteq N_A(B, m)$.

Из леммы 7.1.1 вытекает

**7.1.7. Следствие.** Если $< B, [\ ] >$ – n-арная подгруппа n-арной группы $< A, [\ ] >$, где $n = k(m-1) + 1$, $k \geq 1$, то

$$N_A(B, m) = \{x \in A \mid [x\underbrace{B\ldots B}_{n-1}] = [\underbrace{B\ldots B}_{m-1} x \underbrace{B\ldots B}_{n-m}] = [\underbrace{B\ldots B}_{n-1} x]\}.$$

В частности,

$$N_A(B) = \{x \in A \mid [x\underbrace{B\ldots B}_{n-1}] = [Bx\underbrace{B\ldots B}_{n-1}] = [\underbrace{B\ldots B}_{n-1} x]\}.$$

## §7.2. СВОЙСТВА m-ПОЛУНОРМАЛИЗАТОРА

**7.2.1. Теорема** [12]. Пусть $< B, [\ ] >$ – n-арная подгруппа n-арной группы $< A, [\ ] >$, где $n = s(m-1) + 1$, $s \geq 1$. Тогда $< N_A(B, m), [\ ] >$ – n-арная подгруппа в $< A, [\ ] >$, причём

$$N_A(B, m) \subseteq N_A(B, k),$$

где $m - 1$ делит $k - 1$, $k - 1$ делит $n - 1$.

*Доказательство.* Если $x_1, \ldots, x_n \in N_A(B, m)$, то

$$[[x_1 \ldots x_n] \underbrace{B \ldots B}_{n-1}] = [x_1 \ldots x_{n-1}[x_n \underbrace{B \ldots B}_{n-1}]] =$$

$$= [x_1 \ldots x_{n-1}[\underbrace{B \ldots B}_{m-1} x_n \underbrace{B \ldots B}_{n-m}]] =$$

$$= [x_1 \ldots x_{n-1}[\underbrace{B \ldots B}_{n}]\underbrace{B \ldots B}_{m-2} x_n \underbrace{B \ldots B}_{n-m}] =$$

$$= [x_1 \ldots x_{n-2}[x_{n-1}\underbrace{B \ldots B}_{n-1}]\underbrace{B \ldots B}_{m-1} x_n \underbrace{B \ldots B}_{n-m}] =$$

$$= [x_1 \ldots x_{n-2}[\underbrace{B \ldots B}_{m-1} x_{n-1}\underbrace{B \ldots B}_{n-m}]\underbrace{B \ldots B}_{m-1} x_n \underbrace{B \ldots B}_{n-m}] =$$

$$= [x_1 \ldots x_{n-2}\underbrace{B \ldots B}_{m-1}[x_{n-1}\underbrace{B \ldots B}_{n-1}]x_n \underbrace{B \ldots B}_{n-m}] =$$

$$= [x_1 \ldots x_{n-2}\underbrace{B \ldots B}_{m-1}[\underbrace{B \ldots B}_{n-1} x_{n-1}]x_n \underbrace{B \ldots B}_{n-m}] =$$

$$= [x_1 \ldots x_{n-2}\underbrace{B \ldots B}_{m-1}[\underbrace{B \ldots B}_{n}]x_{n-1} x_n \underbrace{B \ldots B}_{n-m}] =$$

$$= [x_1 \ldots x_{n-2}\underbrace{B \ldots B}_{m-1} x_{n-1} x_n \underbrace{B \ldots B}_{n-m}] = \ldots$$

$$\ldots = [\underbrace{B \ldots B}_{m-1}[x_1 \ldots x_n]\underbrace{B \ldots B}_{n-m}],$$

то есть

$$[[x_1 \ldots x_n]\underbrace{B \ldots B}_{n-1}]] = [\underbrace{B \ldots B}_{m-1}[x_1 \ldots x_n]\underbrace{B \ldots B}_{n-m}]; \quad (1)$$

$$[[x_1 \ldots x_n]\underbrace{B \ldots B}_{n-1}] = [x_1 \ldots x_{n-1}[x_n\underbrace{B \ldots B}_{n-1}]] =$$

$$= [x_1 \ldots x_{n-1}[\underbrace{B \ldots B}_{n-1} x_n]] = \ldots = [\underbrace{B \ldots B}_{n-1}[x_1 \ldots x_n]],$$

то есть

$$[[x_1 \ldots x_n]\underbrace{B\ldots B}_{n-1}] = [\underbrace{B\ldots B}_{n-1}[x_1 \ldots x_n]]. \qquad (2)$$

Из (1), (2) и следствия 7.1.7 вытекает

$$[x_1 \ldots x_n] \in N_A(B, m). \qquad (3)$$

Если теперь $x \in N_A(B, m)$, то учитывая нейтральность последовательностей

$$\underbrace{x\ldots x}_{n-2}\overline{x}, \quad \overline{x}\underbrace{x\ldots x}_{n-2},$$

получим

$$[\overline{x}\underbrace{B\ldots B}_{n-1}] = [\overline{x}\underbrace{B\ldots B}_{n-1}\underbrace{x\ldots x}_{n-2}\overline{x}] = [\overline{x}\,[\underbrace{B\ldots B}_{n-1}x]\underbrace{x\ldots x}_{n-3}\overline{x}] =$$

$$= [\overline{x}\,[x\underbrace{B\ldots B}_{n-1}]\underbrace{x\ldots x}_{n-3}\overline{x}] = [\overline{x}x\underbrace{B\ldots B}_{n-1}\underbrace{x\ldots x}_{n-3}\overline{x}] = \ldots$$

$$\ldots = [\overline{x}\underbrace{x\ldots x}_{n-2}\underbrace{B\ldots B}_{n-1}\overline{x}] = [\underbrace{B\ldots B}_{n-1}\overline{x}],$$

то есть

$$[\overline{x}\underbrace{B\ldots B}_{n-1}] = [\underbrace{B\ldots B}_{n-1}\overline{x}]. \qquad (4)$$

Используя (4), получим

$$[\overline{x}\underbrace{B\ldots B}_{n-1}] = [\overline{x}\underbrace{B\ldots B}_{m-1}\underbrace{x\ldots x}_{n-2}\overline{x}\underbrace{B\ldots B}_{n-m}] =$$

$$= [\overline{x}\underbrace{B\ldots B}_{m-2}[\underbrace{B\ldots B}_{n}]\underbrace{x\ldots x}_{n-2}\overline{x}\underbrace{B\ldots B}_{n-m}] =$$

$$= [\overline{x}\underbrace{B\ldots B}_{m-1}[\underbrace{B\ldots B}_{n-1}x]\underbrace{x\ldots x}_{n-3}\overline{x}\underbrace{B\ldots B}_{n-m}] =$$

$$= [\overline{x}\underbrace{B\ldots B}_{m-1}[\underbrace{B\ldots B}_{m-1}x\underbrace{B\ldots B}_{n-m}]\underbrace{x\ldots x}_{n-3}\overline{x}\underbrace{B\ldots B}_{n-m}] =$$

$$= [\,\overline{x}\underbrace{B\ldots B}_{2(m-1)}x\underbrace{B\ldots B}_{m-n}\underbrace{x\ldots x}_{n-3}\overline{x}\underbrace{B\ldots B}_{n-m}\,] = \ldots$$

$$\ldots = [\,\overline{x}\underbrace{B\,\ldots\,B}_{(n-1)(m-1)}\underbrace{x\ldots x}_{n-2}\underbrace{B\,\ldots\,B}_{(n-2)(n-m)}\overline{x}\underbrace{B\ldots B}_{n-m}\,] =$$

$$= [\,\overline{x}\underbrace{B\,\ldots\,B}_{(n-1)(m-2)+1}\underbrace{B\ldots B}_{n-2}\underbrace{x\ldots x}_{n-2}\underbrace{B\,\ldots\,B}_{(n-2)(n-m)}\overline{x}\underbrace{B\ldots B}_{n-m}\,] =$$

$$= [\,\overline{x}\underbrace{B\ldots B}_{n-1}\underbrace{x\ldots x}_{n-2}\underbrace{B\,\ldots\,B}_{(n-2)(n-m)}\overline{x}\underbrace{B\ldots B}_{n-m}\,] =$$

$$= [[\,\overline{x}\underbrace{B\ldots B}_{n-1}]\underbrace{x\ldots x}_{n-2}\underbrace{B\,\ldots\,B}_{(n-2)(n-m)}\overline{x}\underbrace{B\ldots B}_{n-m}\,] =$$

$$= [[\,\underbrace{B\ldots B}_{n-1}\overline{x}\,]\underbrace{x\ldots x}_{n-2}\underbrace{B\,\ldots\,B}_{(n-2)(n-m)}\overline{x}\underbrace{B\ldots B}_{n-m}\,] =$$

$$= [\,\underbrace{B\ldots B}_{n-1}\overline{x}\underbrace{x\ldots x}_{n-2}\underbrace{B\,\ldots\,B}_{(n-2)(n-m)}\overline{x}\underbrace{B\ldots B}_{n-m}\,] =$$

$$= [\,\underbrace{B\,\ldots\,B}_{(n-1)+(n-2)(n-m)}\overline{x}\underbrace{B\ldots B}_{n-m}\,] =$$

$$= [[\,\underbrace{B\,\ldots\,B}_{(n-1)(n-m)+1}\underbrace{B\ldots B}_{m-2}\overline{x}\underbrace{B\ldots B}_{n-m} = [\,\underbrace{B\ldots B}_{m-1}\overline{x}\underbrace{B\ldots B}_{n-m}\,],$$

то есть

$$[\,\overline{x}\underbrace{B\ldots B}_{n-1}\,] = [\,\underbrace{B\ldots B}_{m-1}\overline{x}\underbrace{B\ldots B}_{n-m}\,]. \qquad (5)$$

Из (4), (5) и следствия 7.1.7 вытекает

$$\overline{x} \in N_A(B, m). \qquad (6)$$

Из (3) и (6), согласно критерию Дёрнте, следует, что $< N_A(B, m), [\,] > -$ n-арная подгруппа в $< A, [\,] >$.

Пусть теперь k − 1 = j(m − 1) и x − произвольный элемент из $N_A(B, m)$. Тогда, ввиду определения 7.1.6, имеем:

$$[x\underbrace{B\ldots B}_{n-1}] = [\underbrace{B\ldots B}_{j(m-1)} x\underbrace{B \ \ldots \ B}_{n-1-j(m-1)}] = [\underbrace{B\ldots B}_{k-1} x\underbrace{B\ldots B}_{n-k}], \ j = 1, \ldots, s;$$

$$[x\underbrace{B\ldots B}_{n-1}] = [x\underbrace{B\ldots B}_{n-1}],$$

и по следствию 7.1.7 $x \in N_A(B, k)$, откуда

$$N_A(B, m) \subseteq N_A(B, k). \qquad \blacksquare$$

**7.2.2. Следствие** [4]. Для любой n-арной подгруппы $< B, [\ ] >$ n-арной группы $< A, [\ ] >$

$$< N_A(B), [\ ] > \text{ и } < HN_A(B), [\ ] >$$

— n-арные подгруппы в $< A, [\ ] >$, причём $N_A(B) \subseteq HN_A(B)$.

**7.2.3. Теорема** [114]. Пусть $< B, [\ ] >$ — n-арная подгруппа n-арной группы $< A, [\ ] >$, причем m − 1 и k − 1 делят n − 1,

$$r - 1 = (m - 1, k - 1).$$

Тогда

$$N_A(B, r) = N_A(B, m) \cap N_A(B, k).$$

*Доказательство.* Включение

$$N_A(B, r) \subseteq N_A(B, m) \cap N_A(B, k) \qquad (1)$$

следует из теоремы 7.2.1.

Так как r − 1 = (m − 1, k − 1), то существуют целые числа α и β такие, что

$$\alpha(m - 1) + \beta(k - 1) = r - 1.$$

Пусть для определенности α > 0, β < 0, то есть

$$\alpha(m - 1) = -\beta(k - 1) + (r - 1), -\beta(k - 1) > 0.$$

Выберем целое t, удовлетворяющее неравенству

$$-\beta(k-1) < t(n-1). \qquad (2)$$

После этого можно выбрать целое s, удовлетворяющее неравенству

$$t(n-1) + \beta(k-1) < s(n-1) - \alpha(m-1). \qquad (3)$$

Из (2) следует

$$t(n-1) + \beta(k-1) > 0,$$

откуда и из (3) получаем

$$\alpha(m-1) < s(n-1). \qquad (4)$$

Если теперь

$$x \in N_A(B, m) \cap N_A(B, k),$$

то, дважды применяя следствие 7.1.2 и учитывая (2) – (4), получим

$$[x\underbrace{B\ldots B}_{n-1}] = [x\underbrace{B\ldots B}_{s(n-1)}] = [\underbrace{B\ldots B}_{\alpha(m-1)} x \underbrace{B \ \ldots \ B}_{s(n-1)-\alpha(m-1)}] =$$

$$= [\underbrace{B\ldots B}_{r-1} [\underbrace{B\ldots B}_{-\beta(k-1)} x \underbrace{B \ \ldots \ B}_{t(n-1)+\beta(k-1)}] \underbrace{B \ \ldots \ B}_{s(n-1)-\alpha(m-1)-\beta(k-1)-t(n-1)}] =$$

$$= [\underbrace{B\ldots B}_{r-1} [x\underbrace{B\ldots B}_{t(n-1)}] \underbrace{B \ \ldots \ B}_{s(n-1)-(r-1)-t(n-1)}] =$$

$$= [\underbrace{B\ldots B}_{r-1} x \underbrace{B \ \ldots \ B}_{s(n-1)-(r-1)}] = [\underbrace{B\ldots B}_{r-1} x \underbrace{B\ldots B}_{n-r}].$$

Применяя следствие 7.1.7 и учитывая верное равенство

$$[x\underbrace{B\ldots B}_{n-1}] = [\underbrace{B\ldots B}_{n-1} x],$$

получим x ∈ $N_A(B, r)$. Таким образом, доказано включение

$$N_A(B, m) \cap N_A(B, k) \subseteq N_A(B, r). \qquad (5)$$

Из (1) и (5) следует требуемое равенство. ∎

Так как $HN_A(B) = N_A(B, n)$, то из теоремы 7.2.3 вытекает

**7.2.4. Следствие** [114]. Пусть $<B, [\ ]>$ – n-арная подгруппа n-арной группы $<A, [\ ]>$. Если m – 1 делит n – 1,

$$r - 1 = (m - 1, n - 1),$$

то

$$N_A(B, r) = N_A(B, m) \cap HN_A(B).$$

**7.2.5. Следствие** [114]. Пусть $<B, [\ ]>$ – n-арная подгруппа n-арной группы $<A, [\ ]>$, и пусть m – 1 и k – 1 делят n – 1, (m – 1, k – 1) = 1. Тогда:

$$N_A(B) = N_A(B, m) \cap N_A(B, k).$$

**7.2.6. Следствие** [114]. Пусть $<B, [\ ]>$ – n-арная подгруппа n-арной группы $<A, [\ ]>$ и пусть m – 1 делит n – 1,

$$(m - 1, n - 1) = 1.$$

Тогда

$$N_A(B) = N_A(B, m) \cap HN_A(B).$$

Если в правой части равенства

$$N_A(B, r) = N_A(B, m) \cap N_A(B, k)$$

операцию ∩ заменить операцией ∨, то есть рассмотреть n-арную подгруппу, порожденную m-полунормализатором $N_A(B, m)$ и k-полунормализатором $N_A(B, k)$, то будет ли эта n-арная подгруппа t-полунормализатором для некоторого t, и как связано число t с числами m и k?

**7.2.7. Предложение.** Если $<B, [\ ]>$ – n-арная подгруппа n-арной группы $<A, [\ ]>$, $m-1$ и $k-1$ делят $n-1$, то

$$N_A(B, m) \vee N_A(B, k) \subseteq N_A(B, t),$$

где $t-1$ – наименьшее общее кратное чисел $m-1$ и $k-1$:

$$t-1 = [m-1, k-1].$$

*Доказательство.* Так как $m-1$ делит $n-1$, $k-1$ делит $n-1$, то $t-1 = [m-1, k-1]$ делит $n-1$, то есть можно рассматривать t-полунормализатор $N_A(B, t)$. Так как $t-1$ кратно $m-1$ и $k-1$, то по теореме 7.2.1

$$N_A(B, m) \subseteq N_A(B, t), N_A(B, k) \subseteq N_A(B, t),$$

откуда

$$N_A(B, m) \vee N_A(B, k) \subseteq N_A(B, t). \blacksquare$$

**7.2.8. Вопрос.** Верно ли, что

$$N_A(B, m) \vee N_A(B, k) = N_A(B, t),$$

где $t-1 = [m-1, k-1]$, то есть верно ли обращение предыдущего предложения?

С предыдущим вопросом связан следующий

**7.2.9. Вопрос.** Будет, ли совокупность

$$\{N_A(B, m) \mid m-1 \text{ делит } n-1\}$$

подрешеткой решетки всех n-арных подгрупп n-арной группы, $<A, [\ ]>$?

Ясно, что при положительном ответе на вопрос 2, $N_A(B)$ и $HN_A(B)$ будут соответственно наименьшим и наибольшим элементами этой подрешетки.

# §7.3. КРИТЕРИИ ИНВАРИАНТНОСТИ И ПОЛУИНВАРИАНТНОСТИ

**7.3.1. Теорема.** n-Арная подгруппа $<B, [\,]>$ n-арной группы $<A, [\,]>$ инвариантна в ней тогда и только тогда, когда подгруппа $B^*(A)$ инвариантна в группе $A^*$.

***Доказательство.*** *Необходимость.* Пусть

$$u = \theta_A(a\underbrace{b\ldots b}_{i-1}),\ v = \theta_A(c\underbrace{b\ldots b}_{j-1})$$

произвольные элементы из $A^*$ и $B^*(A)$ соответственно, где $i, j \in \{1, \ldots, n-1\}$, $a \in A$, $c \in B$, $b$ – фиксированный элемент из $B$.

Так как $<B, [\,]>$ инвариантна в $<A, [\,]>$, $\overline{b} \in B$ и $\overline{a}\underbrace{a\ldots a}_{n-2}$ – нейтральная последовательность, то

$$u^{-1}vu = \theta_A(\overline{b}\underbrace{b\ldots b}_{n-i-1}\overline{a}\underbrace{a\ldots a}_{n-3})\theta_A(c\underbrace{b\ldots b}_{j-1})\theta_A(a\underbrace{b\ldots b}_{i-1}) =$$

$$= \theta_A(\overline{b}\underbrace{b\ldots b}_{n-i-1}\overline{a}\underbrace{a\ldots a}_{n-3}c\underbrace{b\ldots b}_{j-1}a\underbrace{b\ldots b}_{i-1}) =$$

$$= \theta_A(\overline{b}\underbrace{b\ldots b}_{n-i-1}[\overline{a}\underbrace{a\ldots a}_{n-3}ca][\overline{a}\underbrace{a\ldots a}_{n-3}ba]$$

$$\underbrace{[\overline{a}\underbrace{a\ldots a}_{n-3}aba]\ldots[\overline{a}\underbrace{a\ldots a}_{n-3}aba]}_{j-2}\underbrace{b\ldots b}_{i-1}) =$$

$$= \theta_A(\overline{b}\underbrace{b\ldots b}_{n-i-1}c'\underbrace{b'\ldots b'}_{j-1}\underbrace{b\ldots b}_{i-1}) \in B^*(A),$$

где $c', b' \in B$. Таким образом,

$$u^{-1}vu \in B^*(A)$$

для любых $u \in A^*$, $v \in B^*(A)$. Следовательно, $B^*(A)$ инвариантна в $A^*$.

*Достаточность.* Пусть $x$ и $b$ произвольные элементы из $A$ и $B$ соответственно. Так как

$$\theta_A(x) \in A^*,\ \theta_A(b) \in B^*(A),$$

то из инвариантности $B^*(A)$ в $A^*$ следует

$$\theta_A(x)\theta_A(b)\theta_A^{-1}(x) \in B^*(A),$$

$$\theta_A(xb\underbrace{x\ldots x}_{n-3}\bar{x}) \in B^*(A),$$

$$[xb\underbrace{x\ldots x}_{n-3}\bar{x}] = b' \in B,$$

откуда

$$[xB\underbrace{x\ldots x}_{n-3}\bar{x}] \subseteq B. \qquad (1)$$

Аналогично из

$$\theta_A^{-1}(x)\theta_A(b)\theta_A(x) \in B^*(A)$$

получаем

$$[\bar{x}\underbrace{x\ldots x}_{n-3}Bx] \subseteq B,$$

откуда

$$[x[\bar{x}\underbrace{x\ldots x}_{n-3}Bx]\underbrace{x\ldots x}_{n-3}\bar{x}] \subseteq [xB\underbrace{x\ldots x}_{n-3}\bar{x}],$$

$$B \subseteq [xB\underbrace{x\ldots x}_{n-3}\bar{x}]. \qquad (2)$$

Из (1) и (2) следует равенство

$$[xB\underbrace{x\ldots x}_{n-3}\overline{x}] = B.$$

По теореме 2.3.9 $<B, [\,]>$ инвариантна в $<A, [\,]>$. ∎

**7.3.2. Предложение.** Если $<B, [\,]>$ – полуинвариантная n-арная подгруппа n-арной группы $<A, [\,]>$, то $B_o(A)$ – инвариантная подгруппа группы $A^*$.

***Доказательство.*** Зафиксируем $b \in B$ и пусть

$$u = \theta_A(a\underbrace{b\ldots b}_{i-1}),\, a \in A,\, i \in \{1, \ldots, n-1\},$$

$$v = \theta_A(c\underbrace{b\ldots b}_{n-2}),\, c \in B,$$

произвольные элементы из $A^*$ и $B_o(A)$ соответственно.

Так как $<B, [\,]>$ полуинвариантна в $<A, [\,]>$, $\overline{b} \in B$ и $\overline{a}\underbrace{a\ldots a}_{n-2}$ – нейтральная последовательность, то

$$u^{-1}vu = \theta_A(\overline{b}\underbrace{b\ldots b}_{n-i-1}\overline{a}\underbrace{a\ldots a}_{n-3})\theta_A(c\underbrace{b\ldots b}_{n-2})\theta_A(a\underbrace{b\ldots b}_{i-1}) =$$

$$= \theta_A(\overline{b}\underbrace{b\ldots b}_{n-i-1}\overline{a}\underbrace{a\ldots a}_{n-3}[c\underbrace{b\ldots b}_{n-2}a]\underbrace{b\ldots b}_{i-1}) =$$

$$= \theta_A(\overline{b}\underbrace{b\ldots b}_{n-i-1}\overline{a}\underbrace{a\ldots a}_{n-3}[ab_1\ldots b_{n-1}]\underbrace{b\ldots b}_{i-1}) =$$

$$= \theta_A(\overline{b}\underbrace{b\ldots b}_{n-i-1}b b_1\ldots b_{n-1}\underbrace{b\ldots b}_{i-1}) =$$

$$= \theta_A([\overline{b}\underbrace{b\ldots b}_{n-i-1}b b_1\ldots b_i]b_{i+1}\ldots b_{n-1}\underbrace{b\ldots b}_{i-1}) =$$

$$= \theta_A(b' b_{i+1}\ldots b_{n-1}\underbrace{b\ldots b}_{i-1}) \in B_o(A),$$

где $b', b_1 \ldots b_{n-1} \in B$. Таким образом,

$$u^{-1}vu \in B_o(A)$$

для любых $u \in A^*$, $v \in B_o(A)$. Следовательно, $B_o(A)$ инвариантна в $A^*$. ∎

**7.3.3. Теорема.** n-Арная подгруппа $< B, [\ ] >$ n-арной группы $< A, [\ ] >$ полуинвариантна в ней тогда и только тогда, когда подгруппа $B_o(A)$ инвариантна в группе $A_o$.

***Доказательство.*** *Необходимость.* Следует из предложения 7.3.2.

*Достаточность.* Пусть x и c – произвольные элементы из A и B соответственно, b – фиксированный из B. Так как

$$\theta_A(x\underbrace{b\ldots b}_{n-2}) \in A_o,\ \theta_A(c\underbrace{b\ldots b}_{n-2}) \in B_o(A)$$

и $B_o(A)$ инвариантна в $A_o$, то

$$\theta_A(x\underbrace{b\ldots b}_{n-2})\theta_A(c\underbrace{b\ldots b}_{n-2})\theta_A^{-1}(x\underbrace{b\ldots b}_{n-2}) \in B_o(A),$$

$$\theta_A(x\underbrace{b\ldots b}_{n-2})\theta_A(c\underbrace{b\ldots b}_{n-2})\theta_A(\overline{b}\,\overline{x}\underbrace{x\ldots x}_{n-3}) \in B_o(A),$$

$$\theta_A(x\underbrace{b\ldots b}_{n-2}c\underbrace{b\ldots b}_{n-2}\overline{b}\,\overline{x}\underbrace{x\ldots x}_{n-3}) \in B_o(A),$$

$$x\underbrace{b\ldots b}_{n-2}c\underbrace{b\ldots b}_{n-2}\overline{b}\,\overline{x}\underbrace{x\ldots x}_{n-3}\theta_A\ b_1\ \ldots\ b_{n-1},$$

где $b_1, \ldots, b_{n-1} \in B$. Тогда

$$[x\underbrace{b\ldots b}_{n-2}c\underbrace{b\ldots b}_{n-2}\overline{b}\,\overline{x}\underbrace{x\ldots x}_{n-3}x] = [b_1\ \ldots\ b_{n-1}x],$$

$$[x\underbrace{b\ldots b}_{n-2}c\underbrace{b\ldots b}_{n-2}\overline{b}] \in [\underbrace{B\ldots B}_{n-1}x],$$

откуда, учитывая произвольный выбор c ∈ B и то, что b – фиксированный из B, получим

$$[x\underbrace{B\ldots B}_{n-1}] \subseteq [\underbrace{B\ldots B}_{n-1}x]. \qquad (1)$$

Аналогично, так как

$$\theta_A(\underbrace{b\ldots b}_{n-2}x) \in A_o,\ \theta_A(\underbrace{b\ldots b}_{n-2}c) \in B_o(A)$$

и $B_o(A)$ инвариантна в $A_o$, то

$$\theta_A^{-1}(\underbrace{b\ldots b}_{n-2}x)\theta_A(\underbrace{b\ldots b}_{n-2}c)\theta_A(\underbrace{b\ldots b}_{n-2}x) \in B_o(A),$$

$$\theta_A(\overline{x}\underbrace{x\ldots x}_{n-3}\overline{b}\underbrace{b\ldots b}_{n-2}c\underbrace{b\ldots b}_{n-2}x) \in B_o(A),$$

$$[\underbrace{B\ldots B}_{n-1}x] \subseteq [x\underbrace{B\ldots B}_{n-1}]. \qquad (2)$$

Из (1) и (2) следует

$$[x\underbrace{B\ldots B}_{n-1}] = [\underbrace{B\ldots B}_{n-1}x].$$

Следовательно, < B, [ ] > полуинвариантна в < A, [ ] >. ∎

## §7.4. СВЯЗЬ МЕЖДУ n-АРНЫМИ АНАЛОГАМИ И ИХ БИНАРНЫМИ ПРОТОТИПАМИ

**7.4.1. Лемма** [115]. Пусть < B, [ ] > – n-арная подгруппа n-арной группы < A, [ ] >, b ∈ B,

$$u = \theta_A(x\underbrace{b\ldots b}_{n-2}) \in N_{A_o}(B_o(A)).$$

Тогда $x \in HN_A(B)$.

*Доказательство.* По условию

$$u^{-1}vu \in B_o(A)$$

для любого

$$v = \theta_A(b_o\underbrace{b \ldots b}_{n-2}) \in B_o(A),$$

то есть

$$\theta_A(\overline{b}\,\overline{x}\underbrace{x \ldots x}_{n-3})\theta_A(b_o\underbrace{b \ldots b}_{n-2})\theta_A(x\underbrace{b \ldots b}_{n-2}) \in B_o(A),$$

откуда

$$\overline{b}\,\overline{x}\underbrace{x \ldots x}_{n-3}b_o\underbrace{b \ldots b}_{n-2}x\underbrace{b \ldots b}_{n-2}\theta_A b_1 \ldots b_{n-1}$$

для некоторых $b_1, \ldots, b_{n-1} \in B$. Тогда

$$[b_o\underbrace{b \ldots b}_{n-2}x] = [x\underbrace{b \ldots b}_{n-2}b_1 \ldots b_{n-1}\overline{b}] \in [x\underbrace{B \ldots B}_{n-1}].$$

Так как $b_o$ выбран в B произвольно, то доказано включение

$$[\underbrace{B \ldots B}_{n-1}x] \subseteq [x\underbrace{B \ldots B}_{n-1}]. \qquad (1)$$

Снова применяя условие, получаем

$$uvu^{-1} \in B_o(A),$$

то есть

$$\theta_A(x\underbrace{b \ldots b}_{n-2})\theta_A(b_o\underbrace{b \ldots b}_{n-2})\theta_A(\overline{b}\,\overline{x}\underbrace{x \ldots x}_{n-3}) \in B_o(A),$$

откуда

$$x\underbrace{b\ldots b}_{n-2}b_o\underbrace{b\ldots b}_{n-2}\bar{b}\,\bar{x}\underbrace{x\ldots x}_{n-3}\theta_A c_1\ldots c_{n-1}$$

для некоторых $c_1 \ldots c_{n-1} \in B$. Тогда

$$[x\underbrace{b\ldots b}_{n-2}b_o] = [c_1\ldots c_{n-2}x] \in [\underbrace{B\ldots B}_{n-1}x],$$

откуда

$$[x\underbrace{B\ldots B}_{n-1}] \subseteq [\underbrace{B\ldots B}_{n-1}x]. \qquad (2)$$

Из (1) и (2) следует

$$[x\underbrace{B\ldots B}_{n-1}] = [\underbrace{B\ldots B}_{n-1}x].$$

Следовательно, $x \in HN_A(B)$. ∎

**7.4.2. Теорема** [115]. Если $< B, [\ ] > -n$-арная подгруппа $n$-арной группы $< A, [\ ] >$, то

$$(HN_A(B))_o(A) = N_{A_o}(B_o(A)).$$

*Доказательство.* Зафиксируем $h \in HN_A(B)$ и выберем произвольный

$$u = \theta_A(h_o\underbrace{h\ldots h}_{n-2}) \in (HN_A(B))_o(A),\ h_o \in HN_A(B).$$

Если $b_o$ – произвольный, $b$ – фиксированный элементы из $B$, то

$$v = \theta_A(b_o\underbrace{b\ldots b}_{n-2})$$

– произвольный элемент из $B_o(A)$. Так как $h_o, \bar{h} \in HN_A(B)$, то

$$u^{-1}vu = \theta_A(\bar{h}\,\bar{h}_o\underbrace{h_o\ldots h_o}_{n-3})\theta_A(b_o\underbrace{b\ldots b}_{n-2})\theta_A(h_o\underbrace{h\ldots h}_{n-2}) =$$

$$= \theta_A(\bar{h}\,\bar{h}_o\underbrace{h_o\ldots h_o}_{n-3}\,b_o\underbrace{b\ldots b}_{n-2}\,h_o\underbrace{h\ldots h}_{n-2}) =$$

$$= \theta_A(\bar{h}\,\bar{h}_o\underbrace{h_o\ldots h_o}_{n-3}\,[b_o\underbrace{b\ldots b}_{n-2}h_o]\underbrace{h\ldots h}_{n-2}) =$$

$$= \theta_A(\bar{h}\,\bar{h}_o\underbrace{h_o\ldots h_o}_{n-3}\,[h_o b_1 \ldots b_{n-1}]\underbrace{h\ldots h}_{n-2}) =$$

$$= \theta_A([\bar{h}\,b_1 \ldots b_{n-1}]\underbrace{h\ldots h}_{n-2}) =$$

$$= \theta_A([b'_1 \ldots b'_{n-1}\,\bar{h}]\underbrace{h\ldots h}_{n-2}) = \theta_A(b'_1 \ldots b'_{n-1}),$$

где $b_1, \ldots, b_{n-1}, b'_1, \ldots, b'_{n-1} \in B$. Следовательно,

$$u^{-1}vu \in B_o(A),$$

откуда $u \in N_{A_o}(B_o(A))$ и доказано включение

$$(HN_A(B))_o(A) \subseteq N_{A_o}(B_o(A)). \qquad (1)$$

Так как любой элемент $u \in N_{A_o}(B_o(A))$ можно представить в виде

$$u = \theta_A(x\underbrace{b\ldots b}_{n-2}),\, b \in B,$$

то по лемме 7.4.1 $x \in HN_A(B) = N$, откуда, учитывая $B \subseteq HN_A(B)$, получаем

$$u = \theta_A(x\underbrace{b\ldots b}_{n-2}) \in (HN_A(B))_o(A).$$

Следовательно,

$$N_{A_o}(B_o(A)) \subseteq (HN_A(B))_o(A). \qquad (2)$$

Из (1) и (2) следует требуемое равенство. ∎

Согласно замечанию 2.2.20 соответствующая группа $N_o$ n-арной подгруппы $< N, [\,] >$ n-арной группы $< A, [\,] >$ изоморфна подгруппе $N_o(A)$ соответствующей группы $A_o$. Поэтому из теоремы 7.4.2 вытекает

**7.4.3. Следствие** [115]. Соответствующая группа полунормализатора $< HN_A(B), [\,] >$ n-арной подгруппы $< B, [\,] >$ в n-арной группе $< A, [\,] >$ изоморфна нормализатору подгруппы $B_o(A)$ в соответствующей группе $A_o$:

$$(HN_A(B))_o \simeq N_{A_o}(B_o(A)).$$

**7.4.4. Теорема** [12]. Полунормализатор n-арной подгруппы $< B, [\,] >$ в n-арной группы $< A, [\,] >$ совпадает с нормализатором подгруппы $< {}_aB = B_a, @ >$ в группе $< A, @ >$ для любого $a \in HN_A(B)$.

*Доказательство.* Так как $a \in HN_A(B)$, то

$${}_aB = [a\underbrace{B\ldots B}_{n-1}] = [\underbrace{B\ldots B}_{n-1}a] = B_a,$$

откуда, используя определения полунормализатора в n-арной группе и нормализатора в группе, получим

$$HN_A(B) = \{x \in A \mid [x\underbrace{B\ldots B}_{n-1}] = [\underbrace{B\ldots B}_{n-1}x]\} =$$

$$= \{x \in A \mid [[x\alpha a]\underbrace{B\ldots B}_{n-1}] = [\underbrace{B\ldots B}_{n-1}[a\alpha x]]\} =$$

$$= \{x \in A \mid [x\alpha[a\underbrace{B\ldots B}_{n-1}]] = [[\underbrace{B\ldots B}_{n-1}a]\alpha x]\} =$$

$$= \{x \in A \mid x @ {}_aB = B_a @ x\} = \{x \in A \mid x @ {}_aB = {}_aB @ x\} =$$

$$= N_{<A,@>}(<{}_aB, @>),$$

то есть

$$HN_A(B) = N_{<A,@>}(<{_a}B, @>).\qquad\blacksquare$$

**7.4.5. Следствие** [12]. Полунормализатор n-арной подгруппы $<B, [\ ]>$ в n-арной группе $<A, [\ ]>$ совпадает с нормализатором подгруппы $<B, @>$ в группе $<A, @>$ для любого $a \in B$.

**7.4.6. Предложение** [12]. Если $<B, [\ ]>$ – n-арная подгруппа n-арной группы $<A, [\ ]>$, то

$$N_A(B) = \{x \in A \mid [xB\underbrace{x\ldots x}_{n-3}\overline{x}] = B\} =$$

$$= \{x \in A \mid [\overline{x}\underbrace{x\ldots x}_{n-3}Bx] = B\}.$$

***Доказательство.*** Если $x \in N_A(B)$, то по следствию 7.2.2 $\overline{x} \in N_A(B)$. Поэтому

$$[xB\underbrace{x\ldots x}_{n-3}\overline{x}] = [x[\underbrace{B\ldots B}_{n}]\underbrace{x\ldots x}_{n-3}\overline{x}] =$$

$$= [xB[\underbrace{B\ldots B}_{n-1}x]\underbrace{x\ldots x}_{n-4}\overline{x}] = [xB[x\underbrace{B\ldots B}_{n-1}]\underbrace{x\ldots x}_{n-4}\overline{x}] =$$

$$= [x[Bx\underbrace{B\ldots B}_{n-2}]B\underbrace{x\ldots x}_{n-4}\overline{x}] = [x[x\underbrace{B\ldots B}_{n-1}]B\underbrace{x\ldots x}_{n-4}\overline{x}] =$$

$$= [xx[\underbrace{B\ldots B}_{n}]\underbrace{x\ldots x}_{n-4}\overline{x}] = [xxB\underbrace{x\ldots x}_{n-4}\overline{x}] = \ldots$$

$$\ldots[\underbrace{x\ldots x}_{n-2}B\overline{x}] = [\underbrace{x\ldots x}_{n-2}[\underbrace{B\ldots B}_{n}]\overline{x}] =$$

$$= [\underbrace{x\ldots x}_{n-2}B[\underbrace{B\ldots B}_{n-1}\overline{x}]] = [\underbrace{x\ldots x}_{n-2}B[\overline{x}\underbrace{B\ldots B}_{n-1}]] =$$

$$= [\underbrace{x \ldots x}_{n-2}[B\overline{x}\underbrace{B \ldots B}_{n-1}]B] = [\underbrace{x \ldots x}_{n-2}[\overline{x}\underbrace{B \ldots B}_{n-1}]B] =$$

$$= [\underbrace{x \ldots x}_{n-2}\overline{x}[\underbrace{B \ldots B}_{n}]] = B,$$

то есть

$$[xB\underbrace{x \ldots x}_{n-3}\overline{x}] = B.$$

Следовательно,

$$N_A(B) \subseteq \{x \in A \mid [xB\underbrace{x \ldots x}_{n-3}\overline{x}] = B\}. \qquad (1)$$

Если теперь

$$[xB\underbrace{x \ldots x}_{n-3}\overline{x}] = B \qquad (2)$$

для некоторого $x \in A$, то

$$[[xB\underbrace{x \ldots x}_{n-3}\overline{x}]x\underbrace{B \ldots B}_{n-2}] = [Bx\underbrace{B \ldots B}_{n-2}],$$

$$[x\underbrace{B \ldots B}_{n-1}] = [Bx\underbrace{B \ldots B}_{n-2}]. \qquad (3)$$

Из (2) следует также

$$[\underbrace{\overline{x}x \ldots x}_{n-3}[xB\underbrace{x \ldots x}_{n-3}\overline{x}]x] = [\underbrace{\overline{x}x \ldots x}_{n-3}Bx],$$

$$B = [\underbrace{\overline{x}x \ldots x}_{n-3}Bx],$$

откуда

$$[Bx\underbrace{B...B}_{n-2}] = [Bx[\underbrace{\overline{x}x...xBx}_{n-3}]...[\underbrace{\overline{x}x...xBx}_{n-3}]] =$$

$$\underbrace{\phantom{[Bx[\overline{x}x...xBx]...[\overline{x}x...xBx]]}}_{n-2}$$

$$= [B\,x\underbrace{\overline{x}x...xB}_{n-3}...\underbrace{x\overline{x}x...xB}_{n-3}x] = [\underbrace{B...B}_{n-1}x],$$

$$\underbrace{\phantom{[Bx\overline{x}x...xB...x\overline{x}x...xBx]}}_{n-2}$$

то есть

$$[Bx\underbrace{B...B}_{n-2}] = [\underbrace{B...B}_{n-1}x]. \qquad (4)$$

Из (3) и (4), ввиду следствия 7.1.7, следует $x \in N_A(B)$, то есть

$$\{x \in A \mid [xB\underbrace{x...x}_{n-3}\overline{x}] = B\} \subseteq N_A(B). \qquad (5)$$

Из (1) и (5) следует

$$N_A(B) = \{x \in A \mid [xB\underbrace{x...x}_{n-3}\overline{x}] = B\}.$$

Равенство

$$N_A(B) = \{x \in A \mid [\overline{x}\underbrace{x...x}_{n-3}Bx] = B\}$$

доказывается аналогично. ∎

**7.4.7. Следствие.** Если $< B, [\,] > -$ тернарная подгруппа тернарной группы $< A, [\,] >$, то

$$N_A(B) = \{x \in A \mid [xB\overline{x}] = B\} = \{x \in A \mid [\overline{x}Bx] = B\}.$$

Ясно, что n-арная подгруппа $< B, [\,] >$ m-полуинвариантна в n-арной группе $< A, [\,] >$ тогда и только тогда, когда $N_A(B, m) = A$. В частности $< B, [\,] >$ инвариантна (полуинвариантна) в $< A, [\,] >$ тогда и только тогда, когда

$$N_A(B) = A \quad (HN_A(B) = A).$$

Для подмножества B n-арной группы < A, [ ] >, в отличие от n-арной подгруппы < B, [ ] >, множества

$$N_A(B, 2) = \{x \in A \mid [x \overset{n-1}{B}] = [B x \overset{n-1-i}{B}], \forall i = 1, \ldots, n-1\}$$

и

$$N = \{x \in A \mid [xB\underbrace{x \ldots x}_{n-2} \bar{x}] = B\}$$

могут не совпадать. Однако, имеет место включение $N \subseteq N_A(B, 2)$ [12, 116].

**7.4.8. Лемма** [115]. Если $x \in N_A(B)$, то

$$[\underbrace{x \ldots x}_{i-1} B \underbrace{x \ldots x}_{n-i-1} \bar{x}] = B,$$

$$[\bar{x}\underbrace{x \ldots x}_{n-i-1} B \underbrace{x \ldots x}_{i-1}] = B$$

для любого $i = 1, \ldots, n-1$.

*Доказательство.* Докажем второе равенство. Если $i = 1$, то $B = B$. Если $i = 2$, то по предложению 7.4.6

$$[\bar{x}\underbrace{x \ldots x}_{n-3} Bx] = B.$$

Из последнего равенства имеем

$$[\bar{x}\underbrace{x \ldots x}_{n-3}[\bar{x}\underbrace{x \ldots x}_{n-3} Bx]x] = [\bar{x}\underbrace{x \ldots x}_{n-3} Bx]$$

$$[\bar{x}\underbrace{x \ldots x}_{n-4} Bxx] = B,$$

откуда

$$[\bar{x}\underbrace{x \ldots x}_{n-3}[\bar{x}\underbrace{x \ldots x}_{n-4} Bxx]x] = [\bar{x}\underbrace{x \ldots x}_{n-3} Bx],$$

$$[\,\overline{x}\underbrace{x\,...\,x}_{n-5}Bxxx] = B.$$

Продолжая, получим

$$[\,\overline{x}xB\underbrace{x\,...\,x}_{n-3}\,] = B,$$

$$[\,\overline{x}B\underbrace{x\,...\,x}_{n-2}\,] = B.$$

Таким образом, второе равенство верно для любого $i = 1, ..., n - 1$.

Первое равенство доказывается аналогично. ∎

**7.4.9. Теорема** [115]. Если $< B, [\ ] >$ – n-арная подгруппа n-арной группы $< A, [\ ] >$, то

$$(N_A(B))^*(A) = N_{A^*}(B^*(A)).$$

*Доказательство.* Зафиксируем $h \in N_A(B)$ и выберем произвольный

$$u = \theta_A(h_o\underbrace{h\,...\,h}_{i-1}) \in (N_A(B))^*(A),\ h_o \in N_A(B).$$

Если $b_o$ – произвольный, $b$ – фиксированный элементы из $B$, то

$$v = \theta_A(b_o\underbrace{b\,...\,b}_{j-1})$$

– произвольный элемент из $B^*(A)$. Так как $h_o, \overline{h} \in N_A(B)$, то

$$u^{-1}vu = \theta_A(\overline{h}\underbrace{h\,...\,h}_{n-i-1}\overline{h}_o\underbrace{h_o\,...\,h_o}_{n-3})\theta_A(b_o\underbrace{b\,...\,b}_{j-1})\theta_A(h_o\underbrace{h\,...\,h}_{i-1}) =$$

$$= \theta_A(\overline{h}\underbrace{h\,...\,h}_{n-i-1}\overline{h}_o\underbrace{h_o\,...\,h_o}_{n-3}b_o\underbrace{b\,...\,b}_{j-1}h_o\underbrace{h\,...\,h}_{i-1}) =$$

$$= \theta_A(\overline{h}\underbrace{h\ldots h}_{n-i-1}[\overline{h}_o\underbrace{h_o\ldots h_o}_{n-3}b_oh_o]$$

$$\underbrace{[\overline{h}_o\underbrace{h_o\ldots h_o}_{n-3}bh_o]\ldots[\overline{h}_o\underbrace{h_o\ldots h_o}_{n-3}bh_o]}_{j-1}\underbrace{h\ldots h}_{i-1}) =$$

$$= \theta_A(\overline{h}\underbrace{h\ldots h}_{n-i-1}b'_o\underbrace{b'\ldots b'}_{j-1}\underbrace{h\ldots h}_{i-1}) =$$

$$= \theta_A([\overline{h}\underbrace{h\ldots h}_{n-i-1}b'_o\underbrace{h\ldots h}_{i-1}]\underbrace{[\overline{h}\underbrace{h\ldots h}_{n-i-1}b'h\underbrace{\ldots h}_{i-1}]\ldots[\overline{h}\underbrace{h\ldots h}_{n-i-1}b'h\underbrace{\ldots h}_{i-1}]}_{j-1}) =$$

$$= \theta_A(b''_o\underbrace{b''\ldots b''}_{j-1}),$$

где $b'_o$, $b'$, $b''_o$, $b'' \in B$. Следовательно,

$$u^{-1}vu \in B^*(A),$$

откуда $u \in N_{A*}(B^*(A))$ и доказано включение

$$(N_A(B))^*(A) \subseteq N_{A*}(B^*(A)). \qquad (1)$$

Пусть теперь $c \in B$

$$u = \theta_A(x\underbrace{c\ldots c}_{i-1}) = \theta_A(x)\theta_A(\underbrace{c\ldots c}_{i-1})$$

произвольный элемент из $N_{A*}(B^*(A))$. Так как

$$\theta_A(\underbrace{c\ldots c}_{i-1}) \in B^*(A) \subseteq N_{A*}(B^*(A)),$$

то из последнего равенства следует

$$\theta_A(x) \in N_{A*}(B^*(A)). \qquad (2)$$

Тогда

$$\theta_A^{-1}(x)\theta_A(b)\theta_A(x) \in B^*(A)$$

для любого $b \in B$, откуда

$$\theta_A(\overline{x}\underbrace{x \ldots x}_{n-3}bx) \in B^*(A),$$

$$[\overline{x}\underbrace{x \ldots x}_{n-3}bx] = b'$$

для некоторого $b' \in B$. Так как $b$ выбран в $B$ произвольно, то

$$[\overline{x}\underbrace{x \ldots x}_{n-3}Bx] \subseteq B. \qquad (3)$$

Из (2) также следует

$$\theta_A(x)\theta_A(b)\theta_A^{-1}(x) \in B^*(A)$$

для любого $b \in B$, откуда

$$[xB\underbrace{x \ldots x}_{n-3}\overline{x}] \subseteq B.$$

Из последнего равенства вытекает

$$B \subseteq [\overline{x}\underbrace{x \ldots x}_{n-3}Bx]. \qquad (4)$$

Из (3) и (4) следует

$$[\overline{x}\underbrace{x \ldots x}_{n-3}Bx] = B,$$

откуда, ввиду предложения 7.4.6, $x \in N_A(B)$. Тогда

$$u = \theta_A(x\underbrace{c \ldots c}_{i-1}) \in (N_A(B))^*(A),$$

откуда

$$N_{A^*}(B^*(A)) \subseteq (N_A(B))^*(A). \qquad (5)$$

Из (1) и (5) следует требуемое равенство. ∎

Согласно теореме 2.2.19 универсальная обертывающая группа Поста N* n-арной подгруппы < N, [ ] > n-арной группы < A, [ ] > изоморфна подгруппе N*(A) универсальной обертывающей группы Поста A*. Поэтому из теоремы 7.4.9 вытекает

**7.4.10 Следствие** [115]. Универсальная обертывающая группа Поста нормализатора < $N_A(B)$, [ ] > n-арной подгруппы < B, [ ] > в n-арной группе < A, [ ] > изоморфна нормализатору подгруппы B*(A) в универсальной обертывающей группе Поста A*:

$$(N_A(B))^* \simeq N_{A^*}(B^*(A)).$$

## ДОПОЛНЕНИЯ И КОММЕНТАРИИ

**1.** При доказательстве теоремы 7.2.1 существенно использовалось то, что < B, [ ] > – n-арная подгруппа n-арной группы < A, [ ] >. В то же время, следствие 7.2.2 получено С.А. Русаковым для более общего случая, когда B – произвольное подмножество n-арной группы. В связи с этим возникает

**Вопрос.** Можно ли в теореме 7.2.1, а также в полученной с ее помощью теореме 7.2.3, n-арную подгруппу < B, [ ] > n-арной группы < A, [ ] > заменить любым подмножеством этой n-арной группы.

**2.** Аналоги теоремы 7.3.3 есть у Поста [3], Русакова [4], Целакоского и Илич [117]. Обратное утверждение к предложению 7.3.2 следует из теоремы 7.3.3.

**3.** В [26] установлено, что если < H, [ ] > и < N, [ ] > – тернарные подгруппы тернарной группы отражений < $B_n$, [ ] >, причём N ⊂ H, то < N, [ ] > – инвариантна в < H, [ ] > тогда и только тогда, когда |H : N| = 2. Отсюда следует, что всякая тернарная подгруппа из < $B_n$, [ ] > либо совпадает со своим нормализатором в < $B_n$, [ ] >, либо

имеет в нём индекс 2. В частности, если n – нечётное, то все тернарные подгруппы из < $B_n$, [ ] > совпадают со своими нормализаторами в < $B_n$, [ ] >, то есть в этом случае в < $B_n$, [ ] > нет инвариантных тернарных подгрупп.

# Г Л А В А 8

# n-АРНЫЕ АНАЛОГИ ЦЕНТРА ГРУППЫ

В данной главе определяются различные n-арные аналоги центра группе и изучаются их свойства, а также их связь со своими бинарными прототипами.

## §8.1. ОПРЕДЕЛЕНИЯ

Напомним, что центром n-арной группы < A, [ ] > называется [3, 4] множество

$$Z(A) = \{z \in A \mid zx \sim xz \text{ для всех } x \in A\},$$

а централизатором подмножества $B \subseteq A$ в n-арной группе < A, [ ] > называется [118] множество

$$C_A(B) = \{z \in A \mid zx \sim xz \text{ для всех } x \in B\}.$$

Ясно, что центр n-арной группы < A, [ ] > совпадает с централизатором подмножества A в n-арной группе < A, [ ] >, то есть

$$Z(A) = C_A(A).$$

Центр n-арной группы является частным случаем более общего понятия m-полуцентра.

**8.1.1. Определение** [119, 120]. *m-Полуцентром* n-арной группы < A, [ ] >, где n = k(m – 1) + 1, k ≥ 1, называется множество

$$Z(A, m) = \{z \in A \mid z x_1^{m-2} x \sim x x_1^{m-2} z, \forall x, x_1, \ldots, x_{m-2} \in A\}.$$

Если m = 2, то последовательность $x_1^{m-2}$ – пустая, а определение 8.1.1 m-полуцентра совпадает с определением центра, то есть

$$Z(A) = Z(A, 2).$$

Если m = n, то n-полуцентр называется *полуцентром* и обозначается символом HZ(A), то есть

$$HZ(A) = Z(A, n).$$

Ясно, что

$$HZ(A) = \{z \in A \mid [z\,x_1^{n-2}\,x] = [x\,x_1^{n-2}\,z],\ \forall x, x_1, \ldots, x_{n-2} \in A\}.$$

Известно, что все единицы n-арной группы лежат в ее центре. Покажем, что полуцентр n-арной группы, в отличие от её центра, может не содержать единицы этой n-арной группы.

**8.1.2. Пример.** Пусть < A, [ ] > – n-арная группа, производная от неабелевой группы A с единицей e и пусть x ∈ A и x не лежит в центре группы A. Тогда существует y ∈ A такой, что yx ≠ xy, откуда

$$yx = e\underbrace{e \ldots e}_{n-3}yx = [e\underbrace{e \ldots e}_{n-3}yx] \neq [x\underbrace{e \ldots e}_{n-3}ye] = x\underbrace{e \ldots e}_{n-3}ye = xy.$$

Следовательно, элемент e, являющийся единицей в < A, [ ] >, не лежит в её полуцентре.

Еще одним n-арным аналогом центра группы является понятие m-полуцентрализатора.

**8.1.3. Определение** [119, 120]. *m-Полуцентрализатором* подмножества B в n-арной группе < A, [ ] >, n = k(m − 1) + 1, k ≥ 1, называется множество

$$C_A(B, m) = \{z \in A \mid z\,x_1^{m-2}\,x \sim x\,x_1^{m-2}\,z,\ \forall x, x_1, \ldots, x_{m-2} \in B\}.$$

Ясно, что m-полуцентр n-арной группы < A, [ ] > совпадает с m-полуцентрализатором подмножества A в n-арной группе < A, [ ] >, то есть

$$Z(A, m) = C_A(A, m).$$

Если m = 2, то определение 2-полуцентрализатора совпадает с определением централизатора, то есть

$$C_A(B) = C_A(B, 2).$$

Если m = n, то n-полуцентрализатор называется *полуцентрализатором* и обозначается символом $HC_A(B)$, то есть

$$HC_A(B) = C_A(B, n).$$

Ясно, что

$$HC_A(B) = \{z \in A \mid [z\,x_1^{n-2}\,x] = [x\,x_1^{n-2}\,z],\ \forall x, x_1, \ldots, x_{n-2} \in B\}.$$

Таким образом,

$$Z(A) = C_A(A, 2),$$
$$C_A(B) = C_A(B, 2),$$
$$Z(A, m) = C_A(A, m),$$
$$HZ(A) = C_A(A, n),$$
$$HC_A(B) = C_A(B, n),$$

то есть понятие m-полуцентрализатора включает в себя все другие приведенные выше n-арные аналоги центра группы и поэтому является самым широким из них.

Легко проверяется, что

$$Z(A) = \{z \in A \mid [xz\bar{x}\underbrace{x \ldots x}_{n-3}] = z \text{ для всех } x \in A\} =$$

$$= \{z \in A \mid [\bar{x}\underbrace{x \ldots x}_{n-3}zx] = z \text{ для всех } x \in A\},$$

$$C_A(B) = \{z \in A \mid [xz\bar{x}\underbrace{x \ldots x}_{n-3}] = z \text{ для всех } x \in B\} =$$

$$= \{z \in A \mid [\bar{x}\underbrace{x \ldots x}_{n-3}zx] = z \text{ для всех } x \in B\}.$$

Следующее предложение является следствием определений.

**8.1.4. Предложение.** n-Арная группа $< A, [\,] >$ является m-полуабелевой тогда и только тогда, когда она совпадает со своим m-полуцентром. В частности, n-арная группа является абелевой (полуабелевой) тогда и только тогда, когда она совпадает со своим центром (полуцентром).

## §8.2. ОСНОВНАЯ ТЕОРЕМА

**8.2.1. Лемма** [119, 120]. Пусть $< A, [\,] >$ – n-арная группа,

$$n = k(m-1) + 1, k \geq 1.$$

Тогда:

1) если $u_1, \ldots, u_m \in C_A(B, m)$, $x_1, \ldots, x_{m-2}, x \in B$, то

$$u_1 \ldots u_m x_1 \ldots x_{m-2} x \sim x x_1 \ldots x_{m-2} u_1 \ldots u_m;$$

2) если $s \geq 1$, $u, v \in C_A(B, m) \cap B$, $x_1, \ldots, x_{s(m-1)-1} \in B$, то

$$u x_1 \ldots x_{s(m-1)-1} v \sim v x_1 \ldots x_{s(m-1)-1} u.$$

*Доказательство.* 1) Используя определение m-полуцентрализатора, получим

$$u_1 \ldots u_m x_1 \ldots x_{m-2} x \sim u_1 \ldots u_{m-1} x x_1 \ldots x_{m-3} x_{m-2} u_m \sim$$

$$\sim u_1 \ldots u_{m-2} x_{m-2} x x_1 \ldots x_{m-4} x_{m-3} u_{m-1} u_m \sim$$

$$\sim u_1 \ldots u_{m-3} x_{m-3} x_{m-2} x x_1 \ldots x_{m-4} u_{m-2} u_{m-1} u_m \sim \ldots$$

$$\ldots \sim u_1 u_2 x_2 \ldots x_{m-2} x x_1 u_3 \ldots u_m \sim$$

$$\sim u_1 x_1 x_2 \ldots x_{m-2} x u_2 u_3 \ldots u_m \sim x x_1 \ldots x_{m-2} u_1 u_2 \ldots u_m,$$

то есть

$$u_1 \ldots u_m x_1 \ldots x_{m-2} x \sim x x_1 \ldots x_{m-2} u_1 \ldots u_m.$$

2) Положим

$$\alpha_1 = x_1 \ldots x_{m-2},\ \alpha_2 = x_m \ldots x_{2(m-1)-1},\ \ldots$$

$$\ldots,\ \alpha_s = x_{(s-1)(m-1)+1} \ldots x_{s(m-1)-1}.$$

Тогда, учитывая $u \in C_A(B, m) \cap B$, получим

$$u x_1 \ldots x_{s(m-1)-1} v = u \alpha_1 x_{m-1} \alpha_2 x_{2(m-1)} \ldots x_{(s-1)(m-1)} \alpha_s v \sim$$

$$\sim x_{m-1} \alpha_1 u \alpha_2 x_{2(m-1)} \ldots x_{(s-1)(m-1)} \alpha_s v \sim$$

$$\sim x_{m-1} \alpha_1 x_{2(m-1)} \alpha_2 u \ldots x_{(s-1)(m-1)} \alpha_s v \sim \ldots$$

$$\ldots \sim x_{m-1} \alpha_1 x_{2(m-1)} \alpha_2 x_{3(m-1)} \ldots u \alpha_s v \sim$$

$$\sim x_{m-1} \alpha_1 x_{2(m-1)} \alpha_2 x_{3(m-1)} \ldots x_{(s-1)(m-1)} \alpha_{s-1} v \alpha_s u \sim$$

$$\sim x_{m-1} \alpha_1 x_{2(m-1)} \alpha_2 x_{3(m-1)} \ldots v \alpha_{s-1} x_{(s-1)(m-1)} \alpha_s u \sim \ldots$$

$$\ldots \sim x_{m-1} \alpha_1 v \alpha_2 x_{2(m-1)} \ldots \alpha_{s-1} x_{(s-1)(m-1)} \alpha_s u \sim$$

$$\sim v \alpha_1 x_{m-1} \alpha_2 x_{2(m-1)} \ldots x_{(s-1)(m-1)} \alpha_s u = v x_1 \ldots x_{s(m-1)-1} u,$$

то есть

$$u x_1 \ldots x_{s(m-1)-1} v \sim v x_1 \ldots x_{s(m-1)-1} u. \qquad \blacksquare$$

**8.2.2. Теорема** [119, 120]. Пусть $< A, [\ ] > -$ n-арная группа и пусть $n = k(m-1)+1$, $k \geq 1$, $B \subseteq A$. Тогда, если $C_A(B, m) \neq \varnothing$, то $< C_A(B, m), [\ ] > -$ n-арная подгруппа в $< A, [\ ] >$, лежащая в её m-полунормализаторе $N_A(B, m)$.

*Доказательство.* Пусть $z_1, \ldots, z_{k(m-1)+1} -$ произвольные

элементы из $C_A(B, m)$, $x_1, \ldots, x_{m-2}, x$ – произвольные элементы из $B$, и положим

$$\alpha_1 = z_1 \ldots z_{m-1};$$

$$\alpha_2 = z_m \ldots z_{2(m-1)} = z_m \alpha_2', \ \alpha_2' = z_{m+1} \ldots z_{2(m-1)};$$

$$\alpha_3 = z_{2(m-1)+1} \ldots z_{3(m-1)} = z_{2(m-1)+1} \alpha_3', \ \alpha_3' = z_{2(m-1)+2} \ldots z_{3(m-1)};$$

$$\ldots\ldots\ldots\ldots\ldots\ldots\ldots\ldots\ldots\ldots\ldots\ldots\ldots\ldots\ldots\ldots\ldots\ldots\ldots\ldots$$

$$\alpha_{k-2} = z_{(k-3)(m-1)+1} \ldots z_{(k-2)(m-1)} = z_{(k-3)(m-1)+1} \alpha_{k-2}',$$

$$\alpha_{k-2}' = z_{(k-3)(m-1)+2} \ldots z_{(k-2)(m-1)};$$

$$\alpha_{k-1} = z_{(k-2)(m-1)+1} \ldots z_{(k-1)(m-1)} = z_{(k-2)(m-1)+1} \alpha_{k-1}',$$

$$\alpha_{k-1}' = z_{(k-2)(m-1)+2} \ldots z_{(k-1)(m-1)};$$

$$\alpha_k = z_{(k-1)(m-1)+1} \ldots z_{k(m-1)+1} = z_{(k-1)(m-1)+1} \alpha_k',$$

$$\alpha_k' = z_{(k-1)(m-1)+2} \ldots z_{k(m-1)+1}.$$

Используя определение m-полуцентрализатора и лемму 8.2.1, получим

$$[z_1 \ldots z_n]x_1 \ldots x_{m-2}x \sim \alpha_1 \ldots \alpha_{k-1}\underbrace{\alpha_k x_1 \ldots x_{m-2} x}_{\text{лемма}} \sim$$

$$\sim \alpha_1 \ldots \alpha_{k-1} x x_1 \ldots x_{m-2} \alpha_k = \alpha_1 \ldots \alpha_{k-1} \underbrace{x x_1 \ldots x_{m-2} z_{(k-1)(m-1)+1}}_{\text{определение}} \alpha_k' \sim$$

$$\sim \alpha_1 \ldots \alpha_{k-2} \underbrace{\alpha_{k-1} z_{(k-1)(m-1)+1} x_1 \ldots x_{m-2} x}_{\text{лемма}} \alpha_k' \sim$$

$$\sim \alpha_1 \ldots \alpha_{k-2} \, x x_1 \ldots x_{m-2} \alpha_{k-1} z_{(k-1)(m-1)+1} \alpha_k' =$$

$$= \alpha_1 \ldots \alpha_{k-2} \underbrace{x x_1 \ldots x_{m-2} z_{(k-2)(m-1)+1}}_{\text{определение}} \alpha_{k-1}' \alpha_k \sim$$

$$\sim \alpha_1 \ldots \alpha_{k-3} \underbrace{\alpha_{k-2} z_{(k-2)(m-1)+1} x_1 \ldots x_{m-2} x}_{\text{лемма}} \alpha'_{k-1} \alpha_k \sim$$

$$\sim \alpha_1 \ldots \alpha_{k-3} x x_1 \ldots x_{m-2} \alpha_{k-2} z_{(k-2)(m-1)+1} \alpha'_{k-1} \alpha_k =$$

$$= \alpha_1 \ldots \alpha_{k-3} x x_1 \ldots x_{m-2} \alpha_{k-2} \alpha_{k-1} \alpha_k \sim \ldots$$

$$\ldots \sim \alpha_1 x x_1 \ldots x_{m-2} \alpha_2 \alpha_3 \ldots \alpha_k =$$

$$= \alpha_1 \underbrace{x x_1 \ldots x_{m-2} z_m}_{\text{определение}} \alpha'_2 \alpha_3 \ldots \alpha_k \sim \underbrace{\alpha_1 z_m x_1 \ldots x_{m-2} x}_{\text{лемма}} \alpha'_2 \alpha_3 \ldots \alpha_k \sim$$

$$\sim x x_1 \ldots x_{m-2} \alpha_1 z_m \alpha'_2 \alpha_3 \ldots \alpha_k \sim x x_1 \ldots x_{m-2} \alpha_1 \alpha_2 \ldots \alpha_k \sim$$

$$\sim x x_1 \ldots x_{m-2} [z_1 \ldots z_n],$$

то есть

$$[z_1 \ldots z_n] x_1 \ldots x_{m-2} x = x x_1 \ldots x_{m-2} [z_1 \ldots z_n].$$

Следовательно,

$$[z_1 \ldots z_n] \in C_A(B, m),$$

и множество $C_A(B, m)$ замкнуто относительно n-арной операции [ ].

Покажем, что для $\bar{z} \in C_A(B, m)$ для любого $z \in C_A(B, m)$.
Так как

$$z y_1 \ldots y_{m-2} y \sim y y_1 \ldots y_{m-2} z$$

для любых $y, y_1, \ldots, y_{m-2} \in B$, то

$$\underbrace{z \ldots z}_{n} \bar{z} z y_1 \ldots y_{m-2} y \bar{z} \sim \underbrace{z \ldots z}_{n} \bar{z} y y_1 \ldots y_{m-2} z \bar{z},$$

$$z z z y_1 \ldots y_{m-2} y \bar{z} \sim \bar{z} [\underbrace{z \ldots z}_{n}] y y_1 \ldots y_{m-3} y_{m-2} z \bar{z}.$$

К левой части последнего соотношения трижды применим определение m-полуцентрализатора, а к правой части – дока-

занный выше факт

$$[\underbrace{z \ldots z}_{n}] \in C_A(B, m):$$

$$zzyy_1 \ldots y_{m-3}y_{m-2}z\bar{z} \sim \bar{z}y_{m-2}yy_1 \ldots y_{m-3}[\underbrace{z \ldots z}_{n}]z\bar{z},$$

$$zy_{m-2}yy_1 \ldots y_{m-4}y_{m-3}zz\bar{z} \sim \bar{z}y_{m-2}yy_1 \ldots y_{m-3}zzz,$$

$$y_{m-3}y_{m-2}yy_1 \ldots y_{m-4}zzz\bar{z} \sim \bar{z}y_{m-2}yy_1 \ldots y_{m-3}zzz,$$

$$y_{m-3}y_{m-2}yy_1 \ldots y_{m-4}\bar{z}zzz \sim \bar{z}y_{m-2}yy_1 \ldots y_{m-4}y_{m-3}zzz,$$

$$y_{m-3}y_{m-2}yy_1 \ldots y_{m-4}\bar{z} \sim \bar{z}y_{m-2}yy_1 \ldots y_{m-4}y_{m-3}.$$

Из последнего соотношения, сделав переобозначения

$$x = y_{m-3},\ x_1 = y_{m-2},\ x_2 = y,\ x_3 = y_1,\ \ldots,\ x_{m-2} = y_{m-4},$$

получим

$$\bar{z}x_1x_2 \ldots x_{m-2}x \sim xx_1x_2 \ldots x_{m-2}\bar{z}$$

для любых $x, x_1, \ldots, x_{m-2} \in B$. Следовательно, $\bar{z} \in C_A(B, m)$, и по критерию Дёрнте $< C_A(B, m), [\ ] > -$ n-арная подгруппа в $< A, [\ ] >$.

Если $z \in C_A(B, m)$, то

$$[zb_1 \ldots b_{n-1}] = [zb_1 \ldots b_{m-2}b_{m-1}b_m \ldots b_n] =$$

$$= [b_{m-1}b_1 \ldots b_{m-2}zb_m \ldots b_{2m-3}b_{2(m-1)}b_{2m-1} \ldots b_n] =$$

$$= [b_{m-1}b_1 \ldots b_{m-2}b_{2(m-1)}b_m \ldots b_{2m-3}zb_{2m-1} \ldots b_n] = \ldots$$

$$\ldots = [b_{m-1}b_1 \ldots b_{m-2}b_{2(m-1)}b_m \ldots b_{2m-3} \ldots$$

$$\ldots b_{k(m-1)}b_{(k-1)(m-1)+1} \ldots b_{k(m-1)-1}z]$$

для любых $b_1, \ldots, b_{n-1} \in B$, откуда

$$[z \overset{n-1}{B}] = [\overset{m-1}{B} z \overset{n-m}{B}] = [\overset{2(m-1)}{B} z \overset{n-1-2(m-1)}{B}] = \ldots$$

$$\ldots = [\overset{(k-1)(m-1)}{B} z \overset{m-1}{B}] = [\overset{n-1}{B} z].$$

Следовательно, $z \in N_A(B, m)$ и $C_A(B, m) \subseteq N_A(B, m)$. ∎

**8.2.3. Следствие.** Пусть $< A, [\ ] > -$ n-арная группа, $B \subseteq A$. Тогда $< C_A(B), [\ ] >$ и $< HC_A(B), [\ ] > -$ n-арные подгруппы в $< A, [\ ] >$, причём

$$C_A(B) \subseteq N_A(B); \ HC_A(B) \subseteq HN_A(B).$$

**8.2.4. Замечание.** Утверждение следствия 8.2.3 о централизаторе $< C_A(B), [\ ] >$ доказано в [118] для случая, когда $< B, [\ ] > -$ n-арная подгруппа в $< A, [\ ] >$.

## §8.3. СВОЙСТВА m-ПОЛУЦЕНТРАЛИЗАТОРА

**8.3.1. Предложение** [119, 120]. Если $n = k(m - 1) + 1$, то m-полуцентр $< Z(A, m), [\ ] >$ n-арной группы $< A, [\ ] >$ является её m-полуабелевой m-полуинвариантной n-арной подгруппой.

*Доказательство.* То, что $< Z(A, m), [\ ] >$ n-арная подгруппа в $< A, [\ ] >$, доказано в теореме 8.2.2, так как $Z(A, m) = C_A(A, m)$.

m-Полуабелевость $< Z(A, m), [\ ] >$ следует из определений m-полуабелевости и m-полуцентра.

Если $x \in A$, $z_1, \ldots, z_{n-1} \in Z(A, m)$, то

$$[xz_1 \ldots z_{n-1}] = [z_{m-1}z_1 \ldots z_{m-2}xz_m \ldots z_n] =$$

$$= [z_{m-1}z_1 \ldots z_{m-2}z_{2(m-1)}z_m \ldots z_{2m-3}xz_{2m-1} \ldots z_n] = \ldots$$

$$\ldots [z_{m-1}z_1 \ldots z_{m-2}z_{2(m-1)}z_m \ldots z_{k(m-1)}z_{(k-1)(m-1)+1} \ldots z_{k(m-1)-1}x],$$

откуда в силу произвольного выбора $z_1, \ldots, z_{n-1} \in Z(A, m)$, следует

$$[x\underbrace{Z(A,m) \ldots Z(A,m)}_{n-1}] =$$

$$= [\underbrace{Z(A,m)\ldots Z(A,m)}_{i(m-1)} x \underbrace{Z(A,m)\ldots Z(A,m)}_{n-1-i(m-1)}]$$

для любого i = 1, …, k. Следовательно, < Z(A, m), [ ] > m-полуинвариантна в < A, [ ] >. ∎

Полагая в предложении 8.3.1 m = 2, получим

**8.3.2. Следствие** [3, 4]. Центр < Z(A), [ ] > n-арной группы < A, [ ] > является её абелевой инвариантной n-арной подгруппой.

Полагая в предложении 8.3.1 m = n, получим

**8.3.3. Следствие.** Полуцентр < HZ(A), [ ] > n-арной группы < A, [ ] > является её полуабелевой полуинвариантной n-арной подгруппой.

**8.3.4. Замечание.** Полуцентр n-арной группы удовлетворяет более сильному, чем полуинвариантность условию, так как является на самом деле нормальной подгруппой в n-арной группе. Определение нормальной подгруппы дано в §5.2.

**8.3.5. Предложение** [119, 120]. Для n-арной подгруппы < B, [ ] > n-арной группы < A, [ ] > централизатор < $C_A(B)$, [ ] > является инвариантной n-арной подгруппой в нормализаторе < $N_A(B)$, [ ] >.

*Доказательство.* По следствию 8.2.3 < $C_A(B)$, [ ] > – n-арная подгруппа в < A, [ ] >, причём $C_A(B) \subseteq N_A(B)$.

Пусть x – произвольный элемент из $N_A(B)$ и пусть

$$u = [\bar{x}\underbrace{x\ldots x}_{n-3}zx] \in [\bar{x}\underbrace{x\ldots x}_{n-3}C_A(B)x],\ z \in C_A(B).$$

Так как x ∈ $N_A(B)$, то

$$B = [xB\bar{x}\underbrace{x\ldots x}_{n-3}],$$

и следовательно,

$$[xb\bar{x}\underbrace{x\ldots x}_{n-3}] = b' \in B$$

для любого $b \in B$. Из $z \in C_A(B)$ следует

$$zb' \sim b'z,$$

откуда

$$z[xb\bar{x}\underbrace{x\ldots x}_{n-3}] \sim [xb\bar{x}\underbrace{x\ldots x}_{n-3}]z,$$

$$\bar{x}\underbrace{x\ldots x}_{n-3}z[xb\bar{x}\underbrace{x\ldots x}_{n-3}]x \sim \bar{x}\underbrace{x\ldots x}_{n-3}[xb\bar{x}\underbrace{x\ldots x}_{n-3}]zx,$$

$$[\bar{x}\underbrace{x\ldots x}_{n-3}zx]b \sim b[\bar{x}\underbrace{x\ldots x}_{n-3}zx].$$

Следовательно,

$$u = [\bar{x}\underbrace{x\ldots x}_{n-3}zx] \in C_A(B),$$

откуда, ввиду произвольного выбора $u$, следует

$$[\bar{x}\underbrace{x\ldots x}_{n-3}C_A(B)x] \subseteq C_A(B). \qquad (1)$$

Если в приведенном выше доказательстве поменять ролями $\bar{x}\underbrace{x\ldots x}_{n-3}$ и $x$, то получим включение

$$[xC_A(B)\bar{x}\underbrace{x\ldots x}_{n-3}] \subseteq C_A(B).$$

откуда

$$[\bar{x}\underbrace{x\ldots x}_{n-3}[xC_A(B)\bar{x}\underbrace{x\ldots x}_{n-3}]x] \subseteq [\bar{x}\underbrace{x\ldots x}_{n-3}C_A(B)x],$$

$$C_A(B) \subseteq [\,\overline{x}\underbrace{x \ldots x}_{n-3} C_A(B) x\,]. \tag{2}$$

Из (1) и (2) следует равенство

$$[\,\overline{x}\underbrace{x \ldots x}_{n-3} C_A(B) x\,] = C_A(B)$$

для любого $x \in N_A(B)$. Следовательно, $< C_A(B), [\,] >$ инвариантна в $< N_A(B), [\,] >$. ∎

Так как нормализатор инвариантной n-арной подгруппы в n-арной группе совпадает с самой n-арной подгруппой, то из предложения 8.3.5 вытекает

**8.3.6. Следствие.** Если $< B, [\,] >$ инвариантная n-арная подгруппа n-арной группы $< A, [\,] >$, то централизатор $< C_A(B), [\,] >$ также инвариантен в $< A, [\,] >$.

**8.3.7. Замечание.** Так как $Z(A) = C_A(A)$, то следствие 8.3.2 вытекает из следствия 8.3.6.

## §8.4. СВЯЗЬ МЕЖДУ n-АРНЫМИ АНАЛОГАМИ И ИХ БИНАРНЫМИ ПРОТОТИПАМИ

Отметим, что полуцентр $HZ(A)$ n-арной группы $< A, [\,] >$, производной от группы $A$, лежит в центре $Z(A)$ этой группы. Действительно, если $z \in HZ(A)$, то

$$[z\underbrace{e \ldots e}_{n-2} x] = [x\underbrace{e \ldots e}_{n-2} z]$$

для любого $x \in A$, откуда $zx = xz$ и значит $z \in Z(A)$. Следовательно, $HZ(A) \subseteq Z(A)$.

В действительности имеет место более сильное утверждение.

**8.4.1. Предложение** [119]. Пусть $< A, [\,] > -$ n-арная груп-

па, a ∈ HZ(A). Тогда < HZ(A), ⓐ > – инвариантная подгруппа группы < A, ⓐ >, лежащая в её центре < Z(A), ⓐ >.

*Доказательство.* Если z ∈ HZ(A), то

$$[z\bar{a}\underbrace{a\ldots a}_{n-3}x] = [x\bar{a}\underbrace{a\ldots a}_{n-3}z], x \in A,$$

откуда z ⓐ x = x ⓐ z для любого x ∈ A. Следовательно, z ∈ Z(a) и

$$HZ(A) \subseteq Z(A, ⓐ).$$

То, что < HZ(A), ⓐ > – подгруппа в < A, ⓐ >, следует из следствия 8.3.3 и теоремы 2.2.3. Кроме того, для любого x ∈ A верно

$$HZ(A) ⓐ x = [HZ(A)\bar{a}\underbrace{a\ldots a}_{n-3}x] = [x\bar{a}\underbrace{a\ldots a}_{n-3}HZ(A)] = x ⓐ HZ(A),$$

то есть

$$HZ(A) ⓐ x = x ⓐ HZ(A).$$

Следовательно, < HZ(A), ⓐ > – инвариантна в < A, ⓐ >.  ∎

**8.4.2. Замечание.** Инвариантность < HZ(A), ⓐ > в < A, ⓐ > вытекает также из полуинвариантности < HZ(A), [ ] > в < A, [ ] > и следствия 2.3.13.

**8.4.3. Предложение** [119]. Пусть < T = HZ(A), [ ] > – полуцентр n-арной группы < A, [ ] >. Тогда $T_o(A)$ – подгруппа центра $Z(A_o)$ группы $A_o$.

*Доказательство.* Пусть $u = \theta_A(z_1 \ldots z_{n-1})$ – произвольный элемент из $T_o(A)$, где $z_1, \ldots, z_{n-1} \in T$, $v = \theta_A(a_1 \ldots a_{n-1})$ – произвольный элемент из $A_o$.

Так как $z_i \in T$, то

$$uv = \theta_A(z_1 \ldots z_{n-1})\theta_A(a_1 \ldots a_{n-1}) = \theta_A(z_1 \ldots z_{n-1}a_1 \ldots a_{n-1}) =$$

$$= \theta_A(z_1 \ldots z_{n-2}[z_{n-1}a_1 \ldots a_{n-2}a_{n-1}]) =$$

$$= \theta_A(z_1 \ldots z_{n-2}[a_{n-1}a_1 \ldots a_{n-2}z_{n-1}]) =$$

$$= \theta_A(z_1 \ldots z_{n-3}[z_{n-2}a_{n-1}a \ldots a_{n-3}a_{n-2}]z_{n-1}) =$$

$$= \theta_A(z_1 \ldots z_{n-2}[a_{n-2}a_{n-1}a_1 \ldots a_{n-3}z_{n-2}]z_{n-1}) = \ldots$$

$$\ldots = \theta_A([z_1 a_2 \ldots a_{n-1}a_1]z_2 \ldots z_{n-1}) =$$

$$= \theta_A([a_1 a_2 \ldots a_{n-1}z_1]z_2 \ldots z_{n-1}) =$$

$$= \theta_A(a_1 \ldots a_{n-1}z_1 \ldots z_{n-1}) = \theta_A(a_1 \ldots a_{n-1})\theta_A(z_1 \ldots z_{n-1}) = vu,$$

то есть $uv = vu$. Следовательно, $u \in Z(A_o)$ и доказано включение $T_o(A) \subseteq Z(A_o)$. Осталось воспользоваться замечанием 2.2.20, согласно которому, $T_o(A)$ – подгруппа в $A_o$. ∎

**8.4.4. Теорема** [119]. Пусть $< T = C_A(B), [\,] > $ – централизатор подмножества $B$ в n-арной группе $< A, [\,] >$. Тогда $T^*(A)$ совпадает с централизатором $C_{A^*}(B^*(A))$ подмножества $B^*(A)$ в группе $A^*$: $T^*(A) = C_{A^*}(B^*(A))$.

*Доказательство.* По следствию 8.2.3 $< T, [\,] >$ – n-арная подгруппа в $< A, [\,] >$, а по теореме 2.2.19 $T^*(A)$ – подгруппа группы $A^*$.

Согласно 2) теоремы 2.2.19, произвольный элемент $u \in T^*(A)$ может быть представлен в виде

$$u = \theta_A(z_1 \ldots z_i)$$

где $z_1, \ldots, z_i \in T$ $(i = 1, \ldots, n-1)$. Произвольный элемент $v$ подмножества $B^*(A)$ имеет вид

$$v = \theta_A(a_1 \ldots a_j)$$

где $a_1, \ldots, a_j \in B$ $(j = 1, \ldots, n-1)$.

Так как $z_i \in T = C_A(B)$, то $\theta_A(z_i a_j) = \theta_A(a_j z_i)$ для всех $i, j \in \{1, \ldots, n-1\}$, откуда

$$uv = \theta_A(z_1 \ldots z_i)\theta_A(a_1 \ldots a_j) = \theta_A(z_1 \ldots z_i a_1 \ldots a_j) =$$

$$= \theta_A(z_1 \ldots z_{i-1})\theta_A(z_i a_1)\theta_A(a_2 \ldots a_j) =$$

$$= \theta_A(z_1 \ldots z_{i-1})\theta_A(a_1 z_i)\theta_A(a_2 \ldots a_j) =$$

$$= \theta_A(z_1 \ldots z_{i-1} a_1 z_i a_2 \ldots a_j) = \ldots = \theta_A(z_1 \ldots z_{i-1} a_1 \ldots a_j z_i) = \ldots$$

$$\ldots = \theta_A(a_1 \ldots a_j z_1 \ldots z_i) = \theta_A(a_1 \ldots a_j)\theta_A(z_1 \ldots z_i) = vu,$$

то есть $uv = vu$. Следовательно, $u \in C_{A^*}(B^*(A))$ и доказано включение

$$T^*(A) \subseteq C_{A^*}(B^*(A)). \qquad (1)$$

Пусть теперь

$$w = \theta_A(a_1 \ldots a_i), a_i \in A$$

– произвольный элемент из $C_{A^*}(B^*(A))$, $z_2, \ldots, z_i$ – произвольные элементы из T. Так как существует $z_1 \in A$ такой, что

$$\theta_A(a_1 \ldots a_i) = \theta_A(z_1 z_2 \ldots z_i),$$

то

$$w = \theta_A(z_1)\theta_A(z_2 \ldots z_i) \in C_{A^*}(B^*(A)),$$

откуда

$$\theta_A(z_1) \in C_{A^*}(B^*(A)),$$

так как по доказанному выше

$$\theta_A(z_2 \ldots z_i) \in C_{A^*}(B^*(A)).$$

Тогда

$$\theta_A(z_1)\theta_A(b) = \theta_A(b)\theta_A(z_1)$$

для любого $b \in B$, откуда

$$\theta_A(z_1 b) = \theta_A(b z_1),$$

$$z_1 b \sim b z_1,$$

$$z_1 \in C_A(B) = T.$$

Это означает, что

$$w = \theta_A(z_1 z_2 \ldots z_i) \in T^*(A),$$

и следовательно,

$$C_{A^*}(B^*(A)) \subseteq T^*(A). \qquad (2)$$

Из (1) и (2) следует требуемое равенство. ∎

Согласно замечанию 2.2.20 обертывающая группа $T^*$ n-арной подгруппы $<T, [\ ]>$ n-арной группы $<A, [\ ]>$ изоморфна подгруппе $T^*(A)$ обертывающей группы $A^*$. Поэтому из теоремы 8.4.4 вытекает

**8.4.5. Следствие** [119]. Обертывающая группа $T^*$ централизатора $<T = C_A(B), [\ ]>$ n-арной подгруппы $<B, [\ ]>$ в n-арной группе $<A, [\ ]>$ изоморфна централизатору подгруппы $B^*(A)$ в обертывающей группе $A^*$:

$$T^* \simeq C_{A^*}(B^*(A)).$$

Если $B = A$, то

$$C_A(A) = Z(A),\ B^*(A) = A^*,\ C_{A^*} B^*(A) = Z(A^*).$$

Поэтому из теоремы 8.4.4 вытекает

**8.4.6. Следствие.** Пусть $<T = Z(A), [\ ]>$ центр n-арной группы $<A, [\ ]>$. Тогда $T^*(A)$ совпадает с центром $Z(A^*)$ группы $A^*$: $T^*(A) = Z(A^*)$.

**8.4.7. Следствие.** Обертывающая группа $T^*$ центра $<T = Z(A), [\ ]>$ n-арной группы $<A, [\ ]>$ изоморфна центру группы $Z(A^*)$ обертывающей группы $A^*$:

$$T^* \simeq Z(A^*).$$

## §8.5. m-ПОЛУЦЕНТРАЛИЗАТОР ЭЛЕМЕНТА

**8.5.1. Определение.** m-Полуцентрализатор подмножества $M = \{a\}$ в n-арной группе $< A, [\ ] >$, где $n = k(m-1) + 1$, $k \geq 1$, называется *m-полуцентрализатором элемента* a в n-арной группе $< A, [\ ] >$ и обозначается символом $C_A(a, m)$.

Таким образом,

$$C_A(a, m) = \{z \in A \mid z\underbrace{a \ldots a}_{m-1} \sim \underbrace{a \ldots a}_{m-1} z\}.$$

В частности,

$$C_A(a) = \{z \in A \mid za \sim az\}$$

– централизатор элемента a в n-арной группе $< A, [\ ] >$,

$$HC_A(a) = \{z \in A \mid [z\underbrace{a \ldots a}_{n-1}] = [\underbrace{a \ldots a}_{n-1} z]\}$$

– полуцентрализатор элемента a в n-арной группе $< A, [\ ] >$.

Так как

$$a^{[k]}\underbrace{a \ldots a}_{m-1} \sim \underbrace{a \ldots a}_{m-1} a^{[k]}$$

для любой степени элемента a, то m-полуцентрализатор $C_A(a, m)$ вместе с элементом a содержит и всего его степени, в том числе, и косой элемент $\bar{a}$.

Таким образом имеет место

**8.5.2. Предложение** [119]. Пусть $< A, [\ ] >$ – n-арная группа, $a \in A$, $n = k(m-1) + 1$, $k \geq 1$. Тогда:

1) m-полуцентрализатор $C_A(a, m)$ включает в себя циклическую n-арную подгруппу, порождённую элементом a;

2) если $<A, [\,]>$ конечная, то порядок элемента a делит порядок $<C_A(a, m), [\,]>$.

Полагая в предложении 8.5.2 m = 1, получим

**8.5.3. Следствие.** Пусть $<A, [\,]>$ – n-арная группа, $a \in A$. Тогда:

1) централизатор $C_A(a)$ включает в себя циклическую n-арную подгруппу, порожденную элементом a;

2) если $<A, [\,]>$ конечная, то порядок элемента a делит порядок $<C_A(a), [\,]>$.

Полагая в предложении 8.5.2 m = n, получим

**8.5.4. Следствие.** Пусть $<A, [\,]>$ – n-арная группа, $a \in A$. Тогда:

1) полуцентрализатор $HC_A(a)$ включает в себя циклическую n-арную подгруппу, порожденную элементом a;

2) если $<A, [\,]>$ конечная, то порядок элемента a делит порядок $<HC_A(a), [\,]>$.

**8.5.5. Замечание.** Так как

$$C_A(a) \subseteq C_A(a, m) \subseteq HC_A(a),$$

то предложение 8.5.2 и следствие 8.5.4 могут быть получены как следствия, вытекающие из следствия 8.5.3.

**8.5.6. Следствие.** Если порядок элемента a n-арной группы $<A, [\,]>$, где $n = k(m-1) + 1$, $k \geq 1$, равен 2, то n-арная группа $<C_A(a, m), [\,]>$ имеет четный порядок.

**8.5.7. Следствие.** Если порядок элемента a n-арной группы $<A, [\,]>$ равен 2, то n-арные группы $<C_A(a), [\,]>$ и $<HC_A(a), [\,]>$ имеют четный порядок.

## §8.6. СЛАБЫЙ m-ПОЛУЦЕНТРАЛИЗАТОР

**8.6.1. Определение**[119, 121]. Если $n = k(m-1) + 1$, $k \geq 1$,

то *слабым m-полуцентрализатором (m-полуцентрализатором типа D)* подмножества B в n-арной группе $< A, [\ ] >$ называется множество

$$DC_A(B, m) = \{z \in A \mid z\underbrace{x \ldots x}_{m-1} \sim \underbrace{x \ldots x}_{m-1} z, \ \forall x \in B\}.$$

Ясно, $C_A(B, m) \subseteq DC_A(B, m)$.

Если $m = 2$, то определение слабого 2-полуцентрализатора совпадает с определениями 2-полуцентрализатора и централизатора, то есть

$$DC_A(B, 2) = C_A(B, 2) = C_A(B).$$

Если $m = n$, то слабый n-полуцентрализатор называется *слабым полуцентрализатором (полуцентрализатором типа D)* и обозначается символом $HDC_A(B)$, то есть

$$HDC_A(B) = \{z \in A \mid [z\underbrace{x \ldots x}_{n-1}] = [\underbrace{x \ldots x}_{n-1} z], \ \forall x \in B\}.$$

Если в приведенном определении положить $A = B$, то получим определение *слабого m-полуцентра (m-полуцентра типа D)* n-арной группы $< A, [\ ] >$:

$$DZ(A, m) = \{z \in A \mid z\underbrace{x \ldots x}_{m-1} \sim \underbrace{x \ldots x}_{m-1} z, \ \forall x \in A\}.$$

Слабый 2-полуцентр $DZ(A, 2)$ совпадает с 2-полуцентром и центром, то есть

$$DZ(A, 2) = Z(A, 2) = Z(A).$$

Слабый n-полуцентр называется *слабым полуцентром (полуцентром типа D)* и обозначается символом $HDZ(A)$, то есть $DZ(A, n) = HDZ(A)$. Таким образом,

$$HDZ(A) = \{z \in A \mid [z\underbrace{x \ldots x}_{n-1}] = [\underbrace{x \ldots x}_{n-1} z], \ \forall x \in A\}.$$

Легко заметить, что HZ(A) ⊆ HDZ(A).

Ясно, что слабый m-полуцентрализатор элемента a и m-полуцентрализатор этого же элемента a совпадают, то есть

$$DC_A(a, m) = C_A(a, m).$$

В частности,

$$HDC_A(a) = HC_A(a).$$

Имеет место очевидное

**8.6.2. Предложение.** Если < A, [ ] > – n-арная группа, a ∈ A, то

$$DZ(A, m) = \bigcap_{a \in A} C_A(a, m).$$

В частности,

$$HDZ(A) = \bigcap_{a \in A} HC_A(a).$$

Легко проверяется и справедливость следующего утверждения.

**8.6.3. Предложение.** n-Арная группа < A, [ ] > является слабо полуабелевой тогда и только тогда, когда HDZ(A) = A.

**8.6.4. Теорема** [119, 121]. Пусть < A, [ ] > – n-арная группа, B ⊆ A. Тогда:

1) если m – 1 делит n – 1 и $DC_A(B, m) \neq \varnothing$, то тогда < $DC_A(B, m)$, [ ] > – n-арная подгруппа в < A, [ ] >, в частности, < $HDC_A(B)$, [ ] > – n-арная подгруппа в < A, [ ] >;

2) если k – 1 делят n – 1, m – 1 делит k – 1, то

$$DC_A(B, m) \subseteq DC_A(B, k);$$

3) если m – 1 и k – 1 делят n – 1, r – 1 = (m – 1, k – 1), то

$$DC_A(B, r) = DC_A(B, m) \cap DC_A(B, k);$$

4) если m – 1 делит n – 1, r – 1 = (m – 1, n – 1), то

$$DC_A(B, r) = DC_A(B, m) \cap HDC_A(B);$$

5) если m – 1 и k – 1 делят n – 1, (m – 1, k – 1) = 1, то

$$C_A(B) = DC_A(B, m) \cap DC_A(B, k);$$

6) если m – 1 делит n – 1, (m – 1, n – 1) = 1, то

$$C_A(B) = DC_A(B, m) \cap HDC_A(B).$$

***Доказательство.*** 1) Если $z_1, \ldots, z_n \in DC_A(B, m)$, то из определения 8.6.1 следует

$$[z_1 \ldots z_n] \underbrace{x \ldots x}_{m-1} \sim \underbrace{x \ldots x}_{m-1} [z_1 \ldots z_n]$$

для любого $x \in B$. Следовательно, $[z_1 \ldots z_n] \in DC_A(B, m)$.

Если теперь $z \in DC_A(B, m)$, то, учитывая нейтральность последовательностей $\underbrace{z \ldots z}_{n-2} \bar{z}$ и $\bar{z} \underbrace{z \ldots z}_{n-2}$, получим

$$z \underbrace{x \ldots x}_{m-1} \sim \underbrace{x \ldots x}_{m-1} z,$$

$$\bar{z} z \underbrace{x \ldots x}_{m-1} \underbrace{z \ldots z}_{n-3} \bar{z} \sim \bar{z} \underbrace{x \ldots x}_{m-1} z \underbrace{z \ldots z}_{n-3} \bar{z},$$

$$\bar{z} \underbrace{x \ldots x}_{m-1} z \underbrace{z \ldots z}_{n-3} \bar{z} \sim \bar{z} \underbrace{z \ldots z}_{n-3} \underbrace{x \ldots x}_{m-1} \bar{z},$$

$$\bar{z} \underbrace{x \ldots x}_{m-1} \sim \underbrace{x \ldots x}_{m-1} \bar{z}$$

для любого $x \in B$. Следовательно, $\bar{z} \in DC_A(B, m)$. Согласно критерию Дёрнте, $< DC_A(B, m), [\ ] >$ – n-арная подгруппа в $< A, [\ ] >$.

2) Так как k – 1 = t(m – 1), для некоторого целого t то из $z \in DC_A(B, m)$ следует

$$\underbrace{z x \ldots x}_{k-1} = \underbrace{z x \ldots x}_{t(m-1)} = \underbrace{z x \ldots x}_{m-1} \underbrace{x \ldots x}_{(t-1)(m-1)} \sim$$

$$\sim \underbrace{x \ldots x}_{m-1} \underbrace{z x \ldots x}_{(t-1)(m-1)} \sim \ldots \sim \underbrace{x \ldots x}_{t(m-1)} z = \underbrace{x \ldots x}_{k-1} z,$$

то есть

$$\underbrace{z x \ldots x}_{k-1} \sim \underbrace{x \ldots x}_{k-1} z$$

для любого $x \in B$. Следовательно, $z \in DC_A(B, k)$ и верно включение

$$DC_A(B, m) \subseteq DC_A(B, k).$$

3) Включение

$$DC_A(B, r) \subseteq DC_A(B, m) \bigcap DC_A(B, k) \qquad (1)$$

следует из 2).

Так как $r - 1 = (m - 1, k - 1)$, то существуют целые числа $\alpha$ и $\beta$ такие, что

$$\alpha(m - 1) + \beta(k - 1) = r - 1.$$

Пусть для определенности $\alpha > 0$, $\beta < 0$, то есть

$$\alpha(m - 1) = -\beta(k - 1) + (r - 1), -\beta(k - 1) > 0.$$

Если

$$z \in DC_A(B, m) \bigcap DC_A(B, k), x \in B,$$

то

$$\underbrace{z x \ldots x}_{r-1} \underbrace{x \ldots x}_{-\beta(k-1)} = \underbrace{z x \ldots x}_{\alpha(m-1)} \sim \underbrace{x \ldots x}_{\alpha(m-1)} z =$$

$$= \underbrace{x \ldots x}_{r-1} \underbrace{x \ldots x}_{-\beta(k-1)} z \sim \underbrace{x \ldots x}_{r-1} z \underbrace{x \ldots x}_{-\beta(k-1)},$$

откуда следует

$$\underbrace{zx \ldots x}_{r-1} \sim \underbrace{x \ldots x z}_{r-1}.$$

Следовательно, $z \in DC_A(B, r)$ и доказано включение

$$DC_A(B, m) \cap DC_A(B, k) \subseteq DC_A(B, r). \qquad (2)$$

Из (1) и (2) следует требуемое равенство.
4), 5) и 6) следуют из 3). ∎

**8.6.5. Следствие** [119, 121]. Пусть $< A, [\ ] >$ – n-арная группа. Тогда:

1) если $m - 1$ делит $n - 1$, $DZ(A, m) \neq \varnothing$, то тогда $< DZ(A, m), [\ ] >$ – n-арная подгруппа в $< A, [\ ] >$, в частности, $< HDZ(A), [\ ] >$ – n-арная подгруппа в $< A, [\ ] >$;

2) если $k - 1$ делит $n - 1$, $m - 1$ делит $k - 1$, то

$$DZ(A, m) \subseteq DZ(A, k);$$

3) если $m - 1$ и $k - 1$ делят $n - 1$, $r - 1 = (m - 1, n - 1)$, то

$$DZ(A, r) = DZ(A, m) \cap DZ(A, k);$$

4) если $m - 1$ делит $n - 1$, $r - 1 = (m - 1, n - 1)$, то

$$DZ(A, r) = DZ(A, m) \cap HDZ(A);$$

5) если $m - 1$ и $k - 1$ делят $n - 1$, $(m - 1, k - 1) = 1$, то

$$Z(A) = DZ(A, m) \cap DZ(A, k);$$

6) если $m - 1$ делит $n - 1$, $(m - 1, n - 1) = 1$, то

$$Z(A) = DZ(A, m) \cap HDZ(A).$$

### §8.7. m-ПОЛУЦЕНТРАЛИЗАТОР ТИПА T

В теореме 7.2.1 доказано, что если $m - 1$ делит $k - 1$, $k - 1$

делит n − 1, то

$$N_A(B, m) \subseteq N_A(B, k).$$

В связи с этим возникает вопрос: будет ли верным включение

$$C_A(B, m) \subseteq C_A(B, k)$$

при тех же m, k и n, что и выше? Из примера 8.1.2 следует, что для n-арной группы < A, [ ] >

$$C_A(B, m) \not\subset C_A(B, k)$$

при m = 2, k = n. Поэтому ответ на поставленный вопрос является отрицательным.

Ситуация меняется к лучшему, если рассмотреть еще один n-арный аналог централизатора подмножества в группе.

**8.7.1. Определение** [119, 121]. Если n = k(m − 1) + 1, где k ≥ 1, то *m-полуцентрализатором типа T* подмножества B в n-арной группе < A, [ ] > называется множество

$$TC_A(B, m) = \{z \in A \mid zx_1\ldots x_{m-1} \sim x_1\ldots x_{m-1}z, \forall x_1, \ldots, x_{m-1} \in B\}.$$

Ясно, что $TC_A(B, m) \subseteq DC_A(B, m)$.

Если m = 2, то определение 2-полуцентрализатора типа T совпадает с определениями слабого 2-полуцентрализатора, 2-полуцентрализатора и централизатора, то есть

$$TC_A(B, 2) = DC_A(B, 2) = C_A(B, 2) = C_A(B).$$

Если m = n, то n-полуцентрализатор типа T называется *полуцентрализатором типа T* и обозначается символом $HTC_A(B)$, то есть

$$HTC_A(B) = \{z \in A \mid [zx_1\ldots x_{n-1}] = [x_1\ldots x_{n-1}z], \forall x_1,\ldots,x_{n-1} \in B\}.$$

Если в приведенном определении положить A = B, то получим определение *m-полуцентра типа T* n-арной группы

< A, [ ] >:

$$TZ(A, m) = \{z \in A \mid zx_1\ldots x_{m-1} \sim x_1\ldots x_{m-1}z, \quad \forall x_1,\ldots,x_{m-1} \in A\}.$$

2-Полуцентр $TZ(A, 2)$ типа $T$ совпадает со слабым 2-полуцентром, 2-полуцентром и центром, то есть

$$TZ(A, 2) = DZ(A, 2) = Z(A, 2) = Z(A).$$

n-Полуцентр типа $T$ называется *полуцентром типа T* и обозначается символом $HTZ(A)$, то есть $TZ(A, n) = HTZ(A)$. Таким образом,

$$HTZ(A) = \{z \in A \mid [zx_1\ldots x_{n-1}] = [x_1\ldots x_{n-1}z], \forall x_1,\ldots,x_{n-1} \in A\}.$$

Легко заметить, что $HTZ(A) \subseteq HDZ(A)$.

Ясно, что m-полуцентрализатор типа $T$ элемента $a$, слабый m-полуцентрализатор $a$ и m-полуцентрализатор этого же элемента $a$ совпадают, то есть

$$TC_A(a, m) = DC_A(a, m) = C_A(a, m).$$

В частности,

$$HTC_A(a) = HDC_A(a) = HC_A(a).$$

Имеет место очевидное

**8.7.2. Предложение.** n-Арная группа $< A, [ ] >$ удовлетворяет тождеству

$$[xx_1\ldots x_{n-1}] = [x_1\ldots x_{n-1}x]$$

тогда и только тогда, когда $HTZ(A) = A$.

**8.7.3. Теорема** [119, 121]. Пусть $< A, [ ] >$ – n-арная группа, $B \subseteq A$. Тогда:

1) если $m - 1$ делит $n - 1$ и $TC_A(B, m) \neq \varnothing$, то тогда $< TC_A(B, m), [ ] >$ – n-арная подгруппа в $< A, [ ] >$, в частности, $< HTC_A(B), [ ] >$ – n-арная подгруппа в $< A, [ ] >$;

2) если $k - 1$ делят $n - 1$, $m - 1$ делит $k - 1$, то

$$TC_A(B, m) \subseteq TC_A(B, k);$$

3) если $m - 1$ и $k - 1$ делят $n - 1$, $r - 1 = (m - 1, k - 1)$, то

$$TC_A(B, r) = TC_A(B, m) \cap TC_A(B, k);$$

4) если $m - 1$ делит $n - 1$, $r - 1 = (m - 1, n - 1)$, то

$$TC_A(B, r) = TC_A(B, m) \cap HTC_A(B);$$

5) если $m - 1$ и $k - 1$ делят $n - 1$, $(m - 1, k - 1) = 1$, то

$$C_A(B) = TC_A(B, m) \cap TC_A(B, k);$$

6) если $m - 1$ делит $n - 1$, $(m - 1, n - 1) = 1$, то

$$C_A(B) = TC_A(B, m) \cap HTC_A(B).$$

***Доказательство.*** 1) Если $z_1, \ldots, z_n \in TC_A(B, m)$, то из определения 8.7.1 следует

$$[z_1 \ldots z_n]x_1\ldots x_{m-1} \sim x_1\ldots x_{m-1}[z_1 \ldots z_n]$$

для любых $x_1, \ldots, x_{m-1} \in B$. Следовательно,

$$[z_1 \ldots z_n] \in TC_A(B, m).$$

Если теперь $z \in TC_A(B, m)$, то, учитывая нейтральность последовательностей $\underbrace{z \ldots z}_{n-2} \bar{z}$ и $\bar{z} \underbrace{z \ldots z}_{n-2}$, получим

$$zx_1\ldots x_{m-1} \sim x_1\ldots x_{m-1}z,$$

$$\bar{z}zx_1\ldots x_{m-1}\underbrace{z \ldots z}_{n-3}\bar{z} \sim \bar{z}x_1\ldots x_{m-1}z\underbrace{z \ldots z}_{n-3}\bar{z},$$

$$\bar{z}x_1\ldots x_{m-1}z\underbrace{z \ldots z}_{n-3}\bar{z} \sim \bar{z}z\underbrace{z \ldots z}_{n-3}x_1\ldots x_{m-1}\bar{z},$$

$$\bar{z}x_1\ldots x_{m-1} \sim x_1\ldots x_{m-1}\bar{z}$$

для любых $x_1, \ldots, x_{m-1} \in B$. Следовательно, $\bar{z} \in TC_A(B, m)$. Согласно критерию Дернте, $<TC_A(B, m), [\,]> -$ n-арная под-группа в $< A, [\ ] >$.

2) Так как $k - 1 = t(m - 1)$, для некоторого целого $t$ то из $z \in TC_A(B, m)$ следует

$$zx_1\ldots x_{k-1} = zx_1\ldots x_{t(m-1)} = zx_1\ldots x_{m-1}x_m\ldots x_{t(m-1)} \sim$$

$$\sim x_1\ldots x_{m-1}zx_m\ldots x_{t(m-1)} \sim \ldots \sim x_1\ldots x_{t(m-1)}z = x_1\ldots x_{k-1}z,$$

то есть

$$zx_1\ldots x_{k-1} \sim x_1\ldots x_{k-1}z$$

для любых $x_1, \ldots, x_{k-1} \in B$. Следовательно, $z \in TC_A(B, k)$ и верно включение

$$TC_A(B, m) \subseteq TC_A(B, k).$$

3) Включение

$$TC_A(B, r) \subseteq TC_A(B, m) \cap TC_A(B, k) \qquad (1)$$

следует из 2).

Так как $r - 1 = (m - 1, k - 1)$, то существуют целые числа $\alpha$ и $\beta$ такие, что

$$\alpha(m - 1) + \beta(k - 1) = r - 1.$$

Пусть для определенности $\alpha > 0$, $\beta < 0$, то есть

$$\alpha(m - 1) = -\beta(k - 1) + (r - 1), -\beta(k - 1) > 0.$$

Если

$$z \in TC_A(B, m) \cap TC_A(B, k), x_1, \ldots, x_{\alpha(m-1)} \in B,$$

то

$$zx_1\ldots x_{r-1}x_r\ldots x_{-\beta(k-1)+r-1} = zx_1\ldots x_{\alpha(m-1)} \sim x_{\alpha(m-1)}z =$$

$$= x_1\ldots x_{r-1}x_r\ldots x_{-\beta(k-1)+r-1}z \sim x_1\ldots x_{r-1}zx_r\ldots x_{-\beta(k-1)+r-1},$$

откуда следует

$$zx_1\ldots x_{r-1} \sim x_1\ldots x_{r-1}z.$$

Следовательно, $z \in TC_A(B, r)$ и доказано включение

$$TC_A(B, m) \cap TC_A(B, k) \subseteq TC_A(B, r). \qquad (2)$$

Из (1) и (2) следует требуемое равенство.
4) Следует из 3) при $k = n$, так как $TC_A(B, n) = HTC_A(B)$.
5) Следует из 3) при $r = 2$, так как $TC_A(B, 2) = C_A(B)$.
6) Следует из 4) и 5). ∎

**8.7.4. Следствие.** Пусть $< A, [\,] > -$ n-арная группа. Тогда:

1) если $m - 1$ делит $n - 1$, $TZ(A, m) \neq \varnothing$, то тогда $< TZ(A, m), [\,] > -$ n-арная подгруппа в $< A, [\,] >$, в частности, $< HTZ(A), [\,] > -$ n-арная подгруппа в $< A, [\,] >$;

2) если $k - 1$ делит $n - 1$, $m - 1$ делит $k - 1$, то

$$TZ(A, m) \subseteq TZ(A, k);$$

3) если $m - 1$ и $k - 1$ делят $n - 1$, $r - 1 = (m - 1, k - 1)$, то

$$TZ(A, r) = TZ(A, m) \cap TZ(A, k);$$

4) если $m - 1$ делит $n - 1$, $r - 1 = (m - 1, n - 1)$, то

$$TZ(A, r) = TZ(A, m) \cap HTZ(A);$$

5) если $m - 1$ и $k - 1$ делят $n - 1$, $(m - 1, k - 1) = 1$, то

$$Z(A) = TZ(A, m) \cap TZ(A, k);$$

6) если $m - 1$ делит $n - 1$, $(m - 1, n - 1) = 1$, то

$$Z(A) = TZ(A, m) \cap HTZ(A).$$

## §8.8. $(\Sigma, m)$-ПОЛУЦЕНТРАЛИЗАТОР

m-Полуцентрализаторы и m-полуцентрализаторы типа T можно объединить в рамках общего понятия.

Пусть $<A, [\ ]>$ n-арная группа, $B \subseteq A$, $m-1$ делит $n-1$, $\Sigma$ – подмножество множества $S_{m-1}$ всех подстановок на $m-1$ символах.

**8.8.1. Определение** [119, 121]. *($\Sigma$, m)-полуцентрализатором* подмножества B в $<A, [\ ]>$ называется множество

$$C_A(B, \Sigma, m) = \{z \in A \mid zx_1 \ldots x_{m-1} \sim x_{\sigma(1)} \ldots x_{\sigma(m-1)}z,$$

$$\forall x_1, \ldots, x_{m-1} \in B, \forall \sigma \in \Sigma\}.$$

($\Sigma$, 2)-Полуцентрализатор вырождается в централизатор $C_A(B)$.

($\Sigma$, n)-Полуцентрализатор называется *$\Sigma$-полуцентрализатором* и обозначается символом $HC_A(B, \Sigma)$, то есть

$$HC_A(B, \Sigma) = C_A(B, \Sigma, n).$$

Если $A = B$, то ($\Sigma$, m)-полуцентрализатор $C_A(A, \Sigma, m)$ называется *($\Sigma$, m)-полуцентром* n-арной группы $<A, [\ ]>$ и обозначается символом $Z(A, \Sigma, m)$, то есть

$$Z(A, \Sigma, m) = C_A(A, \Sigma, m).$$

($\Sigma$, 2)-Полуцентр вырождается в центр $Z(A)$.

($\Sigma$, n)-Полуцентр называется *$\Sigma$-полуцентром* и обозначается символом $HZ(A, \Sigma)$, то есть

$$HZ(A, \Sigma) = Z(A, \Sigma, n).$$

Если $\Sigma = \{\sigma\}$, то ($\{\sigma\}$, m)-полуцентрализатор называется *($\sigma$, m)-полуцентрализатором* и обозначается следующим образом: $C_A(B, \sigma, m)$.

Аналогично определяются *($\sigma$, m)-полуцентр* $Z(A, \sigma, m)$, *$\sigma$-полуцентрализатор* $HC_A(B, \sigma)$ и *$\sigma$-полуцентр* $HZ(A, \sigma)$.

Ясно, что если

$$\tau = (m - 1 \; m - 2 \ldots 2 \; 1),$$

то определение $(\tau, m)$-полуцентрализатора совпадает с определением m-полуцентрализатора.

Если же $\varepsilon$ – тождественная подстановка на $m - 1$ символах, то определение $(\varepsilon, m)$-полуцентрализатора совпадает с определением m-полуцентрализатора типа T.

Так как $\tau^m = \tau$, $\varepsilon^m = \varepsilon$, то в связи с теоремой 8.2.2 закономерен

**8.8.2. Вопрос.** Если $\sigma^m = \sigma$, то будет ли $(\sigma, m)$-полуцентрализатор $C_A(B, \sigma, m)$ n-арной подгруппой в $< A, [\;] >$?

Следующий вопрос связан с тем, что $\{\varepsilon\}$ – подгруппа в симметрической группе $S_{m-1}$.

**8.8.3. Вопрос.** Если $\Sigma$ – подгруппа в симметрической груп-пе $S_{m-1}$, то будет ли $(\Sigma, m)$-полуцентрализатор $C_A(B, \Sigma, m)$ n-арной подгруппой в $< A, [\;] >$?

Следующий вопрос связан с тем, что $< \{\tau\}, [\;] >$ и $< \{\varepsilon\}, [\;] >$ – m-арные подгруппы в m-арной группе $< S_{m-1}, (\;) >$, производной от симметрической группы $S_{m-1}$.

**8.8.4. Вопрос.** Если $\Sigma$ – m-арная подгруппа в m-арной группе $< S_{m-1}, (\;) >$, где

$$(\sigma_1 \ldots \sigma_m) = \sigma_1 \ldots \sigma_m,$$

то будет ли $(\Sigma, m)$-полуцентрализатор $C_A(B, \Sigma, m)$ n-арной подгруппой в $< A, [\;] >$?

Ясно, что при положительном ответе на вопрос 8.8.4, положительными будут и ответы на вопросы 8.8.2 и 8.8.3.

# ДОПОЛНЕНИЯ И КОММЕНТАРИИ

**1.** Понятие центра n-арной группы $<A, [\,]>$, введенное Постом [3], С.А. Русаков использовал в [122] для определения n-арных групп с центральными рядами.

**Определение** [122, Русаков С.А.]. Пусть $<A, [\,]>$ – n-арная группа с непустым центром. Ряд n-арных подгрупп

$$A_o \subseteq A_1 \subseteq \ldots \subseteq A_i \subseteq A_{i+1} \subseteq \ldots \subseteq A_{k-1} \subseteq A_k = A, \qquad (*)$$

в котором каждая n-арная подгруппа $<A_i, [\,]>$ инвариантна в $<A, [\,]>$ и, кроме того

$$A_o \subseteq Z(A), \ A_{i+1}/A_i \subseteq Z(A/A_i), i = 0, 1, \ldots, k-1,$$

называется *центральным рядом* n-арной группы $<A, [\,]>$.

n-Арные группы, обладающими центральными рядами, С.А. Русаков назвал нильпотентными.

Некоторые свойства таких n-арных групп приведены в следующей теореме.

**Теорема** [122]. Пусть n-арная группа $<A_i, [\,]>$ обладает центральным рядом (*). Тогда:

1) если $<B, [\,]>$ – n-арная подгруппа в $<A, [\,]>$ и $B \cap A_o \neq \varnothing$, то $<B, [\,]>$ обладает центральным рядом;

2) гомоморфный образ n-арной группы $<A, [\,]>$ обладает центральным рядом;

3) если $<B, [\,]>$ – полуинвариантная n-арная подгруппа в $<A, [\,]>$, то n-арная фактор-группа $<A/B, [\,]>$ обладает центральным рядом;

4) любая собственная n-арная подгруппа n-арной группы $<A, [\,]>$ отлична от своего нормализатора в $<A, [\,]>$;

5) любая максимальная n-арная подгруппа n-арной группы $<A, [\,]>$ инвариантна в $<A, [\,]>$.

**2.** Результаты С.А. Русакова из [122] об n-арных группах с центральными рядами продублированы С. Илич в [123]. Там же в [123] есть аналог следствия 8.4.6.

**3.** В. Дудек заметил [40], что m-полуабелева n-арная группа, в частности полуабелева n-арная группа с непустым центром является абелевой.

**4.** В [41] В. Дудек для всякого элемента a n-арной группы $<A, [\ ]>$ определил непустое множество

$$\{x \in A \mid [x\underbrace{a\ldots a}_{n-1}] = [\underbrace{a\ldots a}_{n-1} x]\},$$

которое, как несложно заметить, совпадает с полуцентрализатором $HC_A(a)$ элемента a в n-арной группе $<A, [\ ]>$.

**Теорема** [41]. Пусть $<A, [\ ]>$ – n-арная группа. Тогда
1) множество $HC_A(a)$ совпадает с централизатором элемента a в группе $<A, \text{ⓐ}>$, то есть $HC_A(a) = C_{<A, \text{ⓐ}>}(a)$;
2) $<HC_A(a), [\ ]>$ – n-арная подгруппа в $<A, [\ ]>$;
3) $<A, [\ ]>$ – слабо полуабелева тогда и только тогда, когда $HC_A(a) = A$ для любого $a \in A$.

В [41] утверждение 2) предыдущей теоремы сформулировано для слабо m-полуабелевых n-арных групп, хотя доказательство приведено для произвольной n-арной группы.

Отметим, что утверждение 2) предыдущей теоремы является следствием теоремы 8.2.2.

**5.** n-Арную группу $<A, [\ ]>$ назовем *m-полуабелевой типа T*, где $m - 1$ делит $n - 1$, если для любых $x_1, \ldots, x_m$ последовательности

$$x_1 x_2 \ldots x_m \text{ и } x_2 \ldots x_m x_1$$

эквивалентны.

Из определения эквивалентных последовательностей и теоремы 1.3.3 вытекает

**Предложение 1.** n-Арная группа $<A, [\ ]>$ является m-полуабелевой типа T, если в ней выполняется тождество

$$[xx_1\ldots x_{m-1}\underbrace{x\ldots x}_{n-m}] = [x_1\ldots x_{m-1}\underbrace{x\ldots x}_{n-m+1}].$$

Ясно, что 2-полуабелевость типа T n-арной группы совпадает с ее полуабелевостью.

n-Полуабелевую типа T n-арную группу назовем *полуабелевой типа T*. Таким образом, n-арная группа < A, [ ] > называется полуабелевой типа T, если в ней выполняется тождество

$$[xx_1\ldots x_{n-1}] = [x_1\ldots x_{n-1}x].$$

из предложения 8.7.2.

Всякая m-полуабелева типа T n-арная группа является полуабелевой типа T.

Понятно, что из m-полуабелевости типа T n-арной группы следует ее слабая m-полуабелевость. В частности, из полуабелевости типа T n-арной группы следует ее слабая полуабелевость.

Легко также заметить, что m-полуцентр типа T m-полуабелевой типа T n-арной группы совпадает с ней самой, то есть TZ(A, m) = A. В частности, полуцентр типа T полуабелевой типа T n-арной группы совпадает с ней самой, то есть HTZ(A) = A.

Следствие 8.7.4 позволяет сформулировать

**Предложение 2.** Если n-Арная группа является m-полуабелевой типа T и k-полуабелевой типа T, то она является r-полуабелевой типа T, где

$$r - 1 = (m - 1, k - 1).$$

**Следствие 1.** Любая m-полуабелева типа T n-арная группа является r-полуабелевой типа T, где

$$r - 1 = (m - 1, n - 1).$$

**Следствие 2.** Если (m – 1, k – 1) = 1, то n-арная группа < A, [ ] >, являющаяся одновременно m-полуабелевой типа T и k-полуабелевой типа T, будет абелевой.

**Следствие 3.** Если (m – 1, n – 1) = 1, то n-арная группа < A, [ ] >, являющаяся m-полуабелевой типа T, будет абелевой.

# Г Л А В А  9

# n-АРНЫЕ ГРУППЫ С ИДЕМПОТЕНТАМИ

n-Арные группы являются очень широким обобщением понятия группы. Поэтому в исследованиях по n-арным группам целесообразно выделять n-арные группы в той или иной мере близкие к группам. При этом желательно, чтобы выделяемые n-арные группы были не настолько близкими к группам, как производные n-арные группы, многие свойства которых являются простыми следствиями соответствующих групповых результатов. Последнему требованию удовлетворяет изучающийся в данной главе достаточно широкий класс n-арных групп с идемпотентами, содержащий помимо всех производных n-арных групп и всех идемпотентных n-арных групп, также все конечные n-арные группы, порядки которых взаимно просты с n − 1, и, как будет установлено, все n-арные группы, допускающие автоморфизм, оставляющий неподвижным единственный элемент.

## §9.1. n-АРНАЯ ПОДГРУППА ЕДИНИЦ

Так как в n-арной группе, в отличие от группы, может быть несколько единиц, то возникает задача изучения множества всех единиц произвольной n-арной группы. Множество это до сих пор никем не изучалось. Поэтому отсутствовала какая-либо информация о его строении. Неизвестно было даже, является ли оно в общем случае n-арной подгруппой в n-арной группе. Утвердительный ответ на этот, как оказалось, несложный вопрос получен в самом начале этого параграфа.

Для всякой n-арной группы $\mathscr{A}$ = < A, [ ] > обозначим че-рез $E(\mathscr{A}) = E(A)$ множество всех её единиц.

**9.1.1. Теорема** [124, 125]. Если $E(A) \neq \varnothing$, то $< E(A), [\ ] >$ – характеристическая n-арная подгруппа n-арной группы $< A, [\ ] >$, лежащая в её центре.

*Доказательство.* Из определения единицы вытекает, что для любой единицы $e \in A$ и любого $x \in A$ последовательности $ex$ и $xe$ эквивалентны. Поэтому, если

$$e = [e_1 e_2 \ldots e_n],$$

где $e_1, e_2, \ldots, e_n$ – единицы n-арной группы $< A, [\ ] >$, то для любого $i = 1, 2, \ldots, n$ и любого $x \in A$ имеем

$$[\underbrace{e \ldots e}_{i-1} x \underbrace{e \ldots e}_{n-i}] =$$

$$= [[\underbrace{[e_1 e_2 \ldots e_n] \ldots [e_1 e_2 \ldots e_n]}_{i-1}] x \underbrace{[e_1 e_2 \ldots e_n] \ldots [e_1 e_2 \ldots e_n]}_{n-i}]] =$$

$$= [\underbrace{e_n \ldots e_n}_{i-1} [\ldots [\underbrace{e_2 \ldots e_2}_{i-1} [\underbrace{e_1 \ldots e_1}_{i-1} x \underbrace{e_1 \ldots e_1}_{n-i}] \underbrace{e_2 \ldots e_2}_{n-i}] \ldots] \underbrace{e_n \ldots e_n}_{n-i}] =$$

$$= [\underbrace{e_n \ldots e_n}_{i-1} [\ldots [\underbrace{e_2 \ldots e_2}_{i-1} x \underbrace{e_2 \ldots e_2}_{n-i}] \ldots] \underbrace{e_n \ldots e_n}_{n-i}] = \ldots$$

$$\ldots = [\underbrace{e_n \ldots e_n}_{i-1} x \underbrace{e_n \ldots e_n}_{n-i}] = x,$$

то есть

$$[\underbrace{e \ldots e}_{i-1} x \underbrace{e \ldots e}_{n-i}] = x.$$

Следовательно,

$$[e_1 e_2 \ldots e_n] = e \in E(A).$$

Учитывая, что $\bar{e} = e$ для любого $e \in E(A)$, и применяя критерий Дёрнте, заключаем, что $< E(A), [\ ] >$ – n-арная подгруппа n-арной группы $< A, [\ ] >$.

Если $\alpha \in \operatorname{Aut} A$, $e \in E(A)$, $x \in A$, то

$$[\underbrace{e^{\alpha}...e^{\alpha}}_{i-1} x \underbrace{e^{\alpha}...e^{\alpha}}_{n-i}] = [\underbrace{e^{\alpha}...e^{\alpha}}_{i-1}(x^{\alpha^{-1}})^{\alpha}\underbrace{e^{\alpha}...e^{\alpha}}_{n-i}] =$$

$$= [\underbrace{e...e}_{i-1} x^{\alpha^{-1}} \underbrace{e...e}_{n-i}]^{\alpha} = (x^{\alpha^{-1}})^{\alpha} = x,$$

то есть

$$[\underbrace{e^{\alpha}...e^{\alpha}}_{i-1} x \underbrace{e^{\alpha}...e^{\alpha}}_{n-i}] = x$$

для любого $x \in A$ и любого $i = 1, 2, ..., n$. Следовательно,

$$E^{\alpha}(A) \subseteq E(A)$$

для любого $\alpha \in \operatorname{Aut} A$, в частности,

$$E^{\alpha^{-1}}(A) \subseteq E(A),$$

откуда

$$(E^{\alpha^{-1}}(A))^{\alpha} \subseteq E^{\alpha}(A), \ \ E(A) \subseteq E^{\alpha}(A).$$

Таким образом, $E^{\alpha}(A) = E(A)$, и $< E(A), [\ ] >$ – характеристична в $< A, [\ ] >$.

Так как для любых $e \in E(A)$, $x \in A$ последовательности $ex$ и $xe$ эквивалентны, то $e \in Z(A)$, то есть $E(A) \subseteq Z(A)$. ∎

**9.1.2. Определение.** Если $E(A) \neq \varnothing$, то n-арную подгруппу $< E(A), [\ ] >$ назовем *n-арной подгруппой единиц* n-арной группы $< A, [\ ] >$.

**9.1.3. Следствие.** n-Арная подгруппа единиц n-арной группы инвариантна в ней.

**9.1.4. Замечание.** Так как всякая единица n-арной группы совпадает со своим косым, то, согласно критерию Дёрнте, для доказательства того, что некоторое подмножество n-арной подгруппы единиц является её n-арной подгруппой, достаточно установить замкнутость n-арной операции.

Ясно, что если $e \in E(A)$, то $<\{e\}, [\ ]>$ – n-арная подгруппа в $<E(A), [\ ]>$.

Как показывает следующее предложение, в n-арной подгруппе единиц при $n > 2$ могут существовать n-арные подгруппы, отличные от одноэлементных и от самой n-арной подгруппы $<E(A), [\ ]>$.

**9.1.5. Предложение** [124, 125]. Если $<A, [\ ]>$ – тернарная группа, $e_1, e_2 \in E(A)$, то $<\{e_1, e_2\}, [\ ]>$ – тернарная подгруппа тернарной группы $<E(A), [\ ]>$.

*Доказательство.* Так как

$$[e_1 e_1 e_2] = [e_1 e_2 e_1] = [e_2 e_1 e_1] = e_2,$$

$$[e_2 e_2 e_1] = [e_2 e_1 e_2] = [e_1 e_2 e_2] = e_1,$$

$$[e_1 e_1 e_1] = e_1, [e_2 e_2 e_2] = e_2,$$

то $<\{e_1, e_2\}, [\ ]>$ – тернарная полугруппа. Учитывая замечание 9.1.4, видим, что $<\{e_1, e_2\}, [\ ]>$ – тернарная группа. ∎

**9.1.6. Следствие** [124, 125]. Если конечная тернарная группа содержит более одной единицы, то её n-арная подгруппа единиц, её центр и она сама имеют чётные порядки.

**9.1.7. Предложение** [124, 125]. Если $a, b, c$ – различные единицы тернарной группы $<A, [\ ]>$, то $<\{a, b, c, [abc]\}, [\ ]>$ – тернарная подгруппа четвёртого порядка в $<E(A), [\ ]>$.

*Доказательство.* Положим

$$d = [abc].\ D = \{a, b, c, d\}.$$

Ясно, что $D \subseteq E(A)$. Если $x, y \in D$, то

$$[xxx] = x \in D, \quad [xxy] = y \in D.$$

Рассмотрим теперь все возможные случаи, когда под знаком тернарной операции стоят три различных элемента из D:

$$[abc] = d \in D;$$

$$[abd] = [ab[abc]] = [aa[bbc]] = [aac] = c \in D;$$

$$[acd] = [ac[abc]] = [aa[ccb]] = [aab] = b \in D;$$

$$[bcd] = [bc[abc]] = [bb[cca]] = [bba] = a \in D.$$

Так как множество D замкнуто относительно тернарной операции [ ], то, согласно замечанию 9.1.4, $< D, [\ ] >$ – тернарная подгруппа в $< E(A), [\ ] >$. Предположим, что

$$[abc] = a,$$

откуда

$$[a[abc]c] = [aac], \ [aabcc] = c, \ b = c,$$

что противоречит условию. Следовательно, $[abc] \neq a$. Аналогично доказывается, что $[abc] \neq b$, $[abc] \neq c$, то есть все элементы в D различны ∎

Представляет интерес следующее

**9.1.8. Предложение** [124, 125]. Если n-арная группа $< A, [\ ] >$ является производной от m-арной группы $< A, \lfloor\ \rfloor >$, то любая единица в $< A, \lfloor\ \rfloor >$ является единицей в $< A, [\ ] >$.

*Доказательство.* По условию,

$$[a_1 a_2 \ldots a_n] = \lfloor \lfloor \ldots \lfloor \lfloor a_1 \ldots a_m \rfloor a_{m+1} \ldots a_{2m-1} \rfloor \ldots \rfloor a_{(k-1)(m-1)+2} \ldots a_{k(m-1)+1} \rfloor,$$

для любых $a_1, \ldots, a_n \in A$, где $k > 1$. Если e – единица m-арной группы $< A, [\ ] >$, a – произвольный элемент из A, то

$$[a\underbrace{e\ldots e}_{n-1}] = \lfloor\lfloor\ldots\lfloor\lfloor a\underbrace{e\ldots e}_{m-1}\rfloor\underbrace{e\ldots e}_{m-1}\rfloor\ldots\rfloor\underbrace{e\ldots e}_{m-1}\rfloor =$$

$$= \lfloor\lfloor\ldots\lfloor a\underbrace{e\ldots e}_{m-1}\rfloor\ldots\rfloor\underbrace{e\ldots e}_{m-1}\rfloor = \ldots = \lfloor a\underbrace{e\ldots e}_{m-1}\rfloor = a,$$

$$[ea\underbrace{e\ldots e}_{n-2}] = \lfloor\lfloor\ldots\lfloor\lfloor ea\underbrace{e\ldots e}_{m-2}\rfloor\underbrace{e\ldots e}_{m-1}\rfloor\ldots\rfloor\underbrace{e\ldots e}_{m-1}\rfloor =$$

$$= \lfloor\lfloor\ldots\lfloor a\underbrace{e\ldots e}_{m-1}\rfloor\ldots\rfloor\underbrace{e\ldots e}_{m-1}\rfloor = \ldots = \lfloor a\underbrace{e\ldots e}_{m-1}\rfloor = a,$$

$$[\underbrace{e\ldots e}_{n-1}a] = \lfloor\lfloor\ldots\lfloor\lfloor \underbrace{e\ldots e}_{m}\rfloor\underbrace{e\ldots e}_{m-1}\rfloor\ldots\rfloor\underbrace{e\ldots e}_{m-2}a\rfloor =$$

$$= \lfloor\lfloor\ldots\lfloor \underbrace{e\ldots e}_{m}\rfloor\ldots\rfloor\underbrace{e\ldots e}_{m-2}a\rfloor = \ldots = \lfloor\underbrace{e\ldots e}_{m-1}a\rfloor = a.$$

Следовательно,

$$[a\underbrace{e\ldots e}_{n-1}] = [ea\underbrace{e\ldots e}_{n-2}] = [\underbrace{e\ldots e}_{n-1}a],$$

и по предложению 1.2.10 e является единицей n-арной группы $< A, [\ ] >$. ∎

Пусть $< A, [\ ] >$ – n-арная группа, обладающая по крайней мере одной единицей. Для всякой единицы e n-арной группы $< A, [\ ] >$ определим на A бинарную операцию

$$x \circledcirc y = [x\underbrace{e\ldots e}_{n-2}y].$$

Так как $\underbrace{e\ldots e}_{n-2}$ – обратная последовательность для элемента e, то при a = e операция ⓔ совпадает с операцией ⓐ из теоремы Глускина-Хоссу (§1.5). Легко проверяется (см. например, предложение 1.5.2), что $< A, \circledcirc > $ – группа с единицей e.

**9.1.9. Предложение.** Если e и ε – две единицы n-арной группы $< A, [\ ] >$, то группы $< A, \circledcirc >$ и $< A, \circledcirc >$ – изоморфны.

***Доказательство.*** Определим преобразование $\alpha$ множества A по правилу

$$\alpha: x \to x \circledcirc \varepsilon.$$

Ясно, что $\alpha$ – подстановка на A и, кроме того,

$$(x \circledcirc y)^\alpha = [x\underbrace{e\ldots e}_{n-2}y] \circledcirc \varepsilon = [[x\underbrace{e\ldots e}_{n-2}y]\underbrace{e\ldots e}_{n-2}\varepsilon] =$$

$$= [x\underbrace{e\ldots e}_{n-2}\underbrace{\varepsilon\ldots\varepsilon}_{n-1}y\underbrace{e\ldots e}_{n-2}\varepsilon] = [[x\underbrace{e\ldots e}_{n-2}\varepsilon]\underbrace{\varepsilon\ldots\varepsilon}_{n-2}y\underbrace{e\ldots e}_{n-2}\varepsilon]] =$$

$$= (x \circledcirc \varepsilon) \circledcirc (y \circledcirc \varepsilon) = x^\alpha \circledcirc y^\alpha,$$

то есть

$$(x \circledcirc y)^\alpha = x^\alpha \circledcirc y^\alpha.$$

Следовательно, $\varphi$ - изоморфизм групп $< A, \circledcirc >$ и $< A, \circledcirc >$. ∎

Заметим, что изоморфность групп $< A, \circledcirc >$ и $< A, \circledcirc >$ без явного указания изоморфизма вытекает из следствия 1.6.2.

**9.1.10. Замечание.** Если n-арная группа $< A, [\ ] >$ является производной от группы $< A, * >$ с единицей e, то

$$x * y = x * \underbrace{e * \ldots e}_{n-2} * y = [x\underbrace{e\ldots e}_{n-2}y] = x \circledcirc y,$$

то есть операции $*$ и $\circledcirc$ совпадают.

Обозначим через $Z(A, [\ ])$ – центр n-арной группы $< A, [\ ] >$, а через $Z(A, \circledcirc )$ центр группы $< A, \circledcirc >$.

**9.1.11. Лемма** [124, 125]. Для любой единицы e n-арной группы $< A, [\ ] >$ верно

$$Z(A, [\ ]) = Z(A, \circledcirc ).$$

***Доказательство.*** Если $z \in Z(A, [\ ])$, то для любого $x \in A$ последовательности $zx$ и $xz$ эквивалентны, откуда

$$[zx\underbrace{e\ldots e}_{n-2}] = [xz\underbrace{e\ldots e}_{n-2}]. \qquad (1)$$

Так как e – единица в < A, [ ] >, то xe ~ ex, ze ~ ez. Поэтому из последнего равенства получаем

$$[z\underbrace{e\ldots e}_{n-2}x] = [x\underbrace{e\ldots e}_{n-2}z], \qquad (2)$$

$$z \circlede x = x \circlede z, \qquad (3)$$

то есть $z \in Z(A, \circlede)$. Следовательно,

$$Z(A, [\,]) \subseteq Z(A, \circlede).$$

Если теперь $z \in Z(A, \circlede)$, то последовательно выполняются (3), (2) и (1), то есть последовательности zx и xz эквивалентны, и поэтому $z \in Z(A, [\,])$. Следовательно,

$$Z(A, \circlede) \subseteq Z(A, [\,]).$$

Таким образом, доказано равенство

$$Z(A, [\,]) = Z(A, \circlede). \qquad \blacksquare$$

**9.1.12. Замечание.** Лемма 9.1.11 позволяет ввести обозначение

$$Z(A) = Z(A, [\,]) = Z(A, \circlede)$$

общее для группы < A, $\circlede$ > и n-арной группы < A, [ ] > – производной от неё.

**9.1.13. Следствие** [124, 125]. Для любых единиц e и ε n-арной группы < A, [ ] > центры групп < A, $\circlede$ > и < A, $\circledvarepsilon$ > совпадают.

Так как единица n-арной группы < A, [ ] > лежит в её центре, то справедливо

**9.1.14. Следствие.** Для любых единиц e и ε n-арной группы < A, [ ] > верно

$$\varepsilon \in Z(A, \circlede).$$

Следствие 9.1.14 можно сформулировать иначе: для любого $e \in E(A)$ верно

$$E(A) \subseteq Z(A, \copyright).$$

**9.1.15. Лемма** [124, 125]. Если $e$ и $\varepsilon$ – единицы n-арной группы $< A, [\ ] >$, то

$$\underbrace{\varepsilon \copyright \varepsilon \copyright \ldots \varepsilon}_{n-1} = e.$$

*Доказательство.* Так как $e\varepsilon \sim \varepsilon e$, то

$$\underbrace{\varepsilon \copyright \varepsilon \copyright \ldots \varepsilon}_{n-1} \copyright e = [\underbrace{\underbrace{\varepsilon e \ldots e}_{n-2}\underbrace{\varepsilon e \ldots e}_{n-2}\ldots \underbrace{\varepsilon e \ldots e}_{n-2} e}_{n-1}] =$$

$$= [\underbrace{\varepsilon \ldots \varepsilon}_{n-1}[\underbrace{e \quad \ldots \quad e}_{(n-2)(n-1)+1}]] = [\underbrace{\varepsilon \ldots \varepsilon}_{n-1} e] = e,$$

то есть

$$\underbrace{\varepsilon \copyright \varepsilon \copyright \ldots \varepsilon}_{n-1} \copyright e = e.$$

откуда

$$\underbrace{\varepsilon \copyright \varepsilon \copyright \ldots \varepsilon}_{n-1} = e. \qquad \blacksquare$$

**9.1.16. Лемма** [124, 125]. Пусть $A$ – группа с единицей $e$, $z$ – элемент из её центра, удовлетворяющий условию $z^{n-1} = e$ ($n \geq 2$); $< A, [\ ] >$ – n-арная группа, производная от группы $A$. Тогда $z$ – единица n-арной группы $< A, [\ ] >$.

*Доказательство.* Так как

$$z^{n-1} = e, z \in Z(A),$$

то

$$[\underbrace{z \ldots z}_{i-1} a \underbrace{z \ldots z}_{n-i}] = z^{i-1} a z^{n-i} = z^{n-1} a = ea = a,$$

то есть

$$[\underbrace{z\ldots z}_{i-1}a\underbrace{z\ldots z}_{n-i}] = a$$

для любого a ∈ A и любого i = 1, 2, …, n. Следовательно, z – единица n-арной группы < A, [ ] >. ∎

Следующая теорема является основной в данном параграфе.

**9.1.17. Теорема** [124, 125]. Если e – единица n-арной группы < A, [ ] >, то

$$E(A) = \{z \in Z(A, ⓒ) \mid \underbrace{z ⓒ z ⓒ \ldots z}_{n-1} = e\}.$$

*Доказательство.* Если ε ∈ E(A), то по следствию 9.1.14, ε ∈ Z(A, ⓒ), а по лемме 9.1.15

$$\underbrace{ε ⓒ ε ⓒ \ldots ε}_{n-1} = e.$$

Следовательно,

$$E(A) \subseteq \{z \in Z(A, ⓒ) \mid \underbrace{z ⓒ z ⓒ \ldots z}_{n-1} = e\}.$$

Пусть теперь

$$z \in Z(A, ⓒ), \underbrace{z ⓒ z ⓒ \ldots z}_{n-1} = e.$$

Так как n-арная группа < A, [ ] > является производной от группы < A, ⓒ >, то по лемме 9.1.16

$$\{z \in Z(A, ⓒ) \mid \underbrace{z ⓒ z ⓒ \ldots z}_{n-1} = e\} \subseteq E(A).$$

Из доказанных включений получается требуемое равенство. ∎

Доказанная теорема может быть сформулирована иначе.

**9.1.18. Теорема** [124, 125]. Если $< A, [\ ] > -$ n-арная группа, производная от группы A, то

$$E(A) = \{z \in Z(A) \mid z^{n-1} = e\},$$

где e – единица группы A.

## §9.2. СЛЕДСТВИЯ ИЗ ОСНОВНОЙ ТЕОРЕМЫ И ПРИМЕРЫ

**9.2.1. Следствие** [124, 125]. Если $< A, [\ ] > -$ n-арная группа, производная от абелевой группы A, то

$$E(A) = \{a \in A \mid a^{n-1} = e\},$$

где e – единица группы A.

**9.2.2. Пример** [124, 125]. Пусть $C^*$ – мультипликативная группа комплексных чисел. Так как она абелева, то, согласно следствию 9.2.1, n-арная подгруппа единиц n-арной группы $< C^*, [\ ] >$, производной от группы $C^*$, имеет вид

$$E(C^*) = \{z \in C^* \mid z^{n-1} = 1\} = \{\cos\frac{2k\pi}{n-1} + i\sin\frac{2k\pi}{n-1} \mid k = 0,1,\ldots,n-2\},$$

то есть имеет ровно n – 1 единиц.

**9.2.3. Пример** [124, 125]. Пусть $< C_{p^\infty}, [\ ] > -$ $(p^k+1)$-арная группа, где k = 1, 2, …, производная от квазициклической группы $C_{p^\infty}$. Так как $C_{p^\infty}$ – абелева и содержит единственную циклическую подгруппу $Z_{p^k}$ порядка $p^k$, то

$$E(C_{p^\infty}) = \{a \in C_{p^\infty} \mid a^{p^k} = 1\} = Z_{p^k},$$

то есть в $< C_{p^\infty}, [\ ] >$ ровно $p^k$ – единиц.

**9.2.4. Пример** [124, 125]. Пусть $< Z_k, [\ ] > -$ (m + 1)-арная группа, производная от циклической группы $Z_k$ порядка k, где m делит k. Так как $Z_k$ – абелева и содержит единственную циклическую подгруппу $Z_m$ порядка m, то

$$E(Z_k) = \{a \in Z_k \mid a^m = 1\} = Z_m,$$

то есть в $< Z_k, [\ ] >$ ровно m единиц.

**9.2.5. Следствие** [124, 125]. Пусть подгруппа B центра $Z(A)$ группы A имеет конечный период n – 1, $< A, [\ ] >$ – n-арная группа, производная от группы A. Тогда $< B, [\ ] >$ – n-арная подгруппа в $< A, [\ ] >$, все элементы которой являются единицами.

*Доказательство.* Так как B – подгруппа группы A, то $< B, [\ ] >$ – n-арная подгруппа n-арной группы $< A, [\ ] >$. Если b – произвольный элемент из B, то

$$b^{n-1} = e, b \in Z(A),$$

и по теореме 9.1.18 $b \in E(A)$, то есть $B \subseteq E(A)$. ∎

**9.2.6. Следствие** [124, 125]. Пусть периодическая часть B центра $Z(A)$ группы A имеет конечный период n – 1, $< A, [\ ] >$ – n-арная группа, производная от группы A. Тогда $< B, [\ ] >$ – n-арная подгруппа в $< A, [\ ] >$, совпадающая с n-арной подгруппой $< E(A), [\ ] >$.

*Доказательство.* Так как $Z(A)$ – абелева, то B – подгруппа группы A ([109], с. 90), откуда и из следствия 9.2.5 вытекает, что $< B, [\ ] >$ – n-арная подгруппа в $< A, [\ ] >$, причем $B \subseteq E(A)$.

Если z – единица n-арной группы $< A, [\ ] >$, то по теореме 9.1.18 $z^{n-1} = e$, и поэтому $z \in B$, то есть $E(A) \subseteq B$. Из доказанных включений получаем равенство $B = E(A)$. ∎

Так как порядок и период циклической группы совпадают, то из следствия 9.2.5 вытекает

**9.2.7. Следствие** [124, 125]. Пусть $Z_{n-1}$ – циклическая подгруппа порядка n – 1 группы A, лежащая в её центре,

< A, [ ] > – n-арная группа, производная от группы A. Тогда < $Z_{n-1}$, [ ] > – n-арная подгруппа в < A, [ ] >, все элементы которой являются единицами.

**9.2.8. Следствие** [124, 125]. В n-арной группе < $Z_{n-1}$, [ ] >, производной от циклической группы $Z_{n-1}$ порядка n – 1, все элементы являются единицами.

**9.2.9. Пример.** Так как $A_3$ – циклическая группа третьего порядка, то, согласно следствию 9.2.8, в 4-арной группе < $A_3$, [ ] >, производной от группы $A_3$, все три элемента являются единицами.

**9.2.10. Следствие** [124, 125]. Пусть центр Z(A) группы A имеет период n – 1, < A, [ ] > – n-арная группа, производная от группы A. Тогда Z(A) = E(A).

**9.2.11. Пример.** Пусть < R, [ ] > – тернарная группа, производная от группы кватернионов,

$$R = \{1, a, a^2, a^3, b, ba, ba^2, ba^3\}.$$

Так как Z(R) = {1, $a^2$} – циклическая группа второго порядка, то, согласно следствию 9.2.10, тернарная группа < R, [ ] > имеет ровно две единицы 1 и $a^2$.

**9.2.12. Пример.** Так как при четных n, центр $Z(D_n)$ диэдральной группы $D_n$ включает помимо единицы e ещё и поворот на угол π, то, согласно следствию 9.2.10, тернарная группа < $D_n$, [ ] > при четном n имеет ровно две единицы.

**9.2.13. Следствие** [124, 125]. В n-арной группе < A, [ ] >, производной от абелевой группы периода n – 1, все элементы являются единицами.

**9.2.14. Следствие** [124, 125]. n-Арная группа < A, [ ] >, производная от группы A, имеющей тривиальный центр Z(A) = {e}, обладает единственной единицей e.

**9.2.15. Пример** [124, 125]. Так как $Z(S_n)$ = {e} при n ≥ 3, то, согласно следствию 9.2.14, производная m-арная группа < $S_n$, [ ] > при n ≥ 3, m ≥ 3 обладает единственной единицей e.

**9.2.16. Пример** [124, 125]. Так как $Z(A_n) = \{e\}$ при $n \geq 4$, то, согласно следствию 9.2.14, производная m-арная группа $< A_n, [\ ] >$ при $n \geq 4$, $m \geq 3$ обладает единственной единицей.

**9.2.17. Пример** [124, 125]. Так как $Z(D_n) = \{e\}$ при нечетном n, то, согласно следствию 9.2.14, производная m-арная группа $< D_n, [\ ] >$ при нечетном n и $m \geq 3$ обладает единственной единицей e.

Укажем примеры n-арных групп, не обладающих единицами.

**9.2.18. Пример** [124, 125]. Пусть $< T_n, [\ ] >$ – тернарная группа всех нечётных подстановок степени n (пример 1.1.8), которая является тернарной подгруппой тернарной группы $< S_n, [\ ] >$, производной от группы $S_n$. Предположим, что $\varepsilon$ – единица тернарной группы $< T_n, [\ ] >$. Так как все транспозиции лежат в $T_n$, то $\varepsilon\alpha \sim \alpha\varepsilon$ для любой транспозиции $\alpha$, то есть

$$[\varepsilon\alpha\beta] = [\alpha\varepsilon\beta]$$

для любого $\beta \in T_n$, откуда $\varepsilon\alpha = \alpha\varepsilon$. Всякую подстановку $\sigma \in S_n$ можно представить в виде произведения $\sigma = \alpha_1\alpha_2\ldots\alpha_k$ транспозиций. Поэтому

$$\varepsilon\sigma = \varepsilon\alpha_1\alpha_2\ldots\alpha_k = \alpha_1\alpha_2\ldots\alpha_k\varepsilon = \sigma\varepsilon,$$

то есть $\varepsilon\sigma = \sigma\varepsilon$, откуда $\varepsilon\sigma \sim \sigma\varepsilon$. Следовательно, $\varepsilon$ – единица тернарной группы $< S_n, [\ ] >$, то есть $\varepsilon \in E(S_n) = \{e\}$, что невозможно, так как $\varepsilon$ – нечетная подстановка. Таким образом, $E(T_n) = \varnothing$.

**9.2.19. Пример** [124, 125]. Пусть n – нечетное, $< B_n, [\ ] >$ – тернарная группа отражений правильного n-угольника (пример 1.2.6), которая является тернарной подгруппой тернарной группы $< D_n, [\ ] >$, производной от диэдральной группы $D_n$. Предположим, что $\varepsilon$ – единица тернарной группы $< B_n, [\ ] >$, то есть $\varepsilon\psi \sim \psi\varepsilon$, что равносильно $\varepsilon\psi = \psi\varepsilon$ для любого отражения $\psi \in B_n$. Так как всякий поворот $\varphi \in D_n$ можно представить в виде произведения $\varphi = \psi_1\psi_2$ двух отражений, то

$$\varepsilon\varphi = \varepsilon\psi_1\psi_2 = \psi_1\varepsilon\psi_2 = \psi_1\psi_2\varepsilon = \varphi\varepsilon,$$

то есть $\varepsilon\varphi = \varphi\varepsilon$, откуда $\varepsilon\varphi \sim \varphi\varepsilon$. Следовательно, $\varepsilon$ – единица тернарной группы $< D_n, [\ ] >$, то есть $\varepsilon \in E(D_n) = \{e\}$, что невозможно, так как $\varepsilon$ – отражение. Таким образом, $E(B_n) = \varnothing$.

Так как группы $S_n$ и $D_n$ являются обертывающими группами для тернарных групп из примеров 9.2.18 и 9.2.19 соответственно, то последние два примера обобщаются следующим образом.

**9.2.20. Предложение** [124, 125]. Пусть G– обертывающая группа n-арной группы $< A, [\ ] >$, $< G, [\ ] >$ – n-арная группа, производная от группы G, и пусть $E(A) \cap E(G) = \varnothing$. Тогда $E(A) = \varnothing$.

*Доказательство.* Предположим, что $\varepsilon$ – единица n-арной группы $< A, [\ ] >$. Так как $\varepsilon$ лежит в центре n-арной группы $< A, [\ ] >$, то

$$[\varepsilon a d_1 \ldots d_{n-2}] = [a \varepsilon d_1 \ldots d_{n-2}],$$

$$\varepsilon a d_1 \ldots d_{n-2} = a \varepsilon d_1 \ldots d_{n-2},$$

для любых $a, d_1, \ldots, d_{n-2} \in A$, откуда

$$\varepsilon a = a \varepsilon \qquad (1)$$

Всякий элемент $g \in G$ можно представить в виде произведения $g = a_1 \ldots a_k$ ($k \geq 1$) элементов из A. Используя (1), получим

$$\varepsilon g = \varepsilon a_1 \ldots a_k = a_1 \ldots a_k \varepsilon = g \varepsilon,$$

то есть $\varepsilon g = g \varepsilon$. Следовательно, $\varepsilon$ лежит в центре группы G.

С другой стороны, учитывая, что $\varepsilon$ – единица в $< A, [\ ] >$, получим

$$[\underbrace{\varepsilon \ldots \varepsilon}_{n}] = \varepsilon,\ \varepsilon^n = \varepsilon,\ \varepsilon^{n-1} = e,$$

где e – единица группы G, а значит и n-арной группы $< G, [\ ] >$. Так как $\varepsilon \in Z(G)$, $\varepsilon^{n-1} = e$, то, согласно теореме 9.1.18, $\varepsilon \in E(G)$, что невозможно, так как $E(A) \cap E(G) = \varnothing$. Таким образом, в $< A, [\ ] >$ нет единиц. ∎

Отметим, что конструированию периодических n-арных групп без единицы посвящен §6 монографии С.А. Русакова [4].

## §9.3. n-АРНАЯ ПОДГРУППА ЕДИНИЦ ПРЯМОГО ПРОИЗВЕДЕНИЯ n-АРНЫХ ГРУПП

Пусть $\{\mathscr{A}_i = <A_i, [\ ]_i>\,|\,i \in I\}$ – непустое семейство n-арных групп,
$$A = \prod_{i \in I} A_i = \{a\colon I \to \bigcup_{i \in I} A_i \,|\, a(i) \in A_i\} \ -$$
декартово произведение их носителей, [ ] – n-арная операция, определенная на A покомпонентно:

$$[a_1 a_2 \ldots a_n](i) = [a_1(i) a_2(i) \ldots a_n(i)]_i,\ i \in I$$

для всех $a_1, a_2, \ldots, a_n \in A$. Тогда алгебра $\mathscr{A} = <A, [\ ]> = \prod_{i \in I} \mathscr{A}_i$ называется декартовым или прямым произведением n-арных групп $\mathscr{A}_i$ ($i \in I$), и, согласно предложению 5.1 из [4], является n-арной группой.

**9.3.1. Теорема** [124]. Если $\mathscr{A} = \prod \mathscr{A}_i$ – декартово произведение n-арных групп $\mathscr{A}_i$ ($i \in I$), то справедливы следующие утверждения:

1) $e \in E(\mathscr{A})$ тогда и только тогда, когда $e(i) \in E(\mathscr{A}_i)$ для любого $i \in I$;

2) $E(\mathscr{A}) \neq \varnothing$ тогда и только тогда, когда $E(\mathscr{A}_i) \neq \varnothing$ для любого $i \in I$;

3) если $E(\mathscr{A}_i) \neq \varnothing$ для любого $i \in I$, то $E(\mathscr{A}) = \prod E(\mathscr{A}_i)$.

*Доказательство.* 1) Пусть $e \in E(\mathscr{A})$, $a_i$ – произвольный элемент из $A_i$ и зафиксируем элемент $a \in A$ такой, что $a(i) = a_i$. Так как $e \in E(\mathscr{A})$, то

$$[\underbrace{e\ldots e}_{k-1} a \underbrace{e\ldots e}_{n-k}] = a,\ k = 1, 2, \ldots, n, \tag{1}$$

откуда

$$[\underbrace{e\ldots e}_{k-1} a \underbrace{e\ldots e}_{n-k}](i) = a(i), \qquad (2)$$

$$[\underbrace{e(i)\ldots e(i)}_{k-1} a(i) \underbrace{e(i)\ldots e(i)}_{n-k}]_i = a(i), \qquad (3)$$

$$[\underbrace{e(i)\ldots e(i)}_{k-1} a_i \underbrace{e(i)\ldots e(i)}_{n-k}]_i = a_i.$$

Так как элемент $a_i$ выбран в n-арной группе $\mathscr{A}_i = <A_i, [\ ]>$ произвольно, то из последнего равенства вытекает $e(i) \in E(\mathscr{A}_i)$.

Если теперь a – произвольный элемент из A, и элемент $e \in A$ удовлетворяет условию $e(i) \in E(\mathscr{A}_i)$ для любого $i \in I$, то выполняется равенство (3) для любых $k = 1, 2, \ldots, n$ и $i \in I$, откуда последовательно получаем (2) и (1). Следовательно, $e \in E(\mathscr{A})$.

2) Если $e \in E(\mathscr{A})$, то по доказанному в 1), $e(i) \in E(\mathscr{A}_i)$ для любого $i \in I$, то есть $E(\mathscr{A}_i) \neq \varnothing$.

Если теперь $e_i \in E(\mathscr{A}_i)$, то, положив $e(i) = e_i$ для любого $i \in I$, получим элемент $e \in A$, который, согласно 1), лежит в $E(\mathscr{A})$, то есть $E(\mathscr{A}) \neq \varnothing$.

3) Если $e \in E(\mathscr{A})$, то согласно 1), $e(i) \in E(\mathscr{A}_i)$ для любого $i \in I$, то есть

$$e: I \to \bigcup_{i \in I} E(\mathscr{A}_i).$$

Следовательно, $e \in \prod E(\mathscr{A}_i)$, и, значит, $E(\mathscr{A}) \subseteq \prod E(\mathscr{A}_i)$.

Пусть теперь $e \in \prod E(\mathscr{A}_i)$, то есть

$$e: I \to \bigcup_{i \in I} E(\mathscr{A}_i), \ e(i) \in E(\mathscr{A}_i).$$

Так как $E(\mathscr{A}_i) \subseteq A_i$, то

$$e: I \to \bigcup_{i \in I} A_i, \ e(i) \in A_i,$$

то есть e ∈ ∏$A_i$ = A. Кроме того, согласно 1), e ∈ E(𝒜), откуда, в силу произвольного выбора e, получаем ∏E($𝒜_i$) ⊆ E(𝒜). Из доказанных включений вытекает равенство

$$E(𝒜) = \prod E(𝒜_i).$$  ∎

Все рассмотренные выше примеры n-арных групп таковы, что может сложиться впечатление, будто число единиц в n-арной группе не превосходит её арности. На самом деле это не так. Убедиться в этом позволяет теорема 9.3.1, применённая к n-арным группам с уже известными n-арными подгруппами единиц.

**9.3.2. Пример** [124]. Пусть ℛ = < R, [ ] > – тернарная группа, производная от группы из примера 9.2.11. Так как E(ℛ) = {1, $a^2$}, то согласно 3) теоремы 9.3.1,

$$E(ℛ × ℛ) = E(ℛ) × E(ℛ) = \{(1, 1), (1, a^2), (a^2, 1), (a^2, a^2)\}.$$

Следовательно, тернарная группа ℛ × ℛ имеет ровно четыре единицы. Вообще, тернарная группа $\prod_{i=1}^{k} 𝒞$ имеет $2^k$ единиц. Если же I – бесконечное множество, то тернарная группа $\prod_{i \in I} 𝒞$ имеет бесконечную n-арную подгруппу единиц.

Следующая теорема может быть получена в качестве следствия из теоремы 9.3.1, однако мы проведём её прямое доказательство.

**9.3.3. Теорема** [124]. Если $𝒜_i$ = < $A_i$, [ ]$_i$ > – n-арные группы, производные от групп $A_i$ (i ∈ I), то n-арная группа $𝒜 = \prod 𝒜_i$ = < A, [ ] > является производной от группы ∏$A_i$.

*Доказательство.* Обозначим через ⊙ операцию в группе $A_i$ (i ∈ I), а через ○ – операцию в группе ∏$A_i$. Так как

$$[a_1 a_2 … a_n](i) = [a_1(i) a_2(i) … a_n(i)]_i =$$

$$= a_1(i) ⊙ a_2(i) ⊙ … ⊙ a_n(i) = (a_1 ○ a_2 ○ … ○ a_n)(i),$$

то есть

$$[a_1a_2\ldots a_n](i) = (a_1 \circ a_2 \circ \ldots \circ a_n)(i),$$

то

$$[a_1a_2\ldots a_n] = a_1 \circ a_2 \circ \ldots \circ a_n$$

для любых $a_1, a_2, \ldots, a_n \in A$. Следовательно, n-арная группа $\mathcal{A}$ является производной от группы $\prod A_i$. ∎

## §9.4. МНОЖЕСТВО I(A)

Для всякой n-арной группы $\mathcal{A} = <A, [\ ]>$ обозначим через $I(\mathcal{A}) = I(A)$ множество всех её идемпотентов.

Ясно, что $E(A) \subseteq I(A)$. Если же $<A, [\ ]>$ – абелева, то $E(A) = I(A)$, и, согласно теореме 9.1.1, $<I(A), [\ ]>$ – n-арная подгруппа абелевой n-арной группы $<A, [\ ]>$.

**9.4.1. Пример.** Так как полиадические группы $<C^*, [\ ]>$, $<C_{p^\infty}, [\ ]>$ и $<Z_k, [\ ]>$ из примеров 9.2.2, 9.2.3 и 9.2.4 абелевы, то

$$I(C^*) = \{\cos\frac{2k\pi}{n-1} + i\sin\frac{2k\pi}{n-1} \mid k = 0, 1, \ldots, n–2\},$$

$$I(C_{p^\infty}) = Z_{p^k}, \quad I(Z_k) = Z_m.$$

Приведем примеры, показывающие, что множество всех идемпотентов n-арной группы в общем случае не образует в ней n-арную подгруппу.

**9.4.2. Пример.** Пусть $<S_3, [\ ]>$ - тернарная группа, производная от симметрической группы $S_3$. Легко проверяется, что

$$I(S_3) = \{e, \alpha, \beta, \gamma\},$$

где e – тождественная подстановка; $\alpha$, $\beta$ и $\gamma$ – нечётные подстановки. Так как $|S_3| = 6$, $|I(S_3)| = 4$, и 4 не делит 6, то множество $I(S_3)$ не образует в $<S_3, [\ ]>$ – тернарную подгруппу.

Пример 9.4.2 обобщается следующим образом.

**9.4.3. Пример** [126]. Пусть $< D_n, [\ ] >$ – тернарная группа, производная от диэдральной группы $D_n$. Все отражения правильного n-угольника образуют в $< D_n, [\ ] >$ тернарную подгруппу $< B_n, [\ ] >$, все элементы которой, как установлено в [26], являются идемпотентами. Следовательно, $B_n \subseteq I(D_n)$.

Предположим, что поворот $\varphi \in C_n$ является идемпотентом в $< D_n, [\ ] >$, то есть

$$[\varphi\varphi\varphi] = \varphi,$$

откуда

$$\varphi\varphi\varphi = \varphi;\ \varphi^2 = e,$$

где e – тождественный поворот. Так как $C_n$ – циклическая группа порядка n, то при нечётном n только единица удовлетворяет последнему равенству, а при чётном n в $C_n$ имеется ещё один поворот $\varphi \neq e$ такой, что $\varphi^2 = e$. Таким образом,

$$I(D_{2k+1}) = \{e\} \bigcup B_{2k+1},$$

$$I(D_{2k}) = \{e, \varphi\} \bigcup B_{2k},\ \varphi^2 = e.$$

Так как $|D_n| = 2n$, $|I(D_n)| = n + 1$ при нечётном n, и $|I(D_n)| = n + 2$ при чётном n, то множество всех идемпотентов тернарной группы $< D_n, [\ ] >$ при $n > 2$ не образует в ней тернарную подгруппу.

Следующие леммы понадобятся нам для нахождения условий, при которых множество всех идемпотентов n-арной группы образует в ней n-арную подгруппу.

**9.4.4. Лемма** [126]. Если $\varepsilon$ – идемпотент n-арной группы $< A, [\ ] >$, то $\varepsilon^\alpha$ – идемпотент n-арной группы $< B, [\ ] >$ для любого гомоморфизма

$$\alpha : < A, [\ ] > \to < B = A^\alpha, [\ ] >.$$

В частности, если $I(A) \neq \varnothing$, то $I^\alpha(A) \subseteq I(A)$ для любого эндоморфизма n-арной группы $< A, [\ ] >$.

*Доказательство.* Так как

$$[\underbrace{\varepsilon^\alpha \ldots \varepsilon^\alpha}_{n}] = [\underbrace{\varepsilon \ldots \varepsilon}_{n}]^\alpha = \varepsilon^\alpha,$$

то $\varepsilon^\alpha$ – идемпотент в $< B, [\ ] >$. ∎

**9.4.5. Лемма** [126]. Если $I(A) \neq \varnothing$, то

$$[\underbrace{\varepsilon \ldots \varepsilon}_{i-1} I(A) \underbrace{\varepsilon \ldots \varepsilon}_{n-i}] = I(A)$$

для любого $\varepsilon \in I(A)$ и любого $i = 1, \ldots, n$.

***Доказательство.*** Если $i = 1$ или $i = n$, то доказывать нечего. Так как для любого $i = 2, \ldots, n-1$ последовательности

$$\underbrace{\varepsilon \ldots \varepsilon}_{i-1}, \quad \underbrace{\varepsilon \ldots \varepsilon}_{n-i}$$

являются взаимно обратными, то преобразования

$$a \to [\underbrace{\varepsilon \ldots \varepsilon}_{i-1} a \underbrace{\varepsilon \ldots \varepsilon}_{n-i}], \quad a \to [\underbrace{\varepsilon \ldots \varepsilon}_{n-i} a \underbrace{\varepsilon \ldots \varepsilon}_{i-1}]$$

будут автоморфизмами n-арной группы $< A, [\ ] >$. Тогда по лемме 9.4.4,

$$[\underbrace{\varepsilon \ldots \varepsilon}_{i-1} I(A) \underbrace{\varepsilon \ldots \varepsilon}_{n-i}] \subseteq I(A), \qquad (1)$$

$$[\underbrace{\varepsilon \ldots \varepsilon}_{n-i} I(A) \underbrace{\varepsilon \ldots \varepsilon}_{i-1}] \subseteq I(A). \qquad (2)$$

Из (2) получаем

$$[\underbrace{\varepsilon \ldots \varepsilon}_{i-1} [\underbrace{\varepsilon \ldots \varepsilon}_{n-i} I(A) \underbrace{\varepsilon \ldots \varepsilon}_{i-1}] \underbrace{\varepsilon \ldots \varepsilon}_{n-i}] \subseteq [\underbrace{\varepsilon \ldots \varepsilon}_{i-1} I(A) \underbrace{\varepsilon \ldots \varepsilon}_{n-i}],$$

$$[\underbrace{\varepsilon \ldots \varepsilon}_{n-1} I(A) \underbrace{\varepsilon \ldots \varepsilon}_{n-1}] \subseteq [\underbrace{\varepsilon \ldots \varepsilon}_{i-1} I(A) \underbrace{\varepsilon \ldots \varepsilon}_{n-i}],$$

$$I(A) \subseteq [\underbrace{\varepsilon \ldots \varepsilon}_{i-1} I(A) \underbrace{\varepsilon \ldots \varepsilon}_{n-i}]. \qquad (3)$$

Из включений (1) и (3) следует равенство

$$[\underbrace{\varepsilon\ldots\varepsilon}_{i-1}I(A)\underbrace{\varepsilon\ldots\varepsilon}_{n-i}] = I(A). \qquad \blacksquare$$

**9.4.6. Следствие.** Для любого $\varepsilon \in I(A)$ верно

$$[\varepsilon\, I(A)\underbrace{\varepsilon\ldots\varepsilon}_{n-2}] = I(A), \quad [\underbrace{\varepsilon\ldots\varepsilon}_{n-2}I(A)\,\varepsilon] = I(A).$$

**9.4.7. Теорема** [126, 127]. Если $a \in I(A)$, то $< I(A), [\ ] > -$ n-арная подгруппа в $< A, [\ ] >$ тогда и только тогда, когда $< I(A), @ > -$ подгруппа в $< A, @ >$.

*Доказательство.* Необходимость является следствием предложения 1.5.2.

Достаточность вытекает из первого равенства следствия 9.4.6 и следствия 2.2.7. $\blacksquare$

**9.4.8. Предложение** [126, 127]. $E(A) = I(A)\bigcap Z(A)$.

*Доказательство.* Если $E(A) = \varnothing$, то $I(A)\bigcap Z(A) = \varnothing$.

Пусть теперь $E(A) \neq \varnothing$. Включение $E(A) \subseteq I(A)$ очевидно, а включение $E(A) \subseteq Z(A)$ доказано в теореме 9.1.1. Следовательно,

$$E(A) \subseteq I(A)\bigcap Z(A).$$

Если теперь $e \in I(A)\bigcap Z(A)$, то для любого $x \in A$

$$[\underbrace{e\ldots e}_{n-1}x] = x,$$

и последовательности $ex$ и $xe$ эквивалентны. Поэтому

$$[\underbrace{e\ldots e}_{i-1}x\underbrace{e\ldots e}_{n-i}] = x$$

для любого $i = 1, 2, \ldots, n$, то есть $e \in E(A)$. Следовательно,

$$I(A)\bigcap Z(A) \subseteq E(A).$$

Из доказанных включений следует требуемое равенство. ∎

Прежде чем сформулировать следующие два предложения, напомним, что квадратная матрица называется подстановочной над полем F, если в каждой строке и каждом столбце этой матрицы ровно один элемент совпадает с единицей поля F, а все остальные элементы равны нулю этого же поля.

**9.4.9. Предложение** [126]. Если $<GL_n(F), [\ ]>$ – $(n! + 1)$-арная группа, производная от полной линейной группы $GL_n(F)$, то множество всех подстановочных матриц над F образует в $<GL_n(F), [\ ]>$ $(n! + 1)$-арную подгруппу, лежащую в $I(GL_n(F))$.

*Доказательство.* Известно, что существует мономорфизм $\varphi$ симметрической группы $S_n$ в группу $GL_n(F)$, при котором образ группы $S_n$ совпадает с множеством $\mathcal{B}$ всех подстановочных матриц над F. Так как $\mathcal{B}$ – подгруппа в $GL_n(F)$, то $<\mathcal{B}, [\ ]>$ – $(n! + 1)$-арная подгруппа в $<GL_n(F), [\ ]>$, а так как $|\mathcal{B}| = |S_n| = n!$, то $B^{n!} = E_n$ для любой матрицы $B \in \mathcal{B}$, откуда

$$[\underbrace{B\ldots B}_{n!+1}] = B^{n!}B = E_nB = B.$$

Следовательно, B – идемпотент в $<GL_n(F), [\ ]>$, и значит, $\mathcal{B} \subseteq I(GL_n(F))$. ∎

Дословно повторяя доказательство предыдущего предложения, а также учитывая то, что образом знакопеременной группы $A_n$ при мономорфизме $\varphi$ является множество всех подстановочных матриц с определителем равным единице, получим

**9.4.10. Предложение** [126]. Если $<SL_n(F), [\ ]>$ – $((n!/2)+1)$-арная группа, производная от специальной линейной группы $SL_n(F)$, то множество всех подстановочных матриц над F с определителем равным единице образует в $<SL_n(F), [\ ]>$ $((n!/2) + 1)$-арную подгруппу, лежащую в $I(SL_n(F))$.

Ясно, что для всякого идемпотента a n-арной группы $<A, [\ ]>$ преобразование $\beta$ из теоремы Глускина-Хоссу (§1.5) имеет вид

$$\beta : x \to [ax\underbrace{a\ldots a}_{n-2}].$$

**9.4.11 Теорема** [126, 127]. Если a – идемпотент n-арной группы $<A, [\ ]>$, то

$$I(A) = \{b \in A \mid b \,@\, b^\beta \,@\, \ldots \, b^{\beta^{n-2}} = a\}. \qquad (1)$$

*Доказательство.* Преобразование $\beta$ и элемент

$$a = [\underbrace{a\ldots a}_{n}],$$

согласно теореме Глускина-Хоссу, удовлетворяют условиям

$$a \,@\, x = x^{\beta^{n-1}} \,@\, a, \quad x \in A;$$

$$[x_1 \ldots x_n] = x_1 \,@\, x_2^\beta \,@\, \ldots \, x_{n-1}^{\beta^{n-2}} \,@\, x_n^{\beta^{n-1}} \,@\, a, \; x_1, \ldots, x_n \in A.$$

Так как a – единица группы $<A, @>$, то из первого равенства получаем

$$x = x^{\beta^{n-1}}. \qquad (2)$$

Подставляя (2) во второе равенство, получаем

$$[x_1 \ldots x_n] = x_1 \,@\, x_2^\beta \,@\, \ldots \, x_{n-1}^{\beta^{n-2}} \,@\, x_n. \qquad (3)$$

Положив в (3) $x_1 = \ldots = x_n = b \in I(A)$, получим

$$[\underbrace{b\ldots b}_{n}] = b \circleda b^{\beta} \circledа \ldots b^{\beta^{n-2}} \circledа b,$$

откуда, с учётом $b \in I(A)$, вытекает

$$b = b \circledа b^{\beta} \circledа \ldots b^{\beta^{n-2}} \circledа b, \qquad (4)$$

$$a = b \circledа b^{\beta} \circledа \ldots b^{\beta^{n-2}}, \qquad (5)$$

$$I(A) \subseteq \{b \in A \mid b \circledа b^{\beta} \circledа \ldots b^{\beta^{n-2}} = a\}. \qquad (6)$$

Пусть теперь $b$ – произвольный элемент из $A$, удовлетворяющий (5). Тогда верно (4), откуда, учитывая (3), получаем

$$b = [\underbrace{b\ldots b}_{n}],$$

то есть $b \in I(A)$. Следовательно,

$$\{b \in A \mid b \circledа b^{\beta} \circledа \ldots b^{\beta^{n-2}} = a\} \subseteq I(A). \qquad (7)$$

Из (6) и (7) вытекает (1). ∎

**9.4.12. Замечание.** Из доказательства теоремы 9.4.11 видно, что верно двойственное к (1) равенство

$$I(A) = \{b \in A \mid b^{\beta} \circledа \ldots b^{\beta^{n-2}} \circledа b = a\}.$$

**9.14.13. Следствие** [126, 127]. Если $< A, [\ ] >$ – n-арная группа, производная от группы $A$, то

$$I(A) = \{b \in A \mid b^{n-1} = e\},$$

где $e$ – единица группы $A$.

*Доказательство.* Так как $e$ единица группы $A$, то

$$x \circledе y = [x\underbrace{e\ldots e}_{n-2}y] = x\underbrace{e\ldots e}_{n-2}y = xy,$$

то есть $x \circledе y = xy$. Ясно, что $e \in E(A) \subseteq I(A)$. Поэтому

$$x^\beta = [ex\underbrace{e\ldots e}_{n-2}] = x,$$

то есть β – тождественное преобразование. Следовательно,

$$I(A) = \{b \in A \mid b \odot b^\beta \odot \ldots b^{\beta^{n-2}} = e\} = \{b \in A \mid b^{n-1} = e\}. \quad \blacksquare$$

**9.4.14. Пример** [126, 127]. Пусть $<\Gamma, [\,]>$ – тернарная группа, производная от группы Γ всех движений плоскости. Согласно следствию 9.4.13,

$$I(\Gamma) = \{\gamma \in \Gamma \mid \gamma^2 = e\}.$$

Поэтому все элементы множества I(Γ) исчерпываются всеми отражениями, всеми поворотами на угол π и всеми произведениями отражения на сдвиг, перпендикулярный оси отражения.

**9.4.15. Пример** [126, 127]. Пусть $<GL_n(q), [\,]>$ – q-арная группа (q > 2), производная от полной линейной группы $GL_n(q)$ над конечным полем F(q) из q элементов. Так как $|Z(GL_n(q))| = q - 1$, то по теореме 9.1.18

$$E(GL_n(q)) = Z(GL_n(q)) = \{\alpha E_n \mid \alpha \in F(q), \alpha \neq 0\},$$

то есть $E(GL_n(q))$ состоит из всех скалярных матриц, содержащихся в $GL_n(q)$.

Покажем, что в $<GL_n(q), [\,]>$ имеются идемпотенты, не являющиеся единицами. Для этого обозначим через $\mathcal{D}$ множество всех диагональных матриц с ненулевыми элементами из F(q) на главной диагонали, то есть

$$\mathcal{D} = \{\mathrm{diag}(\alpha_1, \ldots, \alpha_n) \mid \alpha_i \in F(q), \alpha_i \neq 0\}.$$

Так как

$$(\mathrm{diag}(\alpha_1, \ldots, \alpha_n))^{q-1} = \mathrm{diag}(\alpha_1^{q-1}, \ldots, \alpha_n^{q-1}) = \mathrm{diag}(1, \ldots, 1) = E_n,$$

то по следствию 9.4.13, $\mathcal{D} \subseteq I(GL_n(q))$.

Ясно, что $E(GL_n(q)) \subseteq \mathcal{D}$, а при q > 1

$$E(GL_n(q)) \subset \mathcal{D}.$$

Отметим, что в $< GL_n(q), [\ ] >$ могут быть идемпотенты, не являющиеся диагональными матрицами. Например, в 4-арной группе $< GL_3(4), [\ ] >$ элементы

$$\begin{pmatrix} 0 & 0 & 1 \\ 1 & 0 & 0 \\ 0 & 1 & 0 \end{pmatrix}, \begin{pmatrix} 0 & 1 & 0 \\ 0 & 0 & 1 \\ 1 & 0 & 0 \end{pmatrix}$$

являются идемпотентами.

**9.4.16. Пример** [126, 127]. Пусть $< SL_n(q), [\ ] >$ – m-арная группа, где $m = 1 + (n, q - 1)$, производная от специальной линейной группы $SL_n(q)$ над $F(q)$. Так как

$$|Z(SL_n(q))| = (n, q - 1) = m - 1,$$

то по теореме 9.1.18

$$E(SL_n(q)) = Z(SL_n(q)) = \{\alpha E_n \mid \alpha \in F(q), \alpha^n = 1\}.$$

В частности, если q – нечётное, то

$$E(SL_2(q)) = \left\{ \begin{pmatrix} 1 & 0 \\ 0 & 1 \end{pmatrix}, \begin{pmatrix} -1 & 0 \\ 0 & -1 \end{pmatrix} \right\}.$$

Так как определители матриц из примера 9.4.15 равны единице поля, то обе они являются идемпотентами в $< SL_3(4), [\ ] >$, но не являются в ней единицами. Поэтому в общем случае множества $E(SL_n(q))$ и $I(SL_n(q))$ не совпадают.

Отметим, что если $q - 1 < n$, то идемпотентами в $< SL_n(q), [\ ] >$, отличными от единиц, являются также диагональные матрицы, у которых на главной диагонали $q - 1$ элементов равны $\alpha \in F(q)$, где $\alpha \neq 0$, $\alpha \neq 1$, а остальные элементы на главной диагонали равны единице поля.

**9.4.17. Следствие** [126, 127]. Пусть A – группа с тождеством $a^{n-1} = e$, $< A, [\ ] >$ – n-арная группа, производная от группы A. Тогда $I(A) = A$.

**9.4.18. Следствие** [126, 127]. Пусть A – конечная группа порядка $n - 1$, $< A, [\ ] >$ – n-арная группа, производная от группы A. Тогда $I(A) = A$.

## §9.5. ИДЕМПОТЕНТНЫЕ n-АРНЫЕ ГРУППЫ

**9.5.1. Определение.** n-Арная группа $<A, [\ ]>$ называется *идемпотентной*, если каждый её элемент является идемпотентом, то есть $I(A) = A$.

Так как в идемпотентной n-арной группе $<A, [\ ]>$ $Z(A) \subseteq A = I(A)$, то из предложения 9.4.8 вытекает

**9.5.2. Следствие.** Если $<A, [\ ]>$ идемпотентная n-арная группа, то

$$E(A) = Z(A).$$

Примерами идемпотентных n-арных групп могут служить n-арные группы из следствий 9.4.17 и 9.4.18.

Среди идемпотентных n-арных групп с непустой n-арной подгруппой единиц представляют интерес те из них, для которых $|E(A)| = 1$, то есть содержащие только одну единицу.

**9.5.3. Предложение** [126]. n-Арная группа $<A, [\ ]>$, производная от группы $A$ с тождеством $a^{n-1} = e$ и тривиальным центром, является идемпотентной с единственной единицей.

*Доказательство.* Согласно следствию 9.4.17, $<A, [\ ]>$ – идемпотентная n-арная группа, а по следствию 9.2.14 $<A, [\ ]>$ обладает единственной единицей. ∎

**9.5.4. Следствие** [126]. n-Арная группа $<A, [\ ]>$, производная от конечной группы порядка $n-1$ с тривиальным центром, является идемпотентной с единственной единицей.

Так как простая неабелева группа имеет тривиальный центр, то имеет место

**9.5.5. Следствие** [126]. Если $A$ – конечная простая неабелева группа, то $(|A|+1)$-арная группа $<A, [\ ]>$, производная от группы $A$, является идемпотентной с единственной единицей.

**9.5.6. Пример** [126]. Так как $Z(S_n) = \{e\}$ при $n \geq 3$, то, согласно следствию 9.5.4, производная $(n! + 1)$-арная группа $< S_n, [\ ] >$ при $n \geq 3$ является идемпотентной с единственной единицей.

**9.5.7. Пример** [126]. Так как $Z(A_n) = \{e\}$ при $n \geq 4$, то, согласно следствию 9.5.4, производная $(\frac{n!}{2} + 1)$-арная группа $< A_n, [\ ] >$ при $n \geq 4$ является идемпотентной с единственной единицей.

**9.5.8. Пример** [126]. Так как $Z(D_n) = \{e\}$ при нечётном $n$, то, согласно следствию 9.5.4, производная $(2n + 1)$-арная группа $< D_n, [\ ] >$ при нечётном $n$ является идемпотентной с единственной единицей.

**9.5.9. Пример** [126]. Так как $Z(PGL_n(F(q))) = \{e\}$, то, согласно следствию 9.5.4, производная $(m+1)$-арная группа $< PGL_n(F(q)), [\ ] >$, где

$$m = |PGL_n(F(q))| = \frac{1}{q-1}\prod_{i=0}^{n-1}(q^n - q^i),$$

является идемпотентной с единственной единицей.

**9.5.10. Пример** [126]. Так как $Z(PSL_n(F(q))) = \{e\}$, то, согласно следствию 9.5.4, производная $(m+1)$-арная группа $< PSL_n(F(q)), [\ ] >$, где

$$m = |PSL_n(F(q))| = \frac{1}{(q-1)(n,\ q-1)}\prod_{i=0}^{n-1}(q^n - q^i),$$

является идемпотентной с единственной единицей.

Приведём пример бесконечной идемпотентной полиадической группы с единственной единицей.

**9.5.11. Пример** [126]. Пусть $B(m, n)$ – бесконечная группа с $m$ порождающими и тождеством $x^n = e$ ($m \geq 1$, $n \geq 665$, $n$ – нечётное), построенная Адяном [128]. По теореме 3.4 [128] центр группы $B(m, n)$ при указанных $m$ и $n$ тривиален. Поэтому согласно предложению 9.5.3, $(n + 1)$-арная группа $< B(m, n), [\ ] >$, производная от группы $B(m, n)$, при нечётном $n \geq 665$, $m > 1$ является бесконечной идемпотентной с единственной единицей.

Ясно, что все идемпотентные n-арные группы составляют многообразие, которое в многообразии всех n-арных групп выделяется тождеством

$$[\underbrace{x\ldots x}_{n}] = x.$$

Поэтому декартово произведение идемпотентных n-арных групп является идемпотентной n-арной группой. В действительности имеет место более общее утверждение.

**9.5.12. Теорема** [126]. Если $\mathscr{A} = \prod \mathscr{A}_j$ – декартово произведение n-арных групп $\mathscr{A}_j$ ($j \in J$), то справедливы следующие утверждения:

1) $e \in I(\mathscr{A})$ тогда и только тогда, когда $e(j) \in I(\mathscr{A}_j)$ для любого $j \in J$;

2) $I(\mathscr{A}) \neq \varnothing$ тогда и только тогда, когда $I(\mathscr{A}_j) \neq \varnothing$ для любого $j \in J$;

3) если $< I(\mathscr{A}_j), [\,]_j >$ – непустая n-арная подгруппа n-арной группы $\mathscr{A}_j = < A_j, [\,]_j >$ для любого $j \in J$, то $< I(\mathscr{A}), [\,] >$ – n-арная подгруппа n-арной группы $\mathscr{A} = < \prod A_j, [\,] >$, причём

$$I(\mathscr{A}) = \prod I(\mathscr{A}_j).$$

*Доказательство*. 1) Если $e \in I(\mathscr{A})$, то

$$[\underbrace{e\ldots e}_{n}] = e, \qquad (1)$$

откуда

$$[\underbrace{e\ldots e}_{n}](j) = e(j), \qquad (2)$$

$$[\underbrace{e(j)\ldots e(j)}_{n}]_j = e(j), \qquad (3)$$

то есть $e(j) \in I(\mathscr{A}_j)$ для любого $j \in J$.

Если теперь элемент $e \in A = \prod A_j$ удовлетворяет условию $e(j) \in I(A_j)$ для любого $j \in J$, то выполняется равенство (3), откуда последовательно получаем (2) и (1). Следовательно, $e \in I(\mathcal{A})$.

2) Если $e \in I(\mathcal{A})$, то по доказанному в 1), $e(j) \in I(\mathcal{A}_j)$ для любого $j \in J$, то есть $I(\mathcal{A}_j) \neq \varnothing$.

Если теперь $e_j \in I(\mathcal{A}_j)$, то, положив $e(j) = e_j$ для любого $j \in J$, получим элемент $e \in A$, который, согласно 1), лежит в $I(\mathcal{A})$, то есть $I(\mathcal{A}) \neq \varnothing$.

3) По предложению 5.1 [4], $\prod I(\mathcal{A}_j)$ — n-арная группа, которая, очевидно, является n-арной подгруппой в $<\prod A_j, [\,]>$.

Если $e \in I(\mathcal{A})$, то согласно 1), $e(j) \in I(\mathcal{A}_j)$ для любого $j \in J$, то есть

$$e : J \to \bigcup_{j \in J} I(\mathcal{A}_j).$$

Следовательно, $e \in \prod I(\mathcal{A}_j)$, и, значит, $I(\mathcal{A}) \subseteq \prod I(\mathcal{A}_j)$.

Пусть теперь $e \in \prod I(\mathcal{A}_j)$, то есть

$$e : J \to \bigcup_{j \in J} I(\mathcal{A}_j), \ e(j) \in I(\mathcal{A}_j).$$

Так как $I(\mathcal{A}_j) \subseteq A_j$, то

$$e : J \to \bigcup_{j \in J} \mathcal{A}_j, \ e(j) \in \mathcal{A}_j,$$

то есть $e \in \prod A_j = A$. Кроме того, согласно 1), $e \in I(\mathcal{A})$, откуда, в силу произвольного выбора $e$, получаем $\prod I(\mathcal{A}_j) \subseteq I(\mathcal{A})$. Из доказанных включений вытекает требуемое равенство, из которого, в свою очередь, следует, что $< I(\mathcal{A}), [\,] >$ — n-арная подгруппа в $<\prod A_j, [\,]>$. ∎

**9.5.13. Следствие** [126]. Пусть

$$\mathcal{A}_1 = <A_1, [\ ]_1>, \ldots, \mathcal{A}_k = <A_k, [\ ]_k>$$

– n-арные группы. Тогда справедливы следующие утверждения:

1) $(e_1, \ldots, e_k) \in I(\mathcal{A}_1 \times \ldots \times \mathcal{A}_k)$ тогда и только тогда, когда $e_j \in I(\mathcal{A}_j)$ для любого $j = 1, \ldots, k$;

2) $I(\mathcal{A}_1 \times \ldots \times \mathcal{A}_k) \neq \varnothing$ тогда и только тогда, когда $I(\mathcal{A}_j) \neq \varnothing$ для любого $j = 1, \ldots, k$;

3) если $<I(\mathcal{A}_j), [\ ]_j>$ – n-арная подгруппа n-арной группы $\mathcal{A}_j$ для любого $j = 1, \ldots, k$, то $<I(\mathcal{A}_1 \times \ldots \times \mathcal{A}_k), [\ ]>$ – n-арная подгруппа n-арной группы $\mathcal{A}_1 \times \ldots \times \mathcal{A}_k$, причём

$$I(\mathcal{A}_1 \times \ldots \times \mathcal{A}_k) = I(\mathcal{A}_1) \times \ldots \times I(\mathcal{A}_k).$$

**9.5.14. Лемма**. Для всякого идемпотента a n-арной группы $<A, [\ ]>$ отображение

$$\beta : x \to [ax\underbrace{a\ldots a}_{n-2}], \ x \in A$$

является автоморфизмом группы $<A, \circledast>$, удовлетворяющим условию $\beta^{n-1} = \varepsilon$, где $\varepsilon$ – тождественный автоморфизм $<A, \circledast>$.

*Доказательство*. По предложению 1.5.4, $\beta$ – автоморфизм группы $<A, \circledast>$, который, согласно теореме Глускина-Хоссу (седьмое тождество следствия 1.5.6), удовлетворяет условию

$$x^{\beta^{n-1}} \circledast d = d \circledast x, \ x \in A,$$

где

$$d = [\underbrace{a\ldots a}_{n}].$$

Так как a – идемпотент в $< A, [\ ] >$, то d совпадает с элементом a, являющимся единицей в $< A, @ >$. Поэтому из предпоследнего равенства следует

$$x^{\beta^{n-1}} = x.$$ ∎

**9.5.15. Теорема** [126, 127]. Справедливы следующие утверждения:

1) если $< A, [\ ] >$ – идемпотентная n-арная группа, то для любого $a \in A$ группа $< A, @ >$ обладает автоморфизмом $\beta_a$ таким, что $\beta_a^{n-1}$ – тождественное отображение, и

$$b @ b^{\beta_a} @ \ldots b^{\beta_a^{n-2}} = a \qquad (1)$$

для любого $b \in A$;

2) если $< A, \circ >$ – группа, обладающая автоморфизмом $\beta$ таким, что $\beta^{n-1}$ – тождественное отображение и

$$b \circ b^{\beta} \circ \ldots b^{\beta^{n-2}} = e \qquad (2)$$

для любого $b \in A$, где e – единица группы $< A, \circ >$, то $< A, [\ ] >$ – идемпотентная n-арная группа, где

$$[x_1 x_2 \ldots x_n] = x_1 \circ x_2^{\beta} \circ \ldots x_{n-1}^{\beta^{n-2}} \circ x_n, \qquad (3)$$

причём операции $@$ и $\circ$ совпадают.

*Доказательство*. 1) Если $< A, [\ ] >$ – идемпотентная n-арная группа, то для всякого $a \in A$ положим

$$\beta_a : x \to [ax\underbrace{a\ldots a}_{n-2}]$$

Так как $I(A) = A$, то по теореме 9.4.11 для любого $b \in A$ верно (1), а по лемме 9.5.14, $\beta_a^{n-1}$ – тождественное отображение.

2) Так как единица e группы $< A, \circ >$ и её автоморфизм $\beta$ удовлетворяют условиям

$$e^\beta = e, \quad x^{\beta^{n-1}} \circ e = e \circ x,$$

то по обратной теореме Глускина-Хоссу $< A, [\ ] >$ – n-арная группа, где

$$[x_1 x_2 \ldots x_n] = x_1 \circ x_2^\beta \circ \ldots x_{n-1}^{\beta^{n-2}} \circ x_n^{\beta^{n-1}} \circ e,$$

откуда, учитывая тождественность $\beta^{n-1}$, получаем (3).

Подставляя в (3) $x_1 = \ldots = x_n = b \in A$, получим

$$[\underbrace{b \ldots b}_{n}] = b \circ b^\beta \circ \ldots b^{\beta^{n-2}} \circ b,$$

откуда в силу (2) следует

$$[\underbrace{b \ldots b}_{n}] = b,$$

то есть $b$ – идемпотент в $< A, [\ ] >$. Так как элемент $b$ выбран в $A$ произвольно, то $< A, [\ ] >$ – идемпотентная n-арная группа.

Полагая в (3) $x_2 = \ldots = x_{n-1} = e$ и учитывая $e^\beta = e$, получим

$$[x_1 \underbrace{e \ldots e}_{n-2} x_n] = x_1 \circ e^\beta \circ \ldots e^{\beta^{n-2}} \circ x_n = x_1 \circ x_n,$$

то есть

$$[x_1 \underbrace{e \ldots e}_{n-2} x_n] = x_1 \circ x_n$$

для любых $x_1, x_n \in A$. По доказанному все элементы в $< A, [\ ] >$, в том числе и $e$, идемпотенты. Поэтому

$$[x_1 \underbrace{e \ldots e}_{n-2} x_n] = x_1 \circledcirc x_n,$$

откуда, учитывая предыдущее равенство, получаем

$$x_1 \circledcirc x_n = x_1 \circ x_n,$$

то есть операции ⊚ и ∘ совпадают. ∎

Напомним, что автоморфизм β группы A называется расщепляющим, если для любого b ∈ A верно

$$bb^{\beta}\ldots b^{\beta^{n-1}} = 1,$$

где n – порядок автоморфизма β.

Утверждение 2) теоремы 9.5.15 позволяет сформулировать

**9.5.16. Следствие** [126, 127]. Если группа < A, ∘ > допускает расщепляющий автоморфизм β порядка n – 1, то < A, [ ] > – идемпотентная n-арная группа с n-арной операцией (3).

**9.5.17. Следствие** [126, 127]. Если конечная группа < A, ∘ > допускает автоморфизм β порядка n – 1 без неподвижных точек, то < A, [ ] > – идемпотентная n-арная группа с n-арной операцией (3).

*Доказательство*. По теореме V.8.9 из [129], β – расщепляющий автоморфизм порядка n – 1, и применяется предыдущее следствие. ∎

**9.5.18. Следствие** [126, 127]. Если < A, [ ] > – конечная идемпотентная n-арная группа, где n – 1 простое, то для любого a ∈ A группа < A, ⊚ > – нильпотентна.

*Доказательство*. Согласно 1) теоремы 9.5.15, $\beta_a$ – расщепляющий автоморфизм простого порядка n – 1 группы < A, ⊚ >. Тогда по теореме V.8.13 [129] группа < A, ⊚ > – нильпотентна. Следствие доказано.

**9.5.19. Замечание** [126, 127]. Придавая n в следствии 9.5.18 конкретные значения, можно получить большое число новых следствий. В частности, группа < A, ⊚ > будет нильпотентной для любой конечной идемпотентной n-арной группы < A, [ ] >, где n = 3, 4, 6, 8.

**9.5.20. Замечание.** Так как группа $<A, \text{@}>$ изоморфна соответствующей группе Поста $<A_0, *>$ (предложение 1.6.1), то следствие 9.5.18 останется верным, если в нем группу $<A, \text{@}>$ заменить группой $<A_o, *>$.

**9.5.21. Лемма** [126, 127]. Если $\beta$ – автоморфизм n-арной группы $<A, [\ ]>$, и $a^\beta = a$ для некоторого $a \in A$, то $\beta$ – автоморфизм группы $<A, \text{@}>$.

*Доказательство.* Пусть $\varphi$ – автоморфизм n-арной группы $<A, [\ ]>$, оставляющий неподвижным элемент $a \in A$, и пусть $a_1, \ldots, a_{n-2}$ – обратная последовательность для $a$, то есть

$$[aa_1 \ldots a_{n-2}x] = x$$

для любого $x \in A$, откуда, учитывая, что $\varphi$ – автоморфизм $<A, [\ ]>$, а также условие $a^\varphi = a$, получим

$$[aa_1 \ldots a_{n-2}x]^\varphi = x^\varphi,$$

$$[a^\varphi a_1^\varphi \ldots a_{n-2}^\varphi x^\varphi] = x^\varphi,$$

$$[a a_1^\varphi \ldots a_{n-2}^\varphi x^\varphi] = x^\varphi.$$

Из последнего равенства видно, что $a_1^\varphi \ldots a_{n-2}^\varphi$ – обратная последовательность для $a$. Следовательно, последовательности $a_1 \ldots a_{n-2}$ и $a_1^\varphi \ldots a_{n-2}^\varphi$ эквивалентны. Поэтому

$$(x \text{ @ } y)^\varphi = [xa_1 \ldots a_{n-2}y]^\varphi = [x^\varphi a_1^\varphi \ldots a_{n-2}^\varphi y^\varphi] =$$

$$= [x^\varphi a_1 \ldots a_{n-2} y^\varphi] = x^\varphi \text{ @ } y^\varphi,$$

то есть

$$(x \text{ @ } y)^\varphi = x^\varphi \text{ @ } y^\varphi,$$

и $\varphi$ – автоморфизм группы $<A, \text{@}>$. ∎

**9.5.22. Теорема** [126, 127]. Если $<A, [\ ]>$ – конечная n-арная группа, допускающая автоморфизм порядка $n-1$, ос-тавляющий неподвижным единственный элемент, то на A можно определить n-арную операцию $\lfloor \ \rfloor$ так, что $< A, \lfloor \ \rfloor >$ идемпотентная n-арная группа

*Доказательство*. Пусть $\beta$ – автоморфизм порядка $n-1$ n-арной группы $< A, [\ ] >$, оставляющий неподвижным единственный элемент $a \in A$. По лемме 9.5.21, $\beta$ – автоморфизм группы $< A, @ >$, единственным неподвижным элементом которого, согласно условию теоремы, является единица $a$. Поэтому по следствию 9.5.17, $< A, \lfloor \ \rfloor >$ – идемпотентная n-арная группа с n-арной операцией

$$\lfloor x_1 \ldots x_n \rfloor = x_1 \ @ \ x_2^{\beta} \ @ \ \ldots \ x_{n-1}^{\beta^{n-2}} \ @ \ x_n. \qquad \blacksquare$$

**9.5.23. Предложение** [126, 127]. Если n-арная группа $< A, [\ ] >$ допускает автоморфизм $\varphi$, оставляющий неподвижным единственный элемент $a$, то $a$ – идемпотент.

*Доказательство*. Так как $a^{\varphi} = a$, то

$$[\underbrace{a \ldots a}_{n}]^{\varphi} = [\underbrace{a^{\varphi} \ldots a^{\varphi}}_{n}] = [\underbrace{a \ldots a}_{n}],$$

то есть $[\underbrace{a \ldots a}_{n}]$ – неподвижный элемент относительно автоморфизма $\varphi$. А так как $a$ – единственный неподвижный относительно $\varphi$ элемент, то $[\underbrace{a \ldots a}_{n}] = a$, то есть $a$ – идемпотент. $\blacksquare$

**9.5.24. Следствие** [126, 127]. n-Арная группа без идемпотентов не допускает автоморфизм с единственным неподвижным элементом.

В связи с предложением 9.5.23 интересен следующий вопрос: если n-арная группа допускает автоморфизм, оставляющий неподвижным единственный элемент, то будет ли этот элемент единицей?

Отрицательный ответ на этот вопрос даёт следующий

**9.5.25. Пример** [126,127]. Пусть $< T_3 = \{\alpha, \beta, \gamma\}, [\ ] >$ тернарная группа нечётных подстановок с n-арной операцией

$$[xyz] = xyz,$$

где $\alpha = (12)$, $\beta = (13)$, $\gamma = (23)$. В $< T_3, [\ ] >$ все элементы – идемпотенты, не являющиеся единицами. Рассмотрим внутренний автоморфизм

$$\varphi_\alpha : x \to [\alpha x \alpha].$$

Так как $\varphi_\alpha(\alpha) = \alpha$, $\varphi_\alpha(\beta) = \gamma$, $\varphi_\alpha(\gamma) = \beta$, то $\varphi_\alpha$ – автоморфизм тернарной группы $< T_3, [\ ] >$, оставляющий неподвижным единственный элемент $\alpha$, не являющийся единицей. Аналогично определяются автоморфизмы $\varphi_\beta$ и $\varphi_\gamma$, которые также оставляют неподвижными по одному элементу $\beta$ и $\gamma$, которые не являются единицами.

Пример 9.5.25 можно обобщить. Пусть $< B_n, [\ ] >$ – тернарная группа отражений правильного n-угольника, в которой, как известно [26], все элементы являются идемпотентами и нет единиц. Определим на $B_n$ для любого $b \in B_n$ преобразование

$$\varphi_b : x \to [bxb],$$

которое очевидно является внутренним автоморфизмом $< B_n, [\ ] >$.

**9.5.26. Предложение** [126, 127]. Для любого $b \in B_n$ автоморфизм $\varphi_b$ при нечётном n оставляет неподвижным единственный элемент b, а при чётном n неподвижными относительно $\varphi_b$ остаются ровно два элемента.

*Доказательство*. Множество $B_n$ можно представить в виде

$$B_n = \{b, bc, bc^2, \ldots, bc^{n-1}\} = \{b_1, b_2, \ldots, b_n\},$$

где c – образующий поворот циклической группы поворотов $C_n$. Ясно, что $\varphi_b(b) = b = b_1$. Если же $i = 2, \ldots, n$, то

$$\varphi_b(b_i) = [bbc^{i-1}b] = bbc^{i-1}b = c^{i-1}b =$$

$$= bc^{n+1-i} = bc^{n+2-i-1} = b_{n+2-i},$$

то есть

$$\varphi_b(b_i) = b_{n+2-i}, \; i = 2, \ldots, n.$$

Предположим, что $b_i = b_{n+2-i}$, откуда $i = n + 2 - i$, $i = \dfrac{n+2}{2}$.

Последнее равенство возможно только при чётном n. И в этом случае неподвижными относительно $\varphi_b$ остаются только элементы b и $b_{\frac{n+2}{2}}$. ∎

### §9.6. СИЛОВСКОЕ СТРОЕНИЕ ИДЕМПОТЕНТНОЙ n-АРНОЙ ГРУППЫ

**9.6.1 Предложение** [12, 130]. Пусть $< B, [\;] >$ – полуинвариантная n-арная подгруппа n-арной группы $< A, [\;] >$, a – идемпотент из А. Тогда

$$< B_a = [\underbrace{B\ldots B}_{n-1}a], [\;] >$$

– полуинвариантная n-арная подгруппа в $< A, [\;] >$ такая, что $A/B = A/B_a$. Если же $a \notin B$, то $B_a \neq B$.

*Доказательство.* Если

$$u_1 = [b_1^{'}\ldots b_{n-1}^{'}a], \ldots, u_n = [b_1^{(n)}\ldots b_{n-1}^{(n)}a]$$

– произвольные элементы из $B_a$, то, учитывая полуинвариантность $< B, [\;] >$ в $< A, [\;] >$, идемпотентность a, а также то, что $< B, [\;] >$ – n-арная подгруппа в $< A, [\;] >$, получим

$$[u_1\ldots u_n] = [[b_1^{'}\ldots b_{n-1}^{'}a]\ldots[b_1^{(n)}\ldots b_{n-1}^{(n)}a]] =$$

$$= [b_1^{*}\ldots b_{n-1}^{*}[\underbrace{a\ldots a}_{n}]] = [b_1^{*}\ldots b_{n-1}^{*}a] \in B_a,$$

где $b_1^*, \ldots, b_{n-1}^* \in B$. Следовательно, $<B_a, [\ ]>$ – n-арная полугруппа.

Рассмотрим в $<B_a, [\ ]>$ уравнение

$$[u_1 \ldots u_{n-1} x] = u, \qquad (1)$$

где $u = [b_1 \ldots b_{n-1} a] \in B_a$, которое имеет решение $x = u_n \in A$, то есть.

$$[[b_1' \ldots b_{n-1}' a] \ldots [b_1^{(n-1)} \ldots b_{n-1}^{(n-1)} a] u_n] = [b_1 \ldots b_{n-1} a].$$

Из последнего равенства, снова, учитывая полуинвариантность $<B, [\ ]>$ в $<A, [\ ]>$, идемпотентность $a$, а также то, что $<B, [\ ]>$ – n-арная подгруппа в $<A, [\ ]>$, получим

$$[\widetilde{b}_1 \ldots \widetilde{b}_{n-1} [\underbrace{a \ldots a}_{n-1} u_n]] = [b_1 \ldots b_{n-1} a],$$

$$[\widetilde{b}_1 \ldots \widetilde{b}_{n-1} u_n] = [b_1 \ldots b_{n-1} a],$$

где $\widetilde{b}_1, \ldots, \widetilde{b}_{n-1} \in B$. Если $\beta$ – обратная последовательность для последовательности $\widetilde{b}_1 \ldots \widetilde{b}_{n-1}$, то из последнего равенства получаем

$$[\beta[\widetilde{b}_1 \ldots \widetilde{b}_{n-1} u_n]] = [\beta[b_1 \ldots b_{n-1} a]],$$

$$u_n = [[\beta b_1] b_2 \ldots b_{n-1} a] = [b b_2 \ldots b_{n-1} a] \in B_a,$$

где $b = [\beta b_1] \in B$. Следовательно, в $<B_a, [\ ]>$ разрешимо уравнение (1).

Аналогично доказывается разрешимость в $<B_a, [\ ]>$ уравнения

$$[y u_2 \ldots u_n] = u.$$

Таким образом, установлено, что $<B_a, [\ ]>$ – n-арная подгруппа в $<A, [\ ]>$.

Учитывая полуинвариантность $<B, [\ ]>$ в $<A, [\ ]>$, а также идемпотентность a, получим

$$[\underbrace{B_a \ldots B_a}_{n-1} c] = [\underbrace{[\underbrace{B\ldots Ba}_{n-1}]\ldots[\underbrace{B\ldots Ba}_{n-1}]}_{n-1} c] =$$

$$= [\underbrace{B \ldots B}_{(n-1)(n-1)} \underbrace{a\ldots a}_{n-1} c] = [\underbrace{B\ldots B}_{n-1} c] = [c\underbrace{B\ldots B}_{n-1}] =$$

$$= [c\underbrace{B \ldots B}_{(n-1)(n-1)} \underbrace{a\ldots a}_{n-1}] = [c\underbrace{[\underbrace{B\ldots Ba}_{n-1}]\ldots[\underbrace{B\ldots Ba}_{n-1}]}_{n-1}] = [c\underbrace{B_a \ldots B_a}_{n-1}],$$

то есть

$$[\underbrace{B_a \ldots B_a}_{n-1} c] = [c\underbrace{B_a \ldots B_a}_{n-1}]$$

для любого $c \in A$. Следовательно, $<B_a, [\ ]>$ – полуинвариантна в $<A, [\ ]>$.

Так как для любого $c \in A$ верно

$$[\underbrace{B\ldots B}_{n-1} c] = [\underbrace{B\ldots B}_{n-2}[\underbrace{B \ldots B}_{(n-1)(n-2)+1}][\underbrace{a\ldots a}_{n-1} c]] =$$

$$= [\underbrace{\underbrace{B\ldots B}_{n-1}\ldots\underbrace{B\ldots B}_{n-1}}_{n-1} \underbrace{a\ldots a}_{n-1} c] =$$

$$= [\underbrace{[\underbrace{B\ldots Ba}_{n-1}]\ldots[\underbrace{B\ldots Ba}_{n-1}]}_{n-1} c] = [\underbrace{B_a \ldots B_a}_{n-1} c],$$

то есть

$$[\underbrace{B\ldots B}_{n-1} c] = [\underbrace{B_a \ldots B_a}_{n-1} c],$$

то

$$A/B = A/B_a.$$

Так как $a \notin B$, $a \in B_a$, то $B \neq B_a$, точнее $B \cap B_a = \varnothing$. ∎

**9.6.2. Лемма** [126, 130]. Если $< B, [\ ] >$ – полуинвариантная n-арная подгруппа идемпотентной n-арной группы $< A, [\ ] >$, то $< A, [\ ] >$ является объединением непересекающихся полуинвариантных в $< A, [\ ] >$ n-арных подгрупп, мощность каждой из которых совпадает с мощностью множества B.

*Доказательство*. Пусть

$$A = B + [\underbrace{B \ldots B}_{n-1} a] + \ldots$$

– разложение n-арной группы $< A, [\ ] >$ на непересекающиеся правые смежные классы по n-арной подгруппе $< B, [\ ] >$. По предыдущему предложению все смежные классы из указанного разложения являются полуинвариантными n-арными подгруппами в $< A, [\ ] >$. ∎

**9.6.3. Теорема** [126, 130]. Если $< A, [\ ] >$ – конечная идемпотентная n-арная группа порядка $p^k m$, где $(p, m) = 1$, $p$ и $n-1$ – простые, то в $< A, [\ ] >$ существует ровно m полуинвариантных $p$ – силовских n-арных подгрупп $< P_1, [\ ] >, \ldots, < P_m, [\ ] >$, причём

$$A = \bigcup_{i=1}^{m} P_i, \qquad P_i \cap P_j = \varnothing \ (i \neq j).$$

*Доказательство*. По следствию 9.5.18, группа $< A, @ >$ – нильпотентна. Поэтому в $< A, @ >$ – существует инвариантная $p$ – силовская подгруппа $< P_1, @ >$, которая, следовательно, является единственной. По предложению 1.5.4, $\beta_a$ – автоморфизм группы $< A, @ >$, откуда, учитывая единственность $p$ – силовской подгруппы в $< A, @ >$, получаем

$$[aP_1\underbrace{a\ldots a}_{n-2}] = P_1.$$

Тогда по следствию 2.2.7, $<P_1, [\ ]>$ – n-арная подгруппа в $<A, [\ ]>$, которая, очевидно, является p – силовской в $<A, [\ ]>$, а по следствию 2.3.13 и полуинвариантной в $<A, [\ ]>$.

Применяя теперь лемму 9.6.2, получаем разложение

$$A = P_1 + P_2 + \ldots + P_m,$$

n-арной группы $<A, [\ ]>$ в объединение непересекающихся полуинвариантных n-арных подгрупп $<P_i, [\ ]>$, где

$$P_2 = [\underbrace{P_1\ldots P_1}_{n-1} a_2], \ldots, P_m = [\underbrace{P_1\ldots P_1}_{n-1} a_m], a_2, \ldots, a_m \in A.$$

Предположим, что в $<A, [\ ]>$ существует p – силовская n-арная подгруппа $<P, [\ ]>$ – отличная от

$$<P_1, [\ ]>, \ldots, <P_m, [\ ]>.$$

Ясно, что $P_i \cap P \neq \emptyset$ для некоторого $i \in \{1, \ldots, m\}$. Тогда для фиксированного $c \in P_i \cap P$ подгруппы $<P_i, ©>$ и $<P, ©>$ являются различными p – силовскими в группе $<A, ©>$, которая по следствию 1.6.2 изоморфна группе $<A, @>$ с единственной p – силовской подгруппой $<P_1, @>$. Полученное противоречие завершает доказательство. ∎

Так как в конечной нильпотентной группе для любого множества π простых чисел существует единственная π-холловская подгруппа, то, дословно повторяя доказательство теоремы 9.6.3, получим следующее её обобщение.

**9.6.4. Теорема** [126, 130]. Если $<A, [\ ]>$ – конечная идемпотентная n-арная группа порядка km, где $(k, m) = 1$, $n - 1$ – простое, то в $<A, [\ ]>$ существует ровно m полуинвариантных n-арных подгрупп $<P_1, [\ ]>, \ldots, <P_m, [\ ]>$ порядка k, причём

$$A = \bigcup_{i=1}^{m} P_i, \quad P_i \cap P_j = \varnothing \ (i \neq j).$$

**9.6.5. Замечание**. Утверждение о существовании в n-арной группе, удовлетворяющей условиям теоремы 9.6.4, полуинвариантной n-арной подгруппы порядка k может быть получено, как следствие из утверждения теоремы 9.6.3 о существовании силовских n-арных подгрупп и следствия 4.2.1 [4].

**9.6.6. Теорема** [126, 130]. Если $<A, [\ ]>$ – конечная идемпотентная n-арная группа порядка $|A| = p_1^{\alpha_1} \ldots p_m^{\alpha_m}$ ($p_1, \ldots, p_m$, $n - 1$ – простые), то $<A, [\ ]>$ единственным образом разлагается в a-прямое произведение

$$<A, [\ ]> = <A(p_1), [\ ]> \overset{a}{\times} \ldots \overset{a}{\times} <A(p_m), [\ ]> \qquad (1)$$

своих $p_i$ – силовских ($i = 1, \ldots, m$) n-арных подгрупп $<A(P_i), [\ ]>$.

*Доказательство*. По теореме 9.6.3 для любого $p_i$ ($i = 1, \ldots, m$) в $<A, [\ ]>$ существует точно $|A| / p_i^{\alpha_i}$ полуинвариантных $p_i$ – силовских n-арных подгрупп, среди которых только одна, обозначим её через $A(p_i)$, содержит идемпотент a. При этом $<A(p_i), @> - p_i$ – силовская подгруппа группы $<A, @>$, которая по следствию 9.5.18 нильпотентна. Поэтому

$$<A, @> = <A(p_1), @> \times \ldots \times <A(p_m), @>,$$

откуда, учитывая лемму 2.7.21, получаем (1).

Единственность разложения (1) является следствием единственности $p_i$ – силовской n-арной подгруппы в $<A, [\ ]>$, содержащей a. ∎

Согласно утверждению 3) теоремы 5.2 из [4], n-арная группа $<A, [\ ]>$, являющаяся a-прямым произведением своих n-арных подгрупп $<B_1, [\ ]>, \ldots, <B_m, [\ ]>$, изоморфна прямому произведению

$$< B_1, [\ ] > \times \ldots \times < B_m, [\ ] >.$$

Поэтому справедливо

**9.6.7. Следствие** [126, 130]. Конечная идемпотентная n-арная группа $< A, [\ ] >$ порядка

$$|A| = p_1^{\alpha_1} \ldots p_m^{\alpha_m} \ (p_1, \ldots, p_m, n-1 - \text{простые})$$

изоморфна прямому произведению

$$< A(p_1), [\ ] > \times \ldots \times < A(p_m), [\ ] >$$

своих $p_i$-силовских n-арных подгрупп, содержащих один и тот же произвольный элемент $a \in A$.

Дословно повторяя доказательство теоремы 9.6.6, получим следующее её обобщение.

**9.6.8. Теорема** [126, 130]. Если $< A, [\ ] >$ – конечная идемпотентная n-арная группа ($n - 1$ – простое),

$$\pi(A) = \pi_1 \cup \ldots \cup \pi_m, \quad \pi_i \cap \pi_j = \varnothing \ (i \neq j),$$

то $< A, [\ ] >$ единственным образом разлагается в a-прямое произведение

$$< A, [\ ] > = < A(\pi_1), [\ ] > \overset{a}{\times} \ldots \overset{a}{\times} < A(\pi_m), [\ ] >$$

своих $\pi_i$ – холловых $(i = 1, \ldots, m)$ n-арных подгрупп $< A(\pi_i), [\ ] >$.

**9.6.9. Следствие** [126, 130]. Если $< A, [\ ] >$ – конечная идемпотентная n-арная группа ($n - 1$ — простое),

$$\pi(A) = \pi_1 \cup \ldots \cup \pi_m, \quad \pi_i \cap \pi_j = \varnothing \ (i \neq j),$$

то $< A, [\ ] >$ изоморфна прямому произведению

$$< A(\pi_1), [\ ] > \times \ldots \times < A(\pi_m), [\ ] >$$

своих $\pi_i$ – холловых (i = 1, …, m) n-арных подгрупп, содержащих один и тот же произвольный элемент a ∈ A.

**9.6.10. Теорема** [124]. Если n-арная подгруппа единиц < E(A), [ ] > n-арной группы < A, [ ] > является конечной, m – делитель порядка |E(A)|, то < E(A), [ ] > является объединением непересекающихся инвариантных в < A, [ ] > n-арных подгрупп порядка m.

*Доказательство*. Каждый элемент n-арной группы < E(A), [ ] > образует одноэлементную n-арную подгруппу. Поэтому если m = 1, то < E(A), [ ] > является объединением всех своих одноэлементных подгрупп. Если же m = |E(A)|, то доказывать нечего.

Пусть теперь 1 < m < |E(A)|, и зафиксируем элемент e ∈ E(A). Ясно, что < E(A), ⓔ > – подгруппа группы < A, ⓔ >. Из абелевости < E(A), ⓔ > вытекает существование подгруппы < B, ⓔ > порядка m. Так как < B, ⓔ > подгруппа в < E(A), ⓔ >, то < B, [ ] > – n-арная подгруппа в < E(A), [ ] >. Из того, что < E(A), [ ] > лежит в центре < Z(A), [ ] >, следует инвариантность, а значит и полуинвариантность < B, [ ] > в < A, [ ] >. Применяя к < B, [ ] > лемму 9.6.2, получим разложение

$$E(A) = B + [\underbrace{B…B}_{n-1}a] + … , \quad a \in E(A)$$

на непересекающиеся смежные классы, каждый из которых является n-арной подгруппой порядка m. ∎

Как показывает следующий пример, указанное в теореме 9.6.10 разложение n-арной подгруппы единиц в общем случае не является единственным.

**9.6.11. Пример.** Рассмотрим тернарную группу $\mathscr{R} \times \mathscr{R}$ из примера 9.3.2 с n-арной подгруппой единиц

$$E(\mathscr{R} \times \mathscr{R}) = \{(1, 1), (1, a^2), (a^2, 1), (a^2, a^2)\}.$$

Ясно, что

$$E(\mathscr{R} \times \mathscr{R}) = A + B = C + D = F + G,$$

где

$$A = \{(1, 1), (1, a^2)\}, \quad B = \{(a^2, 1), (a^2, a^2)\},$$

$$C = \{(1, 1), (a^2, 1)\}, \quad D = \{(1, a^2), (a^2, a^2)\},$$

$$F = \{(1, 1), (a^2, a^2)\}, \quad G = \{(1, a^2), (a^2, 1)\}.$$

По предложению 9.1.5 любые две единицы тернарной группы образуют тернарную подгруппу. Поэтому тернарная подгруппа единиц $E(\mathscr{R} \times \mathscr{R})$ разлагается тремя различными способами в объединение своих непересекающихся n-арных подгрупп второго порядка.

Следующая теорема обеспечивает единственность разложения тернарной группы единиц.

**9.6.12. Теорема** [124]. Если n-арная подгруппа единиц $< E(A), [\,] >$ n-арной группы $< A, [\,] >$ является конечной порядка $km$, где $(k, m) = 1$, то в $< E(A), [\,] >$ существует ровно $m$ n-арных подгрупп $< P_1, [\,] >, \ldots, < P_m, [\,] >$ порядка $k$. Причем

$$E(A) = \bigcup_{i=1}^{m} P_i, \quad P_i \cap P_j = \varnothing \ (i \neq j).$$

*Доказательство.* По теореме 9.6.10, существует разложение

$$E(A) = \bigcup_{i=1}^{m} P_i, \quad P_i \cap P_j = \varnothing \ (i \neq j)$$

n-арной группы $< E(A), [\,] >$ на непересекающиеся n-арные подгруппы порядка $k$, где $< P_1, \circledcirc > -$ подгруппа порядка $k$ группы $< E(A), \circledcirc >$,

$$P_2 = [\underbrace{P_1 \ldots P_1}_{n-1} a_2], \ldots, P_m = [\underbrace{P_1 \ldots P_1}_{n-1} a_m], a_2, \ldots, a_m \in E(A).$$

Так как $< E(A), \circledcirc > -$ абелева группа, $(k, m) = 1$, то $< P_1, \circledcirc > -$ единственная в $< E(A), \circledcirc >$ подгруппа порядка $k$.

Предположим, что в $<E(A), [\,] >$ существует n-арная подгруппа $< P, [\,] >$ порядка k, отличная от

$$< P_1, [\,] >, \ldots, < P_m, [\,] >.$$

Ясно, что $P_i \cap P \neq \varnothing$ для некоторого $i \in \{1, \ldots, m\}$. Тогда для фиксированного $\varepsilon \in P_i \cap P$ подгруппы $< P_i, \circledcirc >$ и $< P, \circledcirc >$ являются различными подгруппами порядка k группы $< E(A), \circledcirc >$, которая по предложению 9.1.9 изоморфна группе $< E(A), \circledcirc >$ с единственной подгруппой $< P_1, \circledcirc >$ порядка k. Полученное противоречие завершает доказательство. ∎

**9.6.13. Следствие.** Если n-арная подгруппа единиц $< E(A), [\,] >$ n-арной группы $< A, [\,] >$ является конечной порядка $p^k m$, где $(p, m) = 1$, то в $< E(A), [\,] >$ существует ровно m p-силовских n-арных подгрупп $< P_1, [\,] >, \ldots, < P_m, [\,] >$. Причем

$$E(A) = \bigcup_{i=1}^{m} P_i, \quad P_i \cap P_j = \varnothing \ (i \neq j).$$

**9.6.14. Замечание.** Формулировка теоремы 9.6.12 включает в себя тривиальные случаи $k = 1$, $m = 1$.

### §9.7. ПОЛИАДИЧЕСКИЕ ГРУППЫ, ДОПУСКАЮЩИЕ РЕГУЛЯРНЫЙ АВТОМОРФИЗМ

В данном параграфе продолжается начатое в §9.5 изучение n-арных групп, допускающих автоморфизм, называемый в дальнейшем регулярным, с единственным неподвижным элементом. Здесь, в частности, будут получены n-арные аналоги следующих теорем.

**Теорема A** [131, теорема 1.48]. Конечная группа, допускающая регулярный автоморфизм, разрешима.

**Теорема B** [129, теорема V.8.13, Hughes, Kegel, Thompson]. Конечная группа, допускающая расщепляющий автоморфизм простого порядка, нильпотентна.

**Теорема C** [132, Thompson]. Конечная группа, допускающая регулярный автоморфизм простого порядка, нильпотентна.

**Теорема D** [133, Хухро Е.И.]. Разрешимая группа, допускающая регулярный, расщепляющий автоморфизм простого порядка, нильпотентна.

Ясно, что прежде, чем сформулировать n-арные аналоги отмеченных теорем, необходимо определить для n-арных групп понятия, аналогичные понятиям разрешимости и нильпотентности для групп. Среди большого числа существующих n-арных обобщений разрешимости и нильпотентности выберем следующие.

**9.7.1. Определение** [134, Щучкин Н.А.]. n-Арная группа $< A, [\ ] >$ называется *полуразрешимой* (*полунильпотентной*), если некоторая её соответствующая группа разрешима (нильпотентна).

По предложению 1.6.1, для любого $a \in A$ группа $< A, \circledast > $ изоморфна соответствующей группе Поста $< A_o, * >$, а согласно утверждению 2) теоремы 1.4.9, группа $< A_o, * >$ изоморфна любой соответствующей группе n-арной группы $< A, [\ ] >$. Поэтому имеет место

**9.7.2. Определение** [135]. n-Арная группа $< A, [\ ] >$ называется полуразрешимой (полунильпотентной), если для любого $a \in A$ группа $< A, \circledast >$ разрешима (нильпотентна).

Так как, согласно следствию 1.6.2, для любых $a, c \in A$ группы $< A, \circledast >$ и $< A, \copyright >$ изоморфны, то полуразрешимые (полунильпотентные) n-арные группы можно определять следующим образом.

**9.7.3. Определение** [135]. n-Арная группа $< A, [\ ] >$ называется полуразрешимой (полунильпотентной), если для некоторого $a \in A$ группа $< A, @ >$ разрешима (нильпотентна).

Отметим, что в связи со сложившейся в теории n-арных групп терминологией, название – полуразрешимые (полунильпотентные) n-арные группы представляется нам более удачным, чем используемый Щучкиным в оригинале [134] термин – разрешимые (нильпотентные) n-арные группы. Говоря о терминологии, мы имеем ввиду во-первых критерий Поста полуабелевости n-арной группы, согласно которому n-арная группа является полуабелевой тогда и только тогда, когда её соответствующая группа абелева (теорема 2.6.5), и во-вторых критерий полуцикличности n-арной группы, согласно которому n-арная группа $< A, [\ ] >$ является полуциклической тогда и только тогда, когда её соответствующая группа $A_o$ циклическая (теорема 2.5.31).

Веским аргументом в пользу выбранных определений полуразрешимости (полунильпотентности) является их хорошее по сравнению с другими n-арными аналогами разрешимости (нильпотентности) поведение, выражающееся прежде всего в том, что: а) класс всех полуразрешимых n-арных групп содержит все n-арные группы бипримарного порядка и все n-арные группы нечётного порядка; б) всякая n-арная p-группа полунильпотентна. Отметим также, что n-арная подгруппа полуразрешимой (полунильпотентной) n-арной группы сама является полуразрешимой (полунильпотентной). В то же время, например, определённые в [4] нильпотентные n-арные группы могут содержать n-арные подгруппы, не являющиеся нильпотентными. Легко также проверяется, что класс всех полунильпотентных n-арных групп содержит все полуабелевые n-арные группы, а сам входит в класс всех полуразрешимых n-арных групп. Полунильпотентными, согласно следствию 9.5.18, будут и все конечные идемпотентные n-арные группы, где $n - 1$ – простое.

Возможны, как уже отмечалось, и другие, в общем случае неэквивалентные n-арные обобщения разрешимости, как, например, следующие, принадлежащие Русакову [4]: конечная n-арная группа $< A, [\ ] >$ называется разрешимой (полуразрешимой), если она обладает рядом

$$A = A_0 \supseteq A_1 \supseteq A_2 \supseteq \ldots \supseteq A_{k-1} \supseteq A_k \qquad (*)$$

своих n-арных подгрупп таких, что $< A_i, [\ ] >$ – инвариантна (полуинвариантна) в $< A_{i-1}, [\ ] >$ (i = 1, …, k), а числа $|A_{i-1} / A_i|$ и $|A_k|$ – прстые.

Конечные n-арные группы, обладающие субинвариантным (субполуинвариантным) рядом (*) будем в дальнейшем называть *разрешимыми* (*полуразрешимыми*) *по Русакову*.

Ясно, что разрешимая по Русакову n-арная группа будет и полуразрешимой по Русакову.

**9.7.4. Предложение** [135]. Всякая полуразрешимая по Русакову n-арная группа является полуразрешимой.

*Доказательство*. Если $< A, [\ ] >$ – полуразрешимая по Русакову n-арная группа, то существует субполуинвариантный ряд (*). Зафиксируем элемент $a \in A_k$. Тогда по следствию 2.3.13 ряд

$$A = A_0 \supseteq A_1 \supseteq A_2 \supseteq \ldots \supseteq A_{k-1} \supseteq A_k \supseteq \{a\}$$

является субнормальным рядом подгрупп группы $< A, @ >$, порядок каждого фактора которого является простым числом. Следовательно, $< A, @ >$ – разрешимая группа, что влечёт полуразрешимость n-арной группы $< A, [\ ] >$. ∎

Следующий пример показывает, что класс всех полуразрешимых n-арных групп, в общем случае, шире класса всех полуразрешимых по Русакову n-арных групп.

**9.7.5. Пример** [135]. Пусть n – 1 нечётное составное число или n – 1 = $p^\alpha q^\beta$ (p и q – простые). Определим на циклической группе A = < a > n-арную операцию

$$[b_1 \ldots b_n] = b_1 \ldots b_n a.$$

Согласно утверждению 1) леммы 2.5.25 < A, [ ] > – циклическая n-арная группа. Так как для любого a ∈ A группа < A, ⓐ > имеет тот же нечётный или бипримарный порядок, что и группа A, то < A, ⓐ > – разрешима, а значит n-арная группа < A, [ ] > – полуразрешима. В то же время, < A, [ ] > не является полуразрешимой по Русакову, так как, ввиду 2) леммы 2.5.25, в ней нет n-арных подгрупп, отличных от неё самой.

Примером неполуразрешимой тернарной группы может служить тернарная группа < $T_n$, [ ] > всех нечётных подстановок степени n ≥ 5, для которой, как известно (пример 1.4.18), соответствующей группой является знакопеременная группа $A_n$.

**9.7.6. Предложение** [135]. Если < A, [ ] > – конечная n-арная группа, a – фиксированный элемент из A, то справедливы следующие утверждения:

1) если группа < A, ⓐ > допускает регулярный автоморфизм, то < A, [ ] > – полуразрешима;

2) если группа < A, ⓐ > допускает расщепляющий автоморфизм простого порядка, то < A, [ ] > – полунильпотентна;

3) если группа < A, ⓐ > допускает регулярный автоморфизм простого порядка, то < A, [ ] > – полунильпотентна.

*Доказательство*. 1) По теореме A, группа < A, ⓐ > разрешима, что влечёт за собой полуразрешимость n-арной группы < A, [ ] >.

2) Используется теорема B.

3) Используется теорема C. ∎

Предложение 9.7.6 и следствие 9.5.18 позволяют сформулировать

**9.7.7. Предложение** [135]. Конечная идемпотентная n-арная группа, где n − 1 − простое, полунильпотентна.

Согласно лемме 9.5.21, автоморфизм n-арной группы < A, [ ] >, оставляющий неподвижным элемент a, является автоморфизмом группы < A, @ >. Поэтому имеет место

**9.7.8. Лемма** [135]. Автоморфизм n-арной группы < A, [ ] >, оставляющий неподвижным единственный элемент a, является регулярным автоморфизмом группы < A, @ >.

Лемма 9.7.8 и утверждение 1) предложения 9.7.6 позволяют сформулировать следующий n-арный аналог теоремы А.

**9.7.9. Теорема** [135]. Конечная n-арная группа, допускающая регулярный автоморфизм, полуразрешима.

Лемма 9.7.8 и утверждение 3) предложения 9.7.6 позволяют сформулировать следующий n-арный аналог теоремы С.

**9.7.10. Теорема** [135]. Конечная n-арная группа, допускающая регулярный автоморфизм простого порядка, полунильпотентна.

Отметим, что сами утверждения 1), 2) и 3) предложения 9.7.6 являются n-арными аналогами соответственно теорем А, В и С.

**9.7.11. Замечание.** По лемме 9.7.8, всякий регулярный автоморфизм n-арной группы < A, [ ] >, оставляющий неподвижным элемент a, является регулярным автоморфизмом группы < A, @ >, который в свою очередь, в случае конечности A, является расщепляющим [129, теорема V.8.9]. Поэтому утверждение 3) предложения 9.7.6, а также теорема 9.7.10 являются следствиями утверждения 2) предложения 9.7.6.

Если < A, ∘ > − группа, c ∈ A, β − автоморфизм < A, ∘ > такой, что

$$c^\beta = c, \quad b^{\beta^{n-1}} = c \circ b \circ c^{-1}$$

для любого $b \in B$, то согласно обратной теореме Глускина – Хоссу, $< A, [\ ]_{\circ,\beta,c} >$ – n-арная группа с n-арной операцией

$$[a_1 a_2 \ldots a_n]_{\circ,\beta,c} = a_1 \circ a_2^\beta \circ \ldots a_n^{\beta^{n-1}} \circ c.$$

При этом по предложению 2.2.13 операция $\circ$ совпадает с операцией ⊚, где e – единица группы $< A, \circ >$. Это совпадение операций, а также предложение 9.7.6 позволяют сформулировать следующую теорему.

**9.7.12. Теорема** [135]. Пусть $< A, \circ >$ – конечная группа. Тогда справедливы следующие утверждения:

1) если $< A, \circ >$ допускает регулярный автоморфизм, то $< A, [\ ]_{\circ,\beta,c} >$ – полуразрешимая n-арная группа;

2) если $< A, \circ >$ допускает расщепляющий автоморфизм простого порядка, то $< A, [\ ]_{\circ,\beta,c} >$ – полунильпотентная n-арная группа;

3) если $< A, \circ >$ допускает регулярный автоморфизм простого порядка, то $< A, [\ ]_{\circ,\beta,c} >$ – полунильпотентная n-арная группа.

Заметим, что утверждение 3) теоремы 9.7.12 является следствием утверждения 2) этой же теоремы.

Отметим также, что автоморфизмы, которые допускает группа $< A, \circ >$ в утверждениях 1) – 3) теоремы 9.7.12, не обязаны совпадать с автоморфизмом $\beta$.

По теореме 9.5.15, если группа $< A, \circ >$ обладает автоморфизмом $\beta$ таким, что $\beta^{n-1}$ – тождественное преобразование и

$$b \circ b^\beta \circ \ldots b^{\beta^{n-2}} = e$$

для любого b ∈ A, где e – единица группы < A, ∘ >, то < A, ⌊ ⌋ > – идемпотентная n-арная группа с n-арной операцией

$$\lfloor a_1 a_2 \ldots a_n \rfloor = a_1 \circ a_2^\beta \circ \ldots \circ a_{n-1}^{\beta^{n-2}} \circ a_n = [a_1 a_2 \ldots a_n]_{\circ, \beta, c}. \quad (**)$$

причём операции ⊚ и ∘ совпадают. Поэтому, положив в теореме 9.7.12 c = e, и, учитывая, что регулярный автоморфизм конечной группы является расщепляющим [102, теорема V.8.9], получим

**9.7.13. Следствие** [135]. Пусть < A, ∘ > – конечная группа, β – её автоморфизм порядка n – 1. Тогда справедливы следующие утверждения:

1) если β – регулярный, то < A, ⌊ ⌋ > – полуразрешимая идемпотентная n-арная группа;

2) если β – расщепляющий, n - 1 – простое, то < A, ⌊ ⌋ > – полунильпотентная идемпотентная n-арная группа;

3) если β – регулярный, n - 1 – простое, то < A, ⌊ ⌋ > – полунильпотентная идемпотентная n-арная группа.

Следующее следствие дополняет теорему 9.5.22.

**9.7.14. Следствие** [135]. Пусть < A, [ ] > – конечная n-арная группа, допускающая регулярный автоморфизм β порядка n – 1. Тогда на A можно определить n-арную операцию ⌊ ⌋ так, что < A, ⌊ ⌋ > – полуразрешимая идемпотентная n-арная группа. Если же n – 1 – простое, то < A, ⌊ ⌋ > – полунильпотентна.

*Доказательство*. По лемме 9.7.8, β – регулярный автоморфизм группы < A, @ > для некоторого a ∈ A. Обозначим операцию @ символом ∘ и применим следствие 9.7.13. ∎

**9.7.15. Теорема** [135]. Пусть < A, ∘ > – конечная группа порядка |A| = st, где (s,t) = 1; β – её регулярный автоморфизм порядка n – 1, ⌊ ⌋ – операция, определяемая (**). Тогда n-арная группа < A, ⌊ ⌋ > обладает, по крайней мере, одной n-

арной подгруппой порядка s, и любые две n-арные подгруппы порядка s сопряжены в $< A, \lfloor \ \rfloor >$.

***Доказательство***. Согласно теореме A, $< A, \circ >$ – разрешимая группа, в которой по теореме Холла существует, по крайней мере, одна подгруппа порядка s, и любые две подгруппы порядка s сопряжены в $< A, \circ >$. Так как для единицы e группы $< A, \circ >$ операции ⓔ и $\circ$ совпадают, то по предложению 1.6.1, соответствующая группа Поста $< A_o, * >$ для n-арной группы $< A, \lfloor \ \rfloor >$ также обладает, по крайней мере, одной подгруппой порядка s, и любые две подгруппы порядка s сопряжены в $< A_o, * >$.

Согласно замечанию после теоремы V.8.10 [129], $(|A| = |A_o|, n - 1) = 1$, откуда $(t, n - 1) = 1$. Тогда по теореме 3.2.1 [4], n-арная группа $< A, \lfloor \ \rfloor >$ обладает, по крайней мере, одной n-арной подгруппой порядка s, и любые две n-арные подгруппы порядка s сопряжены в $< A, \lfloor \ \rfloor >$. ∎

**9.7.16. Следствие** [135]. Пусть $< A, [ \ ] >$ – конечная n-арная группа порядка $|A| = st$, где $(s, t) = 1$, допускающая регулярный автоморфизм β порядка $n - 1$. Тогда на A можно определить n-арную операцию $\lfloor \ \rfloor$ так, что $< A, \lfloor \ \rfloor >$ – n-арная группа, обладающая, по крайней мере, одной n-арной подгруппой порядка s, и любые две n-арные подгруппы порядка s сопряжены в $< A, \lfloor \ \rfloor >$.

***Доказательство***. По лемме 9.7.8, β – регулярный автоморфизм группы $< A, @ >$ для некоторого $a \in A$. Обозначим операцию @ символом $\circ$ и применим теорему 9.7.15. ∎

**9.7.17. Лемма** [135]. Если для некоторого элемента a n-арной группы $< A, [ \ ] >$ группа $< A, @ >$ – разрешима и допускает регулярный расщепляющий автоморфизм простого порядка, то $< A, [ \ ] >$ полунильпотентна.

***Доказательство***. По теореме D, группа $< A, @ >$ – нильпотентна. Следовательно, n-арная группа $< A, [ \ ] >$ – полунильпотентна. ∎

Леммы 9.7.8 и 9.7.17 позволяют сформулировать следующий n-арный аналог теоремы C.

**9.7.18. Теорема** [135]. Пусть $< A, [\ ] >$ – полуразрешимая n-арная группа с автоморфизмом $\varphi$ простого порядка, оставляющим неподвижным единственный элемент a. Если $\varphi$ – расщепляющий автоморфизм группы $< A, @ >$, то $< A, [\ ] >$ – полунильпотентна.

Следующая теорема является следствием леммы 9.7.17, а также совпадения операций @ и $\circ$ для некоторого $a \in A$.

**9.7.19. Теорема** [135]. Если разрешимая группа $< A, \circ >$ допускает регулярный расщепляющий автоморфизм простого порядка, то $< A, [\ ]_{\circ, \beta, c} >$ – полунильпотентная n-арная группа.

Регулярный расщепляющий автоморфизм и автоморфизм $\beta$ из теоремы 9.7.19 не обязаны совпадать.

Полагая в теореме 9.7.19 c = e, и, учитывая следствие 9.5.16, получим

**9.7.20. Следствие** [135]. Если $< A, \circ >$ – разрешимая группа, $\beta$ – её регулярный расщепляющий автоморфизм простого порядка n – 1, $\lfloor\ \rfloor$ – операция, определяемая (∗∗), то $< A, \lfloor\ \rfloor >$ – полунильпотентная идемпотентная n-арная группа.

**9.7.21. Следствие** [135]. Если для некоторого $a \in A$ автоморфизм

$$\beta_a : x \to [ax\underbrace{a \ldots a}_{n-2}]$$

полуразрешимой идемпотентной n-арной группы $< A, [\ ] >$, где n – 1 – простое, является регулярным, то $< A, [\ ] >$ – полунильпотентна.

*Доказательство*. По теореме 9.5.15 $\beta_a$ – расщепляющий автоморфизм группы $< A, @ >$. Тогда по теореме 9.7.18 $< A, [\ ] >$ – полунильпотентна. ∎

## ДОПОЛНЕНИЯ И КОММЕНТАРИИ

**1.** Идемотентам в n-арных полугруппах посвящена статья В. Дудека [136], в которую он включил некоторые результаты из препринта [124], а именно теорему 9.1.1, утверждающую, что множество всех единиц n-арной группы является её n-арной подгруппой, а также предложение 9.1.5 и следствие 9.1.6.

**2.** В тернарных группах идемпотентность влечет полуабелевость, так как имеет место

**Предложение** [137]. Всякая идемпотентная тернарная группа $< A, [\ ] >$ является полуабелевой.

*Доказательство.* Так как в n-арной группе всякий идемпотент совпадает со своим косым, то, виду предложения 1.2.26,

$$[cba] = [\overline{c}\ \overline{b}\ \overline{a}\ ] = \overline{[abc]} = [abc],$$

то есть

$$[cba] = [abc].$$

для любых a, b, c ∈ A, Следовательно, $< A, [\ ] >$ – полуабелева. ∎

**3.** В [128] приведены различные системы тождеств, определяющих многообразие всех идемпотентных n-арных групп.

**4.** Ясно, что все критерии полуразрешимости n-арной группы, поученные в §9.7 верны по модую классификации конечных простых групп, так как при получении этих критериев использовалась теорема А, являющаяся следствием указанной классификации.

# ЛИТЕРАТУРА

**1. Dörnte, W.** Untersuchungen über einen verallgemeinerten Gruppenbegrieff / W. Dörnte // Math. Z. – 1928. – Bd. 29. – S. 1 – 19.

**2. Prüfer, H.** Theorie der abelshen Gruppen. I. Grundeigenschaften / H. Prüfer // Math. Z. – 1924. – Bd. 20. – S. 165 – 187.

**3. Post, E.L.** Polyadic groups / E.L. Post // Trans. Amer. Math. Soc. – 1940. – Vol. 48, №2. – P.208 – 350.

**4. Русаков, С.А.** Алгебраические n-арные системы / С.А. Русаков. – Мн.: Навука і тэхніка, 1992. – 245 с.

**5. Русаков, С.А.** Некоторые приложения теории n-арных групп / С.А. Русаков. – Минск: Беларуская навука, 1998. – 167 с.

**6. Сушкевич, А.К.** Теория обобщенных групп / А.К. Сушкевич. – Харьков; Киев, 1937. – 176 с.

**7. Курош, А.Г.** Общая алгебра: Лекции 1969/70 учебного года / А.Г. Курош. – М.: Наука, 1974. – 160 с.

**8. Bruck, R.H.** A survey of binary systems / R.H.Bruck. – Berlin; Heldelberg; New York: Springer-Verlad, 1966. – 185 p.

**9. Бурбаки, Н.** Алгебра. Алгебраические структуры, линейная и полилинейная алгебра / Н. Бурбаки. – М.: Физматгиз, 1962.

**10. Артамонов, В.А.** Универсальные алгебры / В.А. Артамонов // Итоги науки и техники. – Сер. Алгебра. Топология. Геометрия. – 1976. – С. 191 – 248.

**11. Glazek, K.** Bibliographi of n-groups (poliadic groups) and same group like n-ary systems / K. Glazek // Proc. of the sympos. n-ary structures. – Skopje, 1982. – P. 259 – 289.

**12. Гальмак, А. М.** Конгруэнции полиадических групп / А.М. Гальмак. – Минск: Беларуская навука, 1999. – 182 с.

**13. Глускин, Л.М.** Позиционные оперативы / Л.М. Глускин // Мат.сборник. – 1965. – Т.68(110), №3. – С.444 – 472.

**14. Hosszu, M.** On the explicit form of n-group operations / M. Hosszu // Publ. Math. – 1963. – V.10, №1 – 4. – P.88 – 92.

**15. Гальмак, А.М.** О приводимости n-арных групп / А.М. Гальмак // Вопросы алгебры. – 1996. – Вып. 10. – С. 164 – 169.

**16. Соколов, Е.И.** О теореме Глускина-Хоссу для n-групп Дёрнте / Е.И. Соколов // Мат. исследования. – Вып.39. – С.187 – 189.

**17. Артамонов, В.А.** Свободные n-арные группы / В.А. Артамонов // Мат. заметки. – 1970. – Т.8, №4. – С. 499 – 507.

**18. Артамонов, В.А.** О шрайеровых многообразиях n-групп и n-полугрупп / В.А. Артамонов // Труды семинара им. И.Г. Петровского. – 1979. – Вып.5. – С. 193 – 202.

**19. Гальмак, А.М.** Трансляции n-арных групп / А.М. Гальмак // Докл. АН БССР. – 1986. – Т.30, №8. – С. 677 – 680.

**20. Гальмак, А.М.** О приводимости n-арных групп / А.М. Гальмак // Препринты ИМ АН БССР. – 1976. – 6(242). – 36 с.

**21. Гальмак, А.М.** Приводимость полиадических групп / А.М. Гальмак // Докл. АН БССР. – 1985. – Т.29, №10. – С. 874 – 877.

**22. Dudek, W.A.** On a generalisation of Hosszu theorem / W.A. Dudek, J. Michalski // Denconstratio Math, 1982. – Vol. 15, №3. – P. 783 – 805.

**23. Дириенко, И.И.** К теореме Глускина-Хоссу об n-группах / И.И. Дириенко, О.В. Колесников. – Харьков, 1980. – 10 с. – Деп. в ВИНИТИ №374 – 80.

**24. Гальмак, А.М.** Теоремы Поста и Глускина-Хоссу / А.М. Гальмак. – Гомель, 1997. – 85 с.

**25. Гальмак, А.М.** Тернарные группы отражений / А.М. Гальмак // Междунар. Мат. Конф. – Тез. докл. – Гомель, 1994. – С. 33.

**26. Гальмак, А. М.** Тернарные группы отражений / А.М. Гальмак, Г.Н. Воробьев. – Минск: Беларуская навука, 1998. – 128с.

**27. Masat Fransis, E.** A useful characterisation of a normal subgroups / E. Masat Fransis // Math. Mag. – 1979. – Vol. 52, № 3. – P. 171 – 173.

**28. Гальмак, А.М.** Инвариантные подгруппы n-арных групп и их обобщения / А.М. Гальмак // Вопросы алгебры. – Мн.: Университетское, 1990. – Вып. 5. – С. 91 – 94.

**29. Воробьев, Г.Н.** О сопряженности n-арных подгрупп / Г.Н. Воробьев // Весці Акадэміі навук Беларусі. Сер. фіз.-мат. навук. – 1996. – №1. – С. 121.

**30. Воробьев, Г.Н.** О полусопряженности n-арных подгрупп / Г.Н. Воробьев // Вопросы алгебры. – Гомель, 1997. – Вып. 10. – С.157 – 163.

**31. Серпинский, В.** 250 задач по элементарной теории чисел / В. Серпинский. – М.: Просвещение, 1968. – 160 с.

**32. Гальмак, А.М.** Абелевы n-арные группы и их обобщения / А.М. Гальмак // Вопросы алгебры. – Минск: Университетское, 1987. – Вып. 3. – С. 86 – 93.

**33. Dudek, W.A.** Remarcs on n-groups / W.A. Dudek // Demonstratio Math. – Vol. 13, №1. – 1980. – P.165 – 181.

**34. Колесников, О.В.** Разложение n-групп / О.В. Колесников // Мат. исслед. – Вып. 51. – Квазигруппы и лупы. Кишинёв: Штиинца, 1979. – С. 88 – 92.

**35. Плоткин, Б.И.** Группы автоморфизмов алгебраических систем / Б.И. Плоткин. – Минск: Наука, 1966. – 603 с.

**36. Glazek, K.** Abelian n-groups / K. Glazek, B. Gleichgewicht // Proc. Congr. Math. Soc. J. Bolyai. – Esztergom. – 1977. – P. 321 – 329.

**37. Гальмак, А.М.** Полуабелевые n-арные группы с идемпотентами / А.М. Гальмак // Веснік ВДУ ім. П.М. Машэрава. – 1999. – № 2(12). – С. 56 – 60.

**38. Воробьев, Г.Н.** Сопряженные n-арные подгруппы и их обобщения / Г.Н. Воробьев // Веснік ВДУ ім. П.М. Машэрава. – 1997. – № 2(4). – С.59 – 64.

**39. Гаврилов, В.В.** О полуциклических n-арных группах / В.В. Гаврилов // Конф. математиков Беларуси. – Тез. докл. Гродно, 1992. – С. 15.

**40. Дудек, В.А.** m-Полуабелевые n-арные группы / В.А. Дудек // Изв. АН ССР Молдова. – Математика. – 1990. – №2. – С. 66 – 70.

**41. Dudek, W.A.** On the class of weakly semiabelian polyadic groups / W.A. Dudek // Discrete Math. – Appl. – Vol. 6, №5. – P. 427 – 433.

**42. Мальцев, А.И.** К общей теории алгебраических систем / А.И. Мальцев // Мат. сб. – 1954. – Т.35, №1. – С.3 – 20.

**43. Monk, J.D.** On the general theory of m-groups / J.D. Monk, F.M. Sioson // Fund. Math. – 1971. – №72. – P. 233 – 244.

**44. Кулаженко, Ю.И.** Критерии полуабелевости n-арной группы / Кулаженко Ю.И. // Веснік ВДУ ім. П.М. Машэрава. 1997. №3(5). С. 61 – 64.

**45. Sioson, F.M.** On Free Abelian n-Groups I / F.M. Sioson // Proc. Japan Acad. – 1967. – Vol. 43. – P. 876 – 879.

**46. Sioson, F.M.** On Free Abelian n-Groups II / F.M. Sioson // Proc. Japan Acad. – 1967. – Vol. 43. – P. 880 – 883.

**47. Sioson, F.M.** On Free Abelian n-Groups III / F.M. Sioson // Proc. Japan Acad. – 1967. – Vol. 43. – P. 884 – 888.

**48. Гальмак, А.М.** Об определении n-арной группы / А.М. Гальмак // Междунар. конф. по алгебре. – Тез. докл. – Новосибирск, 1991. – С. 30.

**49. Тютин, В.И.** К аксиоматике n-арных групп / В.И. Тютин // Докл. АН БССР. – 1985. – Т.29, №8. – С. 691 – 693.

**50. Гальмак, А.М.** О некоторых новых определениях n-арной группы / А.М. Гальмак // Третья междунар. конф. по алгебре. – Тез. докл. – Красноярск, 1993. – С. 33.

**51. Гальмак, А.М.** Новые определения n-арной группы / А.М. Гальмак // Конф. математиков Беларуси. – Тез. докл. Гродно, 1992. – С. 17.

**52. Celakoski, N.** On some axiom systems for n-groups / N. Celakoski // Мат. Бил. Сојуз. друшт. мат. СРМ. – 1997. – Кн. 1. – P. 5 – 14.

**53. Dudek, W.** A note on the axioms of n-groups / W. Dudek, K. Glazek, B. Gleichgewicht // Colloq Math. Soc. J. Bolyai. – 1977. – Vol. 29. – P. 195 – 202.

**54. Русаков, С.А.** К определению n-арной группы / С.А. Русаков // Докл. АН БССР. – 1979. – Т.23, №11. – С. 965 – 967.

**55. Ulshofer, K.** Schlichtere Gruppenaxiome / K. Ulshofer // Praxis Math. 1(1972). – S. 1 – 2.

**56. Galmak, A.M.** Remarks on polyadic groups / A.M. Galmak // Quasigroups and related Systems. – 2000. – Vol. 7. – P. 67 – 70.

**57. Тютин, В.И.** Об условиях, при которых n-арная полугруппа является n-арной группой / В.И. Тютин // Арифметическое и подгрупповое строение конечных групп. – Минск: Навука і тэхніка, 1986. – С.161 – 170.

**58. Гальмак, А.М.** Определения n-арной группы / А.М. Гальмак // Препринт ГГУ им. Ф. Скорины. – 1994. – № 16. – 43 с.

**59. Tvermoes, H.** Om en Generalisation af Gruppenbegrebet / H. Tvermoes. – Kobenhavn. – 1952. – 107 s.

**60. Tvermoes, H.** Über eine verellgemeinerung des Gruppenbegriffs / H. Tvermoes // Math. Scand. – 1953. – Bd. 1. – S. 18 – 30.

**61. Robinson, D.W.** n-Groups with identity elements / D.W. Robinson // Math. Mag. – 1958. – Vol. 31, №.5. – P. 255 – 258.

**62. Слипенко, А.К.** Регулярные оперативы и идеальные эквивалентности / А.К. Слипенко // Докл. АН УССР. – Сер.А. – 1977. – №3. – С. 218 – 221.

**63. Ušan, J.** n-Groups in the light of the neutral operations / J. Ušan // Mathematika Moravica. – 2003. – Special Vol. – 162 p.

**64. Gleichgewicht B.** Remarks on n-groups as abstract algebras / B. Gleichgewicht, K. Glazek // Collq Math. – 1967. – Vol. 17, №2. – P. 209 – 219.

**65. Sokhatski, F.M.** Invertible elements in associates and semigroups 1 / F.M. Sokhatski // Quasigroups and related Systems. – 1998. – Vol. 5. – P. 53 – 68.

**66. Sokhatski, F.M.** Invertible elements in associates and semigroups 2 / F.M. Sokhatski, O.W. Yurevych // Quasigroups and related Systems. – 1999. – Vol. 6. – P. 61 – 70.

**67. Сохоцкий, Ф.Н.** Об ассоциативности многоместных операций / Ф.Н. Сохоцкий // Дискретная математика. – 1992. – №4. – С. 66 – 84.

**68. Юревич, О.В.** Критерії оборотності элементів в асоціатах / О.В. Юревич // Укр. мат. журнал. – 2001. – Т. 53, №11. – С. 1556 – 1563.

**69. Юревич, О.В.** Про скрещену ізотопію поліагруп / О.В. Юревич // Труды института прикладной математики и механики НАН Украины. – 2005. – Вып. 11. – С. 34 – 39.

**70. Ušan, J.** On NP-polyagroups / J. Ušan, R. Galić // Math. Comm. – 2001. – Vol. 6, №2. – P. 153 – 159.

**71. Белоусов, В.Д.** n-Арные квазигруппы / В.Д. Белоусов. – Кишинев: Штиинца, 1972. – 228 с.

**72. Соколов, Е.И.** О приводимости (i, j)-ассоциативных n-квазигрупп / Е.И. Соколов // Изв. АН МССР. – 1968. – №3. – С. 10 – 18.

**73. Čupona, G.** Vector valued semigroups / G. Čupona // Semigroup Forum. – 1983. – № 26. – P. 65 – 74.

**74. Čupona, G.** On topological n-groups / G. Čupona // Билтен. на Друшт. на Мат. и физ. од СРМ. – 1971. – Кн. 22. – P. 5 – 10.

**75.Crombez, G.** On topological n-groups / G. Crombez, G. Six // Abh. Math. Sem. Univ. – Gamburg. – 1974. – №41. – P.115 – 124.

**76. Žižović, M.** Topological analogy of Hosszu-Gluskin's Theorem / M. Žižović // Mat. vesnik. – 1976. – №13. – P. 233 – 235.

**77. Enders, N.** On topological n-groups and their corresponding groups / N. Enders // Discussiones Mathematicae. Algebra and Stohastic Methods. – 1995. – №15. – P. 163 – 169.

**78. Ušan, J.** On topological n-groups / J. Ušan // Mathematika Moravica. – 1998. – №2. – P. 149 – 159.

**79. Мухин, В.В.** О вложении n-арных абелевых топологических полугрупп в n-арные топологические группы / В.В. Мухин, Буржуф Хамза // Вопросы алгебры. Гомель. – 1996. – Вып. 9. – С.57 – 60.

**80. Crombez, G.** On partially ordered n-groups / G. Crombez // Abh. Math. Sem. Univ. – Gamburg. – 1972. – №38. – P.141 – 146.

**81. Ušan, J., Žižović M.** On ordered n-groups / J. Ušan, M. Žižović // Quasigroups and related Systems. – 1997. – Vol. 4. – P. 77 – 87.

**82. Гальмак, А.М.** n-Арные перестановки / А.М. Гальмак // Кн. Некоторые вопросы алгебры и прикладной математики. – Гомель, 2002. – С. 45 – 49.

**83. Гальмак, А.М.** n-Арные морфизмы алгебраических систем / А.М. Гальмак // Международная матем. конференция. – Тез. докл. – Минск, 1993. – С.11 – 12.

**84. Гальмак, А.М.** Обобщённые морфизмы алгебраических систем / А.М. Гальмак // Вопросы алгебры. – Гомель. – 1998. – Вып. 12. – С.36 – 46.

**85. Плоткин, Б.И.** Группа автоморфизмов алгебраических систем / Б.И. Плоткин – М.: Наука, 1966. – 604 с.

**86. Гальмак, А.М.** Полиадические аналоги теорем Кэли и Биркгофа / А.М. Гальмак // Известия ВУЗов. – Математика. – 2001. – №2 (465). – С. 13 – 18.

**87. Birkhof, G.** On groups of automorphisms / G. Birkhof // Revista Union Mat. – Argentina. – 1946. – Vol. 11, №4. – P.155 – 157.

**88. Гальмак, А.М.** n-Арные аналоги теоремы Биркгофа / А.М. Гальмак // Материалы междунар. конф., посвященная памяти академика С. А. Чунихина. – Ч. I. – Гомель, 1995. – С. 49.

**89. Crombez, G.** On (n, m)-rings / G. Crombez // Abh. Math. Sem. Univ. – Hamburg. – 1972. – Vol. 37. – P.180 – 199.

**90. Скорняков, Л.А.** Элементы общей алгебры / Л.А. Скорняков – М.: Наука, 1983. – 273с.

**91. Galmak, A.M.** Generalized morphisms of abelian m-ary groups / A.M. Galmak // Discussiones Mathematicae. General Algebra and Applications. – 2001 – №21. – P. 47 – 55.

**92. Čupona, G.** On [m, n]-rings / G. Čupona // Bull. Soc. math. phys. Mased. – 1965. – Vol. 16. – P. 5 – 10.

**93. Crombez, G.** On (m, n)-quotient rings / G. Crombez, J. Timm // Abh. Math. Semin. Univ. – Hamburg. – 1972. – Vol. 37. – P. 200 – 203.

**94. Dudek, W.** On the divisibility theory in (m, n)-rings / W. Dudek // Dem. Math. – 1981. – Vol. XIV, №1. – P. 19 – 32.

**95. Никитин, А.Н.** Радикал Джекобсона артиновых (2, n)-кольца / А.Н. Никитин // Вестн. Моск. ун-та: Сер. 1. – Математика. Механика. – 1984. – №4. – С. 18 – 22.

**96. Никитин, А.Н.** Полупростые артиновы (2, n)-кольца / А.Н. Никитин // Вестн. Моск. ун-та: Сер. 1. – Математика. Механика. – 1984. – №6. – С. 3 – 7.

**97. Артамонов, В.А.** Артиновы (2, n)-кольца / В.А. Артамонов, А.Н. Никитин // Алгебр. сист. – Волгоград, 1989. – С. 13 – 23.

**98. Celakovski, N.** On (F, G)-rings / Celakovski N.// Год. сбор. Мат. фак. ун-т. – Skopje. – 1977. – Vol. 28. P. 5 – 15.

**99. Кондратова, А.Д.** Об одном свойстве (m, n)-кольца / А.Д. Кондратова // Вопросы алгебры и прикладной математики. – Сб. науч. тр. под ред. С.А. Русакова. – Гомель, 1995. – С. 97 – 101.

**100. Кондратова, А.Д.** Об изоморфном вложении (M, N)-колец / А.Д. Кондратова // Веснік ВДУ ім. П.М. Машэрава. – 1998. – № 1(7). – С. 75 – 78.

**101. Кондратова, А.Д.** Об (M, N)-кольцах, обладающих (M, N)-кольцами частных / А.Д. Кондратова // Вопросы алгебры. – Гомель. 1999. – Вып. 14. – С. 173 – 180.

**102. Суворова, А.Д.** Обобщенные (M, N)-полукольца / А.Д. Суворова // Некоторые вопросы алгебры и прикладной математики. – Сб. науч. тр. – Гомель, 2002. – С. 78 – 100.

**103. Кравченко, Ю.В.** Полупрямые произведения полиадических мультиколец / Ю.В. Кравченко // Вопросы алгебры. – Гомель. 1996. – Вып. 10. – С. 190 – 198.

**104. Кравченко, Ю.В.** О формациях полиадических мультиколец / Ю.В. Кравченко // Вестник БГУ. Сер. физ.-мат. наук. – 1998. – №1. – С. 53 – 57.

**105. Кравченко, Ю.В.** О теореме Жордана-Гельдера для полиадических мультиколец / Ю.В. Кравченко, С.П. Новиков // Некоторые вопросы алгебры и прикладной математики. – Сб. науч. тр. Гомель. – 2002. – С. 101 – 110.

**106. Гальмак, А.М.** n-Арные аналоги нормальных подгрупп / А.М. Гальмак // Веснік МДУ ім. А.А. Кулешова. – 2004. – № 2 – 3(15). – С. 153 – 159.

**107. Воробьев, Г.Н.** О сопряженности полусопряженности n-арных подгрупп в n-арной группе / Г.Н. Воробьев // Известия ГГУ им. Ф.Скорины. – 2006. – № 5. – С. 26 – 32.

**108. Гальмак, А.М.** Представление полиадической группы подстановками на смежных классах / А.М. Гальмак // XVIII Всесоюзная алгебраическая конф. – Тез. докл. – Ч. I. – Кишинев, 1985. – С. 107.

**109. Каргаполов, М.И.** Основы теории групп / М.И. Каргаполов, Ю.И. Мерзляков. – М.: Наука, 1982. – 288 с.

**110. Гальмак, А.М.** О разложениях в обертывающей группе Поста / А.М. Гальмак // Веснік МДУ ім. А.А. Кулешова. – 2006. – № 2 – 3(24). – С. 182 – 189.

**111. Гальмак, А.М.** Разложения обертывающей группы Поста / А.М. Гальмак // Вестник Полоцкого государственного университета. – Серия С. – 2005. – №10. – С. 14 – 18.

**112. Гальмак, А.М.** О теореме Шура для n-арных групп / А.М. Гальмак // Укр. мат. журнал. – 2006. – Т. 58, – № 5. – С. 730 – 741.

**113. Холл, М.** Теория групп / М. Холл. – М.:Изд-во иностранной литературы, 1962. – 468 с.

**114. Гальмак, А.М.** m-Полунормализаторы в n-арной группе / А.М. Гальмак // Известия ГГУ им Ф. Скорины. – 2006. – № 5. – С. 33 – 38.

**115. Galmak, A.M.** Some n-ary analogs of the notion of normalizer of an n-ary subgroup in a group / A.M. Galmak // Bulutinul akademiei de stiente a republicei Moldova. – Matematica. – 2005. – № 3 (49). – P. 63 – 70.

**116. Гальмак, А.М.** К определению инвариантных подмножеств в n-арной группе / А.М. Гальмак // Веснік МДУ ім. А.А. Кулешова. – 1999. – № 2 – 3(3). – С. 88 – 90.

**117. Celakovski, N.** A note on invariant subgroups of n-groups / N. Celakovski, S. Ilic. // Proc. Conf. "Algebra and Logic". – Zagreb, – 1984. – P. 21 – 28.

**118. Тютин, В.И.** n-Арные группы с f-центральными рядами / В.И. Тютин // Вопросы алгебры. – Гомель. – 1987. – Вып. 3. – С. 97 – 116.

**119. Гальмак, А.М.** n-Арные аналоги центра группы / А.М. Гальмак // Препринт института прикладной оптики НАН Беларуси. – 2004. – № 16. – 35с.

**120. Гальмак, А.М.** n-Арные аналоги центра группы / А.М. Гальмак // Веснік МДУ ім. А.А. Куляшова – 2005. – №1 (20). – С. 90 – 97.

**121. Гальмак, А.М.** Полиадические аналоги центра группы / А.М. Гальмак // Весник Полоцкого государственного университета. – Серия С. – 2004. – №11. – С. 24 – 28.

**122. Русаков, С.А.** К теории нильпотентных n-арных групп / С.А. Русаков // Конечные группы. – Минск: Навука і тэхніка, 1978. – С.104 – 130.

**123. Ilic, S.** On nilpotent n-groups / S. Ilic. // Математички весник, 1986. – №38. – Р. 291 – 298.

**124. Гальмак, А.М.** n-Арная подгруппа единиц / А.М. Гальмак // Препринт ГГУ им. Ф.Скорины. – 1998. – № 77. – 23 с.

**125. Гальмак, А.М.** n-Арная подгруппа единиц / А.М. Гальмак // Весці НАН РБ. – 2003. – №2. – С.25 – 30.

**126. Гальмак, А. М.** Идемпотенты в n-арных группах / А.М. Гальмак // Препринт ГГУ им. Ф.Скорины. – 1998. – №81. – 28с.

**127. Гальмак, А.М.** Идемпотентные n-арные группы / А.М. Гальмак // Весці НАН РБ. – 2000. – №2. – С.42 – 45.

**128. Адян, С.И.** Проблема Бернсайда и тождества в группах / С.И. Адян. – М.: Наука, 1975. – 335 с.

**129. Huppert, B.** Endliche Gruppen I / B. Huppert. – Berlin; Heidelberg; New York: Springer, 1967. – 793 p.

**130. Гальмак, А.М.** Силовское строение идемпотентных n-арных групп / А.М. Гальмак // Укр. мат. журнал. – 2001. – №11. – С.1488 – 1494.

**131. Горенстейн, Д.** Конечные простые группы / Д. Горенстейн. – М.: Мир, 1985. – 350с.

**132. Thompson, J.** Finite groups with fixed – point free automorphismus of finite order / J. Thompson // Proc. Nat. Acad. Sci USA. – 1959. – №45. – P.578 – 581.

**133. Хухро, Е.И.** Разрешимая группа, допускающая регулярный расщепляющий автоморфизм простого порядка, нильпотентна / Е.И. Хухро // Алгебра и логика. – 1978. – Т.17, №5. – С.611 – 618.

**134. Щучкин, Н.А.** Разрешимые и нильпотентные n-группы / Н.А. Щучкин // Алгебраические системы. – Волгоград, 1989. – С.133 – 139.

**135. Гальмак, А.М.** Полиадические группы, допускающие регулярный автоморфизм / А.М. Гальмак // Известия ГГУ им. Ф.Скорины. – 2002. – № 5. Вопросы алгебры – 18. – С.104 – 111.

**136. Dudek, W.** Idempotents in n-ary semigroups / W. Dudek // Southeast Asian Bulletin of Mathematics. – 2001. – №25. – P. 97 – 104.

**137. Dudek, W.** Autodistributive n-groups / W. Dudek // Commentationes Math. – Annales Soc. Math. Polonae, Prace Matematyczne. – 1983. – №23. – P. 1 – 11.

**138. Dudek, W.** Varieties of polyadic groups / W. Dudek // Filomat. – 1995. – №9. – P. 657 – 674.

# ПРЕДМЕТНЫЙ УКАЗАТЕЛЬ

# УСЛОВНЫЕ ОБОЗНАЧЕНИЯ

$$a_m^k = \begin{cases} a_m a_{m+1} \ldots a_k, & \text{если } m \leq k, \\ \varnothing, & \text{если } m > k; \end{cases}$$

$$\overset{k}{a} = \begin{cases} \underbrace{a \ldots a}_{k}, & k > 0, \\ \varnothing, & k = 0. \end{cases}$$

$$a^{[s]} = \begin{cases} a, & s = 0, \\ [\overset{s(n-1)+1}{a}], & s > 0, \\ [\overset{-2s}{\overline{a}} \overset{-s(n-3)+1}{a}], & s < 0. \end{cases}$$

$$\overset{m}{B} = \begin{cases} \underbrace{B \ldots B}_{m}, & \text{если } m \geq 1; \\ \varnothing, & \text{если } m \leq 0; \end{cases}$$

$\overline{a}$ – косой элемент для элемента a;

$\alpha^{-1}$ – обратная последовательность для последовательности $\alpha$;

$l(\alpha)$ – длина последовательности $\alpha$;

$F_A$ – свободная полугруппа над алфавитом A;

$\sim$ или $\theta$ – отношение эквивалентности Поста на $F_A$;

$\mathscr{A} = F_A / \theta$;

$A^{(i)} = \{\theta(\alpha a) \mid a \in A\} = \{\theta(a\alpha) \mid a \in A\}$;

$A_n$ – знакопеременная группа степени n;

$< B_n, [\ ] >$ – тернарная группа всех отражений правильного n-угольника;

$C_n$ – циклическая группа порядка n;

$D_n$ – диэдральная группа порядка 2n;

$S_n$ – симметрическая группа степени n;

$[\ ]$ – n-арная операция;

$< A, [\ ] >$ – n-арная группа;

$x@y = [xa^{-1}y]$, где a – элемент n-арной группы $< A, [\ ] >$;

$B_a = \{[b_1 \ldots b_{n-1}a] \mid b_1, \ldots, b_{n-1} \in B\}$;

$_aB = \{[ab_1 \ldots b_{n-1}] \mid b_1, \ldots, b_{n-1} \in B\}$;

$B^{(i)}(A) = \{\theta_A(\alpha) \in A^{(i)} \mid \exists b_1, \ldots, b_i \in B, \alpha \theta_A b_1 \ldots b_i\}$, $i = 1, \ldots, n-1$;

$B_o(A) = B^{(n-1)}(A) = \{\theta_A(\alpha) \in A_o \mid \exists b_1, \ldots, b_{n-1} \in B, \alpha \theta_A b_1 \ldots b_{n-1}\}$;

$B^*(A) = \{\theta_A(\alpha) \in A^* \mid \exists b_1, \ldots, b_i \in B \ (i \geq 1), \alpha \theta_A b_1 \ldots b_i\}$;

$|A : B|$ – индекс n-арной подгруппы $< B, [\ ] >$ в n-арной группе $< A, [\ ] >$;

$< B, [\ ] > \vee < C, [\ ] >$ – n-арная подгруппа, порождённая n-арными подгруппами $< B, [\ ] >$ и $< C, [\ ] >$;

$L(A, [\ ])$ – множество всех n-арных подгрупп n-арной группы $< A, [\ ] >$;

$< B_1, [\ ] > \overset{a}{\times} \ldots \overset{a}{\times} < B_m, [\ ] >$ – a-прямое произведение n-арных подгрупп $< B_1, [\ ] >, \ldots, < B_m, [\ ] >$;

$A(\pi)$ – $\pi$-холлова n-арная подгруппа в $< A, [\ ] >$;

$A(p)$ – p-силовская n-арная подгруппа в $< A, [\ ] >$.

$N_A(B)$ – нормализатор B в $< A, [\ ] >$;

$HN_A(B)$ – полунормализатор B в $< A, [\ ] >$;

$N_A(B, m)$ – m-полунормализатор B в $< A, [\ ] >$.

$Z(A)$ – центр $< A, [\ ] >$;

$HZ(A)$ – полуцентр $< A, [\ ] >$;

$Z(A, m)$ – m-полуцентр $< A, [\ ] >$;

$HDZ(A)$ – слабый полуцентр $< A, [\ ] >$;

$DZ(A, m)$ – слабый m-полуцентр $< A, [\ ] >$;

$HTZ(A)$ – полуцентр типа T $< A, [\ ] >$;

$TZ(A, m)$ – m-полуцентр типа T $< A, [\ ] >$;

$Z(A, \Sigma)$ – $\Sigma$-полуцентр $< A, [\ ] >$;

$Z(A, \Sigma, m)$ – $(\Sigma, m)$-полуцентр $< A, [\ ] >$;

$C_A(B)$ – централизатор B в $< A, [\ ] >$;

$HC_A(B)$ – полуцентрализатор B в $< A, [\ ] >$;

$C_A(B, m)$ – m-полуцентрализатор B в $< A, [\ ] >$;

$HDC_A(B)$ – слабый полуцентрализатор B в $< A, [\ ] >$;

$DC_A(B, m)$ – слабый m-полуцентрализатор B в $< A, [\ ] >$;

$HTC_A(B)$ – полуцентрализатор типа T B в $< A, [\ ] >$;

$TC_A(B, m)$ – m-полуцентрализатор типа T B в $< A, [\ ] >$;

$HC_A(B, \Sigma)$ – $\Sigma$-полуцентрализатор B в $< A, [\ ] >$;

$C_A(B, \Sigma, m)$ – $(\Sigma, m)$-полуцентрализатор B в $< A, [\ ] >$;

$Aut(n, A)$ – множество всех n-арных автоморфизмов алгебры A;

$Aut(A_1, \ldots, A_{n-1})$ – множество всех n-арных автоморфизмов последовательности $\{A_1, \ldots, A_{n-1}, A_1\}$;

$End(n, A)$ – множество всех n-арных эндоморфизмов алгебры A;

$End(A_1, \ldots, A_{n-1})$ – множество всех n-арных эндоморфизмов последовательности $\{A_1, \ldots, A_{n-1}, A_1\}$;

$E(A)$ – n-арная подгруппа единиц n-арной группы $< A, [\ ] >$;

I(A) – множество всех идемпотентов n-арной группы < A, [ ] >;

L(A) – множество всех левых n-арных сдвигов n-арной группы < A, [ ] >;

R(A) – множество всех правых n-арных сдвигов n-арной группы < A, [ ] >;

$S_{A_1, \ldots, A_{n-1}}(\sigma)$ – множество всех n-арных подстановок последовательности $\{A_1, \ldots, A_{n-1}\}$, определяемых подстановкой $\sigma \in S_{n-1}$;

$S_{A_1, \ldots, A_{n-1}}(T)$ – множество всех n-арных подстановок последовательности $\{A_1, \ldots, A_{n-1}\}$, определяемых подстановками $\sigma \in T \subseteq S_{n-1}$;

$S_{A_1, \ldots, A_{n-1}}(S_{n-1})$ – множество всех n-арных подстановок последовательности $\{A_1, \ldots, A_{n-1}\}$.

# СОДЕРЖАНИЕ